\documentclass{memo-l}

\usepackage{amsmath}
\usepackage{amssymb}
\usepackage{verbatim}
\usepackage{epic}
\usepackage{eepic}
\usepackage{pstricks}

\newtheorem{theorem}{Theorem}[section]
\newtheorem{lem}[theorem]{Lemma}
\newtheorem{prop}[theorem]{Proposition}
\newtheorem{cor}[theorem]{Corollary}
\newtheorem{conj}[theorem]{Conjecture}
\newtheorem{claim}{Claim}

\theoremstyle{definition}
\newtheorem{definition}[theorem]{Definition}

\newtheorem{problem}[theorem]{Problem}

\theoremstyle{remark}
\newtheorem{remark}[theorem]{Remark}
\newtheorem{question}[theorem]{Question}

\numberwithin{section}{chapter}
\numberwithin{equation}{chapter}
\numberwithin{figure}{chapter}

\makeindex

\newcommand{\N}{\mathbb{N}}
\newcommand{\Q}{\mathbb{Q}}
\newcommand{\ZZ}{\mathbb{Z}}

\newcommand{\F}{\mathbb{F}}

\newcommand{\Z}{\mathrm{Z}}
\newcommand{\B}{\mathbb{B}}
\newcommand{\pow}{\mathcal{P}}
\newcommand{\dom}{\mathrm{dom}}

\newcommand{\TC}{\mathrm{TC}}
\newcommand{\TCM}{\mathrm{TC}_\mathrm{M}}
\newcommand{\TCF}{\mathrm{TC}_\mathrm{F}}
\newcommand{\TCMF}{\mathrm{TC}_\mathrm{MF}}
\newcommand{\TCP}{\mathrm{TC}_\mathrm{p}}

\newcommand{\Su}{\mathrm{S}}
\newcommand{\SuM}{\mathrm{S}_\mathrm{M}}
\newcommand{\SuF}{\mathrm{S}_\mathrm{F}}
\newcommand{\SuMF}{\mathrm{S}_\mathrm{MF}}

\newcommand{\Stab}{\mathrm{Stab}}

\newcommand{\bDelta}{\bar{\Delta}}
\newcommand{\bF}{\bar{F}}

\newcommand{\ccc}{\text{ccc}\ }
\newcommand{\tDelta}{\tilde{\Delta}}
\newcommand{\bdelta}{\overline{\delta}}
\newcommand{\brho}{{\overline{\rho}}}
\newcommand{\bx}{\overline{x}}
\newcommand{\bg}{\overline{g}}

\newcommand{\sM}{\mathcal{M}}
\newcommand{\sR}{\mathcal{R}}

\newcommand{\sP}{\mathcal{P}}
\newcommand{\defeq}{\doteq}

\newcommand{\ecc}{\mathrm{ecc}}
\newcommand{\lecc}{\mathrm{lecc}}

\newcommand{\sG}{\mathcal{G}}
\newcommand{\sI}{\mathcal{I}}
\newcommand{\res}{\restriction}
\newcommand{\lh}{\mathrm{lh}}

\renewcommand{\emptyset}{\varnothing}
\renewcommand{\:}{\,:\,}

\newcommand{\cf}{c.f.}
\newcommand{\ie}{i.e.}
\newcommand{\ac}{\mathrm{AC}}
\newcommand{\zfc}{\mathrm{ZFC}}
\newcommand{\zf}{\mathrm{ZF}}

\begin{document}
\frontmatter
\title{Group Colorings and Bernoulli Subflows}

\author{Su Gao}
\address{Department of Mathematics, University of North Texas, 1155 Union Circle \#311430, Denton, TX 76203-5017}
\email{sgao@unt.edu}
\thanks{Su Gao's research was supported by the U.S. NSF Grants DMS-0501039 and DMS-0901853.}

\author{Steve Jackson}
\address{Department of Mathematics, University of North Texas, 1155 Union Circle \#311430, Denton, TX 76203-5017}
\email{jackson@unt.edu}
\thanks{Steve Jackson's research was supported by the U.S. NSF Grant DMS-0901853.}

\author{Brandon Seward}
\address{Department of Mathematics, University of Michigan, 530 Church Street, Ann Arbor, MI 48109-1043}
\email{b.m.seward@gmail.com}
\thanks{Brandon Seward's research was supported by two REU supplements to the U.S. NSF Grant DMS-0501039 and an NSF Graduate Research Fellowship.}

\date{December 28, 2011}

\begin{abstract}
In this paper we study the dynamics of Bernoulli flows and their subflows over general countable groups. One of the main themes of this paper is to establish the correspondence between the topological and the symbolic perspectives. From the topological perspective, we are particularly interested in free subflows (subflows in which every point has trivial stabilizer), minimal subflows, disjointness of subflows, and the problem of classifying subflows up to topological conjugacy. Our main tool to study free subflows will be the notion of hyper aperiodic points; a point is {\it hyper aperiodic} if the closure of its orbit is a free subflow. We show that the notion of hyper aperiodicity corresponds to a notion of {\it $k$-coloring} on the countable group, a key notion we study throughout the paper. In fact, for all important topological notions we study, corresponding notions in group combinatorics will be established. Conversely, many variations of the notions in group combinatorics are proved to be equivalent to some topological notions. In particular, we obtain results about the differences in dynamical properties between pairs of points which disagree on finitely many coordinates.

Another main theme of the paper is to study the properties of free subflows and minimal subflows. Again this is done through studying the properties of the hyper aperiodic points and minimal points. We prove that the set of all (minimal) hyper aperiodic points is always dense but meager and null. By employing notions and ideas from descriptive set theory, we study the complexity of the sets of hyper aperiodic points and of minimal points, and completely determine their descriptive complexity. In doing this we introduce a new notion of countable flecc groups and study their properties.
We also obtain the following results for the classification problem of free subflows up to topological conjugacy. For locally finite groups the topological conjugacy relation for all (free) subflows is hyperfinite and nonsmooth. For nonlocally finite groups the relation is Borel bireducible with the universal countable Borel equivalence relation.

The third, but not the least important, theme of the paper is to develop constructive methods for the notions studied. To construct $k$-colorings on countable groups, a fundamental method of construction of multi-layer marker structures is developed with great generality. This allows one to construct an abundance of $k$-colorings with specific properties. Variations of the fundamental method are used in many proofs in the paper, and we expect them to be useful more broadly in geometric group theory. As a special case of such marker structures, we study the notion of ccc groups and prove the ccc-ness for countable nilpotent, polycyclic, residually finite, locally finite groups and for free products.

\end{abstract}

\maketitle

\tableofcontents

\mainmatter
\chapter{Introduction}

In this paper we study Bernoulli flows over arbitrary countable groups (these are also known as Bernoulli shifts, Bernoulli systems, and Bernoulli schemes). The overall focus of this paper is on the development and application of constructive methods, with a particular emphasis on questions surrounding free subflows. The topics, methods, and results presented here should be of interest to at least researchers in descriptive set theory, symbolic dynamics, and topological dynamics, and may be of interest to researchers in C$^*$-algebras, ergodic theory, geometric group theory, and percolation theory. In Section \ref{BERN FLOW} we remind the reader the definitions of Bernoulli flow and subflow and also discuss the importance of Bernoulli flows to various areas of mathematics. In Section \ref{BASIC NOT} we introduce some basic notation and terminology which is needed for this chapter. In Section \ref{SEC EX FREE SUB} we discuss the question of the existence of free subflows. This question has been recently answered and is of importance to this paper. In Sections \ref{SEC INTRO HYP AP}, \ref{SEC INTRO COMPLX}, \ref{INTRO TILE}, and \ref{INTRO AER} we discuss the main results of this paper and at the same time discuss relevance to and motivation from various areas of mathematics, namely descriptive set theory, ergodic theory, geometric group theory, symbolic dynamics, and topological dynamics. A significant aspect of this paper is the invention of some versatile tools which add structure to arbitrary countable groups and offer significant aid in constructing points in Bernoulli flows. These tools are developed in great generality and likely have applications beyond their use here. These constructive methods and their potential utility to various areas of mathematics are discussed in Section \ref{INTRO FUND METH}. Finally, in Section \ref{INTRO OUTLINE} we give a brief outline to the paper and discuss chapter dependencies. We encourage the reader to make use of the detailed index found at the end of the paper which includes both terminology and notation.

\section{Bernoulli flows and subflows} \label{BERN FLOW}

Let us first begin by presenting the most general definition of a Bernoulli flow (also known as Bernoulli shift, Bernoulli system, and Bernoulli scheme). If $G$ is a countable group and $K$ is a set with the discrete topology and with a probability measure $\nu$, then the Bernoulli flow over $G$ with alphabet $K$ is defined to be
$$K^G = \{x \: G \rightarrow K\} = \prod_{g \in G} K$$
together with the product topology, the product measure $\nu^G$, and the following action of $G$: for $x \in K^G$ and $g \in G$, $g \cdot x \in K^G$ is defined by $(g \cdot x)(h) = x(g^{-1} h)$. The set $K$ is always assumed to have at least two elements as otherwise $K^G$ consists of a single point.

The action of $G$ on $K^G$ is quite intuitive. For example, if $G = \ZZ$ and $K = \{0,1\}$ then $K^G = \{0,1\}^\ZZ$ can be viewed as the space of all bi-infinite sequences of $0$'s and $1$'s. $\ZZ$ then acts by shifting these sequences left and right (the action of $5 \in \ZZ$ shifts these sequences $5$ units to the right). Similarly, $\{0,1\}^{\ZZ^2}$ can be visualized as the space of $\{0, 1\}$-labelings of the two dimensional lattice $\ZZ^2 \subseteq \mathbb{R}^2$ with the action of $\ZZ^2$ moving the labels in the obvious fashion. Comprehension of these examples should lead to an intuitive understanding of the action of $G$ on $K^G$.

Under the product topology, the basic open sets of $K^G$ are the sets of the form
$$\{x \in K^G \: \forall 1 \leq i \leq n \ x(h_i) = k_i\}$$
where $h_1, h_2, \ldots, h_n \in G$, $k_1, k_2, \ldots, k_n \in K$, and $n \geq 1$. Thus the action of $G$ on $K^G$ is continuous. It is not difficult to see that the basic open sets of $K^G$ are both open and closed (i.e. clopen\index{clopen}). Since every point is the intersection of a decreasing sequence of basic open sets, it follows that $K^G$ is totally disconnected (meaning that the only connected sets are the one point sets). $K^G$ is also seen to be perfect (meaning that there are no isolated points). Furthermore, $K^G$ is compact if and only if $K$ is finite. Thus by a well known theorem of topology, $K^G$ is homeomorphic to the Cantor set whenever $K$ is finite. On basic open sets the measure $\nu^G$ is given by
$$\nu^G(\{x \in K^G \: \forall 1 \leq i \leq n \ x(h_i) = k_i\}) = \prod_{1 \leq i \leq n} \nu(k_i).$$
Therefore the action of $G$ on $K^G$ is measure preserving. It may not be so clear, but this action is in fact ergodic.

In addition to Bernoulli flows, we are also very interested in their \emph{subflows} (also known as subshifts or subsystems). A subflow of a Bernoulli flow $K^G$ is simply a closed subset of $K^G$ which is stable under the action of $G$. Bernoulli flows and their subflows show up in many areas of mathematics. One reason is that they have a rich diversity of dynamical properties which allow them to model many phenomenon. This ``modeling'' shows up in many contexts, such as in ergodic theory, descriptive set theory, percolation theory, topological dynamics, and symbolic dynamics. In ergodic theory and descriptive set theory, the orbit structures of Bernoulli flows are used to model the orbit structures of measurable group actions on other measure spaces. More generally, they are used to model countable Borel equivalence relations as a well known result of Feldman-Moore states that every countable Borel equivalence relation on a standard Borel space is induced by a Borel action of a countable group (\cite{FM}). In the site percolation model of percolation theory, Bernoulli flows of the form $\{0,1\}^G$ are used to model the flow of liquids through porous materials. In topological dynamics, it is known that if a group $G$ acts continuously on a compact topological space $X$ and the action is expansive, then there is a Bernoulli flow $K^G$ over $G$, a subflow $S \subseteq K^G$, and a continuous surjection $\phi: S \rightarrow X$ which commutes with the action of $G$ (meaning $\phi(g \cdot s) = g \cdot \phi(s)$ for all $g \in G$ and $s \in S$). Furthermore, if $X$ is totally disconnected then $\phi$ can be chosen to be a homeomorphism. Similarly, if $X$ can be partitioned by a collection of clopen sets, then there is a subflow $S$ of a Bernoulli flow $K^G$ and a continuous surjection $\phi: X \rightarrow S$ which commutes with the action of $G$. These types of facts can be used to study Bernoulli flows via topological dynamics (for example, as in \cite{GU}), but more frequently topological dynamical systems are studied via Bernoulli flows. This latter approach led to the invention of symbolic dynamics (\cite{MH}) and its subsequent growth over the past seventy years. A classical example of the use of symbolic dynamics is the modeling of geodesics flows on manifolds by (the suspension of) subflows of Bernoulli flows over $\ZZ$. Traditionally only Bernoulli flows over $\ZZ$ and $\ZZ^n$ are studied in symbolic dynamics, but more recently Bernoulli flows over hyperbolic groups have been used to model the dynamics of hyperbolic groups acting on their boundary (\cite{CP}).

A key aspect of the importance of Bernoulli flows is their modeling capabilities, but there are several other reasons to study them as well. Indeed, Bernoulli flows may be considered interesting in and of themselves. This viewpoint can be seen in at least descriptive set theory, ergodic theory, and symbolic dynamics. Bernoulli flows serve as very natural examples of orbit equivalence relations, of measure preserving ergodic group actions, and of continuous group actions on compact spaces. At the same time, Bernoulli flows have very simple definitions yet their dynamical properties are very difficult to fully understand. A particularly nice and many times useful aspect of Bernoulli flows is that they are susceptible to combinatorial arguments, something which is typically not seen in other dynamical systems. Indeed, combinatorial approaches are a predominant feature both in symbolic dynamics and in this paper. Another source of motivation for studying Bernoulli flows is to understand the relationship between the algebraic properties of the acting group and the dynamical properties of the Bernoulli flow (a research program suggested by Gottschalk in \cite{G}). There are some known results of this type. For example, with complete knowledge of the dynamical properties of a Bernoulli flow $K^G$, one can determine if $G$ is amenable (\cite{CFW}), if $G$ has Kazhdan's property (T) (\cite{GW}), and the rank of $G$ if $G$ is a nonabelian free group (\cite{Ga}), to name a few. This is another aspect of Bernoulli flows which appears on several occasions in this paper. Finally, in topological dynamics Bernoulli flows are also studied in order to reveal properties of the greatest ambit of $G$, since it is known that the greatest ambit of $G$ is the enveloping semigroup of the Bernoulli flow $\{0,1\}^G$ (see \cite{GU}).

In this paper we study the dynamics of Bernoulli flows from the symbolic and topological viewpoints and employ ideas from descriptive set theory to gain further understanding. Although we do not study Bernoulli flows from the ergodic theory perspective, there is a topic we study (tileability properties of groups) which could be of interest to researchers in ergodic theory and geometric group theory.

\section{Basic notions} \label{BASIC NOT}

We study Bernoulli flows from the symbolic and topological perspectives. We therefore only want to consider Bernoulli flows over finite alphabets (these are precisely the compact Bernoulli flows, as mentioned in the previous section). So throughout the paper the term ``Bernoulli flow'' will always mean ``Bernoulli flow over a finite alphabet.'' We will also not make use of any measures (aside from a single lemma). So we will never specify measures on the alphabets or on the Bernoulli flows. Since the alphabet $K$ is always finite and the particular elements of $K$ are unimportant, we will always use $K = \{0, 1, \ldots, k-1\}$ for some positive integer $k > 1$. As is common in logic and descriptive set theory, we let the positive integer $k$ denote the set $\{0, 1, \ldots, k-1\}$. We therefore write $k^G = \{0, 1, \ldots, k-1\}^G$.

Let $G$ be a countable group and let $X$ be a compact Hausdorff space on which $G$ acts continuously (such as the Bernoulli flow $k^G$). A closed subset of $X$ which is stable under the group action is called a \emph{subflow of $X$}. We denote the closure of sets $A \subseteq X$ by $\overline{A}$. If $x \in X$, then the orbit of $x$ is denoted
$$[x] = \{g \cdot x \: g \in G\}.$$
Notice that $\overline{[x]}$ is the smallest subflow of $X$ containing $x \in X$. If $g \in G - \{1_G\}$, $x \in X$, and $g \cdot x = x$ then we call $g$ a \emph{period} of $x$. We call $x \in X$ \emph{periodic} if it has a period and otherwise we call $x$ \emph{aperiodic} (notice that here ``periodic'' and ``aperiodic'' differ from conventional use since most commonly these two terms relate to whether or not the orbit of $x$ is finite). A subflow of $X$ is called \emph{free} if it consists entirely of aperiodic points, and $x \in X$ is called \emph{hyper aperiodic} if $\overline{[x]}$ is free (in \cite{DS} such points are called limit aperiodic). In the specific case where $X$ is the Bernoulli flow $k^G$, we use \emph{$k$-coloring} interchangeably with ``hyper aperiodic.'' Notice that $x \in X$ is hyper aperiodic if and only if $x$ is contained in some free subflow, and furthermore the collection of all hyper aperiodic points is precisely the union of the collection of free subflows. A subflow $S \subseteq X$ is \emph{minimal} if $\overline{[s]} = S$ for all $s \in S$. Similarly, a point $x \in X$ is \emph{minimal} if $\overline{[x]}$ is minimal (this again differs from conventional terminology, since such points $x$ are usually called ``almost periodic''). Two points $x, y \in X$ are called \emph{orthogonal} if $\overline{[x]}$ and $\overline{[y]}$ are disjoint. Finally, two subflows $S_1, S_2 \subseteq X$ are \emph{topologically conjugate} if there is a homeomorphism $\phi: S_1 \rightarrow S_2$ which commutes with the action of $G$ (meaning $\phi(g \cdot s) = g \cdot \phi(s)$ for all $g \in G$ and $s \in S_1$). From the viewpoint of symbolic and topological dynamics, topologically conjugate subflows are essentially identical.

As mentioned in the previous section, a useful property of Bernoulli flows is that many topological and dynamical properties are found to have equivalent combinatorial characterizations. In fact, it is known that hyper aperiodicity, orthogonality, minimality, and topological conjugacy can all be expressed in a combinatorial fashion. We heavily rely on the combinatorial characterizations of these properties within the paper, and as a convenience to the reader we include proofs of these characterizations. Our heavy use of the combinatorial characterization of hyper aperiodicity led us to frequently use the term ``$k$-coloring'' in place of ``hyper aperiodic.'' The term emphasizes the combinatorial condition and is also reminiscent of the term ``coloring'' in graph theory as both roughly mean ``nearby things look different.'' We use the term ``hyper aperiodic'' within this chapter in order to emphasize the dynamical property as well as to avoid the possibility of the reader confusing ``$k$-colorings'' with arbitrary elements of $k^G$.

Now having gone through the basic definitions, let us repeat the second sentence of this introduction. The overall focus of this paper is on the development and application of constructive methods for Bernoulli flows, with a particular emphasis on questions surrounding free subflows.

\section{Existence of free subflows} \label{SEC EX FREE SUB}

The most basic, natural, and fundamental question one can ask about free subflows is:
\begin{quote}
Does every Bernoulli flow contain a free subflow? Equivalently, does every Bernoulli flow contain a hyper aperiodic point?
\end{quote}
This question is an important source of motivation for this paper, so we discuss it here at some length. One may at first hope that this question is answered by an existential measure theory or Baire category argument. Indeed, a promising well known fact is that the collection of aperiodic points in a Bernoulli flow always has full measure and is comeager\index{comeager} (i.e. second category, the countable intersection of dense open sets). However, it is not clear if a comeager set of full measure must contain a subflow, and furthermore a simple argument (included here in Section \ref{SECT SMALLNESS}) shows that the collection of all hyper aperiodic points in a Bernoulli flow is of measure zero and meager\index{meager} (i.e. first category, countable union of nowhere dense sets). The failure of measure theory and Baire category arguments suggests that a constructive approach to this question is needed. This is a bit concerning because even in the case of Bernoulli flows over $\ZZ$ constructions for hyper aperiodic points are not very simple. Nevertheless, we are led to ask: for which groups $G$ can one construct a hyper aperiodic element in at least some Bernoulli flow $k^G$? The two-sided Morse-Thue sequence provides a well known example of a hyper aperiodic point for all Bernoulli flows over $\ZZ$, so the existence of free subflows of Bernoulli flows over $\ZZ$ (and possibly $\ZZ^n$) has been known since at least the 1920's (when the Morse-Thue sequence was introduced). Only very recently were constructions for hyper aperiodic points found for other groups. In 2007, Dranishnikov-Shroeder proved that if $G$ is a torsion free hyperbolic group and $k \geq 9$ then $k^G$ contains a free subflow (\cite{DS}). Their proof essentially used the Morse-Thue sequence along certain geodesic rays of $G$. Shortly after, Glasner-Uspenskij proved in \cite{GU} that if $G$ is abelian or residually finite and $k > 1$ then $k^G$ contains a free subflow. They did this by constructing topological dynamical systems with certain properties and then using the modeling capabilities of Bernoulli flows to conclude that these Bernoulli flows contained free subflows.

The existence question of free subflows was finally resolved in a recent paper by the authors (\cite{GJS}) which provided a method for constructing hyper aperiodic points in $k^G$ for every countable group $G$ and every $k > 1$. In this paper we spend a great deal of time reproving this fact here, and in fact this paper entirely supersedes \cite{GJS}. In addition to presenting a general proof which applies to all Bernoulli flows, we also present alternative specialized proofs in the case of Bernoulli flows over abelian groups, solvable groups, FC groups, residually finite groups, and free groups (a group is FC if every conjugacy class is finite). As mentioned previously, a significant aspect of this paper is the development of powerful tools for constructing elements of Bernoulli flows. A very primitive and obscure form of these tools appeared in \cite{GJS} under a dense and technical presentation. Thankfully here the presentation is much more spread out, the tools are clearly distinguished and greatly generalized, and significant effort was put into making these tools understandable and more widely applicable. It is with the use of these tools that we prove essentially all of the results mentioned in the next four sections. The tools we develop in this paper come from \cite{GJS} and hence come from trying to answer the existence question for free subflows. We therefore place a lot of focus on the existence question in the first half of the paper and use the question as primary motivation for developing our tools.

For those readers who have a background in descriptive set theory, we would like to remark that the first two authors' original motivation for proving the existence of free subflows (in \cite{GJS}) came from the theory of Borel equivalence relations and in particular the theory of hyperfinite equivalence relations. In proving that the orbit equivalence relation on $2^\ZZ = \{0,1\}^\ZZ$ is hyperfinite, a key marker lemma by Slaman and Steel was the following (\cite{SS}).

\begin{lem}[Slaman-Steel]
Let $F(\ZZ)$ be the set of aperiodic points in $2^\ZZ$. Then there is an infinite decreasing sequence of Borel complete sections of $F(\ZZ)$
$$S_0 \supseteq S_1 \supseteq S_2 \supseteq \cdots$$
such that $\bigcap_{n \in \N} S_n = \varnothing$.
\end{lem}

This lemma remains true when $\ZZ$ is replaced by any countably infinite group $G$. The existence of decreasing sequences of complete sections that are relatively closed in $F(\ZZ)$ would allow one to easily construct a continuous embedding of $E_0$ into the orbit equivalence relation on $F(\ZZ)$. However, the existence of a free subflow of $2^G$ immediately implies (by compactness) that for every countable group $G$ there cannot exist a decreasing sequence of relatively closed complete sections of $F(G)$ whose intersection is empty (although a continuous embedding of $E_0$ into the orbit equivalence relation on $F(\ZZ)$ still does exist). The relationship of free subflows to this type of marker question is discussed a bit further in \cite{GJS}.

In the following four sections we discuss the results of this paper followed by a section discussing this paper's methods and tools. Again, we would like to emphasize that although the existence question of free subflows was previously resolved, we use the question here as primary motivation for developing our tools, and these tools in turn are vital to the proofs of nearly all of the results mentioned in the next four sections.

\section{Hyper aperiodic points and $k$-colorings} \label{SEC INTRO HYP AP}

When work began on this paper, one of the original goals was to investigate some of the basic properties of the set of hyper aperiodic points since at the time it was merely known that hyper aperiodic points existed. A natural first question is: How many hyper aperiodic points are there? Of course, this phrasing of the question is rather trivial since if $x$ is hyper aperiodic then $\overline{[x]}$ is uncountable and consists entirely of hyper aperiodic points. However, this question becomes more meaningful when attention is restricted to sets of hyper aperiodic points which are pairwise orthogonal. In this case, the answer to the question is not at all clear. A second natural question is: Is the set of hyper aperiodic points (equivalently, the union of the collection of free subflows) dense? Even more restrictive versions of these questions exist where one considers points which are both hyper aperiodic and minimal. Recall the fact mentioned in the previous section that the set of hyper aperiodic points is always of measure zero and is meager. This tells us that there is a dividing line after which these ``largeness'' questions regarding the set of hyper aperiodic points will have negative answers. Nevertheless, the results mentioned in this section reveal that the set of hyper aperiodic points is surprisingly large in a few respects. The above questions are all answered in succession over two chapters. The following is the crowning theorem of these investigations. 

\begin{theorem} \label{INTRO PERFDENSE}
Let $G$ be a countably infinite group, and let $k > 1$ be an integer. If $U \subseteq k^G$ is open and nonempty, then there exists a perfect (hence uncountable) set $P \subseteq U$ which consists of pairwise orthogonal minimal hyper aperiodic points.
\end{theorem}

We remind the reader that the proofs of all of our main results, including the theorem above, are entirely constructive. This theorem has three nice corollaries, all of which are new results. The first corollary ties in with the descriptive set theory connection mentioned in the previous section and requires further argument after the theorem above. The other two corollaries follow immediately from the theorem above but are also given direct proofs within this paper.

\begin{cor} \label{INTRO SS}
Let $G$ be a countably infinite group, let $k > 1$ be an integer, and let $F(G)$ denote the set of aperiodic points in $k^G$. If $(S_n)_{n \in \N}$ is a decreasing sequence of (relatively) closed complete sections of $F(G)$ (meaning each $S_n$ meets every orbit in $F(G)$) then
$$G \cdot \left( \bigcap_{n \in \N} S_n \right)$$
is dense in $k^G$.
\end{cor}

\begin{cor} \label{INTRO MINDENSE}
If $G$ is a countable group and $k > 1$ is an integer, then the collection of minimal points in $k^G$ is dense in $k^G$.
\end{cor}

\begin{cor} \label{INTRO DENSE}
If $G$ is a countably infinite group and $k > 1$ is an integer, then the collection of hyper aperiodic points in $k^G$ (equivalently the union of the free subflows of $k^G$) is dense in $k^G$.
\end{cor}

This last corollary is equivalent to the statement: if $A \subseteq G$ is finite and $y: A \rightarrow k = \{0, 1, \ldots, k-1\}$ then there is a hyper aperiodic point $x \in k^G$ which extends the function $y$. The above result therefore leads to the question: Which functions $y: A \rightarrow k$ with $A \subseteq G$ can be extended to a hyper aperiodic point $x \in k^G$? While we were unable to answer this question in full generality, we do prove two strong theorems which make substantial progress on resolving the question. The first such theorem is below. It completely characterizes those domains $A \subseteq G$ for which every function $y: A \rightarrow k$ can be extended to a hyper aperiodic point.

\begin{theorem} \label{INTRO SLENDER}
Let $G$ be a countably infinite group, let $A \subseteq G$, and let $k > 1$ be an integer. The following are all equivalent:
\begin{enumerate}
\item[\rm (i)] for every $y: A \rightarrow k$ there exists a perfect (hence uncountable) set $P \subseteq k^G$ consisting of pairwise orthogonal hyper aperiodic points extending $y$;
\item[\rm (ii)] for every $y: A \rightarrow k$ there exists a hyper aperiodic point $x \in k^G$ extending $y$;
\item[\rm (iii)] the function on $A$ which is identically $0$ can be extended to a hyper aperiodic point $x \in k^G$;
\item[\rm (iv)] there is a finite set $T \subseteq G$ so that for all $g \in G$ there is $t \in T$ with $gt \not\in A$.
\end{enumerate}
Furthermore, if $A \subseteq G$ satisfies any of the above equivalent properties then there is a continuous function $f: k^A \rightarrow k^G$ (where $k^A$ has the product topology) whose image consists entirely of hyper aperiodic points and with the property that $f(y)$ extends the function $y$ for each $y \in k^A$.
\end{theorem}

Notice that the set $A$ can be quite large while still satisfying clause (iv). For example, one could take $G = \ZZ^n$ and $A = \ZZ^n - (1000 \ZZ)^n$. Clearly any finite set $A$ satisfies (iv) if $G$ is infinite, so this theorem greatly generalizes the previous corollary. Also, since $H$ satisfies (iv) if $H \leq G$ is a proper subgroup, we have the following.

\begin{cor} \label{INTRO SLENDER2}
Let $G$ be a countably infinite group and let $k > 1$ be an integer. If $H \leq G$ is a proper subgroup of $G$, then every element of $k^H$ can be extended (continuously) to a (perfect set of pairwise orthogonal) hyper aperiodic point(s) in $k^G$.
\end{cor}

The second and final theorem addressing the extendability question stated above is the following. Recall that $k^G = \{0, 1, \ldots, k-1\}^G$.

\begin{theorem} \label{INTRO COFIN}
Let $G$ be a countable group and let $k > 1$ be an integer. If $A \subseteq G$ and $y: A \rightarrow k$, then let $y_* \in (k+1)^G$ be the function satisfying $y_*(a) = y(a)$ for $a \in A$ and $y_*(g) = k$ for $g \in G - A$. Then $y$ can be extended to a hyper aperiodic point in $k^G$ provided $G - A$ is finite and $\overline{[y_*]} \cap k^G$ consists of aperiodic points.
\end{theorem}

We remark that for many groups one thinks of, such as $\ZZ^n$, the above theorem is rather obvious. However, this is not always the case as there are groups whose Bernoulli flows have quite strange behavior. An example somewhat related to the theorem above is that for certain countable groups $G$, there are $x, y \in k^G$ which differ at precisely one coordinate and yet $x$ is hyper aperiodic while $y$ is periodic. This particular phenomenon is carefully studied in this paper and is discussed in Section \ref{INTRO AER}.

Regarding the extendability question stated above, we make the following conjecture.

\begin{conj}
Let $G$ be a countable group, let $A \subseteq G$, let $k > 1$ be an integer, and let $y: A \rightarrow k$. Define $y_* \in (k+1)^G$ by setting $y_*(a) = y(a)$ for $a \in A$ and $y_*(g) = k$ for $g \in G - A$. Then $y$ can be extended to a hyper aperiodic point in $k^G$ if and only if $\overline{[y_*]} \cap k^G$ consists of aperiodic points.
\end{conj}

If $y$ can be extended to a hyper aperiodic point in $k^G$, then it is easy to see that $\overline{[y_*]} \cap k^G$ consists of aperiodic points. The difficult question to resolve is if this condition is sufficient. Clearly this conjecture implies Theorem \ref{INTRO COFIN}. Also, if $A$ satisfies clause (iv) of Theorem \ref{INTRO SLENDER} and $y$ and $y_*$ are as above, then $\overline{[y_*]} \cap k^G$ must be empty. Thus the implication (iv) $\Rightarrow$ (ii) appearing in Theorem \ref{INTRO SLENDER} also follows from the above conjecture. We would like to emphasize that in all of the work we have done studying Bernoulli flows, we have always found the obvious necessary conditions to be sufficient. This is the main reason for us formally making this conjecture.

Related to the extendability question discussed above, one can ask a similar question of which functions $y: A \rightarrow k$ have the property that every point in $k^G$ extending $y$ is hyper aperiodic. There is a combinatorial characterization for this property, but it is rather trivial. However, if $A = H \leq G$ is a subgroup, then one can characterize this property through dynamical conditions on $y$ when $y$ is viewed as an element of the Bernoulli flow $k^H$. It is easy to see that if $y \in k^H$ and every extension of $y$ to $k^G$ is hyper aperiodic, then $y$ must itself be hyper aperiodic (as an element of $k^H$). So the question comes down to: for a subgroup $H \leq G$, which hyper aperiodic $y \in k^H$ have the property that every point in $k^G$ extending $y$ is hyper aperiodic? This is answered by the following theorem.

\begin{theorem} \label{INTRO AUTEXT}
Let $G$ be a countable group and let $k > 1$ be an integer. For a subgroup $H \leq G$, the following are equivalent:
\begin{enumerate}
\item[\rm (i)] there is some hyper aperiodic $y \in k^H$ for which every $x \in k^G$ extending $y$ is hyper aperiodic;
\item[\rm (i)] for every hyper aperiodic $y \in k^H$ and every $x \in k^G$ extending $y$, $x$ is hyper aperiodic;
\item[\rm (ii)] $H$ is of finite index in $G$ and $\langle g \rangle \cap H \neq \{1_G\}$ for every $1_G \neq g \in G$.
\end{enumerate}
Moreover, if every nontrivial subgroup $H \leq G$ satisfies the above equivalent conditions, then $G = \ZZ$.
\end{theorem}

In proving the above theorem, we prove the following interesting proposition. The proof of this proposition is nontrivial, and we do not know if its truth was previously known.

\begin{prop} \label{INTRO BNDSBGRP}
If $G$ is an infinite group and every nontrivial subgroup is of finite index, then $G = \ZZ$.
\end{prop}

\section{Complexity of sets and equivalence relations} \label{SEC INTRO COMPLX}

In this paper we study some complexity questions related to Bernoulli flows, and we approach such questions from the perspective of descriptive set theory. We remark that it is natural to use descriptive set theory as other notions of complexity (such as computability theory) are not generally applicable to Bernoulli flows since, for instance, not all groups have solvable word problem. There are two complexity issues we study here. The first is the descriptive complexities of the set of hyper aperiodic points, the set of minimal points, and the set of minimal hyper aperiodic points. The second is the complexity, under the theory of countable Borel equivalence relations, of the topological conjugacy relation among subflows of a common Bernoulli flow. We do not expect all readers to have previous knowledge of descriptive set theory and so we include a section which briefly introduces the notions and ideas surrounding the theory of countable Borel equivalence relations. The material should be readable to those who are interested. Also, we do review some terminology of these areas here, but only very briefly.

We first recall a bit of terminology from descriptive set theory. A topological space $X$ is Polish if it is separable and if its topology can be generated by a complete metric. A set is $\mathbf{\Sigma}^0_2$ (i.e. $F_\sigma$) if it can be expressed as the countable union of closed sets. Similarly, a set is $\mathbf{\Pi}^0_3$ (i.e. $F_{\sigma\delta}$) if it can be represented as the countable intersection of $\mathbf{\Sigma}^0_2$ sets. A subset $A \subseteq X$ of a Polish space $X$ is $\mathbf{\Sigma}^0_2$-complete if it is $\mathbf{\Sigma}^0_2$ and if for every Polish space $Y$ and every $\mathbf{\Sigma}^0_2$ subset $B \subseteq Y$ there is a continuous function $f: Y \rightarrow X$ with $B = f^{-1}(A)$. A similar definition applies to $\mathbf{\Pi}^0_3$-complete. Intuitively, $\mathbf{\Sigma}^0_2$-complete sets are the most complicated among all $\mathbf{\Sigma}^0_2$ subsets of Polish spaces, and similarly for $\mathbf{\Pi}^0_3$-complete sets.

The study of the descriptive complexity of the set of hyper aperiodic points leads us to define a new class of groups which we call flecc. We provide the simplest definition here. A group $G$ is \emph{flecc} if there is a finite set $A \subseteq G - \{1_G\}$ with the property that for every nonidentity $g \in G$ there is $n \in \ZZ$ and $h \in G$ with $h g^n h^{-1} \in A$. The choice of the name flecc comes from the acronym ``finitely many limit extended conjugacy classes.'' Various properties of flecc groups, other characterizations of flecc groups, and the meaning of the acronym can all be found in Section \ref{sec:flecc}. To the best of our knowledge, the class of flecc groups have never been isolated despite being associated to an interesting dynamical property. This dynamical property is the following. Let $G$ be a countable flecc group, let $A \subseteq G$ be the finite set described above, and let $X$ be any set on which $G$ acts. Then $X$ contains a periodic point if and only if $X$ contains a point having a period in the finite set $A$. To see this, suppose $x \in X$, $g \in G - \{1_G\}$, and $g \cdot x = x$. Then there is $n \in \ZZ - \{0\}$ and $h \in G$ with $h g^n h^{-1} = a \in A$. So we have that $a \in A$ is a period of the point $y = h \cdot x$. This dynamical property of flecc groups leads to a dichotomy in the descriptive complexity of the set of hyper aperiodic points, as seen in the theorem below.

\begin{theorem} \label{INTRO CMPLX1}
Let $G$ be a countable group and $k > 1$ an integer. Then the set of hyper aperiodic points in $k^G$ is closed if $G$ is finite, $\mathbf{\Sigma}^0_2$-complete if $G$ is an infinite flecc group, and $\mathbf{\Pi}^0_3$-complete if $G$ is an infinite nonflecc group.
\end{theorem}

The following theorem restricts attention to sets of minimal points and the dichotomy related to flecc groups disappears.

\begin{theorem} \label{INTRO CMPLX2}
Let $G$ be a countable group and $k > 1$ an integer. Then the set of minimal points in $k^G$ and the set of minimal hyper aperiodic points in $k^G$ are both closed if $G$ is finite and are both $\mathbf{\Pi}^0_3$-complete if $G$ is infinite.
\end{theorem}

Next we discuss the complexity of the topological conjugacy relation on subflows of a common Bernoulli flow. In other words, we study how difficult it is to determine if two subflows are topologically conjugate. For a countable group $G$ and an integer $k > 1$, let $\TC$ denote the topological conjugacy relation on subflows of $k^G$. Specifically, $\TC$ is the equivalence relation on the set of subflows of $k^G$ defined by the rule: $S_1 \ \TC \ S_2$ if and only if $S_1$ and $S_2$ are topologically conjugate. Let $\TCF$, $\TCM$, and $\TCMF$ denote the restriction of $\TC$ to the set of free subflows, minimal subflows, and free and minimal subflows, respectively. Also define an equivalence relation $\TCP$ on $k^G$ by declaring $x \ \TCP \ y$ if and only if $\overline{[x]}$ and $\overline{[y]}$ are topologically conjugate via a homeomorphism sending $x$ to $y$. We show that these five equivalence relations are always countable Borel equivalence relations.

Before stating the theorems, let us quickly introduce a few notions from the theory of Borel equivalence relations. An equivalence relation $E$ on a Polish space $X$ is Borel if it is a Borel subset of $X \times X$ (under the product topology), and the equivalence relation $E$ is countable if every equivalence class is countable. Given Borel equivalence relations $E$ and $F$ on $X$ and $Y$ respectively, $F$ is said to be Borel reducible to $E$ if there is a Borel function $f: Y \rightarrow X$ such that $y_1 \ F \ y_2$ if and only if $f(y_1) \ E \ f(y_2)$. Intuitively, in this situation $E$ is at least as complicated as $F$, or $F$ is no more complicated than $E$. There are countable Borel equivalence relations which all other countable Borel equivalence relations are Borel reducible to (so intuitively they are of maximal complexity), and such equivalence relations are called universal countable Borel equivalence relations. Finally, recall that $E_0$ is the nonsmooth hyperfinite equivalence relation on $2^\N$ defined by: $x \ E_0 \ y$ if and only if there is $n \in \N$ so that $x(m) = y(m)$ for all $m \geq n$.  

\begin{theorem} \label{INTRO E0}
Let $G$ be a countably infinite group and let $k > 1$ be an integer. Then $E_0$ continuously embeds into $\TCP$ and Borel embeds into $\TC$, $\TCF$, $\TCM$, and $\TCMF$.
\end{theorem}

This theorem has two immediate corollaries. We point out that on the space of all subflows of $k^G$ we use the Vietoris topology (see Section \ref{SECT BASIC TC}), or equivalently the topology induced by the Hausdorff metric. In symbolic and topological dynamics there is a lot of interest in finding invariants, and in particular searching for complete invariants, for topological conjugacy, particularly for subflows of Bernoulli flows over $\ZZ$ or $\ZZ^n$. The following corollary says that, up to the use of Borel functions, there are no complete invariants for the topological conjugacy relation on any Bernoulli flow.

\begin{cor} \label{INTRO E01}
Let $G$ be a countably infinite group and let $k > 1$ be an integer. Then there is no Borel function defined on the space of subflows of $k^G$ which computes a complete invariant for any of the equivalence relations $\TC$, $\TCF$, $\TCM$, or $\TCMF$. Similarly, there is no Borel function on $k^G$ which computes a complete invariant for the equivalence relation $\TCP$.
\end{cor}

The above theorem and corollary imply that from the viewpoint of Borel equivalence relations, the topological conjugacy relation on subflows of a common Bernoulli flow is quite complicated as no Borel function can provide a complete invariant. However, the above results do not rule out the possibility of the existence of algorithms for computing complete invariants among subflows described by finitary data, such as subflows of finite type. The above theorem also leads to another nice corollary. We do not know if the truth of the following corollary was previously known.

\begin{cor} \label{INTRO E02}
For every countably infinite group $G$, there are uncountably many pairwise non-topologically conjugate free and minimal continuous actions of $G$ on compact metric spaces.
\end{cor}
The following theorem completely classifies the complexity of $\TC$ and $\TCF$ for all countably infinite groups $G$. Again we see the interplay between group theoretic properties and dynamic properties as this theorem presents a dichotomy between locally finite and nonlocally finite groups. Recall that a group is called locally finite if every finite subset generates a finite subgroup.

\begin{theorem} \label{INTRO UNIEQ}
Let $G$ be a countably infinite group and let $k > 1$ be an integer. If $G$ is locally finite then $\TC$, $\TCF$, $\TCM$, $\TCMF$, and $\TCP$ are all Borel bi-reducible with $E_0$. If $G$ is not locally finite then $\TC$ and $\TCF$ are universal countable Borel equivalence relations.
\end{theorem}

This last theorem generalizes a result of John Clemens which states that for the Bernoulli flow $k^{\ZZ^n}$ the equivalence relation $\TC$ is a universal countable Borel equivalence relation (\cite{JC}).

\section{Tilings of groups} \label{INTRO TILE}

A key aspect of the main constructions of this paper is the use of marker structures. Marker structures can be placed on groups or on sets on which groups act and are a geometrically motivated way of studying groups and their actions. They have been used numerous times in ergodic theory, the theory of Borel equivalence relations (especially the theory of hyperfinite equivalence relations), and even in symbolic dynamics (for studying the automorphism groups of Bernoulli flows over $\ZZ$). In working with Bernoulli flows, it became apparent that solving problems through algebraic methods was cumbersome, placed restrictions on the groups we could consider, and resulted in case-by-case proofs. However, we found that solving problems through geometric methods (specifically through marker structures) relaxed and many times removed restrictions on the groups and resulted in unified arguments. The stark difference between algebraic and geometric methods can clearly be seen in this paper in Chapters \ref{chap:basicconstructions} (algebraic) and \ref{CHAP MARKER} (geometric). It is for this reason that we define and study marker structures on groups. Naturally, better marker structures lead to better proofs. We therefore consider strong types of marker structures such as tilings and \ccc sequences of tilings. 

For a countably infinite group $G$ and a finite set $T \subseteq G$, we call $T$ a \emph{tile} if there is a set $\Delta \subseteq G$ such that the the set $\{\gamma T \: \gamma \in \Delta\}$ partitions $G$. Such a pair $(\Delta, T)$ is a \emph{tiling} of $G$. A sequence of tilings $(\Delta_n, T_n)_{n \in \N}$ is \emph{coherent} if each set $\gamma T_{n+1}$ with $\gamma \in \Delta_{n+1}$ is the union of left $\Delta_n$ translates of $T_n$. A sequence of tilings $(\Delta_n, T_n)_{n \in \N}$ is \emph{centered} if $1_G \in \Delta_n$ for all $n \in \N$. Finally, a centered sequence of tilings $(\Delta_n, T_n)_{n \in \N}$ is \emph{cofinal} if $T_n \subseteq T_{n+1}$ and $G = \bigcup_{n \in \N} T_n$. We abbreviate the three adjectives ``coherent, centered, and cofinal'' to simply \emph{ccc}. We call $G$ a \emph{ccc group} if $G$ admits a \ccc sequence of tilings.

The study of \ccc groups has applications to ergodic theory as it ties in with Rokhlin sets and is closely related with the study of monontileable amenable groups initiated by Chou (\cite{Ch}) and Weiss (\cite{W}). In fact, \ccc groups form a subset of what Weiss called MT groups in \cite{W}. In our terminology, a group is MT if it admits a centered and cofinal sequence of tilings. Ccc groups are also pertinent to the theory of hyperfinite equivalence relations. Progress in the theory of hyperfinite equivalence relations has so far been dependent on finding better and better marker structures on Bernoulli flows over larger and larger classes of groups. The study of \ccc groups, and in fact marker structures on groups in general, gives an upper bound to the types of marker structures which can be constructed on Bernoulli flows and also may give an informal sense of properties such marker structures may have. The notion of \ccc groups is also interesting from the geometric group theory perspective. Groups and their Cayley graphs display such a high degree of symmetry and homogeneity that it is hard to imagine the existence of a group which is not \text{ccc}, or even worse, not MT. This is strongly contrasted by the fact that it seems to be in general very difficult to determine if a group is \text{ccc}. This geometric property is therefore somewhat mysterious.

We prove that a few large classes of groups are \text{ccc}, as indicated in the following theorem.

\begin{theorem}
The following groups are \ccc groups:
\begin{enumerate}
\item[\rm (i)] countable locally finite groups;
\item[\rm (ii)] countable residually finite groups;
\item[\rm (iii)] countable nilpotent groups;
\item[\rm (iv)] countable solvable groups $G$ for which $[G, G]$ is polyclic;
\item[\rm (v)] countable groups which are the free product of a collection of nontrivial groups.
\end{enumerate}
\end{theorem}

Notice that every countable polycyclic group satisfies (iv) and is thus a \ccc group. Abelian groups are nilpotent, and linear groups (complex and real) are residually finite, so these classes of groups are also covered by the above theorem.

We do not know of any countably infinite group which is not \text{ccc}. Ideally, the methods of proof used here would help in finding new classes of groups which are \ccc and in constructing better marker structures on Bernoulli flows and other spaces on which groups act. To this end, this paper includes entirely distinct proofs that the following classes of groups are \text{ccc}: finitely generated abelian groups, nonfinitely generated abelian groups, nilpotent groups, polycyclic groups, residually finite groups, locally finite groups, nonabelian free groups, and free products of nontrivial groups. Despite having a direct proof that polycyclic groups are \ccc (one which does not use the fact that polycyclic groups are residually finite), we were unable to determine if solvable groups are \text{ccc}.

\section{The almost equality relation} \label{INTRO AER}

Finally, in this paper we study the behavior of periodic, aperiodic, and hyper aperiodic points under the almost-equality relation. Points $x, y \in k^G$ are \emph{almost equal}, written $x =^* y$, if as functions on $G$ they differ at only finitely many coordinates. In studying the almost equality relation and in establishing the results mentioned in this section, we had to make substantial use of tools and notions from geometric group theory. To be specific, the proofs of many of the theorems in this section relied heavily upon the notion of the number of ends of a group and on Stallings' theorem regarding groups with more than one end.

We also introduce and study a notion stronger than almost equality. If $x, y \in k^G$ then we write $x =^{**} y$ if $x$ and $y$ disagree on \emph{exactly} one element of $G$ (so $x \neq^{**} x$). We first study the relationship between periodicity and almost equality and obtain the following. Below, $\Stab(x)$ denotes the stabilizer subgroup $\{g \in G \: g \cdot x = x\}$.

\begin{theorem} \label{INTRO AE PERIODIC}
Let $k > 1$ be an integer.
\begin{enumerate}
\item[\rm (i)] Let $G$ be a countable group not containing any subgroup which is a free product of nontrivial groups. Then for every $x \in k^G$ and every $y =^{**} x$, either $x$ is aperiodic or $y$ is aperiodic.
\item[\rm (ii)] Let $G$ be a countable group containing $\ZZ_2 * \ZZ_2$ as a subgroup and with the property that every subgroup of $G$ which is the free product of two nontrivial groups is isomorphic to $\ZZ_2 * \ZZ_2$. Then for every $x \in k^G$ there is an aperiodic $y \in k^G$ with $y =^* x$, but there are periodic $w, z \in k^G$ with $w =^{**} z$.
\item[\rm (iii)] Let $G$ be a countable group containing a subgroup which is the free product of two nontrivial groups one of which has more than two elements. Then there is $x \in k^G$ such that every $y =^*x$ is periodic.
\end{enumerate}
Furthermore for any countable group $G$, if $x =^{**} y \in k^G$ then $\langle \Stab(x) \cup \Stab(y) \rangle \cong \Stab(x) * \Stab(y)$.
\end{theorem}

Clearly every countable group is described by precisely one of the clauses (i), (ii), and (iii). We point out that torsion groups fall into clause (i) and an amenable group will fall into either cluase (i) or clause (ii), depending on whether or not it contains $\ZZ_2 * \ZZ_2$ as a subgroup. As the dynamical properties described are mutually incompatible, we see that for each individual clause, the stated dynamical property characterizes the class of groups described. Thus, for example, if $G$ is a countable group such that for every $x \in k^G$ and every $y =^{**} x$ either $x$ or $y$ is aperiodic, then $G$ does not contain any subgroup which is a free product of nontrivial groups. We show that the groups described in clause (iii) are precisely those countable groups containing nonabelian free subgroups. Thus the dynamical property stated in clause (iii) provides a dynamical characterization of those groups  which contain nonabelian free subgroups.

The above theorem does not require the group $G$ to be infinite, so for finite groups clause (i) leads to the following corollary.

\begin{cor} \label{INTRO FINITE APER}
If $G$ is a finite group and $k > 1$ is an integer then $k^G$ contains at least $(k-1) k^{|G|-1}$ many aperiodic points.
\end{cor}

We also study the behavior of hyper aperiodic points under almost equality. A difficult basic question is if any point almost equal to a hyper aperiodic point must be hyper aperiodic itself. The following theorem says that, suprisingly, this is not always the case.

\begin{theorem} \label{INTRO ACP}
Let $G$ be a countable group and let $k > 1$ be an integer. The following are equivalent:
\begin{enumerate}
\item[\rm (i)] there are $x, y \in k^G$ with $x$ hyper aperiodic and $y$ not hyper aperiodic but $y =^* x$;
\item[\rm (ii)] there are $x, y \in k^G$ with $x$ hyper aperiodic and $y$ periodic but $y =^{**} x$;
\item[\rm (iii)] there is a nonidentity $u \in G$ such that every nonidentity $v \in \langle u \rangle$ has finite centralizer in $G$.
\end{enumerate}
\end{theorem}

We show that abelian groups, nilpotent groups, and nonabelian free groups never satisfy the equivalent conditions listed above. However, we find examples of  groups which are polycyclic (hence solvable) and finite extensions of abelian groups which do satisfy the equivalent conditions above.

As the previous theorem indicates, in general there may be points which are not hyper aperiodic but are almost equal to a hyper aperiodic point. We study the behaviour of such points and arrive at the following theorem.

\begin{theorem} \label{INTRO AE HYP APER}
For a countable group $G$, an integer $k > 1$, and $x \in k^G$, the following are all equivalent:
\begin{enumerate}
\item[\rm (i)] there is a hyper aperiodic $y \in k^G$ with $y =^* x$;
\item[\rm (ii)] there is a hyper aperiodic $y \in k^G$ such that $x$ and $y$ disagree on at most one coordinate;
\item[\rm (iii)] either $x$ is hyper aperiodic or else every $y =^{**} x$ is hyper aperiodic;
\item[\rm (iv)] every limit point of $[x]$ is aperiodic;
\item[\rm (v)] for every nonidentity $s \in G$ there are finite sets $A, T \subseteq G$ so that for all $g \in G - A$ there is $t \in T$ with $x(gt) \neq x(gst)$;
\end{enumerate}
\end{theorem}

We remark that clause (v) is very similar to the combinatorial characterization of being hyper aperiodic. Specifically, a point $x \in k^G$ is hyper aperiodic if and only if it satisfies the condition stated in clause (v) with $A$ restricted to being the empty set.

The method of proof of the previous two theorems leads to the following interesting corollary regarding more general dynamical systems.

\begin{cor} \label{INTRO TOPCOR}
For a countable group $G$, the following are equivalent:
\begin{enumerate}
\item[\rm (i)] for every compact Hausdorff space $X$ on which $G$ acts continuously and every $x \in X$, if every limit point of $[x]$ is aperiodic then $x$ is hyper aperiodic;
\item[\rm (ii)] for every nonidentity $u \in G$ there is a nonidentity $v \in \langle u \rangle$ having infinite centralizer in $G$.
\end{enumerate}
\end{cor}

\section{The fundamental method} \label{INTRO FUND METH}

All of the proofs within this paper are constructive, and nearly all of them rely on a single general, adaptable, and powerful method for constructing points in Bernoulli flows which we call the \emph{fundamental method}. A tremendous amount of time and effort was put into developing the fundamental method as if it were a general theory in itself, and in fact the method was intentionally developed in much greater generality than we make use of here. The fundamental method relies on three objects: a blueprint on the group, a locally recognizable function, and a sequence of functions of subexponential growth. These objects are not fixed but rather each must satisfy a general definition. We will not give precise definitions of these objects at this time, but we will give some indication as to what these objects are.

A blueprint on a group $G$ is a sequence $(\Delta_n, F_n)_{n \in \N}$ which is somewhat similar to a \ccc sequence of tilings. The sets $\Delta_n \subseteq G$ are in some sense evenly spread out within $G$ as there are finite sets $B_n$ with $\Delta_n B_n = G$. The left translates of $F_n$ by $\Delta_n$ are pairwise disjoint, and if a left translate of $F_k$ by $\Delta_k$ meets a left translate of $F_n$ by $\Delta_n$, then the former must be a subset of the latter provided $k \leq n$. Furthermore, for $k < n$ the left translates of $F_k$ by $\Delta_k$ appear in an identical fashion within every left translate of $F_n$ by an element of $\Delta_n$. It is a nontrivial fact that every countably infinite group has a blueprint. In this paper we in fact prove that a very strong type of blueprint exists on every countably infinite group. A locally recognizable function is a function $R: A \rightarrow k$ where $1_G \in A \subseteq G$ is finite and $k > 1$ is an integer. This function must have the property that the identity, $1_G$, is recognizable in the sense that if $a \in A$ and $R(ab) = R(b)$ for all $b \in A$ whenever both are defined, then $a = 1_G$. Again, we show that locally recognizable functions always exist and we give several nontrivial examples. Finally, a sequence of functions of subexponential growth is a sequence $(p_n)_{n \geq 1}$ such that each $p_n: \N \rightarrow \N$ is of subexponential growth.

We present a fixed construction which when given any blueprint, locally recognizable function, and sequence of functions of subexponential growth (under a few restrictions) produces a function $c$ taking values in $k$ and having a large infinite subset of $G$ as domain. This function $c$ has very nice properties related to the blueprint, the locally recognizable function, and the sequence of functions of subexponential growth. One can in fact determine if $g \in \Delta_n$ merely by the behavior of $c$ on the set $g F_n$. Furthermore, each left translate of $F_n$ by $\Delta_n$ has its own proprietary points on which $c$ is undefined. The number of such points is at least $\log p_n(|F_n|)$. If $t$ is the number of undefined points within a translate of $F_n$, then using $k$ values one can extend $c$ on this translate of $F_n$ in $k^t$ many ways. So the logarithm essentially disappears and on each left translate of $F_n$ by $\Delta_n$ one can essentially encode an amount of information which is subexponentially related to the size of $F_n$. The fact that the members of $\Delta_n$ are distinguishable within $c$ allows one to both encode and decode information using the undefined points of $c$. This fact is tremendously useful. Finally, the relationship of $c$ to the locally recognizable function $R$ is that near every member of $\Delta_1$ one sees a ``copy'' of $R$ in $c$.

Each of the three objects which go into the construction have their own strengths and weaknesses in terms of creating points in Bernoulli flows with certain desired properties. In fact, we prove the existence of hyper aperiodic points by using only functions of subexponential growth, we prove the density of hyper aperiodic points by primarily using locally recognizable functions, and we prove the existence of a point which is not hyper aperiodic but almost equal to a hyper aperiodic point by using a special blueprint. The fundamental method refers to the coordinated use of these three objects in achieving a goal of constructing a particular type of element of a Bernoulli flow. To further aid the use of the fundamental method, we develop two general tools which enhance the fundamental method. The first tool is a general method of constructing minimal points in Bernoulli flows. The second is a method of constructing sets of points in Bernoulli flows which have the property that the closures of the orbits of the points are pairwise not topologically conjugate. The fundamental method and these tools have been tremendously useful within this paper as nearly all of our results rely on them, and we hope that they will be similarly useful to other researchers in the future.

In view of some basic questions which were only recently answered in \cite{GU}, it seems that the dynamics of Bernoulli flows over general countable groups has received little investigation from the symbolic and topological points of view. This paper demonstrates that this need not be the case. Although there are problems in symbolic dynamics which seem too difficult even in the case of Bernoulli flows over $\ZZ^n$, there are likely many other interesting questions and properties which can be pursued over a larger class of groups or possibly even all countable groups. The fundamental method, which is entirely combinatorial, provides at least one way of approaching this. Such investigations will also benefit topological dynamics both through the modeling capabilities of Bernoulli flows and through supplying new examples of continuous group actions with various properties. Such examples may lead to quite general results similar to Corollary \ref{INTRO TOPCOR} (in this corollary, showing (ii) implies (i) is quite simple, and to show the negation of (ii) implies the negation of (i) one uses a subflow of a Bernoulli flow). On a final note, we also mention that blueprints provide a nontrivial structure to all countably infinite groups which could be useful in various situations.

\section{Brief outline} \label{INTRO OUTLINE}

In Chapter 2 we reintroduce the main definitions, terminology, and notation of this paper. This is done in a more elaborate and detailed manner than in this introduction. We also state and prove the combinatorial characterizations for dynamical properties such as hyper aperiodicity / $k$-coloring, orthogonality, and minimality (recall that hyper aperiodic and $k$-coloring have identical meanings on the Bernoulli flow $k^G$). We also present other notions, terminology, and ideas which are not present in this introduction. Various simple lemmas related to these concepts are also presented. It is recommended that readers do not skip Chapter 2 as the terminology and notation introduced is important to the rest of the paper.

In Chapter 3 we present various algebraic methods for constructing hyper aperiodic points / $k$-colorings. In particular, we present several general methods for extending $k$-colorings on a smaller group to a larger group. We also give direct (algebraic) proofs that all abelian groups admit $k$-colorings (different proofs are provided for different classes of abelian groups), all solvable groups admit $k$-colorings, all nonabelian free groups admit $k$-colorings, and all residually finite groups admit $k$-colorings. This chapter is independent of all later chapters and can be skipped if desired. The chapter should be of interest to anyone with a strong interest in hyper aperiodic points / $k$-colorings. The chapter also demonstrates the limitations of trying to construct $k$-colorings through algebraic methods.

In Chapter 4 we define a general notion of a marker structure on a group. We then use marker structures to provide a new (geometric) proof that all abelian and all FC groups admit $k$-colorings (a group is FC if all of its conjugacy classes are finite). These marker methods are more geometrically motivated than algebraicly motivated and prove to be much more succesful than the algebraic methods used in Chapter 3. For the rest of the paper our proofs rely on these geometric and marker structure methods. We therefore study groups with particularly strong marker structures - the \ccc groups. The study of \ccc groups comprises a significant portion of Chapter 4. This chapter contains all of the results mentioned in Section \ref{INTRO TILE} above. Technically speaking, the only thing from Chapter 4 which is used in later chapters is the definition of a marker structure. However, the marker structure proof that abelian and FC groups admit $k$-colorings is a good source of motivation for the machinery which is developed in Chapter 5.

Chapter 5 focuses almost entirely on developing the fundamental method and its related machinery. The only result pertaining to Bernoulli flows in this chapter is that the minimal points in a Bernoulli flow are always dense (Corollary \ref{INTRO MINDENSE} from this introduction). This chapter is of great importance to the rest of the paper. The only sections which can be read if Chapter 5 is skipped are Sections \ref{SECT SMALLNESS}, \ref{sec:flecc}, \ref{INTRO RELATNS}, \ref{SECT BASIC TC}, \ref{sec:charext}, \ref{SECT ALMOST EQUAL}, and \ref{sec:autoext}.

To make up for Chapter 5 being nearly devoid of results pertaining to Bernoulli flows, Chapter 6 focuses on presenting short, simple, and satisfying applications of the fundamental method. Each section focuses on one of the objects used in the fundamental method: functions of subexponential growth, locally recognizable functions, and blueprints. Many results are included in this chapter. Those mentioned in this introduction include a weaker version of Theorem \ref{INTRO PERFDENSE} which does not mention minimality, Corollary \ref{INTRO DENSE}, Theorem \ref{INTRO ACP}, and Corollary \ref{INTRO TOPCOR}. The only sections which do not rely on Chapter 6 are those listed in the previous paragraph which do not rely on Chapter 5.

In Chapter 7 we return to developing machinery again. We develop the two tools mentioned in the previous section which enhance the fundamental method. More specifically, we develop tools for using the fundamental method to construct minimal points and to construct sets of points which have the property that the closures of the orbits of the points are pairwise not topologically conjugate. We also investigate properties of fundamental functions - those functions which are constructed through the fundamental method. Additionally, we use the tools we develop to prove Theorem \ref{INTRO PERFDENSE} and Corollary \ref{INTRO SS}. Chapter 10 can be read without reading Chapter 7.

In Chapter 8 we investigate the descriptive complexity of various important subsets of Bernoulli flows. Specifically, we study the descriptive complexities of the sets of hyper aperiodic points, minimal points, and points which are both minimal and hyper aperiodic. We prove Theorems \ref{INTRO CMPLX1} and \ref{INTRO CMPLX2}. We also spend an entire section defining flecc groups and studying their properties. Chapters 9 and 10 are independent of Chapter 8.

In Chapter 9 we study the complexity of the topological conjugacy relation among subflows of a common Bernoulli flow. In other words, we study how difficult it is to determine when two subflows are topologically conjugate. This is done using the theory of countable Borel equivalence relations. A brief introduction to the theory of Borel equivalence relations is included in Section \ref{INTRO RELATNS}. In this chapter we prove Theorems \ref{INTRO E0} and \ref{INTRO UNIEQ} and Corollaries \ref{INTRO E01} and \ref{INTRO E02}. Chapter 10 is independent of Chapter 9.

In Chapter 10 we study both the extendability of partial functions to $k$-colorings and further properties of the almost-equality relation. We prove Theorems \ref{INTRO SLENDER}, \ref{INTRO COFIN}, and \ref{INTRO AUTEXT}, Corollary \ref{INTRO SLENDER2}, and all of the results mentioned in Section \ref{INTRO AER} aside from Theorem \ref{INTRO ACP} and Corollary \ref{INTRO TOPCOR} (Theorem \ref{INTRO ACP} and Corollary \ref{INTRO TOPCOR} are proven in Chapter 6). Somewhat surprisingly, Sections \ref{sec:charext}, \ref{SECT ALMOST EQUAL}, and \ref{sec:autoext} do not rely on any previous material in the paper aside from a few definitions. Furthermore, Section \ref{sec:charext} has an entirely self contained proof of Theorem \ref{INTRO SLENDER}. The dialog present in the proof of this theorem assumes the reader is familiar with the fundamental method, but this is only to aid in comprehension of the proof as no knowledge of the fundamental method is technically required.

Finally, in Chapter 11 we list some open problems.

Chapter dependencies are illustrated in Figure \ref{FIG CHAPFLOW} below. We again would like to remind the reader that there is a detailed index at the end of this paper which includes both terminology and notation.

\begin{figure}[ht]\label{FIG CHAPFLOW}
\begin{center}
\setlength{\unitlength}{3mm}
\begin{picture}(30,43)(5,-.5)

\put(0,26){
\put(0,0){\line(1,0){6}}
\put(0,0){\line(0,1){2}}
\put(6,0){\line(0,1){2}}
\put(0,2){\line(1,0){6}}
\put(3,.4){\makebox(0,0)[b]{Chapter 2}}
}

\put(11,36){
\put(0,0){\line(1,0){6}}
\put(0,0){\line(0,1){2}}
\put(6,0){\line(0,1){2}}
\put(0,2){\line(1,0){6}}
\put(3,.4){\makebox(0,0)[b]{Chapter 3}}
}

\put(11,31){
\put(0,0){\line(1,0){6}}
\put(0,0){\line(0,1){2}}
\put(6,0){\line(0,1){2}}
\put(0,2){\line(1,0){6}}
\put(3,.4){\makebox(0,0)[b]{Chapter 4}}
}

\put(11,26){
\put(0,0){\line(1,0){6}}
\put(0,0){\line(0,1){2}}
\put(6,0){\line(0,1){2}}
\put(0,2){\line(1,0){6}}
\put(3,.4){\makebox(0,0)[b]{Chapter 5}}
}

\put(22,26){
\put(0,0){\line(1,0){6}}
\put(0,0){\line(0,1){2}}
\put(6,0){\line(0,1){2}}
\put(0,2){\line(1,0){6}}
\put(3,.4){\makebox(0,0)[b]{Chapter 6}}
}

\put(33,26){
\put(0,0){\line(1,0){6}}
\put(0,0){\line(0,1){2}}
\put(6,0){\line(0,1){2}}
\put(0,2){\line(1,0){6}}
\put(3,.4){\makebox(0,0)[b]{Chapter 7}}
}

\put(13.5,8){
\put(0,0){\line(1,0){12}}
\put(0,0){\line(0,1){6}}
\put(12,0){\line(0,1){6}}
\put(0,6){\line(1,0){12}}
\put(6,0){\line(0,1){4}}
\put(0,4){\line(1,0){12}}
\put(6,4.4){\makebox(0,0)[b]{Chapter 8}}
\put(3,2.4){\makebox(0,0)[b]{Sections}}
\put(3,.8){\makebox(0,0)[b]{8.1 \& 8.3}}
\put(9,2.4){\makebox(0,0)[b]{Sections}}
\put(9,.8){\makebox(0,0)[b]{8.2 \& 8.4}}
}

\put(13.5,0){
\put(0,0){\line(1,0){12}}
\put(0,0){\line(0,1){6}}
\put(12,0){\line(0,1){6}}
\put(0,6){\line(1,0){12}}
\put(6,0){\line(0,1){4}}
\put(0,4){\line(1,0){12}}
\put(6,4.4){\makebox(0,0)[b]{Chapter 9}}
\put(3,2.4){\makebox(0,0)[b]{Sections}}
\put(3,.8){\makebox(0,0)[b]{9.1 \& 9.2}}
\put(9,2.4){\makebox(0,0)[b]{Sections}}
\put(9,.8){\makebox(0,0)[b]{9.3 \& 9.4}}
}

\put(7,16){
\put(0,0){\line(1,0){14}}
\put(0,0){\line(0,1){6}}
\put(14,0){\line(0,1){6}}
\put(0,6){\line(1,0){14}}

\put(9,0){\line(0,1){4}}

\put(0,4){\line(1,0){14}}
\put(7,4.4){\makebox(0,0)[b]{Chapter 10}}
\put(4.5,2.4){\makebox(0,0)[b]{Sections}}
\put(4.5,.8){\makebox(0,0)[b]{10.2, 10.3, 10.4}}
\put(11.5,2.4){\makebox(0,0)[b]{Section}}
\put(11.5,.8){\makebox(0,0)[b]{10.1}}
}

\put(6,27){\vector(1,0){5}}
\put(3,28){\line(0,1){9}}
\put(3,32){\vector(1,0){8}}
\put(3,37){\vector(1,0){8}}
\multiput(14,33)(0,.2){15}{\line(0,1){.1}}
\put(14,33.2){\vector(0,-1){.2}}
\multiput(14,28)(0,.2){15}{\line(0,1){.1}}
\put(14,28.2){\vector(0,-1){.2}}
\put(17,27){\vector(1,0){5}}
\put(28,27){\vector(1,0){5}}
\put(3,26){\line(0,-1){25}}
\put(3,17){\vector(1,0){4}}
\put(3,9){\vector(1,0){10.5}}
\put(3,1){\vector(1,0){10.5}}
\put(25,26){\line(0,-1){9}}
\put(25,17){\vector(-1,0){4}}
\put(36,26){\line(0,-1){25}}
\put(36,9){\vector(-1,0){10.5}}
\put(36,1){\vector(-1,0){10.5}}

\end{picture}
\caption{Flowchart of dependency between chapters. Solid arrows indicate dependencies; dashed arrows indicate motivation.}
\end{center}
\end{figure}
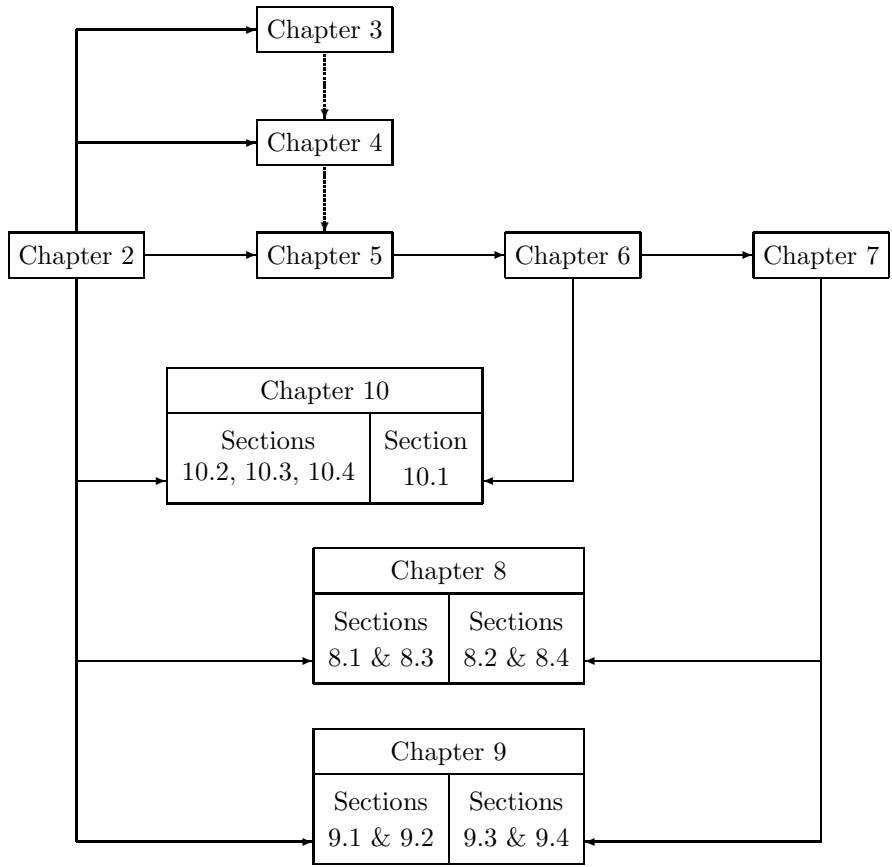 
\chapter{Preliminaries}

In this chapter we go through the definitions presented in the previous chapter in more detail. In the first section we go over some basic definitions and notation. The second, third, and fourth sections discuss the three most central notions: $k$-colorings, minimality, and orthogonality. In these sections we prove that these notions admit equivalent dynamical and combinatorial definitions. We also present a few basic lemmas related to these properties. In the fifth section we tweak the definition of $k$-colorings in various ways to obtain other interesting notions. These notions play an important role in this paper but were not mentioned in the previous chapter. Section five also contains several lemmas and propositions related to these notions. The sixth section discusses further generalizations of the notion of a $k$-coloring, however the discussion in this section is purely speculative as the notions introduced in this section are not studied within the paper. Finally, in the seventh section we discuss generalizations to the action of $G$ on $(2^\N)^G$. This last section has connections with descriptive set theory and topological dynamics.

\section{Bernoulli flows} \label{SEC BERNFLOW2}

For a positive integer $k$, we let $k$ denote the set $\{0, 1, \ldots, k-1\}$. If $G$ is a countable group, then a \emph{Bernoulli flow}\index{Bernoulli flow}\index{Bernoulli $G$-flow} over $G$, or a \emph{Bernoulli $G$-flow}, is a topological space of the form
$$k^G = \{0, 1, \ldots, k-1\}^G = \{x \: G \rightarrow k\} = \prod_{g \in G} \{0, 1, \ldots, k-1\},$$
equipped with the product topology, together with the following action of $G$: for $x \in k^G$ and $g \in G$, $g \cdot x \in k^G$ is defined by $(g \cdot x)(h) = x(g^{-1} h)$ for $h \in G$. Notice that $k^G$ is compact and metrizable. We will often find it convenient to work with the following compatible metric on $k^G$. Fix a countable group $G$, and fix an enumeration $g_0, g_1, g_2, \ldots$ of the group elements of $G$ with $g_0 = 1_G$ (the identity element). For $x, y \in k^G$, we define\index{$d(x,y)$}
$$d(x,y)=\left\{\begin{array}{ll}0, & \mbox{ if $x=y$,} \\
2^{-n}, & \mbox{ if $n$ is the least such that $x(g_n)\neq y(g_n)$.}
\end{array}\right.$$
Notice that the action of $G$ on $k^G$ is continuous.

Since $1^G$ is trivial (it consists of a single point), $2^G$ is in some sense the ``smallest'' Bernoulli flow over $G$. As will be apparent, all of the questions we pursue in this paper are most restrictive (i.e. the most difficult to answer) in the setting of $2^G$. We therefore work nearly exclusively with $2^G$, however all of our results hold for $k^G$ ($k > 1$) by making obvious modifications to the proofs. While our main results were stated in terms of $k^G$ in the previous chapter, within the body of this paper we will only state our results in terms of $2^G$. Nevertheless, many definitions will be given in terms of $k^G$.

Although we will work primarily with Bernoulli flows, there are times when we wish to discuss more general dynamical systems. To accommodate this we introduce some notation and definitions in a more general setting. Let $G$ be a countable group and let $X$ be a compact metrizable space on which $G$ acts continuously. If $x \in X$, then the orbit of $x$ is denoted\index{$[x]$}
$$[x] = \{g \cdot x \: g \in G\}.$$
A \emph{subflow}\index{subflow}\index{Bernoulli flow!subflow}\index{Bernoulli $G$-flow!subflow} of $X$ is a closed subset of $X$ which is stable under the group action. If $A \subseteq X$, then we denote the closure of $A$ by $\overline{A}$. Notice that $\overline{[x]}$ is the smallest subflow of $X$ containing $x$.

We call $g \in G - \{1_G\}$ a \emph{period}\index{period} of $x \in X$ if $g \cdot x = x$. We call $x$ \emph{periodic}\index{periodic} if it has a period, and otherwise we call $x$ \emph{aperiodic}. As a word of caution, we point out that our use of the word periodic differs from conventional use (usually it means that the orbit of $x$ is finite). A subflow of $X$ is \emph{free}\index{free subflow} if it consists entirely of aperiodic points.

In this paper, aperiodic points and free subflows of $2^G$ play a particularly important role. We denote by $F(G)$\index{$F(G)$} the collection of all aperiodic points of $2^G$. $F(G)$ is called the \emph{free part}\index{free part}. It is a dense $G_\delta$ subset of $2^G$, is stable under the group action, and is closed in $2^G$ if and only if $G$ is finite. An important achievement of both this paper and the authors' previous paper \cite{GJS} is showing that while $F(G)$ is not compact (for infinite groups), it does display some compactness types of properties. Notice that in the case of $2^G$, a subflow is free if and only if it is contained in $F(G)$.

If $Y$ is another compact metrizable space on which $G$ acts continuously, then $X$ and $Y$ are \emph{topologically conjugate}\index{topologically conjugate} if there is a homeomorphism $\phi: X \rightarrow Y$ which \emph{commutes with the action of $G$}\index{commutes with the group action}, meaning $\phi(g \cdot x) = g \cdot \phi(x)$ for all $g \in G$ and $x \in X$.

Much of this paper is concentrated on constructing elements of $2^G$ with special properties. These functions $G \rightarrow 2 = \{0,1\}$ will be mostly defined by induction. In the middle of a construction we will be only working with partial functions from $G$ to $2$. This motivates the following notations and definitions. A {\it partial function}\index{partial function} $c$ from $G$ to $2$, denoted $c:G\rightharpoonup 2$, is a function $c:\dom(c)\to 2$ with $\dom(c)\subseteq G$. The set of all partial functions from $G$ to $2$ is denoted \index{$2^{\subseteq G}$}$2^{\subseteq G}$. We also denote the set of all partial functions from $G$ to $2$ with finite domains by \index{$2^{<G}$}$2^{<G}$. The action of $G$ on $2^G$ induces a natural action of $G$ on $2^{\subseteq G}$ as follows. For $g \in G$ and$ c \in 2^{\subseteq G}$, let $g \cdot c$ be the function with domain $g \cdot \dom(c)$ given by $(g \cdot c)(h) = c(g^{-1} h)$ for $h \in \dom(g \cdot c) = g \cdot \dom(c)$. Alternatively, $2^{\subseteq G}$ may be viewed simply as $3^G$. The bijection $\phi:2^{\subseteq G} \rightarrow 3^G$ is given by
$$\phi(c)(g) = \begin{cases}
c(g) & \text{if } g \in \dom(c) \\
2 & \text{otherwise}.
\end{cases}$$
It is easy to see that $\phi$ is a bijection and that $\phi$ commutes with the action of $G$. Therefore, $2^{\subseteq G}$ may be considered as $3^G$ if desired. In particular, this provides us with a nice compact metrizable topology on $2^{\subseteq G}$.

There still remains three more definitions which are very central to this paper. These definitions are introduced and discussed in each of the three next sections.

\section{$2$-colorings}

The related notions of $2$-coloring, $k$-coloring, and hyper aperiodic points are the most central notions of this paper. We first define these notions.

\begin{definition}\index{hyper aperiodic}
Let $G$ be a countable group, and let $X$ be a compact metrizable space on which $G$ acts continuously. A point $x \in X$ is called \emph{hyper aperiodic} if every point in $\overline{[x]}$ is aperiodic. Equivalently, $x$ is hyper aperiodic if $\overline{[x]}$ is a free subflow of $X$.
\end{definition}

In the particular context of Bernoulli flows we have the following specialized terminology.

\begin{definition}\index{$k$-coloring}\index{$2$-coloring}
Let $G$ be a countable group and $k\geq 2$ an integer. A point $x \in k^G$ is called a \emph{$k$-coloring} if it is hyper aperiodic. That is, $x$ is a $k$-coloring if every point in $\overline{[x]}$ is aperiodic. Equivalently, $x$ is a $k$-coloring if for every $s\in G$ with $s\neq 1_G$ there is a finite set $T\subseteq G$ such that
$$ \forall g\in G\ \exists t\in T\ x(gt)\neq x(gst).$$
\end{definition}

We will very shortly prove the equivalence of the two statements given in the previous definition.

In the context of Bernoulli flows, the terms $k$-coloring and hyper aperiodic are interchangeable. The term $k$-coloring, or to be more precise, $2$-coloring, is used with much greater frequency within the paper than the term hyper aperiodic. The reason is that $2$-coloring was the original term and the term hyper aperiodic was adopted much later on in order to discuss the concept in a more general setting. Still, a nice feature of of the term $k$-coloring is that it emphasizes the combinatorial characterization and is also reminiscent of the term coloring in graph theory as both notions roughly mean that nearby things look different.

Before proving that the two conditions in the previous definition are equivalent, we introduce one more definition.

\begin{definition}\index{block}
Let $G$ be a countable group, $k\geq 2$ an integer, $x\in k^G$,
and $s\in G$ with $s\neq 1_G$. We say that $x$ {\it blocks} $s$ if no element of $\overline{[x]}$ has period $s$. Equivalently, $x$ blocks $s$ if there is a finite set $T\subseteq G$ such that
$$ \forall g\in G\ \exists t\in T\ x(gt)\neq x(gst).$$
\end{definition}

Notice that $x \in k^G$ is a $k$-coloring if and only if $x$ blocks all non-identity $s \in G$.

The following lemma, which proves the equivalence of the combinatorial and dynamical conditions found in the previous two definitions, originally appeared, independently, in \cite{GJS} and \cite{GU}. For the convenience of the reader we include the proof below.

\begin{lem}[\cite{GJS}; Pestov, c.f. \cite{GU}]\label{lem:basiccoloringlemma}
Let $G$ be a countable group, $k \geq 2$ an integer, $x \in k^G$, and $s \in G$ with $s \neq 1_G$. Then $s \cdot y \neq y$ for all $y \in \overline{[x]}$ if and only if there is a finite set $T \subseteq G$ so that
$$\forall g \in G \ \exists t \in T \ x(g t) \neq x(g s t).$$
\end{lem}

\begin{proof}
$(\Leftarrow)$ Assume $x$ has the combinatorial property. Suppose $y\in\overline{[x]}$, that is,
there are $h_m\in G$ with $h_m\cdot x\to y$ as $m\to\infty$. We show that $s \cdot y \neq y$.
Assume not and suppose $s \cdot y = y$. Then by the continuity of the action we have that
$s^{-1}h_m\cdot x\to s^{-1}\cdot y = y$. Let $T\subseteq G$ be a finite set such that for any $g\in G$ there is $t\in T$ with $x(gt)\neq x(gst)$.
Let $n$ be large enough so that $T\subseteq \{g_0,\dots, g_n\}$, where $g_0, g_1, \ldots$ is the enumeration of $G$ used in defining the metric on $k^G$. Let $m\geq n$ be such that
$d(h_m\cdot x, y), d(s^{-1}h_m\cdot x, y)<2^{-n}$. Now fix $t\in T$ with $(h_m\cdot x)(t)=x(h_m^{-1}t)\neq x(h_m^{-1}st)=(s^{-1}h_m\cdot x)(t)$.
Then $y(t)=(h_m\cdot x)(t)\neq (s^{-1}h_m\cdot x)(t)=y(t)$, a contradiction.

$(\Rightarrow)$ Assume $s \cdot y \neq y$ for all $y \in \overline{[x]}$. Denote $C=\overline{[x]}$. Then for any $y\in C$, $s^{-1}\cdot y\neq y$,
and hence there is $t\in G$ with $(s^{-1}\cdot y)(t)\neq y(t)$. Define a function $\tau: C\to G$ by letting $\tau(y)=g_n$ where $n$ is the least so that
$(s^{-1}\cdot y)(g_n)\neq y(g_n)$. Then $\tau$ is a continuous function. Since $C$ is compact we get that $\tau(C)\subseteq G$ is finite. Let $T=\tau(C)$.
Then for any $g\in G$, there is $t\in T$ such that $x(gt)=(g^{-1}\cdot x)(t)\neq (s^{-1}g^{-1}\cdot x)(t)=x(gst)$. This proves that $x$ has the combinatorial property.
\end{proof}

\begin{cor}
Let $G$ be a countable group, $k \geq 2$ an integer, and $x \in k^G$. Then $x$ is hyper aperiodic, i.e. each $y \in \overline{[x]}$ is aperiodic, if and only if for every non-identity $s \in G$ there is a finite set $T \subseteq G$ such that
$$\forall g \in G \ \exists t \in T \ x(g t) \neq x(g s t).$$
\end{cor}

In view of the previous corollary, the problem of constructing free Bernoulli subflows is
reduced to the problem of constructing elements $x\in 2^G$ with a particularly combinatorial property. The combinatorial characterization of $2$-colorings is therefore vital to our constructions.

Under the dynamical definition of blocking, the following corollary is rather obvious. However, we will tend to focus mostly on the combinatorial definition of blocking, and from the combinatorial viewpoint the statement of the following corollary is not so expected. It is therefore worth pointing out for future reference.

\begin{cor}\label{cor:blockinglemma}
Let $G$ be a countable group, $k\geq 2$ an integer, $x\in k^G$, and $s\in G$ with $s\neq 1_G$. If $x$ blocks $g^{-1} s^n g$ for any integer $n$ and $g\in G$ with $s^n \neq 1_G$, then $x$ blocks $s$. In particular, if $x$ blocks $s^n$ for any integer $n$ with $s^n \neq 1_G$, then $x$ blocks $s$.
\end{cor}

\begin{proof}
We will use the dynamical characterization of blocking. If $x$ does not block $s$ then there is $z\in \overline{[x]}$ with $s\cdot z=z$. Then $g^{-1}s^ng\cdot (g^{-1}\cdot z)=g^{-1}\cdot z\in \overline{[x]}$. Hence $x$ does not block $g^{-1}s^ng$ for any $g\in G$ and $n$ with $s^n \neq 1_G$.
\end{proof}

We will also find it useful to discuss blockings for partial functions on $G$.

\begin{definition}\label{def:block}\index{block}
Let $G$ be a countable group, $c\in 2^{\subseteq G}$, and $s\in G$ with $s\neq 1_G$. We say that $c$ {\it blocks} $s$ if for any $x\in 2^G$ with $c\subseteq x$, $x$ blocks $s$.
\end{definition}

If $G$ is infinite no element of $2^{<G}$ can block any $s\in G$. There are, however, partial functions with coinfinite domains that can block all $s\in G$ with $s\neq 1_G$. As an example, suppose $y\in 2^{\mathbb Z}$ is a $2$-coloring on ${\mathbb Z}$. Consider a partial function $c: {\mathbb Z}\rightharpoonup 2$ defined by $c(2n)=y(n)$ for all $n\in {\mathbb Z}$. Then $c$ blocks all $s\in {\mathbb Z}$ with $s\neq 0$. This is because, for any $x\in 2^{\mathbb Z}$ with $c\subseteq x$, $x$ blocks all $2s$ with $s\neq 0$, and therefore by Corollary~\ref{cor:blockinglemma} $x$ blocks all $s$ with $s\neq 0$. It follows that any $x\in 2^G$ with $c\subseteq x$ is also a $2$-coloring on ${\mathbb Z}$.

Before closing this section we remark that in a countably infinite group, a single finite set cannot witness the blocking of all shifts.

\begin{lem}
Let $G$ be a countably infinite group. Then there are no finite sets $T\subseteq G$
such that for all $s\in G$ with $s\neq 1_G$, there is $t\in T$ with $x(gt)\neq x(gst)$.
\end{lem}

\begin{proof}
Assume not, and let $T\subseteq G$ be such a finite set. By induction we can define an infinite sequence $(h_n)$ of elements of $G$ so that $(h_nT)$ are pairwise disjoint. In fact, let $h_0\in G$ be arbitrary. In general, suppose $h_0,\dots, h_n$ have been defined so that $h_0T,\dots, h_nT$ are pairwise disjoint. Let $h_{n+1}\in G-(h_0T\cup\dots\cup h_nT)T^{-1}$. Since $G$ is infinite and $(h_0T\cup\dots\cup h_nT)T^{-1}$ is finite, such $h_{n+1}$ exists. Then $h_{n+1}T$ is disjoint from $h_0T,\dots, h_nT$. Now consider the partial functions $c_n\in 2^{\subseteq G}$ with $\dom(c_n)=T$ defined by $c_n(t)=x(h_nt)$. Since there are only finitely many partial functions with domain $T$, there are $n\neq m$ such that $c_n=c_m$. Thus $x(h_nt)=x(h_mt)$ for all $t\in T$. Thus if we let $s=h_m^{-1}h_n$, $T$ fails to witness that $x$ blocks $s$, a contradiction.
\end{proof}

\section{Orthogonality}

The notion of orthogonality is another central notion to this paper. On the one hand, a pair of points being orthogonal says that they are distinct from one another in a strong sense, and on the other hand the notion of orthogonality carries along with it some nice properties which we will make use of frequently.

\begin{definition}\index{orthogonal}\index{$x\, \bot\, y$}
Let $G$ be a countable group, let $X$ be a compact metrizable space on which $G$ acts continuously, and let $x_0, x_1 \in X$. We say that $x_0$ and $x_1$ are {\em orthogonal}, denoted $x_0\,\bot\, x_1$, if $\overline{[x_0]}$ and $\overline{[x_1]}$ are disjoint. If $X$ is a Bernoulli flow, then this is equivalent to the existence of a finite set $T\subseteq G$ such that
$$\forall g_0, g_1 \in G\ \exists t\in T\ x_0(g_0 t)\neq x_1(g_1 t).$$
\end{definition}

The following lemma implies that within the context of Bernoulli flows, the two conditions given in the previous definition are equivalent.

\begin{lem}\label{lem:orthogonallemma} Let $G$ be a countable group, $k\geq 2$ an integer, and $x_0, x_1\in k^G$. Then $\overline{[x_0]}$ and $\overline{[x_1]}$ are disjoint if and only if there is a finite set $T \subseteq G$ such that
$$\forall g_0, g_1 \in G\ \exists t\in T\ x_0(g_0 t)\neq x_1(g_1 t).$$
\end{lem}

\begin{proof}
$(\Rightarrow)$ Suppose $\overline{[x_0]}\cap \overline{[x_1]}=\emptyset$. Since they are both compact it follows that there is some $\delta>0$
such that for any $z_0\in \overline{[x_0]}$ and $z_1\in \overline{[x_1]}$, $d(z_0,z_1)\geq \delta$. Recall that the metric on $k^G$ is defined in terms of an enumeration $1_G=g_0, g_1,\dots$ of $G$. Let $n$ be large enough such that $\delta\geq 2^{-n}$.
Then in particular for any $g_0, g_1 \in G$, $d(g_0^{-1} \cdot x_0, g_1^{-1} \cdot x_1)\geq 2^{-n}$. This implies that there is $t\in \{g_0,\dots, g_n\}$ such that
$x_0(g_0 t) = (g_0^{-1} \cdot x_0)(t) \neq (g_1^{-1} \cdot x_1)(t) = x_1(g_1 t)$.

$(\Leftarrow)$ Conversely, suppose that $T$ is a finite set such that for every pair $g_0, g_1 \in G$ there is $t \in T$ with $x_0(g_0 t) \neq x_1(g_1 t)$. Let $n$ be large enough such that $T\subseteq \{g_0,\dots, g_n\}$. Then for any $y_0\in [x_0]$ and $y_1\in [x_1]$, there is $t\in T$ such that $y_0(t)\neq y_1(t)$, and thus $d(y_0,y_1)\geq 2^{-n}$. It follows that for any $z_0\in \overline{[x_0]}$ and $z_1\in \overline{[x_1]}$, $d(z_0,z_1)\geq 2^{-n}$,
and therefore $\overline{[x_0]}\cap \overline{[x_1]}=\emptyset$.
\end{proof}

We will frequently work with infinite sets of pairwise orthogonal elements. In this situation we remark that the pairwise orthogonality cannot be witnessed by a single finite set.

\begin{lem} Let $G$ be a countable group, $k\geq 2$ an integer, and $x_0, x_1, \dots$ be infinitely many pairwise orthogonal elements of $k^G$. Then there are no finite sets $T\subseteq G$ such that for any $n\neq m$ there is $t\in T$ with $x_n(t)\neq x_m(t)$.
\end{lem}

\begin{proof}
Assume not, and let $T\subseteq G$ be such a finite set. Consider the partial functions $c_n\in 2^{\subseteq G}$
with $\dom(c_n)=T$ defined by $c_n(t)=x_n(t)$. Since there are only finitely many partial functions with domain $T$, there are $n\neq m$ such that $c_n=c_m$. Thus for all $t\in T$, $x_n(t)=x_m(t)$, a contradiction.
\end{proof}

$2$-colorings were constructed on every countable group in \cite{GJS}, and they are also constructed in this paper (with a much improved construction). The known methods for constructing $2$-colorings on general countable groups involve purely geometric and combinatorial methods, and the constructions are rather long and technical. There is therefore motivation to develop simpler constructions for more restricted classes of groups. We do this in Chapter \ref{chap:basicconstructions} and Section \ref{sec:abmar}. The construction in Section \ref{sec:abmar} uses geometric methods, as in the general setting. However, in Chapter \ref{chap:basicconstructions} we construct $2$-colorings on all solvable groups, all free groups, and all residually finite groups by using algebraic methods. The notion of orthogonality plays a key role in these constructions. The following definition will be used in Chapter \ref{chap:basicconstructions}.

\begin{definition}\index{$(\lambda,k)$-coloring property}\index{coloring property}\index{$2$-coloring!coloring property}
Let $G$ be a countable group, $k\geq 2$ an integer, and $\lambda\geq 1$ a cardinal number.
$G$ is said to have the {\em $(\lambda,k)$-coloring property}, if there exist $\lambda$ many pairwise orthogonal $k$-colorings of $G$.
$G$ is said to
have the {\em coloring property} if $G$ has the $(1,2)$-coloring property, i.e.,
there is a $2$-coloring on $G$.
\end{definition}

We point out that it has already been proven in \cite{GJS} that every countable group has the coloring property, i.e. admits a $2$-coloring, and moreover that every countably infinite group has the $(2^{\aleph_0},2)$-coloring property. So for any cardinal $\lambda\leq 2^{\aleph_0}$ and integer $k\geq 2$, every countably infinite group has the $(\lambda,k)$-coloring property. The above definition is therefore rather trivial, but nonetheless it will be useful for studying algebraic constructions of $2$-colorings in Chapter \ref{chap:basicconstructions}.

It is easy that finite groups have the coloring property. However, it is not clear how many orthogonal $k$-colorings there can be.

\begin{lem}\label{lem:finite} Every finite group has the coloring property. Every finite group with at least 3 elements has the $(2,2)$-coloring property. The two element group ${\mathbb Z}_2$ does not have the
$(2,2)$-coloring property.
\end{lem}

\begin{proof}
For any finite group $G$ let $c(1_G)=0$ and $c(g)=1$ for all $g\neq 1_G$. Then $c$ is a $2$-coloring on $G$.
Let $\bar{c}(g)=1-c(g)$ for all $g\in G$. If $|G|>2$ then $c$ and $\bar{c}$ are both $2$-colorings and $c\,\bot\, \bar{c}$. ${\mathbb Z}_2$ has only two $2$-colorings, but they are
not orthogonal (they are in the same orbit).
\end{proof}

\section{Minimality}

We now discuss the classical notion of minimality. We again see that in the context of Bernoulli flows this dynamical notion has a combinatorial characterization.

\begin{definition}\index{minimal}
Let $G$ be a countable group, and let $X$ be a compact metrizable space on which $G$ acts continuously. A subflow $Y \subseteq X$ is \emph{minimal} if $\overline{[y]} = Y$ for all $y \in Y$. A point $x \in X$ is \emph{minimal} if $\overline{[x]}$ is minimal. If $X$ is a Bernoulli flow, then $x \in X$ is minimal if it satisfies the following: for every finite $A \subseteq G$ there exists a finite $T \subseteq G$ such that
$$\forall g\in G\ \exists t\in T\ \forall a\in A\ x(g t a) = x(a).$$
\end{definition}

As a word of caution, we point out that our definition of minimality of a point is not standard; it is relatively common to call $x$ almost-periodic if $\overline{[x]}$ is minimal. However, for us almost-periodic will have a different meaning.

In a moment we will prove that in the context of Bernoulli flows the two stated conditions for minimality of a point are equivalent. However, first we prove the standard fact that minimal subflows always exist. We first give a standard well-known
argument for their existence assuming $\ac$.

\begin{lem} [$\zfc$] \label{EXIST MIN}
Let $X$ be a compact Hausdorff space, and let $G$ be a group which acts on $X$. Then $X$ contains a minimal subflow.
\end{lem}

\begin{proof}
Consider the collection of all subflows of $X$, ordered by reverse inclusion. This collection is nonempty because it contains $X$ itself. If $(X_n)_{n \in \N}$ is a chain, then each $X_n$ is a subflow and hence is closed, thus compact. So $Y= \bigcap_{n \in \N} X_n$ is a nonempty compact set. Since $X$ is Hausdorff, $Y$ is closed. $Y$ is also clearly $G$-invariant. The claim now follows by applying Zorn's Lemma.
\end{proof}

In fact, $\ac$ is not needed to prove Lemma~\ref{EXIST MIN},
at least in the case when $X$ is Polish. We are not sure if this has been observed before so
we give the proof in the following lemma.

\begin{lem} [$\zf$]\label{existmin2}
Let $X$ be a compact Polish space on which the group $G$ acts continuously. Then there is a minimal
subflow.
\end{lem}

\begin{proof}
Let $\{ U_n\}_{n \in \N}$ enumerate a base for $X$.
Let $F(X)$ be the standard Borel space of closed (so compact) non-empty subsets of $X$ with the usual Effros Borel
structure (which, since $X$ is compact, is generated by the Vietoris topology on $F(X)$). By the Borel selection theorem in
descriptive set theory (c.f. \cite[Theorem 12.13]{KechrisBook} or \cite[Theorem 1.4.6]{Gao}) there is a Borel function $s \colon F(X)\to X$ which is a
selector, that is, $s(F)\in F$ for all $F\in F(X)$. We define inductively closed invariant sets
$F_\alpha$ satisfying $F_\beta \subsetneq F_\alpha$ for all $\alpha<\beta$. Let $F_0=X$.
If $\alpha$ is a limit ordinal, let $F_\alpha=\bigcap_{\beta<\alpha} F_\beta$ (which is non-empty by compactness). If $\alpha=\beta+1$, stop
the construction if $F_\beta$ is minimal. Otherwise, let $n$ be least such that $F_\beta \cap U_n \neq \emptyset$
and $F_\beta-G\cdot U_n=F_\beta-\bigcup_{g\in G}g\cdot U_n\neq \emptyset$. Such an $n$ clearly exists if
$F_\beta$ is an invariant but not  minimal closed set. Let $A=F_\beta-U_n$. Let $F=\{ x \in A \colon
\overline{[x]} \subseteq A\}$. Note that $X-F=(X-A)\cup \{ x \colon \exists g \in G\ g \cdot x \in X-A\}$. Since
$G$ acts continuously on $X$, this shows that $X-F$ is open, so $F$ is closed. Let $x_\alpha=s(F)$, and let
$F_\alpha= \overline{[x_\alpha]}$. Clearly $F_\alpha$ is a closed invariant set which is properly contained
in $F_\beta$. The above transfinite recursion defining the $F_\alpha$ is clearly done in $\zf$. The construction must
stops at some ordinal $\theta$, and we are done as $F_\theta$ is then a minimal subflow.
\end{proof}

In fact, using a little more descriptive set theory one can prove more. We state this in the next lemma.

\begin{lem} [$\zf$] \label{existmin3}
Let $X$ be a compact Polish space and let $G$ be a Polish group acting in a Borel way on $X$. Then there is a minimal subflow.
\end{lem}

\begin{proof}
We proceed as in Lemma~\ref{existmin2}, defining by transfinite recursion a sequence $F_\alpha$ of (non-empty)
closed, invariant subsets of $X$. At limit stages we again take intersections. If $\alpha=\beta+1$ and $F_\beta$ is not
a minimal flow, we again let $n$ be least such that $F_\beta\cap U_n\neq\emptyset$ and $F_\beta-G\cdot U_n\neq\emptyset$. Again let $A=F_\beta-U_n$ and

\begin{equation*}
F= \{ x \colon \overline{[x]} \subseteq A\} =
\{ x \colon \forall g \in G\ (g \cdot x \in A)\}.
\end{equation*}
Note that $F$ is a non-empty $\boldsymbol{\Pi}^1_1$ set, using that the action of $G$ on $X$ is Borel.
In fact, consider the relation $R \subseteq F(X) \times X$ defined by
$R(A,x) \leftrightarrow \forall g \in G\ (g \cdot x \in A)$. This is a $\boldsymbol{\Pi}^1_1$ relation in the Polish space
$F(X)\times X$. It is a theorem of $\zf$ that $\boldsymbol{\Pi}^1_1$ subsets of products of Polish spaces
admit $\boldsymbol{\Pi}^1_1$ uniformizations (recall a uniformization $R'$ of a relation $R \subseteq X \times Y$
means $R' \subseteq R$, $\dom(R')=\dom(R)$, and $R'$ is the graph of a (partial) function). Here we do not care about the complexity
of the uniformization, only that the relation $R$ has a uniformizing function, call it $s$, provably in $\zf$. The proof
then finishes as in Lemma~\ref{existmin2}, letting $x_\alpha=s(F)$ and $F_\alpha=\overline{[x_\alpha]}$ as before.
\end{proof}

Now we prove the equivalence of the dynamical and combinatorial characterizations of minimality (of a point) in the context of Bernoulli flows.

\begin{lem}\label{lem:minimallemma}
Let $G$ be a countable group and $x\in k^G$. Then $\overline{[x]}$ is a minimal subflow
iff for every finite $A \subseteq G$ there exists a finite $T \subseteq G$ such that
$$\forall g\in G\ \exists t\in T\ \forall a\in A\ x(g t a) = x(a).$$
\end{lem}

\begin{proof}
($\Rightarrow$) Assume $\overline{[x]}$ is a minimal subflow. Let $A \subseteq G$ be arbitrary but finite, and let $n$ be large enough such that $A \subseteq \{g_0, g_1, \ldots, g_n\}$, where $g_0, g_1, \ldots$ is the enumeration of $G$ used in defining the metric on $k^G$. For every $z \in \overline{[x]}$ there exists $h \in G$ with $d(h \cdot z, x) < 2^{-n}$ since $[z]$ is dense in $\overline{[x]}$. Define $\phi(z)=g_m$, where $m$ is the least integer such that $d(g_m^{-1} \cdot z, x) < 2^{-n}$. Then $\phi: \overline{[x]}\to G$ is continuous. Since $\overline{[x]}$ is compact, it follows that $\phi(\overline{[x]})$ is finite. Set $T = \phi(\overline{[x]})$. In particular, we have that for any $g \in G$ there is $t \in T$ with $d(t^{-1} \cdot (g^{-1} \cdot x), x) < 2^{-n}$. Therefore, for all $a \in A$, $\ x(g t a) = (t^{-1} \cdot g^{-1} \cdot x)(a) = x(a)$.

($\Leftarrow$) Now assume $x$ has the stated combinatorial property. Let $z \in \overline{[x]}$. It suffices to show that $x\in \overline{[z]}$. For this we fix an arbitrary $\epsilon > 0$ and show that $d([z],x)<\epsilon$. Then since $\epsilon$ is arbitrary, we would actually have $d([z],x)=0$ and so $x\in\overline{[z]}$. For this let $n$ be large enough such that $2^{-n} < \epsilon$, and set $A = \{g_0, g_1, \ldots, g_n\}$. By our assumption, there is a finite $T \subseteq G$ such that for all $g \in G$ there is $t \in T$ with $x(g t a) = x(a)$ for all $a \in A$. Let $h_i$ be a sequence in $G$ with $h_i \cdot x \rightarrow z$ as $i \rightarrow \infty$. Let $m$ be large enough such that $T A \subseteq \{g_0, g_1, \ldots, g_m\}$, and fix $i$ with $d(h_i \cdot x, z) < 2^{-m}$. Then for some $t \in T$, $\ x(h_i^{-1} t a) = x(a)$ for all $a \in A$.  Thus $z(t a) = (h_i \cdot x)(t a) = x(h_i^{-1} t a) = x(a)$ for all $a \in A$. This implies that $d([z],x)<2^{-n}<\epsilon$, as promised.
\end{proof}

The combinatorial characterization of minimality allows us to explicitly construct minimal elements of $2^G$ without appealing to Zorn's lemma. It also has the following immediate corollary about the descriptive complexity of the set of all minimal elements.

\begin{cor} Let $G$ be a countable group. The set of all minimal elements of $2^G$ is ${\bf\Pi}^0_3$.
\end{cor}

We also note the following basic fact, which provides a useful way to obtain orthogonal elements through minimality.

\begin{lem}
Let $G$ be a countable group and $x\in 2^G$. If $y\in 2^G-\overline{[x]}$ is minimal then $y\,\bot\,x$.
\end{lem}

\begin{proof}
By Lemma~\ref{lem:orthogonallemma} it suffices to show that $\overline{[y]}\cap\overline{[x]}=\emptyset$.
Assume not, and let $z\in \overline{[y]}\cap\overline{[x]}$. Then $\overline{[z]}\subseteq \overline{[x]}$. Moreover, by minimality of $y$, $\overline{[z]}=\overline{[y]}$. Thus $y\in \overline{[y]}=\overline{[z]}\subseteq \overline{[x]}$, contradicting our assumption.
\end{proof}

\section{\label{sec:2.5}Strengthening and weakening of $2$-colorings}

In this section we introduce some natural strengthening and weakening of $2$-colorings which will be further studied in later chapters. We first give their definitions.

\begin{definition}\index{almost equal}\index{$=^*$} Let $G$ be a countable group and $x,y\in 2^G$. We call $x$ and $y$ \emph{almost equal}, denoted $x =^* y$, if the set $\{g\in G\,:\, x(g)\neq y(g)\}$ is finite.
\end{definition}

\begin{definition} Let $G$ be a countable group and $x\in 2^G$.
\begin{enumerate}
\item[(1)] For $s\in G$ with $s\neq 1_G$, we say that $x$ {\it nearly blocks} $s$ if there are finite sets $S, T\subseteq G$ such that
$$ \forall g\not\in S\ \exists t\in T\ x(gt)\neq x(gst). $$
$x$ is called a {\it near $2$-coloring}
on $G$ if $x$ nearly blocks $s$ for all $s\in G$ with $s\neq 1_G$.
\item[(2)] $x$ is called an {\it almost $2$-coloring} on $G$ if there is a $2$-coloring $y$ on $G$ such that $x=^*y$.
\item[(3)] For $s\in G$ with $s\neq 1_G$, we say that $x$ {\it strongly blocks} $s$ if $x$ blocks $s$ and there are infinitely many $g\in G$ such that $x(sg)\neq x(g)$. $x$ is called a {\it strong $2$-coloring} on $G$ if $x$ strongly blocks $s$ for all $s\in G$ with $s\neq 1_G$,
\end{enumerate}
\end{definition}
\index{near $2$-coloring}\index{$2$-coloring!near $2$-coloring}
\index{nearly block}\index{block!nearly block}
\index{almost $2$-coloring}\index{$2$-coloring!almost $2$-coloring}
\index{strong $2$-coloring}\index{$2$-coloring!strong $2$-coloring}
\index{strongly block}\index{block!strongly block}

We first mention that there is an equivalent dynamical characterization for near $2$-colorings. We remind the reader that if $A$ is a subset of a topological space $X$ and $x \in X$, then $x$ is said to be a limit point of $A$ if $x$ lies in the closure of $A - \{x\}$.

\begin{lem} \label{NEAR EQUIV}
Let $G$ be a countable group and let $x \in 2^G$. The following are equivalent:
\begin{enumerate}
\item[\rm (i)] $x$ is a near $2$-coloring;
\item[\rm (ii)] for every non-identity $s \in G$ there are finite sets $S, T \subseteq G$ so that for all $g \in G - S$ there is $t \in T$ with $x(g t) \neq x(g s t)$;
\item[\rm (iii)] every limit point of $[x]$ is aperiodic.
\end{enumerate}
\end{lem}

\begin{proof}
The equivalence of (i) and (ii) is by definition.

(ii) $\Rightarrow$ (iii). Let $y \in \overline{[x]}$ be a limit point of $[x]$. Then there is a non-repeating sequence $(g_n)_{n \in \N}$ of group elements of $G$ with $y = \lim g_n \cdot x$. Fix a non-identity $s \in G$. It suffices to show that $s^{-1} \cdot y \neq y$. Let $S, T \subseteq G$ be finite and such that for all $g \in G - S$ there is $t \in T$ with $x(g t) \neq x(g s t)$. Let $m \in \N$ be such that for all $n \geq m$ and all $t \in T$
$$y(t) = (g_n \cdot x)(t) \text{ and } y(s t) = (g_n \cdot x)(s t).$$
Since $(g_n)_{n \in \N}$ is non-repeating and $S$ is finite, there is $n \geq m$ with $g_n^{-1} \not\in S$. Let $t \in T$ be such that $x(g_n^{-1} t) \neq x(g_n^{-1} s t)$. Then we have
$$y(t) = (g_n \cdot x)(t) = x(g_n^{-1} t) \neq x(g_n^{-1} s t) = (g_n \cdot x)(s t) = y(s t) = (s^{-1} \cdot y)(t).$$
Therefore $s^{-1} \cdot y \neq y$.

(iii) $\Rightarrow$ (ii). Fix a non-identity $s \in G$. We must find sets $S$ and $T$ satisfying (ii). Let $C$ be the set of limit points of $[x]$. Then $C$ is closed, compact, and nonempty (since $2^G$ is compact). Let $g_0, g_1, g_2, \ldots$ be the enumeration of $G$ used in defining the metric $d$ on $2^G$. Define $\phi: C \rightarrow \N$ by letting $\phi(y)$ be the least $n$ with $y(g_n) \neq y(s g_n) = (s^{-1} \cdot y)(g_n)$. Then $\phi$ is continuous. Since $C$ is compact, $\phi$ has finite image. So there is a finite $T \subseteq G$ containing the image of $\phi$. So for every $y \in C$ there is $t \in T$ with $y(t) \neq y(s t)$. Let $S$ be the set of $g \in G$ for which $x(g t) = x(g s t)$ for all $t \in T$. Towards a contradiction, suppose $S$ is infinite. By compactness of $2^G$, we can pick a non-repeating sequence $(g_n)_{n \in \N}$ of elements of $S$ such that $g_n^{-1} \cdot x$ converges to some $y \in 2^G$. Now $y \in C$ and for $t \in T$ and sufficiently large $n \in \N$ we have
$$y(t) = (g_n^{-1} \cdot x)(t) = x(g_n t) = x(g_n s t) = (g_n^{-1} \cdot x)(s t) = y(s t).$$
So $y(t) = y(s t)$ for all $t \in T$, a contradiction. We conclude that $S$ is finite.
\end{proof}

\begin{lem}\label{lem:almostcoloringlemma} Let $G$ be a countable group.
Then the following hold:
\begin{enumerate}
\item[(a)] Every strong $2$-coloring on $G$ is a $2$-coloring on $G$.
\item[(b)] Every $2$-coloring on $G$ is an almost $2$-coloring on $G$.
\item[(c)] Every almost $2$-coloring on $G$ is a near $2$-coloring on $G$.
\item[(d)] Every aperiodic near $2$-coloring on $G$ is a $2$-coloring on $G$.
\item[(e)] $x$ is a strong $2$-coloring on $G$ iff for all $y=^*x$, $y$ is a $2$-coloring on $G$.
\end{enumerate}
\end{lem}

\begin{proof} (a) and (b) are immediately obvious and (c) follows from the previous lemma. We only show (d) and (e).

For (d) assume that $x$ is an aperiodic near $2$-coloring on $G$. By the previous lemma all of the limit points of $[x]$ are aperiodic. Since $\overline{[x]}$ is the union of $[x]$ with the limit points of $[x]$, it follows that $\overline{[x]}$ is free (i.e. consists entirely of aperiodic points). Thus $x$ is a $2$-coloring.

Now for (e) we first show $(\Rightarrow)$. Assume $x$ is a strong $2$-coloring on $G$. Let $y=^*x$ and $A=\{g\in G\,:\, x(g)\neq y(g)\}$.
Then $y$ is an almost $2$-coloring, and in particular a near $2$-coloring. By (d), it suffices to show that
$y$ is aperiodic. Let $s\in G$ with $s\neq 1_G$. Let $g\in G-(A\cup s^{-1}A)$ be such that $x(sg)\neq x(g)$. Then $g, sg\not\in A$, and $y(g)=x(g)\neq x(sg)=y(sg)$. Hence $y$ is aperiodic.

For $(\Leftarrow)$, assume that for all $y=^* x$, $y$ is a $2$-coloring on $G$. In particular $x$ is a
$2$-coloring on $G$. We show that $x$ strongly blocks $s$ for all $s\in G$ with $s\neq 1_G$. Fix such an $s$. Consider two cases.

Case 1: $s$ has infinite order, i.e., $\langle s\rangle$ is infinite. Let $T\subseteq G$ be a finite set witnessing that $x$ blocks $s$. Since $TT^{-1}\cap \langle s\rangle$ is finite, there is $m\in\N$ such that for all $k$ with $|k|\geq m$, $s^k\not\in TT^{-1}$. Fix such an $m$. Then we have that for all distinct $n,k\in\N$, $s^{nm}T\cap s^{km}T=\emptyset$.
By blocking we have that for all $n\in\N$ there is $t_n\in T$ such that $x(s^{nm}t_n)\neq x(s^{nm}st_n)$. Thus $x(s^{nm}t_n)\neq x(ss^{nm}t_n)$ for all $n\in\N$. Since the set $\{s^{nm}t_n\,:\, n\in\N\}$ is infinite, we have that $x$ strongly blocks $s$.

Case 2: $s$ has finite order. Toward a contradiction, assume that $A=\{ t\in G\,:\, x(t)\neq x(st)\}$ is finite. Then for all $t\not\in A$, $x(st)=x(t)$. Now define $y=^*x$ so that $\{g\in G\,:\, y(g)\neq x(g)\}\subseteq \langle s\rangle A$ and $y$ is constant on $\langle s\rangle A$. Then $y(st)=y(t)$ for all $t\in G$. Thus $y$ is not a $2$-coloring, contradiction.
\end{proof}

Thus we have the following implications:
$$\begin{array}{c}
 \mbox{strong $2$-coloring} \\
  \left\Downarrow \displaystyle\frac{{}}{{}}\mbox{(a)}\right.  \\
\ \ \ \ \ \ \ \ \ \ \ \ \ \ \ \ \ \ \ \ \ \ \ \ \ \ \ \ \ \ \ \ \ \ \ \ \mbox{$2$-coloring}\Longleftrightarrow \mbox{aperiodic near $2$-coloring } \\
\left\Downarrow \displaystyle\frac{{}}{{}}\mbox{(b)}\right. \\
\mbox{almost $2$-coloring}  \\
\left\Downarrow \displaystyle\frac{{}}{{}}\mbox{(c)}\right. \\
\mbox{near $2$-coloring}
\end{array}
$$
We will show in Section~\ref{SEC ACP} that the converses of (a) and (b) are false. On the other hand, in Section \ref{SECT ALMOST EQUAL} we will prove that the converse of (c) is true. We are particularly interested in the following property for
countable groups.

\begin{definition} A countable group $G$ is said to have the {\it almost $2$-coloring property} ({\it ACP})
if every almost $2$-coloring on $G$ is a $2$-coloring on $G$.
\end{definition}\index{almost $2$-coloring!almost $2$-coloring property}\index{ACP}\index{$2$-coloring!almost $2$-coloring property}

The following lemma is easy to prove.

\begin{lem}\label{lem:ACPsimplechar} Let $G$ be a countable group. Then the following are equivalent:
\begin{enumerate}
\item[(i)] $G$ has the ACP;
\item[(ii)] Every $2$-coloring on $G$ is a strong $2$-coloring on $G$.
\item[(iii)] Every almost $2$-coloring on $G$ is a strong $2$-coloring on $G$.
\end{enumerate}
\end{lem}

\begin{proof}
It is immediate that (iii) is equivalent
to the combination of (i) and (ii), thus it suffices to show the equivalence of (i) and (ii).
For (i)$\Rightarrow$(ii), suppose $G$ has the ACP. Let $x$ be a $2$-coloring on $G$. Let $y=^*x$. Then $y$ is an almost $2$-coloring. By the ACP $y$ is a $2$-coloring. Thus we have shown that every $y=^*x$ is a $2$-coloring on $G$. By Lemma~\ref{lem:almostcoloringlemma} (e) $x$ is a strong $2$-coloring on $G$. The converse (ii)$\Rightarrow$(i) is similar.
\end{proof}

We consider the notion of centralizer in a group $G$ in the following proposition. For $g\in G$, the {\it centralizer} of $g$ in $G$ is defined as\index{$\Z_G(g)$}
$$ \Z_G(g)=\{h\in G\,:\, gh=hg\}. $$

\begin{prop} \label{PROP HALF ACP}
Let $G$ be a countably infinite group. If for every $1_G \neq u \in G$ there is $1_G \neq v \in \langle u \rangle$ with $|\Z_G(v)| = \infty$, then every near $2$-coloring on $G$ is a $2$-coloring on $G$. In particular, a group $G$ has the ACP if for every $1_G \neq u \in G$ there is $1_G \neq v \in \langle u \rangle$ with $|\Z_G(v)| = \infty$.
\end{prop}

\begin{proof}
Let $G$ be a group with the stated property. Let $x \in 2^G$ be a near $2$-coloring. We will show that $x$ is a $2$-coloring by showing that $x$ is aperiodic and then applying clause (d) of Lemma \ref{lem:almostcoloringlemma}. Towards a contradiction, suppose $x$ is not aperiodic. So there is $1_G \neq u \in G$ with $u \cdot x = x$. Let $1_G \neq v \in \langle u \rangle$ be such that $|\Z_G(v)| = \infty$. Notice $v \cdot x = x$. Let $g_1, g_2, \ldots$ be any non-repeating sequence of elements in $\Z_G(v)$. By compactness of $2^G$ and by passing to a subsequence if necessary, we may suppose that $(g_n \cdot x)_{n \in \N}$ is a convergent sequence. Set $y = \lim g_n \cdot x$. Since each $g_n \in \Z_G(v)$, we have
$$v \cdot y = v \cdot ( \lim g_n \cdot x) = \lim v \cdot g_n \cdot x = \lim g_n \cdot v \cdot x = \lim g_n \cdot x = y.$$
Thus $y$ is a limit point of $[x]$ and is periodic. This contradicts Lemma \ref{NEAR EQUIV}. We conclude that $x$ must be aperiodic and is thus a $2$-coloring.
\end{proof}

The condition above is in fact both necessary and sufficient for $G$ to have the ACP. The proof of necessity will be provided in Section \ref{SEC ACP}. From the previous proposition we can easily list a few classes of groups which have the ACP. Recall the following definition of FC groups.

\begin{definition} \label{DEFN FCGROUP}
If $G$ is a group in which every conjugacy class is finite then $G$ is
called an {\it FC group}. Specifically, an FC group $G$ is a group such that
for all $g\in G$, $\{ h g h^{-1}\,:\, h\in G\}$ is finite.
\end{definition}\index{FC group}

\begin{cor}\label{cor:ACPexamples}
Let $G$ be a countably infinite group. Then $G$ has the ACP if any of the following is true:
\begin{enumerate}
\item[\rm (i)] every non-identity element of $G$ has infinite order;
\item[\rm (ii)] $G$ is a free abelian or free non-abelian group;
\item[\rm (iii)] $G$ is nilpotent;
\item[\rm (iv)] $G$ is an FC group.
\end{enumerate}
\end{cor}

\begin{proof}
(i). For any $1_G \neq v \in G$ we have that $\langle v \rangle \subseteq \Z_G(v)$. So if every non-identity group element has infinite order, then every group element has infinite centralizer. So by the previous proposition $G$ has the ACP.

(ii). This follows immediately from (i).

(iii). Set $G_0 = G$ and in general define $G_{n+1} = [G, G_n]$. Since $G$ is nilpotent, $G_n$ is trivial for sufficiently large $n$. Let $n$ be such that $G_n$ is infinite and $G_{n+1}$ is finite. Fix $1_G \neq v \in G$. We have that for all $g \in G_n$, $[g, v] \in G_{n+1}$. If $g, h \in G_n$ satisfy $[g, v] = [h, v]$ then
$$ g^{-1} v^{-1}g v = [g, v] = [h, v] =  h^{-1} v^{-1}h v$$
so
$$h g^{-1} v g h^{-1} = v$$
and therefore $h g^{-1} \in \Z_G(v)$. Since $G_n$ is infinite and $G_{n+1}$ is finite, it immediately follows that infinitely many elements of $G_n$ lie in $\Z_G(v)$. By the previous proposition $G$ has the ACP.

(iv). Fix $1_G \neq v \in G$. If $g, h \in G$ satisfy $g v g^{-1} = h v h^{-1}$ then it follows $h^{-1} g \in \Z_G(v)$. Since $G$ is infinite and the conjugacy class of $v$ is finite, it follows that $\Z_G(v)$ is infinite. So by the previous proposition $G$ has the ACP.
\end{proof}

In Section \ref{SEC ACP}, we will show that solvable, polycyclic, and virtually abelian groups in general do not have the ACP.
In contrast, in Section \ref{SEC UNI 2 COL} we will show that every countably infinite group has a strong $2$-coloring.

\section{\label{sec:2.6}Other variations of $2$-colorings}

In this section we introduce some further concepts related to $2$-colorings. These will not be
our main subjects of investigation. However, we will note from time to time that our methods
for constructing $2$-colorings can also be applied to obtain these variations.

First we consider the dual notion of a {\it right action} of $G$ on $2^G$:\index{right action}
$$ (g\cdot x)(h)=x(hg). $$
This induces a dual version of all the concepts that we have defined and considered throughout this
chapter. Consequently we obtain the notion of right $2$-colorings.

\begin{definition} Let $G$ be a countable group. An element $x\in 2^G$ is called a {\it right $2$-coloring} if for any $s\in G$ with $s\neq 1_G$ there is a finite set $T\subseteq G$
such that
$$ \forall g\in G\ \exists t\in T\ x(tg)\neq x(tsg). $$
\end{definition}\index{right action!right $2$-coloring}\index{$2$-coloring!right $2$-coloring}

It is now natural to ask whether the concepts of $2$-colorings and right $2$-colorings can be combined.

\begin{definition} Let $G$ be a countable group. An element $x\in 2^G$ is called a {\it two-sided $2$-coloring} if it is both a $2$-coloring and a right $2$-coloring.
\end{definition}\index{right action!two-sided $2$-coloring}\index{two-sided $2$-coloring}\index{$2$-coloring!two-sided $2$-coloring}

Of course, for abelian groups $2$-colorings and right $2$-colorings coincide, hence also with two-sided $2$-colorings. For non-abelian groups, very little is known for two-sided $2$-colorings. In Section~\ref{sec:3.3}, we will give for non-abelian free groups examples of two-sided $2$-colorings and of $2$-colorings that are not
right (or two-sided) $2$-colorings.

Next we note that the definition of $2$-colorings does not explicitly mention the inverse operation in
a group, and therefore can be similarly defined for any semigroup.

\begin{definition} Let $S$ be a countable semigroup. An element $2^S$ is called a {\it $2$-coloring} on the semigroup $S$ if for any $s\in S$ there is a finite set $T\subseteq S$
such that
$$ \forall g\in S\ [\,g\neq gs\rightarrow \exists t\in T\ x(gt)\neq x(gst)\,]. $$
\end{definition}\index{$2$-coloring!$2$-coloring on a semigroup}

We will not systematically explore $2$-colorings on semigroups in this paper. Instead, we will just consider
some $2$-colorings on ${\mathbb N}$. These are intrinsically related to $2$-colorings on ${\mathbb Z}$.

\begin{definition} A $2$-coloring $x\in 2^{\mathbb Z}$ is {\it unidirectional} if for all $s\in {\mathbb Z}$ there is a finite $T\subseteq {\mathbb N}$ such that
$$ \forall g\in {\mathbb Z}\ \exists t\in T\ x(g+t)\neq x(g+s+t). $$
\end{definition}\index{unidirectional}

Thus for unidirectional $2$-colorings on ${\mathbb Z}$ one can always search for distinct colors by shifting to the right. It is clear that, if a $2$-coloring on ${\mathbb Z}$ is unidirectional, then
its restriction on ${\mathbb N}$ is a $2$-coloring on ${\mathbb N}$. However, we have the following observation.

\begin{lem}\label{lem:unidirectional} Any $2$-coloring on ${\mathbb Z}$ is unidirectional.
\end{lem}

\begin{proof}
Suppose $x$ is a $2$-coloring on ${\mathbb Z}$. Fix $s\in {\mathbb Z}$ with $s\neq 0$. Let
$T\subseteq Z$ be the finite set witnessing that $x$ blocks $s$. Let $m$ be the least element of $T$.
Then we claim that $|m|+T\subseteq{\mathbb N}$ also witnesses that $x$ blocks $s$. To see this let
$g\in {\mathbb Z}$ be arbitrary. Consider the element $|m|+g$. By blocking there is $t\in T$ such
that $x(|m|+g+t)\neq x(|m|+g+s+t)$. Therefore $x(g+(|m|+t))\neq x(g+s+(|m|+t))$ with $|m|+t\in |m|+T$, as required.
\end{proof}

Thus indeed the restriction to ${\mathbb N}$ of any $2$-coloring on ${\mathbb Z}$ is a $2$-coloring on
${\mathbb N}$. Conversely, it is also easy to check that if $y\in 2^{\mathbb N}$ is a $2$-coloring on
${\mathbb N}$ then the element $x\in 2^{\mathbb Z}$ defined by
$$ x(n)=\left\{\begin{array}{ll}
y(n), & \mbox{ if $n\geq 0$,} \\
y(-n), & \mbox{ otherwise.}
\end{array}\right. $$
is a $2$-coloring on ${\mathbb Z}$.

\section{Subflows of $(2^{\N})^G$} \label{SEC 2NG}

Some of our results in this paper about Bernoulli subflows can be directly generalized to more general dynamical systems. Among continuous actions of $G$ the shift action on $(2^{\N})^G$ plays an important role. Let us recall the following basic fact from \cite{DJK} about this dynamical system. Again for the convenience of the reader we
give the proof below.

\begin{lem}\label{lem:universalBorelGspace} Let $G$ be a countable group with a Borel action on a standard Borel space $X$. Then there is a Borel embedding $\theta: X\to (2^{\N})^G$ such that for all $g\in G$ and $x\in X$, $\theta(g\cdot x)=g\cdot \theta(x)$.
\end{lem}

\begin{proof}
Let $U_0, U_1,\dots$ be a sequence of Borel sets in $X$ separating points. Define $\theta: X\to (2^{\N})^G$ by
$$ \theta(x)(g)(i)=1\iff g^{-1}\cdot x\in U_i. $$
Then $\theta$ is as required.
\end{proof}

Thus $(2^{\N})^G$ contains a $G$-invariant Borel subspace that is Borel isomorphic to the Borel $G$-space $X$.
In this sense $(2^{\N})^G$ is a {\it universal} Borel $G$-space among all standard Borel $G$-spaces. In the case that the space $X$ is a zero-dimensional Polish space and the action of $G$ on $X$ is continuous, we can improve the embedding $\theta$ to be continuous.

\begin{lem}\label{lem:universalPolishGspace} Let $G$ be a countable group with a continuous action on a zero-dimensional Polish space $X$. Then there is a continuous embedding $\theta: X\to (2^{\N})^G$ such that for all $g\in G$ and $x\in X$, $\theta(g\cdot x)=g\cdot \theta(x)$.
\end{lem}

\begin{proof}
In the proof of Lemma~\ref{lem:universalBorelGspace} if we take the $U_i$'s from a countable clopen base of $X$ the resulting $\theta$ is continuous.
\end{proof}

If $X$ is compact in addition, then the resulting $\theta$ is a homeomorphic embedding.

Because of these universality properties of $(2^\N)^G$ we are especially interested in establishing results about its subflows. Note that $(2^{\N})^G$ is isomorphic to the space $2^{\N\times G}$, and it is more convenient to
consider this latter space when we consider combinatorial properties of elements. The following lemmas are analogous to their counterparts, Lemmas~\ref{lem:basiccoloringlemma}, \ref{lem:orthogonallemma} and \ref{lem:minimallemma}, for Bernoulli flows. We state them without proof.

\begin{lem}\label{lem:analogbasic} Let $G$ be a countable group and $x\in 2^{\N\times G}$. Then $x$ is hyper aperiodic iff for any $s\in G$ there is $N\in \N$ and finite $T\subseteq G$ such that
$$ \forall g\in G\ \exists n<N\ \exists t\in T\ x(n,gt)\neq x(n, gst). $$
\end{lem}

\begin{lem}\label{lem:analogorthogonal} Let $G$ be a countable group and $x_0, x_1\in 2^{\N\times G}$. Then $x_0$ and $x_1$ are orthogonal iff there is $N\in\N$ and finite $T\subseteq G$ such that
$$ \forall g_0, g_1\in G\ \exists n<N\ \exists t\in T\ x_0(n,g_0t)\neq x_1(n,g_1t). $$
\end{lem}

\begin{lem}\label{lem:analogminimal} Let $G$ be a countable group and $x\in 2^{\N\times G}$. Then $x$ is minimal iff for all $N\in\N$ and finite $A\subseteq G$ there is a finite $T\subseteq G$ such that
$$ \forall g\in G\ \exists t\in T\ \forall n<N\ \forall a\in A\ x(n, gta)=x(n,a). $$
\end{lem}

Note that by Lemma~\ref{lem:analogbasic} if $x\in 2^{\N\times G}$ is such that $x(0,\cdot)$ is a $2$-coloring on $G$ then $x$ is hyper aperiodic. Hence the existence of hyper aperiodic points is an immediate corollary of the existence of $2$-colorings on $G$. In Chapter \ref{CHAP STUDY} we will show, among other facts, that every non-empty open subset of $(2^\N)^G$ contains a perfect set of pairwise orthogonal minimal hyper aperiodic points.
\chapter{\label{chap:basicconstructions}Basic Constructions of $2$-Colorings}

In this chapter we give some basic constructions of $2$-colorings on groups. The methods introduced here
are not as powerful as the one explored later in this paper. But they are simple and intuitive,
and using these methods we are able to construct $2$-colorings on all solvable groups, all
free groups and some of their extensions. In fact, other than constructing $2$-colorings on free groups (including ${\mathbb Z}$),
this chapter focuses primarily on constructing $2$-colorings on group extensions.

\section{$2$-Colorings on supergroups of finite index}

In this section we consider two constructions to obtain $2$-colorings on a countable group from $2$-colorings on a subgroup of finite index.

Let $G$ be a countable group and $H\leq G$ with $1<|G:H|=m<\infty$. Let $\alpha_1=1_G, \alpha_2, \dots, \alpha_m$ enumerate a set of representatives for all left cosets of $H$ in $G$. Given $x, y\in 2^H$, we define a function $\kappa_H(x;y)\in 2^G$ by
$$ \kappa_H(x;y)(g)=\left\{\begin{array}{ll} x(g), & \mbox{ if $g\in H$,} \\ y(h), & \mbox{ if $g\not\in H$ and $g=\alpha_ih$ for $1<i\leq m$.} \end{array}\right.
$$

\begin{figure}[h]
\begin{center}
\setlength{\unitlength}{1mm}
\begin{picture}(55,26)(0,0)

\put(0,0){\line(1,0){50}}
\put(50,0){\line(0,1){25}}
\put(0,0){\line(0,1){25}}
\put(0,25){\line(1,0){50}}
\put(53,22){\makebox(0,0)[b]{$G$}}
\put(5,21){\makebox(0,0)[b]{$H$}}
\put(5,10){\makebox(0,0)[b]{$x$}}
\put(10,0){\line(0,1){25}}
\put(15,21){\makebox(0,0)[b]{$\alpha_2H$}}
\put(15,9){\makebox(0,0)[b]{$y$}}
\put(20,0){\line(0,1){25}}
\put(40,0){\line(0,1){25}}
\put(45,21){\makebox(0,0)[b]{$\alpha_mH$}}
\put(44,9){\makebox(0,0)[b]{$y$}}
\put(30,10){\makebox(0,0)[b]{$\cdots\cdots$}}

\end{picture}
\caption{The function $\kappa_H(x;y)$.}
\end{center}
\label{fig:extensioncoloring1}
\end{figure}
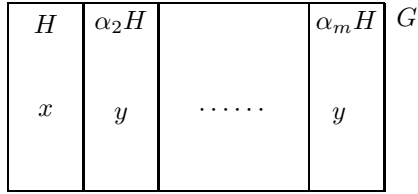

Thus $\kappa_H(x;y)$ is obtained by imposing $x$ on $H$ and $y$ on every other left coset of $H$ viewed as a copy of $H$ (see Figure~\ref{fig:extensioncoloring1}). Apparently the definition of $\kappa_H(x;y)$ depends on the particular choice of left coset representatives, and they are omitted in the notation just for simplicity. However, the results we prove below about $\kappa_H(x;y)$ will not depend on this choice. We first observe the following fact.

\begin{lem}\label{lem:Hxycoloring} Let $G$ be a countable group and $H\leq G$ with $|G:H|<\infty$. If $x$ is a $2$-coloring on $H$ and $y\in 2^H$ is such that $y\,\bot\, x$, then $\kappa_H(x;y)$ is a $2$-coloring on $G$.
\end{lem}

\begin{proof} Let $T_0=\{\alpha_1=1_G, \alpha_2, \dots, \alpha_m\}$. Since $x\,\bot\, y$, there is a finite set $T_1\subseteq H$ such that for all $h_0, h_1\in H$ there is $\tau\in T_1$ such that $x(h_0\tau)\neq y(h_1\tau)$. Given $s\in G$ with $s\neq 1_G$, let
$$I_s=\{1\leq i\leq m\: \alpha_i^{-1}s\alpha_i\in H\}.$$
Since $x$ is a $2$-coloring, for each $i\in I_s$, there is a finite set $T_{s,i}\subseteq H$ such that for all
$h\in H$ there is $\tau\in T_{s,i}$ such that $x(h\tau)\neq x(h\alpha_i^{-1}s\alpha_i\tau)$. Let
$$ T=T_0\left( T_1\cup\bigcup_{i\in I_s} T_{s,i}\right). $$
We verify that $T$ witnesses that $s$ is blocked by $\kappa_H(x;y)$. For this let $g\in G$. First there is some $1\leq i\leq m$ such that $g^{-1}\in \alpha_iH$. Then $g\alpha_i\in H$. If $gs\alpha_i\not\in H$, say $gs\alpha_i=\alpha_jh$ for $1<j\leq m$ and $h\in H$, then there is $\tau\in T_1$ such that
$$ \kappa_H(x;y)(g\alpha_i\tau)=x(g\alpha_i\tau)\neq y(h\tau)=\kappa_H(x;y)(gs\alpha_i\tau), $$
which finishes the proof by taking $t=\alpha_i\tau\in T_0T_1\subseteq T$.

If $gs\alpha_i\in H$, then $i\in I_s$ since $\alpha_i^{-1}s\alpha_i=(g\alpha_i)^{-1}(gs\alpha_i)\in H$. In this case there
is $\tau\in T_{s,i}$ such that
$$ \kappa_H(x;y)(g\alpha_i\tau)=x(g\alpha_i\tau)\neq x(g\alpha_i(\alpha_i^{-1}s\alpha_i)\tau)=x(gs\alpha_i\tau)=\kappa_H(x;y)(gs\alpha_i\tau). $$
Again, by letting $t=\alpha_i\tau\in T_0T_{s,i}\subseteq T$, we have that $\kappa_H(x;y)(gt)\neq \kappa_H(x;y)(gst)$, and our proof is complete.
\end{proof}

The idea of the above proof can be informally summarized as the following procedure. Given $s$ and $g$ we first transfer $g$ back to the ``standard" set $H$. If the corresponding element $gs$ is transferred to the same set, then we note that they are related by one of finitely many conjugates of $s$, and use the $2$-coloring property of $x$. If $gs$ stays out of $H$, then we use the orthogonality of $y$ and $x$ to finish the proof.

If we assume instead that $y$ is a $2$-coloring (and $x$ is not), then we can use the compliment of $H$ as our standard set, but this idea encounters a difficulty when $g$ and $gs$ are transferred to different left cosets (by the right multiplication of the same element) outside $H$. In this case we note that, if we assume that $H$ is a normal subgroup of $G$, then the difficulty disappears. Thus we have the following corollary of the above proof.

\begin{cor}\label{cor:Hxycoloring} Let $G$ be a countable group and $H\unlhd G$ with $1<|G:H|<\infty$. If $y$ is a $2$-coloring on $H$ and $x\in 2^H$ is such that $x\,\bot\, y$, then $\kappa_H(x;y)$ is a $2$-coloring on $G$.
\end{cor}

\begin{proof} Given $s\in G$ with $s\neq 1_G$, the witnessing set $T$ for $\kappa_H(x;y)$ blocking $s$ is the same as in the proof of Lemma~\ref{lem:Hxycoloring}. In fact, let $\alpha_i$ be such that $g\alpha_i\in H$. If $gs\alpha_i\not\in H$ then the proof is finished as before since $x\,\bot\, y$. If $gs\alpha_i\in H$, then $s=\alpha_i((g\alpha_i)^{-1}gs\alpha_i)\alpha_i^{-1}\in\alpha_i H\alpha_i^{-1}=H$. In this case let $j\neq i$, so that $g\alpha_j\not\in H$. Then $gs\alpha_j=g\alpha_j(\alpha_j^{-1}s\alpha_j)\not\in H$. Let $h\in H$ be such that $g\alpha_j=\alpha_kh$ for some $k$. Then $gs\alpha_j=\alpha_kh(\alpha_j^{-1}s\alpha_j)$, and there is $\tau\in T_{s,j}$ such that
$$ \kappa_H(x;y)(g\alpha_j\tau)=y(h\tau)\neq y(h(\alpha_j^{-1}s\alpha_j)\tau)=\kappa_H(x;y)(gs\alpha_j\tau) $$
by the assumption that $y$ is a $2$-coloring. Letting $t=\alpha_j\tau\in T_0T_{s,j}\subseteq T$, we have that
$\kappa_H(x;y)(gt)\neq \kappa_H(x;y)(gst)$ as required.
\end{proof}

In particular, this corollary applies when $H\leq G$ and $|G:H|=2$.

The same idea of the proof of Lemma~\ref{lem:Hxycoloring} can also be used to study when $\kappa_H(x_0;y_0)\,\bot\,\kappa_H(x_1;y_1)$. For instance, it can be shown that, if either $x_0$ or $y_0$ is orthogonal to both $x_1$ and $y_1$, then $\kappa_H(x_0;y_0)\,\bot\,\kappa_H(x_1;y_1)$ (note that this holds independently from the choice of left coset representatives in the definitions of $\kappa_H(x_0;y_0)$ and $\kappa_H(x_1;y_1)$). It follows that if $\{x_0,y_0, x_1, y_1\}$ is a set of pairwise orthogonal elements and $\{x_0,y_0\}\neq\{x_1,y_1\}$, then $\kappa_H(x_0;y_0)\,\bot\,\kappa_H(x_1;y_1)$. Below we state without proof a simple fact that can be justified with similar arguments.

\begin{lem}\label{lem:Hxyorthogonality} Let $G$ be a countable group and $H\leq G$ with $|G:H|<\infty$. If $X,Y\subseteq 2^H$ are disjoint such that $X\cup Y$ is a set of pairwise orthogonal elements of $2^H$, then the set
$$\{\kappa_H(x;y)\: x\in X,\ y\in Y\} $$
is a set of pairwise orthogonal elements of $2^G$.
\end{lem}

Throughout the rest of the paper we use $0$ to denote the constant $0$ function on a group and $1$ to denote the constant $1$ function.\index{$0$}\index{$1$} It follows immediately from Lemma~\ref{lem:orthogonallemma} that for any $2$-coloring $x$ on $H$, $x\,\bot\, 0,1$.

Recall that $H$ is said to have the $(\lambda,2)$-coloring property (where $\lambda\geq 1$ is a cardinal number)
if there exist $\lambda$ many pairwise orthogonal $2$-colorings on $H$. We thus have the following corollary.

\begin{cor}\label{cor:finiteindexCP} Let $G$ be a countable group, $H\leq G$ with $1<|G:H|<\infty$, and $\lambda\geq 1$ a cardinal number.
Suppose $H$ has the $(\lambda,2)$-coloring property. Then the following hold:
\begin{enumerate}
\item[(i)] If $\lambda$ is infinite then $G$ has the $(\lambda, 2)$-coloring property.
\item[(ii)] If $\lambda$ is finite then $G$ has the $(\frac{1}{2}\lambda(\lambda+3),2)$-coloring property.
\end{enumerate}
\end{cor}

\begin{proof}
Let $X$ be a set of pairwise orthogonal $2$-colorings on $H$ with $|X|=\lambda$. If $\lambda$ is infinite, then
note that $\{\kappa_H(x;0)\: x\in X\}$ is a set of pairwise orthogonal $2$-colorings on $G$ by Lemmas~\ref{lem:Hxycoloring} and \ref{lem:Hxyorthogonality}. Since $|\{\kappa_H(x;0)\: x\in X\}|=|X|=\lambda$, $G$ has the $(\lambda, 2)$-coloring property. If $\lambda$ is finite, we enumerate the elements of $X$ as $x_1,\dots, x_{\lambda}$. Consider the collection
$$ \left\{\kappa_H(x_i;x_j)\: 1\leq i<j\leq \lambda\right\}\cup\left\{\kappa_H(x_i;y)\: 1\leq i\leq \lambda,\ y\in\{0,1\}\right\}. $$
By Lemmas~\ref{lem:Hxycoloring}, \ref{lem:Hxyorthogonality} and the remark preceding Lemma~\ref{lem:Hxyorthogonality}, this is a set of pairwise orthogonal $2$-colorings on $G$. Its cardinality is $\frac{1}{2}\lambda(\lambda-1)+2\lambda=\frac{1}{2}\lambda(\lambda+3)$.
\end{proof}

For the rest of this section we consider a generalization of $\kappa_H(x;y)$ defined as follows. For $x_1,\dots, x_m\in 2^H$ define
$$\kappa_H(x_1,\dots, x_m)(g)=\kappa_H(x_1,\dots, x_m)(\alpha_ih)=x_i(h) $$
for $g=\alpha_ih$, where $1\leq i\leq m$, $h\in H$, and $\alpha_1=1_G, \dots, \alpha_m$ enumerate a set of representatives for all left cosets of $H$ in $G$.

\begin{figure}[h]
\begin{center}
\setlength{\unitlength}{1mm}
\begin{picture}(55,26)(0,0)

\put(0,0){\line(1,0){50}}
\put(50,0){\line(0,1){25}}
\put(0,0){\line(0,1){25}}
\put(0,25){\line(1,0){50}}
\put(53,22){\makebox(0,0)[b]{$G$}}
\put(5,21){\makebox(0,0)[b]{$H$}}
\put(5,10){\makebox(0,0)[b]{$x_1$}}
\put(10,0){\line(0,1){25}}
\put(15,21){\makebox(0,0)[b]{$\alpha_2H$}}
\put(15,10){\makebox(0,0)[b]{$x_2$}}
\put(20,0){\line(0,1){25}}
\put(40,0){\line(0,1){25}}
\put(45,21){\makebox(0,0)[b]{$\alpha_mH$}}
\put(44,10){\makebox(0,0)[b]{$x_m$}}
\put(30,10){\makebox(0,0)[b]{$\cdots\cdots$}}

\end{picture}
\caption{The function $\kappa_H(x_1,\dots, x_m)$.}
\end{center}
\label{fig:extensioncoloring2}
\end{figure}
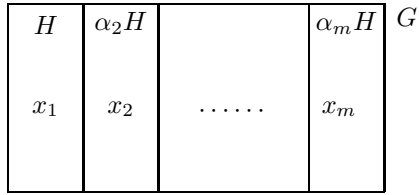

Figure~\ref{fig:extensioncoloring2} illustrates the definition of $\kappa_H(x_1,\dots, x_m)$. Clearly, $\kappa_H(x;y)=\kappa_H(x,y,\dots,y)$. In the following we prove a generalization of Lemma~\ref{lem:Hxycoloring} which guarantees that $\kappa_H(x_1,\dots, x_m)$ is a $2$-coloring by assuming one of the $x_i$ is a $2$-coloring orthogonal to all other $x_j$'s. In the proof we will use a well known lemma of Poincar\'{e}, which we recall below.

\begin{lem}\label{lem:Poincare} Let $G$ be a group and $H\leq G$ with $|G:H|<\infty$. Then there is $K\unlhd G$ such that $K\leq H$ and $|H:K|<\infty$.
\end{lem}

\begin{proof}
Let $\Sigma$ be the collection of all left cosets of $H$ in $G$. For each $g\in G$, let $\varphi(g)$ be a permutation of $\Sigma$ given by $\varphi(g)(\alpha H)=g\alpha H$. Then $\varphi: G\to S(\Sigma)$ is a group
homomorphism, where $S(\Sigma)$ is the group of all permutations of $\Sigma$. Here $S(\Sigma)$ is finite since $\Sigma$ is finite. Let $K=\ker(\varphi)$. Then $K\unlhd G$. It follows from the finiteness of $S(\Sigma)$ that $G/K$ is finite. Hence to finish the proof it suffices to verify that $K\leq H$. For this let $g\in K$, then $\varphi(g)=1_{S(\Sigma)}$, and in particular $\varphi(g)(H)=gH=H$, so $g\in H$.
\end{proof}

\begin{theorem}\label{thm:Hxxcoloring} Let $G$ be a countable group and $H\leq G$ with $|G:H|=m<\infty$. Let $x_1,\dots, x_m\in 2^H$. If there is $1\leq i\leq m$ such that $x_i$ is a $2$-coloring on $H$ and $x_i\,\bot\, x_j$ for any $1\leq j\leq m$ with $j\neq i$, then $\kappa_H(x_1,\dots, x_m)$ is a $2$-coloring on $G$.
\end{theorem}

\begin{proof}
Let $K\unlhd G$ be given by the preceding lemma. Then $K\leq H$ and $|G:K|<\infty$. Let $\gamma_1=1_G, \dots, \gamma_n$ enumerate a set of representatives for all cosets of $K$ in $G$. Let $1\leq i\leq m$ be such that
$x_i$ is a $2$-coloring on $H$ and that $x_i\,\bot\, x_j$ for all $j\neq i$, $1\leq j\leq m$.
 Since $K\leq H$ there is $1\leq p\leq n$ such that $K\gamma_p=\gamma_pK\subseteq \alpha_iH$.
 Let $T_0=\{\gamma_q^{-1}\gamma_p\: 1\leq q\leq n\}$. By the orthogonality assumptions there is a finite set $T_1\subseteq H$ such that for all $j\neq i$, $1\leq j\leq m$, and $h, h'\in H$ there is $\tau\in T_1$ such that $x_i(h\tau)\neq x_j(h'\tau)$.

Given $s\in G$ with $s\neq 1_G$, let
$$I_s=\{1\leq q\leq n\: \gamma_p^{-1}\gamma_q s\gamma_q^{-1}\gamma_p\in H\}.$$
Since $x_i$ is a $2$-coloring, for each $q\in I_s$, there is a finite set $T_{s,q}\subseteq H$ such that for all
$h\in H$ there is $\tau\in T_{s,q}$ such that
$$x_i(h\tau)\neq x_i(h\gamma_p^{-1}\gamma_q s\gamma_q^{-1}\gamma_p\tau).$$
Let
$$ T=T_0\left( T_1\cup\bigcup_{q\in I_s} T_{s,q}\right). $$
We claim that $T$ witnesses that $s$ is blocked by $\kappa_H(x_1, \dots, x_m)$. For this let $g\in G$. First there is some $1\leq q\leq n$ such that $g\in K\gamma_q$. Then $g\gamma_q^{-1}\in K$ and $g\gamma_q^{-1}\gamma_p\in K\gamma_p\subseteq \alpha_iH$. Let $h\in H$ be such that $g\gamma_q^{-1}\gamma_p=\alpha_ih$. If $gs\gamma_q^{-1}\gamma_p\not\in \alpha_iH$, say $gs\gamma_q^{-1}\gamma_p=\alpha_jh'$ for $j\neq i$, $1\leq j\leq m$, and $h'\in H$, then there is $\tau\in T_1$ such that
$$ \kappa_H(x_1,\dots,x_m)(g\gamma_q^{-1}\gamma_p\tau)=x_i(h\tau)\neq x_j(h'\tau)=\kappa_H(x_1,\dots, x_m)(gs\gamma_q^{-1}\gamma_p\tau), $$
which finishes the proof by taking $t=\gamma_q^{-1}\gamma_p\tau\in T_0T_1\subseteq T$.

If $gs\gamma_q^{-1}\gamma_p\in \alpha_iH$, then $q\in I_s$ since
$$\gamma_p^{-1}\gamma_qs\gamma_q^{-1}\gamma_p=(g\gamma_q^{-1}\gamma_p)^{-1}(gs\gamma_q^{-1}\gamma_p)\in (\alpha_i^{-1}H)(\alpha_iH)=H.$$
In this case there
is $\tau\in T_{s,q}$ such that
$$\begin{array}{l} \kappa_H(x_1,\dots, x_m)(g\gamma_q^{-1}\gamma_p\tau)=x_i(h\tau) \\
\neq x_i(h(\gamma_p^{-1}\gamma_qs\gamma_q^{-1}\gamma_p)\tau)=\kappa_H(x_1,\dots, x_m)(gs\gamma_q^{-1}\gamma_p\tau).
\end{array}
$$
Again, by letting $t=\gamma_q^{-1}\gamma_p\tau\in T_0T_{s,q}\subseteq T$ our proof is complete.
\end{proof}

Despite the tedious notation the idea of the above proof is quite simple: we use the underlying normal subgroup to transfer the elements to a standard set just as we did in the proof of Lemma~\ref{lem:Hxycoloring}, and then use
the assumptions of $2$-coloring and orthogonality to finish the proof. The same idea can be applied again to investigate when $\kappa_H(x_1,\dots, x_m)\,\bot\,\kappa_H(y_1,\dots,y_m)$. We state the following observation without proof.

\begin{lem}\label{lem:Hxxorthogonality} Let $G$ be a countable group and $H\leq G$ with $|G:H|=m<\infty$. Let
$x_1,\dots, x_m,y_1,\dots, y_m\in 2^H$. If there is $1\leq i\leq m$ such that $x_i\,\bot\,y_j$ for all $1\leq j\leq m$, then
$$\kappa_H(x_1,\dots, x_m)\,\bot\,\kappa_H(y_1,\dots, y_m).$$
\end{lem}

Using Theorem~\ref{thm:Hxxcoloring} and Lemma~\ref{lem:Hxxorthogonality} one can improve Corollary~\ref{cor:finiteindexCP} with the
general $\kappa_H(x_1,\dots, x_m)$ in place of $\kappa_H(x;y)$.

\section{$2$-Colorings on group extensions}  \label{SECT COLOR GRP EXT}

We begin by defining a natural map $2^G \times 2^H \rightarrow 2^{G \times H}$.

\begin{definition} Let $G$ and $H$ be countable groups, $x\in 2^G$, and $y\in 2^H$. Then the
{\it product} $xy$ is an element of $2^{G\times H}$ defined by
$$ (xy)(g,h)=x(g)y(h). $$
\end{definition}\index{product}

\begin{figure}[h]
\begin{center}
\setlength{\unitlength}{1mm}
\begin{picture}(70,28)(0,0)

\put(3,5){\line(1,0){30}}
\put(3,5){\line(0,1){21}}
\put(32,0){\makebox(0,0)[b]{$H$}}
\put(3,0){\makebox(0,0)[b]{$y$}}
\put(15,5){\line(0,1){21}}
\put(15,0){\makebox(0,0)[b]{$1$}}
\put(20,0){\makebox(0,0)[b]{$0$}}
\put(20,5){\line(0,1){21}}
\put(17,15){\makebox(0,0)[b]{$x$}}
\put(22,15){\makebox(0,0)[b]{$0$}}
\put(0,23){\makebox(0,0)[b]{$G$}}

\put(48,5){\line(1,0){30}}
\put(48,5){\line(0,1){21}}
\put(77,0){\makebox(0,0)[b]{$H$}}
\put(45,4){\makebox(0,0)[b]{$x$}}
\put(48,13){\line(1,0){30}}
\put(45,13){\makebox(0,0)[b]{$1$}}
\put(45,18){\makebox(0,0)[b]{$0$}}
\put(48,18){\line(1,0){30}}
\put(65,14){\makebox(0,0)[b]{$y$}}
\put(65,19){\makebox(0,0)[b]{$0$}}
\put(45,23){\makebox(0,0)[b]{$G$}}

\end{picture}
\caption{The product $xy$ viewed from two different perspectives.}
\end{center}
\label{fig:productcoloring}
\end{figure}
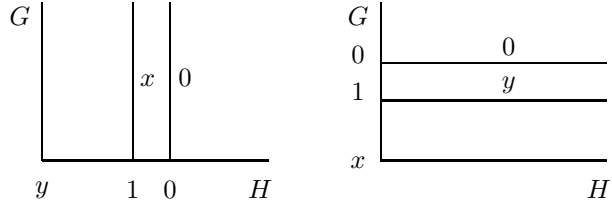

One way to view the product $xy$ is to regard $y$ as labeling the cosets of $G$ in $G\times H$ and impose the function $x$ on the copy of $G$ when the $y$ label is $1$ and $0$ when the $y$ label is $0$ (see Figure~\ref{fig:productcoloring}). Of course, by symmetry $x$ could be viewed as labeling cosets of $H$ as well.
The following proposition collects some elementary facts about the product. In the statement we use $0$ to
denote the constant zero element in $2^G$, $2^H$, or $2^{G\times H}$.

\begin{prop}\label{prop:productcoloring} Let $G$ and $H$ be countable groups, $x, x_1, x_2\in 2^G$, and $y, y_1, y_2\in 2^H$. Then the following hold:
\begin{enumerate}
\item[\rm (i)] $xy$ is a $2$-coloring iff both $x$ and $y$ are $2$-colorings.
\item[\rm (ii)] If $x_1\bot\, x_2$ and $0\not\in\overline{[y_1]}\cup\overline{[y_2]}$, then
$x_1y_1\bot\, x_2y_2$.
\item[\rm (iii)] If $0\not\in\overline{[y]}$, then $x_1\bot\, x_2$ iff $x_1y\bot\, x_2y$.
\item[\rm (iv)] If $xy\neq 0$, then $xy$ is minimal iff both $x$ and $y$ are minimal.
\end{enumerate}
\end{prop}

\begin{proof}
(i) First assume that $x$ and $y$ are both $2$-colorings. Note that $0\not\in\overline{[x]}$, and therefore there is a finite set $A\subseteq G$ such that
$$ \forall g\in G\ \exists a\in A\ x(ga)=1. $$
Similarly there is a finite set $B\subseteq H$ such that
$$ \forall h\in H\ \exists b\in B\ y(hb)=1. $$
To show that $xy$ is a $2$-coloring, fix a nonidentity $(s,u)\in G\times H$.
Without loss of generality assume $s\neq 1_G$. Thus we can find a finite set $T\subseteq G$
witnessing that $x$ blocks $s$. If $u\neq 1_H$ then we also have a finite set $S\subseteq H$
witnessing that $y$ blocks $u$. If $u=1_H$ we set $S=\emptyset$.
Then we claim that $(T\times B)\cup(A\times S)$ witnesses
that $xy$ blocks $(s,u)$ in $G\times H$. To see this, let $(g,h)\in G\times H$ be arbitrary.
We consider two cases. Case 1: $u=1_H$. Then we may find $b\in B$ such that $y(hb)=1$ and $t\in T$ such
that $x(gt)\neq x(gst)$. Note that also $y(hub)=1$. We have
$$ (xy)(gt,hb)=x(gt)y(hb)=x(gt)\neq x(gst)=x(gst)y(hub)=(xy)(gst,hub). $$
Since $(t,b)\in T\times B$ we are done. Case 2: $u\neq 1_H$. In this case we find $v\in S$ such that
$y(hv)\neq y(huv)$. If $y(hv)=1$ we find $a\in A$ such that $x(ga)=1$; if $y(huv)=1$ we find $a\in A$
such that $x(gsa)=1$. Either way we have
$$ (xy)(ga,hv)=x(ga)y(hv)=y(hv)\neq y(huv)=x(gsa)y(huv)=(xy)(gsa,huv). $$
Since $(a,v)\in A\times S$, our proof is completed.

For the converse assume without loss of generality that $x$ is not a $2$-coloring. Then there is
some $z\in\overline{[x]}$ with a nontrivial period $s\neq 1_G$. It follows that $zy\in \overline{[xy]}$
and that $(s,1_H)$ is a period of $zy$. Thus $xy$ is not a $2$-coloring.

The proof for (ii) is similar. For (iii) it suffices to note that, if $T\times S\subseteq G\times H$
is a finite set witnessing $x_1y\bot\, x_2y$, then $T$ witnesses $x_1\bot\, x_2$.

To prove (iv) we first assume that both $x$ and $y$ are minimal. Let $A \subseteq G \times H$ be finite. Without loss of generality we may assume $A=B\times C$ for $B\subseteq G$ and $C\subseteq H$. Let $T_B \subseteq G$ be finite with the property that for all $g \in G$ there is $t \in T_B$ with $x(g t b) = x(b)$ for all $b \in B$. Similarly, let $T_C\subseteq H$ be finite such that for all $h\in H$ there is $\tau\in T_C$ with $y(h\tau c)=y(c)$ for all $c\in C$. We claim that $T = T_B\times T_C$ works for $A$. For
this let $(g,h) \in G \times H$ be arbitrary. Let $t \in T_B$ be such that $x(gtb) = x(g)$ for all $b \in B$, and let $\tau \in T_C$ be such that $y(h\tau c) = y(c)$ for all $c \in C$. Then $(t,\tau) \in T$ and for all $(b,c) \in A$,
$$ (xy)(g tb, h\tau c) = x(gtb)y(h\tau c)= x(b)y(c)=(xy)(b,c).$$
This shows that $xy$ is minimal.

For the converse we assume $xy\neq 0$ is minimal. Note that we have both $x\neq 0$ and $y\neq 0$. We show that $x$ is minimal, and by symmetry it follows that $y$ is minimal too. For this fix $h_0\in H$ with $y(h_0)=1$. Let $A\subseteq G$ be finite. Without loss of generality we assume that there is $g_0\in A$ with $x(g_0)=1$. Since
$A\times\{h_0\}$ is a finite subset of $G\times H$, by the minimality of $xy$, there is a finite $T\subseteq G\times H$ such that for all $(g,h)\in G\times H$ there is $(t,\tau)\in T$ with
$(xy)(gta,h\tau h_0)=(xy)(a, h_0)$ for all $a\in A$. Let $T_G=\{t\in G\,:\, \exists \tau\in H\ (t,\tau)\in T\}$.
We claim that $T_G$ works for $A$. For this let $g\in G$ be arbitrary. There is $(t,\tau)\in T$ such that
$(xy)(gta, \tau h_0)=(xy)(a, h_0)$ for all $a\in A$. In particular, $t\in T_G$ and $(xy)(gtg_0, \tau h_0)=(xy)(g_0,h_0)=x(g_0)y(h_0)=1$. It follows that $y(\tau h_0)=1$ and therefore $x(gta)=(xy)(gta, \tau h_0)=(xy)(a, h_0)=x(a)$ for all $a\in A$. This shows that $x$ is minimal as required.
\end{proof}

\begin{cor}\label{cor:productcoloring} Let $G$ and $H$ be countable groups, and $\lambda, \kappa\geq 1$ cardinal numbers.
If $G$ has the $(\lambda,2)$-coloring property and $H$ has the $(\kappa,2)$-coloring property, then
$G\times H$ has the $(\lambda\cdot\kappa, 2)$-coloring property.
\end{cor}

One can also consider a slightly more general construction on the product group $G\times H$ as follows. For $x,z\in 2^G$ and $y\in 2^H$, define $xy_z\in 2^{G\times H}$ by
$$ (xy_z)(g,h)=x(g)y(h)+z(g)(1-y(h)). $$

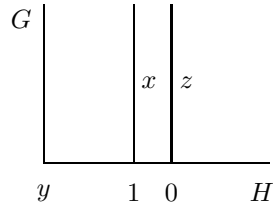
\begin{figure}[h]
\begin{center}
\setlength{\unitlength}{1mm}
\begin{picture}(30,28)(0,0)

\put(3,5){\line(1,0){30}}
\put(3,5){\line(0,1){21}}
\put(32,0){\makebox(0,0)[b]{$H$}}
\put(3,0){\makebox(0,0)[b]{$y$}}
\put(15,5){\line(0,1){21}}
\put(15,0){\makebox(0,0)[b]{$1$}}
\put(20,0){\makebox(0,0)[b]{$0$}}
\put(20,5){\line(0,1){21}}
\put(17,15){\makebox(0,0)[b]{$x$}}
\put(22,15){\makebox(0,0)[b]{$z$}}
\put(0,23){\makebox(0,0)[b]{$G$}}

\end{picture}
\caption{The function $xy_z$.}
\end{center}
\label{fig:semiproductcoloring}
\end{figure}

Here again $y$ is used to label the cosets of $G$ in $G\times H$. When the label is $1$ the coset is imposed the function $x$, and when the label is $0$ the coset is imposed the function $z$ (see Figure~\ref{fig:semiproductcoloring}).
Then similar to Proposition~\ref{prop:productcoloring} (i) one can show that, if both $x$ and $y$
 are $2$-colorings, and $z\bot\,x$, then $xy_z$ is a $2$-coloring on $G\times H$. Conversely, if $xy_z$
 is a $2$-coloring, we can only conclude that $y$ is a $2$-coloring due to asymmetry in this
 construction. In fact, both $x$ and $z$ can be periodic in this case. For example,
 let $x_0, z_0, y$ be $2$-colorings on ${\mathbb Z}$, and let $1$ denote the constant $1$ element
 in $2^{\mathbb Z}$. Let $x=x_01$ and $z=1z_0$. Then $x$ and $z$ are both periodic elements in
 $2^{{\mathbb Z}\times {\mathbb Z}}$, and $x\bot\,z$. It is easy to check that
$xy_z$ is a $2$-coloring on ${\mathbb Z}^3$.

We next consider general group extensions. Recall that in the preceding section we have considered
the case where $H$ is a normal subgroup of finite index in a countable group $G$. The constructions there fail to work when $H$ has infinite index in $G$, because the witnessing sets are no longer finite. In the next theorem we get around this problem by making use of $k$-colorings on the quotient.

\begin{theorem}\label{thm:transitivecoloring} Let $m,k\geq 2$ be integers and $\lambda\geq 1$ a cardinal number. Let $G$
be a countable group and $H\unlhd G$. Suppose $G/H$ has the $(\lambda,m)$-coloring property and
$H$ has the $(m,k)$-coloring property. Then
 $G$ has the $(\lambda,k)$-coloring property.
\end{theorem}

\begin{proof} We first define a $k$-coloring $x$ on $G$ assuming that
$z$ is an $m$-coloring on $G/H$ and $y_0, \dots, y_{m-1}$ are pairwise orthogonal $k$-colorings on $H$.
Let $R$ be a transversal for the cosets of $H$, i.e., $R$ contains exactly one element of each coset of $H$.
Let $\sigma: G\to R$ be such that for every $g\in G$, $\sigma(g)\in Hg=gH$.
Then define
$x: G\to k$ by
letting $$x(g)=y_{z(Hg)}(\sigma(g)^{-1}g).$$

We check that $x$ is a $k$-coloring on $G$.
For this fix $s\in G$ with $s\neq 1_G$. First assume $s\in H$. Since $y_0, \dots, y_{m-1}$ are all $k$-colorings
there are finite subsets $T_0, \dots, T_{m-1}\subseteq H$ such that for all
$h\in H$ and $i<m$ there are $t_i\in T_i$ such that $y_i(ht_i)\neq y_i(hst_i)$.
Let $T=\bigcup_{i<m}T_i$. We check that for any $g\in G$ there is $t\in T$ such that $x(gt)\neq x(gst)$.
Let $g\in G$. If $z(Hg)=i$ then for any $t\in T$, since $s,t\in H$, we have
$$x(gt)=y_i(\sigma(gt)^{-1}gt)=y_i(\sigma(g)^{-1}gt)$$ and
$$ x(gst)=y_i(\sigma(gst)^{-1}gst)=y_i(\sigma(g)^{-1}gst). $$
Thus if we let $t=t_i$ so that $y_i(ht_i)\neq y_i(hst_i)$ where $h=\sigma(g)^{-1}g$, then $x(gt)\neq x(gst)$.

Now we assume that $s\not\in H$. From the assumption that $z$ is an $m$-coloring we obtain a finite set
$F\subseteq R$ such that for any $g\in G$ there is $f\in F$ such that
$z(Hgf)\neq z(Hgsf)$. Let $\Gamma\subseteq H$ witness the orthogonality of $y_i$ and $y_j$ for all pairs $i,j<m$, $i\neq j$. That is,
for any $i,j<m$, $i\neq j$, and any $g_i, g_j\in H$, there is $\gamma\in\Gamma$ such that $y_i(g_i\gamma)\neq y_j(g_j\gamma)$.
Let $T=F\Gamma$. We again check that for any $g\in G$ there is $t\in T$ such that $x(gt)\neq x(gst)$.
First fix an $f\in F$ such that $z(Hgf)\neq z(Hgsf)$. For definiteness let $z(Hgf)=i$ and $z(Hgsf)=j$.
Then for any $\gamma\in\Gamma\subseteq H$,
$$ x(gf\gamma)=y_i(\sigma(gf)^{-1}gf\gamma) $$
and
$$ x(gsf\gamma)=y_j(\sigma(gsf)^{-1}gsf\gamma). $$
Thus letting $h_i=\sigma(gf)^{-1}gf$, $h_j=\sigma(gsf)^{-1}gsf$ and applying the orthogonality we obtain
a $\gamma\in\Gamma$ such that $y_i(h_i\gamma)\neq y_j(h_j\gamma)$. Letting $t=f\gamma$, we have thus verified that $x$ is a $k$-coloring.

Now we assume $z_0$ and $z_1$ are two orthogonal $m$-colorings on $G/H$.
Let $x_0$ and $x_1$ be defined similarly as above.
We verify that $x_0\bot\, x_1$, i.e., there is a finite set $\Phi\subseteq G$
such that for any $g_0, g_1\in G$ there is $\varphi\in\Phi$ such that $x_0(g_0\varphi)\neq x_1(g_1\varphi)$.
Let $F\subseteq R$ be finite such that for all $g_0, g_1\in G$ there is $f\in F$ such that $z_0(Hg_0f)\neq z_1(Hg_1f)$.
Let $\Gamma\subseteq H$ witness the orthogonality of all pairs $y_i$ and $y_j$ for $i,j<m$ and $i\neq j$.
Let $\Phi=F\Gamma$. Then for any $g_0, g_1\in G$, letting $f\in F$ be fixed as above,
$h_0=\sigma(g_0f)^{-1}g_0f$, $h_1=\sigma(g_1f)^{-1}g_1f$, $i=z_0(Hg_0f)$, $j=z_0(Hg_1f)$, and $\gamma$ such
that $y_i(h_0\gamma)\neq y_j(h_1\gamma)$, then
$$ x_0(g_0f\gamma)=y_{z_0(Hg_0f)}(\sigma(g_0f\gamma)^{-1}g_0f\gamma)=y_i(h_0\gamma), \mbox{ and} $$
$$ x_1(g_1f\gamma)=y_j(h_1\gamma). $$
Thus $x_0(g_0f\gamma)\neq x_1(g_1f\gamma)$.
\end{proof}

The following approach is an alternative way to obtain $2$-colorings on an extension from those on
a normal subgroup. Instead of assuming the existence of any $2$-coloring on the quotient we consider a
strong notion of a uniform $2$-coloring property on the normal subgroup.

\begin{definition}\label{def:uniformproperty}
 Let $G$ be a countable group. We say that $G$ has the {\it uniform $2$-coloring property} if there exists a perfect set $\{x_\sigma \,:\, \sigma \in 2^\N\}$ of pairwise orthogonal 2-colorings on $G$ such that
 \begin{enumerate}
 \item[(i)] for any $s\in G$ with $s\neq 1_G$, there is a finite set $T\subseteq G$ such that for any $x\in \{x_\sigma \,:\, \sigma \in 2^\N\}$, we have
     $$ \forall g\in G\ \exists t\in T\ x(g t) \neq x(g s t); $$
 \item[(ii)] for each $n\in \N$ there is a finite set $A_n \subseteq G$ such that
 for any $\sigma, \tau \in 2^\N$ with $\sigma(n) \neq \tau(n)$,
 $$ \forall g_0,g_1 \in G\ \exists a\in A_n\ x_\sigma(g_0 a) \neq x_\tau(g_1 a).$$
 \end{enumerate}
\end{definition}\index{$2$-coloring!uniform $2$-coloring property}\index{uniform $2$-coloring property}

\begin{theorem}\label{thm:uniformextension} Let $G$ be a countable group and $H \unlhd G$.
If $H$ has the uniform $2$-coloring property then so does $G$.
\end{theorem}

\begin{proof} Suppose $H$ has the uniform $2$-coloring property. As in the above definition, there
is a collection of $2$-colorings $\{y_\sigma \,:\,\sigma\in 2^\N \}$ on $H$, for each $s\in H$
with $s\neq 1_H$ there is $T_H(s)\subseteq H$, and for each $n\in\N$
 there is $A_n\subseteq H$  satisfying (i) and (ii).

We first deal with the case $|G:H|=\infty$.
Let $1_G=r_0, r_1, \ldots$ enumerate a transversal of the cosets of $H$ in $G$. Let
$\phi: {\mathbb N}\times{\mathbb N}\to \{0,1\}$ be a function with the following property:
\begin{quote} for any $i, n\in\N$, letting $j, k\in\N$ be the unique integers satisfying $r_jH=r_ir_{n+1}H$ and $r_kH=r_ir_{n+1}^{-1}H$, we have
either $\phi(i,n)\neq \phi(j,n)$ or $\phi(i,n)\neq\phi(k,n)$.
\end{quote}
To see that such a function exists, note that for any fixed $n\in{\mathbb N}$, the right multiplication
by $r_{n+1}$ induces a permutation $\pi_n$ on ${\mathbb N}$ such that $r_{\pi_n(i)}H=r_ir_{n+1}H$. Note that
$\pi_n$ has no fixed points. Thus in the statement of the property $j=\pi_n(i)$ and $k=\pi_n^{-1}(i)$.
The permutation $\pi_n$ can be decomposed into basic cycles of either finite or infinite length. In either case it is easy to assign values to indices so that no three consecutive indices in each cycle
are assigned the same value. Since $k,i,j$ are consecutive indices, we must have $\phi(i,n)\neq \phi(j,n)$ or $\phi(i,n)\neq \phi(k,n)$.

We then define infinitely many elements $\tau_i\in 2^\N$ for $i\in \N$ by letting $\tau_i(n)=\phi(i,n)$.  We will also use a coding function $\langle\cdot,\cdot\rangle: 2^\N \times 2^\N \rightarrow 2^\N$ defined by
$\langle\tau,\sigma\rangle(2n)=\tau(n)$ and $\langle\tau,\sigma\rangle(2n+1)=\sigma(n)$ for all $n\in\N$.

We are now ready to construct a collection $\{x_\sigma\,:\, \ \sigma \in 2^\N\}$ of pairwise orthogonal 2-colorings on $G$.
For each $\sigma \in 2^\N$ define $x_\sigma$ by
$$ x_\sigma(r_ih)=y_{\langle\tau_i,\sigma\rangle}(h). $$

We verify that each $x_\sigma$ is a 2-coloring on $G$. Let $s\in G$ with $s\neq 1_G$ and let $g\in G$ be arbitrary. If $s\in H$ then there exists $t\in T_H(s)$ such that $x_\sigma(g t)\neq x_\sigma(g s t)$. If $s \not\in H$ let $s\in r_{n+1}H$, $g\in r_iH$, $gs\in r_jH$ and $gs^{-1}\in r_kH$. Then by the property of $\phi$ either $\phi(i,n)\neq \phi(j,n)$ or $\phi(i,n)\neq \phi(k,n)$. Therefore either $\langle\tau_i,\sigma\rangle(2n) \neq \langle\tau_j,\sigma\rangle(2n)$ or
$\langle\tau_i,\sigma\rangle(2n) \neq \langle\tau_k,\sigma\rangle(2n)$. It follows that if we let $T= A_{2n} \cup s^{-1} A_{2n}$ then there exists $t\in T$ such that $x_\sigma(g t)\neq x_\sigma(g s t)$. Note that the choice of $T$ does not depend on $\sigma$ so our collection of 2-colorings on $G$ satisfies property (i) in Definition~\ref{def:uniformproperty}.

For property (ii) in Definition~\ref{def:uniformproperty} it is clear that the set
$B_n=A_{2n+1}$ works for $n\in\N$. This finishes
the proof in the case $H$ has infinite index in $G$.

As for the case $|G:H|=m<\infty$, we can use an easy adaptation of the above construction. In this case
the function $\phi$ would be only defined on a finite domain $(m-1)\times m$. We then extend its definition to
${\mathbb N}\times{\mathbb N}$ using value $0$ and proceed as above. The resulting functions are as required.
\end{proof}

\section{\label{sec:3.2}$2$-Colorings on ${\mathbb Z}$}

For the rest of this chapter we construct concrete $2$-colorings on concrete groups. In this section we
focus on the group ${\mathbb Z}$. We show that ${\mathbb Z}$ has the uniform $2$-coloring property.

We will use the following notation. Let
$$ 2^{\prec{\mathbb Z}}=\bigcup_{l\leq r\in{\mathbb Z}} 2^{[l,r]}. $$
For $p\in 2^{\prec{\mathbb Z}}$ we let $|p|=l-r+1$ if $p\in 2^{[l,r]}$.
For $p\in 2^{\prec{\mathbb Z}}$ we let $\bar{p}(i)=1-p(i)$ for all $i\in \dom(p)$ and call it the {\it conjugate} of $p$. \index{conjugate} Thus $\dom(\bar{p})=\dom(p)$.
For $p, q\in 2^{\prec{\mathbb Z}}$, we write $p\subseteq q$ if $\dom(p)\subseteq\dom(q)$ and $q\upharpoonright \dom(p)=p$.
The group ${\mathbb Z}$ acts on $2^{\prec{\mathbb Z}}$ naturally: for $s\in {\mathbb Z}$ and $p\in 2^{\prec{\mathbb Z}}$, let
$$ (s+p)(i)=p(i-s). $$
Thus $\dom(s+p)=s+\dom(p)$. We write $p\sim q$ if there is $s\in {\mathbb Z}$ such that $s+p=q$.
For $p_0, p_1\in 2^{\prec{\mathbb Z}}$, say $p_0\in 2^{[l_0,r_0]}$, $p_1\in 2^{[l_1,r_1]}$, we write $p_0^\smallfrown p_1$ or simply $p_0p_1$
for the unique $q\in 2^{[l_0, r_0+1+r_1-l_1]}$ such that $q\!\upharpoonright\! [l_0,r_0]=p_0$ and $q\upharpoonright [r_0+1, r_0+1+r_1-l_1]\sim p_1$.
By iteration we can define the notation $p_0^\smallfrown p_1^\smallfrown \cdots ^\smallfrown p_n$ or $p_0p_1\cdots p_n$.

We also let
$$ P=\bigcup_{k\in{\mathbb N}} 2^{[-k,k]}. $$
For $p, q\in P$, say $p\in 2^{[-k,k]}$ and $q\in 2^{[-l,l]}$, we write $p\sqsubseteq q$ if $2k+1\mid l-k$ and for all $i\in{\mathbb Z}$,
if $D=[(2k+1)i+k+1, (2k+1)(i+1)+k]\subseteq\dom(q)$ then ${q\!\upharpoonright\! D}\sim p$ or ${q\!\upharpoonright\! D}\sim\bar{p}$.

\begin{figure}[h]
\begin{center}
\setlength{\unitlength}{1mm}
\begin{picture}(100,28)(0,0)

\put(0,20){
\put(0,4){\line(1,0){100}}
\put(2,0){\makebox(0,0)[b]{${\mathbb Z}$}}
\put(15,-0.5){\makebox(0,0)[b]{$-l$}}
\put(15,4){\line(0,1){1}}
\put(45,-0.5){\makebox(0,0)[b]{$-k$}}
\put(45,4){\line(0,1){1}}
\put(55,0){\makebox(0,0)[b]{$k$}}
\put(55,4){\line(0,1){1}}
\put(85,0){\makebox(0,0)[b]{$l$}}
\put(85,4){\line(0,1){1}}
\put(50,5){\makebox(0,0)[b]{$\overbrace{\ \ \ \ \ \ \ }^{2k+1}$}}
}

\put(0,3){
\put(45,10){\line(1,0){10}}
\put(45,10){\line(0,1){1}}
\put(55,10){\line(0,1){1}}
\put(50,11){\makebox(0,0)[b]{$p$}}
}

\put(0,-18){
\put(15,20){\line(1,0){70}}
\put(15,20){\line(0,1){1}}
\put(25,20){\line(0,1){1}}
\put(35,20){\line(0,1){1}}
\put(45,20){\line(0,1){1}}
\put(55,20){\line(0,1){1}}
\put(65,20){\line(0,1){1}}
\put(75,20){\line(0,1){1}}
\put(85,20){\line(0,1){1}}
\put(50,16){\makebox(0,0)[b]{$q$}}
\put(20,21){\makebox(0,0)[b]{$\bar{p}$}}
\put(30,21){\makebox(0,0)[b]{$p$}}
\put(40,21){\makebox(0,0)[b]{$p$}}
\put(50,21){\makebox(0,0)[b]{$\bar{p}$}}
\put(60,21){\makebox(0,0)[b]{$\bar{p}$}}
\put(70,21){\makebox(0,0)[b]{$p$}}
\put(80,21){\makebox(0,0)[b]{$\bar{p}$}}
}

\end{picture}
\caption{An illustration of $p\sqsubseteq q$.}
\end{center}
\label{fig:sqsubseteq}
\end{figure}
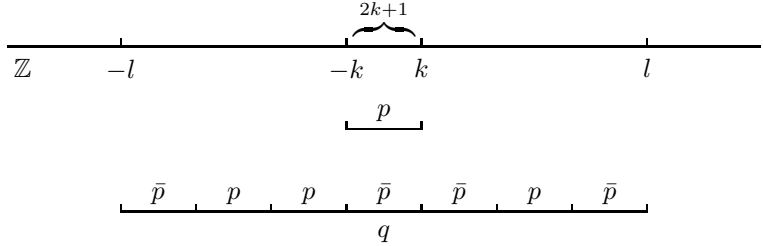

Note that $\sqsubseteq$ is a transitive relation, i.e., if $p_0\sqsubseteq p_1$ and $p_1\sqsubseteq p_2$ then $p_0\sqsubseteq p_2$.
Also for $p\in 2^{[-k,k]}$ and $x\in 2^{\mathbb Z}$, we write $p\sqsubseteq x$ if for all $i\in{\mathbb N}$, $p\sqsubseteq x\upharpoonright[-i(2k+1)-k,i(2k+1)+k]$.
We now define two operations on $P$, $\Phi_0$ and $\Phi_1$. For $p\in P$, let $\Phi_0(p)$ and $\Phi_1(p)$ be the unique elements of $P$ so that
$$\Phi_0(p)\sim pp\bar{p}pp\bar{p}p \ \ \mbox{ and }\ \ \Phi_1(p)\sim \bar{p}p\bar{p}p\bar{p}p\bar{p}. $$
Note that $p\sqsubseteq \Phi_0(p),\Phi_1(p)$ and $|\Phi_0(p)|=|\Phi_1(p)|=7|p|$. Also for $i=0,1$, $\overline{\Phi_i(p)}=\Phi_i(\bar{p})$.

\begin{lem} \label{lem:bot} Let $p, q\in P$ and $x,y\in 2^{\mathbb Z}$. If $|p|=|q|$ and $\Phi_0(p)\sqsubseteq x$, $\Phi_1(q)\sqsubseteq y$, then $x\bot\, y$.
\end{lem}

\begin{proof}
Let $T=\{ i|p|\,:\, 0\leq i\leq 7\}$. Let $$p_0\in\{ \Phi_0(p)\Phi_0(p), \Phi_0(\bar{p})\Phi_0(\bar{p}), \Phi_0(\bar{p})\Phi_0(p), \Phi_0(p)\Phi_0(\bar{p})\}$$
and $$p_1\in\{ \Phi_1(q)\Phi_1(q), \Phi_1(\bar{q})\Phi_1(\bar{q}), \Phi_1(\bar{q})\Phi_1(q), \Phi_1(q)\Phi_1(\bar{q})\}.$$ By direct inspection it can be shown
that for any $0\leq g_0, g_1< 7|p|$ there is $t\in T$ such that $p_0(p)(g_0+t)\neq p_1(q)(g_1+t)$.
In fact, the $0,1$-sequence $\langle p_0(g_0+t)\: t\in T\rangle$ consists of at least two occurrences of
$00$ and $11$ which are separated by at most four digits in between, but this property fails for the sequence
$\langle p_1(g_1+t)\: t\in T\rangle$.
Now for any $g_0, g_1\in {\mathbb Z}$ there are $p_0, p_1$ as above and $p_0'\subseteq x$, $p_1'\subseteq y$ such that $g_0\in \dom(p_0')$, $g_1\in \dom(p_1')$,
$p_0\sim p_0'$ and $p_1\sim p_1'$. Thus there is $t\in T$ such that $p_0'(g_0+t)\neq p_1'(g_1+t)$. In particular $x(g_0+t)\neq y(g_1+t)$.
\end{proof}

In the sequel we use the notation $2^{<\N}$ to denote the set of all finite binary sequences. I.e., $2^{<\N}=\bigcup_{n\in\N} 2^n$. For $u\in 2^{<\N}$, let $|u|$ denote the length of $u$.

\begin{theorem} \label{thm:coloringZ} ${\mathbb Z}$ has the uniform $2$-coloring property.
\end{theorem}

\begin{proof}
We define a system $(p_u)_{u\in 2^{<\N}}$ of elements of $P$ by induction on $|u|$ so that the following conditions are satisfied:
\begin{enumerate}
\item[(i)] for all $u, v\in 2^{<\N}$ with $|u|=|v|$, $|p_u|=|p_v|$;
\item[(ii)] for all $u\in 2^{<\N}$, $\Phi_0(p_u)\sqsubseteq p_{u^\smallfrown 0}$ and $\Phi_1(p_u)\sqsubseteq p_{u^\smallfrown 1}$;
\item[(iii)] for every $u\in 2^{<\N}$, there is $i\in \dom(p_u)$ such that $i+|u|\in \dom(p_u)$ and $p_u(i)\neq p_u(i+|u|)$.
\end{enumerate}
To begin the definition, let $\dom(p_\emptyset)=\{0\}$ and $p_\emptyset(0)=0$. In general suppose all $p_u$ where $|u|\leq n$ have been defined.
We first define $q_{u^\smallfrown 0}$, $q_{u^\smallfrown 1}$ to satisfy the conditions (ii) and (iii). For this let
$i\in \dom(\Phi_0(p_u))$ so that $i+n+1\not\in \dom(\Phi_0(p_u))$. Let $q_{u^\smallfrown 0}\sqsupseteq \Phi_0(p_u)$ be such that
$q_{u^\smallfrown 0}(i+n+1)\neq \Phi_0(p_u)(i)$. Then define $q_{u^\smallfrown 1}\sqsupseteq \Phi_1(p_u)$ similarly. After all $q_v$ where $|v|=n+1$ have
been defined this way let $l=\max\{|q_v|\,:\, |v|=n+1\}$ and define $p_{u\smallfrown i}$ so that $|p_{u\smallfrown i}|=l$, $q_{u\smallfrown i}\subseteq p_{u\smallfrown i}$ and $\Phi_i(p_u)\sqsubseteq p_{u\smallfrown i}$.
This finishes the definition of $(p_u)_{u\in 2^{<\N}}$.

Now for $\alpha\in 2^\N$ we let $x_\alpha=\bigcup_{n\in\N} p_{\alpha\upharpoonright n}$. We claim that each $x_\alpha$ is a $2$-coloring
on ${\mathbb Z}$. To verify this let $s\in{\mathbb Z}$. Let $n=|s|$ and $u=\alpha\upharpoonright n$. Let $i\in \dom(p_u)$ such
that $i+n\in\dom(p_u)$ and $p_u(i)\neq p_u(i+n)$. Let $T=[0,2|p_u|]$. Now let $g\in {\mathbb Z}$ be arbitrary. Then noting that $p_u\sqsubseteq x_\alpha$,
there is $q\subseteq x_\alpha$ with $g\in\dom(q)$ and $q\sim p_u$ or $q\sim \bar{p}_u$. Letting $j$ to be the least integer greater than $g$ with $j\not\in \dom(q)$, we have that
$x_\alpha(i+j+\frac{1}{2}(|p_u|-1))\neq x_\alpha(i+n+j+\frac{1}{2}(|p_u|-1))$. Thus if we let $t=i+j+\frac{1}{2}(|p_u|-1)-g\in T$ (if $s>0$)
or $t=i+n+j+\frac{1}{2}(|p_u|-1)\in T$ (if $s<0$), we must have that $x_{\alpha}(g+t)\neq x_\alpha(g+s+t)$.
Note that the set $T$ only depends on $s$ and not on $\alpha$, since $|p_u|$ only depends on $|u|$ by property (i). This shows that the set $\{x_\alpha\,:\, \alpha\in 2^\N\}$ satisfies
Definition~\ref{def:uniformproperty} (i).

Finally, suppose $\alpha,\beta\in 2^\N$ with $\alpha(n)\neq\beta(n)$. Let $u=\alpha\upharpoonright n$ and $v=\beta\upharpoonright n$. Without loss of generality assume $u^\smallfrown 0\subseteq \alpha$ and
$v^\smallfrown 1\subseteq \beta$. Then $\Phi_0(p_u)\sqsubseteq p_{u^\smallfrown 0}\sqsubseteq x_\alpha$ and $\Phi_1(p_v)\sqsubseteq p_{v^\smallfrown 1}\sqsubseteq x_\beta$. By Lemma~\ref{lem:bot}, $x_\alpha\bot\, x_\beta$. Moreover, the proof of Lemma~\ref{lem:bot} shows that the witnessing set
can be taken as $\{in\,:\, 0\leq i\leq 7\}$, which depends only on $n$ and not on $\alpha$ and $\beta$. This shows that the set $\{x_\alpha\,:\, \alpha\in 2^\N\}$ satisfies
Definition~\ref{def:uniformproperty} (ii).
\end{proof}

We remark that, using Lemma~\ref{lem:minimallemma}, it is easy to check that all $2$-colorings constructed in the above proof are minimal. By an obvious modification of the above proof, we have the following corollary.

\begin{cor} Let $U$ be any given open subset of $2^{\mathbb Z}$. Then there is a perfect set of
pairwise orthogonal minimal $2$-colorings in $U$.
\end{cor}

The following corollary follows immediately from Theorem~\ref{thm:uniformextension}.

\begin{cor}\label{cor:normaloverZ} Let $G$ be a countable group. If ${\mathbb Z}\unlhd G$, then $G$ has the uniform $2$-coloring property. In particular, for any countable group $G$, $G\times {\mathbb Z}$ has the uniform $2$-coloring property.
\end{cor}

Before closing this section we briefly turn to a curious question about constructing $2$-colorings on ${\mathbb Z}$ that are orthogonal to their conjugates.
Note that our construction above  does not produce $2$-colorings on ${\mathbb Z}$
orthogonal to their own conjugates. One needs a slightly different
construction to achieve this.

For any $x\in 2^{\mathbb Z}$ let $\overline{x}\in 2^{\mathbb Z}$ be defined by $\overline{x}(n)=1-x(n)$ for all $n\in{\mathbb Z}$.
If $x$ is a $2$-coloring then so is $\bar{x}$. For any $x\in 2^{\mathbb Z}$ we also define $x'\in 2^{\mathbb Z}$ as follows:
$$ x'(n)=\left\{\begin{array}{ll}
x(n/3), & \mbox{ if $3\mid n$,} \\
0, & \mbox{ otherwise.}
\end{array}\right.
$$
The following facts are easy to see. For any $x\in 2^{\mathbb Z}$, $x'\bot\,\overline{x'}$,
since $x'$ does not contain two consecutive 1s. Also, if  $x$ is a $2$-coloring on ${\mathbb Z}$,
then so is $x'$. This is because $x'$ blocks $3n$ for any $n\neq 0$, and therefore it blocks $n$ for all $n\neq 0$ by Corollary~\ref{cor:blockinglemma}. Using Lemma~\ref{lem:minimallemma} it is clear that $x$ is minimal iff $\overline{x}$ is minimal iff $x'$ is minimal. Finally if $x,y\in 2^{\mathbb Z}$, then $x\bot\, y$ iff $\overline{x}\bot\, \overline{y}$ iff
 $x'\bot\, y'$.

Thus we have the following corollary.

\begin{cor} There is a perfect set $X$ of pairwise orthogonal minimal $2$-colorings on ${\mathbb Z}$ such that for any $x\in X$, $x\bot\, \overline{x}$.
\end{cor}

One can also modify the construction in an obvious way to obtain such families of $2$-colorings inside
any given open set.

\section{\label{sec:3.3}$2$-Colorings on nonabelian free groups}

In this section we show that all nonabelian free groups have the uniform $2$-coloring property.
We will need the following observation.

\begin{definition} Two elements $x_0,x_1\in 2^{\mathbb Z}$ are {\it positively orthogonal}, denoted $x_0{\bot\!}^+ x_1$, if there is a finite $T\subseteq\N$ such that
$$ \forall g_0,g_1\in {\mathbb Z}\ \exists t\in T\ x_0(g_0+t)\neq x_1(g_1+t). $$
\end{definition}\index{positively orthogonal}\index{orthogonal!positively orthogonal}

\begin{lem}\label{lem:positiveorthogonal} For $x_0, x_1\in 2^{\mathbb Z}$, if $x_0\bot\, x_1$ then $x_0{\bot\!}^+x_1$.
\end{lem}

\begin{proof}
This is similar to the proof of Lemma~\ref{lem:unidirectional}. Let $T\subseteq {\mathbb Z}$ witness
$x_0\bot\,x_1$. Let $m$ be the least element of $T$. Then $|m|+T\subseteq \N$ witnesses $x_0\bot\,x_1$
as well.
\end{proof}

We are now ready to consider free groups. Let ${\mathbb F}_n$ be the free group with $n$ generators, where $n\geq 2$ is an integer. For notational uniformity we use $F_\omega$ to denote the free group with countably infinitely many generators, where $\omega$ denotes the first infinite ordinal.  We will combine the two cases by considering ${\mathbb F}_n$ with $n$ generators, where $2\leq n\leq \omega$.\index{$\omega$}\index{${\mathbb F}_n$}\index{${\mathbb F}_\omega$}

Fix $2\leq n\leq \omega$. For any $x\in 2^{\mathbb Z}$, we define
$x^*\in 2^{{\mathbb F}_n}$ by
$x^*(w)=x(|w|)$, where $|w|$ is the length of the reduced word $w$.

\begin{theorem}\label{thm:free} Let $2\leq n\leq \omega$. If $x$ is a $2$-coloring on ${\mathbb Z}$, then $x^*$ is a $2$-coloring on ${\mathbb F}_n$. In addition, for $x_0,x_1\in 2^{\mathbb Z}$, if $x_0\bot\,x_1$, then $x_0^*\bot\,x_1^*$.
\end{theorem}

\begin{proof}  Let $A=\{ a_m\,:\, m< n\}$ be a generating set of elements of
${\mathbb F}_n$. Let $s\in {\mathbb F}_n$ with $s\neq 1_{{\mathbb F}_n}$.
For each integer $i\in [-2|s|,2|s|]$ with $i\neq 0$, let $L_i\subseteq {\mathbb N}$ be a finite set
such that for any $j\in{\mathbb Z}$ there is $l\in L_i$ with $x(j+l)\neq x(j+i+l)$.  Let
$$T=\{ t\in {\mathbb F}_n\,:\,
\exists i, i'\in [-|s|,|s|]\ \ |t|\in L_i+i'\}. $$
We check that for any $g\in {\mathbb F}_n$ there is $t\in T$ such that
$x^*(gt)\neq x^*(gst)$. For this let $g\in {\mathbb F}_n$. We consider two cases. Case 1: $|g|\neq|gs|$.
In this case let $i=|gs|-|g|$. Then $0<|i|\leq |s|$. Let $l\in L_i\subseteq {\mathbb N}$ be such that $x(|g|+l)\neq x(|g|+i+l)$.
There is $t$ with $|t|=l$ such that $|gt|=|g|+l$ and $|gst|=|gs|+l$.
Now $t\in T$ (with $i'=0$)
and $x^*(gt)=x(|g|+l)\neq x(|g|+i+l)=x(|gs|+l)=x^*(gst)$.
Case 2: $|g|=|gs|$. Then from the structure of the free group we get $u, v\in {\mathbb F}_n$ such that
$s=u^{-1}v$, $|s|=2|u|=2|v|$ and $|gu^{-1}|=|g|-|u|$. Note that $|gu^{-1}|\neq |gu^{-1}vu^{-1}|$ and their difference
$i\leq |s|$. Thus by a similar construction as that in Case 1 there is $t_0$ with $|t_0|\in L_i$ such that
$x^*(gu^{-1}t_0)\neq x^*(gu^{-1}vu^{-1}t_0)$. Now let $t=u^{-1}t_0$, then $x^*(gt)\neq x^*(gst)$ and $t\in T$
with $i'=|t|-|t_0|\leq |u^{-1}|\leq |s|$. This finishes the proof that $x^*$ is a $2$-coloring on ${\mathbb F}_n$.

Now suppose $x_0\bot\,x_1$. Then by Lemmas~\ref{lem:unidirectional} and \ref{lem:positiveorthogonal}
$x_0$ and $x_1$ are positively orthogonal unidirectional $2$-colorings on ${\mathbb Z}$. Let $L\subseteq {\mathbb N}$ be
such that for any $j_0, j_1\in {\mathbb Z}$ there is $l\in L$ with $x_0(j_0+l)\neq x_1(j_1+l)$. Let $a_0, a_1\in A$ be arbitrary and
$T=\{a_i^l, a_i^{-l}\,:\, i=0,1,\ l\in L\}$.
Let $g_0, g_1\in F_n$ be arbitrary. Let $j_0=|g_0|$ and $j_1=|g_1|$. Let $l\in L$ be such that $x_0(j_0+l)\neq x_1(j_1+l)$.
Then there is $t\in T$ such that $|t|=l$, $|g_0t|=|g_0|+|t|$ and $|g_1t|=|g_1|+|t|$. Then $x_0^*(g_0t)=x_0(|g_0|+|t|)=x_0(j_0+l)\neq x_1(j_1+l)=x_1(|g_1|+|t|)=x_1^*(g_1t)$.
This shows that $x_0^*\bot\, x_1^*$.
\end{proof}

\begin{theorem}\label{thm:uniformfree} For any $1\leq n\leq\omega$ the free group ${\mathbb F}_n$ has the uniform $2$-coloring property.
\end{theorem}

\begin{proof}
Let $\{x_\alpha\,:\,\alpha\in 2^\N\}$ be a collection of $2$-colorings on ${\mathbb Z}$ witnessing
the uniform $2$-coloring property for ${\mathbb Z}$ from Theorem~\ref{thm:coloringZ}. Then for $2\leq n\leq \omega$, the collection $\{x_\alpha^*\,:\, \alpha\in 2^\N\}$ witnesses the uniform $2$-coloring
property for ${\mathbb F}_n$. This is because, by the above proof, the set $T$ witnessing that
 $x_\alpha^*$ blocks $s$ depends on $s$ only and does not depend on $\alpha$; in addition, if
 the set $L$ witnessing
 the orthogonality of $x_\alpha$ and $x_\beta$ depends only on the index $n$ where $\alpha(n)\neq\beta(n)$, then the set $T$ witnessing the orthogonality of $x_\alpha^*$ and $x_\beta^*$
 depends only on $n$.
\end{proof}

Since the free groups have the ACP (Corollary~\ref{cor:ACPexamples} (ii)), we have the following immediate corollary.

\begin{cor}\label{cor:perfectfree} Let $2\leq n\leq\omega$. Let $U$ be any given open subset of $2^{{\mathbb F}_n}$. Then there is a perfect set of pairwise orthogonal $2$-colorings in $U$.
\end{cor}

We also have the following immediate corollary from Theorem~\ref{thm:uniformextension}.

\begin{cor} Let $G$ be a countable group. If for some $1\leq n\leq\omega$, ${\mathbb F}_n\unlhd G$, then $G$ has the uniform $2$-coloring property.
\end{cor}

Note that the definition of $x^*$ makes sense even for $n=1$. And in this case the proofs of the theorems still
work and give another collection witnessing the uniform $2$-coloring property for ${\mathbb Z}$.

Moreover, when only restrictions of $2$-colorings on the semigroup ${\mathbb N}$ are considered, we obtain a collection of $2$-colorings on $\N$ witnessing the uniform $2$-coloring property for $\N$. Thus in particular, $\N$ has the uniform $2$-coloring property.

Finally we remark that if $x$ is a $2$-coloring on ${\mathbb Z}$, the $2$-coloring $x^*$ is actually a two-sided $2$-coloring. This is because, the dual definition of $x^*$ would be the same for right actions of ${\mathbb F}_n$ and the dualized proof of Theorem~\ref{thm:free} would show that $x^*$ is
a right $2$-coloring.

Before closing this section we give a construction of a $2$-coloring on ${\mathbb F}_n$, $n>1$, that is not a two-sided $2$-coloring.

\begin{theorem} For $n>1$ there exists a $2$-coloring on ${\mathbb F}_n$ that is not a right $2$-coloring.
\end{theorem}

\begin{proof}
It suffices to construct a $2$-coloring on ${\mathbb F}_n$ that is right-periodic. Fix $n>1$. Let $a$ be one of the generators of ${\mathbb F}_n$, and let $F$ be the free subgroup generated by the other $n-1$ generators. Let
$x_0\,\bot\,x_1\in 2^F$ be $2$-colorings on $F$ (they can be obtained by Theorem~\ref{thm:uniformfree}). Let $y\in 2^{\mathbb Z}$ be any $2$-coloring on ${\mathbb Z}$, and $y^*$ be the word-length $2$-coloring of ${\mathbb F}_n$ (Theorem~\ref{thm:free}).

We now construct a $2$-coloring $z$ on ${\mathbb F}_n$ so that $z(1_{{\mathbb F}_n})=0$ and
$z(wa)=z(w)$ for all $w\in {\mathbb F}_n$. To define such a $z$ it is clearly sufficient to define the values of $z(w)$ for all nonempty words $w\in {\mathbb F}_n$ that do not end in $a$ or $a^{-1}$. Such a word can be uniquely written as $w=u_0a^{p_0}u_1a^{p_1}\dots u_k$, where $k\geq 0$, $u_0,\dots, u_k\in F$, $p_0,\dots,p_{k-1}\in{\mathbb Z}-\{ 0\}$ and $u_1,\dots,u_k\neq 1_F$ if $k>0$. Let $w_1=wu_k^{-1}$. We define $z$ by
$$z(w)= x_{y^*(w_1)}(u_k).
$$
$z$ is clearly right-periodic, hence is not a right $2$-coloring.

We verify that $z$ is a $2$-coloring. Fix
a nonidentity $s\in {\mathbb F}_n$. Let $T_0$ be a finite subset of $F$ so that
for any $h, h'\in F$ there is $t\in T_0$ such that $x_0(ht)\neq x_1(h't)$. Let
$M=\max\{|u|\,:\, u\in T_0\}$. Let
$N$ be a large enough positive integer so that
for any $0<k\leq |s|$ and for any $m\in {\mathbb Z}$, there is $0< l\leq N$ such that
$y(m+l)\neq y(m+k+l)$. Such $N$ exists since $y$ is unidirectional (Lemma~\ref{lem:unidirectional}).
Let
$$T=\{ t\in{\mathbb F}_n\,:\, |t|\leq 2|s|+N+M\}. $$

We claim that $T$ witnesses that $z$ blocks $s$. Let $g\in {\mathbb F}_n$. First notice that there is $s'\in {\mathbb F}_n$ with $|s'|\leq |s|$ such that $|gs'|\neq |gss'|$. In fact, there is such an $s'$ among the initial segments of $s$. Then note that there is a generator $b$ of ${\mathbb F}_n$ (not necessarily distinct from $a$) so that $|gs'b^\epsilon|=|gs'|+1$ and $|gss'b^\epsilon|=|gss'|+1$ for some
$\epsilon\in\{-1,1\}$. Thus for $t_0=s'b^\epsilon$ we have that $|t_0|\leq 2|s|$ and $|gt_0|\neq |gst_0|$.
Next we consider $t_1=t_0a^k$ where $0<|k|\leq N$. There is such a $k$ so that $|t_1|=|t_0|+|k|$ and
$y^*(gt_1)\neq y^*(gst_1)$. Let $i=y^*(gt_1)$ and $i'=y^*(gst_1)$. Then by the
orthogonality of $x_0$ and $x_1$ there is $t=t_1u$ for some nonidentity $u\in F$ such that
$$ z(gt)=x_i(u)\neq x_{i'}(u)=z(gst). $$
Obviously $|t|\leq 2|s|+N+M$.
\end{proof}

\section{$2$-Colorings on solvable groups}

In this section we establish the uniform $2$-coloring property for all countably infinite solvable groups.
We first do this for all countably infinite abelian groups.

If an abelian
group contains at least one element of infinite order then we are done by Corollary~\ref{cor:normaloverZ}. Thus we only need to deal with countably infinite abelian torsion
groups here. There are two concrete situations we need to discuss before coming back to the general
argument.

The first situation concerns a direct sum of infinitely many finite groups. Let $H_0, H_1, \dots, H_n,\dots$ be nontrivial finite groups
and $H=\oplus_n H_n$. We show that $H$ has the uniform $2$-coloring property.

\begin{lem} \label{lem:productfinite}
Let $\pi\in 2^\N$ be such that $0,1\not\in\overline{[\pi]}$. For any $h\in H$ define
$$c_\pi(h)=\left\{\begin{array}{ll} 0, & \mbox{ if $h=1_H$,} \\
\pi(n), & \mbox{ if $h\neq 1_H$ and $n\in\N$ is the least such that $h_n\neq 1_{H_n}$.}
\end{array}\right.
$$
Then $c_\pi$ is a $2$-coloring on $H$. Moreover, if $\pi_0\neq\pi_1$ and $0,1\not\in \overline{[\pi_0]}, \overline{[\pi_1]}$, then $c_{\pi_0}\bot\, c_{\pi_1}$.
\end{lem}

\begin{proof} First it is easily seen that $0,1\not\in\overline{[\pi]}$ iff there is $b\in\N$ such that for any $n\in\N$ there is $m<b$ with $\pi(n)\neq\pi(n+m)$. We will use this equivalence below without elaboration.

Let $s\in H$ with $s\neq 1_H$. Let $n_s$ be the least $n$ such that $s_n\neq 1_{H_n}$. Let $T=\oplus_{n\leq n_s+b}H_n$.
Now suppose $h\in H$ is arbitrary. Let $t_0=\oplus_{n\leq n_s} h_n^{-1}$. Then for all $n\leq n_s$, $(ht_0)_n=1_{H_n}$.
Similarly, for all $n<n_s$, $(hst_0)_n=h_ns_n(t_0)_n=h_nh_n^{-1}=1_{H_n}$. However, $(hst_0)_{n_s}=h_{n_s}s_{n_s}h_{n_s}^{-1}\neq 1_{H_{n_s}}$. Note that for any $t_1\in H$ with $(t_1)_n=1_{H_n}$ for all $n\leq n_s$, $c_\pi(hst_0t_1)=c_\pi(hst_0)$. Now if $c_\pi(ht_0)\neq c_\pi(hst_0)$ we are done since $t_0\in T$. Suppose $c_\pi(ht_0)=c_\pi(hst_0)$. By the assumption
on $\pi$ there is $m<b$ such that $\pi(n_s+1)\neq \pi(n_s+1+m)$. We consider two cases. Case 1: $\pi(n_s+1)\neq c_\pi(ht_0)$. In this case let $k_{n_s+1}\in H_{n_s+1}$ be any nonidentity element and let $t_1=k_{n_s+1}$. Then
$n_s+1$ is the least $n$ so that $(ht_0t_1)_n\neq 1_{H_n}$. Hence $c_\pi(ht_0t_1)=\pi(n_s+1)\neq c_\pi(hst_0)=c_\pi(hst_0t_1)$. Thus $t=t_0t_1$ is as required. Case 2: $\pi(n_s+1+m)\neq c_\pi (ht_0)$. Let $t_1=\oplus_{n_s+1\leq n< n_s+1+m}
h_n^{-1}\oplus k_{n_s+1+m}$ where $k_{n_s+1+m}\in H_{n_s+1+m}$ is an arbitrary element $\neq h_{n_s+1+m}^{-1}$. Then
$c_\pi(ht_0t_1)=\pi(n_s+1+m)\neq c_\pi(hst_0)=c_\pi(hst_0t_1)$. Note that $n_s+1+m\leq n_s+b$, thus $t=t_0t_1\in T$
is as required. This shows that $c_\pi$ is a $2$-coloring.

Now suppose $\pi_0\neq\pi_1$ and $0,1\not\in\overline{[\pi_0]},\overline{[\pi_1]}$, and let the witness be $b_0$ and $b_1$. Let $b_2$ be the least $n$ such that $\pi_0(n)\neq \pi_1(n)$.
Let $b=b_0+b_1+b_2$ and $T=\oplus_{n\leq b} H_n$. Then we claim that for any $g_0, g_1\in H$ there is $t\in T$ such that $c_{\pi_0}(g_0t)\neq c_{\pi_1}(g_1t)$.
Let $g_0, g_1\in H$. We consider two cases. Case 1: $(g_0)_i=(g_1)_i$ for all $i\leq b_2$. Then let $t\in \oplus_{n\leq b_2}H_n\subseteq T$ be such that
$(g_0t)_i=1_{H_i}$ for all $i<b_2$ and $(g_0t)_{b_2}\neq 1_{H_{b_2}}$. Then the same is true for $g_1t$, and thus $c_{\pi_0}(g_0t)=\pi_0(b_2)\neq \pi_1(b_2)=
c_{\pi_1}(g_1t)$. Case 2: $(g_0)_i\neq (g_1)_i$ for some $i\leq b_2$. Then let $t_0\in \oplus_{n\leq b_2}H_n\subseteq T$ be such that for some $i\leq b_2$,
$(g_0t_0)_i=1_{H_i}\neq (g_1t_0)_i$ and for all $j<i$, $(g_0t_0)_j=1_{H_j}=(g_1t_0)_j$. If $c_{\pi_0}(g_0t_0)\neq c_{\pi_1}(g_1t_0)$ there is nothing more to prove.
Otherwise, note that $c_{\pi_1}(g_1t_0)=\pi_1(i)$ for the above mentioned $i\leq b_2$ and $c_{\pi_0}(g_0t_0)=\pi_0(k)$ for some $k>i$. Since $0,1\not\in\overline{[\pi_0]}$,
there is $m<b_0$ such that $\pi_0(k)\neq \pi_0(k+m)$. Thus there is $t_1\in \oplus_{k\leq n\leq k+m}H_n\subseteq T$ such that $c_{\pi_0}(g_0t_0t_1)=\pi_0(k+m)$.
But then $c_{\pi_1}(g_1t_0t_1)=\pi_1(i)\neq \pi_0(k+m)=c_{\pi_0}(g_0t_0t_1)$, so $t=t_0t_1$ is as required. This completes the proof of the lemma.
\end{proof}

\begin{theorem} \label{thm:productfinite} Let $H_0, H_1, \dots, H_n, \dots$ be nontrivial finite groups and $H=\oplus_n H_n$. Then $H$ has the uniform $2$-coloring property.
\end{theorem}

\begin{proof}
Let $\{x_\alpha\,:\,\alpha\in 2^\N\}$ be the collection of $2$-colorings on ${\mathbb Z}$
constructed in the proof of Theorem~\ref{thm:coloringZ}.
Then each $\pi_\alpha=x_{\alpha}\upharpoonright\N$ satisfies $0,1\not\in\overline{[\pi_\alpha]}$.

By the proof of the above lemma,
for any $s\in H$ with $s\neq 1_H$, the witnessing set $T$ for the blocking of $s$ by $c_{\pi_\alpha}$ only depends on $s$ and not on $\alpha$. This shows that the collection $\{c_{\pi_\alpha}\,:\, \alpha\in 2^\N\}$ satisfies Definition~\ref{def:uniformproperty} (i).

To check that $\{c_{\pi_\alpha}\,:\, \alpha\in 2^\N\}$ satisfies Definition~\ref{def:uniformproperty} (ii), we let $A_n\subseteq\N$ be
 given by Theorem~\ref{thm:coloringZ} such that for all $\alpha,\beta\in 2^\N$ with $\alpha(n)\neq\beta(n)$, we have
 $$ \forall g_0,g_1\in {\mathbb Z}\ \exists a\in A_n\ x_\alpha(g_0+a)\neq\beta(g_1+a).$$
  Note that we could take $A_n\subseteq\N$ because of Lemma~\ref{lem:positiveorthogonal}.
  Let $b_n=\max A_n$. Then in particular for any $\alpha,\beta$ as above, there is some $m<b_n$ such that
  $x_\alpha(m)\neq x_\beta(m)$. By the proof of the above lemma, if we let $T_n=\oplus_{m\leq b_n+8}H_m$, then
  $$ \forall h_0,h_1\in H\ \exists t\in T_n\ c_{\pi_\alpha}(h_0t)\neq c_{\pi_\beta}(h_1t). $$
  Since $T_n$ does not depend on $\alpha$ and $\beta$, our proof is complete.
\end{proof}

Next we consider the quasicyclic group ${\mathbb Z}(p^\infty)$ for any prime $p$.

\begin{theorem} \label{thm:quasicyclic} Let $p$ be a prime number. Then ${\mathbb Z}(p^\infty)$
 has the uniform $2$-coloring property.
\end{theorem}

\begin{proof} Every element $g$ of ${\mathbb Z}(p^\infty)$ can be expressed as
$$ \gamma(a_0,\dots,a_{N-1})=\displaystyle\frac{a_0}{p}+\displaystyle\frac{a_1}{p^2}+\displaystyle\frac{a_2}{p^3}+\cdots+
\displaystyle\frac{a_{N-1}}{p^N} $$
for some $N\geq 0$ and $0\leq a_n<p$ for $n=0,\dots,N-1$. For notational convenience we denote $g(n)=a_n$ for $n=0,\dots,N-1$, and more generally, for $n\geq N$, let $g(n)=0$. Now for $g\in {\mathbb Z}(p^\infty)$ let $n_g$ be the least $n$
such that $g(n_g)\neq 0$. Then similar to the proof of Lemma~\ref{thm:productfinite}, we have the following claim.
\begin{quote}
Let $\pi\in 2^\N$ be such that $0,1\not\in\overline{[\pi]}$.
For any $g\in {\mathbb Z}(p^\infty)$ define $c_\pi(g)=\pi(n_g)$. Then $c_\pi$ is a $2$-coloring on ${\mathbb Z}(p^\infty)$.
Moreover, if $\pi_0\neq\pi_1$ and $0,1\not\in\overline{[\pi_0]},\overline{[\pi_1]}$, then $c_{\pi_0}\bot\, c_{\pi_1}$.
\end{quote}
The proof is also similar. In fact, let $s\in {\mathbb Z}(p^\infty)$ so that $s\neq 0$. Let $T=\{ t\in {\mathbb Z}(p^\infty)\,:\, t(n)=0 \mbox{ for all $n>n_s+b$} \}$. Then for all $g\in {\mathbb Z}(p^\infty)$, let $t_0=-\gamma(g\upharpoonright (n_s+1))$. We have that $n_{g+t_0}>n_s$. Thus for any $t_1$ with $n_{t_1}>n_s$, $c_\pi(g+s+t_0+t_1)=c_\pi(g+s+t_0)$,
whereas for some such $t_1$ with $n_{t_1}\leq n_s+b$, we can arrange that $c_\pi(g+t_0+t_1)\neq c_\pi(g+s+t_0)$.
Hence if we let $t=t_0+t_1$ then $c_\pi(g+t)\neq c_\pi(g+s+t)$.

The rest of the proof is similar to that of Theorem~\ref{thm:productfinite}.
\end{proof}

Now we are ready to establish the uniform $2$-coloring property for all countably infinite abelian
groups. As noted before we only need to deal with the torsion case. Also recall that any abelian group can be written as the direct sum of its maximal divisible subgroup and a reduced subgroup.
In the case of a divisible group there is at least one prime $p$ such that the quasicyclic group ${\mathbb Z}(p^\infty)$ is contained in the group.

\begin{theorem}\label{thm:abelian} Let $G$ be a countably infinite abelian group. Then $G$ has the
uniform $2$-coloring property.
\end{theorem}

\begin{proof}
Assume that $G$ is a torsion group. If $G$ has a nontrivial divisible subgroup then there is some prime $p$ such
that ${\mathbb Z}(p^\infty)\unlhd G$. In this case we are done by the preceding theorem and Theorem~\ref{thm:uniformextension}. Suppose $G$ is reduced. We consider
two cases. Case 1: There are infinitely many prime $p$ for which there exist elements of order $p$. In this case
let $p_0, p_1,\dots, p_n,\dots$ be distinct prime numbers and $g_0,g_1,\dots, g_n,\dots$ be nonzero elements so that
$p_ng_n=0$. Then $H=\langle g_0,g_1,\dots,g_n,\dots\rangle$ is isomorphic to the direct sum $\oplus_n {\mathbb Z}_{p_n}$.
Since $H\unlhd G$, by Theorem~\ref{thm:productfinite} and Theorem~\ref{thm:uniformextension} we have that $G$ has the uniform $2$-coloring property.
Case 2: There are only finitely many primes $p$ so that $G$ has a nontrivial $p$-component. Let $G_p$ be the $p$-component
of $G$, i.e., the subgroup of all elements of $G$ whose order is a power of $p$. Let $p_0,\dots, p_n$ be all primes such
that $G_{p_i}$ is nontrivial. Then $G=\oplus_{i\leq n} G_{p_i}$. Thus at least one of $G_{p_i}$ is infinite. Fix such a
$p$. Since we assume that $G$ is reduced, we claim that there are infinitely many elements in $G_p$ with order $p$.
In fact, define a partial order $<$ defined on $G_p$ by $h<g$ iff there is $k\geq 1$ such that $p^kh=g$.
Then since $G_p$ is reduced, $<$ is a wellfounded tree on $G_p$, i.e., there is no infinite $<$-descending sequence in $G_p$. If there are only finitely many elements of order $p$ in
$G_p$, then the tree is finite splitting. For this, just note that if $g_0, \dots, g_n,\dots $ are infinitely many distinct
elements with $pg_0=pg_1=\dots=pg_n=\dots$, then for any $n\geq 1$, $p(g_0-g_n)=0$, and thus $g_0-g_1, \dots, g_0-g_n,\dots$
are infinitely many distinct elements of order $p$. It follows by K\"{o}nig's lemma that a finite splitting wellfounded
tree is finite, and thus $G_p$ would be finite if there are only finitely many elements of order $p$.

Finally, suppose there are infinitely many elements of order $p$ in $G_p$. We define by induction a sequence $h_n$
of elements in $G_p$ as follows. Let $h_0$ be any nonzero element of order $p$ in $G_p$. In general, if $h_0,\dots,h_n$ have been defined, then note that $\langle h_0,\dots, h_n\rangle$ is isomorphic to ${\mathbb Z}_p^{n+1}$, hence finite,
and let $h_{n+1}$ be any nonzero element of order $p$ not in $\langle h_0,\dots, h_n\rangle$. Our assumption guarantees
that this construction will not stop at any finite stage. Also, when the infinite sequence $h_0,\dots, h_n,\dots$ is
defined, we have that $\langle h_0,\dots, h_n,\dots\rangle$ is isomorphic to the direct sum $\oplus_n {\mathbb Z}_p$. Now by
Theorem~\ref{thm:productfinite} and Theorem~\ref{thm:uniformextension}, $G$ has the uniform $2$-coloring property, and our theorem is proved.
\end{proof}

Finally we expand the result to all countably infinite solvable groups.

\begin{theorem}\label{thm:solvableuniform} Let $G$ be a countably infinite solvable group. Then $G$ has the uniform $2$-coloring property.
\end{theorem}

\begin{proof}
Suppose $G$ has rank $n\geq 1$ and its derived series are as follows:
$$ G\unrhd G' \unrhd G'' \unrhd \cdots \unrhd G^{(n)}=\{1_G\}. $$
Then for each $i<n$, $G^{(i)}/G^{(i+1)}$ is abelian. Let $n_0$ be the smallest such that $G^{(n_0)}$ is finite. Then $0<n_0\leq n$.
By Theorem~\ref{thm:uniformextension} it suffices to show that $G^{(n_0-1)}$ has the uniform $2$-coloring property. By assumption, $G^{(n_0-1)}/G^{(n_0)}$ is an infinite
abelian group, thus it has the uniform $2$-coloring property by Theorem~\ref{thm:abelian}. If $n_0=n$ then we have that $G^{(n_0-1)}\cong
G^{(n_0-1)}/G^{(n_0)}$ and we are done. If $n_0<n$, we must have that $|G^{(n_0)}|>2$, since otherwise $G^{(n_0-1)}$ is in fact abelian;
thus by Lemma~\ref{lem:finite} $G^{(n_0)}$ has the $(2,2)$-coloring property. Let $y_0$ and $y_1$ be
orthogonal $2$-colorings on $G^{(n_0)}$. Let $\{z_\alpha\,:\, \alpha\in 2^\N\}$ be a collection of $2$-colorings on $G^{(n_0-1)}/G^{(n_0)}$ witnessing its uniform $2$-coloring property. Using the construction
in the proof of Theorem~\ref{thm:transitivecoloring} to define a collection of $2$-colorings on $G^{(n_0-1)}$. It can be easily verified that the resulting $2$-colorings witness
the uniform $2$-coloring property for $G^{(n_0-1)}$.
\end{proof}

\section{$2$-Colorings on residually finite groups}

In this short section we present a final method of constructing $2$-colorings through algebraic methods. We show how to construct a $2$-coloring on any countable residually finite group.

\begin{theorem}
If $G$ is a countable residually finite group, then $G$ has the coloring property.
\end{theorem}

\begin{proof}
If $G$ is finite then it clearly has the coloring property. So suppose that $G$ is countably infinite. Let $(K_n)_{n \in \N}$ be a decreasing sequence of finite index normal subgroups of $G$ with $\bigcap K_n = \{1_G\}$. Such a sequence exists by the definition of residual finiteness. Define $x \in 2^G$ by setting $x(g) = n \mod 2$ where $n$ satisfies $g \in K_n - K_{n+1}$ (and define $x(1_G)$ arbitrarily). We claim that $x$ is a $2$-coloring on $G$. Fix a nonidentity $s \in G$ and let $n$ satisfy $s \in K_n - K_{n+1}$. Let $T_0$ be a set of representatives for the cosets of $K_{n+1}$ in $G$, let $T_1$ be a set of representatives for the cosets of $K_{n+3}$ in $K_{n+1}$, and let $T = T_0 T_1$. Let $g \in G$ be arbitrary. Since $s \not\in K_{n+1}$, $g$ and $gs$ are not in the same coset of $K_{n+1}$. Consequently, by considering the group $G / K_{n+1}$ we see that there is $t_0 \in T_0$ with $g t_0 \in K_{n+1}$ and $g s t_0 \not\in K_{n+1}$. Notice that if $m$ satisfies $g s t_0 \in K_m - K_{m+1}$, then $m < n+1$. By considering the group $K_{n+1} / K_{n+3}$, we see that there are $t_1, t_2 \in T_1 \subset K_{n+1}$ with $g t_0 t_1 \in K_{n+1} - K_{n+2}$ and $g t_0 t_2 \in K_{n+2} - K_{n+3}$. So clearly $x(g t_0 t_1) \neq x(g t_0 t_2)$. Also, since $g s t_0 \in K_m - K_{m+1}$ and the sequence $(K_r)_{r \in \N}$ is decreasing, we have $g s t_0 t_1, g s t_0 t_2 \in K_m - K_{m+1}$ as well (since $t_1, t_2 \in K_{n+1} \subset K_{m+1}$). So $x(g s t_0 t_1) = x(g s t_0 t_2)$. It follows that either $x(g t_0 t_1) \neq x(g s t_0 t_1)$ or else $x(g t_0 t_2) \neq x(g s t_0 t_2)$. Since $t_0 t_1, t_0 t_2 \in T$, we conclude that $x$ is a $2$-coloring on $G$.
\end{proof}

\begin{cor}
All finitely generated abelian groups, all finitely generated nilpotent groups, all polycyclic groups, all countable (real or complex) linear groups, and all countable nonabelian free groups admit a $2$-coloring.
\end{cor}

\begin{proof}
All of these groups are residually finite.
\end{proof}

The discovery of the above construction occurred very late in the developement of this paper. In effect, we did not investigate if one can use constructions similar to the one above to establish that residually finite groups have the uniform $2$-coloring property. Of course, we do prove in Section \ref{SEC UNI 2 COL} that every countably infinite group has the uniform $2$-coloring property. It may be nice though if the methods of this section could show this fact directly for residually finite groups. We leave the resolution of this question to interested readers.

\chapter{Marker Structures and Tilings} \label{CHAP MARKER}

In this chapter we will introduce a general notion of marker structures on countable groups
and study some of their general properties. As an immediate application of this notion
we give in Section~\ref{sec:abmar} another proof that all abelian groups admit a $2$-coloring,
and in fact the proof will be generalized to establish that all FC groups admit a $2$-coloring.
The concept of marker region will be one of the main tools we use for our main results in
future chapters. In the remainder of this chapter we then introduce and study the  related notion
of a \ccc group. We will show that ccc groups include, among others, all nilpotent groups, all polycyclic groups, all residually finite groups, all locally finite groups, and all groups which are free products of nontrivial groups.
The results of this chapter are relatively independent and will not be needed for the
rest of the paper.

\section{Marker structures on groups}
We introduce the general notion of a {\em marker structure} on a countable group $G$,
and introduce also several specializations of this notion. This point of view is crucial
for the main results of this paper to appear in the following chapters.
In Chapter~\ref{chap:basicconstructions}
we gave certain more algebraic arguments which showed that every countable solvable group
has a $2$-coloring, in fact, we showed this in a strong form already
(\cf\  Theorem~\ref{thm:solvableuniform}).
However, we have not been able to push the methods used in those proofs further. In particular,
we have not been able to use them to show that every countable group $G$ admits a $2$-coloring.
For this we seem to need arguments involving a certain more
geometric nature, which leads to the concept of
a marker structure. The concept of marker structures on various types of groups is certainly
not new to this paper, and indeed has been a central notion in ergodic theory and the theory
of equivalence relations for a long time. In particular, it plays a key role in most hyperfiniteness
proofs. The theorem of Weiss (\cf\ \cite{DJK} \cite{JKL})
that all of the Borel actions of the group $\mathbb{Z}^n$ are hyperfinite
uses these concepts in the proof.
Recall a Borel equivalence relation on a Polish space is {\em hyperfinite}
if it can be written as  an increasing union of finite
Borel equivalence relations. The more recent proof of
Gao and Jackson \cite{GJ} that all Borel actions of any countable abelian group are also hyperfinite
makes use of even better marker structures on these groups. Indeed, the best known results
on the hyperfiniteness problem (determining which groups have only hyperfinite Borel actions)
involve carefully examining the nature of the marker structures that can be put on such groups.
Although this is an interesting connection, the arguments of this paper do not require familiarity
with the notion or theory of hyperfinite equivalence relations.

In the main results of this paper the existence of certain carefully
controlled marker structures is also of central importance. However,
for many of our results the point of view is somewhat different. We
are often interested now in what marker structure can be put on {\em
  arbitrary} countable groups. Of course, as we restrict the class of
groups, we expect marker structures with better properties. The point
we wish to emphasize is that there is a common thread between many of these
other arguments (such as hyperfiniteness proofs) and the arguments of
this paper, and this is what we abstract into the notion of a marker
structure. Various specializations of this notion result in
interesting concepts which have been studied on their own, such as the
class of MT groups defined independently by Chou \cite{Ch} and Weiss \cite{W} in their study of monotileable amenable groups.  Although the marker structures we use
in our main results can be put on any group, it is nonetheless
interesting to ask exactly which types of structures can be put on
various groups. For example, we will introduce the concept of a
ccc tiling of a group. Some very basic questions about which groups
admit such marker structures remain open.

We next give the general notion of a marker structure and various specializations
of the concept. We first use this concept to give a completely different proof that
all of the abelian, and then all of the FC groups admit a $2$-coloring. The proofs of this section
have a decidedly more geometric flavor than the previous arguments; this seems to be inherent in the
concept of a marker structure. Part of the reason for presenting these proofs is that they
foreshadow the more involved arguments necessary for general groups. Indeed, the short proofs for abelian
and FC groups to follow can be seen as a rough outline of the procedure for general groups,
with some of the key technical difficulties removed. We then introduce the strongest notion of
marker structure which seems relevant for the type of constructions one might do along these lines,
and this leads to the notion of a ccc tiling. Again, we will not prove these exist on arbitrary
groups (nor do we need to for our main results), but it becomes an interesting independent
question as to when these exist. As we said above, this is likely related to other questions
such as the hyperfiniteness problem.

We point out two technical distinctions before giving the actual definitions. First,
for the kinds of arguments we do we are mainly interested in not a single
marker structure (defined below), but a sequence of such structures. This is generally also the case
in arguments from ergodic theory as well as hyperfiniteness theory. Second, for the results of this
paper we are interested in the marker structures on the groups themselves as opposed to
on some Polish space on which the groups act. This is in contrast to many of the arguments in
ergodic theory and descriptive set theory where marker arguments occur. In putting marker structures
on the groups themselves, there is no issue of definability that enters in as in the Polish
space case. Thus it becomes easier, at least in theory, to put such structures on the group.
So, the inability to put a type of marker structure on a group puts an upper-bound
on what one can do with the equivalence relation defined by a Borel action of the group
(at least if the action of the group is free). Again, this gives an independent interest to
questions about marker structures on groups.

\begin{definition} \label{def:markerstructure}
Let $G$ be a countable group. A {\em marker structure} on $G$ is a pair
$(\Delta, \sR)$ where $\Delta \subseteq G$, $\sR \subseteq \sP(G)$
satisfying:
\begin{enumerate}
\item
$\sR$ is a pairwise disjoint collection.
\item
For every $R \in \sR$, $|\Delta \cap R|=1$.
\item
Every $\delta \in \Delta$ lies in some $ R \in \sR$.
\item
The set $\bigcup_{\delta \in \Delta} \delta^{-1} R_\delta$ is finite, where
$R_\delta$ denotes the unique $R \in \sR$ with $\delta \in R$.
\end{enumerate}
\end{definition}\index{marker structure}\index{marker region}\index{marker point}

We call the elements $\delta \in \Delta$ the {\em marker points}, and the
sets $R \in \sR$ the {\em marker regions}.

The definition of marker structure in Definition~\ref{def:markerstructure} is quite general.
It encompasses
all of the marker constructions of this paper as well as all of the known hyperfiniteness proofs.
For the constructions of this paper, however, we are usually interested in marker structures with
additional properties. The next definition records some of these additional properties,

\begin{definition}
A marker structure $(\Delta,\sR)$ is {\em regular}
if there is a single (necessarily finite) $F \subseteq G$ such that for
all $\delta \in \Delta$ we have $\delta^{-1} R_\delta=F$, where $R_\delta$ is the unique element of
$\sR$ which contains $\delta$. A marker structure $(\Delta, \sR)$ is {\em centered}
if $1_G \in \Delta$.
A marker structure $(\Delta,\sR)$ is  {\em total} if $G=\bigcup \sR$.
A marker structure is a {\em tiling} if it is regular and total.
\end{definition}\index{marker structure!regular}\index{marker structure!centered}\index{marker structure!total}\index{tiling}\index{tile}\index{regular}\index{centered}\index{total}

In the case of a regular marker
structure $(\Delta,\sR)$, we usually present the marker structure as
$(\Delta, F)$, where $F$ is the common value of $\delta^{-1}R_\delta$
for $\delta \in \Delta$. Thus, the marker regions are the sets of the form
$R=\delta F$ for some $\delta \in \Delta$.
Note that for a regular marker structure $(\Delta,F)$ we necessarily have
$1_G \in F$ since $1_G= \delta^{-1} \delta$ and $\delta \in R_\delta$.
Conversely, given a $\Delta \subseteq G$ and a finite $F \subseteq G$ with
$1_G \in F$, if we set $\sR=\{ \delta F \colon \delta \in \Delta\}$,
then $(\Delta,\sR)$ is a regular marker structure iff
whenever $\delta_1 \neq \delta_2$ are in $\Delta$, then
$\delta_1 F \cap \delta_2 F =\emptyset$. Notice that this is stronger than just requiring
that $\{ \delta F \colon \delta \in \Delta\}$ forms a pairwise disjoint
collection (since the latter condition allows for the possibility that
$\delta_1 F=\delta_2 F$ for some $\delta_1\neq \delta_2$ in $\Delta$).

In the case of a tiling, the $\delta
F=\delta(\delta^{-1}R_\delta)=R_\delta$ partition the group $G$, where
again $F$ is the common value of $\delta^{-1} R_\delta$. In fact, a regular
marker structure $(\Delta,F)$ is a tiling iff $\bigcup_{\delta \in \Delta}
\delta F=G$.
We call $F$ the
{\em tile} for the tiling $\sR$. In particular, we usually also present tilings as
$(\Delta,F)$, where $F$ is the tile.  Total marker structures occur in
hyperfiniteness proofs, but for most of the main results of this paper
we must work with marker structures which are not total. Nevertheless,
we establish certain strong tiling properties for classes of groups in
this chapter.

In the main arguments of this paper, and also in hyperfiniteness
proofs, one needs to consider not just a single marker structure, but a
sequence $(\Delta_n, \sR_n)$ of marker structures. We introduce some
more terminology for such sequences.

\begin{definition}
A sequence of marker structures $(\Delta_n,\sR_n)$ is {\em coherent}
if for $k\leq n$ and marker regions $R_k\in \sR_k$, $R_n \in \sR_n$, we have
$R_k \cap R_n \neq \emptyset$ implies $R_k \subseteq R_n$.
A sequence of marker structures $(\Delta_n, \sR_n)$ is {\em cofinal}
if for every finite $A \subseteq \bigcup_n \bigcup \sR_n$ there is an $n \in \N$ such that
for all $m \geq n$ we have $A \subseteq R_m$ for some $R_m \in \sR_m$.
A sequence of  marker structures $(\Delta_n,\sR_n)$ is {\em centered}
if  for each $n$ the marker structure $(\Delta_n,\sR_n)$ is centered
(that is, $1_G \in \Delta_n$ for all $n$).
\end{definition}\index{sequence of marker structures!coherent}\index{sequence of marker structures!cofinal}\index{sequence of marker structures!centered}\index{coherent}\index{cofinal}\index{centered}

As with single marker structures, we usually present sequences $(\Delta_n, \sR_n)$ of regular
marker structures as $(\Delta_n, F_n)$ where $F_n$ is the common value
of $\delta^{-1} R_\delta$ for $\delta \in \Delta_n$ and $R_\delta$ the
unique $R \in \sR_n$ containing $\delta$.

Note that a sequence of tilings $(\Delta_n,F_n)$ is coherent iff every
$n$th level marker region $\delta_n F_n$ is contained in a (unique)
$n+1$st level marker region $\delta_{n+1}F_{n+1}$.  Also, if
$(\Delta_n,F_n)$ is a coherent sequence of tilings then each $n+1$st
level marker region $\delta_{n+1} F_{n+1}$ is a disjoint union of $n$th
level marker regions $\delta F_n$, for some finite set of $\delta \in
\Delta_n$.  Finally, note that for a sequence of total
marker structures $(\Delta_n, \sR_n)$, being cofinal is just saying
that for every finite $A \subseteq G$, for large enough $n$ we have
that $A$ is contained in a single $n$th level marker region $R_n
\in\sR_n$.  Moreover, a sequence of centered regular marker structures
$(\Delta_n,F_n)$ is cofinal iff every finite $A \subseteq \bigcup_n
\bigcup \sR_n$ (where $\sR_n=\{ \delta F_n \colon \delta \in \Delta_n\}$)
is contained in $F_n$ for large enough $n$.  This is
because we may assume $A$ contains $1_G$, and each $F_n$ contains
$1_G$ for a centered tiling. In particular, a sequence $(\Delta_n, F_n)$
of centered tilings is cofinal iff $\bigcup_n F_n=G$.

The next definition gives a name to groups admitting the strongest form of tilings we will consider.

\begin{definition}
A countable group $G$ is a
{\em \ccc group} if $G$ has a coherent, cofinal, centered  sequence of tilings $( \Delta_n, F_n)$.
\end{definition}\index{tiling!ccc}\index{ccc group}

The significance of  \ccc tilings is that they are in some sense the most highly
controlled marker structure we can get on a group. It is easy to see that various simple
groups such as $\mathbb{Z}$ or $\mathbb{Z}^n$ are \ccc groups. The general situation, however, is not
clear, which leads  to the following question.

\begin{question} \label{cccquestion}
Which groups are \ccc groups?
\end{question}\index{MT group}

Recall from \cite{Ch} and \cite{W} the concept of MT groups defined independently by Chou and Weiss. A countable group is called an {\it MT group} if it admits cofinal tilings. Chou and Weiss independently proved that
the class of MT groups is closed under group extensions and that all countable residually finite groups and all countable solvable groups are MT. Chou further proved that any free product of nontrivial groups is MT. Chou and Weiss raised the following question.

\begin{question}[Chou \cite{Ch}, Weiss \cite{W}]\label{MTquestion}
Which groups are MT groups?
\end{question}

The above questions are important for several reasons. The results of this
paper depend heavily on being able to construct sufficiently good
marker regions for a general group.  We suspect that in future
applications of these methods, it may become important to identify
even better classes of marker regions for groups (or some special
families of groups). Aside from the applications to the current paper,
these general questions also arise in other considerations.  For
example, suppose $G$ is a countable group acting in a Borel way on a
standard Borel space $X$. Recall the equivalence relation $E$ on $X$
generated by the action is said to be {\em hyperfinite} if $E$ is the
increasing union of finite Borel sub-equivalence relations $E_n$ on
$X$. That is, $E$ is hyperfinite if we can find (Borel) marker regions
on $X$ whose equivalence classes union to all of $X$. Here the marker
regions are on the Polish space $X$, and not the group $G$. However,
marker regions $R$ for $X$, assuming the action of $G$ on $X$ is free,
easily induce marker regions for $G$ by simply fixing a particular
equivalence class $[x]$ and considering the relation $g \sim h$ iff
$(g \cdot x) \, R \, (h \cdot x)$. The other direction does not go
through, so having marker regions with certain properties on a group
$G$ is in general a weaker assertion that having marker regions with
these properties on $X$. The two questions, though, are certainly
related, and having the regions on $G$ is a necessary condition for
having them on $X$.

We will consider the question of which groups are \ccc groups later in
this chapter. For now we note the simple observation that the
centeredness requirement is mainly for convenience as it can always be
achieved.

\begin{prop} \label{prop:centered}
Every tiling  $(\Delta, F)$ of a group $G$ has a presentation
as a centered tiling. That is, there is a centered
tiling  $(\Delta', F')$ having the same marker regions (i.e., $\{ \delta F \colon \delta \in \Delta\}
=\{ \delta' F' \colon \delta' \in \Delta'\}$). In particular,
if $G$ admits a coherent, cofinal sequence of tilings, then $G$ is a \ccc group.
\end{prop}

\begin{proof}
Suppose $(\Delta, F)$ is a tiling for $G$. Let $\delta \in \Delta$ be such that $1_G \in \delta F$.
Let $\Delta'=\Delta \delta^{-1}$ and $F'= \delta F$. This clearly works.
\end{proof}

\section{$2$-Colorings on abelian and FC groups by markers} \label{sec:abmar}

In this section we use the notion of marker structure to give a proof that all abelian groups admit
a $2$-coloring. This proof is quite different from that of Theorem~\ref{thm:abelian}. This rather simple
proof foreshadows the proof for general groups to be given in Chapters~\ref{CHAP FM} and \ref{chap:7}, and
will serve to motivate some of the later constructions. We then extend the argument slightly to show that
every FC group also admits a $2$-coloring (the definition of an FC group is given below).
This result does not seem to follow from the
methods of the previous chapters. It will also show some of the difficulties associated with the group
being nonabelian, and will give further motivation for the general constructions later.

We remark that we use multiplicative notation throughout, even when the group is abelian.

We will first introduce some notation. For any graph $\Gamma$ let
$V(\Gamma)$ and $E(\Gamma)$ denote the vertices and edges of $\Gamma$
respectively. For any two graphs $\Gamma_1$ and $\Gamma_2$ let
$\Gamma_{1}\cup\Gamma_{2}=(V(\Gamma_{1})\cup
V(\Gamma_{2}),E(\Gamma_{1})\cup E(\Gamma_{2}))$.
\begin{theorem}
If $G$ is a countable abelian group then $G$ has a $2$-coloring.
\end{theorem}

\begin{proof}
When $G$ is finite this is clear, so we may assume $G$ is countably
infinite and let $1_G=g_0, g_1, \ldots, g_n, \ldots$ be an enumeration
of the elements of $G$.

We begin by constructing a sequence
$(F_n)_{n\in\N^+}$ of finite subsets of $G$. First choose $F_1$ such
that $|F_1|\geq 3$ and for some $a_1\in\ F_1$, $a_1  g_1\in F_1$. We
will continue the construction inductively and assume that $F_1, F_2,
\ldots, F_{k-1}$ have been defined for some $k>1$. Choose any
$\lambda_1, \lambda_2, \lambda_3 \in G$ with $\lambda_i  F_{k-1} F^{-1}_{k-1} \cap
\lambda_j  F_{k-1}  F^{-1}_{k-1}=\varnothing$ for
$i\neq j$ $\ i, j\in \{1,2,3\}$ and choose a finite $F_k\subseteq G$
such that $$\{\lambda_1, \lambda_2, \lambda_3, g_k  \lambda_1, g_k
\lambda_2, g_k  \lambda_3\}\,   F^2_{k-1}  F^{-1}_{k-1} \subseteq F_k.$$
\indent Now, for each $n\in\N$ fix $\Delta_n \subseteq G$ such that
$\{\gamma  F_n \,:\, \gamma\in \Delta_n\}$ is a collection of
maximally disjoint translates of $F_n$. Thus, we have defined a sequence
$(\Delta_n,F_n)$ of regular marker structures on $G$. Note that
this sequence is neither coherent nor cofinal.

We claim that for every $n>1$ and $\gamma\in
\Delta_n$ there exist distinct $z_1, z_2, z_3 \in \Delta_{n-1}$ with
$\{z_i,z_i  g_n\}  F_{n-1} \subseteq \gamma  F_n$ for each $i\in
\{1,2,3\}$. By construction there exist $\lambda_1, \lambda_2,
\lambda_3 \in G$ such that for $i, j\in \{1, 2, 3\}$ with $i\neq j$,
$\lambda_i  F_{n-1}  F^{-1}_{n-1} \cap \lambda_j  F_{n-1}
F^{-1}_{n-1}=\varnothing$ and $\{\lambda_i,\lambda_i  g_n\}  F^2_{n-1}
F^{-1}_{n-1} \subseteq F_n$. Since the $\Delta_{n-1}$-translates of
$F_{n-1}$ are maximally disjoint, for some $z_1\in \Delta_{n-1}$, $z_1
 F_{n-1} \cap \gamma \lambda_1  F_{n-1} \neq \varnothing$ and
therefore $z_1 \in \gamma  \lambda_1  F_{n-1}  F_{n-1}^{-1}$. Similarly
we find there exists $z_2\in \Delta_{n-1}\cap \gamma \lambda_2
F_{n-1}  F^{-1}_{n-1}$ and $z_3\in \Delta_{n-1}\cap \gamma
\lambda_3  F_{n-1}  F^{-1}_{n-1}$. Since $\lambda_i  F_{n-1}
F^{-1}_{n-1} \cap \lambda_j  F_{n-1}  F^{-1}_{n-1}=\varnothing$ for
$i\neq j$ $\ i, j\in \{1, 2, 3\}$, $z_1,\ z_2,$ and $z_3$ must be
distinct. Finally, for $i\in \{1, 2, 3\}$, $\{z_i,z_i  g_n\}
F_{n-1} \subseteq \{\gamma  \lambda_i,\gamma  \lambda_i  g_n\}
F^2_{n-1}  F^{-1}_{n-1} \subseteq \gamma  F_n$.


\begin{figure}[ht]
\begin{center}
\setlength{\unitlength}{6mm}
\begin{picture}(20,18)(-1,-8)

\put(16,2){\makebox(0,0)[b]{$\gamma F_n$}}
\qbezier(0,4.5)(0,9)(9,9)
\qbezier(0,4.5)(0,0)(9,0)
\qbezier(9,0)(18,0)(18,4.5)
\qbezier(18,4.5)(18,9)(9,9)

\put(0,0){

\put(-1,-0.2){
\put(3.5,5){\makebox(0,0)[b]{$z_1F_{n-1}$}}
\qbezier(3,4)(3,3)(4.5,3)
\qbezier(3,4)(3,5)(4.5,5)
\qbezier(4.5,3)(6,3)(6,4)
\qbezier(4.5,5)(6,5)(6,4)
\put(4,4){\circle*{0.1}}
\put(5,3.3){\makebox(0,0)[b]{$z_1\!\in\! \Delta_{n-1}$}}
}

\put(0.6,1.0){
\put(3.5,5){\makebox(0,0)[b]{$z_1g_nF_{n-1}$}}
\qbezier(3,4)(3,3)(4.5,3)
\qbezier(3,4)(3,5)(4.5,5)
\qbezier(4.5,3)(6,3)(6,4)
\qbezier(4.5,5)(6,5)(6,4)
\put(4,4){\circle*{0.1}}
\put(4.7,3.8){\makebox(0,0)[b]{$z_1g_n$}}
}

\put(3,3.8){\vector(4,3){1.6}}
}

\put(8,-2){

\put(-1,-0.2){
\put(3.5,5){\makebox(0,0)[b]{$z_3F_{n-1}$}}
\qbezier(3,4)(3,3)(4.5,3)
\qbezier(3,4)(3,5)(4.5,5)
\qbezier(4.5,3)(6,3)(6,4)
\qbezier(4.5,5)(6,5)(6,4)
\put(4,4){\circle*{0.1}}
\put(5,3.3){\makebox(0,0)[b]{$z_3\!\in\! \Delta_{n-1}$}}
}

\put(0.6,1.0){
\put(3.5,5){\makebox(0,0)[b]{$z_3g_nF_{n-1}$}}
\qbezier(3,4)(3,3)(4.5,3)
\qbezier(3,4)(3,5)(4.5,5)
\qbezier(4.5,3)(6,3)(6,4)
\qbezier(4.5,5)(6,5)(6,4)
\put(4,4){\circle*{0.1}}
\put(4.7,3.8){\makebox(0,0)[b]{$z_3g_n$}}
}

\put(3,3.8){\vector(4,3){1.6}}
}

\put(4.3,1.5){

\put(-1,-0.2){
\put(3.5,5){\makebox(0,0)[b]{$z_2F_{n-1}$}}
\qbezier(3,4)(3,3)(4.5,3)
\qbezier(3,4)(3,5)(4.5,5)
\qbezier(4.5,3)(6,3)(6,4)
\qbezier(4.5,5)(6,5)(6,4)
\put(4,4){\circle*{0.1}}
\put(5,3.3){\makebox(0,0)[b]{$z_2\!\in\! \Delta_{n-1}$}}
}

\put(0.6,1.0){
\put(3.5,5){\makebox(0,0)[b]{$z_2g_nF_{n-1}$}}
\qbezier(3,4)(3,3)(4.5,3)
\qbezier(3,4)(3,5)(4.5,5)
\qbezier(4.5,3)(6,3)(6,4)
\qbezier(4.5,5)(6,5)(6,4)
\put(4,4){\circle*{0.1}}
\put(4.7,3.8){\makebox(0,0)[b]{$z_2g_n$}}
}

\put(3,3.8){\vector(4,3){1.6}}
}

\put(13,6){\makebox(0,0)[b]{$\cdots\cdots$}}

\end{picture} \vspace{-130pt}
\caption{\label{fig:A1}The proof of Theorem 4.2.1.}
\end{center}
\end{figure}
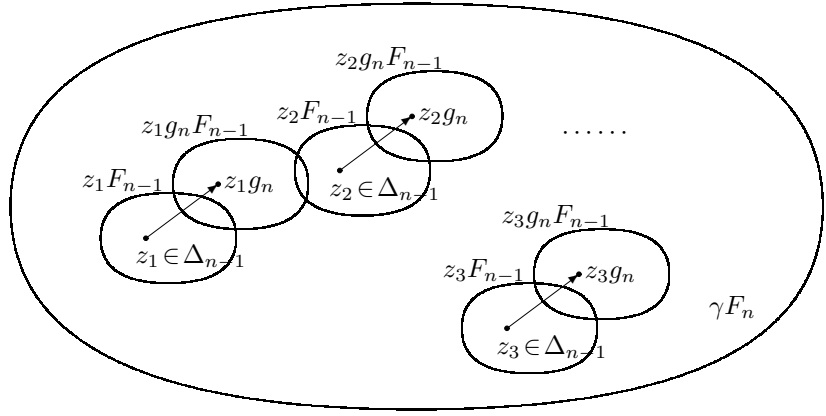

We will now create an increasing  sequence of graphs
$(\Gamma_n)_{n\in\N}$ which we will use to construct a 2-coloring on
$G$. By construction there exists $a_1 \in F_1$ with $a_1  g_1\in
F_1$. Define $\Gamma_1$ to be the graph with edge set
the set of all (undirected) edges between $\gamma a_1$ and
$\gamma a_1 g_1$ for $\gamma \in \Delta_1$. We will write this as
$\bigcup_{\gamma\in\Delta_1} \{(\gamma
a_1,\gamma  a_1  g_1)\}$. Note
that since the $\Delta_1$-translates of $F_1$ are disjoint and
$a_1,a_1  g_1\in F_1$, $\Gamma_1$ is composed of an infinite number
of disconnected components, each of which contains only one edge. On
this note it is clear that $\Gamma_1$ has no cycles. Additionally
since $|F_1|\geq 3$, for all $\gamma \in \Delta_1$ we have  $|(\gamma  F_1)
- V(\Gamma_1)|\geq 1$.

We will continue the construction inductively and
assume that for some $k>1$ $\Gamma_1, \Gamma_2, \ldots, \Gamma_{k-1}$
have been defined such that

\begin{enumerate}
\item[(i)] $\Gamma_{k-1}$ has no cycles;
\item[(ii)] For all $\gamma\in \Delta_{k-1}$, $|(\gamma
F_{k-1}) - V(\Gamma_{k-1})|\geq 1$;
\item[(iii)] For all $i<k$ and $\gamma\in\Delta_i$ there
exists $a\in F_i$ such that $\{\gamma  a,\gamma  a
g_i\}\subseteq V(\Gamma_{k-1})$ and
$(\gamma  a, \gamma  a  g_i)\in E(\Gamma_{k-1})$.
\end{enumerate}

We know that for any $\gamma\in\Delta_k$ there
exist distinct $z_1,z_2,z_3\in\Delta_{k-1}$ with $\{z_i,z_i  g_k\}
F_{k-1} \subseteq \gamma  F_k$ for $i\in\{1,2,3\}$. Therefore for
every $\gamma\in\Delta_k$, $|(\gamma  F_k)-V(\Gamma_{k-1})|\geq 3$ and there exists $a_k^\gamma \in F_k$ such
that $\gamma  a_k^\gamma \in \gamma  F_k -
V(\Gamma_{k-1})$ and $\gamma  a_k^\gamma  g_k\in \gamma  F_k$. We
then define $\Gamma_k$ to be $\Gamma_{k-1}$ together with the edges
in $\bigcup_{\gamma\in\Delta_k}  \{(\gamma  a_k^\gamma,\gamma
a_k^\gamma  g_k)\}$.
Since in each
$\Delta_k$-translate of $F_k$ at most two vertices are being appended
to $\Gamma_{k-1}$ in constructing $\Gamma_k$, we see that for all
$\gamma\in\Delta_k$ that $|\gamma F_k - V(\Gamma_k)|\geq
1$. Additionally, as $\Gamma_{k-1}$ has no cycles and for each
$\gamma\in\Delta_k$ we have $a_k^\gamma\not\in V(\Gamma_{k-1})$, $\Gamma_k$
cannot have any cycles either. Figure~\ref{fig:A1} illustrates our construction.

We let $\Gamma=\bigcup_{n\in\N^+} \Gamma_n$ and claim
that $\Gamma$ has no cycles. Towards a contradiction suppose $\Gamma$
has a cycle. Then the cycle would traverse a finite number of edges,
and therefore for some sufficiently large $n\in\N$ $\ \Gamma_n$ would
contain this cycle. But this is a contradiction since $\Gamma_n$ has
no cycles. Since $\Gamma$ has no cycles it can be 2-colored in the
graph theoretic sense, that is, any two vertices joined by an edge have
different colors. Let $\mu:V(\Gamma)\rightarrow 2$ be a 2-coloring of
$\Gamma$ and let $c:G\rightarrow 2$ be any function such that for all $g\in
V(\Gamma)$ $c(g)=\mu(g)$. We will now show $c$ is a 2-coloring on
$G$. Let $g_i$ be given, let $T=F^2_i  F^{-1}_i$, and let $g\in G$ be
arbitrary. Since the $\Delta_i$-translates of $F_i$ are maximally
disjoint there exists $f_1 \in F_i$ such that $g  f_1 \in \Delta_i
F_i$. It follows there exists $f_2 \in F^{-1}_i$ with $\gamma = g  f_1
f_2 \in \Delta_i$. Finally, there exists $a \in F_i$ such that
$\{\gamma  a, \gamma  a g_i\}\subseteq V(\Gamma_i) \subseteq
V(\Gamma)$ and $(\gamma  a, \gamma  a  g_i)\in E(\Gamma_i)
\subseteq E(\Gamma)$. Therefore $c(g  f_1  f_2  a)=c(\gamma
a)=\mu(\gamma  a) \neq \mu(\gamma  a g_i)=c(\gamma a  g_i)=c(g
f_1  f_2  a  g_i)$. This completes the proof as $f_1  f_2  a \in
F^2_i  F^{-1}_i = T$.
\end{proof}

We now extend the previous argument to FC groups (see Definition \ref{DEFN FCGROUP}).

\begin{theorem}
If $G$ is a countable FC group then $G$ has a $2$-coloring.
\end{theorem}

\begin{proof}
When $G$ is finite this is clear, so we may assume $G$ is countably
infinite and let $1_G=g_0, g_1, \ldots, g_n, \ldots$ be an enumeration
of the elements of $G$. For each $n\in \N$ let $C_n$ denote the finite
conjugacy class containing $g_n$.

We begin by constructing a sequence $(F_n)_{n\in\N^+}$
of finite subsets of $G$. First choose $F_1$ such that $|F_1|\geq
|C_1|+2$ and for some $a_1\in\ F_1$ we have $a_1 C_1 \subseteq F_1$. We will
continue the construction inductively and assume that $F_1, F_2,
\ldots, F_{k-1}$ have been defined for some $k>1$. Let $m=|C_k|$ and
choose any $\lambda_1, \lambda_2, \ldots, \lambda_{m+1} \in G$ such
that for $i, j\in \{1,2,\ldots, m+1\}$ with $i\neq j$
\begin{equation}\label{A1}
\lambda_i F_{k-1} F_{k-1}^{-1} \cap \lambda_j F_{k-1} F_{k-1}^{-1}
= \varnothing,
\end{equation}
\begin{equation}\label{A2}
\lambda_i F_{k-1} F_{k-1}^{-1} F_{k-1} C_k \cap \lambda_j F_{k-1}
F_{k-1}^{-1} F_{k-1} = \varnothing, \text{ and}
\end{equation}
\begin{equation}\label{A3}
\lambda_i F_{k-1} F_{k-1}^{-1} F_{k-1} C_k \cap \lambda_j F_{k-1}
F_{k-1}^{-1} F_{k-1} C_k = \varnothing.
\end{equation}
Finally choose a finite $F_k\subseteq G$ such that for all $i\in
\{1,2,\ldots,m+1\}$
\begin{equation}\label{A4}
\lambda_i F_{k-1} F_{k-1}^{-1} F_{k-1}\cup
\lambda_i F_{k-1} F_{k-1}^{-1} F_{k-1} C_k \subseteq F_k.
\end{equation}
\indent Now, for each $n\in\N^+$ fix $\Delta_n \subseteq G$ such that
$\{\gamma F_n \ |\ \gamma\in \Delta_n\}$ is a maximally disjoint collection
of translates of $F_n$. Thus, we have defined a sequence
$(\Delta_n,F_n)$ of regular marker structures on $G$.

We claim that for every $n\in\N^+$ and $\gamma\in
\Delta_n$ there exist distinct elements $z_1, z_2, \ldots, z_{m+1}\in
\Delta_{n-1}$ where $m=|C_n|$ such that for each $i,j\in
\{1,2,\ldots,m+1\}$ with $i\neq j$
\begin{equation}\label{B1}
z_i F_{n-1} C_n \cap z_j F_{n-1}=\varnothing,
\end{equation}
\begin{equation}\label{B2}
z_i F_{n-1} C_n \cap z_j F_{n-1} C_n=\varnothing, \text{ and}
\end{equation}
\begin{equation}\label{B3}
z_i F_{n-1}\cup z_i F_{n-1} C_n \subseteq \gamma F_n.
\end{equation}
By construction there exist $\lambda_1, \lambda_2, \ldots,
\lambda_{m+1} \in G$ satisfying (\ref{A1}) through (\ref{A4}). Since
the $\Delta_{n-1}$-translates of $F_{n-1}$ are maximally disjoint, for
some $z_1\in \Delta_{n-1}$ $z_1 F_{n-1} \cap \gamma\lambda_1
F_{n-1} \neq \varnothing$ and therefore $z_1 \in \gamma\lambda_1
F_{n-1} F_{n-1}^{-1}$. Similarly we find there exists $z_i\in
\Delta_{n-1}\cap \gamma\lambda_i F_{n-1} F_{n-1}^{-1}$ for each
$i\in \{2,3,\ldots,m+1\}$. It follows that for each $i, j\in \{1, 2,
\ldots, m+1\}$ with $i\neq j$ that $z_i$ and $z_j$ are distinct since
$\lambda_i F_{n-1} F_{n-1}^{-1} \cap \lambda_j F_{n-1}
F_{n-1}^{-1}=\varnothing$. Finally, for $i\in \{1, 2, \ldots, m+1\}$,
$z_i F_{n-1} \subseteq \gamma\lambda_i F_{n-1} F_{n-1}^{-1}
F_{n-1}$ and $z_i F_{n-1} C_n \subseteq \gamma\lambda_i F_{n-1}
F_{n-1}^{-1} F_{n-1} C_n$, so properties (\ref{B1}) through (\ref{B3})
follow from $\lambda_1, \lambda_2, \ldots, \lambda_{m+1}$ satisfying
(\ref{A2}) through (\ref{A4}).

\begin{figure}[ht]
\begin{center}
\setlength{\unitlength}{6mm}
\begin{picture}(30,18)(-1,-8)

\put(18,2){\makebox(0,0)[b]{$\gamma F_n$}}
\qbezier(0,4.5)(0,9)(10,9)
\qbezier(0,4.5)(0,0)(10,0)
\qbezier(10,0)(20,0)(20,4.5)
\qbezier(20,4.5)(20,9)(10,9)

\put(0,0){

\put(-1,-0.2){
\put(3.5,5){\makebox(0,0)[b]{$z_1F_{n-1}$}}
\qbezier(3,4)(3,3)(4.5,3)
\qbezier(3,4)(3,5)(4.5,5)
\qbezier(4.5,3)(6,3)(6,4)
\qbezier(4.5,5)(6,5)(6,4)
\put(4,4){\circle*{0.1}}
\put(5,3.3){\makebox(0,0)[b]{$z_1$}}
}

\put(3,3.8){\vector(4,3){1.6}}
\put(4.6, 5){\circle*{0.1}}
\put(4.7,5.2){\makebox(0,0)[b]{$z_1h_1$}}

\put(3,3.8){\vector(3,1){3.0}}
\put(6,4.8){\circle*{0.1}}
\put(6.2,4.9){\makebox(0,0)[b]{$z_1h_2$}}

\put(3,3.8){\vector(1,0){3.3}}
\put(6.3,3.8){\circle*{0.1}}
\put(7,3.6){\makebox(0,0)[b]{$z_1h_3$}}

\put(6.5,2.8){\makebox(0,0)[b]{$\cdots$}}

\put(3,3.8){\vector(2,-1){2}}
\put(5,2.8){\circle*{0.1}}
\put(5.2,2.1){\makebox(0,0)[b]{$z_1h_m$}}
}

\put(8,-1){

\put(-1,-0.2){
\put(3.5,5){\makebox(0,0)[b]{$z_iF_{n-1}$}}
\qbezier(3,4)(3,3)(4.5,3)
\qbezier(3,4)(3,5)(4.5,5)
\qbezier(4.5,3)(6,3)(6,4)
\qbezier(4.5,5)(6,5)(6,4)
\put(4,4){\circle*{0.1}}
\put(5,3.3){\makebox(0,0)[b]{$z_i$}}
}

\put(3,3.8){\vector(4,3){1.6}}
\put(4.6, 5){\circle*{0.1}}
\put(4.7,5.2){\makebox(0,0)[b]{$z_ih_1$}}

\put(3,3.8){\vector(3,1){3.0}}
\put(6,4.8){\circle*{0.1}}
\put(6.2,4.9){\makebox(0,0)[b]{$z_ih_2$}}

\put(3,3.8){\vector(1,0){3.3}}
\put(6.3,3.8){\circle*{0.1}}
\put(7,3.6){\makebox(0,0)[b]{$z_ih_3$}}

\put(6.5,2.8){\makebox(0,0)[b]{$\cdots$}}

\put(3,3.8){\vector(2,-1){2}}
\put(5,2.8){\circle*{0.1}}
\put(5.2,2.1){\makebox(0,0)[b]{$z_ih_m$}}
}

\put(13,6){\makebox(0,0)[b]{$\cdots\cdots$}}

\end{picture}  \vspace{-130pt}
\caption{\label{fig:A2}The proof of Theorem 4.2.2.}
\end{center}
\end{figure}

We will now create an increasing  sequence of graphs
$(\Gamma_n)_{n\in\N^+}$ which we will use to construct a 2-coloring on
$G$. By construction there exists $a_1 \in F_1$ with $a_1 C_1\subseteq
F_1$. Define $\Gamma_1$ to be the graph with (undirected) edges
$\bigcup_{\gamma\in\Delta_1}\bigcup_{h\in\
C_1} \{ (\gamma a_1, \gamma h a_1)\}$.
Note that $\gamma h a_1 = \gamma a_1 a_1^{-1} h a_1 \in
\gamma a_1 C_1 \subseteq \gamma F_1$ and since the
$\Delta_1$-translates of $F_1$ are disjoint, $\Gamma_1$ is composed of
an infinite number of disconnected components. It is clear from the
construction that $\Gamma_1$ has no cycles. Additionally since
$|F_1|\geq |C_1|+2$, for all $\gamma \in \Delta_1$ we have $|\gamma F_1 -
V(\Gamma_1)|\geq 1$.

We will continue the construction inductively and
assume that for some $k>1$ $\Gamma_1, \Gamma_2, \ldots, \Gamma_{k-1}$
have been defined such that
\begin{enumerate}
\item[(i)]
$\Gamma_{k-1}$ has no cycles;
\item[(ii)] For all $\gamma\in \Delta_{k-1}$, $|\gamma
F_{k-1}-V(\Gamma_{k-1})|\geq 1$;
\item[(iii)] For all $i<k$, $\gamma\in\Delta_i$, and $h\in
C_i$ there exists $a\in F_i$ such that $\{\gamma a,\gamma h
a\}\subseteq V(\Gamma_{k-1})$ and
$(\gamma a, \gamma h a)\in E(\Gamma_{k-1})$.
\end{enumerate}

Enumerate the members of $C_k$ as
$h_1,h_2,\ldots,h_m$. We know that for each $\gamma\in\Delta_k$ there
exist distinct $z_1^\gamma,z_2^\gamma, \ldots, z_{m+1}^\gamma
\in\Delta_{k-1}$ satisfying (\ref{B1}) through (\ref{B3}). For each
$i\in \{1,2,\ldots,m\}$ and each $\gamma\in\Delta_k$ fix $a_i^\gamma
\in F_{k-1}$ such that $z_i^\gamma a_i^\gamma \not\in
V(\Gamma_{k-1})$. Then define $\Gamma_k$ to be the graph $\Gamma_{k-1}$
together with the edges between the points $z_i^\gamma a_i^\gamma$ and
$\gamma h_i \gamma^{-1} z_i^\gamma a_i^\gamma$. That is, we set
$$E(\Gamma_k)=E(\Gamma_{k-1}) \cup \bigcup_{\gamma\in\Delta_k} \bigcup_{1\leq i \leq
m} \{(z_i^\gamma a_i^\gamma, \gamma h_i \gamma^{-1} z_i^\gamma
a_i^\gamma)\}.$$
Since $\gamma
h_i \gamma^{-1} z_i^\gamma a_i^\gamma = z_i^\gamma a_i^\gamma ((z_i^\gamma a_i^\gamma)^{-1}
\gamma h_i \gamma^{-1} z_i^\gamma a_i^\gamma )\in z_i^\gamma F_{k-1} C_k$, it follows
from the definition of the $z_i^\gamma$'s that
\[ z_1^\gamma a_1^\gamma,\,
\gamma h_1 \gamma^{-1} z_1^\gamma a_1^\gamma, \, z_2^\gamma a_2^\gamma, \gamma
h_2 \gamma^{-1} z_2^\gamma a_2^\gamma, \ldots, z_m^\gamma a_m^\gamma, \, \gamma
h_m \gamma^{-1} z_m^\gamma a_m^\gamma\] are all distinct and lie in $\gamma
F_k$ for each $\gamma\in\Delta_k$. Since $\Gamma_{k-1}$ has no cycles
and for each $\gamma\in\Delta_k$ and $i\in \{1,2,\ldots,m\}$
$z_i^\gamma a_i^\gamma \not\in V(\Gamma_{k-1})$ it follows that
$\Gamma_k$ has no cycles either. Additionally for each
$\gamma\in\Delta_k$ and $i\in\{1,2,\ldots,m\}$ $z_i^\gamma a_i^\gamma,
\gamma h_i \gamma^{-1} z_i a_i^\gamma \not\in z_{m+1} F_{k-1}$ by
(\ref{B1}) therefore $|\gamma F_k-V(\Gamma_k)|\geq 1$. The third requirement
on the induction follows as well when we set $a=\gamma^{-1} z_i^\gamma
a_i^\gamma$. Figure~\ref{fig:A2} illustrates our construction.

We let $\Gamma=\bigcup_{n\in\N^+} \Gamma_n$. As before, $\Gamma$ has no cycles
and so can be $2$-colored in the graph theoretic sense.
Let $\mu\colon    V(\Gamma)\rightarrow 2$ be a 2-coloring of
$\Gamma$ and let $c:G\rightarrow 2$ be any function such that for all $g\in
V(\Gamma)$ $c(g)=\mu(g)$. We will now show $c$ is a 2-coloring on
$G$. Let $g_i$ be given, let $T=F_i F_i^{-1} F_i$, and let $g\in G$ be
arbitrary. Since the $\Delta_i$-translates of $F_i$ are maximally
disjoint there exists $f_1 \in F_i$ such that $g f_1 \in \Delta_i
F_i$. It follows there exists $f_2 \in F_i^{-1}$ with $\gamma = g f_1
f_2 \in \Delta_i$. Then we have $g g_i f_1 f_2 = \gamma (f_1 f_2)^{-1}
g_i (f_1 f_2) = \gamma h$ for some $h\in C_i$. Finally, by
construction there exists $a \in F_i$ such that $\{\gamma a, \gamma h
a\} \subseteq V(\Gamma_i) \subseteq V(\Gamma)$ and $(\gamma a, \gamma
h a)\in E(\Gamma_i) \subseteq E(\Gamma)$. Therefore $c(g f_1 f_2
a)=c(\gamma a)=\mu(\gamma a) \neq \mu(\gamma h a)=c(\gamma h a)=c(g
g_i f_1 f_2 a)$. This completes the proof as $f_1 f_2 a \in F_i
F_i^{-1} F_i = T$.
\end{proof}

\section{Some properties of \ccc groups}

In the remainder of this chapter we study \ccc groups. In this section we first establish some basic properties of \ccc groups.
Recall that a countable group $G$ is a \ccc group if it has a sequence of tilings
$(\Delta_n, F_n)$ which are coherent, cofinal, and centered. As we noted
in Proposition~\ref{prop:centered} there is no loss of generality in assuming
centeredness.

We first prove a general lemma which shows that
starting from a \ccc  tiling $(\Delta_n, F_n)$
of the group $G$ we may modify it to get another \ccc tiling $(\tDelta_n, F_n)$
(using the same tiles $F_n$) with some additional uniformity properties.  Recall that in a \ccc
tiling $(\Delta_n,F_n)$, every $F_n$ is a finite disjoint union of
translates $\delta F_{n-1}$ for $\delta \in \Delta_{n-1}$ (since $F_n =1_G F_n$ is
one of the sets in the partition corresponding to $(\Delta_n$, $F_n$)).

\begin{lem} \label{imtiling}
Let $(\Delta_n, F_n)$ be a \ccc tiling of the group $G$.
Then there is a \ccc tiling $(\tDelta_n,F_n)$ of $G$ satisfying the following.
For each $n$, let
$\tDelta^n_{n-1}=\{ \delta \in \tDelta_{n-1} \colon \delta F_{n-1} \subseteq F_n\}$.
Then
\begin{enumerate}
\item
$\tDelta_n= \bigcup_{m>n} \tDelta^m_{m-1} \tDelta^{m-1}_{m-2} \cdots
\tDelta^{n+1}_n$;
\item
$F_n=\tDelta^n_{n-1} \tDelta^{n-1}_{n-2} \cdots \tDelta^1_0 F_0$;
\item
$\tDelta_{n+1} \subseteq \tDelta_n$.
\end{enumerate}
\end{lem}

\begin{proof}
Let $\Delta^n_{n-1}=\{ \delta \in \Delta_{n-1} \colon \delta F_{n-1} \subseteq
F_n\}$. So, $\Delta^n_{n-1}$ is a finite subset of $\Delta_{n-1}$.
For $m>n$ define $\Delta^m_n= \Delta^m_{m-1} \Delta^{m-1}_{m-2} \cdots \Delta^{n+1}_n$.
Note that by centeredness that $\Delta^m_n \subseteq \Delta^{m+1}_n$.
Set $\tDelta_n =\bigcup_{m>n} \Delta^m_n$.
Since $1_G \in \Delta_{n-1}$ and $F_{n} \subseteq F_{n+1}$,
we have $1_G \in \Delta^{n+1}_{n}$ and thus $1_G \in \tDelta_n$
for each $n$. Since we have not changed the $F_n$, we still have $\bigcup_i F_i=G$
and so the $(\tDelta_n,F_n)$ are centered and cofinal.

To see that $(\tDelta_n, F_n)$ is a tiling, we first show that the distinct translates of
$F_n$ by $\tDelta_n$ are disjoint.  Suppose $\delta F_n \cap \eta F_n \neq \emptyset$
with $\delta,\eta \in \tDelta_n$.
Say $\delta= \delta^m_{m-1} \cdots \delta^{n+1}_n$, $\eta= \eta^m_{m-1} \cdots \eta^{n+1}_n$
where $\delta^{i+1}_i$, $\eta^{i+1}_i \in \Delta^{i+1}_i$ (we may assume a common value for $m$
by centeredness). By an immediate induction on $i$ we have that for all $i>n$
that $\tDelta^i_n F_n =F_i$. In particular,
$\delta^{m-1}_{m-2} \cdots \delta^{n+1}_n F_n \subseteq F_{m-1}$ and similarly
$\eta^{m-1}_{m-2} \cdots \eta^{n+1}_n F_n \subseteq F_{m-1}$. By definition, the
distinct $\Delta^m_{m-1}$ translates of $F_{m-1}$ are disjoint, and so $\delta^m_{m-1}=
\eta^m_{m-1}$. Continuing in this manner gives that $\delta^{i+1}_i =\eta^{i+1}_i$
for all $i$. Thus, the two translates of $F_n$ are the same. We next show that
the $\tDelta_n$ translates of $F_n$ cover $G$. This follows from
$F_m= \tDelta^m_n F_n$ and the fact that $G=\bigcup_m F_m$. So,
$G=\bigcup_m \Delta^m_n F_n= \bigcup \tDelta_n F_n$.

To see the tilings $(\tDelta_n,F_n)$ are coherent, consider $\delta F_n$ for $\delta
\in \tDelta_n$. Say $\delta =\delta^m_{m-1} \cdots \delta^{n+1}_n$ where again $\delta^{i+1}_i
\in \Delta_i$. Since $F_n=\Delta^n_{n-1} F_{n-1}$ we have
$\delta F_n= \bigcup_{\delta^n_{n-1} \in \Delta^n_{n-1}}
\delta^m_{m-1} \cdots \delta^{n+1}_n \delta^n_{n-1} F_{n-1}$. Each element
$\delta^m_{m-1} \cdots \delta^{n+1}_n \delta^n_{n-1}$ lies in $\tDelta_{n-1}$,
and this is a disjoint union as we already showed. This shows the tilings are coherent.

Finally, note that the sets $\tilde{\tDelta}_n$ defined starting from the \ccc tiling
$(\tDelta_n,F_n)$ are the same as the sets $\tDelta_n$. This is because the
corresponding sets ${\tDelta}^{i+1}_i$ and $\Delta^{i+1}_i$ are equal.
This in turn follows from the fact that we use the same sets $F_{i+1}$, $F_i$
in these definitions and that $F_{i+1}=\Delta^{i+1}_i F_i$ which shows that
${\tDelta}^{i+1}_i=\Delta^{i+1}_i$ as $\Delta^{i+1}_i \subseteq \tDelta_i$.

The second claim of the lemma is a particular case of $F_n=\tDelta^n_k F_k$
which was noted above, using $k=0$. The third claim is immediate by
centeredness since $\delta^m_{m-1} \cdots \delta^{n+2}_{n+1}=
\delta^m_{m-1} \cdots \delta^{n+2}_{n+1} 1_G\in \tDelta_n$ as $1_G \in \Delta^{n+1}_{n}$.
\end{proof}

Note that the tilings constructed in Lemma~\ref{imtiling} have the following uniformity property.
For any $i<n$, any $\gamma_1$, $\gamma_2 \in \Delta_n$, and for any $g \in F_n$,
we have $\gamma_1 g \in \Delta_i$ iff $\gamma_2 g \in \Delta_i$.  This is because
$\gamma F_n \cap \Delta_i= \gamma\Delta^n_{n-1} \Delta^{n-1}_{n-2} \cdots \Delta^{i+1}_i$
for any $\gamma \in \Delta_n$, which follows easily from the properties of Lemma~\ref{imtiling}.
This uniformity will be an important ingredient for blueprints, which will be
defined and studied in the next chapter (\cf\ Definition~\ref{def:blueprint}).

In the rest of this section we show that the class of \ccc groups is closed under taking direct products, more generally direct sums, and finite index extensions.

The next basic lemma is used in the proof of Theorem~\ref{thm:abeliantiling}.
We will prove an  extension of it also in Lemma~\ref{sumccc}.

\begin{lem} \label{prodccc}
Suppose $G$, $H$ are \ccc groups, and let $(\Delta_n,F_n)$, $(\Delta'_n,F'_n)$ be \ccc tilings
for $G$, $H$ respectively. Then $(\Delta_n\times \Delta'_n, F_n\times F'_n)$ is a \ccc tiling for $G \times H$.
\end{lem}

\begin{proof}
For fixed $n$, every element of $g$ can be written uniquely in the form $\delta f_n$ where $\delta
\in \Delta_n$ and $f \in F_n$, and likewise every element of $H$ can be written uniquely as
$\delta' f'$ for $\delta' \in \Delta'_n$, $f' \in F'_n$. So, $(g,h)$ can be written uniquely as
$(\delta, \delta')\cdot (f,f')$ where $(\delta,\delta') \in \Delta\times \Delta'$ and $(f,f') \in
F_n\times F'_n$. This shows that every $(\Delta_n\times \Delta'_n, F_n\times F'_n)$ is a tiling of $G \times H$.
Since each $\delta_n F_n$ ($\delta_n \in \Delta_n$) is contained in some $\delta_{n+1} F_{n+1}$
($\delta_{n+1} \in \Delta_{n+1}$), and likewise for $H$, it follows immediately
that every $(\delta_n, \delta'_n) \cdot (F_n\times F'_n)$ is contained in some
 $(\delta_{n+1}, \delta'_{n+1}) \cdot (F_{n+1}\times F'_{n+1})$ and so our tilings on $G \times H$ are coherent.
The cofinality of the $G\times H$ tilings follows from that of the $G$, $H$ tilings, and centeredness
follows immediately from $1_G \in \Delta_n$, $1_H \in \Delta'_n$ and so $(1_G,1_H) \in \Delta_n\times \Delta'_n$.
\end{proof}

Of course, Lemma~\ref{prodccc} holds also in the case where a group $\tilde{G}$ is the ``internal''
product of subgroups $G$ and $H$ with the appropriate change of notation.

\begin{lem} \label{sumccc}
Suppose $G=\sum_i G_i$ where each $G_i$ is a \ccc group. Then $G$ is a \ccc group.
\end{lem}

\begin{proof}
Let $(\Delta^i_n, F^i_n)$ be a \ccc tiling for $G_i$. Let $F_n=F^1_n F^2_n \cdots F^n_n$.
Let $\Delta_n= \Delta^1_n \cdots \Delta^n_n \cdot H_n$ where $H_n=\sum_{i=n+1}^\infty G_i$.
Clearly $1_G \in \Delta_n$ for each $n$. Fix an $n$ and we show $(\Delta_n,F_n)$ is a tiling of $G$.
Every $g \in G$ can be written in the form $g= g_1 \cdots g_n h$ where $h \in H_n$.
Then

\begin{equation*}
\begin{split}
g& =(\delta^1_n f^1_n) \cdots (\delta^n_n f^n_n) h\\ & =
(\delta^1_n \cdots \delta^n_n)(f^1_n \cdots f^n_n)h,
\end{split}
\end{equation*}

\noindent
where $\delta^i_n \in \Delta^i_n$
and $f^i_n \in F^i_n$. Since $H_n$ commutes with $G_1,\dots,G_n$, we have
$g=(\delta^1_n \cdots \delta^n_n h) (f^1_n \cdots f^n_n) \in \Delta_n F_n$.
So, the $\Delta_n$ translates of $F_n$ cover $G$.

Suppose next that $\delta F_n \cap \eta F_n \neq \emptyset$, where $\delta, \eta \in \Delta_n$.
So, letting $\delta=\delta^1_n \cdots \delta^n_n h$, $\eta=\eta^1_n \cdots \eta^n_nh'$
where $\delta^i_n, \eta^i_n \in \Delta^i_n$, $h,h' \in H_n$,
we have

$$
\delta^1_n \cdots \delta^n_n h f^1_n \cdots f^n_n= \eta^1_n \cdots \eta^n_n h' k^1_n \cdots k^n_n,$$

\noindent
where $f^i_n,k^i_n \in F^i_n$. Rearranging this gives

$$
(\delta^1_n f^1_n) \cdots (\delta^n_n f^n_n) h=
(\eta^1_n k^1_n) \cdots (\eta^n_n k^n_n) h'.$$

\noindent

Since $G$ is a direct sum of $\sum_{i=1}^n G_i$ and $H_n$ we must have
$h=h'$, and then $\delta^i_n f^i_n=\eta^i_n k^i_n$ for all $i=1,\dots,n$.
Since $(\Delta^i_n,F^i_n)$ is a tiling of $G_i$ this gives
$\delta^i_n=\eta^i_n$ and $f^i_n=k^i_n$. In particular,
$\delta=\eta$. So, $(\Delta_n,F_n)$ is a tiling of $G$.

Finally, consider a translate $T=\delta F_{n+1}$ of $F_{n+1}$ where $\delta \in \Delta_{n+1}$.
Then $T=\delta^1_{n+1} \cdots \delta^n_{n+1} \delta^{n+1}_{n+1} h F_{n+1}$, where $h \in H_{n+1}$.
We have:

\begin{equation*}
\begin{split}
T&=\delta^1_{n+1} \cdots \delta^n_{n+1} \delta^{n+1}_{n+1} h F^1_{n+1} \cdots F^n_{n+1}F^{n+1}_{n+1}
\\ &
= (\delta^1_{n+1} F^1_{n+1}) \cdots (\delta^n_{n+1} F^n_{n+1})
(\delta^{n+1}_{n+1} F^{n+1}_{n+1}) h
\\ &
=\bigcup_{\delta_1 \in A_1} \cdots \bigcup_{\delta_n \in A_n}
\bigcup_{a_{n+1} \in \delta^{n+1}_{n+1} F^{n+1}_{n+1} }
(\delta_1F^1_n) \cdots (\delta_nF^n_n) a_{n+1} h
\\ &
= \bigcup_{\delta_1 \in A_1} \cdots \bigcup_{\delta_n \in A_n}
\bigcup_{a_{n+1} \in \delta^{n+1}_{n+1} F^{n+1}_{n+1} }
(\delta^1_n \cdots \delta^n_n a_{n+1} h) F^1_n \cdots F^n_n
\end{split}
\end{equation*}

\noindent
For some finite sets $A_i \subseteq \Delta^i_n$. This shows that
$T$ is a (disjoint) union of translates of $F_n$ by elements of
$\Delta_n$.
\end{proof}

\begin{lem}\label{lem:finiteindexccc}
Suppose $G$ has a finite index subgroup $H$ which is a \ccc group. Then $G$ is a \ccc group.
\end{lem}

\begin{proof}
Let $(\Delta_n,F_n)$ be a \ccc tiling for $H$ and let $Hx_1, Hx_2,\dots, Hx_m$ be the distinct
right cosets of $H$ in $G$, where $x_1=e$. Let $F'_n= \bigcup_{i=1}^m F_n x_i$. We claim that
$(\Delta_n, F'_n)$ is a \ccc tiling of $G$. For each $n$, every $g \in G$ is of the form
$g= h x_i$ for some $h \in H$, and thus of the form $g= \delta_n f_n x_i$ where
$\delta_n \in \Delta_n$ and $f_n \in F_n$. So, $g \in \delta_n F'_n$. So, the $\Delta_n$
translates of the $F'_n$ cover $G$. Suppose next that $\delta^1_n F'_n \cap \delta^2_n F'_n \neq
\emptyset$, where $\delta^1_n, \delta^2_n \in \Delta_n$. Then,
$\delta^1_n f^1_n x_i= \delta^2_n f^2_n x_j$ for some $f^1_n, f^2_n \in F_n$.
By disjointness of the distinct
cosets we have $i=j$, and therefore $\delta^1_n f^1_n=\delta^2_n f^2_n$. Since $(\Delta_n,F_n)$
is a tiling, it follows that $\delta^1_n=\delta^2_n$ (and $f^1_n=f^2_n$). Thus,
the distinct translates of $F'_n$ by $\Delta_n$ are disjoint. So, each $(\Delta_n,F_n)$
is a tiling of $G$.
By assumption, $1_G=1_H  \in \Delta_n$
for each $n$. It is also easy to see that they are cofinal.

Finally, to show coherence consider a translate $T=\delta_{n+1}F'_{n+1}$, where $\delta_{n+1}
\in \Delta_{n+1}$. So, $T= \bigcup_{i=1}^m \delta_{n+1} F_{n+1} x_i$. Let $A \subseteq
\Delta_n$ be finite such that $\delta_{n+1}F_{n+1}= \bigcup_{\delta \in A} \delta F_n$.
So, $T= \bigcup_{\delta \in A} \bigcup_{i=1}^m \delta F_n x_i=\bigcup_{\delta \in A}
\delta F'_n$. So, each $\Delta_{n+1}$ translate of $F'_{n+1}$ is a (disjoint by above) union of $\Delta_n$
translates of $F'_n$.
\end{proof}

\section{Abelian, nilpotent, and polycyclic groups are \ccc} \label{SECT ABNILPOL}

In this section we show that all abelian, nilpotent, and polycyclic groups are ccc.

\begin{theorem} \label{thm:abeliantiling}
Every countable abelian group is a \ccc group.
\end{theorem}

\begin{proof}
Let $G=\{ 1_G=g_0,g_1,\dots\}$ be abelian.
Suppose that for $i\leq m$ we have constructed tilings $(\Delta^i_n,F^i_n)$
satisfying the following:
\begin{enumerate}
\item \label{tiling1}
Each $(\Delta^i_n,F^i_n)$ is a \ccc tiling of $G_i=\langle g_0,\dots,g_i \rangle$.
\item \label{tiling3}
If $i<m$ each marker region $\delta^i_n F^i_n$ (for $\delta^i_n \in \Delta^i_n$) is contained in a (unique)
region $\delta_n^{i+1} F^{i+1}_n$ (for some $\delta^{i+1}_n\in \Delta^{i+1}_n$).
\end{enumerate}

We proceed to construct the $\Delta_{n}^{m+1}$, $F^{m+1}_{n}$.
Suppose first that $g_{m+1}$ has infinite order in $G/G_{m}$.
In this case, $G_{m+1}\cong G_m \oplus \langle g_{m+1} \rangle$. Since $\langle
g_{m+1} \rangle \cong \mathbb{Z}$, we may easily get a \ccc sequence of tilings
$(\Delta'_n, F'_n)$ for $\langle g_{m+1}\rangle$.
Let $F_{n}^{m+1}=F_n^m \cdot F'_n$, and let $\Delta_{n}^{m+1}=\Delta_n^m \cdot \Delta'_{n}$.
From Lemma~\ref{prodccc} we have that
$(\Delta^{m+1}_{n}, F^{m+1}_{n})$ is a \ccc sequence of tilings of $G_{m+1}$.
Property~(\ref{tiling3}) is clear from the proof of Lemma~\ref{prodccc}
(since $\delta^m_{n} F^m_n =(\delta^m_n \cdot 1_G ) F^m_n \subseteq (\delta^m_n \cdot 1_G) F^{m+1}_{n}
=\delta^{m+1}_{n} F^{m+1}_{n  }$ as
$\delta^{m+1}_n=(\delta^m_n \cdot 1_G) \in \Delta^{m+1}_n=\Delta^m_n \Delta'_n$
since $1_G \in \Delta'_n$ by centeredness, and $F^{m+1}_{n}=F^m_n F'_n$ contains $F^m_n$
as $1_G \in F'_n$ by centeredness).

If $g_{m+1}$ has finite order $k$ in $G/G_m$, Let
$\Delta_{n}^{m+1}=\Delta_n^m$ and $F_{n}^{m+1}=F_n^m\cdot \{ g_{m+1}^0,\dots, g_{m+1}^{k-1} \}$.
This again defines the tiling $(\Delta^{m+1}_{n}, F^{m+1}_{n})$.
Clearly $\bigcup_n F_{n}^{m+1}=G_{m+1}$.
The coherence property that for all $\delta^{m+1}_{n} \in \Delta^{m+1}_{n}$
there is a $\delta^{m+1}_{n+1} \in \Delta^{m+1}_{n+1}$ with
$\delta^{m+1}_{n} F^{m+1}_{n} \subseteq \delta^{m+1}_{n+1}F^{m+1}_{n+1}$
follows immediately from the coherence property for the tilings $(\Delta^m_n, F^m_n)$.
Namely, $\delta^{m+1}_{n}F^{m+1}_{n}=(\delta^{m+1}_{n} F^m_n) \{g_{m+1}^0,\dots, g_{m+1}^{k-1} \}
\subseteq (\delta^{m+1}_{n+1} F^{m}_{n+1})  \{g_{m+1}^0,\dots, g_{m+1}^{k-1} \}
=\delta^{m+1}_{n+1} F^{m+1}_{n+1}$, for some $\delta^{m+1}_{n+1} \in \Delta^{m+1}_{n+1}$.
Next we show that the $\Delta^{m+1}_{n}$ translates of $F^{m+1}_{n}$ cover $G$.
Since $g_{m+1}^k \in G_m$, every element $x$ of $G_{m+1}$
is of the form $x=g \cdot  g_{m+1}^i$ for some $i<k$, where $g \in G_m$.
So, $x= \gamma  f  g_{m+1}^i$ where $\gamma \in \Delta_n^m$ and $f \in F_n^m$.
Thus, $\gamma \in \Delta_{n}^{m+1}$ and $f g_{n+1}^i \in F_{n}^{m+1}$. To see
disjointness of the translates, suppose $\gamma_1 ( f_1 g_{m+1}^i)= \gamma_2 ( f_2 g_{m+1}^j)$,
where $\gamma_1,\gamma_2 \in \Delta_{n}^{m+1}=\Delta_n^m$ and $f_1,f_2 \in
F_n^m$ (without loss of generality, $i \leq j<k$).
So, $g_{m+1}^{j-i} \in G_m$. As $g_{m+1}$ has order $k$ in $G/G_m$,
$i=j$. So, $\gamma_1 f_1=\gamma_2 f_2$, and thus $\gamma_1=\gamma_2$ and $f_1=f_2$.
Property~(\ref{tiling3}) is again clear (since $g_{m+1}^0$ is the identity).

We have now defined $(\Delta^m_n,F^m_n)$ for all $n$, $m$ which satisfy the above properties.
We may assume that each $g_m \notin \langle g_0,\dots,g_{m-1}\rangle$.
Let $k_m$ be the order of $g_m$ in $G/G_{m-1}$ if this order is finite and otherwise
$k_m=\infty$.

\begin{claim}
For each $m$, let $B_m=\{ g_{m+1}^{i_{m+1}} g_{m+2}^{i_{m+2}} \cdots g_j^{i_j}
\colon i_{m+1}<k_{m+1},\dots, i_j<k_j\}$. Then $B_m$
is a set of coset representatives for $G/G_m$. Also, $B_m \supseteq B_{m+1}$.
\end{claim}

\begin{proof}
This is easily checked.
\end{proof}

Notice from the above construction that property~(\ref{tiling3}) was actually established in the following
stronger sense: if $x,y \in G_m$ and $x$, $y$ are in the same marker region
$\delta^m_n F^m_n$, then for any
$i<k_{m+1}$ we have that $xg_{m+1}^i $, $yg_{m+1}^i$ are in the same region $\delta^{m+1}_{n} F^{m+1}_{n}$.

Let $F_n=F_n^n$ and $\Delta_n= B_n\Delta_n^n$. This defines the sequence $(\Delta_n, F_n)$ for $G$.
If $x \in G$, then for each $n$, there is a $\gamma \in B_n$ and a $y \in G_n$
such that $x=\gamma y$. Then is a $\delta \in \Delta_n^n$ and $f \in F_n^n=F_n$
such that $y=\delta f$. So, $x=\gamma \delta f$ and $\gamma \delta \in \Delta_n$.
So, the $\Delta_n$ translates of $F_n$ cover $G$. To see disjointness
suppose $\gamma_1 f_1=\gamma_2 f_2$ where $\gamma_1,\gamma_2 \in \Delta_n$ and
$f_1,f_2 \in F_n$. Say $\gamma_1=b_1  \delta_1$ where $b_1 \in B_n$ and $\delta_1
\in \Delta_n^n$. Likewise, write $\gamma_2=b_2 \delta_2$.
So, $b_1 \delta_1 f_1=b_2 \delta_2 f_2$. Since $\delta_1 f_1$ and $\delta_2 f_2$ are in $G_n$,
this implies $b_1=b_2$, and thus $\delta_1=\delta_2$ and $f_1=f_2$. So, each $(\Delta_n,F_n)$
is a tiling of $G$. Suppose finally that $x,y \in G$ and $x$, $y$
are in the same $n$-level marker region for $G$.
Say $x=b  \delta f_1$, $y=b \delta f_2$ where $b \in B_n$, $\delta \in \Delta_n^n$,
and $f_1,f_2 \in F_n^n$.  So, $\delta f_1$ and $\delta f_2$ are in the same $n$-level
region for $G_n$.
Say $b=g_{n+1}^{i_{n+1}} b'$ where $b' \in B_{n+1}$. By our observation in the first paragraph
after the claim,
$\delta f_1 g_{n+1}^{i_{n+1}}$ and $\delta f_2 g_{n+1}^{i_{n+1}}$
are in the same $n$-level region for $G_{n+1}$ and by coherence are in the same
$n+1$-level region for $G_{n+1}$. So, $x=(\delta b')(f_1 g_{n+1}^{i_{n+1}})$ is in the same
$n+1$-level region for $G$ as $y=(\delta b')(f_2 g_{n+1}^{i_{n+1}})$.
\end{proof}

We now extend Theorem~\ref{thm:abeliantiling} to nilpotent groups.
We use the following lemma, in which $\Z(G)$ denotes the center of $G$.\index{$\Z(G)$}

\begin{lem} \label{centerccc}
Let $H \leq \Z(G)$ and suppose that both $H$ and $G/H$ are \ccc groups. Then
$G$ is a \ccc group.
\end{lem}

\begin{proof}
Let  $(\Delta'_n, F'_n)$ be a \ccc tiling for $H$ and fix also
a \ccc tiling $(\bDelta_n, \bF_n)$ for $G/H$.
For notational convenience we assume (without loss of generality) that
$\bF_0$ and $F'_0$ both consist of just the identity element of $G/H$ and $H$
respectively. Let the $\bDelta^m_n$, ${\Delta'}^m_n$ be defined as in Lemma~\ref{imtiling}.
So,  $\bDelta^n_{n-1}=\{ \bdelta \in \bDelta_{n-1} \colon \bdelta \bF_{n-1} \subseteq \bF_n\}$
and $\bDelta^m_n= \bDelta^m_{m-1} \bDelta^{m-1}_{m-2} \cdots \bDelta^{n+1}_n$, and
likewise for the ${\Delta'}^m_n$. From Lemma~\ref{imtiling} we may assume that
$\bDelta_n= \bigcup_m \bDelta^m_n$ and $\Delta'_n=\bigcup_m {\Delta'}^m_n$.

Fix coset representatives $\Delta^n_{n-1}$ for the elements $\bdelta^n_{n-1}$
in $\bDelta^n_{n-1}$. We may assume that $1_G \in \Delta^n_{n-1}$ for all $n$.
Let $\Delta^m_n= \Delta^m_{m-1} \Delta^{m-1}_{m-2} \cdots \Delta^{n+1}_n$,
and $\Delta_n=\bigcup_{m>n} \Delta^m_n$.
Thus, $\Delta_n$ is a set of coset representatives for $\bDelta_n$.
Note that $\bF_n=\bDelta^n_{n-1} \cdots \bDelta^2_1 \bDelta^1_0$, and so
if we let $F_n =\Delta^n_{n-1} \cdots \Delta^2_1 \Delta^1_0$, then $F_n$
is a set of coset representatives for $\bF_n$.
For the group $H$ we have also defined the finite sets ${\Delta'}^n_m$.
From Lemma~\ref{imtiling} we also have that $F'_n={\Delta'}^n_{n-1}\cdots {\Delta'}^1_0$.

We claim that $(\Delta_n \Delta'_n, F_n F'_n)$ is a \ccc tiling of $G$.
We first show that for each $n$ that $(\Delta_n \Delta'_n, F_n F'_n)$ is a tiling of $G$.
We must show that every element $g$ of $G$ can be written uniquely in the form
$g= d d' f f'$ where $d \in \Delta_n$, $d' \in \Delta'_n$, $f \in F_n$, and $f' \in F'_n$.
To see uniqueness, suppose $d_1 d'_1f_1 f'_1=d_2 d'_2 f_2 f'_2$. Since the $d'$ terms are in
$H\leq \Z(G)$, this can be rewritten as $d_1f_1 d'_1 f'_1= d_2 f_2 d'_2 f'_2$.
In $G/H$ this becomes $d_1 f_1=d_2 f_2$. This implies that in $G/H$ that $d_1=d_2$ and $f_1=f_2$
since $(\bDelta_n, \bF_n)$ is a tiling of $G/H$. Since the distinct points of $\Delta_n$ and $F_n$
are in distinct cosets by $H$, it follows that $d_1=d_2$ and $f_1=f_2$.
We therefore have that $d'_1 f'_1=d'_2 f'_2$. Since $(\bDelta'_n, \bF'_n)$ is a tiling of $H$ we have
$d'_1=d'_2$ and $f'_1=f'_2$. To show existence, let $g \in G$.
Let $d \in \Delta_n$, $f \in F_n$ be such that $g= df$ in $G/H$, that is, $g=dfh$ where $h \in H$.
Since $(\Delta'_n,F'_n)$ is a tiling of $H$, we may write $h=d'f'$ where $d' \in \Delta'_n$ and
$f' \in F'_n$. So, $g=dfd'f'=dd' ff' $ where $dd' \in \Delta_n \Delta'_n$ and
$f f' \in F_n F'_n$.

The tiling are clearly centered as $1_G \in \Delta_n$, $1_G \in \Delta'_n$
and so $1_G \in \Delta_n \Delta'_n$.
Also, the tilings are easily cofinal since each of the tilings on $G/H$ and on $H$ are.
To show coherence, consider an $n$ level tile, which is of the form
$T=d d' F_{n} F'_{n}$, where $d \in \Delta_{n}$, $d' \in \Delta'_{n}$.
We have $F_{n} =\bigcup_{\delta \in \Delta^n_{n-1} } \delta F_{n-1}$ and
$F'_{n}=\bigcup_{\delta' \in {\Delta'}^n_{n-1}} \delta'F'_{n-1}$. So,

\begin{equation*}
\begin{split}
T& = \bigcup_{ \delta \in \Delta^n_{n-1}} \bigcup_{\delta' \in {\Delta'}^n_{n-1}}
d d' \delta F_{n-1} \delta' F'_{n-1}\\ & =
\bigcup_{ \delta \in \Delta^n_{n-1}} \bigcup_{\delta' \in {\Delta'}^n_{n-1}}
d \delta d' \delta' F_{n-1} F'_{n-1}.
\end{split}
\end{equation*}
This shows that $T$ is a union of translates of $F_{n-1}F'_{n-1}$ by elements of
$\Delta_{n-1} {\Delta'}_{n-1}$ (note that $\Delta_n \Delta^n_{n-1} \subseteq \Delta_{n-1}$)
and we are done.
\end{proof}

\begin{theorem} \label{ch05_thm:nilpotenttiling}
Every countable nilpotent group is a \ccc group.
\end{theorem}

\begin{proof}
This follows immediately by an induction on the nilpotency rank of $G$ using
Lemma~\ref{centerccc}.
\end{proof}

A closer scrutiny of the proof of Lemma~\ref{centerccc} gives us the following more general fact. It will be useful when the subgroup $H$ is no longer abelian.\index{$\Z_G(H)$}

\begin{lem} Let $G$ be a countable group, $H\leq G$, and $\Z_G(H)$ be the centralizer of $H$ in $G$, i.e., $\Z_G(H)=\{g\in G\,:\, gh=hg \mbox{ for all } h\in H\}$.
If $G=\Z_G(H)H$ and both $H$ and $G/H$ are \ccc groups, then $G$ is a \ccc group.
\end{lem}

\begin{proof}
Note that the assumption $G=\Z_G(H)H$ implies that $H\unlhd G$. The proof of Lemma~\ref{centerccc} can be repeated verbatim, except when picking the coset representatives for $\Delta^n_{n-1}$, we require them to be chosen from the set $\Z_G(H)$; this is possible since $G=\Z_G(H)H$.
\end{proof}

We next extend Theorem~\ref{thm:abeliantiling} to a class of solvable groups.
Unfortunately we are unable to extend the result to all solvable groups, but we must restrict
to those with finitely generated quotients in the derived series. We recall the following definition.

\begin{definition}
A countable group $G$ is said to be {\em polycyclic} if $G$ has a finite derived series
$G=G^0 \unrhd G^1 \unrhd \cdots \unrhd G^k=\{ 1_G\}$, where $G^{i+1}=(G^i)' = [G_i, G_i]$, and each quotient group
$G^{i-1}/G^{i}$ is finitely generated.
\end{definition}\index{polycyclic}

Being polycyclic is equivalent to having a subnormal series
$G=G_0 \unrhd G_1 \cdots \unrhd G_m=\{ 1_G\}$ (\ie, each $G_{i+1}$ is a normal subgroup
of $G_i$) with all quotient groups $G_{i+1}/G_i$ cyclic.

The collection of polycyclic groups  include some finitely generated groups
which are not of polynomial growth, and therefore not virtually nilpotent (by Wolf \cite{Wo}). Thus Lemma~\ref{lem:finiteindexccc} and Theorem~\ref{ch05_thm:nilpotenttiling}
do not prove that polycyclic groups are ccc. On the other hand, all polycyclic groups are residually finite, and we will prove in the next section that all residually finite groups are ccc.
Here we present a direct proof of the \ccc property for a class of solvable groups containing the class of polycyclic groups. We hope that this proof may be useful in future generalizations of this result to solvable groups.

\begin{theorem} \label{ch05_thm:solvable}
If $G$ is a countable solvable group and $[G, G]$ is polycyclic, then $G$ is a \ccc group. In particular, if $G$ is polycyclic then $G$ is a \ccc group.
\end{theorem}

\begin{proof}
Every polycyclic group is countable, and subgroups of polycyclic groups are polycyclic. So the second sentence in the statement of the theorem follows from the first. So we prove the first sentence. We first show the following simple fact about finitely generated abelian groups
which we will need for the proof.

\begin{lem} \label{ch05_lem:ablem}
Every finitely generated abelian group $A$ has a \ccc tiling
$(\Delta_n, F_n)$ satisfying:

\begin{enumerate}
\item [\rm (1)]
Each $\Delta_n$ is a subgroup of $A$.
\item [\rm (2)]
Each $\Delta_n$ is invariant under any automorphism of $A$.
\item [\rm (3)]
Each $(\Delta_n,F_n)$ satisfies the property of Lemma~\ref{imtiling},
that is, $$\Delta_n=\bigcup_{m>n} \Delta^m_{m-1} \Delta^{m-1}_{m-2} \cdots \Delta^{n+1}_n,$$
where the $\Delta^i_{i-1}$ are as in Lemma~\ref{imtiling}.
\end{enumerate}
\end{lem}

\begin{proof}
Write $A= \mathbb{Z}^k \oplus F$, where $F$ is a finite subgroup. Let $m>1$ be such that
$F^{m} =\{ e\}$. Let $\Delta_n= A^{m^n}$, for all $n \geq 1$. Clearly,
(1) and (2) are satisfied. Let
$F_n=\{ (z_1^{i_1},z_2^{i_2},\dots,z_k^{i_k},f) \colon 0 \leq i_1,\dots i_k <m^n,
f \in F\}$, where $(z_1,\dots, z_k)$ generates $\mathbb{Z}^k$.
Easily, $(\Delta_n,F_n)$ is a \ccc tiling of $A$. Property (3)
is easily checked.
\end{proof}

Consider the derived series of $G$, $G=G^0 \unrhd G^1 \unrhd \cdots \unrhd G^k=\{ 1_G\}$, where $G^{i}=(G^{i-1})' = [G_{i-1}, G_{i-1}]$. By assumption, for each $i > 1$ the quotient group
$G^{i-1}/G^{i}$ is finitely generated.
For each quotient group $G^{i-1}/G^i$, fix a \ccc tiling
$(\bDelta_n^i,\bF_n^i)$ satisfying clause (3) (this can always be done by Lemma \ref{imtiling} and Theorem \ref{thm:abeliantiling}). Additionally, for $i > 1$ require this \ccc tiling to be as in Lemma~\ref{ch05_lem:ablem}.
As in Lemma~\ref{imtiling}, we let
$(\bDelta^i)^n_{n-1} \subseteq \bDelta_{n-1}^i$ be finite such that
$\bF_n^i$ is the disjoint union of translates of $\bF_{n-1}^i$ by
the elements of $(\bDelta^i)^n_{n-1}$. We choose coset representatives
$(\Delta^i)^n_{n-1} \subseteq G^{i-1}$ for $(\bDelta^i)^n_{n-1}$ in $G^{i-1}/G^i$.
Let $\Delta_n^i= \bigcup_{m>n} (\Delta^i)^m_{m-1} \cdots (\Delta^i)^{n+1}_n$.
Let $F_n^i= (\Delta^i)^n_{n-1} \cdots (\Delta^i)^1_0$. Then, the $G^i$ cosets
of the $(\Delta_n^i, F_n^i)$ also give a representation of the given tilings.
Note that the $G^i$ cosets of the elements of $F_n^i$ are exactly the elements of the given sets $\bF_n^i$. Similarly, from property (3) we have that
the $G^i$ cosets of the elements in $\Delta^i_n$ are exactly the elements
of $\bDelta_n^i$.

We summarize some of the properties of these sets.

\begin{enumerate}
\item [\rm (4)]
$\Delta_{n+1}^i \subseteq \Delta_{n}^i$ for all $n$.
\item [\rm (5)]
$\Delta_n^i= \Delta_{n+1}^i (\Delta^i)^{n+1}_n$.
\item [\rm (6)]
$F_n^i \subseteq F_{n+1}^i$.
\item [\rm (7)]
For $i > 1$, $\bDelta_n^i$ is a subgroups of $G^{i-1}/G^i$ and is invariant under every automorphism
of $G^{i-1}/G^i$.
\end{enumerate}

Let now

\begin{equation*}
\begin{split}
F_n& = [(\Delta^1)^n_{n-1} \cdots ({\Delta^k})^n_{n-1}]
[(\Delta^1)^{n-1}_{n-2} \cdots (\Delta^k)^{n-1}_{n-2}]\cdots
[(\Delta^1)^1_0 \cdots (\Delta^k)^{1}_{0}]
\\
\Delta_n &=\bigcup_{m>n} [(\Delta^1)^m_{m-1} \cdots (\Delta^k)^m_{m-1})]\cdots
[ (\Delta^1)^{n+1}_n \cdots (\Delta^k)^{n+1}_n]
\end{split}
\end{equation*}

We show that $(\Delta_n, F_n)$ is a \ccc tiling for $G$.

By definition, $F_{n+1}$ is a union of translates of $F_n$, namely by points of the set
$(\Delta^1)^n_{n-1} \cdots ({\Delta^k})^n_{n-1}$.
We show that these translates are disjoint. More generally, suppose $m>n$ and

\begin{equation} \label{ch05_eqn:uniq}
\begin{split}
& (\delta_m^1 \cdots \delta_m^k)(\delta_{m-1}^1 \cdots\delta_{m-1}^k) \cdots
(\delta_1^1 \cdots \delta_1^k)=
 (\rho_m^1 \cdots \rho_m^k)(\rho_{m-1}^1 \cdots\rho_{m-1}^k) \cdots (\rho_1^1 \cdots \rho_1^k)
\end{split}
\end{equation}

\noindent
where $\delta_j^i$, $\rho_j^i \in (\Delta^i)^j_{j-1}$.
We show that $\delta_j^i=\rho_j^i$, that is, all the corresponding terms in the two expressions
above are equal.

\begin{lem} \label{ch05_lem:iden}
For any $\delta_m^1,\dots, \delta_m^k, \dots, \delta^1_1,\dots \delta_1^k$ in $G$ we have:

\begin{equation*}
\begin{split}
& (\delta_m^1 \delta_m^2\cdots \delta_m^k)(\delta_{m-1}^1 \delta_{m-1}^2
\cdots\delta_{m-1}^k) \cdots (\delta_1^1 \delta_1^2 \cdots \delta_1^k)
\\ & =
(\delta_m^1 \delta_{m-1}^1 \cdots \delta_1^1) \\ & \cdot
({\delta_m^2}^{\delta_{m-1}^1 \cdots \delta_1^1} {\delta_{m-1}^2}^{\delta_{m-2}^1 \cdots \delta_1^1}
\cdots {\delta_1^2} )\\ & \cdot
({\delta_m^3}^{\delta_{m-1}^1 \cdots \delta_1^1 {\delta_{m-1}^2}^{\delta_{m-2}^1 \cdots \delta_1^1}
{\delta_{m-2}^2}^{\delta_{m-3}^1 \cdots \delta_1^1} \cdots \delta_1^2} \cdots)\\ &  \cdots
\\&
=({\delta_m^1}^{c_m^1} \cdots {\delta_1^1}^{c_1^1})({\delta_m^2}^{c_m^2} \cdots {\delta_1^2}^{c_1^2})
\cdots ({\delta_m^k}^{c_m^k} \cdots {\delta_1^k}^{c_1^k})
\end{split}
\end{equation*}

\noindent
where $c_m^1=1_G$ for all $m$, and (inductively on $i$)
$$c_j^i=\delta_{j-1}^1 \cdots \delta^1_1 {\delta_{j-1}^2}^{c_{j-1}^2}
\cdots {\delta_{1}^2}^{c_{1}^2} \cdots {\delta_{j-1}^{i-1}}^{c_{j-1}^{i-1}} \cdots
{\delta_1^{i-1}}^{c_1^{i-1}}.$$
\end{lem}

\begin{proof}
Repeatedly use the identity $xy= y x^y$. First move all the terms
$\delta_m^1, \delta_{m-1}^1, \dots, \delta^1_1$ to the left of the equation. Then move all of the terms
$\delta_m^2$, $\delta_{m-1}^2,\dots$, $\delta_1^2$ to the left to immediately
follow these terms (the terms $\delta_j^2$ have actually become ${\delta_j^2}^{\delta_{j-1}^1
\cdots \delta^1_1}={\delta_j^2}^{c_j^2}$ from the first step).
Continuing in this manner gives the above equation.
\end{proof}

Apply now Lemma~\ref{ch05_lem:iden} to both sides of equation~\ref{ch05_eqn:uniq}.
This gives:

\begin{equation} \label{ch05_eqn:dr}
\begin{split}
& (\delta_m^1 \cdots \delta^1_1)({\delta_m^2}^{c^m_2} \cdots {\delta_1^2}^{c_1^2})\cdots
({\delta_m^k}^{c_m^k} \cdots {\delta_1^k}^{c_1^k})
\\ &
=
(\rho_m^1 \cdots \rho^1_1)({\rho_m^2}^{e_m^2} \cdots {\rho_1^2}^{e_1^2})\cdots
({\rho_m^k}^{e_m^k} \cdots {\rho_1^k}^{e_1^k}).
\end{split}
\end{equation}

\noindent
where the $e_j^i$ are defined as the $c_j^i$ using the $\rho$'s instead of the $\delta$'s.

Considering this equation in $G_0/G_1$ gives
$\bdelta_m^1 \cdots \bdelta^1_1= \brho_m^1 \cdots \brho^1_1$.
Since the sets $\bDelta$ satisfy (3)
we have that $\bdelta_m^1 \cdots \bdelta_2^1 \in \bDelta_1^1$, and also
$\brho_m^1 \cdots \brho_2^1 \in \bDelta_1^1$. Since $\bdelta^1_1, \brho^1_1
\in \bF^1_1$, it follows that $\bdelta^1_1=\brho^1_1$ and
$\bdelta_m^1 \cdots \bdelta^2_1= \brho_m^1\cdots \brho_2^1$. Continuing,
we get that $\bdelta_j^1=\brho_j^1$ for all $j=1,\dots,m$.
It then follows that $\delta_j^1=\rho_j^1$ for all $j$ as well.

From this and the above equation we then have that

$$
({\delta_m^2}^{c_m^2} \cdots {\delta_1^2}^{c_1^2})\cdots
({\delta_m^k}^{c_m^k} \cdots {\delta_1^k}^{c_1^k})=
({\rho_m^2}^{c_m^2} \cdots {\rho_1^2}^{c_1^2})\cdots
({\rho_m^k}^{c_m^k} \cdots {\rho_1^k}^{c_1^k}).
$$

More generally, suppose that after $i-1$ steps we have show that

\begin{equation*}
\begin{split}
& \delta_m^1=\rho_m^1, \dots, \delta^1_1=\rho^1_1
\\ &
\qquad \qquad \vdots
\\ &
\delta_m^{i-1}=\rho_m^{i-1},\dots \delta_1^{i-1}=\rho_1^{i-1}.
\end{split}
\end{equation*}

In particular, from the equation for $c_m^i$ we see that
$c_j^\ell=e_j^\ell$ for all $1 \leq j \leq m$ and all $\ell \leq i$.
From equation~(\ref{ch05_eqn:dr}) we therefore have:

$$
({\delta_m^i}^{c_m^i}\cdots {\delta_1^i}^{c_1^i})\cdots
({\delta_m^k}^{c_m^k} \cdots {\delta_1^k}^{c_1^k})=
({\rho_m^i}^{c_m^i}\cdots {\rho_1^i}^{c_1^i})\cdots
({\rho_m^k}^{e_m^k} \cdots {\rho_1^k}^{e_1^k}).
$$

Considering this equation in $G^{i-1}/G^i$ gives:

$$
{{\bdelta}_m^i}^{c_m^i}\cdots {\bdelta_1^i}^{c_1^i}=
{\brho_m^i}^{c_m^i}\cdots {\brho_1^i}^{c_1^i}.
$$

Conjugating both sides by $(c_1^i)^{-1}$ gives:

$$
{\bdelta_m^i}^{c_m^i (c_1^i)^{-1}} \cdots {\bdelta_2^i}^{c_2^i (c_1^i)^{-1}}{\bdelta_1^i}=
{\brho_m^i}^{c_m^i (c_1^i)^{-1}} \cdots {\brho_2^i}^{c_2^i (c_1^i)^{-1}}{\brho_1^i}.
$$

\noindent
From properties (4) and (7) of the $\bDelta$ sets
we have that ${\bdelta_m^i}^{c_m^i (c_1^i)^{-1}} \cdots {\bdelta_2^i}^{c_2^i  (c_1^i)^{-1}}$
and ${\brho_m^i}^{c_m^i (c_1^i)^{-1}} \cdots {\brho_2^i}^{c_2^i (c_1^i)^{-1} }$ are in $\bDelta_1^i$.
Since $\bdelta_1^i$ and $\brho_1^i$ are in $\bF_1^i$, and the distinct $\bDelta_1^i$
translates of $\bF_1^i$ are disjoint, it follows that
$\bdelta_1^i=\brho_1^i$ and ${\bdelta_m^i}^{c_m^i (c_1^i)^{-1}} \cdots
{\bdelta_2^i}^{c_2^i (c_1^i)^{-1}}=
{\brho_m^i}^{c_m^i (c_1^i)^{-1}} \cdots {\brho_2^i}^{c_2^i (c_1^i)^{-1}}$.
Continuing in this manner yields
$\bdelta_j^i=\brho_j^i$ for all $j=1,\dots,m$. It follows that $\delta_j^i=\rho_j^i$
for all $j=1,\dots,m$.

Starting from equation~\ref{ch05_eqn:uniq} we have shown that all the corresponding terms on
the two sides of the equation are equal, that is $\delta_j^i=\rho_j^i$ for all
$i=1,\dots,k$ and $j=1,\dots,m$. From this we have the following.

\begin{lem}
For each $n$, the $\Delta_n$ translates of $F_n$ are disjoint.
Also, each translate of $F_n$ by an element of $\Delta_n$ is a disjoint union
of translates of $F_{n-1}$ by elements of $\Delta_{n-1}$.
\end{lem}

\begin{proof}
Suppose $x \in \delta F_n \cap \rho F_n \neq \emptyset$, where $\delta, \rho \in \Delta_n$.
By definition of $\Delta_n$ we have that

\begin{equation*}
\begin{split}
\delta& =(\delta_m^1\cdots \delta_m^k) \cdots (\delta_{n+1}^{1}\cdots \delta_{n+1}^k)
\\
\rho & = (\rho_m^1\cdots \rho_m^k) \cdots (\rho_{n+1}^{1}\cdots \rho_{n+1}^k)
\end{split}
\end{equation*}

\noindent
where $\delta_j^i, \rho_j^i \in (\Delta^i)^{j}_{j-1}$ and
we may assume a common value of $m>n$ for both equations since $1_G \in (\Delta^i)^j_{j-1}$
for all $i$ and $j$. Then from the definition of $F_n$ we have that $x$ can be written as:

\begin{equation*}
\begin{split}
x& =(\delta_m^1\cdots \delta_m^k) \cdots (\delta_{n+1}^{1}\cdots \delta_{n+1}^k)
(\delta_n^1 \cdots \delta_n^k)(\delta^1_1 \cdots \delta_1^k)
\\ &
= (\rho_m^1\cdots \rho_m^k) \cdots (\rho_{n+1}^{1}\cdots \rho_{n+1}^k)
(\rho_n^1 \cdots \rho_n^k)(\rho^1_1 \cdots \rho_1^k)
\end{split}
\end{equation*}

All of the corresponding terms of both expressions are equal, and in particular
$\delta=\rho$.

To see the second claim, note that by definition
$$F_n= (\Delta^1)^n_{n-1} \cdots (\Delta^k)^n_{n-1} F_{n-1}.$$
So, $F_n= \bigcup \{ (\delta_n^1 \cdots \delta_n^k) F_{n-1}
\colon \delta_n^i \in (\Delta^i)^n_{n-1}\}$.
Since $\delta_n^1 \cdots \delta_n^k \in \Delta_{n-1}$, the first claim shows that
these translates of $F_{n-1}$ are disjoint.
\end{proof}

The next lemma show that the $\Delta_n$ translates of $F_n$ cover $G$.

\begin{lem}  \label{ch05_lem:existence}
For every $n$ and every $g \in G$, we may write $g$ in the form
$g=(\delta_m^1 \cdots \delta_m^k) \cdots (\delta_{n+1}^1 \cdots \delta_{n+1}^k)f$,
where $f =(\delta_n^1 \cdots \delta_n^k)\cdots (\delta^1_1 \cdots \delta_1^k) \in F_n$,
for some $m>n$.
\end{lem}

\begin{proof}
Fix $n$ and $g \in G$. Every element $\bx$ of a quotient $G^{i-1}/G^i$
can be written in the form $\bx=\brho_m^i \cdots \brho_1^i$ for some $m=m(x)$
(which depends on $x$), where $\brho_j^i \in (\bDelta^i)^j_{j-1}$,
From this and the fact that the identity element is in all
the $(\bDelta_i)^j_{j-1}$ it follows that we may write $g$ as:

$$
g=(\delta_{m_1}^1 \delta_{m_1-1}^1 \cdots \delta^1_1) g_1$$
where $g_1 \in G_1$. This defines the $\delta_j^1$ for $j \leq m_1$. For convenience in the
following argument, we set $\delta_j^1=1_G$ for $j>m_1$. Note from the definition of the $c_j^i$ that
$c_j^2$ is now defined for all $j$. Recall that $c_j^2=\delta_{j-1}^1\cdots \delta^1_1$.
Thus, $c_j^2=c_p^2$ for all $j, p \geq m_1+1$. That is, the $c_j^2$
are eventually constant.

Assume in general that we have defined
$$\delta_{m_1}^1,\dots, \delta^1_1,\ \delta_{m_2}^2,\dots, \delta_1^2,\dots,\
\delta_{m_{\ell-1}}^{\ell-1},\dots, \delta_1^{\ell-1}$$
with
$$
g=(\delta_{m_1}^1,\dots, \delta^1_1)({\delta_{m_2}^2}^{c_{m_2}^2},\dots, {\delta_1^2}^{c_1^2}),
\cdots ({\delta_{m_{\ell-1}}^{\ell-1}}^{c_{m_{\ell-1}}^{\ell-1}},\dots,
{\delta_1^{\ell-1}}^{c_1^{\ell-1}}) g_\ell
$$
where $g_\ell \in G_\ell$ and all the $\delta_j^i$ lie in $(\Delta^i)^j_{j-1}$.
Again for convenience set $\delta_j^i=1_G$ for $j > m_i$, for
$i=1,\dots, \ell-1$. Note that $c_j^i$ is defined for $i=1,\dots, \ell-1$
and all $j$.  Also, inspecting the formula for $c_j^i$ shows that
for all $i=1,\dots, \ell-1$ that $c_j^i$ is eventually
constant for large enough $j$. Let us call this eventual value $c_\infty^i$.
To finish the inductive step in the proof of the lemma it suffices to show that
in the quotient group $G^\ell/G^{\ell+1}$
the (arbitrary) element $\bg_\ell \in G^\ell/G^{\ell+1}$ can be written in the form
$$
\bg_\ell= {\bdelta_{m}^\ell}^{c_m^\ell} \cdots {\bdelta_1^\ell}^{c_1^\ell}
$$
for some $m$, where as usual $\delta_j^\ell \in (\Delta^\ell)^j_{j-1}$.
For the rest of the proof we work in the quotient group $G_\ell/G_{\ell+1}$
and we suppress writing the bars in the notation.
Conjugating by $(c_\infty^\ell)^{-1}$,
it suffices to show that an arbitrary element of the quotient group
$h={g_\ell}^{(c_\infty^\ell)^{-1}}$ can be written in the form
$$
h= {\delta_m^\ell}^{d_m^\ell} \cdots {\delta_1^\ell}^{d_1^\ell}$$
where the $d_j^\ell=1_{G_\ell/G_{\ell+1}}$ for large enough $j$, say for $j > j_0$.
To begin, write $h^{(d_1^\ell)^{-1}}=y \delta_1^\ell$
for some $\delta_1^\ell \in (\Delta^\ell)^1_0$ and $y \in \Delta_1^\ell$.
We can do this since the $\Delta_1^\ell$ translates of $F^\ell_1$ cover
$G_\ell/G_{\ell+1}$ (and recall $F^\ell_1=(\Delta^\ell)^1_0$). Conjugating this gives
$h=z {\delta_1^\ell}^{d_1^\ell}$ where $z= y^{d_1^\ell} \in \Delta_1^\ell$
by the invariance of $\Delta_1^\ell$ under automorphisms.
For the next step, since $z^{(d_2^\ell)^{-1}} \in \Delta_1^\ell$ we may write it as
$z^{(d_2^\ell)^{-1}}= w \delta_2^\ell$ for some $\delta_2^\ell \in (\Delta^\ell)^2_1$
and $w \in \Delta_2^\ell$. Conjugating gives $z= u {\delta_2^\ell}^{d_2^\ell}$.
So, $g= u {\delta_2^\ell}^{d_2^\ell} {\delta_1^\ell}^{d_1^\ell}$ where $u \in \Delta_2^\ell$.
Continuing in this manner we may write $g$ as
$g= v {\delta_{j_0}^\ell}^{d_{j_0}^\ell} \cdots {\delta_1^\ell}^{d_1^\ell}$,
where $v \in \Delta^\ell_{j_0}$. By the property of Lemma~\ref{imtiling}
of the $\Delta_j^\ell$, $v$ is of the form
$v= \delta_m^\ell \cdots \delta_{j_0+1}^\ell$ for some large enough $m$. Therefore,
$g= (\delta_m^\ell \cdots \delta_{j_0+1}^\ell)
({\delta_{j_0}^\ell}^{d_{j_0}^\ell} \cdots {\delta_1^\ell}^{d_1^\ell})$ and we are done
(since $d_j^\ell=1_{G_\ell/G_{\ell+1}}$ for $j> j_0$).
\end{proof}

This completes the proof of Lemma~\ref{ch05_lem:existence}
and of Theorem~\ref{ch05_thm:solvable}.
\end{proof}

The above method of proof does not immediately extend to solvable groups. The obstacle is that the statement of Lemma~\ref{ch05_lem:ablem} does not hold for abelian groups in general.
If $G=\bigoplus_{i=1}^\infty \mathbb{Z}$ is the infinite direct sum of copies of $\mathbb{Z}$, then it
is easy to check that the only sets $\Delta \subseteq G$ which are invariant under all automorphisms
are the sets $G^n$. However, aside from the trivial case $n=1$, none of these sets
can be the set of center points for a tiling of $G$ as $G/G^n$ is not finite.

We close this section with the following open question.

\begin{question}
Is every countable solvable group a \ccc group?
\end{question}

\section{Residually finite and locally finite groups and free products are ccc} \label{SECT RESFIPROD}

Recall the following definition for residually finite groups.

\begin{definition}
A group $G$ is \emph{residually finite} if the intersection of all finite index normal subgroups of $G$ is trivial.
\end{definition}\index{residually finite}

\begin{theorem}\label{thm:residuallyfiniteccc}
If $G$ is a countably infinite residually finite group, then $G$ is ccc.
\end{theorem}

\begin{proof}
Fix a strictly decreasing sequence $G = H_0 \geq H_1 \geq \cdots$ of finite index normal subgroups of $G$ with $\bigcap_{n \in \N} H_n = \{1_G\}$. Now fix an enumeration $g_0, g_1, \ldots$ of the nonidentity group elements of $G$ in a manner such that $g_i \not\in H_{i+1}$.

We will now construct, for each $i \in \N$, a complete set $\Delta^{i+1}_i$ of coset representatives for the right cosets of $H_{i+1}$ inside of $H_i$. In other words, $H_i = H_{i+1} \Delta^{i+1}_i$, and $H_{i+1} \lambda_1 \cap H_{i+1} \lambda_2 = \varnothing$ for distinct $\lambda_1, \lambda_2 \in \Delta^{i+1}_i$. The sets $\Delta^{i+1}_i$ will have some additional properties and must be constructed inductively. Let $\Delta^1_0$ be a complete set of coset representatives for the right cosets of $H_1$ inside of $H_0$ such that $1_G, g_0 \in \Delta^1_0$. Now suppose that $\Delta^1_0$ through $\Delta^{k+1}_k$ have been defined. An easy inductive argument shows that
$$H_{k+1} \Delta^{k+1}_k \Delta^k_{k-1} \cdots \Delta^1_0 = G.$$
Notice also that $|\Delta^{i+1}_i| = [H_i : H_{i+1}]$ and
$$|\Delta^{k+1}_k \Delta^k_{k-1} \cdots \Delta^1_0| = [G: H_{k+1}] < [G: H_{k+2}].$$
Therefore
$$H_{k+2} \Delta^{k+1}_k \Delta^k_{k-1} \cdots \Delta^1_0 \neq G.$$
Find the least $j \in \N$ such that
$$g_j \not\in H_{k+2} \Delta^{k+1}_k \Delta^k_{k-1} \cdots \Delta^1_0.$$
Since $H_{k+1} \Delta^{k+1}_k \Delta^k_{k-1} \cdots \Delta^1_0 = G$, we can find $\gamma \in \Delta^{k+1}_k \Delta^k_{k-1} \cdots \Delta^1_0$ such that $g_j \in H_{k+1} \gamma$. Finally, let $\Delta^{k+2}_{k+1}$ be a complete set of coset representatives for the right cosets of $H_{k+2}$ inside of $H_{k+1}$ such that $1_G, g_j \gamma^{-1} \in \Delta^{k+2}_{k+1}$. This completes the construction of the sets $\Delta^{i+1}_i$. Notice that each $\Delta^{i+1}_i$ is finite.

Now define $F_0 = \{1_G\}$, $\Delta_0 = G$, and for $n > 0$
$$F_n = \Delta^n_{n-1} \Delta^{n-1}_{n-2} \cdots \Delta^1_0,$$
$$\Delta_n = H_n.$$
We claim that $(\Delta_n, F_n)_{n \in \N}$ is a ccc sequence of tilings of $G$. Since $H_i = H_{i+1} \Delta^{i+1}_i$, by induction we easily have $\Delta_n F_n = H_0 = G$. So the $\Delta_n$-translates of $F_n$ cover $G$. If $h, h' \in H_n$, $\lambda_i, \lambda_i' \in \Delta^{i+1}_i$ and
$$h \lambda_{n-1} \lambda_{n-2} \cdots \lambda_0 = h' \lambda_{n-1}' \lambda_{n-2}' \cdots \lambda_0'$$
then after viewing this equation in $G/H_1$ we see that $\lambda_0 = \lambda_0'$ (equality in $G/H_1$ and hence equality in $G$). After canceling $\lambda_0$ and $\lambda_0'$ from both sides and viewing the new equation in $G/H_2$, we see that $\lambda_1 = \lambda_1'$. Continuing in this manner we find that $\lambda_i = \lambda_i'$ for each $0 \leq i < n$ and $h = h'$. Therefore the $\Delta_n$-translates of $F_n$ are disjoint. We conclude that for each $n \in \N$ $(\Delta_n, F_n)$ is a tiling of $G$.

Clearly $1_G \in H_n = \Delta_n$, so the tilings are centered. Coherency is also clear since for $n > 1$ and $h \in H_n$ we have
$$h F_n = h \Delta^n_{n-1} F_{n-1}$$
and $h \Delta^n_{n-1} \subseteq H_{n-1} = \Delta_{n-1}$.

All that remains is to check cofinality. Notice that $F_n \subseteq F_{n+1}$ since $1_G \in \Delta^{n+1}_n$, and also notice $1_G \in F_0$. So we only have to show that for every $j \in \N$ there is $n \in \N$ with $g_j \in F_n$. Towards a contradiction, suppose the sequence of tilings is not cofinal. Let $j \in \N$ be least such that $g_j \not\in F_n$ for all $n \in \N$. Let $k \in \N$ be such that $g_i \in F_k$ for all $i < j$ ($k=0$ if $j=0$). By the way the sets $\Delta^{i+1}_i$ were constructed, we must have for all $n > k$
$$g_j \in H_{n+1} F_n = \Delta_{n+1} F_n$$
(otherwise for some $\gamma \in F_n$, $g_j \gamma^{-1} \in \Delta^{n+1}_n$ and $g_j \in F_{n+1}$). Let $\gamma \in \Delta_k$ be such that $g_j \in \gamma F_k$. We have $g_j \in \Delta_{k+1} F_k$, so there is $\sigma \in \Delta_{k+1}$ with $g_j \in \sigma F_k$. As $\Delta_{k+1} = H_{k+1} \subseteq H_k = \Delta_k$, we have that $\sigma \in \Delta_k$. Then $g_j \in \gamma F_k \cap \sigma F_k$. Therefore $\sigma = \gamma$ and $\gamma \in \Delta_{k+1}$. Repeating this argument, we find that $\gamma \in \Delta_n$ for all $n \geq k$. Thus
$$\gamma \in \bigcap_{n \geq k} \Delta_n = \bigcap_{n \geq k} H_n = \{1_G\}.$$
Now we have $g_j \in \gamma F_k = F_k$, a contradiction. We conclude that the sequence of tiles are cofinal.
\end{proof}

By a theorem of Gruenberg (\cf\ \cite{M}) free products of residually finite groups are residually finite, hence they are also \ccc by Theorem~\ref{thm:residuallyfiniteccc}. All finitely generated nilpotent groups and all polycyclic groups are residually finite. Hence Theorem~\ref{thm:residuallyfiniteccc} proves again that these groups are ccc. Also, all free groups are residually finite, hence are \ccc also by Theorem~\ref{thm:residuallyfiniteccc}. Finally, by the theorem of Gruenberg mentioned above, free products of finite groups are residually finite, and hence are also ccc.

Now we show that countable locally finite groups are \ccc groups.

\begin{theorem}
If $G$ is a countably infinite locally finite group then $G$ is a \ccc group.
\end{theorem}

\begin{proof}
Fix an increasing sequence $A_0 \subseteq A_1 \subseteq \cdots$ of finite subsets of $G$ with $G = \bigcup_{n \in \N} A_n$. For each $n \in \N$ set $F_n = \langle A_n \rangle$. Then we have $G = \bigcup_{n \in \N} F_n$, and since $G$ is locally finite each $F_n$ is finite. For each $n \geq 1$ let $\Delta^n_{n-1}$ be a complete set of coset representatives for the left cosets of $F_{n-1}$ in $F_n$. In other words, $F_n = \Delta^n_{n-1} F_{n-1}$ and $\delta F_{n-1} \cap \delta' F_{n-1} = \varnothing$ for $\delta \neq \delta' \in \Delta^n_{n-1}$. We further require that $1_G \in \Delta^n_{n-1}$ for all $n \geq 1$. Now we set
$$\Delta_n = \bigcup_{m > n} \Delta^m_{m-1} \Delta^{m-1}_{m-2} \cdots \Delta^{n+1}_n.$$
Notice that the members of the above union are increasing since $1_G$ is in each $\Delta^{k+1}_k$. We claim that $(\Delta_n, F_n)_{n \in \N}$ is a \ccc sequence of tilings of $G$.

Clearly $F_{n+1} = \Delta^{n+1}_n F_n$, and it easily follows by induction that
$$\Delta^m_{m-1} \Delta^{m-1}_{m-2} \cdots \Delta^{n+1}_n F_n = \Delta^m_{m-1} F_{m-1} = F_m.$$
Therefore
$$\Delta_n F_n = \bigcup_{m > n} \Delta^m_{m-1} \Delta^{m-1}_{m-2} \cdots \Delta^{n+1}_n F_n = \bigcup_{m > n} F_m = \bigcup_{m \in \N} F_m = G.$$
By definition we have that the $\Delta^{n+1}_n$-translates of $F_n$ are disjoint and are contained in $F_{n+1}$. Since the $\Delta^{n+2}_{n+1}$-translates of $F_{n+1}$ are disjoint (and contained in $F_{n+2}$), it follows that the $\Delta^{n+2}_{n+1} \Delta^{n+1}_n$-translates of $F_n$ are disjoint (and contained in $F_{n+2}$). Inductively assume that the $\Delta^{m-1}_{m-2} \Delta^{m-2}_{m-3} \cdots \Delta^{n+1}_n$-translates of $F_n$ are disjoint and contained in $F_{m-1}$. Since the $\Delta^m_{m-1}$-translates of $F_{m-1}$ are disjoint and contained in $F_m$, it follows that the $\Delta^m_{m-1} \Delta^{m-1}_{m-2} \cdots \Delta^{n+1}_n$-translates of $F_n$ are disjoint and contained in $F_m$. So by induction and by the definition of $\Delta_n$, it follows that the $\Delta_n$-translates of $F_n$ are disjoint. Thus $(\Delta_n, F_n)$ is a tiling of $G$ for each $n \in \N$. If $\gamma \in \Delta_{n+1}$ then
$$\gamma F_{n+1} = \gamma \Delta^{n+1}_n F_n.$$
Since $\gamma \in \Delta_{n+1}$, there is $m > n$ with
$$\gamma \in \Delta^m_{m-1} \Delta^{m-1}_{m-2} \cdots \Delta^{n+2}_{n+1}$$
and therefore $\gamma \Delta^{n+1}_n \subseteq \Delta_n$. Thus, every $\Delta_{n+1}$-translate of $F_{n+1}$ is the union of $\Delta_n$-translates of $F_n$. So $(\Delta_n, F_n)_{n \in \N}$ is a coherent sequence of tilings of $G$. Since $1_G$ is in each $\Delta^{k+1}_k$, we have $1_G \in \Delta_n$ for all $n \in \N$. Also, $G = \bigcup_{n \in \N} F_n$. Thus $(\Delta_n, F_n)_{n \in \N}$ is a \ccc sequence of tilings of $G$ and $G$ is a \ccc group.
\end{proof}

In the rest of this section we consider arbitrary countable free products of nontrivial countable groups and show that they are ccc. We give a proof that is combinatorial by nature, unlike the algebraic construction in the proof of Theorem~\ref{thm:residuallyfiniteccc}. To illustrate the combinatorial argument, we will first give a direct proof that all free groups are ccc.

Recall that, for $n\geq 2$ an integer, $\mathbb{F}_n$ denotes the free group on $n$ generators,
and $\mathbb{F}_\omega$ the free group on countably infinitely many generators.

\begin{theorem} \label{ch05_thm:free}
Every free group $\mathbb{F}_n$ or $\mathbb{F}_\omega$ is a \ccc group.
\end{theorem}

Before turning to the proof of Theorem~\ref{ch05_thm:free} we fix some notation.
If $G=\mathbb{F}_n$ or $G=\mathbb{F}_\omega$ is a free group, let $T$ be the Cayley graph which in this case
is a tree. Recall that if $G=\mathbb{F}_n$ then the nodes of the graph are the elements of $G$
and two elements $y$ and $z$ of $G$ are linked by an edge if either $y=z x_i$ or $y=z (x_i)^{-1}$,
where $x_1,\dots x_n$ are the generators of $G$. Every node has degree $2n$ ($\infty$ if
$G=\mathbb{F}_\omega$). We view $T$ as a rooted tree with root node corresponding to the identity element
$1_G$ of $G$. For $x \in G$ we let the {\em depth} of $x$ be the distance from $x$ to $1_G$
in $T$, and denote it by $d(x)$. This, of course, is just the length of $x$ when expressed as a reduced word in the
generators. This also can be thought of as the depth of $x$ as a node in the rooted tree $T$.
The children of a node $x \in G$ are the nodes adjacent to $x$ with depth $d(x)+1$. For $x \neq 1_G$, the parent of
$x$ is the unique node adjacent to $x$ with depth $d(x)-1$.
The root node $1_G$ has $2n$ children and no parent, every other node has one parent and
$2n-1$ children. We say $S \subseteq T$ is a {\em subtree} if $1_G \in S$ and
$S$ is closed under the parent relation, that is, if $x \in S$ and $y$ is a parent of $x$,
then $y \in S$.
Viewing $T$ as a rooted tree as above, this corresponds to the usual notion of subtree of a tree.
We will identify $G$ and $T$, and so speak of subtrees of $G$ as well.

The following simple lemma is the main point.

\begin{lem} \label{ch05_lem:fl}
Let $S \subseteq G$ be a subtree of the free group $G$, Then $G$ can be tiled
by copies of $S$.
\end{lem}

\begin{proof}
Let $\sM$ be a maximal pairwise disjoint collection of translates of $S$ subject to
the condition that $M \defeq \bigcup \sM$ is a subtree of $T$. It suffices to
show that $M=T$. If not, let $z\in T-M$ have minimal depth. Clearly $z \neq 1_G$,
and so the unique parent $y$ of $z$ is in $M$. We must have $z=y x_i$
or $z=y (x_i)^{-1}$ for some generator $x_i$, and where $y$, when written as a
reduced word, does not end in the term canceling this last term. Suppose
to be specific $z=yx_i$, the other case being similar. Let $k$ be such that
$(x_i^{-1})^k \in S$ but $(x_i^{-1})^{k+1} \notin S$. Let $w= z \cdot (x_i)^k$. Note that $w$ is
in reduced form as written. Consider $\sM'=\sM \cup \{ w S\}$. We claim that
$\sM'$ is still pairwise disjoint and $M'\defeq \bigcup \sM'$ is still a tree, which then
contradicts the maximality of $\sM$.

An element $h$ of $wS$ is of the form $z (x_i)^k s$, where $s \in S$.
Note that $s$ cannot begin with $(x_i^{-1})^{k+1}$ as otherwise, since $S$ is a tree,
$(x_i^{-1})^{k+1}$ would be in $S$. Since $z$ also ends with $x_i$, it follows that
the reduced form for $h$ is of the form $h= z u$ for some (possibly empty) word $u$.
So, $h$ is a descendant of $z$, and since $z \notin M$, it follows that
$h \notin M$ as well. So, $wS \cap M=\emptyset$.

Note that the path from $w$ to $y$ lies entirely in $M'$ since $(x_i^{-1})^\ell \in S$
for all $\ell \leq k$. So, if $h \in wS$ lies on this branch, then any initial segment
of $h$ lies in $M'$. If $h$ is not on this path, then $h$ is of the form
$h= y \cdot x_i \cdot (x_i)^k \cdot (x_i^{-1})^{-\ell} u$ where
$\ell \leq k$, $u$ does not start with $x_i^{-1}$,
and $(x_i^{-1})^{-\ell} u \in S$. Any initial segment $h'$ of $h$ is either on the path
from $w$ to $z$ or else is of the form
$h'= y \cdot x_i \cdot (x_i)^k \cdots (x_i^{-1})^{-\ell} u'$ where $u'$ is an initial segment of $u$.
In the latter case, $(x_i^{-1})^{-\ell} \cdot u' \in S'$ as $S$ is a tree, and so
$h' \in wS \subseteq M'$.
\end{proof}

\begin{proof}[Proof of Theorem~\ref{ch05_thm:free}]
Let $G$ be a free group, and $T$ the corresponding tree as above. Let
$g_0,g_1,\dots$ enumerate $G$. Let $S_0 \subseteq T$ be an arbitrary subtree of $T$.
Suppose after step $i$ the subtree $S_i$ of $T$ has been defined with
$S_0 \subseteq S_1 \subseteq \cdots \subseteq S_i$. Suppose also that each $S_j$ for
$j=1,\dots, i$ is a disjoint union of $S_{j-1}$ and another translate $w_{j-1} S_{j-1}$
of $S_{j-1}$.

To define $S_{i+1}$, let $g_i$ be least in the enumeration such that $g_i \notin S_i$ but
the parent of $g_i$ is in $S_i$. Let $w_i$ be the element constructed in the
proof of Lemma~\ref{ch05_lem:fl} using $S=S_i$ and $z=g_i$. So,
$S_i$ and  $w_i S_i$ are disjoint,  $g_i \in w_iS_i$, and $S_i \cup w_i S_i$ is a subtree of $T$. We then let $S_{i+1}=S_i\cup w_iS_i$.
This finishes the inductive definition of $S_i$ and $w_i$.

From the definition of $g_i$ it follows easily that
$G=\bigcup_i S_i$. This in turn gives our \ccc tiling of $G$.
Namely, the $i$-th level tiling  will
be $(\Delta_i, S_i)$ where
$\Delta_i=\bigcup_{m>i} \Delta^m_{m-1} \Delta^{m-1}_{m-2} \cdots \Delta^{i+1}_i$
where $\Delta^j_{j-1}=\{1_G, w_{j-1}\}$. Note that $S_m=\Delta^m_{m-1} S_{m-1}=\cdots
= \Delta^m_{m-1} \Delta^m_{m-1} \cdots \Delta^{i+1}_i S_i$. Thus, $G=\Delta_i S_i$ and the $\Delta_i$ translates of $S_i$ are pairwise disjoint,
so $(\Delta_i,S_i)$ is a tiling of $G$. Clearly $1_G \in \Delta_i$ as $1_G \in \Delta^{i+1}_i$.
Also, the tilings are coherent since every $\delta S_{i+1}$, for $\delta
\in \Delta_{i+1}$ is of the form
$\delta^m_{m-1} \delta^{m-1}_{m-2} \cdots \delta^{i+2}_{i+1} S_{i+1}$, where $\delta^j_{j-1} \in
\Delta^j_{j-1}$. This is then equal to
$\delta^m_{m-1} \delta^{m-1}_{m-2} \cdots \delta^{i+2}_{i+1} \Delta^{i+1}_i S_i$,
which is the disjoint union of
$\delta^m_{m-1} \delta^{m-1}_{m-2} \cdots \delta^{i+2}_{i+1} S_i$ and
$\delta^m_{m-1} \delta^{m-1}_{m-2} \cdots \delta^{i+2}_{i+1} w_i S_i$. Thus, the sequence $(\Delta_i,S_i)$
gives a \ccc tiling of $G$.
\end{proof}

We now turn to free products. Recall the following definition.

\begin{definition}
Let $\sG=\{ G_i\}_{i \in \sI} $ be a collection of groups.
We assume the $G_i$ are pairwise disjoint as sets.
The {\em free product} $\ast \sG=\ast_i  G_i$ of the collection is the group with generators
$\bigcup_i G_i$ and relations $g \cdot h =k$ for all $g,h,k$ in some common $G_i$
with $g \cdot h=k$ in $G_i$.
\end{definition}\index{free product}\index{free product!trivial}

If all but one group $G$ in $\sG$ are trivial (i.e., contain only one element), then we say that $\sG$ is {\em trivial} and the free product $\ast \sG$
is just isomorphic to $G$. In this case, $\ast \sG$ is \ccc iff $G$ is.
If $\sG$ is nontrivial, the next theorem says that $\ast \sG$
always is a \ccc group.

\begin{theorem}  \label{freeproduct}
If $\sG$ is a countable nontrivial collection of countable
groups, then the free product $\ast \sG$ is a \ccc group.
\end{theorem}

\begin{proof}
To ease notation we consider the case where $\sG=\{ G,H\}$ with both $G$ and $H$ nontrivial, and we denote the free product
in this case by $G\ast H$. The general case is entirely similar.

We first consider the case that one of the groups, say $G$, is infinite.
Every nonidentity element $x$ of $G *H$ can be written uniquely in the form
$g_1 h_1 g_2 h_2 \cdots k_i$ or $h_1 g_1 h_2 g_2 \cdots k_i$ where the $g_i, h_i$
are nonidentity elements of $G$, $H$ respectively, and $k_i$ is in $G$ or $H$ depending
on whether the word length is odd or even. We denote the word length of $x$ by $\lh(x)$. When $l\leq \lh(x)$ we use $x\upharpoonright l$ to denote
the initial segment of $x$ of word length $l$.\index{$G$-tail}


If $F \subseteq G*H$, we say $x \in G*H$ is a $G$-{\em tail} (with respect to $F$) if $x$ is not the identity, $x$ ends
with a term $k_i \in G$, and $x$ is the unique element of $F$ which extends
$x \res (\lh(x)-1)$. We likewise define the notion of an $H$-tail.
We say $F$ is a {\em tree} if $1 \in F$ and any initial segment of an element of $F$ (when written
in one of the above two forms) is also an element of $F$.\index{tree}

The next lemma says that we may create as many new tails as we like.

\begin{lem} \label{tail}
Suppose $F \subseteq G*H$ is a finite tree with at least one $G$-tail and one $H$-tail.
Then for any $k \in \N$, there is an $F' \subseteq G*H$ satisfying the following.

\begin{enumerate}
\item
$F' \supseteq F$ is a finite tree.
\item
$F'$ is a disjoint union of translates $\delta F$ of $F$.
\item
$F'$ has at least $k$ many $G$-tails and at least $k$ many $H$-tails.
\end{enumerate}
\end{lem}

\begin{proof}
Let $A \subseteq G$ be the set of all $g \in G$ such that some $x \in F$,
when written in its reduced form, begins with $g$. Since $F$ is finite, so is $A$.
Let $g_1,\dots, g_k \in G$ be such that $g_iA\cap A=\emptyset$ and $g_i A \cap g_j A= \emptyset$
for all $i \neq j$. We can do this as $A$ is finite and $G$ is infinite.
Consider the collection of translates $F, g_1 F, \cdots, g_k F$. Since
$g_ iA \cap g_j A=\emptyset$ for all $i \neq j$, it follows that for any elements
$f_i, f_j \in F$ that $g_i f_i$ and $g_j f_j$ begin with different elements
of $G$ and so are not equal. So, $g_i F \cap g_j F=\emptyset$. Similarly,
$F \cap g_iF=\emptyset$ for all $i$. Since $F$ is a tree, it is easy to see
that $F \cup g_1F \cup \cdots \cup g_k F$ is also a tree (note that
each $g_i$ is in this union as $1_G \in F$). Finally,
each $g_i F$ has at least one $G$-tail and at least one $H$-tail since $F$
does (if say $g_1 h_1\dots g_k$ is a $G$-tail in $F$, then easily
$g_i g_1 h_1 \dots g_k$ is a $G$-tail in $g_iF$).
\end{proof}

Turning to the proof of Theorem~\ref{freeproduct}, let $z_0,z_1,\dots$
enumerate $G*H$. Assume inductively we have constructed
$F_0 \subseteq F_1 \subseteq \cdots \subseteq F_n$ satisfying:

\begin{enumerate}
\item
Each $F_i$ is a finite tree.
\item
Each $F_{i+1}$ is a disjoint union of translates $\delta F_i$ of $F_i$.
\item
Each $F_i$ has at least one $G$-tail and one $H$-tail.
\end{enumerate}

To construct $F_{n+1}$, Let $z$ be the least element of $G*H$ (in the enumeration
$z_0, z_1,\dots$) such that $z \notin F_n$ but $y \defeq z \res (\lh(z)-1) \in F_n$.
Say to be specific $z=y g$ where $g \in G$ (the case $z=y h$ is similar). From
Lemma~\ref{tail} we may assume that $F_n$ has at least $2$ $G$-tails and
at least $2$ $H$-tails. Let $t \in F_n$ be an $H$-tail relative to $F_n$.
Consider the translate $z t^{-1} F_n$. Since $t$ ends with an $H$ term,
$t^{-1}$ begins with an $H$ term, so $z t^{-1}$ is in reduced form.
Since $t \in F_n$, $z \in z t^{-1} F_n$. Since $t$ is the unique element of
$F_n$ which extends $t \res (\lh(t)-1)$, every element of
$z t^{-1} F_m$ when written in reduced form is of the form
$z u$ for some possibly empty term $u$. Since $F_n$ is a tree, none
of these elements can be in $F_n$. So, $zt^{-1}F_n \cap F_n=\emptyset$.
We set $F_{n+1}=F_n \cup zt^{-1} F_n$. To see that $F_{n+1}$ is a tree,
first note that the path from $z$ to $z t^{-1}$ lies entirely in $F_{n+1}$
This is because every element of this path is of the form
$z t^{-1} v$ where $v$ is an initial segment of $t$, and hence lies in $F_n$.
Every element of $z t^{-1}F_n$ not on this path is of the form
$x=z t^{-1} u w$ where $u w \in F_n$ and $u$ is the maximal initial segment of
$u w$ which cancels an end segment of $t^{-1}$. That is, the reduced form for $x$ is
$x= z s w$ where $s$ is an initial segment of $t^{-1}$. So, an initial segment
of $x$ is either on the path from $z$ to $z t^{-1}$ or else it is of the form
$z s w'= z t^{-1} u w'$ where $w'$ is an initial segment of $w$.
Since $F_n$ is a tree, $u w' \in F_n$, and so in either case we have
$x \in z t^{-1} F_n \subseteq F_{n+1}$. This shows $F_{n+1}$ is a tree.

Finally, $F_{n+1}$ has at most one fewer tail than $F_n$. Namely, every
$G$- or $H$-tail of $F_n$ except possibly $y$ is still a tail of $F_{n+1}$.
So, $F_{n+1}$ satisfies all of our inductive assumptions and also
contains $z$. As in previous arguments, this gives a \ccc set of tilings
for $G*H$. Namely, $(\Delta_n, F_n)$, where
$\Delta_n= \bigcup_{m>n} \Delta^m_{m-1} \Delta^{m-1}_{m-2} \cdots \Delta^{n+1}_n$,
where $\Delta^{i+1}_i$ is the finite set such that $F_{i+1}=\Delta^{i+1}_i F_i$ and
the $\Delta^{i+1}_i$ translates of $F_i$ are pairwise disjoint.

The case where both $G$ and $H$ are finite follows from Theorem~\ref{thm:residuallyfiniteccc},
since $G\ast H$ is residually finite. But here we point out that the above argument can be easily modified to
give a direct proof also, as follows. We modify assumption (3) of the $F_n$ to now be that
each $F_n$ has at least two $G$-tails and two $H$-tails. Note that we may always start with
a finite subtree $F_0$ of $G*H$ with at least two $G$ tails and at least two $H$ tails except in the
case where both $G$ and $H$ have size $2$. For the group $\mathbb{Z}_2 * \mathbb{Z}_2$
it is easy to directly construct a \ccc tiling. For example, one can let $F_n$ be the set of size $3^n$
consisting of all reduced words of length $\leq (3^n-1)/2$. Easily $F_{n+1}$ is a disjoint
union of three translates of $F_n$
(in fact, this group
has an index two subgroup isomorphic to $\mathbb{Z}$ and is also polycyclic, so both
Lemma~\ref{lem:Hxycoloring} and Theorem~\ref{ch05_thm:solvable} apply).
To get $F_{n+1}$, as in the above proof let
$z$ be the least element of $G*H$ such that $z \notin F_n$ but $y \defeq z \res (\lh(z)-1) \in F_n$.
Consider the case $z=yg$ where $g \in G$, the other case being similar.
Let $t \in F_n$ be an $H$-tail relative to $F_n$. We again consider the translate
$zt^{-1}F_n$ and let $F_{n+1}=F_n \cup zt^{-1}F_n$. The argument above shows that $F_{n+1}$
is a tree which is the union of two disjoint translates of $F_n$. Also, at most one of the
tails of $F_n$ is not a tail of $F_{n+1}$. So, it suffices to observe that $F_{n+1}$
has at least one $G$-tail, and also at least one $H$-tail, which are not in $F_n$.
Let $z_1$, $z_2$ be two $H$-tails of $F_n$. Then at least one of
$zt^{-1} z_1$, $zt^{-1}z_2$ is an $H$-tail of $F_{n+1}$ which is not in $F_n$.
This is because any two distinct tails must be incompatible (i.e., when written in reduced form neither
word is an initial segment of the other), and so at least one of $z_1$, $z_2$ is incompatible with
$t$. If, say, $z_1$ is incompatible with $t$, then it is easy to check that $zt^{-1}z_1$ is an $H$-tail of
$F_{n+1}$ which is not in $F_n$. The argument for $G$-tails is similar.

This completes the proof of Theorem~\ref{freeproduct}.
\end{proof}

\chapter{\label{CHAP FM}Blueprints and Fundamental Functions}

This chapter is the backbone to a variety of results we prove in the rest of this paper. In this chapter we  present a powerful and customizable method for constructing elements of $2^G$ with special properties. The concept of a blueprint, a special sequence of regular marker structures, is introduced in the first section. Blueprints simply organize the group theoretic structure of the group $G$ and allows one to carry out sophisticated constructions of elements of $2^G$. In our case, all of these constructions stem from one main construction which appears in the second section. In the third section, it is shown that every countably infinite group admits a blueprint. Finally, in the fourth section we study the growth rate of blueprints and how this impacts the main construction. This entire chapter is very abstract and simply develops tools for later use. It may initially be difficult to appreciate, understand, and see the motivation for what we do in this chapter, however the reader will be greatly rewarded in the next chapter.

\section{Blueprints}

Fundamental functions were originally developed for constructing 2-colorings on arbitrary groups and therefore they have their roots in the methods appearing in Section 4.2. Section 4.2 is not a prerequisite, but it would certainly aid in understanding on an intuitive level what our course of approach is. In this section, we study countable groups themselves under the notion of a blueprint. Blueprints are sequences of regular marker structures which have much more structure than those constructed in Section 4.2. Blueprints have enough structure to allow us to create partial functions on $G$ with very nice properties. This construction appears in the next section and is easily used in the next chapter to construct a $2$-coloring. Our first definition will be central to our studies for the rest of the paper.

\begin{definition} \label{DEFN MAX DIS} \index{maximally disjoint}
Let $G$ be a group and let $A, B, \Delta \subseteq G$. We say that the $\Delta$-translates of $A$ are \emph{maximally disjoint within} $B$ if the following properties hold:
\begin{enumerate}
\item[(i)] for all $\gamma, \psi \in \Delta$, if $\gamma \neq \psi$ then $\gamma A \cap \psi A = \varnothing$;
\item[(ii)] for every $g \in G$, if $g A \subseteq B$ then there exists $\gamma \in \Delta$ with $g A \cap \gamma A \neq \varnothing$.
\end{enumerate}
When property (i) holds we say that the $\Delta$-translates of $A$ are \emph{disjoint}. Furthermore, we say that the $\Delta$-translates of $A$ are \emph{contained and maximally disjoint within} $B$ if the $\Delta$-translates of $A$ are maximally disjoint within $B$ and $\Delta A \subseteq B$.
\end{definition}

Notice that in the definition above we were referring to the left translates of $A$ by $\Delta$ but never explicitly used the term left translates. Throughout the rest of this paper when we use the word translate(s) it will be understood that we are referring to left translate(s). When referring to right translates we will explicitly write out ``right-translate(s)''. Note that in the definition above there is no restriction on $\Delta$ being nonempty. So at times it may be that the $\varnothing$-translates of $A$ are contained and maximally disjoint within $B$.

\begin{definition} \label{def:blueprint} \index{blueprint} \index{blueprint!maximally disjoint} \index{blueprint!centered} \index{blueprint!directed} \index{pre-blueprint} \index{pre-blueprint!maximally disjoint} \index{pre-blueprint!centered} \index{pre-blueprint!directed}\index{coherent}\index{maximally disjoint}\index{directed}\index{uniform}\index{growth}\index{dense}
\index{centered}
Let $G$ be a countably infinite group. A \emph{blueprint} is a sequence $(\Delta_n, F_n)_{n \in \N}$ of regular marker structures satisfying the following conditions:
\begin{enumerate}
\item[\rm (i)] (disjoint) for every $n \in \N$ and distinct $\gamma, \psi \in \Delta_n$, $\gamma F_n \cap \psi F_n = \varnothing$;
\item[\rm (ii)] (dense) for every $n \in \N$ there is a finite $B_n \subseteq G$ with $\Delta_n B_n = G$;
\item[\rm (iii)] (coherent) for $k \leq n$, $\gamma \in \Delta_n$, and $\psi \in \Delta_k$, $\psi F_k \cap \gamma F_n \neq \varnothing \Longleftrightarrow \psi F_k \subseteq \gamma F_n$;
\item[\rm (iv)] (uniform) for $k < n$ and $\gamma, \sigma \in \Delta_n$, $\gamma^{-1} (\Delta_k \cap \gamma F_n) = \sigma^{-1} (\Delta_k \cap \sigma F_n)$;
\item[\rm (v)] (growth) for every $n > 0$ and $\gamma \in \Delta_n$, there are distinct $\psi_1, \psi_2, \psi_3 \in \Delta_{n-1}$ with $\psi_i F_{n-1} \subseteq \gamma F_n$ for each $i = 1,2,3$.
\end{enumerate}
If a sequence $(\Delta_n, F_n)_{n \in \N}$ satisfies all requirements for being a blueprint except for (ii), then it is called a \emph{pre-blueprint}. Furthermore, a (pre-)blueprint $(\Delta_n, F_n)_{n \in \N}$ is
\begin{enumerate}
\item[\rm (1)] \emph{maximally disjoint} if the $\Delta_n$-translates of $F_n$ are maximally disjoint within $G$;
\item[\rm (2)] \emph{centered} if $1_G \in \Delta_n$ for every $n \in \N$;
\item[\rm (3)] \emph{directed} if for every $k \in \N$ and $\psi_1, \psi_2 \in \Delta_k$ there is $n > k$ and $\gamma \in \Delta_n$ with $\psi_1 F_k \cup \psi_2 F_k \subseteq \gamma F_n$.
\end{enumerate}
\end{definition}


\begin{figure}[ht]
\begin{center}
\setlength{\unitlength}{6mm}
\begin{picture}(30,11)(-2,-0.5)

\put(14,9){\makebox(0,0)[b]{$F_1$}}
\qbezier(0,4.5)(0,11)(10.5,11)
\qbezier(0,4.5)(0,0)(7.5,0)
\qbezier(7.5,0)(16,0)(16,7.5)
\qbezier(10.5,11)(16,11)(16,7.5)

\put(-1,-0.2){
\put(5,5){\makebox(0,0)[b]{$\alpha_1F_0$}}
\qbezier(3,4.5)(3,3)(4.5,3)
\qbezier(3,4.5)(3,6)(4.5,6)
\qbezier(4.5,3)(6,3)(6,4.5)
\qbezier(4.5,6)(6,6)(6,4.5)
\put(4,4){\circle*{0.1}}
\put(5.1,3.8){\makebox(0,0)[b]{$\alpha_1\!=\!a_1$}}
}

\put(2.5,-2){
\put(5,5){\makebox(0,0)[b]{$\beta_1 F_0$}}
\qbezier(3,4.5)(3,3)(4.5,3)
\qbezier(3,4.5)(3,6)(4.5,6)
\qbezier(4.5,3)(6,3)(6,4.5)
\qbezier(4.5,6)(6,6)(6,4.5)
\put(4,4){\circle*{0.1}}
\put(5.1,3.8){\makebox(0,0)[b]{$\beta_1\!=\!b_1$}}
}

\put(4.2,1.3){
\put(5,5){\makebox(0,0)[b]{$\gamma_1 F_0$}}
\qbezier(3,4.5)(3,3)(4.5,3)
\qbezier(3,4.5)(3,6)(4.5,6)
\qbezier(4.5,3)(6,3)(6,4.5)
\qbezier(4.5,6)(6,6)(6,4.5)
\put(4,4){\circle*{0.1}}
\put(4.5,3.8){\makebox(0,0)[b]{$\gamma_1$}}
}

\put(1.1,2.9){
\put(5,5){\makebox(0,0)[b]{$\lambda F_0$}}
\qbezier(3,4.5)(3,3)(4.5,3)
\qbezier(3,4.5)(3,6)(4.5,6)
\qbezier(4.5,3)(6,3)(6,4.5)
\qbezier(4.5,6)(6,6)(6,4.5)
\put(4,4){\circle*{0.1}}
\put(4.4,3.8){\makebox(0,0)[b]{$\lambda$}}
}

\put(5.5,4.5){
\put(5,5){\makebox(0,0)[b]{$\lambda F_0$}}
\qbezier(3,4.5)(3,3)(4.5,3)
\qbezier(3,4.5)(3,6)(4.5,6)
\qbezier(4.5,3)(6,3)(6,4.5)
\qbezier(4.5,6)(6,6)(6,4.5)
\put(4,4){\circle*{0.1}}
\put(4.4,3.8){\makebox(0,0)[b]{$\lambda$}}
}

\put(6.5,-1.5){
\put(5,5){\makebox(0,0)[b]{$\lambda F_0$}}
\qbezier(3,4.5)(3,3)(4.5,3)
\qbezier(3,4.5)(3,6)(4.5,6)
\qbezier(4.5,3)(6,3)(6,4.5)
\qbezier(4.5,6)(6,6)(6,4.5)
\put(4,4){\circle*{0.1}}
\put(4.4,3.8){\makebox(0,0)[b]{$\lambda$}}
}

\put(13,6){\makebox(0,0)[b]{$\cdots\cdots$}}
\put(5.8,4.5){\circle{0.1}}
\put(6.3,4.3){\makebox(0,0)[b]{$1_G$}}

\end{picture}
\caption{\label{fig:F1}The composition of $F_1$. In the figure $\lambda$ denotes a generic element of $\Lambda_1$ (see Remark~\ref{remark513}).
All the solid points form the set $D^1_0$. The position of $1_G$ in $F_1$ can be arbitrary.}
\end{center}
\end{figure}
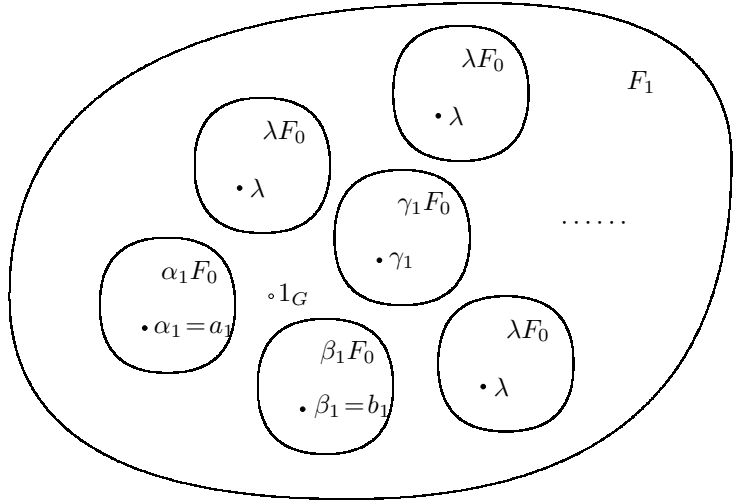


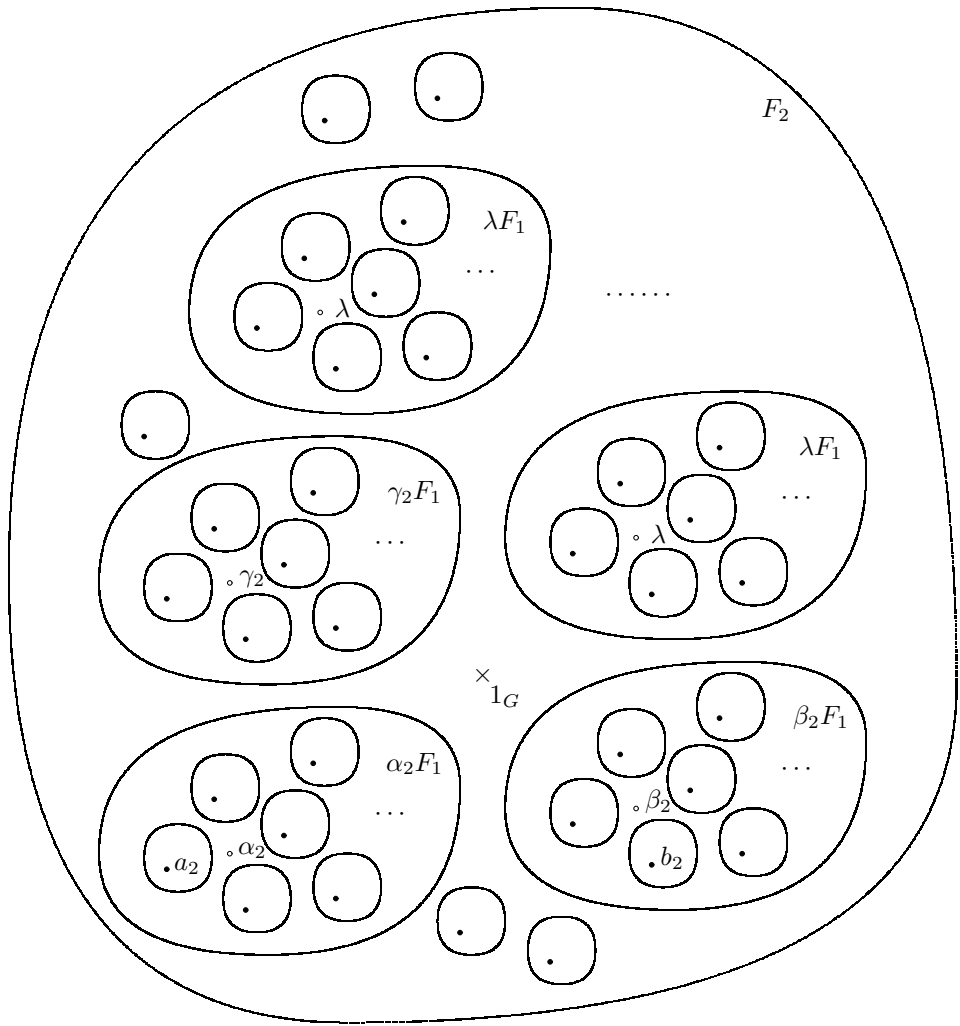
\begin{figure}[ht]
\begin{center}
\setlength{\unitlength}{3mm}
\begin{picture}(30,45)(5,-.5)

\put(34,40){\makebox(0,0)[b]{$F_2$}}
\qbezier(0,20)(0,45)(25,45)
\qbezier(0,20)(0,0)(15,0)
\qbezier(15,0)(42,0)(42,15)
\qbezier(25,45)(42,45)(42,15)

\put(4,3){
\put(14,8){\makebox(0,0)[b]{$\alpha_2F_1$}}
\qbezier(0,4.5)(0,11)(10.5,11)
\qbezier(0,4.5)(0,0)(7.5,0)
\qbezier(7.5,0)(16,0)(16,7.5)
\qbezier(10.5,11)(16,11)(16,7.5)

\put(-1,-0.2){
\qbezier(3,4.5)(3,3)(4.5,3)
\qbezier(3,4.5)(3,6)(4.5,6)
\qbezier(4.5,3)(6,3)(6,4.5)
\qbezier(4.5,6)(6,6)(6,4.5)
\put(4,4){\circle*{0.2}}
\put(4.9,3.8){\makebox(0,0)[b]{$a_2$}}
}

\put(2.5,-2){
\qbezier(3,4.5)(3,3)(4.5,3)
\qbezier(3,4.5)(3,6)(4.5,6)
\qbezier(4.5,3)(6,3)(6,4.5)
\qbezier(4.5,6)(6,6)(6,4.5)
\put(4,4){\circle*{0.2}}
}

\put(4.2,1.3){
\qbezier(3,4.5)(3,3)(4.5,3)
\qbezier(3,4.5)(3,6)(4.5,6)
\qbezier(4.5,3)(6,3)(6,4.5)
\qbezier(4.5,6)(6,6)(6,4.5)
\put(4,4){\circle*{0.2}}
}

\put(1.1,2.9){
\qbezier(3,4.5)(3,3)(4.5,3)
\qbezier(3,4.5)(3,6)(4.5,6)
\qbezier(4.5,3)(6,3)(6,4.5)
\qbezier(4.5,6)(6,6)(6,4.5)
\put(4,4){\circle*{0.2}}
}

\put(5.5,4.5){
\qbezier(3,4.5)(3,3)(4.5,3)
\qbezier(3,4.5)(3,6)(4.5,6)
\qbezier(4.5,3)(6,3)(6,4.5)
\qbezier(4.5,6)(6,6)(6,4.5)
\put(4,4){\circle*{0.2}}
}

\put(6.5,-1.5){
\qbezier(3,4.5)(3,3)(4.5,3)
\qbezier(3,4.5)(3,6)(4.5,6)
\qbezier(4.5,3)(6,3)(6,4.5)
\qbezier(4.5,6)(6,6)(6,4.5)
\put(4,4){\circle*{0.2}}
}

\put(13,6){\makebox(0,0)[b]{$\cdots$}}
\put(5.8,4.5){\circle{0.2}}
\put(6.8,4.3){\makebox(0,0)[b]{$\alpha_2$}}
}

\put(22,5){
\put(14,8){\makebox(0,0)[b]{$\beta_2F_1$}}
\qbezier(0,4.5)(0,11)(10.5,11)
\qbezier(0,4.5)(0,0)(7.5,0)
\qbezier(7.5,0)(16,0)(16,7.5)
\qbezier(10.5,11)(16,11)(16,7.5)

\put(-1,-0.2){
\qbezier(3,4.5)(3,3)(4.5,3)
\qbezier(3,4.5)(3,6)(4.5,6)
\qbezier(4.5,3)(6,3)(6,4.5)
\qbezier(4.5,6)(6,6)(6,4.5)
\put(4,4){\circle*{0.2}}
}

\put(2.5,-2){
\qbezier(3,4.5)(3,3)(4.5,3)
\qbezier(3,4.5)(3,6)(4.5,6)
\qbezier(4.5,3)(6,3)(6,4.5)
\qbezier(4.5,6)(6,6)(6,4.5)
\put(4,4){\circle*{0.2}}
\put(4.9,3.8){\makebox(0,0)[b]{$b_2$}}
}

\put(4.2,1.3){
\qbezier(3,4.5)(3,3)(4.5,3)
\qbezier(3,4.5)(3,6)(4.5,6)
\qbezier(4.5,3)(6,3)(6,4.5)
\qbezier(4.5,6)(6,6)(6,4.5)
\put(4,4){\circle*{0.2}}
}

\put(1.1,2.9){
\qbezier(3,4.5)(3,3)(4.5,3)
\qbezier(3,4.5)(3,6)(4.5,6)
\qbezier(4.5,3)(6,3)(6,4.5)
\qbezier(4.5,6)(6,6)(6,4.5)
\put(4,4){\circle*{0.2}}
}

\put(5.5,4.5){
\qbezier(3,4.5)(3,3)(4.5,3)
\qbezier(3,4.5)(3,6)(4.5,6)
\qbezier(4.5,3)(6,3)(6,4.5)
\qbezier(4.5,6)(6,6)(6,4.5)
\put(4,4){\circle*{0.2}}
}

\put(6.5,-1.5){
\qbezier(3,4.5)(3,3)(4.5,3)
\qbezier(3,4.5)(3,6)(4.5,6)
\qbezier(4.5,3)(6,3)(6,4.5)
\qbezier(4.5,6)(6,6)(6,4.5)
\put(4,4){\circle*{0.2}}
}

\put(13,6){\makebox(0,0)[b]{$\cdots$}}
\put(5.8,4.5){\circle{0.2}}
\put(6.8,4.3){\makebox(0,0)[b]{$\beta_2$}}
}

\put(4,15){
\put(14,8){\makebox(0,0)[b]{$\gamma_2F_1$}}
\qbezier(0,4.5)(0,11)(10.5,11)
\qbezier(0,4.5)(0,0)(7.5,0)
\qbezier(7.5,0)(16,0)(16,7.5)
\qbezier(10.5,11)(16,11)(16,7.5)

\put(-1,-0.2){
\qbezier(3,4.5)(3,3)(4.5,3)
\qbezier(3,4.5)(3,6)(4.5,6)
\qbezier(4.5,3)(6,3)(6,4.5)
\qbezier(4.5,6)(6,6)(6,4.5)
\put(4,4){\circle*{0.2}}
}

\put(2.5,-2){
\qbezier(3,4.5)(3,3)(4.5,3)
\qbezier(3,4.5)(3,6)(4.5,6)
\qbezier(4.5,3)(6,3)(6,4.5)
\qbezier(4.5,6)(6,6)(6,4.5)
\put(4,4){\circle*{0.2}}
}

\put(4.2,1.3){
\qbezier(3,4.5)(3,3)(4.5,3)
\qbezier(3,4.5)(3,6)(4.5,6)
\qbezier(4.5,3)(6,3)(6,4.5)
\qbezier(4.5,6)(6,6)(6,4.5)
\put(4,4){\circle*{0.2}}
}

\put(1.1,2.9){
\qbezier(3,4.5)(3,3)(4.5,3)
\qbezier(3,4.5)(3,6)(4.5,6)
\qbezier(4.5,3)(6,3)(6,4.5)
\qbezier(4.5,6)(6,6)(6,4.5)
\put(4,4){\circle*{0.2}}
}

\put(5.5,4.5){
\qbezier(3,4.5)(3,3)(4.5,3)
\qbezier(3,4.5)(3,6)(4.5,6)
\qbezier(4.5,3)(6,3)(6,4.5)
\qbezier(4.5,6)(6,6)(6,4.5)
\put(4,4){\circle*{0.2}}
}

\put(6.5,-1.5){
\qbezier(3,4.5)(3,3)(4.5,3)
\qbezier(3,4.5)(3,6)(4.5,6)
\qbezier(4.5,3)(6,3)(6,4.5)
\qbezier(4.5,6)(6,6)(6,4.5)
\put(4,4){\circle*{0.2}}
}

\put(13,6){\makebox(0,0)[b]{$\cdots$}}
\put(5.8,4.5){\circle{0.2}}
\put(6.8,4.3){\makebox(0,0)[b]{$\gamma_2$}}
}

\put(22,17){
\put(14,8){\makebox(0,0)[b]{$\lambda F_1$}}
\qbezier(0,4.5)(0,11)(10.5,11)
\qbezier(0,4.5)(0,0)(7.5,0)
\qbezier(7.5,0)(16,0)(16,7.5)
\qbezier(10.5,11)(16,11)(16,7.5)

\put(-1,-0.2){
\qbezier(3,4.5)(3,3)(4.5,3)
\qbezier(3,4.5)(3,6)(4.5,6)
\qbezier(4.5,3)(6,3)(6,4.5)
\qbezier(4.5,6)(6,6)(6,4.5)
\put(4,4){\circle*{0.2}}
}

\put(2.5,-2){
\qbezier(3,4.5)(3,3)(4.5,3)
\qbezier(3,4.5)(3,6)(4.5,6)
\qbezier(4.5,3)(6,3)(6,4.5)
\qbezier(4.5,6)(6,6)(6,4.5)
\put(4,4){\circle*{0.2}}
}

\put(4.2,1.3){
\qbezier(3,4.5)(3,3)(4.5,3)
\qbezier(3,4.5)(3,6)(4.5,6)
\qbezier(4.5,3)(6,3)(6,4.5)
\qbezier(4.5,6)(6,6)(6,4.5)
\put(4,4){\circle*{0.2}}
}

\put(1.1,2.9){
\qbezier(3,4.5)(3,3)(4.5,3)
\qbezier(3,4.5)(3,6)(4.5,6)
\qbezier(4.5,3)(6,3)(6,4.5)
\qbezier(4.5,6)(6,6)(6,4.5)
\put(4,4){\circle*{0.2}}
}

\put(5.5,4.5){
\qbezier(3,4.5)(3,3)(4.5,3)
\qbezier(3,4.5)(3,6)(4.5,6)
\qbezier(4.5,3)(6,3)(6,4.5)
\qbezier(4.5,6)(6,6)(6,4.5)
\put(4,4){\circle*{0.2}}
}

\put(6.5,-1.5){
\qbezier(3,4.5)(3,3)(4.5,3)
\qbezier(3,4.5)(3,6)(4.5,6)
\qbezier(4.5,3)(6,3)(6,4.5)
\qbezier(4.5,6)(6,6)(6,4.5)
\put(4,4){\circle*{0.2}}
}

\put(13,6){\makebox(0,0)[b]{$\cdots$}}
\put(5.8,4.5){\circle{0.2}}
\put(6.8,4.3){\makebox(0,0)[b]{$\lambda$}}
}

\put(8,27){
\put(14,8){\makebox(0,0)[b]{$\lambda F_1$}}
\qbezier(0,4.5)(0,11)(10.5,11)
\qbezier(0,4.5)(0,0)(7.5,0)
\qbezier(7.5,0)(16,0)(16,7.5)
\qbezier(10.5,11)(16,11)(16,7.5)

\put(-1,-0.2){
\qbezier(3,4.5)(3,3)(4.5,3)
\qbezier(3,4.5)(3,6)(4.5,6)
\qbezier(4.5,3)(6,3)(6,4.5)
\qbezier(4.5,6)(6,6)(6,4.5)
\put(4,4){\circle*{0.2}}
}

\put(2.5,-2){
\qbezier(3,4.5)(3,3)(4.5,3)
\qbezier(3,4.5)(3,6)(4.5,6)
\qbezier(4.5,3)(6,3)(6,4.5)
\qbezier(4.5,6)(6,6)(6,4.5)
\put(4,4){\circle*{0.2}}
}

\put(4.2,1.3){
\qbezier(3,4.5)(3,3)(4.5,3)
\qbezier(3,4.5)(3,6)(4.5,6)
\qbezier(4.5,3)(6,3)(6,4.5)
\qbezier(4.5,6)(6,6)(6,4.5)
\put(4,4){\circle*{0.2}}
}

\put(1.1,2.9){
\qbezier(3,4.5)(3,3)(4.5,3)
\qbezier(3,4.5)(3,6)(4.5,6)
\qbezier(4.5,3)(6,3)(6,4.5)
\qbezier(4.5,6)(6,6)(6,4.5)
\put(4,4){\circle*{0.2}}
}

\put(5.5,4.5){
\qbezier(3,4.5)(3,3)(4.5,3)
\qbezier(3,4.5)(3,6)(4.5,6)
\qbezier(4.5,3)(6,3)(6,4.5)
\qbezier(4.5,6)(6,6)(6,4.5)
\put(4,4){\circle*{0.2}}
}

\put(6.5,-1.5){
\qbezier(3,4.5)(3,3)(4.5,3)
\qbezier(3,4.5)(3,6)(4.5,6)
\qbezier(4.5,3)(6,3)(6,4.5)
\qbezier(4.5,6)(6,6)(6,4.5)
\put(4,4){\circle*{0.2}}
}

\put(13,6){\makebox(0,0)[b]{$\cdots$}}
\put(5.8,4.5){\circle{0.2}}
\put(6.8,4.3){\makebox(0,0)[b]{$\lambda$}}
}

\put(28,32){\makebox(0,0)[b]{$\cdots\cdots$}}
\put(21,15){\makebox(0,0)[b]{$\times$}}
\put(22,14){\makebox(0,0)[b]{$1_G$}}

\put(16,0){
\qbezier(3,4.5)(3,3)(4.5,3)
\qbezier(3,4.5)(3,6)(4.5,6)
\qbezier(4.5,3)(6,3)(6,4.5)
\qbezier(4.5,6)(6,6)(6,4.5)
\put(4,4){\circle*{0.2}}
}

\put(20,-1.3){
\qbezier(3,4.5)(3,3)(4.5,3)
\qbezier(3,4.5)(3,6)(4.5,6)
\qbezier(4.5,3)(6,3)(6,4.5)
\qbezier(4.5,6)(6,6)(6,4.5)
\put(4,4){\circle*{0.2}}
}

\put(2,22){
\qbezier(3,4.5)(3,3)(4.5,3)
\qbezier(3,4.5)(3,6)(4.5,6)
\qbezier(4.5,3)(6,3)(6,4.5)
\qbezier(4.5,6)(6,6)(6,4.5)
\put(4,4){\circle*{0.2}}
}

\put(15,37){
\qbezier(3,4.5)(3,3)(4.5,3)
\qbezier(3,4.5)(3,6)(4.5,6)
\qbezier(4.5,3)(6,3)(6,4.5)
\qbezier(4.5,6)(6,6)(6,4.5)
\put(4,4){\circle*{0.2}}
}

\put(10,36){
\qbezier(3,4.5)(3,3)(4.5,3)
\qbezier(3,4.5)(3,6)(4.5,6)
\qbezier(4.5,3)(6,3)(6,4.5)
\qbezier(4.5,6)(6,6)(6,4.5)
\put(4,4){\circle*{0.2}}
}

\end{picture}
\caption{\label{fig:F2}The composition of $F_2$. In the figure $\lambda$ denotes a generic element of $\Lambda_2$. All the circled points form the set $D^2_1$. All the solid points form the set $D^2_0$. The position of $1_G$ in $F_2$ can be arbitrary.
For details of each translates of $F_1$ see Figure~\ref{fig:F1}.}
\end{center}
\end{figure}

\begin{remark}\label{remark513} \index{$D^n_k$} \index{$\Lambda_n$} \index{$\alpha_n$} \index{$\beta_n$} \index{$\gamma_n$} \index{$a_n$} \index{$b_n$} \index{blueprint!$D^n_k$} \index{blueprint!$\Lambda_n$} \index{blueprint!$\alpha_n$} \index{blueprint!$\beta_n$} \index{blueprint!$\gamma_n$} \index{blueprint!$a_n$} \index{blueprint!$b_n$} \index{pre-blueprint!$D^n_k$} \index{pre-blueprint!$\Lambda_n$} \index{pre-blueprint!$\alpha_n$} \index{pre-blueprint!$\beta_n$} \index{pre-blueprint!$\gamma_n$} \index{pre-blueprint!$a_n$} \index{pre-blueprint!$b_n$}
In the context of Borel equivalence relations, properties (i) and (ii) of the above definition are commonly referred to as the \emph{marker property}, however this term will not be used in this paper. Recall that if $(\Delta, F)$ is a regular marker structure, then necessarily $1_G \in F$. So any pre-blueprint $(\Delta_n, F_n)_{n \in \N}$ must have $1_G \in F_n$ for every $n \in \N$. The set $\gamma^{-1} (\Delta_k \cap \gamma F_n)$ appearing in (iv) will be denoted by $D^n_k$ (as the set only depends on $n$ and $k$). Thus, $\gamma D^n_k = \Delta_k \cap \gamma F_n$ for all $\gamma \in \Delta_n$. Clause (iv) of the above definition may seem mysterious at first, but all it says is that $\Delta_k$ meets each $\Delta_n$-translate of $F_n$ in the same manner (the intersection being a left translate of $D^n_k$). With this thought, note that the growth property is equivalent to $|D^n_{n-1}| \geq 3$. For our purposes we will have to distinguish three distinct members from $D^n_{n-1}$. Thus, whenever discussing a pre-blueprint $(\Delta_n, F_n)_{n \in \N}$ the symbols $\alpha_n, \beta_n, \gamma_n$ will be reserved to denote three distinct elements of $D^n_{n-1}$ (for each $n > 0$). The following symbols will also have a reserved meaning:
$$\Lambda_n = D^n_{n-1} - \{\alpha_n, \beta_n, \gamma_n\};$$
$$a_0 = b_0 = 1_G;$$
$$a_n = \alpha_n \alpha_{n-1} \cdots \alpha_1 \ (\text{for } n \geq 1);$$
$$b_n = \beta_n \beta_{n-1} \cdots \beta_1 \ (\text{for } n \geq 1).$$
See Figures \ref{fig:F1} and \ref{fig:F2} for an illustration of what blueprints generally look like.
\end{remark}

Pre-blueprints are not difficult to construct. The notation involved in the construction may be cumbersome, but the conceptual idea is relatively simple. Pre-blueprints can be constructed one step at a time, with each $F_n$ being a union of translates of $F_k$'s for $k < n$. All of the blueprints we construct will have this property, namely $F_n = \bigcup_{k < n} D^n_k F_k$. This simple construction of pre-blueprints will be presented in detail in the third section.

The fact that every countably infinite group admits a blueprint is quite nontrivial and is postponed to the third section of this chapter. The reason for postponing this is two-fold. First, the construction of a blueprint is rather technical and the reader may value the proof more after having seen why blueprints are important. Second, our construction will show that a very strong type of blueprint exists, in particular one that is maximally disjoint, centered, and directed. The constructed blueprint will have many nice properties and we don't want the reader to have an oversimplified view of blueprints when going through proofs involving arbitrary blueprints.

The purpose of this section is to bring to light some of the useful properties of pre-blueprints. Many of these properties will be vital in the next section. The following lemma consists of direct consequences of the definition of pre-blueprints. These facts will be used with high frequency throughout the rest of the paper.

\begin{lem} \label{BP LIST}
Let $G$ be a countably infinite group, and let $(\Delta_n, F_n)_{n \in \N}$ be a pre-blueprint. Then
\begin{enumerate}
\item[\rm (i)] $\Delta_n D^n_k \subseteq \Delta_k$, for all $k < n$;
\item[\rm (ii)] $D^n_k F_k \subseteq F_n$ for all $k < n$;
\item[\rm (iii)] $\lambda_1 F_k \cap \lambda_2 F_k = \varnothing$ for all $k < n$ and distinct $\lambda_1, \lambda_2 \in D^n_k$;
\item[\rm (iv)] $D^n_m D^m_k \subseteq D^n_k$ for all $k < m < n$;
\item[\rm (v)] $D^n_k \neq \varnothing$ for all $k < n$;
\item[\rm (vi)] $\Delta_k f \cap \gamma F_n = \gamma D^n_k f = (\Delta_k \cap \gamma F_n) f$, for all $k < n$, $\gamma \in \Delta_n$, and $f \in F_k$;
\item[\rm (vii)] $\gamma f \in \Delta_k B \Longleftrightarrow \sigma f \in \Delta_k B$, for all $k < n$, $\gamma, \sigma \in \Delta_n$, $f \in F_n$, and $B \subseteq F_k$;
\item[\rm (viii)] both $(\Delta_n a_n)_{n \in \N}$ and $(\Delta_n b_n)_{n \in \N}$ are decreasing sequences;
\item[\rm (ix)] $a_n, b_n \in F_n$ for all $n \in \N$;
\item[\rm (x)] $a_n \neq b_n$ for all $n \geq 1$;
\item[\rm (xi)] $\Delta_n a_n \cap \Delta_k b_k = \varnothing$ for all $n, k > 0$;
\item[\rm (xii)] for $n > k \in \N$
$$\Delta_n D^n_{n-1} a_{n-1} \cap \Delta_k F_k \subseteq \Delta_k a_k$$
and
$$\Delta_n D^n_{n-1} b_{n-1} \cap \Delta_k F_k \subseteq \Delta_k b_k;$$
\item[\rm (xiii)] for $n > k \in \N$
$$\Delta_n D^n_{n-1} a_{n-1} \cap \Delta_k D^k_{k-1} a_{k-1}$$
$$\subseteq \Delta_k \alpha_k a_{k-1} = \Delta_k a_k$$
and
$$\Delta_n D^n_{n-1} b_{n-1} \cap \Delta_k D^k_{k-1} b_{k-1}$$
$$\subseteq \Delta_k \beta_k b_{k-1} = \Delta_k b_k.$$
\end{enumerate}
\end{lem}

\begin{proof}
(i). Let $\gamma \in \Delta_n$ and $\lambda \in D^n_k$. Then
$$\gamma \lambda \in \gamma D^n_k = \Delta_k \cap \gamma F_n \subseteq \Delta_k.$$

(ii). Pick any $\gamma \in \Delta_n$ and $\lambda \in D^n_k$. Then $\gamma \lambda \in \gamma D^n_k = \Delta_k \cap \gamma F_n$. So $\gamma \lambda \in \gamma F_n$ and $\gamma \lambda \in \Delta_k$. By definition, $1_G \in F_k$, so that $\gamma \lambda F_k \cap \gamma F_n \neq \varnothing$. By the coherent property of pre-blueprints, this gives $\gamma \lambda F_k \subseteq \gamma F_n$. Now cancel $\gamma$.

(iii). Pick any $\gamma \in \Delta_n$. Then $\gamma \lambda_1$ and $\gamma \lambda_2$ are distinct elements of $\Delta_k$ by (i). The claim now follows from the disjoint property of pre-blueprints.

(iv). Pick $\psi \in D^n_m$ and $\lambda \in D^m_k$. Let $\gamma \in \Delta_n$ be arbitrary. By (i) $\gamma \psi \in \Delta_m$ and hence $\gamma \psi \lambda \in \Delta_k$. As $1_G \in F_k$, by (ii) we have
$$\gamma \psi \lambda \in \gamma \psi \lambda F_k \subseteq \gamma \psi F_m \subseteq \gamma F_n.$$
Thus, $\gamma \psi \lambda \in \Delta_k \cap \gamma F_n = \gamma D^n_k$.

(v). By the growth property of pre-blueprints, $D^n_{n-1} \neq \varnothing$. Now suppose $D^n_{k+1} \neq \varnothing$. As $D^{k+1}_k \neq \varnothing$, we have $\varnothing \neq D^n_{k+1} D^{k+1}_k \subseteq D^n_k$ by (iv).

(vi). The second equality is clear from the definition of $D^n_k$. We verify the first equality. Let $\psi \in \Delta_k$ be such that $\psi f \in \gamma F_n$. Then $\psi F_k \cap \gamma F_n \neq \varnothing$, so by the coherent property of pre-blueprints $\psi F_k \subseteq \gamma F_n$. By definition, $1_G \in F_k$, so $\psi \in \gamma F_n$. It follows $\psi \in \Delta_k \cap \gamma F_n = \gamma D^n_k$ and hence $\psi f \in \gamma D^n_k f$. On the other hand, $\gamma D^n_k \subseteq \Delta_k$ by (i). Hence $\gamma D^n_k f \subseteq \Delta_k f$. Also, by (ii) $D^n_k f \subseteq F_n$. So $\gamma D^n_k f \subseteq \gamma F_n$. Thus we have $\gamma D^n_k f \subseteq \Delta_k f \cap \gamma F_n$.

(vii). By (vi) we have
$$\gamma f \in \Delta_k B \Longleftrightarrow \gamma f \in \Delta_k B \cap \gamma F_n \Longleftrightarrow \gamma f \in \gamma D^n_k B \Longleftrightarrow f \in D^n_k B $$
$$\Longleftrightarrow \sigma f \in \sigma D^n_k B \Longleftrightarrow \sigma f \in \Delta_k B \cap \sigma F_n \Longleftrightarrow \sigma f \in \Delta_k B.$$

(viii). By (i), $\Delta_n a_n = \Delta_n \alpha_n a_{n-1} \subseteq \Delta_{n-1} a_{n-1}$. The same argument applies to $\Delta_n b_n$.

(ix). By definition $a_0 = 1_G \in F_0$. If we assume $a_{n-1} \in F_{n-1}$, then by (ii)
$$a_n = \alpha_n a_{n-1} \in \alpha_n F_{n-1} \subseteq F_n.$$
By induction, and by a similar argument, we have $a_n, b_n \in F_n$ for all $n \in \N$.

(x). From (ix) we have that $a_n = \alpha_n a_{n-1} \in \alpha_n F_{n-1}$. Similarly, $b_n \in \beta_n F_{n-1}$. The claim then follows from (iii) since $\alpha_n \neq \beta_n$.

(xi). Suppose $0 < k \leq n$. Since $a_k \neq b_k \in F_k$, and since the $\Delta_k$-translates of $F_k$ are disjoint, from (viii) we have
$$\Delta_n a_n \cap \Delta_k b_k \subseteq \Delta_k a_k \cap \Delta_k b_k = \varnothing.$$
The case $0 < n \leq k$ is identical.

(xii) and (xiii). $\Delta_n D^n_{n-1} \subseteq \Delta_{n-1}$ by (i). So by (viii)
$$\Delta_n D^n_{n-1} a_{n-1} \subseteq \Delta_{n-1} a_{n-1} \subseteq \Delta_k a_k.$$
This gives us (xii) and part of (xiii). The equality at the end of (xiii) follows from the definition of $a_k$. The argument for $(b_n)_{n \in \N}$ is identical.
\end{proof}

Hopefully the proof of the previous lemma helps demonstrate to the reader the important role of the coherent and uniform properties of pre-blueprints.

The next lemma consists of properties of stronger types of pre-blueprints.

\begin{lem} \label{STRONG BP LIST}
Let $G$ be a countably infinite group and let $(\Delta_n, F_n)_{n \in \N}$ be a pre-blueprint.
\begin{enumerate}
\item[\rm (i)] If the pre-blueprint is centered, then $(F_n)_{n \in \N}$ is an increasing sequence and $(\Delta_n)_{n \in \N}$ is a decreasing sequence.
\item[\rm (ii)] If the pre-blueprint is maximally disjoint, then it is a blueprint and for all $n\in\N$,  $\Delta_n F_n F_n^{-1} = G$.
\item[\rm (iii)] If the pre-blueprint is directed, then for any $r(1), r(2), \ldots, r(m)$ and $\psi_1 \in \Delta_{r(1)}, \ldots, \psi_m \in \Delta_{r(m)}$, there is $n \in \N$ and $\gamma \in \Delta_n$ so that for every $1 \leq i \leq m$ $\psi_i F_{r(i)} \subseteq \gamma F_n$.
\item[\rm (iv)] If the pre-blueprint is directed and centered, then for any $r(1), r(2), \ldots$, $r(m)$ and $\psi_1 \in \Delta_{r(1)}, \ldots, \psi_m \in \Delta_{r(m)}$, there is $n \in \N$ so that for every $1 \leq i \leq m$ $\psi_i F_{r(i)} \subseteq F_n$.
\item[\rm (v)] If $n > k \geq t$, $\sigma \in \Delta_n$, $A \subseteq G$ is finite, $\Delta_t \cap A F_k F_t F_t F_t^{-1} \subseteq \sigma D^n_t$, and the $\Delta_t$-translates of $F_t$ are maximally disjoint, then for all $\gamma \in \Delta_n$
$$\Delta_k \cap \gamma \sigma^{-1} A = \gamma \sigma^{-1} (\Delta_k \cap A) \subseteq \gamma D^n_k;$$
\item[\rm (vi)] If the pre-blueprint is a directed blueprint and the $\Delta_0$-translates of $F_0$ are maximally disjoint, then it is minimal in the following sense: for every finite $A \subseteq G$ and $N \in \N$ there is a finite set $T \subseteq G$ so that for any $g \in G$
$$\exists t \in T \ \forall 0 \leq k \leq N \ gt(\Delta_k \cap A) = \Delta_k \cap gtA.$$
\item[\rm (vii)] If the pre-blueprint is directed, then
$$\left|\bigcap_{n \in \N} \Delta_n \right|, \ \left|\bigcap_{n \in \N} \Delta_n a_n \right|, \ \left|\bigcap_{n \in \N} \Delta_n b_n \right| \leq 1.$$
\item[\rm (viii)] If the pre-blueprint is directed, centered, and for infinitely many $n$ and infinitely many $k$ $\alpha_n \neq 1_G \neq \beta_k$, then
$$\bigcap_{n \in \N} \Delta_n a_n = \bigcap_{n \in \N} \Delta_n b_n = \varnothing.$$
\end{enumerate}
\end{lem}

\begin{proof}
(i). By definition, $1_G \in F_n$ for each $n \in \N$. Thus $F_n \cap F_{n+1} \neq \varnothing$. If the pre-blueprint is centered, then $1_G \in \Delta_n \cap \Delta_{n+1}$, so it follows from the coherent property of pre-blueprints that $F_n \subseteq F_{n+1}$. By (i) of Lemma \ref{BP LIST}, $\Delta_{n+1} D^{n+1}_n \subseteq \Delta_n$. However, $1_G \in \Delta_n \cap 1_G F_{n+1} = 1_G D^{n+1}_n$, so $\Delta_{n+1} \subseteq \Delta_n$.

(ii). Suppose the pre-blueprint is maximally disjoint, and let $g \in G$. Then the $\Delta_n$-translates of $F_n$ are maximally disjoint within $G$, so $g F_n \cap \Delta_n F_n \neq \varnothing$. Hence there is $f_1, f_2 \in F_n$ and $\gamma \in \Delta_n$ with $g f_1 = \gamma f_2$. It follows $g = \gamma f_2 f_1^{-1} \in \Delta_n F_n F_n^{-1}$. So the dense property is satisfied and the pre-blueprint is a blueprint.

(iii). It suffices to prove the claim for the maximal elements with respect to inclusion among $\{\psi_i F_{r(i)} \: 1 \leq i \leq m\}$. By the coherent property, distinct maximal members of this collection are disjoint. So without loss of generality we may assume $\psi_i F_{r(i)} \cap \psi_j F_{r(j)} = \varnothing$ for $i \neq j$. Also, without loss of generality we may assume $r(1) \leq r(2) \leq \cdots \leq r(m)$. For each $i > 1$ pick $\lambda_i \in D^{r(i)}_{r(1)}$. Then $\psi_i \lambda_i \in \Delta_{r(1)}$ for each $i > 1$. For each $i > 1$ pick $n(i)$ and $\sigma_i \in \Delta_{n(i)}$ with
$$\psi_1 F_{r(1)} \cup \psi_i \lambda_i F_{r(1)} \subseteq \sigma_i F_{n(i)}.$$
So for $i, j > 1$ we have
$$\psi_1 F_{r(1)} \subseteq \sigma_i F_{n(i)} \cap \sigma_j F_{n(j)}.$$
Thus, by the coherent property of pre-blueprints, it must be that one of $\sigma_i F_{n(i)}$ and $\sigma_j F_{n(j)}$ contains the other. If $\sigma F_n$ is the largest member of the $\sigma_i F_{n(i)}$'s with respect to containment, then we have
$$\psi_1 F_{r(1)} \cup \psi_2 \lambda_2 F_{r(1)} \cup \cdots \cup \psi_m \lambda_m F_{r(1)} \subseteq \sigma F_n.$$
As $\psi_i \lambda_i F_{r(1)} \subseteq \psi_i F_{r(i)}$ for each $i > 1$, we have that each $\psi_i F_{r(i)}$ meets $\sigma F_n$ nontrivially. Since the $\psi_i F_{r(i)}$'s are pairwise disjoint, none of them can contain $\sigma F_n$. Thus by the coherent property $\sigma F_n$ must contain each $\psi_i F_{r(i)}$.

(iv). By (iii) there is $n \in \N$ and $\gamma \in \Delta_n$ with $\psi_i F_{r(i)} \subseteq \gamma F_n$. So it will be enough to show $\gamma F_n \subseteq F_m$ for some $m \in \N$. As $1_G \in \Delta_n$, the pre-blueprint being directed implies there is $m \in \N$ and $\sigma \in \Delta_m$ with $\gamma F_n \cup F_n \subseteq \sigma F_m$. By (i), this gives $F_m \cap \sigma F_m \neq \varnothing$. Since $1_G \in \Delta_m$, by the disjoint property of pre-blueprints we must have $\sigma = 1_G$.

(v). Fix $\gamma \in \Delta_n$. We have
$$(\Delta_k \cap A) D^k_t \subseteq (\Delta_k D^k_t) \cap (A D^k_t) \subseteq \Delta_t \cap A F_k F_t F_t F_t^{-1}$$
$$\subseteq \sigma D^n_t \subseteq \sigma F_n.$$
So if $\psi \in \Delta_k \cap A$ then $\psi F_k \cap \sigma F_n \neq \varnothing$ and by the coherent property of pre-blueprints it follows that $\psi \in \sigma D^n_k$. So $\Delta_k \cap A \subseteq \sigma D^n_k$. Therefore
$$\gamma \sigma^{-1} (\Delta_k \cap A) \subseteq \gamma \sigma^{-1} \sigma D^n_k \subseteq \gamma D^n_k \subseteq \Delta_k.$$
Also $\gamma \sigma^{-1} (\Delta_k \cap A) \subseteq \gamma \sigma^{-1} A$. Thus
$$\gamma \sigma^{-1} (\Delta_k \cap A) \subseteq \Delta_k \cap (\gamma \sigma^{-1} A).$$

To show the reverse inclusion, pick $\lambda \in \Delta_k \cap (\gamma \sigma^{-1} A)$. Fix any $\delta \in D^k_t$. Then $\sigma \gamma^{-1} \lambda \delta \in A F_t$ so $\sigma \gamma^{-1} \lambda \delta F_t \subseteq A F_t F_t$. Notice that the $\Delta_t \cap A F_k F_t F_t F_t^{-1}$-translates of $F_t$ are maximally disjoint within $A F_t F_t$ (though not necessarily contained in $A F_t F_t$; see Definition \ref{DEFN MAX DIS}). So there is $\psi \in \Delta_t \cap A F_k F_t F_t F_t^{-1}$ with
$$\psi F_t \cap \sigma \gamma^{-1} \lambda \delta F_t \neq \varnothing.$$
However, $\psi \in \sigma D^n_t$, so $\gamma \sigma^{-1} \psi \in \Delta_t$, $\lambda \delta \in \Delta_t$, and
$$\gamma \sigma^{-1} \psi F_t \cap \lambda \delta F_t \neq \varnothing.$$
Therefore $\lambda \delta = \gamma \sigma^{-1} \psi$. Since $\psi \in \sigma D^n_t \subseteq \sigma F_n$, we have $\gamma \sigma^{-1} \psi \in \gamma F_n$. Thus $\lambda \delta \in \gamma F_n$ so $\lambda F_k \cap \gamma F_n \neq \varnothing$. By the coherency property of blueprints, we have that $\lambda F_k \subseteq \gamma F_n$ and $\lambda \in \gamma D^n_k$. It follows that $\sigma \gamma^{-1} \lambda \in \sigma D^n_k \subseteq \Delta_k$. Thus
$$\sigma \gamma^{-1} \lambda \in \Delta_k \cap A$$
and therefore
$$\lambda \in \gamma \sigma^{-1} (\Delta_k \cap A).$$

(vi). Let $A \subseteq G$ be finite and let $N \in \N$. For each $0 \leq k \leq N$ let $C_k = \Delta_0 \cap A F_k F_0 F_0 F_0^{-1}$. By (iii) there is $n \in \N$ and $\sigma \in \Delta_n$ with $C_k F_0 \subseteq \sigma F_n$ for every $0 \leq k \leq N$. In particular, $C_k \subseteq \sigma D^n_0$. Since we are assuming $(\Delta_n F_n)_{n \in \N}$ is a blueprint, there is a finite $B \subseteq G$ with $\Delta_n B = G$. Set $T = B^{-1} \sigma^{-1}$ and let $g \in G$ be arbitrary. Since $\Delta_n B = G$, there is $b \in B^{-1}$ with $gb = \gamma \in \Delta_n$. We will show that the stated condition is satisfied for $t = b \sigma^{-1} \in T$. So $g t = \gamma \sigma^{-1}$. This follows from (v): for $0 \leq k \leq N$ we have
$$\gamma \sigma^{-1} (\Delta_k \cap A) = \Delta_k \cap (\gamma \sigma^{-1} A).$$
We conclude that the blueprint satisfies the stated minimal condition.

(vii). Let $(f_n)_{n\in\N}$ be a sequence such that $f_n\in F_n$ for all $n\in\N$. We show that $\left|\bigcap_{n \in \N} \Delta_n f_n\right|\leq 1$.
Suppose $g, h \in \bigcap_{n \in \N} \Delta_n f_n$. Then $g, h \in \Delta_0 f_0 \subseteq \Delta_0 F_0$. If our pre-blueprint is directed, then there is $n > 0$ and $\gamma \in \Delta_n$ with $g, h \in \gamma F_n$. However, $g, h \in \Delta_n f_n$, so there are $\sigma_1, \sigma_2 \in \Delta_n$ with $g = \sigma_1 f_n$ and $h = \sigma_2 f_n$. By conclusion (ix) of Lemma \ref{BP LIST} $\sigma_1 F_n \cap \gamma F_n \neq \varnothing$. By the disjoint property of pre-blueprints, we must have $\sigma_1 = \gamma$ and similarly $\sigma_2 = \gamma$. Thus $\sigma_1 = \sigma_2$ and it follows $g = h$.

(viii). Assume that our pre-blueprint is centered and directed. Suppose $g \in \bigcap_{n \in \N} \Delta_n a_n$. We will show $\alpha_n = 1_G$ for all but finitely many $n \in \N$. We have $g \in \Delta_0 a_0 = \Delta_0$ and $1_G \in \Delta_0$. Hence there is $n > 0$ and $\gamma \in \Delta_n$ with $g, 1_G \in \gamma F_n$. As $1_G \in \Delta_n$ and $1_G F_n \cap \gamma F_n \supseteq \{1_G\} \neq \varnothing$, we must have $\gamma = 1_G$. Thus $g \in F_n$. By (i), $(F_m)_{m \in \N}$ is an increasing sequence, so $g \in F_m$ for all $m \geq n$. Fix $m \geq n$ and let $\sigma \in \Delta_{m+1}$ be such that $g = \sigma a_{m+1}$. As $a_{m+1} \in F_{m+1}$, we have $\sigma F_{m+1} \cap F_{m+1} \supseteq \{g\} \neq \varnothing$. We then must have $\sigma = 1_G$. So $g = a_{m+1} = \alpha_{m+1} a_m \in \alpha_{m+1} F_m$. Then $g \in F_m \cap \alpha_{m+1} F_m$ so $\alpha_{m+1} = 1_G$. The case of $\cap \Delta_n b_n$ is similar, so this completes the proof.
\end{proof}

Note that (iii) reveals why the word ``directed'' was chosen. Consider the set $\mathcal{C} = \{\gamma F_n \: n \in \N \wedge \gamma \in \Delta_n\}$ with the partial ordering
$$\psi F_k \prec \gamma F_n \Longleftrightarrow k \leq n \wedge \psi F_k \subseteq \gamma F_n.$$
Conclusion (iii) says that if the pre-blueprint $(\Delta_n, F_n)_{n \in \N}$ is directed, then this partially ordered set is a directed set.

The importance of clause (v) will be better appreciated after the next section. This is because the behavior fundamental functions (partial functions on $G$ constructed in the next section) will be highly dependent on the sets $(\Delta_n)_{n \in \N}$. Knowing how subsets $A \subseteq G$ intersect $\Delta_k$ will be very useful. Clause (v) should also be recognized as being closely related to the uniform property of pre-blueprints. If we set $A = \sigma F_n$ and ignore the assumptions of this clause, then we see that the conclusion is precisely the uniform property appearing in the definition of pre-blueprints.

We point out that the minimal property mentioned in (vi) actually does relate to the minimality of a certain dynamical system. Fix $N \in \N$, and define $x: G \rightarrow 2^{N+1}$ so that for $g, h \in G$
$$x(g) = x(h) \Longleftrightarrow (\forall 0 \leq k \leq N \ g \in \Delta_k \Leftrightarrow h \in \Delta_k).$$
Then one can check via Lemma \ref{lem:minimallemma} that $x$ is minimal if and only if the pre-blueprint satisfies the stated minimal condition for $N$.

\section{Fundamental functions} \label{SECT FM}

Now we will get to see how pre-blueprints are used in constructing well behaved partial functions on $G$. As the reader will see, one reason sequences of marker structures are useful is that it endows organization to the group which allows one to work with the group at the small scale at first and step by step work at larger and larger scales tweaking what has been done previously. Of course, different types of sequences of marker structures are needed for different constructions. Pre-blueprints seem to be precisely the type needed in the main construction of this section. We will construct a partial function on $G$, and the important feature of the constructed function is that one will be able to recognize the structure of the pre-blueprint from the behavior of the function alone. In other words, the organization endowed to the group by this sequence of marker structures, the pre-blueprint, will essentially become an intrinsic feature of the partial function. The members of $\Delta_n$ for each $n > 0$ will be identifiable using what we call a membership test.

\begin{definition} \index{membership test}\index{membership test!simple}\index{test region}\index{test functions}\index{membership test!test region}\index{membership test!test functions}
Let $G$ be a group, $c \in 2^{\subseteq G}$, and $\Delta \subseteq G$. We say $c$ \emph{admits a $\Delta$ membership test} if there is a finite $V \subseteq G$ and $S \subseteq 2^V$ so that $\Delta V \subseteq \dom(c)$ and for all $x \in 2^G$ with $x \supseteq c$
$$g \in \Delta \Longleftrightarrow (g^{-1} \cdot x) \res V \in S.$$
The set $V$ is called a \emph{test region} and the elements of $S$ are called \emph{test functions}. If $S$ can be taken to be a singleton, we may say that $c$ admits a \emph{simple} $\Delta$ membership test. In this case, if $S = \{f\}$ then for all $x \in 2^G$ with $x \supseteq c$
$$g \in \Delta \Longleftrightarrow \forall v \in V \ x(gv) = f(v).$$
\end{definition}

Equivalently, $c$ admits a $\Delta$ membership test if there is a finite $V \subseteq G$ satisfying $\Delta V \subseteq \dom(c)$ and with the property that for any $g \not\in \Delta$ and $\gamma \in \Delta$, there is $v \in V$ with $gv \in \dom(c)$ and $c(gv) \neq c(\gamma v)$. This equivalent characterization is not used in this paper and so we do not include a proof.

In the upcoming theorem, we will create a single function with a simple $\Delta_n$ membership test for each $n > 0$. The membership test will be constructed inductively; the membership test of $\Delta_{n+1}$ will be reliant on the membership test for $\Delta_n$. Establishing the base case of the induction seems to be achieved most easily through the use of a locally recognizable function.

\begin{definition} \label{DEFN LR} \index{locally recognizable function}\index{locally recognizable function!trivial}
Let $G$ be a group, let $A \subseteq G$ be finite with $1_G \in A$, and let $R: A \rightarrow 2$. We call $R$ \emph{locally recognizable} if for every $1_G \neq a \in A$ there is $b \in A$ so that $a b \in A$ and $R(ab) \neq R(b)$. $R$ is called \emph{trivial} if $|\{a \in A \: R(a) = R(1_G)\}| = 1$.
\end{definition}

The lemma below gives an equivalent characterization of locally recognizable functions. The property used in the definition above is the easiest to verify, while the property given in the lemma below is the most useful property in terms of applications.

\begin{lem} \label{LEM LR EQUIV}
Let $G$ be a group, let $A \subseteq G$ be finite with $1_G \in A$, and let $R: A \rightarrow 2$. The function $R$ is locally recognizable if and only if for every $x \in 2^G$ with $x \res A = R$
$$\forall a \in A \ (\forall b \in A \ x(ab) = x(b) \Longrightarrow a = 1_G).$$
\end{lem}

\begin{proof}
($\Rightarrow$). Suppose $R$ is locally recognizable. If $1_G \neq a \in A$, then by definition there is $b \in A$ with $a b \in A$ and $R(ab) \neq R(b)$. So if $x \in 2^G$ satisfies $x \res A = R$, then $x(ab) \neq x(b)$.

($\Leftarrow$). Assume that $R$ has the property stated above. Let $a \in A$ and suppose that for every $b \in A$ either $a b \not\in A$ or else $R(ab) = R(b)$. It suffices to show $a = 1_G$. We may define $R':A \cup a A \rightarrow 2$ by requiring $R'$ to extend $R$ and satisfy $R'(ab) = R(b)$ for every $b \in A$. Then $R'$ is well defined. If $x \in 2^G$ is any extension of $R'$, then $x(ab) = x(b)$ for every $b \in A$. Thus, by assumption this implies $a = 1_G$.
\end{proof}

Every group with more than two elements admits a nontrivial locally recognizable function. If $G$ contains a nonidentity $g$ with $g^2 \neq 1_G$, then set $A = \{1_G, g, g^2\}$ and define $R(1_G) = R(g) = 1$ and $R(g^2) = 0$. If every element of $G$ has order two, then pick distinct nonidentity $g, h \in G$, set $A = \{1_G, g, h, gh\}$, and define $R(1_G) = R(g) = R(h) = 1$ and $R(gh) = 0$ (keep in mind $h g = g h$ as $G$ must be abelian). Also note that if $R: A \rightarrow 2$ is locally recognizable and $B \supseteq A$, then $R': B \rightarrow 2$ is locally recognizable, where
$$R'(b) = \begin{cases}
R(b) & \text{if } b \in A \\
1 - R(1_G) & \text{if } b \in B - A
\end{cases}.$$
Thus, nontrivial locally recognizable functions exist on a large multitude of domains. More advanced examples of locally recognizable functions will be presented in the next chapter where we will see that they are useful for more than just creating a membership test.

\begin{definition} \index{$m$-uniform}\index{uniform!$m$-uniform}
Let $(\Delta_n, F_n)_{n \in \N}$ be a pre-blueprint. A set $A \subseteq G$ is said to be \emph{$m$-uniform} with respect to this pre-blueprint if
$$\forall n \geq m \ \forall \gamma, \sigma \in \Delta_n \ \gamma^{-1} ( A \cap \gamma F_n) = \sigma^{-1} (A \cap \sigma F_n).$$
A partial function $c \in 2^{\subseteq G}$ is said to be $m$-\emph{uniform} with respect to $(\Delta_n, F_n)_{n \in \N}$ if each of the three sets $\dom(c)$, $c^{-1}(0)$, and $c^{-1}(1)$ are $m$-uniform with respect to $(\Delta_n, F_n)_{n \in \N}$.
\end{definition}

The uniform property of pre-blueprints asserts that for every $k \in \N$ $\Delta_k$ is $k$-uniform relative to $(\Delta_n, F_n)_{n \in \N}$.

We are now ready for the construction. It may help to recall some of the fixed notation related to pre-blueprints ($\alpha_n$, $\beta_n$, $\gamma_n$, $a_n$, $b_n$, $D^n_k$, $\Lambda_n$) before continuing.

\begin{theorem} \label{FM}
Let $G$ be a countably infinite group, let $(\Delta_n, F_n)_{n \in \N}$ be a pre-blueprint, and let $R: F_0 \rightarrow 2$ be a nontrivial locally recognizable function. Then there exists a function $c \in 2^{\subseteq G}$ with the following properties:
\begin{enumerate}
\item[\rm (i)] $c(\gamma \gamma_1 f) = R(f)$ for all $\gamma \in \Delta_1$ and $f \in F_0$;
\item[\rm (ii)] $c$ admits a simple $\Delta_n$ membership test with test region a subset of $\gamma_n F_{n-1}$ for each $n \geq 1$;
\item[\rm (iii)] $G - \dom(c)$ is the disjoint union $\bigcup_{n \geq 1} \Delta_n \Lambda_n b_{n-1}$;
\item[\rm (iv)] $c(g) = 1 - R(1_G)$ for all $g \in G - \Delta_1 (\gamma_1 F_0 \cup D^1_0)$;
\item[\rm (v)] $(\gamma F_n - \{\gamma b_n\}) \cap \dom(c) = \gamma (F_n - \{b_n\} - \bigcup_{1 \leq k \leq n} D^n_k \Lambda_k b_{k-1})$ for all $n \geq 1$ and $\gamma \in \Delta_n$;
\item[\rm (vi)] $c(\gamma f) = c(\sigma f)$ for all $n \geq 1$, $\gamma, \sigma \in \Delta_n$, and
$$f \in F_n - \{a_n, b_n\} - \bigcup_{1 \leq k \leq n} D^n_k \Lambda_k b_{k-1};$$
\item[\rm (vii)] for all $n \geq 1$ $c \res (G - \Delta_n \{a_n, b_n\})$ is $n$-uniform.
\end{enumerate}
\end{theorem}

\begin{proof}
We wish to construct a sequence of functions $(c_n)_{n \geq 1}$ satisfying for each $n \geq 1$:
\begin{enumerate}
\item[(1)] $\dom(c_n) = G - \Delta_n a_n - \Delta_n b_n - \bigcup_{1 \leq k \leq n} \Delta_k \Lambda_k b_{k-1}$
\item[(2)] $c_{n+1} \supseteq c_n$;
\item[(3)] $c_n$ admits a simple $\Delta_n$ membership test with test region a subset of $\gamma_n F_{n-1}$.
\end{enumerate}


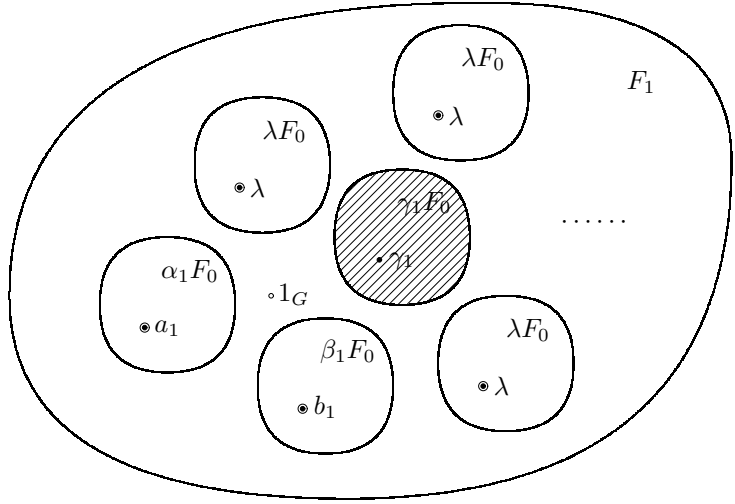
\begin{figure}[ht]
\begin{center}
\setlength{\unitlength}{6mm}
\begin{picture}(30,12)(-2,-.5)

\put(14,9){\makebox(0,0)[b]{$F_1$}}
\qbezier(0,4.5)(0,11)(10.5,11)
\qbezier(0,4.5)(0,0)(7.5,0)
\qbezier(7.5,0)(16,0)(16,7.5)
\qbezier(10.5,11)(16,11)(16,7.5)

\put(-1,-0.2){
\put(5,5){\makebox(0,0)[b]{$\alpha_1F_0$}}
\qbezier(3,4.5)(3,3)(4.5,3)
\qbezier(3,4.5)(3,6)(4.5,6)
\qbezier(4.5,3)(6,3)(6,4.5)
\qbezier(4.5,6)(6,6)(6,4.5)
\put(4,4){\circle*{0.1}}
\put(4,4){\circle{0.2}}
\put(4.5,3.8){\makebox(0,0)[b]{$a_1$}}
}

\put(2.5,-2){
\put(5,5){\makebox(0,0)[b]{$\beta_1 F_0$}}
\qbezier(3,4.5)(3,3)(4.5,3)
\qbezier(3,4.5)(3,6)(4.5,6)
\qbezier(4.5,3)(6,3)(6,4.5)
\qbezier(4.5,6)(6,6)(6,4.5)
\put(4,4){\circle*{0.1}}
\put(4,4){\circle{0.2}}
\put(4.5,3.8){\makebox(0,0)[b]{$b_1$}}
}

\put(4.2,1.3){
\put(5,5){\makebox(0,0)[b]{$\gamma_1 F_0$}}
\qbezier(3,4.5)(3,3)(4.5,3)
\qbezier(3,4.5)(3,6)(4.5,6)
\qbezier(4.5,3)(6,3)(6,4.5)
\qbezier(4.5,6)(6,6)(6,4.5)
\put(4,4){\circle*{0.1}}
\put(4.5,3.8){\makebox(0,0)[b]{$\gamma_1$}}

\put(3.4,3.4){\line(1,1){2.22}}
\put(3.5,3.3){\line(1,1){2.21}}
\put(3.6,3.2){\line(1,1){2.2}}
\put(3.75,3.13){\line(1,1){2.12}}
\put(3.88,3.08){\line(1,1){2.04}}
\put(4.03,3.04){\line(1,1){1.93}}
\put(4.21,3.01){\line(1,1){1.77}}
\put(4.4,3){\line(1,1){1.6}}
\put(4.62,3.02){\line(1,1){1.37}}
\put(4.82,3.02){\line(1,1){1.14}}
\put(5.07,3.08){\line(1,1){0.88}}
\put(5.4,3.2){\line(1,1){0.4}}

\put(3.3,3.5){\line(1,1){2.21}}
\put(3.2,3.6){\line(1,1){2.2}}
\put(3.13,3.75){\line(1,1){2.12}}
\put(3.08,3.88){\line(1,1){2.04}}
\put(3.04,4.03){\line(1,1){1.93}}
\put(3.01,4.21){\line(1,1){1.77}}
\put(3,4.4){\line(1,1){1.6}}
\put(3.02,4.62){\line(1,1){1.37}}
\put(3.02,4.82){\line(1,1){1.14}}
\put(3.08,5.07){\line(1,1){0.88}}
\put(3.2,5.4){\line(1,1){0.4}}
}

\put(1.1,2.9){
\put(5,5){\makebox(0,0)[b]{$\lambda F_0$}}
\qbezier(3,4.5)(3,3)(4.5,3)
\qbezier(3,4.5)(3,6)(4.5,6)
\qbezier(4.5,3)(6,3)(6,4.5)
\qbezier(4.5,6)(6,6)(6,4.5)
\put(4,4){\circle*{0.1}}
\put(4,4){\circle{0.2}}
\put(4.4,3.8){\makebox(0,0)[b]{$\lambda$}}
}

\put(5.5,4.5){
\put(5,5){\makebox(0,0)[b]{$\lambda F_0$}}
\qbezier(3,4.5)(3,3)(4.5,3)
\qbezier(3,4.5)(3,6)(4.5,6)
\qbezier(4.5,3)(6,3)(6,4.5)
\qbezier(4.5,6)(6,6)(6,4.5)
\put(4,4){\circle*{0.1}}
\put(4,4){\circle{0.2}}
\put(4.4,3.8){\makebox(0,0)[b]{$\lambda$}}
}

\put(6.5,-1.5){
\put(5,5){\makebox(0,0)[b]{$\lambda F_0$}}
\qbezier(3,4.5)(3,3)(4.5,3)
\qbezier(3,4.5)(3,6)(4.5,6)
\qbezier(4.5,3)(6,3)(6,4.5)
\qbezier(4.5,6)(6,6)(6,4.5)
\put(4,4){\circle*{0.1}}
\put(4,4){\circle{0.2}}
\put(4.4,3.8){\makebox(0,0)[b]{$\lambda$}}
}

\put(13,6){\makebox(0,0)[b]{$\cdots\cdots$}}
\put(5.8,4.5){\circle{0.1}}
\put(6.3,4.3){\makebox(0,0)[b]{$1_G$}}

\end{picture}
\caption{\label{fig:F3} An illustration of the set $\gamma_1F_0\cup D^1_0$, the shaded area together with all the highlighted (circled-solid) points in the figure. Compare with Figure~\ref{fig:F1}.}
\end{center}
\end{figure}


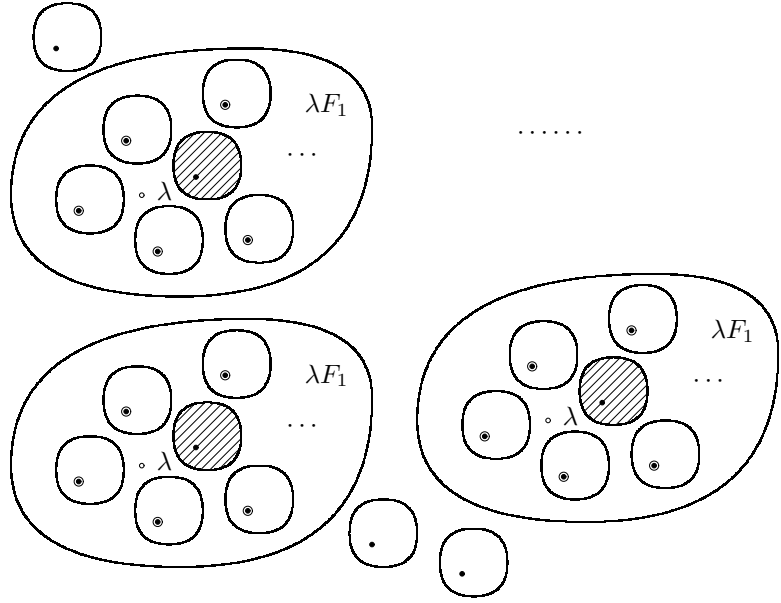
\begin{figure}[ht]
\begin{center}
\setlength{\unitlength}{3mm}
\begin{picture}(30,30)(5,-.5)


\put(4,3){
\put(14,8){\makebox(0,0)[b]{$\lambda F_1$}}
\qbezier(0,4.5)(0,11)(10.5,11)
\qbezier(0,4.5)(0,0)(7.5,0)
\qbezier(7.5,0)(16,0)(16,7.5)
\qbezier(10.5,11)(16,11)(16,7.5)

\put(-1,-0.2){
\qbezier(3,4.5)(3,3)(4.5,3)
\qbezier(3,4.5)(3,6)(4.5,6)
\qbezier(4.5,3)(6,3)(6,4.5)
\qbezier(4.5,6)(6,6)(6,4.5)
\put(4,4){\circle*{0.2}}
\put(4,4){\circle{0.4}}
}

\put(2.5,-2){
\qbezier(3,4.5)(3,3)(4.5,3)
\qbezier(3,4.5)(3,6)(4.5,6)
\qbezier(4.5,3)(6,3)(6,4.5)
\qbezier(4.5,6)(6,6)(6,4.5)
\put(4,4){\circle*{0.2}}
\put(4,4){\circle{0.4}}
}

\put(4.2,1.3){
\qbezier(3,4.5)(3,3)(4.5,3)
\qbezier(3,4.5)(3,6)(4.5,6)
\qbezier(4.5,3)(6,3)(6,4.5)
\qbezier(4.5,6)(6,6)(6,4.5)
\put(4,4){\circle*{0.2}}

\put(3.4,3.4){\line(1,1){2.22}}
\put(3.6,3.2){\line(1,1){2.2}}
\put(3.88,3.08){\line(1,1){2.04}}
\put(4.21,3.01){\line(1,1){1.77}}
\put(4.62,3.02){\line(1,1){1.37}}
\put(5.07,3.08){\line(1,1){0.88}}

\put(3.2,3.6){\line(1,1){2.2}}
\put(3.08,3.88){\line(1,1){2.04}}
\put(3.01,4.21){\line(1,1){1.77}}
\put(3.02,4.62){\line(1,1){1.37}}
\put(3.08,5.07){\line(1,1){0.88}}
}

\put(1.1,2.9){
\qbezier(3,4.5)(3,3)(4.5,3)
\qbezier(3,4.5)(3,6)(4.5,6)
\qbezier(4.5,3)(6,3)(6,4.5)
\qbezier(4.5,6)(6,6)(6,4.5)
\put(4,4){\circle*{0.2}}
\put(4,4){\circle{0.4}}
}

\put(5.5,4.5){
\qbezier(3,4.5)(3,3)(4.5,3)
\qbezier(3,4.5)(3,6)(4.5,6)
\qbezier(4.5,3)(6,3)(6,4.5)
\qbezier(4.5,6)(6,6)(6,4.5)
\put(4,4){\circle*{0.2}}
\put(4,4){\circle{0.4}}
}

\put(6.5,-1.5){
\qbezier(3,4.5)(3,3)(4.5,3)
\qbezier(3,4.5)(3,6)(4.5,6)
\qbezier(4.5,3)(6,3)(6,4.5)
\qbezier(4.5,6)(6,6)(6,4.5)
\put(4,4){\circle*{0.2}}
\put(4,4){\circle{0.4}}
}

\put(13,6){\makebox(0,0)[b]{$\cdots$}}
\put(5.8,4.5){\circle{0.2}}
\put(6.8,4.3){\makebox(0,0)[b]{$\lambda$}}
}

\put(4,15){
\put(14,8){\makebox(0,0)[b]{$\lambda F_1$}}
\qbezier(0,4.5)(0,11)(10.5,11)
\qbezier(0,4.5)(0,0)(7.5,0)
\qbezier(7.5,0)(16,0)(16,7.5)
\qbezier(10.5,11)(16,11)(16,7.5)

\put(-1,-0.2){
\qbezier(3,4.5)(3,3)(4.5,3)
\qbezier(3,4.5)(3,6)(4.5,6)
\qbezier(4.5,3)(6,3)(6,4.5)
\qbezier(4.5,6)(6,6)(6,4.5)
\put(4,4){\circle*{0.2}}
\put(4,4){\circle{0.4}}
}

\put(2.5,-2){
\qbezier(3,4.5)(3,3)(4.5,3)
\qbezier(3,4.5)(3,6)(4.5,6)
\qbezier(4.5,3)(6,3)(6,4.5)
\qbezier(4.5,6)(6,6)(6,4.5)
\put(4,4){\circle*{0.2}}
\put(4,4){\circle{0.4}}
}

\put(4.2,1.3){
\qbezier(3,4.5)(3,3)(4.5,3)
\qbezier(3,4.5)(3,6)(4.5,6)
\qbezier(4.5,3)(6,3)(6,4.5)
\qbezier(4.5,6)(6,6)(6,4.5)
\put(4,4){\circle*{0.2}}

\put(3.4,3.4){\line(1,1){2.22}}
\put(3.6,3.2){\line(1,1){2.2}}
\put(3.88,3.08){\line(1,1){2.04}}
\put(4.21,3.01){\line(1,1){1.77}}
\put(4.62,3.02){\line(1,1){1.37}}
\put(5.07,3.08){\line(1,1){0.88}}

\put(3.2,3.6){\line(1,1){2.2}}
\put(3.08,3.88){\line(1,1){2.04}}
\put(3.01,4.21){\line(1,1){1.77}}
\put(3.02,4.62){\line(1,1){1.37}}
\put(3.08,5.07){\line(1,1){0.88}}
}

\put(1.1,2.9){
\qbezier(3,4.5)(3,3)(4.5,3)
\qbezier(3,4.5)(3,6)(4.5,6)
\qbezier(4.5,3)(6,3)(6,4.5)
\qbezier(4.5,6)(6,6)(6,4.5)
\put(4,4){\circle*{0.2}}
\put(4,4){\circle{0.4}}
}

\put(5.5,4.5){
\qbezier(3,4.5)(3,3)(4.5,3)
\qbezier(3,4.5)(3,6)(4.5,6)
\qbezier(4.5,3)(6,3)(6,4.5)
\qbezier(4.5,6)(6,6)(6,4.5)
\put(4,4){\circle*{0.2}}
\put(4,4){\circle{0.4}}
}

\put(6.5,-1.5){
\qbezier(3,4.5)(3,3)(4.5,3)
\qbezier(3,4.5)(3,6)(4.5,6)
\qbezier(4.5,3)(6,3)(6,4.5)
\qbezier(4.5,6)(6,6)(6,4.5)
\put(4,4){\circle*{0.2}}
\put(4,4){\circle{0.4}}
}

\put(13,6){\makebox(0,0)[b]{$\cdots$}}
\put(5.8,4.5){\circle{0.2}}
\put(6.8,4.3){\makebox(0,0)[b]{$\lambda$}}
}

\put(22,5){
\put(14,8){\makebox(0,0)[b]{$\lambda F_1$}}
\qbezier(0,4.5)(0,11)(10.5,11)
\qbezier(0,4.5)(0,0)(7.5,0)
\qbezier(7.5,0)(16,0)(16,7.5)
\qbezier(10.5,11)(16,11)(16,7.5)

\put(-1,-0.2){
\qbezier(3,4.5)(3,3)(4.5,3)
\qbezier(3,4.5)(3,6)(4.5,6)
\qbezier(4.5,3)(6,3)(6,4.5)
\qbezier(4.5,6)(6,6)(6,4.5)
\put(4,4){\circle*{0.2}}
\put(4,4){\circle{0.4}}
}

\put(2.5,-2){
\qbezier(3,4.5)(3,3)(4.5,3)
\qbezier(3,4.5)(3,6)(4.5,6)
\qbezier(4.5,3)(6,3)(6,4.5)
\qbezier(4.5,6)(6,6)(6,4.5)
\put(4,4){\circle*{0.2}}
\put(4,4){\circle{0.4}}
}

\put(4.2,1.3){
\qbezier(3,4.5)(3,3)(4.5,3)
\qbezier(3,4.5)(3,6)(4.5,6)
\qbezier(4.5,3)(6,3)(6,4.5)
\qbezier(4.5,6)(6,6)(6,4.5)
\put(4,4){\circle*{0.2}}

\put(3.4,3.4){\line(1,1){2.22}}
\put(3.6,3.2){\line(1,1){2.2}}
\put(3.88,3.08){\line(1,1){2.04}}
\put(4.21,3.01){\line(1,1){1.77}}
\put(4.62,3.02){\line(1,1){1.37}}
\put(5.07,3.08){\line(1,1){0.88}}

\put(3.2,3.6){\line(1,1){2.2}}
\put(3.08,3.88){\line(1,1){2.04}}
\put(3.01,4.21){\line(1,1){1.77}}
\put(3.02,4.62){\line(1,1){1.37}}
\put(3.08,5.07){\line(1,1){0.88}}
}

\put(1.1,2.9){
\qbezier(3,4.5)(3,3)(4.5,3)
\qbezier(3,4.5)(3,6)(4.5,6)
\qbezier(4.5,3)(6,3)(6,4.5)
\qbezier(4.5,6)(6,6)(6,4.5)
\put(4,4){\circle*{0.2}}
\put(4,4){\circle{0.4}}
}

\put(5.5,4.5){
\qbezier(3,4.5)(3,3)(4.5,3)
\qbezier(3,4.5)(3,6)(4.5,6)
\qbezier(4.5,3)(6,3)(6,4.5)
\qbezier(4.5,6)(6,6)(6,4.5)
\put(4,4){\circle*{0.2}}
\put(4,4){\circle{0.4}}
}

\put(6.5,-1.5){
\qbezier(3,4.5)(3,3)(4.5,3)
\qbezier(3,4.5)(3,6)(4.5,6)
\qbezier(4.5,3)(6,3)(6,4.5)
\qbezier(4.5,6)(6,6)(6,4.5)
\put(4,4){\circle*{0.2}}
\put(4,4){\circle{0.4}}
}

\put(13,6){\makebox(0,0)[b]{$\cdots$}}
\put(5.8,4.5){\circle{0.2}}
\put(6.8,4.3){\makebox(0,0)[b]{$\lambda$}}
}

\put(28,22){\makebox(0,0)[b]{$\cdots\cdots$}}

\put(16,0){
\qbezier(3,4.5)(3,3)(4.5,3)
\qbezier(3,4.5)(3,6)(4.5,6)
\qbezier(4.5,3)(6,3)(6,4.5)
\qbezier(4.5,6)(6,6)(6,4.5)
\put(4,4){\circle*{0.2}}
}

\put(20,-1.3){
\qbezier(3,4.5)(3,3)(4.5,3)
\qbezier(3,4.5)(3,6)(4.5,6)
\qbezier(4.5,3)(6,3)(6,4.5)
\qbezier(4.5,6)(6,6)(6,4.5)
\put(4,4){\circle*{0.2}}
}

\put(2,22){
\qbezier(3,4.5)(3,3)(4.5,3)
\qbezier(3,4.5)(3,6)(4.5,6)
\qbezier(4.5,3)(6,3)(6,4.5)
\qbezier(4.5,6)(6,6)(6,4.5)
\put(4,4){\circle*{0.2}}
}

\end{picture}
\caption{\label{fig:F4} An illustration of the set $\Delta_1(\gamma_1F_0\cup D^1_0)$. In the figure $\lambda$ denotes a generic element of $\Delta_1$. The set consists
of the shaded regions together with all highlighted (circled-solid) points.
For details of each translates of $F_1$ see Figure~\ref{fig:F3}.}
\end{center}
\end{figure}

Let us first dwell for a moment on (1). Condition (1) is consistent with condition (2) because $\Delta_n a_n$ and $\Delta_n b_n$ are decreasing sequences and $\Delta_{n+1} \Lambda_{n+1} b_n \subseteq \Delta_n b_n$ (conclusions (i) and (viii) of Lemma \ref{BP LIST}). By conclusions (xi) and (xii) of Lemma \ref{BP LIST}, we have $\Delta_n a_n$ is disjoint from $\Delta_n b_n \cup \bigcup_{1 \leq k \leq n} \Delta_k \Lambda_k b_{k-1}$. Therefore for $n > 1$ we desire $\dom(c_n)$ to be
$$\dom(c_{n-1}) \cup (\Delta_{n-1} a_{n-1} - \Delta_n a_n) \cup (\Delta_{n-1} b_{n-1} - \Delta_n [\Lambda_n \cup \{\beta_n\}] b_{n-1}).$$
It is important to note that these unions are disjoint. This tells us that given $c_{n-1}$, we can define $c_n \supseteq c_{n-1}$ to have whichever values on $\Delta_{n-1} a_{n-1} - \Delta_n a_n$ and $\Delta_{n-1} b_{n-1} - \Delta_n [\Lambda_n \cup \{\beta_n\}] b_{n-1}$ without worry of a contradiction between the two or with $c_{n-1}$.

Define
$$c_1: (G - \Delta_1 a_1 - \Delta_1 b_1 - \Delta_1 \Lambda_1) \rightarrow \{0, 1\}$$
by
$$c_1 (g) = \begin{cases}
R(f) & \text{if } g = \gamma \gamma_1 f \text{ where } \gamma \in \Delta_1 \text{ and } f \in F_0 \\
1 - R(1_G) & \text{otherwise}
\end{cases}$$
for $g \in \dom(c_1)$. The function $c_1$ satisfies (1) since $b_0 = 1_G$. Notice that the set referred to in (iv), $G - \Delta_1 (\gamma_1 F_0 \cup D^1_0)$, is a subset of the domain of $c_1$. Clearly each element of this set is mapped to $1 - R(1_G)$ by $c_1$. So as long as the final function $c$ extends $c_1$ clause (iv) will be satisfied. See Figures \ref{fig:F3} and \ref{fig:F4} for an illustration of the set $G - \Delta_1 (\gamma_1 F_0 \cup D^1_0)$.

We claim $c_1$ satisfies (3) with test region $\gamma_1 F_0$. Since $\Delta_1 \gamma_1$, $\Delta_1 a_1$, $\Delta_1 b_1$, and $\Delta_1 \Lambda_1$ are pairwise disjoint subsets of $\Delta_0$, we have that $\Delta_1 \gamma_1 F_0$ is disjoint from $\Delta_1 a_1 \cup \Delta_1 b_1 \cup \Delta_1 \Lambda_1$ (since $1_G \in F_0$). Thus $\Delta_1 \gamma_1 F_0 \subseteq \dom(c_1)$ as required.

Let $x \in 2^G$ be an arbitrary extension of $c_1$. To finish verifying (3), we will show $g \in \Delta_1$ if and only if for all $f \in F_0$ $x(g \gamma_1 f) = R(f)$. If $\gamma \in \Delta_1$, then $\gamma \gamma_1 F_0 \subseteq \dom(c_1)$ and hence $x(\gamma \gamma_1 f) = R(f)$ for all $f \in F_0$. Now suppose $g \in G$ satisfies $x(g \gamma_1 f) = R(f)$ for all $f \in F_0$. Note that $x(h) = c_1(h) = 1 - R(1_G)$ for all $h \in \dom(c_1) - \Delta_1 \gamma_1 F_0$. As $x(g \gamma_1 1_G) = R(1_G)$, either $g \gamma_1 \in \Delta_1 \gamma_1 F_0$ or $g \gamma_1 \not\in \dom(c_1)$. But $g \gamma_1$ cannot be in $G - \dom(c_1) \subseteq \Delta_0$, for then we would have
$$g \gamma_1 F_0 - \{g \gamma_1\} \subseteq \dom(c_1) - \Delta_1 \gamma_1 F_0$$
and hence $R = (\gamma_1^{-1} g^{-1} \cdot x) \res F_0$ would be trivial, a contradiction. So $g \gamma_1 \in \Delta_1 \gamma_1 F_0$. Let $\gamma \in \Delta_1$ and $a \in F_0$ be such that $g \gamma_1 = \gamma \gamma_1 a$. By construction, $x(\gamma \gamma_1 f) = R(f)$ for all $f \in F_0$, and we have that for all $b \in F_0$
$$(\gamma_1^{-1} \gamma^{-1} \cdot x)(a b) = x(\gamma \gamma_1 a b) = x(g \gamma_1 b) = R(b) = (\gamma_1^{-1} \gamma^{-1} \cdot x)(b).$$
Now it follows from Lemma \ref{LEM LR EQUIV} that $a = 1_G$. Thus, $g \gamma_1  = \gamma \gamma_1$ and $g = \gamma \in \Delta_1$.

Now suppose that $c_1$ through $c_{k-1}$ have been constructed and satisfy $(1)$ through $(3)$. We pointed out earlier that we desire $c_k$ to have domain
$$\dom(c_{k-1}) \cup (\Delta_{k-1} a_{k-1} - \Delta_k a_k) \cup (\Delta_{k-1} b_{k-1} - \Delta_k [\Lambda_k \cup \{\beta_k\}] b_{k-1}).$$
We define $c_k$ to satisfy $c_k \supseteq c_{k-1}$ and:
$$c_k (\Delta_{k-1} a_{k-1} - \Delta_k \{\gamma_k, \alpha_k\} a_{k-1}) = \{0\};$$
$$c_k (\Delta_k \gamma_k a_{k-1}) = \{1\};$$
$$c_k (\Delta_k \gamma_k b_{k-1}) = \{1\};$$
$$c_k (\Delta_k \alpha_k b_{k-1}) = \{0\};$$
$$c_k (\Delta_{k-1} b_{k-1} - \Delta_k D^k_{k-1} b_{k-1}) = \{0\}.$$
From our earlier remarks on (1), we know $c_k$ is well defined. It is easily checked that $c_k$ satisfies (1) and (2) (recall that $\Delta_k a_k = \Delta_k \alpha_k a_{k-1}$ and $\Lambda_k = D^k_{k-1} - \{\alpha_k, \beta_k, \gamma_k\}$). See Figure \ref{fig:F5} for an illustration of $c_k$.


\begin{figure}[ht]
\begin{center}
\setlength{\unitlength}{1cm}
\begin{picture}(15,8)(-.3,0)

\put(6,3.7){\oval(12,7.5)}
\put(11, 6.5){\makebox(0,0)[b]{$\gamma F_k$}}

\put(0.4,4.4){
\put(1.5,1.0){\oval(2.8,2.6)}
\put(1.8,.8){\circle*{.1}}
\put(1.8,1.5){\circle*{.1}}
\put(1,.7){\makebox(0,0)[b]{$\gamma \alpha_k b_{k-1}$}}
\put(2.1, .7){\makebox(0,0)[b]{$0$}}
\put(1.4,1.4){\makebox(0,0)[b]{$\gamma a_k$}}
\put(2.1, 1.4){\makebox(0,0)[b]{$?$}}
\put(1.5, -0.8){\makebox(0,0)[b]{$\gamma\alpha_kF_{k-1}$}}
}

\put(3.8,4.4){
\put(1.5,1.0){\oval(2.8,2.6)}
\put(1.8,.8){\circle*{.1}}
\put(1.8,1.5){\circle*{.1}}
\put(1.4,.7){\makebox(0,0)[b]{$\gamma b_k$}}
\put(2.1, .7){\makebox(0,0)[b]{$?$}}
\put(1,1.4){\makebox(0,0)[b]{$\gamma \beta_ka_{k-1}$}}
\put(2.1, 1.4){\makebox(0,0)[b]{$0$}}
\put(1.5, -0.8){\makebox(0,0)[b]{$\gamma\beta_kF_{k-1}$}}
}

\put(7.2,4.4){
\put(1.5,1.0){\oval(2.8,2.6)}
\put(1.8,.8){\circle*{.1}}
\put(1.8,1.5){\circle*{.1}}
\put(1,.7){\makebox(0,0)[b]{$\gamma \gamma_k b_{k-1}$}}
\put(2.1, .7){\makebox(0,0)[b]{$1$}}
\put(1,1.4){\makebox(0,0)[b]{$\gamma \gamma_ka_{k-1}$}}
\put(2.1, 1.4){\makebox(0,0)[b]{1}}
\put(1.5, -0.8){\makebox(0,0)[b]{$\gamma\gamma_kF_{k-1}$}}
}

\put(1.8,1.0){
\put(1.5,1.0){\oval(2.8,2.6)}
\put(1.8,.8){\circle*{.1}}
\put(1.8,1.5){\circle*{.1}}
\put(1,.7){\makebox(0,0)[b]{$\gamma \lambda b_{k-1}$}}
\put(2.1, .7){\makebox(0,0)[b]{$?$}}
\put(1,1.4){\makebox(0,0)[b]{$\gamma \lambda a_{k-1}$}}
\put(2.1, 1.4){\makebox(0,0)[b]{$0$}}
\put(2.0, -0.8){\makebox(0,0)[b]{$\gamma\lambda F_{k-1}$, $\lambda\in \Lambda_k$}}
}

\put(6.8, 2.0){\makebox(0,0)[b]{$\cdots\cdots$}}

\end{picture}
\caption{\label{fig:F5} The definition of $c_k$ ensures a simple membership test for $\Delta_k$}
\end{center}
\end{figure}
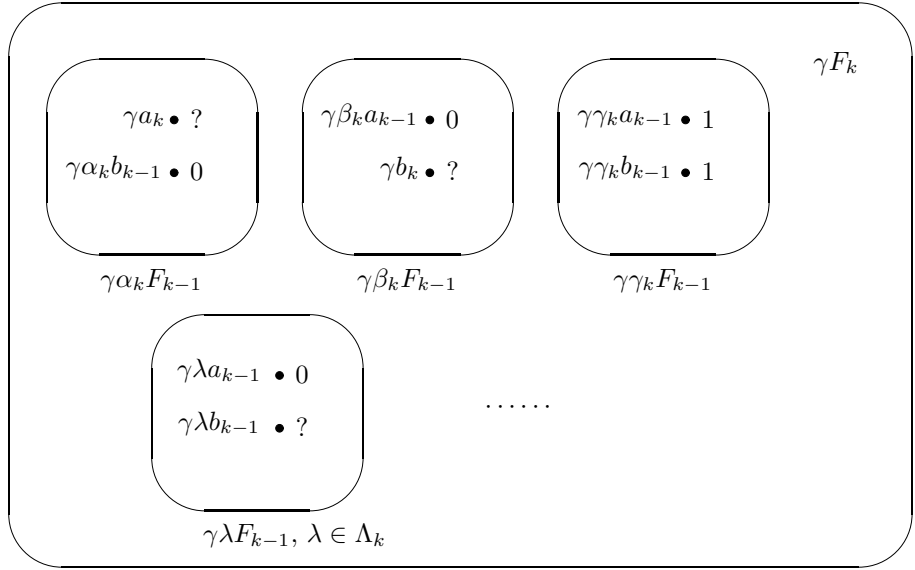

Let $V \subseteq \gamma_{k-1} F_{k-2}$ be the test region referred to in (3) for $n = k-1$, and let $v \in 2^V$ witness the membership test. Set $W = \gamma_k (V \cup \{a_{k-1}, b_{k-1}\})$ and define $w \in 2^W$ so that $w$ extends $\gamma_k \cdot v$ and $w(\gamma_k a_{k-1}) = w(\gamma_k b_{k-1}) = 1$. We claim that $c_k$ satisfies (3) with test region $W$ and witnessing function $w$. Clearly $W \subseteq \gamma_k F_{k-1}$. Let $x \in 2^G$ be an arbitrary extension of $c_k$. If $\gamma \in \Delta_k$, then $\gamma \gamma_k \in \Delta_{k-1}$ and it is then clear from the definition of $c_k$ that $x(\gamma a) = w(a)$ for all $a \in W$. Now suppose $g \in G$ satisfies $x(ga) = w(a)$ for all $a \in W$. Then in particular $x(g \gamma_k a) = v(a)$ for all $a \in V$ and thus $g \gamma_k \in \Delta_{k-1}$. Also, $x(g \gamma_k a_{k-1}) = w(a_{k-1}) = 1$. From how we defined $c_k$, we have for $h \in G$
$$h \in \Delta_{k-1} \text{ and } x(h a_{k-1}) = 1 \Longrightarrow h \in \Delta_k \gamma_k \text{ or } h \in \Delta_k \alpha_k.$$
However, $g \gamma_k \not\in \Delta_k \alpha_k$, for otherwise we would have
$$1 = w(\gamma_k b_{k-1}) = x(g \gamma_k b_{k-1}) = c_k(g \gamma_k b_{k-1}) = 0.$$
We conclude $g \gamma_k \in \Delta_k \gamma_k$ and $g \in \Delta_k$. Thus $c_k$ satisfies (3).

Finally, take the function $c' = \bigcup_{n \geq 1} c_n$ and if necessary extend $c'$ arbitrarily to $\bigcap_{n \in \N} \Delta_n a_n$ and $\bigcap_{n \in \N} \Delta_n b_n$ to get a function $c \in 2^{\subseteq G}$. Properties (i) and (iv) clearly hold due to how $c_1$ was defined. Property (ii) holds since $c \supseteq c_n$ for each $n \geq 1$, and property (iii) follows from (1) and conclusion (xiii) of Lemma \ref{BP LIST} (for the disjointness of the union). We proceed to verify properties (v), (vi), and (vii).

(v). Fix $n \geq 1$ and $\gamma \in \Delta_n$. By conclusion (xii) of Lemma \ref{BP LIST}, $\gamma F_n - \{\gamma b_n\}$ is disjoint from $\Delta_k \Lambda_k b_{k-1}$ for all $k > n$. Also, $\gamma b_n \in \Delta_k b_k =\Delta_k \beta_k b_{k-1}$ for all $k \leq n$ by conclusion (viii) of Lemma \ref{BP LIST}. Since $\Delta_k \beta_k$ and $\Delta_k \Lambda_k$ are disjoint subsets of $\Delta_{k-1}$ and $b_{k-1} \in F_{k-1}$, it follows that $\gamma b_n \not\in \Delta_k \Lambda_k b_{k-1}$ for $k \leq n$. Thus for $k \leq n$
$$(\gamma F_n - \{\gamma b_n\}) \cap \Delta_k \Lambda_k b_{k-1} = \Delta_k \Lambda_k b_{k-1} \cap \gamma F_n = \gamma D^n_k \Lambda_k b_{k-1}$$
by conclusion (vi) of Lemma \ref{BP LIST}. The claim now follows from (iii).

(vi). Fix $n \geq 1$, $\gamma, \sigma \in \Delta_n$, and $f \in F_n - \{a_n, b_n\} - \bigcup_{1 \leq k \leq n} D^n_k \Lambda_k b_{k-1}$. Then by (v) $\gamma f, \sigma f \in \dom(c)$. Also, $f \not\in \{a_n, b_n\}$ and hence $\gamma f, \sigma f \not\in \Delta_n \{a_n, b_n\}$ since the $\Delta_n$-translates of $F_n$ are disjoint and $a_n, b_n \in F_n$. However, $\dom(c) - \dom(c_n) \subseteq \Delta_n \{a_n, b_n\}$, and since $\gamma f, \sigma f \in \dom(c)$ there must be a $m \leq n$ with $\gamma f, \sigma f \in \dom(c_m)$. Let $k \leq n$ be minimal with either $\gamma f$ or $\sigma f$ in $\dom(c_k)$. If $k = 1$, then conclusion (vii) of Lemma \ref{BP LIST} implies that both $\gamma f, \sigma f \in \dom(c_1)$ and $c(\gamma f) = c_1(\gamma f) = c_1(\sigma f) = c(\sigma f)$. Similarly, if $k > 1$ then, after recalling the five equations used to define $c_k$, conclusion (vii) of Lemma \ref{BP LIST} again implies that both $\gamma f, \sigma f \in \dom(c_k)$ and $c(\gamma f) = c_k(\gamma f) = c_k( \sigma f) = c(\sigma f)$.

(vii). For $n \geq 1$ set $c^n = c \res(G - \Delta_n \{a_n, b_n\})$. Fix $n \geq 1$ and $m \geq n$. Let $\gamma, \psi \in \Delta_m$ and let $f \in F_m$. We must show $\gamma f \in \dom(c^n)$ if and only if $\psi f \in \dom(c^n)$ and furthermore $c^n(\gamma f) = c^n(\psi f)$ when these are defined. Since $\Delta_n \{a_n, b_n\} \subseteq G - \dom(c^n)$, we are done in the case $f$ is $a_m$ or $b_m$ (since $\Delta_m \{a_m, b_m\} \subseteq \Delta_n \{a_n, b_n\}$ by conclusion (viii) of Lemma \ref{BP LIST}). So we assume $f$ is neither $a_m$ nor $b_m$. Then by (v) $\gamma f \in \dom(c)$ if and only if $\psi f \in \dom(c)$, and by conclusion (vii) of Lemma \ref{BP LIST} $\gamma f \in \Delta_n \{a_n, b_n\}$ if and only if $\psi f \in \Delta_n \{a_n, b_n\}$. Therefore, $\gamma f \in \dom(c^n)$ if and only if $\psi f \in \dom(c^n)$. Finally, by (vi) we have that $c^n(\gamma f) = c^n(\psi f)$ whenever these are defined.
\end{proof}

Although the previous theorem applies to pre-blueprints, we restrict ourselves to blueprints for the following two definitions.

\begin{definition} \index{canonical function} \index{canonical function!compatible with locally recognizable function}
A function $c \in 2^{\subseteq G}$ is called \emph{canonical} if for some blueprint $(\Delta_n, F_n)_{n \in \N}$ and some nontrivial locally recognizable function $R: F_0 \rightarrow 2$ the conclusions of Theorem \ref{FM} are satisfied. If we wish to emphasize the blueprint, we say $c$ is \emph{canonical with respect to} $(\Delta_n, F_n)_{n \in \N}$. To emphasize the locally recognizable function, we say $c$ is \emph{compatible} with $R$.
\end{definition}

\begin{definition} \label{DEF FUNDF} \index{fundamental function}
Let $(\Delta_n, F_n)_{n \in \N}$ be a blueprint. A function $c \in 2^{\subseteq G}$ is called \emph{fundamental with respect to} $(\Delta_n, F_n)_{n \in \N}$ if some function $c' \subseteq c$ is canonical with respect to $(\Delta_n, F_n)_{n \in \N}$ and there are sets $\Theta_n \subseteq \Lambda_n$ for each $n \geq 1$ such that
$$\dom(c) = G - \bigcup_{n \geq 1} \Delta_n \Theta_n b_{n-1}.$$
In this case, if $R$ is a nontrivial locally recognizable function and $c'$ is compatible with $R$, then we say $c$ is \emph{compatible} with $R$ as well. We simply call $c \in 2^{\subseteq G}$ \emph{fundamental} if it is fundamental with respect to some blueprint.
\end{definition}

Notice that every canonical function is fundamental: set $\Theta_n = \Lambda_n$.

\begin{remark} \index{$\Theta_n$} \index{fundamental function!$\Theta_n$}
When dealing with a fundamental function $c \in 2^{\subseteq G}$, the symbols $\Theta_n$ for each $n \geq 1$ will be reserved. $\Theta_n$ will necessarily be as to satisfy the above definition.
\end{remark}

Clause (v) of Theorem \ref{FM} can be adapted for fundamental functions.

\begin{lem} \label{FUND DOM}
Let $G$ be a countably infinite group, let $(\Delta_n, F_n)_{n \in \N}$ be a blueprint, and let $c \in 2^{\subseteq G}$ be fundamental with respect to this blueprint. Then for all $n \geq 1$ and $\gamma \in \Delta_n$
$$(\gamma F_n - \{\gamma b_n\}) \cap \dom(c) = \gamma \left( F_n - \{b_n\} - \bigcup_{1 \leq k \leq n} D^n_k \Theta_k b_{k-1} \right).$$
\end{lem}

\begin{proof}
Fix $n \geq 1$ and $\gamma \in \Delta_n$. By conclusion (xii) of Lemma \ref{BP LIST}, $\gamma F_n - \{\gamma b_n\}$ is disjoint from $\Delta_k \Theta_k b_{k-1}$ for all $k > n$. Also, $\gamma b_n \in \Delta_k b_k =\Delta_k \beta_k b_{k-1}$ for all $k \leq n$ by conclusion (viii) of Lemma \ref{BP LIST}. Since $\Delta_k \beta_k$ and $\Delta_k \Theta_k$ are disjoint subsets of $\Delta_{k-1}$ and $b_{k-1} \in F_{k-1}$, it follows that $\gamma b_n \not\in \Delta_k \Theta_k b_{k-1}$ for $k \leq n$. Thus for $k \leq n$
$$(\gamma F_n - \{\gamma b_n\}) \cap \Delta_k \Theta_k b_{k-1} = \Delta_k \Theta_k b_{k-1} \cap \gamma F_n = \gamma D^n_k \Theta_k b_{k-1}$$
by conclusion (vi) of Lemma \ref{BP LIST}. The claim now follows from the fact that
$$\dom(c) = G - \bigcup_{n \geq 1} \Delta_n \Theta_n b_{n-1}.$$
\end{proof}

\section{Existence of blueprints} \label{SECT EXIST BP}

In this section we show that every countably infinite group admits a blueprint. All of our future results in the paper rely on blueprints and therefore their existence is vitally important. It is not difficult to construct sequences $(\Delta_n, F_n)_{n \in \N}$ which are dense, nor is it difficult to construct sequences which are coherent. However, the truth is that the key difficulty in constructing a blueprint is simultaneously achieving the coherent property and the dense property. First we outline a simple way to construct pre-blueprints.

\begin{lem} \label{MAKE PREBP}
Let $G$ be a countably infinite group. Let $(F_n)_{n \in \N}$ be a sequence of finite subsets of $G$, and let $\{\delta^n_k \: n, k \in \N, \ k < n\}$ be a collection of finite subsets of $G$ satisfying
\begin{enumerate}
\item[\rm (i)] $1_G \in \delta^n_{n-1}$ for each $n \geq 1$;
\item[\rm (ii)] $|\delta^n_{n-1}| \geq 3$ for each $n \geq 1$;
\item[\rm (iii)] the $\delta^n_k$-translates of $F_k$ are disjoint for all $n, k \in \N$ with $k < n$;
\item[\rm (iv)] $\delta^n_m F_m \cap \delta^n_k F_k = \varnothing$ for all $m \neq k < n$;
\item[\rm (v)] $F_n = \bigcup_{0 \leq k < n} \delta^n_k F_k$ for all $n \geq 1$.
\end{enumerate}
Then there is a sequence $(\Delta_n)_{n \in \N}$ of subsets of $G$ with $\delta^n_k \subseteq \Delta_k$ for every $n, k \in \N$ with $k < n$ and such that $(\Delta_n, F_n)_{n \in \N}$ is a centered and directed pre-blueprint.
\end{lem}

\begin{proof}
The basic idea is that setting $\Delta_k = \bigcup_{n > k} \delta^n_k$ nearly works, except that we have to enlarge this set in order to satisfy the uniform property of pre-blueprints. Since in the end we want $\delta^n_k \subseteq \Delta_k$ and $\delta^n_k \subseteq F_n$ (so $\delta^n_k \subseteq \Delta_k \cap F_n$), in order to achieve the uniform property it is necessary that $\delta^n_k$ be copied wherever any $\delta^m_n$-translates of $F_n$ are located. In other words, for each $n \in \N$ we want to view $F_n$ as carrying all of the sets $\delta^n_k$ ($k < n$) with it. With this mind set, we want to recognize all of the translates of $F_k$ which explicitly or implicitly make up a part of $F_n$. For example, for $k < m < n$ we have $\delta_k^m F_k \subseteq F_m$ and $\delta_m^n F_m \subseteq F_n$ so $\delta_m^n \delta_k^m F_k \subseteq F_n$. Thus, informally we would say the $\delta_m^n \delta_k^m$-translates of $F_k$ are implicitly a part of $F_n$. On the other hand, if for $g \in F_n$ we only have $g F_k \subseteq F_n$ we would not necessarily want to say the $g$-translate of $F_k$ is a part of $F_n$. We will momentarily define sets $D^n_k$. The choice of notation is no mistake. Later when we define the $\Delta_n$'s we will show that the $D^n_k$'s carry the usual meaning for pre-blueprints. Informally, we want $D_k^n$ to be the set of all $g$'s in $F_n$ such that the $g$-translate of $F_k$ either explicitly or implicitly makes up a part of $F_n$. We now give the formal definition for this. For $k \in \N$, define $D_k^k = \{1_G\}$, $D_k^{k+1}= \delta_k^{k+1}$, and in general for $n>k$
$$D_k^n= \delta_{n-1}^n D_k^{n-1} \cup \delta_{n-2}^n D_k^{n-2} \cup \cdots \cup \delta_{k+1}^n D_k^{k+1} \cup \delta_k^n = \bigcup_{k \leq m < n} \delta_m^n D_k^m.$$

We first verify that the $D_k^n$'s possess the following properties:
\begin{enumerate}
\item[\rm (1)] $D_k^n F_k \subseteq F_n$ for all $k, n \in \N$ with $k \leq n$;
\item[\rm (2)] $D_m^n D_k^m \subseteq D_k^n$ for all $k, m, n \in \N$ with $k \leq m \leq n$;
\item[\rm (3)] the $D_k^n$-translates of $F_k$ are disjoint for all $k, n \in \N$ with $k \leq n$;
\end{enumerate}

(Proof of 1). Clearly $D_k^k F_k = F_k$. If we assume $D_k^m F_k \subseteq F_m$ for all $k \leq m < n$, then by (v)
$$D_k^n F_k = \bigcup_{k\leq m<n} \delta_m^n D_k^m F_k \subseteq \bigcup_{k\leq m<n} \delta_m^n F_m \subseteq F_n.$$
The claim now immediately follows from induction.

(Proof of 2). Clearly when $n = m$ we have $D_m^n D_k^m = D_n^n D_k^n = D_k^n$. If we assume $D_m^t D_k^m \subseteq D_k^t$ for all $m \leq t < n$, then
$$D_m^n D_k^m = \bigcup_{m \leq t < n} \delta_t^n D_m^t D_k^m \subseteq \bigcup_{m \leq t < n} \delta_t^n D_k^t \subseteq \bigcup_{k \leq t < n} \delta_t^n D_k^t = D_k^n.$$
The claim now immediately follows from induction.

(Proof of 3). The $D_k^n$-translates of $F_k$ are disjoint when $n = k$ and when $n = k+1$ (by (iii)). Assume the $D_k^m$ translates of $F_k$ are disjoint for all $k \leq m < n$. Recall $D_k^n = \bigcup_{k \leq m < n} \delta_m^n D_k^m$. If $k \leq r < s < n$, then by (iv) we have $\delta_r^n F_r \cap \delta_s^n F_s = \varnothing$. It then follows from (1) that $\delta_r^n D_k^r F_k \cap \delta_s^n D_k^s F_k = \varnothing$. Additionally, if $k \leq m < n$ and $\gamma, \psi \in \delta_m^n$ are distinct, then $\gamma F_m \cap \psi F_m = \varnothing$ by (iii). Again by (1) we have $\gamma D_k^m F_k \cap \psi D_k^m F_k = \varnothing$. Finally, by assumption the $D_k^m$-translates of $F_k$ are disjoint for every $k \leq m < n$. So in particular, for each $k \leq m < n$ and $\gamma \in \delta^n_m$ the $\gamma D^m_k$-translates of $F_k$ are disjoint. It follows that the $D_k^n$-translates of $F_k$ must be disjoint. The claim now follows from induction.

We point out that $D_k^n \subseteq D_k^{n+1}$ since $\delta_n^{n+1} D_k^n \subseteq D_k^{n+1}$ and $1_G \in \delta_n^{n+1}$ by (i). For $k \in \N$ we define $\Delta_k = \bigcup_{n \geq k} D_k^n$. We now check that $(\Delta_n, F_n)_{n \in \N}$ is a pre-blueprint.

(Disjoint). Let $n \in \N$ and $\gamma \neq \psi \in \Delta_n$. Then for some $m > n$ $\gamma, \psi \in D^m_n$. From (3) we then have $\gamma F_n \cap \psi F_n = \varnothing$.

(Coherent). Suppose $k < n$, $\psi \in \Delta_k$, and $\gamma \in \Delta_n$ satisfy $\psi F_k \cap \gamma F_n \neq \varnothing$. Let $m \geq n$ be large enough so that $\psi \in D^m_k$ and $\gamma \in D^m_n$. We will prove $\psi \in \gamma D^n_k$ by induction on $m$. By (1) this will give us $\psi F_k \subseteq \gamma F_n$. Clearly, if $m = n$ then $\gamma = 1_G$ and $\psi \in D_k^n = \gamma D_k^n$. Now suppose the claim is true for all $n \leq i < m$. By the definition of $D^m_k$ and $D^m_n$, there are $k \leq s < m$ and $n \leq t < m$ with $\psi \in \delta^m_s D^s_k$ and $\gamma \in \delta^m_t D^t_n$. However, if $s \neq t$ then by (iv) we have
$$\psi F_k \cap \gamma F_n \subseteq \delta^m_s D^s_k F_k \cap \delta^m_t D^t_n F_n \subseteq \delta^m_s F_s \cap \delta^m_t F_t = \varnothing.$$
So it must be that $s = t$. Let $\lambda, \sigma \in \delta^m_t$ be such that $\psi \in \lambda D^t_k$ and $\gamma \in \sigma D^t_n$. If $\lambda \neq \sigma$ then by (iii) we would have
$$\psi F_k \cap \gamma F_n \subseteq \lambda D^t_k F_k \cap \sigma D^t_n F_n \subseteq \lambda F_t \cap \sigma F_t = \varnothing.$$
So we must have $\lambda = \sigma$. Then $\lambda^{-1} \psi \in D^t_k \subseteq \Delta_k$, $\lambda^{-1} \gamma \in D^t_n \subseteq \Delta_n$, and $\lambda^{-1} \psi F_k \cap \lambda^{-1} \gamma F_n \neq \varnothing$. By the induction hypothesis we conclude $\lambda^{-1} \psi \in \lambda^{-1} \gamma D^n_k$ and hence $\psi \in \gamma D^n_k$.

(Uniform). It suffices to show that $\Delta_k \cap \gamma F_n = \gamma D^n_k$ for $k < n$ and $\gamma \in \Delta_n$. In particular, this will show that $D^n_k$ has its usual meaning. For sufficiently large $m$ $\gamma \in D^m_n$ so $\gamma D^n_k \subseteq D^m_k \subseteq \Delta_k$ by (2). Since $\gamma D^n_k F_k \subseteq \gamma F_n$, we have $\gamma D^n_k \subseteq \Delta_k \cap \gamma F_n$. Conversely, if $\psi \in \Delta_k \cap \gamma F_n$, then in particular $\psi F_k \cap \gamma F_n \neq \varnothing$ since $1_G \in F_k$. Under this assumption, it was shown in the previous paragraph that $\psi \in \gamma D^n_k$.

(Growth). By (ii) $|D^n_{n-1}| = |\delta^n_{n-1}| \geq 3$.

(Centered). By (i) we have $1_G \in \delta^{n+1}_n \subseteq \Delta_n$.

(Directed). Let $n, k \in \N$ and let $\gamma \in \Delta_n$, $\psi \in \Delta_k$. Then for large enough $m$ we have $\gamma \in D^m_n$ and $\psi \in D^m_k$. So by (1) $\gamma F_n, \psi F_k \subseteq F_m = 1_G \cdot F_m$.
\end{proof}

Notice that we provided an explicit construction of the sets $(\Delta_n)_{n \in \N}$ satisfying the lemma.

The previous lemma provides one with an easy way to construct many pre-blueprints. Once $F_{n-1}$ has been defined, one simply chooses sets $\delta^n_k$ for $0 \leq k < n$ satisfying the assumptions of the lemma and then defines $F_n = \bigcup_{0 \leq k < n} \delta^n_k F_k$. In the end one will have collections of sets satisfying the assumptions of the lemma.

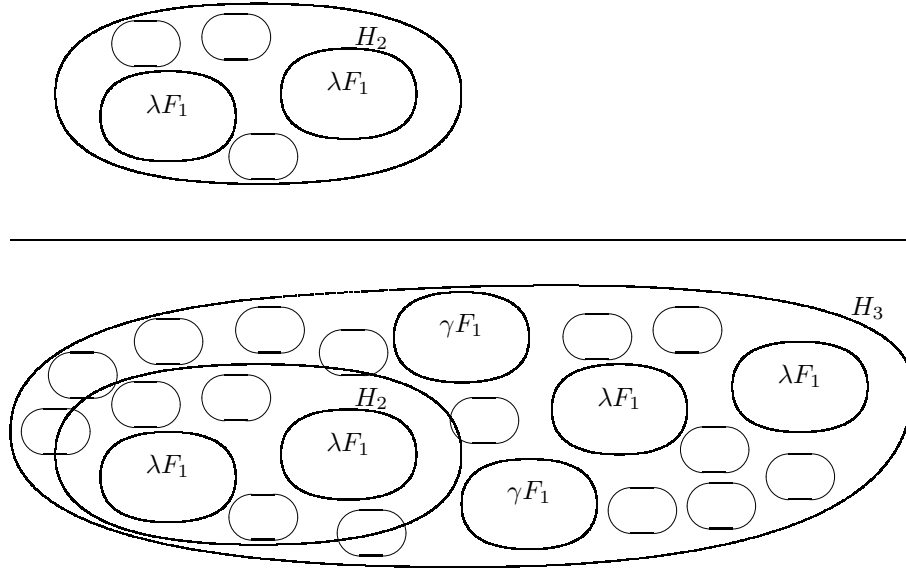
\begin{figure}[ht]
\begin{center}
\setlength{\unitlength}{3mm}
\begin{picture}(42,28)(0,-0.5)

\put(0,14.5){\line(1,0){42}}

\put(0,15){
\put(2,2){

\put(4,6.2){\oval(3,2)}
\put(8,6.5){\oval(3,2)}
\put(9.2,1.2){\oval(3,2)}

\put(14,6){\makebox(0,0)[b]{$H_2$}}
\qbezier(0,4)(0,0)(9,0)
\qbezier(0,4)(0,8)(9,8)
\qbezier(9,0)(18,0)(18,4)
\qbezier(9,8)(18,8)(18,4)}

\put(4,3){
\put(3,2){\makebox(0,0)[b]{$\lambda F_1$}}
\qbezier(0,2)(0,0)(3,0)
\qbezier(0,2)(0,4)(3,4)
\qbezier(3,0)(6,0)(6,2)
\qbezier(3,4)(6,4)(6,2)}

\put(12,4){
\put(3,2){\makebox(0,0)[b]{$\lambda F_1$}}
\qbezier(0,2)(0,0)(3,0)
\qbezier(0,2)(0,4)(3,4)
\qbezier(3,0)(6,0)(6,2)
\qbezier(3,4)(6,4)(6,2)}
}

\put(28,2.5){\oval(3,2)}
\put(31.5,2.7){\oval(3,2)}
\put(35,4){\oval(3,2)}
\put(21,6.5){\oval(3,2)}
\put(16,1.5){\oval(3,2)}
\put(2,6){\oval(3,2)}
\put(3.2,8.5){\oval(3,2)}
\put(7,10){\oval(3,2)}
\put(11.5, 10.5){\oval(3,2)}
\put(15.2,9.5){\oval(3,2)}

\put(17,8.2){
\put(3,2){\makebox(0,0)[b]{$\gamma F_1$}}
\qbezier(0,2)(0,0)(3,0)
\qbezier(0,2)(0,4)(3,4)
\qbezier(3,0)(6,0)(6,2)
\qbezier(3,4)(6,4)(6,2)}

\put(20,0.8){
\put(3,2){\makebox(0,0)[b]{$\gamma F_1$}}
\qbezier(0,2)(0,0)(3,0)
\qbezier(0,2)(0,4)(3,4)
\qbezier(3,0)(6,0)(6,2)
\qbezier(3,4)(6,4)(6,2)}

\put(0,-1){
\put(2,2){

\put(4,6.2){\oval(3,2)}
\put(8,6.5){\oval(3,2)}
\put(9.2,1.2){\oval(3,2)}

\put(14,6){\makebox(0,0)[b]{$H_2$}}
\qbezier(0,4)(0,0)(9,0)
\qbezier(0,4)(0,8)(9,8)
\qbezier(9,0)(18,0)(18,4)
\qbezier(9,8)(18,8)(18,4)}

\put(4,3){
\put(3,2){\makebox(0,0)[b]{$\lambda F_1$}}
\qbezier(0,2)(0,0)(3,0)
\qbezier(0,2)(0,4)(3,4)
\qbezier(3,0)(6,0)(6,2)
\qbezier(3,4)(6,4)(6,2)}

\put(12,4){
\put(3,2){\makebox(0,0)[b]{$\lambda F_1$}}
\qbezier(0,2)(0,0)(3,0)
\qbezier(0,2)(0,4)(3,4)
\qbezier(3,0)(6,0)(6,2)
\qbezier(3,4)(6,4)(6,2)}
}

\put(20,2){
\put(2,2){

\put(4,6.2){\oval(3,2)}
\put(8,6.5){\oval(3,2)}
\put(9.2,1.2){\oval(3,2)}

}

\put(4,3){
\put(3,2){\makebox(0,0)[b]{$\lambda F_1$}}
\qbezier(0,2)(0,0)(3,0)
\qbezier(0,2)(0,4)(3,4)
\qbezier(3,0)(6,0)(6,2)
\qbezier(3,4)(6,4)(6,2)}

\put(12,4){
\put(3,2){\makebox(0,0)[b]{$\lambda F_1$}}
\qbezier(0,2)(0,0)(3,0)
\qbezier(0,2)(0,4)(3,4)
\qbezier(3,0)(6,0)(6,2)
\qbezier(3,4)(6,4)(6,2)}
}

\qbezier(0,6)(0,0)(20,0)
\qbezier(20,0)(40,0)(40,8)
\qbezier(40,8)(40,14)(12,12)
\qbezier(12,12)(0,11)(0,6)
\put(38,11){\makebox(0,0)[b]{$H_3$}}

\end{picture}
\caption{\label{fig:F7}A scenario when the greedy algorithm fails to produce a maximally disjoint family. The upper half of the figure illustrates the construction of $F_2$ by the greedy algorithm: first fill $H_2$ with a maximally disjoint family of translates of $F_1$ (generically marked as $\lambda F_1$), and then fill the remaining part of $H_2$ with a maximally disjoint family of translates of $F_0$ (unmarked). In the lower half of the figure, $F_3$ is constructed similarly, starting with a maximally disjoint family of translates of $F_2$ in $H_3$ (note the translate of $F_2$ on the right). Apparently the resulting collection of translates of $F_1$ is not maximally disjoint.}

\end{center}
\end{figure}

Pre-blueprints are easy to construct, but a nontrivial question is how to modify these methods to construct a blueprint. By conclusion (ii) of Lemma \ref{STRONG BP LIST} we know that it suffices to make each of the $\Delta_n$-translates of $F_n$ maximally disjoint within $G$. It can be seen that in the previous lemma $\Delta_k F_k \subseteq \bigcup_{n \geq k} F_n$ for every $k \in \N$. Since the $F_n$'s are increasing (in the construction we have) we need $F_n$ to come close to exhausting the entire group as $n \rightarrow \infty$. One way to do this is to fix an increasing sequence $(H_n)_{n \in \N}$ with $\bigcup_{n \in \N} H_n = G$ and for each $n \in \N$ try to construct $F_n$ so that it comes close to filling up all of $H_n$. A likely belief is that in order to make the pre-blueprint maximally disjoint we need to not only use the $H_n$'s, but also when constructing $F_n$ the set $\delta^n_{n-1}$ should be chosen so that the $\delta^n_{n-1}$-translates of $F_{n-1}$ are contained and maximally disjoint within $H_n$, then $\delta^n_{n-2}$ should be chosen so that the $\delta^n_{n-2}$-translates of $F_{n-2}$ are contained and maximally disjoint within the space that remains, and carry this on all the way to $\delta^n_0$.
We refer to this approach as the {\it greedy algorithm}. The idea might be that since $\Delta_k$ is all of the explicit and implicit translates of $F_k$ used during the process and since these translates were always chosen to be maximally disjoint within the available space, the $\Delta_k$-translates of $F_k$ should be maximally disjoint. However, this is not the case. With the greedy algorithm, the $\Delta_0$-translates of $F_0$ will definitely by maximally disjoint within $G$, but this is not necessarily the case for the $\Delta_n$-translates of $F_n$ for $n > 0$ (see Figure~\ref{fig:F7} for an illustration of the potential problem). In fact, the situation can be so bad that for all finite sets $A \subseteq G$, $\Delta_1 A \neq G$.

This approach to constructing a blueprint is salvagable with a more careful implementation. Choosing $\delta^n_{n-1}$ so that its translates of $F_{n-1}$ are contained and maximally disjoint within $H_n$ is the right thing to do, but with $\delta^n_{n-2}$ we should be more careful. It is likely that in order for the $\Delta_{n-1}$-translates of $F_{n-1}$ to be maximally disjoint, a translate of $F_{n-1}$ must be used which intersects $H_n$ but is not contained in it (as is the case of Figure~\ref{fig:F7} for $n=2$). If $\delta^n_{n-2}$-translates of $F_{n-2}$ fill up too much of the ``boundary'' of $H_n$, then this will be a problem. Also, we could have the $\delta^n_{n-2}$-translates of $F_{n-2}$ fill up a lot of the ``boundary'' of $F_n$, so the problem could be made worse when constructing $F_{n+1}$. (Translates of $F_{n-2}$ could again fill up the boundary of $H_{n+1}$, and just after these $F_{n-2}$'s could be translates of $F_n$ which therefore also have translates of $F_{n-2}$ making up their boundary. The translates of $F_{n-2}$ could therefore fill up an even thicker portion of the boundary of $H_{n+1}$.) So the idea is we should make sure there is a buffer between the $\delta^n_{n-2}$-translates of $F_{n-2}$ and the complement of $H_n$. By similar reasoning, we should keep the $\delta^n_{n-3}$-translates of $F_{n-3}$ away from the boundary of the region where we were placing the $\delta^n_{n-2}$-translates of $F_{n-2}$. We put this idea in place after the following definition.

\begin{definition} \label{DEFN GROWTH} \index{growth sequence}
Let $G$ be a countably infinite group. A \emph{growth sequence} is a sequence $(H_n)_{n \in \N}$ of finite subsets of $G$ satisfying:
\begin{enumerate}
\item[\rm (i)] $1_G \in H_0$;
\item[\rm (ii)] $\bigcup_{n \in \N} H_n = G$;
\item[\rm (iii)] $H_{n-1} (H_0^{-1} H_0) (H_1^{-1} H_1) \cdots (H_{n-1}^{-1} H_{n-1}) \subseteq H_n$ for each $n \geq 1$;
\item[\rm (iv)] for each $n \geq 1$, if $\Delta \subseteq H_n$ has the property that $g H_{n-1} \cap \Delta H_{n-1} \neq \varnothing$ whenever $g H_{n-1} \subseteq H_n$, then $|\Delta| \geq 3$.
\end{enumerate}
\end{definition}

It is easy to construct a sequence $(H_n)_{n \in \N}$ satisfying (i), (ii), and (iii). Condition (iv) is not difficult to satisfy either, but might not be as obvious. Condition (iv) will be studied more in the next section.

\begin{theorem} \label{EXIST STRONG BP}
Let $G$ be a countably infinite group and let $(H_n)_{n \in \N}$ be a growth sequence. Then there is a maximally disjoint, centered, directed blueprint $(\Delta_n, F_n)_{n \in \N}$ satisfying
\begin{enumerate}
\item[\rm (i)] $F_0 = H_0$;
\item[\rm (ii)] $F_n \subseteq H_n$ for all $n \geq 1$;
\item[\rm (iii)] for all $n \geq 1$ the $D^n_{n-1}$-translates of $F_{n-1}$ are contained and maximally disjoint within $H_n$;
\item[\rm (iv)] for all $n \geq 1$ and $0 \leq k < n$ the $D^n_k$-translates of $F_k$ are maximally disjoint within $H_{n-1}$.
\end{enumerate}
\end{theorem}

\begin{proof}
Set $F_0 = H_0$ so (i) is satisfied. Suppose $F_0$ through $F_{n-1}$ have been constructed with each $F_i \subseteq H_i$. Choose $\delta^n_{n-1}$ so that $1_G \in \delta^n_{n-1}$ and the $\delta^n_{n-1}$-translates of $F_{n-1}$ are contained and maximally disjoint within $H_n$. Note that by the definition of a growth sequence we must have $|\delta^n_{n-1}| \geq 3$. Once $\delta^n_{n-1}$ through $\delta^n_{k+1}$ have been defined with $0 \leq k < n-1$, let $\delta^n_k$ be such that the $\delta^n_k$-translates of $F_k$ are contained and maximally disjoint within
$$B^n_k - \bigcup_{k < m < n} \ \bigcup_{\gamma\in \delta_m^n} \gamma F_m = B^n_k - \bigcup_{k < m < n} \delta_m^n F_m$$
where
$$B^n_k = \{g \in G \: \{g\} (F_{k+1}^{-1} F_{k+1}) (F_{k+2}^{-1} F_{k+2}) \cdots (F_{n-1}^{-1} F_{n-1}) \subseteq H_n\}.$$
Note that $H_{n-1} \subseteq B^n_k$, so $B^n_k \neq \varnothing$. Finally, define
$$F_n = \bigcup_{0 \leq k < n} \delta^n_k F_k.$$
Clearly $F_n \subseteq H_n$. See Figure \ref{fig:F6} for an illustration of the construction of $F_n$.


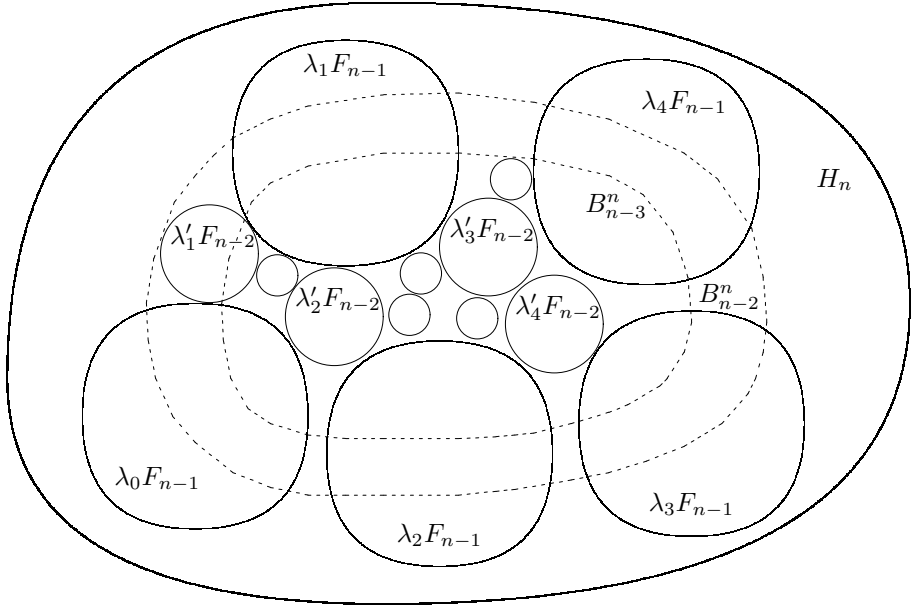
\begin{figure}[ht]
\begin{center}
\setlength{\unitlength}{5mm}
\begin{picture}(28,17)(0,0)

\put(22,11){\makebox(0,0)[b]{$H_n$}}
\qbezier(0,6)(0,16)(10,16)
\qbezier(0,6)(0,0)(10,0)
\qbezier(10,0)(24,0)(24,8)
\qbezier(10,16)(24,16)(24,8)

\linethickness{.08mm}
\put(-4,-4){
\put(8,7){\makebox(0,0)[b]{$\lambda_0 F_{n-1}$}}
\qbezier(6,9)(6,6)(9,6)
\qbezier(6,9)(6,12)(9,12)
\qbezier(9,6)(12,6)(12,9)
\qbezier(9,12)(12,12)(12,9)}

\put(0,3){
\put(9,11){\makebox(0,0)[b]{$\lambda_1 F_{n-1}$}}
\qbezier(6,9)(6,6)(9,6)
\qbezier(6,9)(6,12)(9,12)
\qbezier(9,6)(12,6)(12,9)
\qbezier(9,12)(12,12)(12,9)}

\put(8,2.5){
\put(10,10.5){\makebox(0,0)[b]{$\lambda_4 F_{n-1}$}}
\qbezier(6,9)(6,6)(9,6)
\qbezier(6,9)(6,12)(9,12)
\qbezier(9,6)(12,6)(12,9)
\qbezier(9,12)(12,12)(12,9)}

\put(2.5,-5){
\put(9,6.5){\makebox(0,0)[b]{$\lambda_2 F_{n-1}$}}
\qbezier(6,9)(6,6)(9,6)
\qbezier(6,9)(6,12)(9,12)
\qbezier(9,6)(12,6)(12,9)
\qbezier(9,12)(12,12)(12,9)}

\put(9.2,-4.2){
\put(9,6.5){\makebox(0,0)[b]{$\lambda_3 F_{n-1}$}}
\qbezier(6,9)(6,6)(9,6)
\qbezier(6,9)(6,12)(9,12)
\qbezier(9,6)(12,6)(12,9)
\qbezier(9,12)(12,12)(12,9)}

\begin{dashjoin}{0.1}
\jput(5,4.3){}
\jput(6,3.5){}
\jput(7,3.1){}
\jput(8,2.9){}
\jput(9,2.9){}
\jput(10,2.9){}
\jput(12,2.9){}
\jput(13,3){}
\jput(13.5,3.05){}
\jput(15,3.3){}
\jput(16,3.5){}
\jput(17,3.8){}
\jput(18,4.2){}
\jput(19,4.8){}
\jput(19.7, 5.6){}
\jput(19.9,6){}
\jput(20.1,6.7){}
\jput(20.2,7.5){}
\jput(20.1,8.7){}
\jput(19.8,10){}
\jput(18.9,11.3){}
\jput(18,12.){}
\jput(16,12.9){}
\jput(14,13.35){}
\jput(12,13.6){}
\jput(10,13.55){}
\jput(8,13.25){}
\jput(7,12.9){}
\jput(6,12.4){}
\jput(5.5,11.95){}
\jput(5,11.4){}
\jput(4.45,10.7){}
\jput(3.9,9.2){}
\jput(3.7,8){}
\jput(3.75,7){}
\jput(3.95,6){}
\jput(4.4,5){}
\jput(5,4.3){}
\end{dashjoin}

\put(19.2,7.8){\makebox(0,0)[b]{$B^n_{n-2}$}}
\put(16.2,10.2){\makebox(0,0)[b]{$B^n_{n-3}$}}

\put(7.7,6.65){
\put(1.1, 1.1){\makebox(0,0)[b]{$\lambda'_2 F_{n-2}$}}
\put(1,1){\circle{2.6}}}

\put(4.38,8.32){
\put(1.1, 1.1){\makebox(0,0)[b]{$\lambda'_1 F_{n-2}$}}
\put(1,1){\circle{2.6}}}

\put(11.8,8.5){
\put(1.1, 1.1){\makebox(0,0)[b]{$\lambda'_3 F_{n-2}$}}
\put(1,1){\circle{2.6}}}

\put(13.55,6.45){
\put(1.1, 1.1){\makebox(0,0)[b]{$\lambda'_4 F_{n-2}$}}
\put(1,1){\circle{2.6}}}

\begin{dashjoin}{0.1}
\jput(7,4.8){}
\jput(8,4.5){}
\jput(9,4.4){}
\jput(10,4.4){}
\jput(12,4.4){}
\jput(13,4.45){}
\jput(13.5,4.5){}
\jput(14,4.55){}
\jput(15,4.8){}
\jput(16,5.1){}
\jput(16.5,5.3){}
\jput(17.35, 5.8){}
\jput(17.55,6){}
\jput(18,6.7){}
\jput(18.2,7.5){}
\jput(18.1,8.7){}
\jput(17.7,9.9){}
\jput(16.9,10.9){}
\jput(16,11.4){}
\jput(14,11.85){}
\jput(12,12){}
\jput(10,12.){}
\jput(8,11.65){}
\jput(7,11.2){}
\jput(6.5,10.7){}
\jput(5.9,9.2){}
\jput(5.7,8){}
\jput(5.75,7){}
\jput(5.95,6){}
\jput(6.4,5.2){}
\jput(7,4.8){}
\end{dashjoin}

\put(12.5,7.6){\circle{1.1}}
\put(7.18,8.75){\circle{1.1}}
\put(13.4, 11.3){\circle{1.1}}
\put(11,8.8){\circle{1.1}}
\put(10.7,7.7){\circle{1.1}}

\end{picture}
\caption{\label{fig:F6} The construction of $F_n$. The figure shows the first three steps of the construction. In the first step a maximal disjoint family of translates of $F_{n-1}$ within $H_n$ is selected. In the second step a maximal disjoint family of translates of $F_{n-2}$ within the available part of $B^n_{n-2}$ is selected.
In the third step a maximal disjoint family of translates of $F_{n-3}$ (smallest circles without labels in the figure) within the available part of $B^n_{n-3}$ is selected. This induction process is repeated $n$ times.}
\end{center}
\end{figure}

The $F_n$'s and $\delta^n_k$'s satisfy the assumptions of Lemma \ref{MAKE PREBP}. So if $(\Delta_n)_{n \in \N}$ is as defined in the proof of that lemma, then $(\Delta_n, F_n)_{n \in \N}$ is a centered and directed pre-blueprint. Clearly this pre-blueprint satisfies (i), (ii), and (iii).

The use of the $B^n_k$'s is the key ingredient in this proof. Intuitively, they create the buffer we discussed prior to this theorem. In other words, $B^n_k$ buffers the $\delta^n_k$-translates of $F_k$ from ``the boundary'' of $F_n$. The smaller the value of $k$, the larger the buffer. If a translate of $F_k$ meets $F_n$, then the $\delta^n_t$-translates of $F_t$ for $t < k$ are kept so close to the center of $F_n$ that this translate of $F_k$ cannot meet $\delta^n_t F_t$ for $t < k$ without also meeting $\delta^n_m F_m$ for some $k \leq m < n$. This is formalized in (1) below.

We proceed to verify the following three facts.
\begin{enumerate}
\item[\rm (1)] If $n > k$, $g \in G$, and $g F_k \cap F_n \neq \varnothing$, then $g F_k \cap \delta^n_m F_m \neq \varnothing$ for some $k \leq m < n$;
\item[\rm (2)] $g F_k \cap F_n \neq \varnothing \Longrightarrow g F_k \cap D^n_k F_k \neq \varnothing$ for all $g \in G$ and $k \leq n$;
\item[\rm (3)] the $D_k^n$-translates of $F_k$ are maximally disjoint within $B^n_k$ for all $n, k \in \N$ with $k < n$.
\end{enumerate}

(Proof of 1). It is important to note we require $m \geq k$ as otherwise the claim would be trivial. Let $n > k$ and $g \in G$ satisfy $g F_k \cap F_n \neq \varnothing$. It suffices to show that if $g F_k \cap \delta_m^n F_m = \varnothing$ for all $k < m < n$ then $g F_k \cap \delta_k^n F_k \neq \varnothing$ (since this will validate the claim with $m = k$). As $F_n = \bigcup_{0 \leq t < n} \delta_t^n F_t$, there is $0 \leq t \leq k$ with $g F_k \cap \delta_t^n F_t \neq \varnothing$. If $t = k$, then we are done. So suppose $t < k$. We have
$$g F_k \subseteq \delta_t^n F_t  F_k^{-1} F_k \subseteq \delta_t^n F_t (F_{t+1}^{-1} F_{t+1}) (F_{t+2}^{-1} F_{t+2}) \cdots (F_k^{-1} F_k).$$
Hence
$$g F_k (F_{k+1}^{-1} F_{k+1}) \cdots (F_{n-1}^{-1} F_{n-1}) \subseteq \delta_t^n F_t (F_{t+1}^{-1} F_{t+1}) \cdots (F_{n-1}^{-1} F_{n-1}).$$
However, by definition $\delta_t^n F_t \subseteq B^n_t$. So the right hand side of the expression above is contained within $H_n$, and therefore $g F_k \subseteq B^n_k$. Thus
$$g F_k \subseteq B^n_k - \bigcup_{k < m < n} \delta_m^n F_m.$$
It now follows from the definition of $\delta_k^n$ that $g F_k \cap \delta_k^n F_k \neq \varnothing$. This substantiates our claim.

(Proof of 2). Fix $k \in \N$. If $n = k$ then the claim is clear. Now assume the claim is true for all $k \leq m < n$. Let $g \in G$ satisfy $g F_k \cap F_n \neq \varnothing$. By (1) we have that $g F_k \cap \delta^n_m F_m \neq \varnothing$ for some $k \leq m < n$. Let $\gamma \in \delta^n_m$ be such that $g F_k \cap \gamma F_m \neq \varnothing$. Then $\gamma^{-1} g F_k \cap F_m \neq \varnothing$, so by the induction hypothesis $\gamma^{-1} g F_k \cap D^m_k F_k \neq \varnothing$. By the definition of $D^n_k$ we have $\gamma D^m_k \subseteq \delta^n_m D^m_k \subseteq D^n_k$. Thus, $g F_k \cap D^n_k F_k \neq \varnothing$. The claim now follows from induction.

(Proof of 3). Fix $k < n$ and let $g \in G$ be such that $g F_k \subseteq B^n_k$. We must show $g F_k \cap D^n_k F_k \neq \varnothing$. We are done if $g F_k \cap \delta^n_k F_k \neq \varnothing$ since $\delta^n_k = \delta^n_k D^k_k \subseteq D^n_k$. So suppose $g F_k \cap \delta^n_k F_k = \varnothing$. Recall that in the construction of $F_n$ we defined $\delta_{n-1}^n$ through $\delta_{k+1}^n$ first and then chose $\delta_k^n$ so that its translates of $F_k$ would be maximally disjoint within $B^n_k-\bigcup_{k<m<n} \delta_m^n F_m$. Thus we cannot have $g F_k \subseteq B^n_k-\bigcup_{k<m<n} \delta_m^n F_m$ as this would violate the definition of $\delta_k^n$. Since $g F_k \subseteq B^n_k$, we must have $g F_k \cap (\bigcup_{k<m<n} \delta_m^n F_m) \neq \varnothing$. Let $k < m < n$ and $\gamma \in \delta^n_m$ be such that $g F_k \cap \gamma F_m \neq \varnothing$. Then $\gamma^{-1} g F_k \cap F_m \neq \varnothing$ and thus $\gamma^{-1} g F_k \cap D^m_k F_k \neq \varnothing$ by (2). Now we have $\gamma D^m_k \subseteq \delta^n_m D^m_k \subseteq D^n_k$ so that $g F_k \cap D^n_k F_k \neq \varnothing$. This completes the proof of (3).

Considering (3), we have that in particular the $D^n_k$-translates of $F_k$ are maximally disjoint within (though likely not contained in) $H_{n-1}$ since $H_{n-1} \subseteq B^n_k$. This establishes (iv). Since $\bigcup_{n \in \N} H_n = G$, it follows that the $\Delta_k$-translates of $F_k$ are maximally disjoint within $G$. By clause (ii) of Lemma \ref{STRONG BP LIST}, $(\Delta_n, F_n)_{n \in \N}$ is a maximally disjoint blueprint.
\end{proof}

The previous theorem motivates the following definition.

\begin{definition} \label{DEFN BP GUIDE}
Let $(\Delta_n, F_n)_{n \in \N}$ be a blueprint and let $(H_n)_{n \in \N}$ be a growth sequence. We say the blueprint $(\Delta_n, F_n)_{n \in \N}$ is \emph{guided} by the growth sequence $(H_n)_{n \in \N}$ if the numbered clauses of Theorem \ref{EXIST STRONG BP} are satisfied. Specifically, if:
\begin{enumerate}
\item[\rm (i)] $F_0 = H_0$;
\item[\rm (ii)] $F_n \subseteq H_n$ for all $n \geq 1$;
\item[\rm (iii)] for all $n \geq 1$ the $D^n_{n-1}$-translates of $F_{n-1}$ are contained and maximally disjoint within $H_n$;
\item[\rm (iv)] for all $n \geq 1$ and $0 \leq k < n$ the $D^n_k$-translates of $F_k$ are maximally disjoint within $H_{n-1}$.
\end{enumerate}
\end{definition}

Notice that the blueprint in the previous definition is not required to be centered, maximally disjoint, nor directed. However, we do have the following.

\begin{lem} \label{LEM GUIDED BP}
Let $G$ be a countably infinite group and let $(\Delta_n, F_n)_{n \in \N}$ be a blueprint guided by a growth sequence $(H_n)_{n \in \N}$. Then
\begin{enumerate}
\item[\rm (i)] If $(\Delta_n, F_n)_{n \in \N}$ is centered, then it is directed and maximally disjoint within $G$;
\item[\rm (ii)] $H_n \subseteq F_{n+2} F_0^{-1}$ for all $n \in \N$;
\item[\rm (iii)] $\psi F_k \cap \gamma H_n \neq \varnothing \Longrightarrow \psi F_k \subseteq \gamma H_{n+1} \Longrightarrow \psi F_k \subseteq \gamma F_{n+2}$, for all $n \geq k$, $\psi \in \Delta_k$, and $\gamma \in \Delta_{n+2}$;
\item[\rm (iv)] $\gamma h \in \Delta_k B \Longleftrightarrow \sigma h \in \Delta_k B$, for all $n \geq k$, $h \in H_n$, $\gamma, \sigma \in \Delta_{n+2}$, and $B \subseteq F_k$.
\end{enumerate}
\end{lem}

\begin{proof}
(i). Since the blueprint is centered $D^n_k = 1_G D^n_k \subseteq \Delta_k$ by conclusion (i) of Lemma \ref{BP LIST}. Therefore $\Delta_k$ is maximally disjoint within $H_{n-1}$ for all $n > k$ by clause (iv) of Definition \ref{DEFN BP GUIDE}. Since $\bigcup_{n > k} H_{n-1} = G$ by clauses (ii) and (iii) of Definition \ref{DEFN GROWTH}, it follows that the $\Delta_k$-translates of $F_k$ are maximally disjoint within $G$. Now let $\psi_1, \psi_2 \in \Delta_k$. Then $\psi_1 F_k \cup \psi_2 F_k \subseteq H_{n-1}$ for some $n > k$. So $\psi_1 F_k$ and $\psi_2 F_k$ must meet $D^n_k F_k$ by clause (iv) of Definition \ref{DEFN BP GUIDE}. However $D^n_k \subseteq \Delta_k$, so it must be that $\psi_1, \psi_2 \in D^n_k = 1_G D^n_k \subseteq 1_G F_n$. Thus $\psi_1 F_k \cup \psi_2 F_k \subseteq 1_G F_n$ and $(\Delta_n, F_n)_{n \in \N}$ is directed.

(ii). If $g \in H_n$ then $g F_0 \subseteq H_n H_0 \subseteq H_{n+1}$. So by clause (iv) of Definition \ref{DEFN BP GUIDE}, $g F_0 \cap D^{n+2}_0 F_0 \neq \varnothing$. By conclusion (ii) of Lemma \ref{BP LIST}, $D^{n+2}_0 F_0 \subseteq F_{n+2}$. Thus $g F_0 \cap F_{n+2} \neq \varnothing$ and $g \in F_{n+2} F_0^{-1}$.

(iii). If $\psi F_k \cap \gamma H_n \neq \varnothing$ then $\psi F_k \subseteq \gamma H_n F_k^{-1} F_k \subseteq \gamma H_{n+1}$. By clause (iv) of Definition \ref{DEFN BP GUIDE}, we have $\psi F_k \cap \gamma D^{n+2}_k F_k \neq \varnothing$. However, $\gamma D^{n+2}_k \subseteq \Delta_k$ by conclusion (i) of Lemma \ref{BP LIST} and so it must be that $\psi \in \gamma D^{n+2}_k$. It then follows that $\psi F_k \subseteq \gamma D^{n+2}_k F_k \subseteq \gamma F_{n+2}$ by conclusion (ii) of Lemma \ref{BP LIST}.

(iv). Suppose $\psi \in \Delta_k$ and $\gamma h \in \psi B$. Then $\gamma h \in \psi F_k \cap \gamma H_n$, so by (iii) $\psi F_k \subseteq \gamma F_{n+2}$ and hence $\psi \in \gamma D^{n+2}_k$. It follows $\sigma \gamma^{-1} \psi \in \sigma D^{n+2}_k \subseteq \Delta_k$ and
$$\sigma h = \sigma \gamma^{-1} \gamma h \in \sigma \gamma^{-1} \psi B \subseteq \Delta_k B.$$
\end{proof}

Centered blueprints guided by a growth sequence are centered, maximally disjoint, directed, and on top of this the close relationship between the blueprint and the growth sequence is quite useful as well. These blueprints are the strongest type of blueprints which we know exist for every countably infinite group.

We end this section with a quick application of blueprints. We do not know a proof of this theorem which does not use blueprints. The theorem therefore appears to be nontrivial.

\begin{theorem} \label{MIN DENSE}
Let $G$ be a countable group. Then the set of minimal elements of $2^G$ is dense.
\end{theorem}

\begin{proof}
If $G$ is finite then every element of $2^G$ is minimal. So suppose $G$ is countably infinite and let $1_G = g_0, g_1, g_2, \ldots$ be the enumeration of $G$ used in defining the metric $d$ on $2^G$. Let $x \in 2^G$ and let $\epsilon > 0$. Let $r \in \N$ be such that $2^{-r} < \epsilon$, and set $A = \{g_0, g_1, \ldots, g_r\}$. Let $(\Delta_n, F_n)_{n \in \N}$ be a directed maximally disjoint blueprint with $A \subseteq F_0$ (use Theorem \ref{EXIST STRONG BP}).

Define $y \in 2^G$ by
$$y(g) = \begin{cases}
x(a) & \text{if } \gamma \in \Delta_0, \ a \in A, \text{ and } g = \gamma a \\
0 & \text{otherwise}
\end{cases}.$$
Since the $\Delta_0$-translates of $F_0$ are disjoint and $A \subseteq F_0$, $y$ is well defined. Also, we have $d(x, \gamma^{-1} \cdot y) < \epsilon$ for any $\gamma \in \Delta_0$. It remains to show that $y$ is minimal (in which case $\gamma^{-1} \cdot y$ is minimal as well). Let $B \subseteq G$ be finite. By conclusion (vi) of Lemma \ref{STRONG BP LIST} there is a finite $T \subseteq G$ so that for any $g \in G$ there is $t \in T$ such that
$$\forall b \in B A^{-1} \ (gtb \in \Delta_0 \Longleftrightarrow b \in \Delta_0).$$
Let $g \in G$ be arbitrary, and let $t \in T$ be such that $gtb \in \Delta_0$ if and only if $b \in \Delta_0$ for every $b \in B A^{-1}$. If $b \in B$ and $gtb = \gamma a$ for some $\gamma \in \Delta_0$ and $a \in A$, then $gtba^{-1} = \gamma \in \Delta_0$. Hence $ba^{-1} \in \Delta_0$ and $b \in \Delta_0 a$. Similarly, if $b \in \Delta_0 a$ then $gtb \in \Delta_0 a$. It follows that $y(gtb) = y(b)$ for all $b \in B$. Thus $y$ is minimal by Lemma \ref{lem:minimallemma}.
\end{proof}

\section{Growth of blueprints}

We will soon see that fundamental functions are highly useful. In fact, all forthcoming results rely on these functions. Recall that canonical functions are only partial functions. Their ``free points'' are precisely $\Delta_n \Lambda_n b_{n-1}$ for $n \geq 1$. In order for these functions to be useful, we need to be able to ensure that they have many free points. In other words, we want to be able to make $|\Lambda_n| = |D^n_{n-1}| - 3$ large for every $n \geq 1$. In this section we will achieve this goal in the best possible way. Specifically, we will show that each $\Lambda_n$ can be made as large as possible relative to the size of $F_n$.

Let $G$ be a group and let $A, B \subseteq G$ be finite with $1_G \in A$. Define \index{$\rho(A; B)$}
$$\rho(B; A) = \min \{|D| \: D \subseteq B \text{ and } \forall g \in B \ (g A \subseteq B \Rightarrow g A \cap D A \neq \varnothing)\}.$$
This is well defined since $B$ is finite. We tailored the definition of $\rho$ so that it would have the following properties.

\begin{lem} \label{LEM RHO}
Let $A, B \subseteq G$ be finite with $1_G \in A$.
\begin{enumerate}
\item[(i)] If $\Delta \subseteq B$ and the $\Delta$-translates of $A$ are contained and maximally disjoint within $B$, then $|\Delta| \geq \rho(B; A)$;
\item[(ii)] If $1_G \in A' \subseteq A$ then $\rho(B; A') \geq \rho(B; A)$;
\item[(iii)] If $C \subseteq G$ is finite and $C A^{-1} A \subseteq B$ then $\rho(B;A) \leq \rho(B-C;A) + \rho(C A^{-1} A; A)$;
\end{enumerate}
\end{lem}

\begin{proof}
For finite $A', B' \subseteq G$ with $1_G \in A'$ define
$$S(B';A') = \{ D \: D \subseteq B' \text{ and } \forall g \in B' \ (g A' \subseteq B' \Rightarrow g A' \cap D A' \neq \varnothing)\}.$$
So $\rho(B';A') = \min\{ |D| \: D \in S(B';A')\}$.

(i). If $g \in B$ and $g A \subseteq B$, then since the $\Delta$-translates of $A$ are maximally disjoint within $B$ we have $g A \cap \Delta A \neq \varnothing$. Since the $\Delta$-translates of $A$ are contained in $B$ and $1_G \in A$ we have $\Delta \subseteq B$. Therefore $\Delta \in S(B;A)$ so we have $|\Delta| \geq \rho(B;A)$.

(ii). Let $D \in S(B; A')$ be such that $|D| = \rho(B;A')$. If $g \in B$ and $g A \subseteq B$ then $g A' \subseteq g A \subseteq B$ so $g A' \cap D A' \neq \varnothing$. So
$$\varnothing \neq g A' \cap D A' \subseteq g A \cap D A.$$
Therefore $D \in S(B; A)$ and $\rho(B;A') = |D| \geq \rho(B;A)$.

(iii). Let $D_1 \in S(B-C;A)$ and $D_2 \in S(C A^{-1} A; A)$ be such that $|D_1| = \rho(B-C;A)$ and $|D_2| = \rho(C A^{-1} A; A)$. Set $D = D_1 \cup D_2$. Then $D \subseteq B$. Let $g \in B$ be such that $g A \subseteq B$. We proceed by cases. \underline{Case 1:} $g A \subseteq B-C$. Then we must have $g A \cap D_1 A \neq \varnothing$. In particular, $g A \cap D A \neq \varnothing$. \underline{Case 2:} $g A \cap C \neq \varnothing$. Then $g A \subseteq C A^{-1} A$, so $g A \cap D_2 A \neq \varnothing$ and hence $g A \cap D A \neq \varnothing$. Therefore $D \in S(B;A)$. So we have
$$\rho(B;A) \leq |D| \leq \rho(B-C;A) + \rho(C A^{-1} A; A).$$
\end{proof}

Note that clause (iv) of Definition \ref{DEFN GROWTH} is equivalent to $\rho(H_n; H_{n-1}) \geq 3$. Clauses (i) and (ii) listed in the lemmma above were implicitly used in verifying the growth property of the blueprint constructed in Theorem \ref{EXIST STRONG BP}.

\begin{lem} \label{RATIO}
Let $G$ be an infinite group, and let $A, B \subseteq G$ be finite with $1_G \in A$. For any $\epsilon > 0$ there exists a finite $C \subseteq G$ containing $B$ such that $\rho(C ; A) > \frac{|C|}{|A|}(1 - \epsilon)$.
\end{lem}

\begin{proof}
Let $\Delta \subseteq G$ be countably infinite and such that the $\Delta$-translates of $A A^{-1}$ are disjoint and $\Delta A A^{-1} A \cap B = \varnothing$. Let $\lambda_1, \lambda_2, \ldots$ be an enumeration of $\Delta$. For each $n \geq 1$, define
$$B_n = B \cup \left( \bigcup_{1 \leq k \leq n} \lambda_k A \right).$$
Fix $n \geq 1$, and let $D \subseteq B_n$ be such that $g A \cap D A \neq \varnothing$ whenever $g \in B_n$ with $g A \subseteq B_n$. It follows that for each $1 \leq i \leq n$ there is $d_i \in D$ with $d_i A \cap \lambda_i A \neq \varnothing$. Then
$$d_i \in \lambda_i A A^{-1}.$$
Since the $\Delta$-translates of $A A^{-1}$ are disjoint, the $d_i$'s are all distinct. Additionally, $d_i A \cap B \subseteq \Delta A A^{-1} A \cap B = \varnothing$ so that $\rho(B_n ; A) - n \geq \rho(B ; A)$. Therefore we have
$$\rho(B_n ; A) \frac{|A|}{|B_n|} \geq \frac{n |A| + \rho(B ; A) |A|}{n |A| + |B|}.$$
Clearly as $n$ goes to infinity the fraction on the right goes to 1. So there is $n \geq 1$ with $\rho(B_n ; A) \frac{|A|}{|B_n|} > 1 - \epsilon$ and $\rho(B_n ; A) > \frac{|B_n|}{|A|}(1 - \epsilon)$. Setting $C = B_n$ completes the proof.
\end{proof}

\begin{definition} \index{subexponential growth} \index{polynomial growth}
A function $f : \N \rightarrow \N$ is said to have \emph{subexponential growth} if for every $u > 1$ there is $N \in \N$ so that $f(n) < u^n$ for all $n \geq N$. Similarly, $f: \N \rightarrow \N$ is said to have \emph{polynomial growth} if there are $a, b, d \in \N$ so that for all $n \in \N$ we have $f(n) \leq a \cdot n^d + b$.
\end{definition}

\begin{lem} \label{LEM SUBEXPGRTH}
If $f, g : \N \rightarrow \N$ have subexponential growth, then their product has subexponential growth and the function $h: \N \rightarrow \N$ defined by $h(n) = \max\{f(k) \: k \leq n\}$ has subexponential growth.
\end{lem}

\begin{proof}
If $u > 1$, then there is $N \in \N$ with $f(n) < (\sqrt{u})^n$ and $g(n) < (\sqrt{u})^n$ for all $n \geq N$. It follows $f(n)\cdot g(n) < u^n$ for all $n \geq N$ so $f \cdot g$ has subexponential growth. If the function $h$ is bounded, then the claim is trivial. So suppose $h$ is not bounded. Let $u > 1$ and let $N \in \N$ be such that $f(n) < u^n$ for all $n \geq N$. Since $h$ is not bounded, there is $M > N$ with $h(M) > h(N)$. It follows that for every $n \geq M$ there is $N < k(n) \leq n$ with $h(n) = f(k(n))$. Thus, for $n \geq M$ we have $h(n) = f(k(n)) < u^{k(n)} \leq u^n$. We conclude $h$ has subexponential growth.
\end{proof}

\begin{lem}\label{FAST GROWTH}
Let $G$ be an infinite group, and let $A, B \subseteq G$ be finite with $1_G \in A$. If $f : \N \rightarrow \N$ has subexponential growth, then there exists a finite $C \subseteq G$ containing $B$ such that $2^{\rho(C; A)} > f(|C|)$.
\end{lem}

\begin{proof}
Let $N \in \N$ be such that $2^{\frac{n}{2 |A|}} > f(n)$ for all $n \geq N$. Let $B' \subseteq G$ be a finite set containing $B$ with $|B'| \geq N$. By Lemma \ref{RATIO} there exists a finite $C \subseteq G$ containing $B'$ with $\rho(C ; A) > \frac{1}{2} \frac{|C|}{|A|}$. Then $C \supseteq B$ and as $|C|$ is at least $N$,
$$2^{\rho(C ; A)} > 2^{\frac{|C|}{2 |A|}} > f(|C|).$$
\end{proof}

\begin{definition} \index{blueprint!any prescribed growth} \index{fundamental function!any prescribed growth} \index{canonical function!any prescribed growth} \index{any prescribed growth} \index{blueprint!any prescribed polynomial growth} \index{fundamental function!any prescribed polynomial growth} \index{canonical function!any prescribed polynomial growth} \index{any prescribed polynomial growth}
Fix a countably infinite group $G$, and let $P$ be a property of blueprints on $G$. We say the blueprints with property $P$ \emph{can have any prescribed growth} (or \emph{can have any prescribed polynomial growth}) if for any sequence $(p_n)_{n \geq 1}$ of functions of subexponential growth (respectively polynomial growth) there is a blueprint $(\Delta_n, F_n)_{n \in \N}$ with property $P$ satisfying for each $n \geq 1$
$$|\Lambda_n| \geq \log_2 \ \max(p_n(|F_n|), p_n(|B_n|))$$
where $B_n$ is a finite set satisfying $\Delta_n B_n B_n^{-1} = G$. Similarly, if $P$ is a property of fundamental (or canonical) functions, then we say the collection of fundamental (canonical) functions with property $P$ \emph{can have any prescribed growth} (or \emph{can have any prescribed polynomial growth}) if for any sequence $(p_n)_{n \geq 1}$ of functions of subexponential growth (respectively polynomial growth) there is a function $c$ fundamental (canonical) with respect to a blueprint $(\Delta_n, F_n)_{n \in \N}$ such that $c$ has property $P$ and for each $n \geq 1$
$$|\Theta_n| \geq \log_2 \ \max(p_n(|F_n|), p_n(|B_n|))$$
where $B_n$ is a finite set satisfying $\Delta_n B_n B_n^{-1} = G$.
\end{definition}

In the previous definition, requiring $\Delta_n B_n B_n^{-1} = G$ instead of $\Delta_n B_n = G$ is a significant detail. The reason is that it is possible for $|B_n B_n^{-1}| = |B_n|^2$ and therefore $p_n(|B_n B_n^{-1}|) = p_n(|B_n|^2)$. However, even though $p_n$ has subexponential growth, the function $q_n$ defined by $q_n(k) = p_n(k^2)$ may not (for example $p_n(k) = 2^{\sqrt{k}}$). There is nothing formally wrong with this, but this is the reason why our proofs do not work if the above definition is changed so that $\Delta_n B_n = G$.

\begin{lem} \label{GROWTH SEQ}
Let $G$ be a countably infinite group. If $(p_n)_{n \geq 1}$ is a sequence of functions of subexponential growth and $1_G \in A \subseteq G$ is finite, then there exists a growth sequence $(H_n)_{n \in \N}$ with $H_0 = A$ and $\rho(H_n; H_{n-1}) \geq \log_2 \ p_n(|H_n|)$ for all $n \geq 1$.
\end{lem}

\begin{proof}
Fix a sequence $(A_n)_{n \in \N}$ with $A_0 = A$ and $\bigcup_{n \in \N} A_n = G$. Set $H_0 = A_0 = A$. Now assume $H_0$ through $H_{n-1}$ have been constructed for $n > 0$. Apply the previous lemma to find a finite $H_n \subseteq G$ satisfying
$$H_n \supseteq A_n \cup H_{n-1} (H_0^{-1} H_0) (H_1^{-1} H_1) \cdots (H_{n-1}^{-1} H_{n-1})$$
and $\rho(H_n; H_{n-1}) \geq \max( \log_2 \ p_n(|H_n|), \ 3)$. It is easily checked that $(H_n)_{n \in \N}$ is a growth sequence with the desired property.
\end{proof}

\begin{cor} \label{GROW BP}
Let $G$ be a countably infinite group and let $1_G \in A \subseteq G$ be finite. The blueprints $(\Delta_n, F_n)_{n \in \N}$ on $G$ which are centered, guided by a growth sequence, and have $A \subseteq F_0$ can have any prescribed growth. In particular, the blueprints on $G$ which are centered, directed, maximally disjoint, and have $A \subseteq F_0$ can have any prescribed growth.
\end{cor}

\begin{proof}
Let $(p_n)_{n \geq 1}$ be a sequence of functions of subexponential growth. By Lemma \ref{LEM SUBEXPGRTH}, we may assume that each $p_n$ is nondecreasing. For $n \geq 1$ and $k \in \N$, define $q_n(k) = 8 \cdot p_n(k)$. Then $(q_n)_{n \geq 1}$ is a sequence of functions of subexponential growth. By Lemma \ref{GROWTH SEQ} there is a growth sequence $(H_n)_{n \in \N}$ with $H_0 = A$ and $\rho(H_n; H_{n-1}) \geq \log_2 \ q_n(|H_n|)$ for all $n \geq 1$. Apply Theorem \ref{EXIST STRONG BP} to get a maximally disjoint, centered, directed, blueprint $(\Delta_n, F_n)_{n \in \N}$ guided by $(H_n)_{n \in \N}$. Then $A = F_0$. Set $B_n = F_n$ and notice that $\Delta_n B_n B_n^{-1} = G$. Fix $n \geq 1$. By clause (iii) of Theorem \ref{EXIST STRONG BP} and by clause (i) of Lemma \ref{LEM RHO} we have
$$|D^n_{n-1}| \geq \rho(H_n; F_{n-1}).$$
By clause (ii) of Lemma \ref{LEM RHO} we have
$$\rho(H_n; F_{n-1}) \geq \rho(H_n; H_{n-1}).$$
Therefore
$$|\Lambda_n| = |D^n_{n-1}| - 3 \geq - 3 + \log_2 \ q_n(|H_n|)$$
$$= \log_2 \ p_n(|H_n|) \geq \log_2 \ \max(p_n(|F_n|), p_n(|B_n|)).$$
\end{proof}

\begin{cor} \label{GEN SUBEXP FREE}
Let $G$ be a countably infinite group, let $1_G \in A \subseteq G$ be finite, and let $R: A \rightarrow 2$ be a locally recognizable function. The functions which are canonical with respect to a centered blueprint guided by a growth sequence and which are also compatible with $R$ can have any prescribed growth. In particular, the collection of all fundamental functions can have any prescribed growth.
\end{cor}

\begin{proof}
If $(\Delta_n, F_n)_{n \in \N}$ is a blueprint and $c \in 2^{\subseteq G}$ is canonical with respect to this blueprint, then for each $n \geq 1$ $\Theta_n = \Lambda_n$.
Therefore the claim immediately follows from Corollary \ref{GROW BP} and Theorem \ref{FM}.
\end{proof}  
\chapter{\label{chap:7}Basic Applications of the Fundamental Method}

In this chapter, we finally get to reap some of the benefits of all the hard work which went into the previous chapter. In this chapter we present quick and easy yet satisfying applications of the tools we have developed. This chapter places emphasis on the wide variety of constructions, properties, and proofs which can be created using the tools from the previous chapter. Each section focuses on a specific object from the previous chapter and relies primarily on this object to prove an important and nontrivial result. All of our work will be in the spirit of a general and recurrent procedure in this paper which we refer to as \emph{the fundamental method}. The fundamental method refers to the coordinated use of functions of subexponential growth, locally recognizable functions, blueprints, and fundamental functions in achieving a goal of constructing certain special elements of $2^G$. This chapter will be a first step in convincing the reader that the fundamental method provides tremendous control in constructing special elements of $2^G$.

\section{The uniform $2$-coloring property} \label{SEC UNI 2 COL}

This section focuses on the use of functions of subexponential growth. We begin by proving that every countably infinite group has a $2$-coloring. Shortly afterwards we strengthen this to show that all countably infinite groups have the uniform $2$-coloring property.

\begin{theorem}\label{GEN COL}
Let $G$ be a countably infinite group, let $(\Delta_n, F_n)_{n \in \N}$ be a blueprint, and for each $n \geq 1$ let $B_n$ be finite with $\Delta_n B_n B_n^{-1} = G$. If $c \in 2^{\subseteq G}$ is fundamental with respect to this blueprint and $|\Theta_n| \geq \log_2 \ (2|B_n|^4+1)$ for each $n \geq 1$, then $c$ can be extended to a function $c'$ with $|\Theta_n(c')| > |\Theta_n(c)| - \log_2 \ (2|B_n|^4+1) - 1$ such that every $x \in 2^G$ extending $c'$ is a $2$-coloring. In particular, every countable group has a $2$-coloring.
\end{theorem}

\begin{proof}
For each $i \geq 1$, define $\B_i : \N \rightarrow \{0,1\}$ so that $\B_i(k)$ is the $i^\text{th}$ digit from least to most significant in the binary representation of $k$ when $k \geq 2^{i-1}$ and $\B_i(k) = 0$ when $k < 2^{i-1}$. Also, for each $n \geq 1$, let $s(n)$ be the smallest integer greater than or equal to $\log_2 \ (2 |B_n|^4+1)$ and fix any distinct $\theta_1^n, \theta_2^n, \ldots, \theta_{s(n)}^n \in \Theta_n$.

Fix an enumeration $s_1, s_2, \ldots$ of the nonidentity elements of $G$. For each $n \geq 1$, let $\Gamma_n$ be the graph with vertex set $\Delta_n$ and edge relation given by
$$(\gamma, \psi) \in E(\Gamma_n) \Longleftrightarrow \gamma^{-1} \psi \in B_n B_n^{-1} s_n B_n B_n^{-1} \text{ or } \psi^{-1} \gamma \in B_n B_n^{-1} s_n B_n B_n^{-1}$$
for distinct $\gamma, \psi \in \Delta_n$. Then $\deg_{\Gamma_n} (\gamma) \leq 2 |B_n|^4$ for each $\gamma \in \Delta_n$. We can therefore find, via the usual greedy algorithm, a graph-theoretic $(2 |B_n|^4 + 1)$-coloring of $\Gamma_n$, say $\mu_n : \Delta_n \rightarrow \{0, 1, \ldots, 2|B_n|^4\}$.

Define $c' \supseteq c$ by setting
$$c'(\gamma \theta_i^n b_{n-1}) = \B_i (\mu_n(\gamma))$$
for each $n \geq 1$, $\gamma \in \Delta_n$, and $1 \leq i \leq s(n)$. Since $2^{s(n)} \geq 2 |B_n|^4 +1$, all integers $0$ through $2 |B_n|^4$ can be written in binary using $s(n)$ digits. Thus no information is lost between the $\mu_n$'s and $c'$. Setting $\Theta_n(c') = \Theta_n(c) - \{\theta_1^n, \ldots, \theta_{s(n)}^n\}$ we clearly have that $c'$ is fundamental and
$$|\Theta_n(c')| = |\Theta_n(c)| - s(n) > |\Theta_n(c)| - \log_2 \ (2 |B_n|^4+1) - 1.$$

Fix $1_G \neq s \in G$. Then $s = s_n$ for some $n \geq 1$. Let $V$ be the test region for the $\Delta_n$ membership test admitted by $c$. Set $T =  B_n B_n^{-1}(V \cup \Theta_n(c)b_{n-1})$. Now let $x \in 2^G$ be an arbitrary extension of $c'$, and let $g \in G$ be arbitrary. Since $\Delta_n B_n B_n^{-1} = G$, there is $b \in B_n B_n^{-1}$ with $g b = \gamma \in \Delta_n$. We proceed by cases.

\underline{Case 1}: $g s b \not\in \Delta_n$. Since $x \supseteq c$, $g b \in \Delta_n$, and $g s b \not\in \Delta_n$, there is $v \in V$ such that $x(g b v) \neq x(g s b v)$. This completes this case since $b v \in T$.

\underline{Case 2}: $g s b \in \Delta_n$. Then
$$(g b)^{-1} (g s b) = b^{-1} s b \in B_n B_n^{-1} s B_n B_n^{-1}.$$
Thus $(g b, g s b) \in E(\Gamma_n)$ so $\mu_n(g b) \neq \mu_n(g s b)$. Consequently, there is $1 \leq i \leq s(n)$ with $x(g b \theta_i^n b_{n-1}) \neq x(g s b \theta_i^n b_{n-1})$. This completes this case since $b \theta_i^n b_{n-1} \in T$. We conclude $x$ is a 2-coloring.

Now we show that every countable group has a $2$-coloring. As mentioned previously, this is immediate for finite groups. So suppose $G$ is a countably infinite group. Define $p_n(k) = 2 k^4+1$ for each $n \geq 1$ and $k \in \N$. Then $(p_n)_{n \geq 1}$ is a sequence of functions of subexponential growth. By Corollary \ref{GEN SUBEXP FREE}, there is a fundamental function $c \in 2^{\subseteq G}$ with $|\Theta_n| \geq \log_2 \ p_n(|B_n|)$ for all $n \geq 1$. Applying the above construction leads to the conclusion that there is a 2-coloring on $G$.
\end{proof}

Notice that this proof shows that for every $1_G \neq s \in G$ there is a finite set $T \subseteq G$ so that for all $x \in 2^G$ extending $c'$ and all $g \in G$ we have
$$\exists t \in T \ x(gst) \neq x(gt)$$
(the main point here is that $T$ did not depend on the extension $x \in 2^G$). This is actually a general phenomenon as the following proposition shows.

\begin{prop} \label{EXTENSION FREE}
Let $G$ be a countably infinite group, and let $c \in 2^{\subseteq G}$ have the property that every $x \in 2^G$ extending $c$ is a $2$-coloring. Then
\begin{enumerate}
\item[\rm (i)] for every nonidentity $s \in G$ there is a finite $T \subseteq G$ so that for all $g \in G$ there is $t \in T$ with $gt, gst \in \dom(c)$ and $c(gt) \neq c(gst)$;
\item[\rm (ii)] if $E(c) = \{x \in 2^G \: c \subseteq x\}$ is the set of full extensions of $c$, then $\overline{G \cdot E(c)}$ is a free subflow of $2^G$.
\end{enumerate}
\end{prop}

\begin{proof}
(i). Towards a contradiction, suppose there is a nonidentity $s \in G$ such that no finite set $T$ exists satisfying (i). First suppose that there is $g \in G$ such that for all $h \in G$ $c(gh) = c(gsh)$ whenever $gh, ghs \in \dom(c)$. Define $x \in 2^G$ by setting $x(gh) = c(gsh)$ when $gsh \in \dom(c)$, $x(gsh) = c(gh)$ when $gh \in \dom(c)$, and $x(gh) = x(gsh) = 0$ if $gh, gsh \not\in \dom(c)$. Then $x$ is well defined and $s^{-1} \cdot (g^{-1} \cdot x) = g^{-1} \cdot x$. This is a contradiction since $x$ extends $c$ and hence must be a $2$-coloring, in particular must be aperiodic. Now suppose that for every $g \in G$ there is $h \in G$ with $gh, gsh \in \dom(c)$ and $c(gh) \neq c(gsh)$. Let $(A_n)_{n \in \N}$ be an increasing sequence of finite subsets of $G$ with $\bigcup_{n \in \N} A_n = G$. For each $n \in \N$, $s$ and $A_n$ do not satisfy (i) so there is $g_n \in G$ with $c(g_n a) = c(g_n s a)$ whenever $a \in A_n$ and $g_n a, g_n s a \in \dom(c)$. The set $\{g_n \: n \in \N\}$ cannot be finite as otherwise we would be in the first case treated above. Therefore, since the $A_n$'s are increasing we can replace the $g_n$'s with a subsequence if necessary and assume that for $n \neq k \in \N$
$$(g_n A_n \cup g_n s A_n) \cap (g_k A_k \cup g_k s A_k) = \varnothing.$$
Define $y \in 2^{\subseteq G}$ as follows. If $n \in \N$ and $a \in A_n$, set $y(g_n a) = c(g_n s a)$ if $g_n s a \in \dom(c)$ and $y(g_n s a) = c(g_n a)$ if $g_n a \in \dom(c)$. Then $y(h) = c(h)$ whenever $h \in \dom(y) \cap \dom(c)$. So $y \cup c \in 2^{\subseteq G}$. Define $x \in 2^G$ by setting $x(h) = (y \cup c)(h)$ if $h \in \dom(y) \cup \dom(c)$ and $x(h) = 0$ otherwise. For any $h \in G$ there is $n \in \N$ with $h \in A_n$ and hence
$$[s^{-1} \cdot (g_n^{-1} \cdot x)](h) = x(g_n s h) = x(g_n h) = (g_n^{-1} \cdot x)(h).$$
Since the action of $G$ on $2^G$ is continuous, it follows that if $z \in 2^G$ is any limit point of $(g_n^{-1} \cdot x)_{n \in \N}$ then $s^{-1} \cdot z = z$. However, this is a contradiction since $x$ extends $c$ and hence must be a $2$-coloring.

(ii). Fix a nonidentity $s \in G$. Let $T \subseteq G$ be as in (i). Let $x \in E(c)$, let $g \in G$, and let $y = g \cdot x \in G \cdot E(c)$. Then by (i) there is $t \in T$ with
$$y(t) = (g \cdot x)(t) = x(g^{-1} t) \neq x(g^{-1} s t) = (g \cdot x)(st) = y(st) = (s^{-1} \cdot y)(t).$$
By considering the metric $d$ on $2^G$, it follows that there is an $\epsilon > 0$ depending only on $s$ such that for all $y \in G \cdot E(c)$ $d(y, s^{-1} \cdot y) > \epsilon$. By the continuity of the action of $G$ on $2^G$ it follows that if $z \in \overline{G \cdot E(c)}$ then $d(z, s^{-1} \cdot z) \geq \epsilon$. In particular, $s^{-1} \cdot z \neq z$. Since $s \in G - \{1_G\}$ was arbitrary, we conclude that $\overline{G \cdot E(c)}$ is a free subflow.
\end{proof}

\begin{prop} \label{PRE STRONG COL}
Let $G$ be a countably infinite group, $(\Delta_n, F_n)_{n \in \N}$ a blueprint, and $c \in 2^{\subseteq G}$ a fundamental function with $\Theta_n \neq \varnothing$ for all $n \geq 1$. Then $c$ can be extended to a function $x \in 2^G$ with the property that for every nonidentity $s \in G$ there are infinitely many $g \in G$ with $x(g) \neq x(sg)$.
\end{prop}

\begin{proof}
Fix any $\theta_n \in \Theta_n$ for each $n \geq 1$. Enumerate $G-\{1_G\}$ as $s_1, s_2,\dots$ so that every nonidentity group element is enumerated infinitely many times. Inductively
define an increasing sequence $k_n$ of natural numbers as follows. For $n=1$ let $k_1=1$. In
general suppose $k_m$, $1\leq m\leq n$, have all been defined. Then let $k_{n+1}>k_n$ be the least such that
$$\theta_{k_{n+1}} b_{k_{n+1}-1},s_{n+1}\theta_{k_{n+1}}b_{k_{n+1}-1}\not\in \{\theta_{k_m}b_{k_m-1},s_m\theta_{k_m}b_{k_m-1} \: 1\leq m\leq n\}.$$
Such $k_{n+1}$ exists since the set of all $\theta_n b_{n-1}$, $n\geq 1$, is infinite.
This finishes the definition of all $k_n$. As a result, the elements
$$ \theta_{k_1} b_{k_1 -1}, s_1\theta_{k_1} b_{k_1 -1}, \theta_{k_2} b_{k_2 -1}, s_2\theta_{k_2} b_{k_2 -1}, \dots, \theta_{k_n}b_{k_n -1}, s_n\theta_{k_n} b_{k_n -1}, \dots $$
are all distinct.

Since $\theta_{k_n} b_{k_n-1}\not\in \dom(c)$ for all $n\geq 1$, we can clearly extend $c$ to an $x \in 2^G$ such that for all $n\geq 1$, $x(s_n\theta_{k_n} b_{k_n -1})\neq x(\theta_{k_n} b_{k_n-1})$.
\end{proof}

An important theorem will be drawn from the previous proposition momentarily, but first we consider orthogonality.

\begin{prop}\label{GEN ORTH}
Let $G$ be a countably infinite group, and let $c \in 2^{\subseteq G}$ be fundamental with $\Theta_n \neq \varnothing$ for each $n \geq 1$. Then for each $\tau \in 2^\N$ there is a fundamental $c_\tau \supseteq c$ with $|\Theta_n(c_\tau)| = |\Theta_n(c)| - 1$ for each $n \geq 1$ and with the property that if $\tau \neq \sigma \in 2^\N$, $x, y \in 2^G$, $x \supseteq c_\tau$, and $y \supseteq c_\sigma$, then $x$ and $y$ are orthogonal. Additionally, for each $\tau \in 2^\N$ there is $x_\tau \in 2^G$ extending $c_\tau$ such that $\{x_\tau \: \tau \in 2^\N\}$ is a perfect set.
\end{prop}

\begin{proof}
Let $(\Delta_n, F_n)_{n \in \N}$ be the blueprint corresponding to $c$. For each $n \geq 1$ pick $\theta_n \in \Theta_n$. For $\tau \in 2^\N$, we define $c_\tau \supseteq c$ by setting
$$c_\tau (\gamma \theta_n b_{n-1}) = \tau(n-1)$$
for each $n \geq 1$ and $\gamma \in \Delta_n$. If we define $x_\tau \supseteq c_\tau$ by letting $x_\tau$ be zero on $G - \dom(c_\tau)$, then the map $\tau \mapsto x_\tau$ is one-to-one and continuous. Therefore $\{x_\tau \: \tau \in 2^\N\}$ is a perfect set.

Let $B_n$ be finite with $\Delta_n B_n B_n^{-1} = G$ and let $V_n$ be the test region for the $\Delta_n$ membership test admitted by $c$. Set $T_n = B_n B_n^{-1}(V_n \cup \{\theta_n b_{n-1}\})$. Now suppose $\tau \neq \sigma \in 2^\N$, and let $n \geq 1$ satisfy $\tau(n-1) \neq \sigma(n-1)$. Let $x, y \in 2^G$ with $x \supseteq c_\tau$ and $y \supseteq c_\sigma$, and let $g_1, g_2 \in G$ be arbitrary. We will show that there is $t \in T_n$ with $x(g_1 t) \neq y(g_2 t)$. There is $b \in B_n B_n^{-1}$ with $g_1 b \in \Delta_n$. We proceed by cases.

\underline{Case 1}: $g_2 b \not\in \Delta_n$. Since $g_1 b \in \Delta_n$ and $g_2 b \not\in \Delta_n$, there is $v \in V_n$ with $x(g_1 b v) \neq y(g_2 b v)$. This completes this case since $bv \in T_n$.

\underline{Case 2}: $g_2 b \in \Delta_n$. Then $x(g_1 b \theta_n b_{n-1}) = \tau(n-1) \neq \sigma(n-1) = y(g_2 b \theta_n b_{n-1})$. This completes this case as $b \theta_n b_{n-1} \in T_n$.
\end{proof}

Notice that in the previous proof, the set witnessing the orthogonality, $T_n$, depended only on the $n \geq 1$ satisfying $\tau(n-1) \neq \sigma(n-1)$. We will need this fact shortly.

The fact that functions of subexponential growth are closed under multiplication together with the abstract nature of the definition of fundamental functions allows one to easily stack constructions on top of one another, as the next three results demonstrate. 

\begin{theorem}\label{thm:strongcoloring}
Every countably infinite group has a strong $2$-coloring.
\end{theorem}

\begin{proof}
For $n \geq 1$ and $k \in \N$ define $p_n(k) = 2 \cdot 2k^4$. Then $(p_n)_{n \geq 1}$ is a sequence of functions of subexponential growth. By Corollary \ref{GEN SUBEXP FREE}, there is a fundamental $c \in 2^{\subseteq G}$ with
$$|\Theta_n| \geq \log_2 \ (4|B_n|^4+2) = 1 + \log_2 \ (2|B_n|^4+1)$$
for each $n \geq 1$, where $B_n$ satisfies $\Delta_n B_n B_n^{-1} = G$. Now apply Theorem \ref{GEN COL} and Proposition \ref{PRE STRONG COL}, in that order.
\end{proof}

\begin{theorem} \label{THM PSET ORTH COL}
If $G$ is a countably infinite group, then $G$ has the uniform $2$-coloring property. In particular, there is a perfect set of pairwise orthogonal $2$-colorings on $G$.
\end{theorem}

\begin{proof}
The proof is nearly identical to that of the previous theorem. The only difference is to apply Theorem \ref{GEN COL} and Proposition \ref{GEN ORTH}, in that order. This immediately shows that there is a perfect set of pairwise orthogonal $2$-colorings on $G$. The collection of functions constructed, together with the comments immediately following the proofs of Theorem \ref{GEN COL} and Proposition \ref{GEN ORTH}, directly demonstrate that $G$ has the uniform $2$-coloring property.
\end{proof}

\begin{theorem}
If $G$ is a countably infinite group, then there is an uncountable collection of pairwise orthogonal strong $2$-colorings on $G$.
\end{theorem}

\begin{proof}
Same proof as the previous two theorems, except use the functions $p_n(k) = 4(2k^4+1)$. At the end, apply Theorem \ref{GEN COL}, Proposition \ref{GEN ORTH}, and Proposition \ref{PRE STRONG COL}, in that order.
\end{proof}

\section{Density of $2$-colorings}

This section focuses on applications of locally recognizable functions. We begin by revealing just how plentiful these functions are.

\begin{prop}\label{LR EXT}
If $G$ is a countably infinite group, $B \subseteq G$ is finite, and $Q : B \rightarrow 2$ is any function, then there exists a nontrivial locally recognizable function $R : A \rightarrow 2$ extending $Q$.
\end{prop}

\begin{proof}
By defining $Q(1_G) = 0$ if necessary, we may assume $1_G \in B$. Set $B_1 = B$. Choose any $a \neq b \in G - B_1$ and set $B_2 = B_1 \cup \{a, b\}$. Next chose any $c \in G - (B_2 B_2 \cup B_2 B_2^{-1})$ and set $B_3 = B_2 \cup \{c\} = B_1 \cup \{a, b, c\}$. Let $A = B_3 B_3$ and define $R: A \rightarrow 2$ by
$$R(g) = \begin{cases}
Q(g) & \text{if } g \in B_1 \\
Q(1_G) & \text{if } g \in \{a, b, c\} \\
1 - Q(1_G) & \text{if } g \in A - B_3.
\end{cases}$$

We claim $R$ is a locally recognizable function (it is clearly nontrivial). Towards a contradiction, suppose there is $y \in 2^G$ extending $R$ such that for some $1_G \neq g \in A$ $y(g h) = y(h)$ for all $h \in A$. In particular, $y(g) = y(1_G) = R(1_G)$ so $g \in B_3$. We first point out that at least one of $a$, $b$, or $c$ is not an element of $g B_3$. We prove this by cases. \underline{Case 1}: $g \in B_2$. Then $c \not\in g B_2 \subseteq B_2 B_2$ and $c \neq g c$ since $g \neq 1_G$. Thus $c \not\in g B_3$. \underline{Case 2}: $g \in B_3 - B_2 = \{c\}$. Then $g = c$. Since $c \not\in B_2 B_2^{-1}$, $a, b \not\in c B_2$. If $a, b \in c B_3$ then we must have $a = c^2 = b$, contradicting $a \neq b$. We conclude $\{a, b\} \not\subset c B_3 = g B_3$.

The key point now is that $\{a, b, c\} \subseteq \{h \in A \: y(h) = y(1_G)\} \subseteq B_3$ but $\{a, b, c\} \not\subseteq g B_3 \subseteq A$. Therefore
$$|\{h \in B_3 \: y(g h) = y(1_G)\}| < |\{h \in A \: y(h) = y(1_G)\}|$$
$$= |\{h \in B_3 \: y(h) = y(1_G) \}| = |\{h \in B_3 \: y(gh) = y(1_G)\}|.$$
This is clearly a contradiction.
\end{proof}

\begin{cor}\label{GEN DENSE}
Let $G$ be a countably infinite group, $x \in 2^G$, and $\epsilon > 0$. Then there exists a nontrivial locally recognizable function $R: A \rightarrow 2$ such that for any fundamental $c \in 2^{\subseteq G}$ compatible with $R$ and any $y \in 2^G$ extending $c$ there is $g \in G$ with $d(x, g \cdot y) < \epsilon$.
\end{cor}

\begin{proof}
Let $k \in \N$ be such that $2^{-k} < \epsilon$, and let $B = \{g_0, g_1, \ldots, g_k\}$, where $g_0, g_1, \ldots$ is the fixed enumeration of $G$ used in defining the metric $d$. Set $Q = x|_B$ and apply the previous proposition to get a nontrivial locally recognizable function $R: A \rightarrow 2$ extending $Q$.

Now let $c$ be fundamental with respect to some $(\Delta_n, F_n)_{n \in \N}$ and be compatible with $R$. Let $y \in 2^G$ extend $c$ and set $g = (\gamma \gamma_1)^{-1}$ for any $\gamma \in \Delta_1$. Then for every $a \in A$, $(g \cdot y)(a) = y(g^{-1} a) = R(a)$. So in particular, for every $b \in B$ $(g \cdot y)(b) = y(g^{-1} b) = R(b) = x(b)$. Therefore $d(x, g \cdot y) < \epsilon$.
\end{proof}

\begin{theorem}\label{thm:density}
If $G$ is a countably infinite group, then the collection of 2-colorings on $G$ is dense in $2^G$.
\end{theorem}

\begin{proof}
Arbitrarily fix $x \in 2^G$ and $\epsilon > 0$. First apply Corollary~\ref{GEN DENSE} to obtain a nontrivial locally recognizable function $R$ such that for any fundamental $c$ compatible with $R$ and any $y\in 2^G$ extending $c$, there is $g\in G$ with
$d(x,g\cdot y)<\epsilon$. Next follow the final argument in the proof of Theorem \ref{GEN COL} to obtain a fundamental $c$ 
compatible with $R$ such
that every $y\in 2^G$ extending $c$ is a $2$-coloring. Let $y$ be an arbitrary element of $2^G$ extending $c$. Let $g\in G$ be  as promised above. Then $y$ is a $2$-coloring, and so is $g\cdot y$. We now have $d(x,g\cdot y)<\epsilon$. This completes the proof of the theorem. 
\end{proof}

\begin{theorem} \label{THM PERF DENSE COL}
Let $G$ be a countably infinite group, $x \in 2^G$, and $\epsilon > 0$. Then there is a perfect set of pairwise orthogonal 2-colorings within the $\epsilon$-ball about $x$.
\end{theorem}

\begin{proof}
Arbitrarily fix $x\in 2^G$ and $\epsilon>0$. First appy Corollary \ref{GEN DENSE} to obtain a locally recognizable function $R$, and then apply the final argument in the proof of Theorem \ref{GEN COL} to obtain a fundamental $c$ compatible with $R$. Next apply Proposition \ref{GEN ORTH} to obtain a perfect set $\{ x_{\tau}\,:\, \tau\in 2^\N\}$ of pairwise orthogonal $2$-colorings extending $c$. Now by the proof of Proposition \ref{GEN DENSE}, if we let $g=(\gamma\gamma_1)^{-1}$ for any $\gamma\in\Delta_1$, we have $d(x,g\cdot x_{\tau})<\epsilon$. Note that this $g$ only depends on the blueprint inducing $c$, and in particular does not depend on $\tau$. We thus obtained a set $\{g\cdot x_\tau\,:\, \tau \in 2^\N\}$ of pairwise orthogonal $2$-colorings within the $\epsilon$-ball about $x$. By the continuity of the group action, $\{g\cdot x_\tau\,:\, \tau\in 2^\N\}=g\cdot\{ x_\tau\,:\, \tau\in 2^\N\}$ is still perfect. 
\end{proof}

We use the following notation. If $c \in 2^{\subseteq G}$, we let $\overline{c} \in 2^G$ denote the {\it conjugate} of $c$:
\index{conjugate}
$$\overline{c}(g) = 1 - c(g), \mbox{ for all $g\in \dom(c)$}.$$

\begin{prop}\label{prop:625}
Let $G$ be a countably infinite group, let $B \subseteq G$ be finite but nonempty, and let $Q : B \rightarrow 2$ be an arbitrary function. Then there exists a nontrivial locally recognizable function $R : A \rightarrow 2$ extending $Q$ with the property that if $c \in 2^{\subseteq G}$ is fundamental and compatible with $R$, $x, y \in 2^G$, $x \supseteq c$, $y \supseteq \overline{c}$, then $x$ is orthogonal to $y$. In particular, if $c \in 2^{\subseteq G}$ is canonical and compatible with $R$ and $x \in 2^G$ extends $c$, then $x$ is orthogonal to its conjugate $\overline{x}$.
\end{prop}

\begin{proof}
By applying Proposition \ref{LR EXT} if necessary, we may assume $Q$ is a nontrivial locally recognizable function. Choose a finite $C \subseteq G$ disjoint from $B$ and having cardinality strictly greater than $B$. In particular $|C|\geq 2$. Let $R: B C^{-1} C \cup C \rightarrow 2$ extend $Q$ and have value $1-Q(1_G)$ on $(B C^{-1} C \cup C) - B$. Then $R$ is a nontrivial locally recognizable function.

Let $c \in 2^{\subseteq G}$ be fundamental with respect to some $(\Delta_n, F_n)_{n \in \N}$ and compatible with $R$. Let $x, y \in 2^G$ with $x \supseteq c$ and $y \supseteq \overline{c}$. Let $H$ be finite with $\Delta_1 H = G$. Set $T = H^{-1} \gamma_1 C$ and let $g_1, g_2 \in G$ be arbitrary. There is $h \in H^{-1}$ with $g_1 h \in \Delta_1$. Then $h \gamma_1 C \subseteq T$. Towards a contradiction, suppose $x(g_1 h \gamma_1 \chi) = y(g_2 h \gamma_1 \chi)$ for all $\chi  \in C$. Then $y(g_2 h \gamma_1 \chi ) = 1 - Q(1_G)$ for all $\chi  \in C$. However, it is easy to see that $\overline{c}$ is fundamental with respect to $(\Delta_n, F_n)_{n\in\N}$ and compatible with $\overline{R}$. So
$$g_2 h \gamma_1 C \subseteq \{u \in G \: y(u) = \overline{R}(1_G)\}. $$
By clause (iv) of Theorem \ref{FM} and the definition of fundamental functions,
$$ \{u \in G \: y(u) = \overline{R}(1_G)\} \subseteq \Delta_1 (\gamma_1 F_0 \cup D^1_0)\ \cup\ \bigcup_{n\geq 1} \Delta_n\Lambda_nb_{n-1}.$$
However, it is easy to see that for all $n\geq 1$, $\Delta_n\Lambda_nb_{n-1}\subseteq \Delta_1D^1_0$. Hence we actually have that $g_2h\gamma_1 C\subseteq\Delta_1(\gamma_1F_0\cup D^1_0)$.

If $g_2 h \gamma_1 \chi  = \psi \in \Delta_1 (D^1_0 - \{\gamma_1\})$, then since $C = \chi  (\chi^{-1} C) \subseteq \chi  F_0$ we have
$$(g_2 h \gamma_1 C - \{g_2 h \gamma_1 \chi \}) \subseteq (\psi F_0 - \{\psi\}) \subseteq G - \Delta_1 (\gamma_1 F_0 \cup D^1_0).$$
Since $|C|\geq 2$, this is a contradiction to our previous statement. So it must be that $g_2 h \gamma_1 C \subseteq \Delta_1 \gamma_1 F_0$. Let $\psi \in \Delta_1$ be such that $g_2 h \gamma_1 C \cap \psi \gamma_1 F_0 \neq \varnothing$. For $f \in F_0$ we have $y(\psi \gamma_1 f) = \overline{R}(1_G) = 1 - Q(1_G)$ only if $f \in B$. So $g_2 h \gamma_1 C \cap \psi \gamma_1 B \neq \varnothing$, and therefore $g_2 h \gamma_1 C \subseteq \psi \gamma_1 B C^{-1} C \subseteq \psi \gamma_1 F_0$. It follows $g_2 h \gamma_1 C \subseteq \psi \gamma_1 B$, but then
$$|C| = |g_2 h \gamma_1 C| \leq |\psi \gamma_1 B| = |B| < |C|$$
which is a contradiction. Therefore there is $\chi  \in C$ with $x(g_1 h \gamma_1 \chi ) \neq y(g_2 h \gamma_1 \chi )$. We conclude $x$ and $y$ are orthogonal.
\end{proof}

\begin{cor}
Let $G$ be a countably infinite group, $x \in 2^G$, and $\epsilon > 0$. Then there is a conjugation invariant perfect set of pairwise orthogonal 2-colorings contained in the union of the balls of radius $\epsilon$ about $x$ and $\overline{x}$.
\end{cor}

\begin{proof}
Fix $x\in 2^G$ and $\epsilon>0$. Apply Proposition \ref{GEN DENSE} to obtain a locally recognizable function $Q$, and then use Proposition \ref{prop:625} to get a locally recognizable function $R$ extending $Q$. 
The rest of the proof follows that of Theorem \ref{THM PERF DENSE COL}.
\end{proof}

\section{Characterization of the ACP} \label{SEC ACP}

This section focuses on the uses of blueprints. We begin by constructing a blueprint which will be needed in Sections \ref{SEC RIGID} and \ref{SECT ISO}. Afterwards, we construct another blueprint and use it to characterize those groups which have the ACP.

\begin{prop}\label{SCHEME 1}
Let $G$ be a countably infinite group. Then there is a centered blueprint $(\Delta_n, F_n)_{n \in \N}$ guided by a growth sequence such that for every $g \in G - \Z(G)$ and every $n \geq 1$ there are infinitely many $\gamma \in \Delta_n$ with $\gamma g \neq g \gamma$. Furthermore, if $(p_n)_{n \geq 1}$ and $(q_n)_{n \geq 1}$ are functions and each $p_n$ has subexponential growth, then there is a blueprint having the properties listed in the previous sentence and with
$$|\Lambda_n| \geq q_n(|F_{n-1}|) + \log_2 \ (p_n(|F_n|))$$
for all $n \geq 1$.
\end{prop}

\begin{proof}
Let $(p_n)_{n \geq 1}$ and $(q_n)_{n \geq 1}$ be sequences of functions with each $p_n$ of subexponential growth. We may assume that each $p_n$ and $q_n$ are nondecreasing. Let $R: A \rightarrow 2$ be any nontrivial locally recognizable function. Without loss of generality $1_G\in A$. Let $(A_n)_{n \in \N}$ be a sequence of finite subsets of $G$ with $A_0 = A$ and $\bigcup_{n \in \N} A_n = G$. Set $H_0 = A_0$. In general, once $H_0$ through $H_{n-1}$ have been constructed, define $H_n$ as follows. Let $C_n$ be a finite set such that the $C_n$-translates of $H_{n-1}$ are disjoint, $C_n H_{n-1} \cap H_{n-1} = \varnothing$, and for every $h \in H_{n-1} - \Z(G)$ there is $c \in C_n$ with $ch \neq hc$. Then choose $H_n$ so that
$$H_n \supseteq C_n H_{n-1} \cup A_n \cup H_{n-1} (H_0^{-1} H_0) \cdots (H_{n-1}^{-1} H_{n-1})$$
and $\rho(H_n; H_{n-1}) \geq 3 + q_n(|H_{n-1}|) + \log_2 \ (p_n(|H_n|))$. The sequence $(H_n)_{n \in \N}$ is then a growth sequence.

Recall that in the proof of Theorem \ref{EXIST STRONG BP} each $\delta^n_{n-1}$ was chosen arbitrarily aside from the requirement that $1_G \in \delta^n_{n-1}$ and the $\delta^n_{n-1}$-translates of $F_{n-1}$ be maximally disjoint and contained within $H_n$. We may therefore require that $C_n \subseteq \delta^n_{n-1}$ for every $n \geq 1$. Let $(\Delta_n, F_n)_{n \in \N}$ be the blueprint constructed from Theorem \ref{EXIST STRONG BP} with this change. Then this blueprint is centered and guided by $(H_n)_{n \in \N}$. Notice $\delta^n_{n-1} = D^n_{n-1} \subseteq \Delta_{n-1}$. Suppose $g \in G - \Z(G)$. Let $n \geq 1$ be such that $g \in H_{n-1}$. Then by the definition of $C_m$ for $m \geq n$, we have that there is $\gamma_m \in C_m \subseteq \delta^m_{m-1} \subseteq \Delta_{m-1}$ with $\gamma_m g \neq g \gamma_m$. If $k > m \geq n$, then $\gamma_k \neq \gamma_m$ since $\gamma_m \in F_m \subseteq F_k$ and $\gamma_k \neq 1_G$. Since $(\Delta_n)_{n \in \N}$ is a decreasing sequence (clause (i) of Lemma \ref{STRONG BP LIST}), it follows that for every $m \geq 1$ there are infinitely many elements of $\Delta_m$ which do not commute with $g$ (namely $\gamma_k$ for all $k \geq \max(n, m)$).

If we set $B_n = F_n$, then the last claim is satisfied as well since each $p_n$ and each $q_n$ is nondecreasing and
$$|\Lambda_n| \geq \rho(H_n; H_{n-1}) - 3$$
(see the proof of Corollary \ref{GROW BP}).
\end{proof}

The proposition below constructs a blueprint which is essential in characterizing which groups have the ACP. Recall the notation
$\Z_G(g) = \{h \in G \: h g = g h\}$.

\begin{prop} \label{SCHEME 2}
Let $G$ be a countably infinite group with an element $u \neq 1_G$ satisfying $|\Z_G(u^i)| < \infty$ whenever $i\in{\mathbb Z}$ and $u^i \neq 1_G$, and let $1_G \in A \subseteq G$ be finite with $u \cdot A = A$. Then there is a blueprint $(\Delta_n, F_n)_{n \in \N}$ such that $u \cdot \Delta_n = \Delta_n$ for every $n \in \N$ and with $F_0 = A$. Furthermore, the blueprints with this property can have any prescribed polynomial growth.
\end{prop}

\begin{proof}
Since $\langle u \rangle \subseteq \Z_G(u)$, the order of $u$ must be finite. So finding a finite set $A$ with $u \cdot A = A$ is not an obstacle to applying this proposition. Notice for $i, j \in \mathbb{Z}$, $g \in G$, and $F \subseteq G$ we have
$$u^i g F \cap u^j g F \neq \varnothing \Longleftrightarrow g^{-1} u^{i-j} g \in F F^{-1}.$$
Also for $g, h \in G$ and $i \in \mathbb{Z}$
$$g^{-1} u^i g = h^{-1} u^i h \Longleftrightarrow h g^{-1} \in \Z_G(u^i).$$
Thus, it follows that if $F \subseteq G$ is finite, then for all but finitely many $g \in G$ the $\langle u \rangle$-translates of $g F$ are disjoint. For finite subsets $F \subseteq G$, define $V(F)$ to be the finite (possibly empty) set consisting of all $g \in G$ with the property that the $\langle u \rangle$-translates of $g F$ are not disjoint. Notice that $u \cdot V(F) = V(F)$. By the above remarks, we have that
$$V(F) = \{g \in G \: \exists i \in \ZZ \ u^i \neq 1_G \text{ and } g^{-1} u^i g \in F F^{-1}\}.$$
So if $M = \max\{ |Z_G(u^i)| \: i \in \ZZ \text{ and } u^i \neq 1_G\}$ then
$$|V(F)| \leq |\langle u \rangle | \cdot M \cdot |F F^{-1}|.$$

Notice that if $g, h \in G$ and $h F \cap \langle u \rangle g F = \varnothing$, then immediately we have $\langle u \rangle h F \cap \langle u \rangle g F = \varnothing$. Of particular importance, if $F, H \subseteq G$ are finite, $u \cdot H = H$, and $H \cap V(F) = \varnothing$, then there exists a set $\delta$ such that $u \cdot \delta = \delta$ and the $\delta$-translates of $F$ are contained and maximally disjoint within $H$.

By considering the function $V$, it is easy to modify the proof of Lemma \ref{RATIO} to arrive at the following conclusion. If $A, B \subseteq G$ are finite, $1_G \in A$, and $\epsilon > 0$, then there is a finite $C \subseteq G$ containing $B$ with $u \cdot C = C$ and $\rho(C; A) > \frac{|C|}{|A|} (1 - \epsilon)$. The changes to the proof of Lemma \ref{RATIO} are the following. Replace $B$ with $\langle u \rangle B$ if necessary so that $u \cdot B = B$. By avoiding the finite set $V(A A^{-1})$, one can choose $\Delta$ so that $u \cdot \Delta = \Delta$. The computation in the proof shows that
$$C = B \cup \Lambda A$$
satisfies $\rho(C; A) > \frac{|C|}{|A|} (1 - \epsilon)$ as long as $\Lambda$ is a sufficiently large finite subset of $\Delta$. We can of course choose a sufficiently large $\Lambda \subseteq \Delta$ with $u \cdot \Lambda = \Lambda$, and hence we obtain a $C$ satisfying the inequality and with $u \cdot C = C$.

An immediate consequence to the previous paragraph is the following. If $A, B \subseteq G$ are finite, $1_G \in A$, and $f: \N \rightarrow \N$ is a function of subexponential growth, then there is a finite $C \subseteq G$ containing $B$ with $u \cdot C = C$ and $2^{\rho(C;A)} > f(|C|)$ (see Lemma \ref{FAST GROWTH}).

Let $(p_n)_{n \geq 1}$ be a sequence of functions of polynomial growth. We may suppose that each $p_n$ is nondecreasing. For $n \geq 1$ and $k \in \N$ define
$$q_n(k) = 8 p_n(2 \cdot |\langle u \rangle | \cdot M \cdot k^4).$$
Then $(q_n)_{n \geq 1}$ is a sequence of functions of polynomial growth. Let $(A_n)_{n \in \N}$ be an increasing sequence of finite subsets of $G$ with $G = \bigcup_{n \in \N} A_n$ and $A_0 = A$. Set $H_0 = A_0$. Once $H_0$ through $H_{n-1}$ have been defined, use the previous paragraph to find a finite $H_n \subseteq G$ satisfying $u \cdot H_n = H_n$,
$$H_n \supseteq A_n \cup V(H_{n-1}) H_{n-1} H_{n-1}^{-1} H_{n-1} \cup H_{n-1} (H_0^{-1} H_0) (H_1^{-1} H_1) \cdots (H_{n-1}^{-1} H_{n-1}),$$
and
$$\rho(H_n; H_{n-1}) \geq \log_2 \ q_n(|H_n|) + \rho(V(H_{n-1}) H_{n-1} H_{n-1}^{-1} H_{n-1}; H_{n-1}).$$
The sequence $(H_n)_{n \in \N}$ is easily checked to be a growth sequence. Notice that by clause (iii) of Lemma \ref{LEM RHO}
$$\rho(H_n; H_{n-1}) \geq \log_2 \ q_n(|H_n|) + \rho(V(H_{n-1}) H_{n-1} H_{n-1}^{-1} H_{n-1}; H_{n-1})$$
implies
$$\rho(H_n - V(H_{n-1})H_{n-1}; H_{n-1}) \geq \log_2 \ q_n(|H_n|).$$

Set $F_0 = H_0 = A_0 = A$. Then $1_G \in F_0$ and $u \cdot F_0 = F_0$. So $1_G \in V(F_0)$ and for any finite $F \subseteq G$ containing $F_0$ we have $1_G \in V(F)$. Suppose $F_0$ through $F_{n-1}$ have been defined with the property that for all $0 \leq m < n$, $u \cdot F_m  = F_m$.
Let $\delta \subseteq G$ be such that $u \cdot \delta = \delta$ and the $\delta$-translates of $F_{n-1}$ are contained and maximally disjoint within $H_n - V(H_{n-1}) H_{n-1}$. Notice that
$$\begin{array}{rcl}|\delta| &\geq & \rho(H_n - V(H_{n-1})H_{n-1}; F_{n-1}) \\
&\geq & \rho(H_n - V(H_{n-1})H_{n-1}; H_{n-1}) \geq \log_2 \ q_n(|H_n|) \geq 3.
\end{array}$$
Set $\delta^n_{n-1} = \delta \cup \{1_G\}$. Then $\delta^n_{n-1}F_{n-1}=\delta F_{n-1}\cup F_{n-1}$, and so 
$$u\cdot (\delta^n_{n-1}F_{n-1})=\delta^n_{n-1}F_{n-1}.$$
Now suppose that $\delta^n_{n-1}$ through $\delta^n_{k+1}$ have been defined for some $0 \leq k < n-1$. Inductively assume
that for all $n+1\geq j\geq k+1$, $u\cdot (\delta^n_jF_j)=\delta^n_jF_j$. Define
$$B^n_k = \{g \in G \: \{g\} (F_{k+1}^{-1} F_{k+1}) (F_{k+2}^{-1} F_{k+2}) \cdots (F_{n-1}^{-1} F_{n-1}) \subseteq H_n\}$$
and notice that $u \cdot B^n_k = B^n_k$ and $H_{n-1} \subseteq B^n_k$. Let $\delta^n_k \subseteq G$ be such that $u \cdot \delta^n_k = \delta^n_k$ and the $\delta^n_k$-translates of $F_k$ are contained and maximally disjoint within
$$ B^n_k - V(H_k) H_k - \bigcup_{k < m < n} \delta_m^n F_m.$$
Denoting the set displayed above by $S$, then $S\cap V(F_k)=\emptyset$ since $1_G\in F_k\subseteq H_k$ and $V(F_k)\subseteq V(H_k)$. This, together with $u\cdot S=S$, guarantees that $\delta^n_k$ exists. Finally, define
$$F_n = \bigcup_{0 \leq k < n} \delta^n_k F_k.$$
Then $u\cdot F_n=F_n$.

We now apply Lemma \ref{MAKE PREBP} to get a pre-blueprint $(\tilde{\Delta}_n, F_n)_{n \in \N}$ (conditions (i) through (v) 
of Lemma \ref{MAKE PREBP} are clearly satisfied). To be specific, for $k \in \N$ set $\tilde{D}^k_k = \{1_G\}$ and for $n > k$ define
$$\tilde{D}^n_k = \bigcup_{k \leq m < n} \delta^n_m \tilde{D}^m_k.$$
Then define $\tilde{\Delta}_k = \bigcup_{n \geq k} \tilde{D}^n_k$.

We claim that $u \cdot (\tilde{D}^n_k - \{1_G\}) = \tilde{D}^n_k - \{1_G\}$ for all $n, k \in \N$ with $n \geq k$. Fix $k \in \N$. The claim is obvious when $n = k$. Now let $n > k$ and suppose the claim is true for $n-1$. Recall from our construction that $u \cdot \delta^n_m = \delta^n_m$ for all $0 \leq m < n-1$ and that $u \cdot (\delta^n_{n-1} - \{1_G\}) = \delta^n_{n-1} - \{1_G\}$. We have
$$\tilde{D}^n_k - \{1_G\} = \left( \bigcup_{k \leq m < n-1} \delta^n_m \tilde{D}^m_k \right) \cup (\delta^n_{n-1} - \{1_G\}) \tilde{D}^{n-1}_k \cup (\tilde{D}^{n-1}_k - \{1_G\})$$
so the claim follows by induction.

We claim that $(\tilde{\Delta}_n, F_n)_{n \in \N}$ is a blueprint. For $k \in \N$ define $T_k = V(H_k) H_k H_k^{-1}$. It will suffice to show that $\tilde{\Delta}_k T_k = G$. We proceed to verify the following three facts.
\begin{enumerate}
\item[\rm (1)] If $n > k$, $g \in G$, and $g F_k \cap F_n \neq \varnothing$, then either $g \in T_k$ or else $g F_k \cap \delta^n_m F_m \neq \varnothing$ for some $k \leq m < n$;
\item[\rm (2)] $g F_k \cap F_n \neq \varnothing \Longrightarrow g \in \tilde{D}^n_k T_k$ for all $g \in G$ and $k \leq n$;
\item[\rm (3)] $g F_k \subseteq B^n_k \Longrightarrow g \in \tilde{D}^n_k T_k$ for all $g \in G$ and $k < n$.
\end{enumerate}

(Proof of 1). Let $n > k$ and $g \in G$ satisfy $g F_k \cap F_n \neq \varnothing$. It suffices to show that if $g \not\in T_k$ and $g F_k \cap \delta_m^n F_m = \varnothing$ for all $k < m < n$ then $g F_k \cap \delta_k^n F_k \neq \varnothing$ (since this will validate the claim with $m = k$). As $F_n = \bigcup_{0 \leq t < n} \delta_t^n F_t$, there is $0 \leq t \leq k$ with $g F_k \cap \delta_t^n F_t \neq \varnothing$. If $t = k$, then we are done. So suppose $t < k$. We have
$$g F_k \subseteq \delta_t^n F_t  F_k^{-1} F_k \subseteq \delta_t^n F_t (F_{t+1}^{-1} F_{t+1}) (F_{t+2}^{-1} F_{t+2}) \cdots (F_k^{-1} F_k).$$
Hence
$$g F_k (F_{k+1}^{-1} F_{k+1}) \cdots (F_{n-1}^{-1} F_{n-1}) \subseteq \delta_t^n F_t (F_{t+1}^{-1} F_{t+1}) \cdots (F_{n-1}^{-1} F_{n-1}).$$
However, by definition $\delta_t^n F_t \subseteq B^n_t$. So the right hand side of the expression above is contained within $H_n$, and therefore $g F_k \subseteq B^n_k$. Also $g \not\in T_k$ implies $g F_k$ is disjoint from $V(H_k) H_k$. Thus
$$g F_k \subseteq B^n_k - V(H_k) H_k - \bigcup_{k < m < n} \delta_m^n F_m.$$
It now follows from the definition of $\delta_k^n$ that $g F_k \cap \delta_k^n F_k \neq \varnothing$. This substantiates our claim.

(Proof of 2). Fix $k \in \N$. We prove the claim by an induction on $n\geq k$. If $n = k$ then the claim is clear. Now assume the claim is true for all $k \leq m < n$. Let $g \in G$ satisfy $g F_k \cap F_n \neq \varnothing$. If $g \in T_k$ then we are done since $1_G \in \tilde{D}^n_k$. So we may assume $g \not\in T_k$. By (1) we have that $g F_k \cap \delta^n_m F_m \neq \varnothing$ for some $k \leq m < n$. Let $\gamma \in \delta^n_m$ be such that $g F_k \cap \gamma F_m \neq \varnothing$. Then $\gamma^{-1} g F_k \cap F_m \neq \varnothing$, so by the induction hypothesis $\gamma^{-1} g \in \tilde{D}^m_k T_k$. By the definition of $\tilde{D}^n_k$ we have $\gamma \tilde{D}^m_k \subseteq \delta^n_m \tilde{D}^m_k \subseteq \tilde{D}^n_k$. Thus, $g \in \tilde{D}^n_k T_k$. 

(Proof of 3). Fix $k < n$ and let $g \in G$ be such that $g F_k \subseteq B^n_k$. We must show $g \in \tilde{D}^n_k T_k$. We are done if $g F_k \cap \delta^n_k F_k \neq \varnothing$ since $F_k F_k^{-1} \subseteq T_k$ and $\delta^n_k = \delta^n_k \tilde{D}^k_k \subseteq \tilde{D}^n_k$. So suppose $g F_k \cap \delta^n_k F_k = \varnothing$. Recall that in the construction of $F_n$ we defined $\delta_{n-1}^n$ through $\delta_{k+1}^n$ first and then chose $\delta_k^n$ so that its translates of $F_k$ would be maximally disjoint within $B^n_k - V(H_k) H_k - \bigcup_{k<m<n} \delta_m^n F_m$. Thus we cannot have $g F_k \subseteq B^n_k - V(H_k) H_k - \bigcup_{k<m<n} \delta_m^n F_m$ as this would violate the definition of $\delta_k^n$. Since $g F_k \subseteq B^n_k$, we must have either $g F_k \cap V(H_k) H_k \neq \varnothing$ or $g F_k \cap (\bigcup_{k<m<n} \delta_m^n F_m) \neq \varnothing$. If $g F_k \cap V(H_k) H_k \neq \varnothing$ then $g \in V(H_k)H_kF_k^{-1}\subseteq T_k \subseteq \tilde{D}^n_k T_k$. So we may suppose $g F_k \cap (\bigcup_{k<m<n} \delta_m^n F_m) \neq \varnothing$. Let $k < m < n$ and $\gamma \in \delta^n_m$ be such that $g F_k \cap \gamma F_m \neq \varnothing$. Then $\gamma^{-1} g F_k \cap F_m \neq \varnothing$ and thus $\gamma^{-1} g \in \tilde{D}^m_k T_k$ by (2). Now we have $\gamma \tilde{D}^m_k \subseteq \delta^n_m \tilde{D}^m_k \subseteq \tilde{D}^n_k$ so that $g \in \tilde{D}^n_k T_k$. This completes the proof of (3).

Recall that we have $H_{n-1} \subseteq B^n_k$. To see $\tilde{\Delta}_k T_k = G$,
fix $g\in G$. Then for sufficiently large $n>k$ we have $gF_k\subseteq H_{n-1}$ since $\{H_n\}_{n\in\N}$ is increasing.
Thus $gF_k\subseteq B^n_k$, and by (3), $g\in \tilde{D}^n_kT_k$. We thus conclude $(\tilde{\Delta}_n, F_n)_{n \in \N}$ is a blueprint.

Pick a nonidentity $\gamma_n \in \tilde{\Delta}_n$ and define $T_n' = T_n \cup \gamma_n^{-1} T_n$. For each $n \in \N$ define $\Delta_n = \tilde{\Delta}_n - \{1_G\}$. Then $(\Delta_n, F_n)_{n \in \N}$ is a blueprint satisfying $u \cdot \Delta_n = \Delta_n$ for each $n \in \N$. To see this is a blueprint, just notice that $\Delta_n T_n' = G$ for each $n \in \N$. We also have that
$$|T_n'| \leq 2|T_n| = 2|V(H_n) H_n H_n^{-1}| \leq 2 \cdot |V(H_n)| \cdot |H_n|^2
 \leq 2 \cdot |\langle u \rangle | \cdot M \cdot |H_n|^4.$$
Therefore
$$|\Lambda_n| = |D^n_{n-1}| - 3 = |\delta^n_{n-1} - \{1_G\}| - 3 \geq -3 + \log_2 \ q_n(|H_n|)$$
$$= -3 + \log_2 \ (8p_n(2 \cdot |\langle u \rangle | \cdot M \cdot |H_n|^4))$$
$$= \log_2 \ p_n(2 \cdot |\langle u \rangle | \cdot M \cdot |H_n|^4) \geq \log_2 \ p_n(|T_n'|).$$
So the blueprint satisfies the growth condition.
\end{proof}

We are now ready to characterize which groups have the ACP (almost $2$-coloring property).

\begin{theorem}\label{thm:ACP}
Let $G$ be a countably infinite group. Then $G$ has the ACP if and only if for every $1_G \neq u \in G$ there is $1_G \neq v \in \langle u \rangle$ with $|\Z_G(v)| = \infty$.
\end{theorem}

\begin{proof}
If $G$ has the stated property, then it was proved in Proposition \ref{PROP HALF ACP} that $G$ has the ACP. So it will suffice to assume that there is $1_G \neq u \in G$ with $|\Z_G(u^i)| < \infty$ whenever $u^i \neq 1_G$ and then show that $2^G$ contains an almost $2$-coloring which is not a $2$-coloring. Fix $u \in G$ with this property.

For $n \geq 1$ and $k \in \N$ define $p_n(k) = 8 |\langle u \rangle | k^4+1$. Then $(p_n)_{n \geq 1}$ is a sequence of functions of polynomial growth. Pick any nontrivial locally recognizable function $R: A \rightarrow 2$. By the previous proposition, there is a blueprint $(\Delta_n, F_n)_{n \in \N}$ and finite sets $(B_n)_{n \in \N}$ with $F_0 = A$ and such that for every $n \in \N$ $u \cdot \Delta_n = \Delta_n$, $\Delta_n B_n B_n^{-1} = G$, and
$$|\Lambda_{n+1}| \geq \log_2 \ p_{n+1}(|B_{n+1}|).$$
Apply Theorem \ref{FM} to get a function $c$ canonical with respect to $(\Delta_n, F_n)_{n \in \N}$ and compatible with $R$. Recall that the value of $c$ on $\bigcap_{n \in \N} \Delta_n a_n$ and $\bigcap_{n \in \N} \Delta_n b_n$ is arbitrary. We may therefore assume that $c$ is constant on $\bigcap_{n \in \N} \Delta_n a_n$ and $\bigcap_{n \in \N} \Delta_n b_n$. Notice that
$$u \cdot \bigcap_{n \in \N} \Delta_n a_n = \bigcap_{n \in \N} \Delta_n a_n$$
and similarly with the $a_n$'s replaced with $b_n$'s.

We claim that $u \cdot c = c$. Notice that
$$G - \dom(c) = \bigcup_{n \geq 1} \Delta_n \Lambda_n b_{n-1}$$
is invariant under left multiplication by $u$. Now fix $g \in \dom(c)$. We will show that $c(g) = c(ug)$. We proceed by cases.

\underline{Case 1:} $g \in \bigcap_{n \in \N} \Delta_n \{a_n, b_n\}$. Then by our earlier comment we immediately have $c(g) = c(ug)$.

\underline{Case 2:} There is $n \geq 1$ with $g \in \Delta_n \{a_n, b_n\} - \Delta_{n+1} \{a_{n+1}, b_{n+1}\}$. \underline{Subcase A:} $g \in \Delta_{n+1} F_{n+1}$. Let $\gamma \in \Delta_{n+1}$ and $f \in F_{n+1}$ be such that $g = \gamma f$. Then $f \neq a_{n+1}, b_{n+1}$. Since $u \cdot \Delta_{n+1} = \Delta_{n+1}$, we have that $u \gamma \in \Delta_{n+1}$. So $c(g) = c(\gamma f) = c(u \gamma f) = c(ug)$ by conclusion (vi) of Theorem \ref{FM}. \underline{Subcase B:} $g \not\in \Delta_{n+1} F_{n+1}$. Then we also have $ug \in \Delta_n \{a_n, b_n\} - \Delta_{n+1} F_{n+1}$. From the proof of Theorem \ref{FM} it can be seen that
$$c(\Delta_n \{a_n, b_n\} - \Delta_{n+1} F_{n+1}) = \{0\}.$$
So $c(g) = c(ug)$.

\underline{Case 3:} $g \not\in \Delta_1 \{a_1, b_1\}$. \underline{Subcase A:} $g \in \Delta_1 F_1$. Let $\gamma \in \Delta_1$ and $f \in F_1$ be such that $g = \gamma f$. Then $f \neq a_1, b_1$ and $u \gamma \in \Delta_1$. So $c(g) = c(\gamma f) = c(u \gamma f) = c(ug)$ by conclusion (vi) of Theorem \ref{FM}. \underline{Subcase B:} $g \not\in \Delta_1 F_1$. Then $ug \not\in \Delta_1 F_1$. By conclusion (iv) of Theorem \ref{FM} we have $c(g) = 1-R(1_G) = c(ug)$.

We conclude that $u \cdot c = c$ as claimed. Fix an enumeration $s_1, s_2, \ldots$ of the nonidentity group elements of $G$. Let $V_n$ be the test region for the $\Delta_n$-membership test admitted by $c$. Set $T_n = B_n B_n^{-1} (V_n \cup F_n)$. Pick any $h_n \in G - \{1_G, s_n^{-1}\} \{1_G, s_n\} T_n T_n^{-1}$. For each $n \geq 1$ let $\Gamma_n$ be the graph with vertex set $\{ \langle u \rangle \gamma \: \gamma \in \Delta_n\}$ and edge relation given by
$$(\langle u \rangle \gamma, \langle u \rangle \psi) \in E(\Gamma) \Longleftrightarrow \exists i \in \ZZ$$
$$\gamma^{-1} u^i \psi \in B_n B_n^{-1} \{1_G, h_n^{-1}\} \{s_n, s_n^{-1}\} \{1_G, h_n\} B_n B_n^{-1}$$
for disinct $\langle u \rangle \gamma$ and $\langle u \rangle \psi$. Notice that this edge relation is well defined. We have that
$$\deg_{\Gamma_n} (\langle u \rangle \gamma) \leq 8 |\langle u \rangle | |B_n|^4.$$
Therefore we can find a graph-theoretic $(8 |\langle u \rangle | |B_n|^4 + 1)$-coloring of $\Gamma_n$, say $\mu_n : V(\Gamma_n) \rightarrow \{0, 1, \ldots, 8 |\langle u \rangle | |B_n|^4\}$.

For each $i \geq 1$, define $\B_i : \N \rightarrow \{0,1\}$ so that $\B_i(k)$ is the $i^\text{th}$ digit from least to most significant in the binary representation of $k$ when $k \geq 2^{i-1}$ and $\B_i(k) = 0$ when $k < 2^{i-1}$.

For $n \geq 1$ set $s(n) = |\Lambda_n|$ and let $\lambda^n_1, \lambda^n_2, \ldots, \lambda^n_{s(n)}$ be an enumeration of $\Lambda_n$. Define $y \supseteq c$ by setting
$$y(\gamma \lambda^n_i b_{n-1}) = \B_i (\mu_n(\langle u \rangle \gamma))$$
for each $n \geq 1$, $\gamma \in \Delta_n$, and $1 \leq i \leq s(n)$. Since $2^{s(n)} \geq 8 |\langle u \rangle | |B_n|^4+1$, all integers $0$ through $8 |\langle u \rangle | |B_n|^4$ can be written in binary using $s(n)$ digits. Thus no information is lost between the $\mu_n$'s and $y$. Since
$$G - \dom(c) = \bigcup_{n \geq 1} \Delta_n \Lambda_n b_{n-1},$$
we have that $y \in 2^G$. Also, it is easily checked from the definition of $y$ and the fact that $u \cdot c = c$ that $u \cdot y = y$. Thus $y$ is periodic. We will show that $y$ is an almost $2$-coloring.

Define $x \in 2^G$ by setting $x(1_G) = 1-y(1_G)$ and $x(g) = y(g)$ for $1_G \neq g \in G$. Clearly $x=^*y$. We claim that $x$ is a $2$-coloring of $G$. Fix $1_G \neq s \in G$. Then for some $n \geq 1$ we have $s = s_n$. Define
$$W = \{g \in G \: \exists i \in \ZZ \ g u^i g^{-1} = s_n\}.$$
Notice that $W=V(\{s_n\})^{-1}$ is finite. Set $T = W \cup T_n \cup h_n T_n$ and let $g \in G$ be arbitrary.

If $g^{-1} \in W$ and $g^{-1} u^i g = s_n$ then we have
$$x(g g^{-1}) = x(1_G) \neq y(1_G) = y(u^i) = x(u^i) = x(g s_n g^{-1}).$$
In this case we are done since $g^{-1} \in W \subseteq T$. So we may suppose that $g^{-1} \not\in W$. It follows that $\langle u \rangle g \neq \langle u \rangle g s_n$ and furthermore
$$\langle u \rangle g t \neq \langle u \rangle g s_n t$$
for all $t \in T$.

Notice that by our choice of $h_n$,
$$1_G \not\in T_n^{-1} \{1_G, s_n^{-1}\} \{1_G, s_n\} h_n T_n.$$
So if $g T_n$ or $g s_n T_n$ contains $1_G$ then $g$ is an element of $T_n^{-1} \{1_G, s_n^{-1}\}$ and therefore $1_G$ is neither an element of $g h_n T_n$ nor $g s_n h_n T_n$. If $1_G \not\in g T_n \cup g s_n T_n$ then set $k = 1_G$ and otherwise set $k = h_n$. In any case we have that $1_G \not\in g k T_n \cup g s_n k T_n$. In particular, for all $t \in T_n$
$$x(g k t) = y(g k t) \text{ and } x(g s_n k t) = y(g s_n k t).$$

Since $\Delta_n B_n B_n^{-1} = G$, there is $b \in B_n B_n^{-1}$ with $g k b \in \Delta_n$. We proceed by cases. \underline{Case 1:} $g s_n k b \not\in \Delta_n$. Since $c$ admits a $\Delta_n$-memebership test with test region $V_n$, there must be $v \in V_n$ with
$$x(g k b v) = y(g k b v) \neq y(g s_n k b v) = x(g s_n k b v).$$
The equalities hold because $b\in B_nB_n^{-1}$ and $v\in V_n$, and therefore $bv\in T_n$ and $kbv\in T$.
This case is completed as we have shown $x(gt)\neq x(gs_nt)$ for $t=kbv\in T$. \underline{Case 2:} $g s_n k b \in \Delta_n$. Then we have
$$(g k b)^{-1} (g s_n k b) = b^{-1} k^{-1} s_n k b \in B_n B_n^{-1} \{1_G, h_n^{-1}\} s_n \{1_G, h_n\} B_n B_n^{-1}.$$
It follows that $(\langle u \rangle g k b, \langle u \rangle g s_n k b) \in E(\Gamma_n)$ since $\langle u \rangle g k b \neq \langle u \rangle g s_n k b$. Thus from the definition of $y$ we have that there is $1 \leq i \leq s(n)$ with
$$x(g k b \lambda^n_i b_{n-1}) = y(g k b \lambda^n_i b_{n-1}) \neq y(g s_n k b \lambda^n_i b_{n-1}) = x(g s_n k b \lambda^n_i b_{n-1}).$$
The equalities hold because $b\in B_nB_n^{-1}$ and $\lambda^n_kb_{n-1}\in F_n$, and therefore $b\lambda^n_ib_{n-1}\in T_n$. Also we have $k b \lambda^n_i b_{n-1} \in T$. This completes the proof.
\end{proof}

We give an example of a group without the ACP.  Consider the group $G = \ZZ_2 * \ZZ_2$. The group has the presentation 
$$G = \langle a, b \ | \ a^2 = b^2 = 1_G \rangle.$$ 
One can check that $\Z_G(a) = \langle a \rangle$. So by the proof of Theorem~\ref{thm:ACP} there is an almost $2$-coloring $y \in 2^G$ with $a \cdot y = y$. This gives the promised counterexample to the converse of Lemma~\ref{lem:almostcoloringlemma} (b). Let $x=^*y$ be a $2$-coloring on $G$. Then $x$ is not a strong $2$-coloring by Lemma~\ref{lem:almostcoloringlemma} (e). Thus $x$ is a counterexample to the converse of Lemma~\ref{lem:almostcoloringlemma} (a). 

Notice that the group $G$ considered above is polycyclic and virtually abelian. Since all polycyclic groups are solvable, we have that in general solvable, polycyclic, and virtually abelian groups do not necessarily 
have the ACP.

The above theorem has a very nice general corollary.

\begin{cor}
For a countable group $G$, the following are equivalent:
\begin{enumerate}
\item[\rm (i)] for every compact Hausdorff space $X$ on which $G$ acts continuously and every $x \in X$, if every limit point of $[x]$ is aperiodic then $x$ is hyper aperiodic;
\item[\rm (ii)] for every nonidentity $u \in G$ there is a nonidentity $v \in \langle u \rangle$ having infinite centralizer in $G$.
\end{enumerate}
\end{cor}

\begin{proof}
Very minor modifications to the proof of Proposition \ref{PROP HALF ACP} give the implication (ii) $\Rightarrow$ (i). On the other hand, suppose $G$ does not satisfy (ii). Then by the previous theorem there is a periodic almost $2$-coloring $x$ on $G$. Set $X = \overline{[x]}$. Then $X$ is a compact Hausdorff space on which $G$ acts continuously. By clause (c) of Lemma \ref{lem:almostcoloringlemma} and by Lemma \ref{NEAR EQUIV}, every limit point of $[x]$ is aperiodic. However, $x$ is periodic and is therefore not hyper aperiodic (not a $2$-coloring).
\end{proof}

\chapter{Further Study of Fundamental Functions} \label{CHAP STUDY}

In this chapter, we will focus on developing general tools which aid in implementing the fundamental method. These tools are developed primarily because we need them in later chapters, however we will develop these tools in more generality than they will be used. The tools developed in this chapter will help with three tasks: understanding the relationship between a fundamental function and the points in the closure of its orbit; understanding how minimality relates to fundamental functions and building minimal fundamental functions; and controlling when two fundamental functions generate topologically conjugate subflows. The first section focuses on the closure of the orbit of fundamental functions. The next three sections deal with minimality, and the final section focuses on topological conjugacy among the subflows generated by fundamental functions.

\section{Subflows generated by fundamental functions} \label{SECT SUBFLOW FUND FNCTN}

In this section we will go through some basic observations regarding the closure of the orbit of a fundamental function. We will see that if $x \in 2^G$ is fundamental with respect to a blueprint $(\Delta_n, F_n)_{n \in \N}$ and $y \in \overline{[x]}$, then there are sets $\Delta^y_n$ such that $y$ is fundamental with respect to the blueprint $(\Delta^y_n, F_n)_{n \in \N}$. Finding the sets $\Delta^y_n$ is made easy by the $\Delta_n$ membership test admitted by $x$ for $n \geq 1$.

Fix a blueprint $(\Delta_n, F_n)_{n \in \N}$ and let $x \in 2^G$ be fundamental with respect to this blueprint. For each $n \geq 1$, let $V_n \subseteq \gamma_n F_{n-1} \subseteq F_n$ be the test region for the simple $\Delta_n$ membership test admitted by $x$ (clause (ii) of Theorem \ref{FM}), and let $P_n \in 2^{V_n}$ be the corresponding test function. So for $n \geq 1$ and $g \in G$, we have $g \in \Delta_n$ if and only if $x(gv) = P_n(v)$ for all $v \in V_n$. Now if $y \in \overline{[x]}$, we define $\Delta_n^y$ for $n \geq 1$ by \index{$\Delta_n^y$}
$$\Delta_n^y = \{g \in G \: \forall v \in V_n \ y(gv) = P_n(v)\}.$$
Notice that $\Delta_n^x = \Delta_n$ and $\Delta_n^{g \cdot y} = g \cdot \Delta_n^y$ for all $n \geq 1$ and all $g \in G$.

\begin{prop} \label{SUBFLOW BP}
Let $G$ be a countably infinite group, let $(\Delta_n, F_n)_{n \in \N}$ be a blueprint, and let $x \in 2^G$ be fundamental with respect to this blueprint. Then for every $y \in \overline{[x]}$ we have:
\begin{enumerate}
\item[\rm (i)] if $y = \lim h_m \cdot x$ then $\Delta_n^y = \bigcup_{k \in \N} \bigcap_{m \geq k} h_m \Delta_n$ for all $n \geq 1$;
\item[\rm (ii)] if $n \geq 1$, $B \subseteq G$ is finite, and $\Delta_n B = G$ then $\Delta_n^y B = G$;
\item[\rm (iii)] $\gamma^{-1} (\Delta_k^y \cap \gamma F_n) = D^n_k$, for all $n \geq k \geq 1$ and $\gamma \in \Delta_n^y$;
\item[\rm (iv)] $(\Delta_n^y, F_n)_{n \in \N}$ is a blueprint, where $\Delta^y_0$ is defined by the formula in (i).
\end{enumerate}
\end{prop}

\begin{proof}
Let $1_G = g_0, g_1, \ldots$ be the fixed enumeration of $G$ used in defining the metric $d$ on $2^G$. For each $n \geq 1$, let $V_n \subseteq \gamma_n F_{n-1} \subseteq F_n$ be the test region for the simple $\Delta_n$ membership test admitted by $x$ (clause (ii) of Theorem \ref{FM}), and let $P_n \in 2^{V_n}$ be the corresponding test function. Let $(h_n)_{n \in \N}$ be such that $y = \lim_{n \rightarrow \infty} h_n \cdot x$.

(i). Let $n \geq 1$, and let $\gamma \in \Delta_n^y$. Let $r \in \N$ be such that $\gamma V_n \subseteq \{g_0, g_1, \ldots, g_r\}$. Let $k \in \N$ be such that $d(h_m \cdot x, y) < 2^{-r}$ for all $m \geq k$. Then for all $m \geq k$ and all $v \in V_n$
$$x(h_m^{-1} \gamma v) = (h_m \cdot x)(\gamma v) = y(\gamma v) = P_n(v).$$
Therefore $h_m^{-1} \gamma \in \Delta_n^x$ and $\gamma \in h_m \Delta_n^x$ for all $m \geq k$.

Now let $n \geq 1$ and suppose $\gamma \in \bigcup_{k \in \N} \bigcap_{m \geq k} h_m \Delta_n^x$. Let $r \in \N$ be such that $\gamma V_n \subseteq \{g_0, g_1, \ldots, g_r\}$, let $k \in \N$ be such that $\gamma \in \bigcap_{m \geq k} h_m \Delta_n^x$, and let $m \geq k$ be such that $d(h_m \cdot x, y) < 2^{-r}$. Then for all $v \in V_n$
$$y(\gamma v) = (h_m \cdot x)(\gamma v) = x(h_m^{-1} \gamma v) = P_n(v).$$
The last equality follows from the fact that $h_m^{-1} \gamma \in \Delta_n^x$ since $\gamma \in h_m \Delta_n^x$. It follows that $\gamma \in \Delta_n^y$.

(ii). Fix $g \in G$. Let $r \in \N$ be such that $g B^{-1} V_n \subseteq \{g_0, g_1, \ldots, g_r\}$. Let $m \in \N$ be such that $d(h_m \cdot x, y) < 2^{-r}$. Clearly $h_m \Delta_n^x B = G$, so there is $b \in B^{-1}$ with $g b \in h_m \Delta_n^x$, or equivalently $h_m^{-1} g b \in \Delta_n^x$. Then for all $v \in V_n$
$$y(g b v) = (h_m \cdot x)(g b v) = x(h_m^{-1} g b v) = P_n(v).$$
Hence $g b \in \Delta_n^y$ and $g \in \Delta_n^y B$. We conclude that $\Delta_n^y B = G$.

(iii). Fix $n \geq k \in \N$ and $\gamma \in \Delta_n^y$. If $\lambda \in D^n_k$, then $\Delta_n^x \lambda \subseteq \Delta_k^x$ by clause (i) of Lemma \ref{BP LIST}. Right multiplication by $\lambda$ is a bijection of $G$, so
$$\Delta_n^y \lambda = \bigcup_{i \in \N} \bigcap_{m \geq i} h_m \Delta_n^x \lambda \subseteq \bigcup_{i \in \N} \bigcap_{m \geq i} h_m \Delta_k^x = \Delta_k^y.$$
Thus $\gamma D^n_k \subseteq \Delta_k^y \cap \gamma F_n$. On the other hand, since $F_n$ is finite we can find $m \in \N$ by (i) with $\gamma \in h_m \Delta_n^x$ and
$$\Delta_k^y \cap \gamma F_n \subseteq h_m \Delta_k^x \cap \gamma F_n = h_m [\Delta_k^x \cap h_m^{-1} \gamma F_n] = \gamma ( \Delta_k^x \cap F_n) = \gamma D^n_k.$$

(iv). Define $\Delta_0^y = \bigcup_{k \in \N} \bigcap_{m \geq k} h_m \Delta_0$. One can easily check that (ii) and (iii) remain true when $ n > k = 0$ and when $n = k = 0$. We verify the conditions listed in Definition \ref{def:blueprint}. (Disjoint). Let $n \in \N$ and $\gamma \neq \psi \in \Delta_n^y$. Then by (i) there is $m \in \N$ with $\gamma, \psi \in h_m \cdot \Delta_n^x$. Then $h_m^{-1} \gamma \neq h_m^{-1} \psi \in \Delta_n^x$ so $h_m^{-1} \gamma F_n \cap h_m^{-1} \psi F_n = \varnothing$ and hence $\gamma F_n \cap \psi F_n = \varnothing$. (Dense). This follows immediately from (ii). (Coherent). Suppose $n \geq k \in \N$, $\gamma \in \Delta_n^y$, $\psi \in \Delta_k^y$, and $\psi F_k \cap \gamma F_n \neq \varnothing$. By (i) there is $m \in \N$ with $\gamma \in h_m \Delta_n^x$ and $\psi \in h_m \Delta_k^x$. So $h_m^{-1} \gamma \in \Delta_n^x$, $h_m^{-1} \psi \in \Delta_k^x$, and $h_m^{-1} \psi F_k \cap h_m^{-1} \gamma F_n \neq \varnothing$. It follows that $h_m^{-1} \psi F_k \subseteq h_m^{-1} \gamma F_n$ and hence $\psi F_k \subseteq \gamma F_n$. (Uniform). This follows immediately from (iii). (Growth). The growth condition on blueprints is equivalent to the property $|D^n_{n-1}| \geq 3$. Therefore this follows immediately from (iii). We conclude that $(\Delta_n^y, F_n)_{n \in \N}$ is a blueprint.
\end{proof}

Technically, the definition of $\Delta_0^y$ in clause (iv) above depends on the sequence $(h_m)$ chosen. However, the set $\Delta_0^y$ is essentially unimportant and the non uniqueness of $\Delta_0^y$ is not a problem for us. It is only important to fix a single choice of $\Delta_0^y$ which satisfies the equation above with respect to some sequence $(h_m)$ with $y = \lim h_m \cdot x$. Notice though that for $n \geq 1$ $\Delta_n^y = \bigcup_{k \in \N} \bigcap_{m \geq k} h_m \Delta_n$ for every sequence $(h_m)$ with $y = \lim h_m \cdot x$. Thus, if we have a particular $\Delta_0^y$ in mind, we can always choose to work with the sequence $(h_m)$ with $y = \lim h_m \cdot x$ and $\Delta_0^y = \bigcup_{k \in \N} \bigcap_{m \geq k} h_m \Delta_0$.

Recall that in general for a blueprint $\alpha_n$, $\beta_n$, and $\gamma_n$ are assumed only to be distinct members of $D^n_{n-1}$, and these group elements are used to define $\Lambda_n$, $a_n$, and $b_n$. Therefore the objects $\alpha_n$, $\beta_n$, $\gamma_n$, $a_n$, $b_n$, and $\Lambda_n$ can all be used with the blueprint $(\Delta_n^y, F_n)_{n \in \N}$ and all the conclusions of Lemmas \ref{BP LIST} and \ref{STRONG BP LIST} will hold. This is a very important observation.

\begin{lem} \label{SUBFLOW BP LIST}
Let $G$ be a countably infinite group, let $(\Delta_n, F_n)_{n \in \N}$ be a blueprint, and let $x \in 2^G$ be fundamental with respect to this blueprint. Then for every $y \in \overline{[x]}$ we have:
\begin{enumerate}
\item[\rm (i)] if $(\Delta_n)_{n \in \N}$ is decreasing then so is $(\Delta_n^y)_{n \in \N}$;
\item[\rm (ii)] if $(\Delta_n, F_n)_{n \in \N}$ is maximally disjoint then so is $(\Delta_n^y, F_n)_{n \in \N}$;
\item[\rm (iii)] if $(\Delta_n, F_n)_{n \in \N}$ is guided by a growth sequence $(H_n)_{n \in \N}$, then $(\Delta_n^y, F_n)_{n \in \N}$ is guided by the same growth sequence;
\item[\rm (iv)] if $(\Delta_n, F_n)_{n \in \N}$ is centered, directed, and the $\Delta_k$-translates of $F_k$ are maximally disjoint for some $k \in \N$ then
$$\left| \bigcap_{n \in \N} \Delta_n^y \right| \leq 1;$$
\item[\rm (v)] if $(\Delta_n, F_n)_{n \in \N}$ is centered and directed and $g \in \bigcap_{n \in \N} \Delta_n^y$ for some $g \in G$, then $g \Delta_n \subseteq \Delta_n^y$ for all $n \in \N$;
\item[\rm (vi)] if $(\Delta_n, F_n)_{n \in \N}$ is centered, directed, and maximally disjoint and $g \in G$ satisfies $g \in \bigcap_{n \in \N} \Delta_n^y$, then $g \Delta_n = \Delta_n^y$ for all $n \in \N$.
\end{enumerate}
\end{lem}

\begin{proof}
(i). Since $(\Delta_n^x)_{n \in \N}$ is a decreasing sequence, we have
$$\Delta_{n+1}^y = \bigcup_{k \in \N} \bigcap_{m \geq k} h_m \Delta_{n+1}^x \subseteq \bigcup_{k \in \N} \bigcap_{m \geq k} h_m \Delta_n^x = \Delta_n^y.$$

(ii). This follows immediately from clause (ii) of the previous proposition (with $B = F_n F_n^{-1}$).

(iii). By referring back to Definition \ref{DEFN BP GUIDE} we see that the property of being guided by a growth sequence only depends on the sets $(F_n)_{n \in \N}$, $(H_n)_{n \in \N}$, and $(D^n_k)_{n \geq k \in \N}$.

(iv). Let $g, h \in \bigcap_{n \in \N} \Delta_n^y$. Since the $\Delta_k^x$-translates of $F_k$ are maximally disjoint within $G$, there is $\psi \in \Delta_k$ with $h^{-1} g \in \psi F_k F_k^{-1}$. Since $(\Delta_n^x, F_n)_{n \in \N}$ is centered and directed, by clause (iv) of Lemma \ref{STRONG BP LIST} there is $n \geq k$ with $\psi F_k \subseteq F_n$. Then $h^{-1} g \in F_n F_k^{-1}$. By clause (i) of Lemma \ref{STRONG BP LIST}, $F_k \subseteq F_n$ and therefore $h^{-1} g \in F_n F_n^{-1}$. It follows that $g F_n \cap h F_n \neq \varnothing$. Since $g, h \in \Delta_n^y$ we must have $g = h$.

(v). Suppose $g \in \Delta_n^y$ for all $n \in \N$. Fix $k \in \N$. If $\gamma \in \Delta_k^x$, then since $(\Delta_n^x, F_n)_{n \in \N}$ is centered and directed, there is $n > k$ with $\gamma \in F_n$. In particular, $\gamma \in D^n_k$. Thus $\Delta_k^x \subseteq \bigcup_{n \geq k} D^n_k$. On the other hand, $\bigcup_{n \geq k} D^n_k \subseteq \Delta_k^x$ by clause (i) of Lemma \ref{BP LIST} (since $(\Delta_n^x, F_n)_{n \in \N}$ is centered). Thus $\Delta_k^x = \bigcup_{n \geq k} D^n_k$. Again by clause (i) of Lemma \ref{BP LIST} we have
$$g \Delta_k^x = g \bigcup_{n \geq k} D^n_k \subseteq \Delta_k^y.$$

(vi). Suppose $g \in \Delta_n^y$ for all $n \in \N$. Fix $n \in \N$. By (v) we have $g \Delta_n^x \subseteq \Delta_n^y$. Let $\gamma \in \Delta_n^y$. Since the $\Delta_n^x$-translates of $F_n$ are maximally disjoint, there is $\psi \in \Delta_n^x$ with $\psi F_n \cap g^{-1} \gamma F_n \neq \varnothing$. Thus $g \psi \in g \Delta_n^x \subseteq \Delta_n^y$ and $g \psi F_n \cap \gamma F_n \neq \varnothing$. By the disjoint property of blueprints we must have that $\gamma = g \psi \in g \Delta_n^x$. We conclude that $g \Delta_n^x = \Delta_n^y$.
\end{proof}

We now show that every function in the closure of the orbit of a fundamental function is fundamental.

\begin{prop}
Let $G$ be a countably infinite group, let $(\Delta_n, F_n)_{n \in \N}$ be a blueprint, and let $x \in 2^G$ be fundamental with respect to this blueprint. Then every $y \in \overline{[x]}$ is fundamental with respect to $(\Delta_n^y, F_n)_{n \in \N}$.
\end{prop}

\begin{proof}
Let $1_G = g_0, g_1, \ldots$ be the fixed enumeration of $G$ used in defining the metric $d$ on $2^G$. For each $n \geq 1$, let $V_n \subseteq \gamma_n F_{n-1} \subseteq F_n$ be the test region for the simple $\Delta_n$ membership test admitted by $x$ (clause (ii) of Theorem \ref{FM}), and let $P_n \in 2^{V_n}$ be the corresponding test function. Let $(h_n)_{n \in \N}$ be such that $y = \lim_{n \rightarrow \infty} h_n \cdot x$.

Since $x$ is fundamental with respect to $(\Delta_n, F_n)_{n \in \N}$, the restriction of $x$ to $G - \bigcup_{n \geq 1} \Delta_n \Lambda_n b_{n-1}$, call this function $x'$, is canonical with respect to $(\Delta_n, F_n)_{n \in \N}$. By the definition of canonical, this means that there is a locally recognizable function $R: A \rightarrow 2$ such that $x'$, $R$, and $(\Delta_n, F_n)_{n \in \N}$ satisfy the conclusions of Theorem \ref{FM}. Let $y'$ be the restriction of $y$ to $G - \bigcup_{n \geq 1} \Delta_n^y \Lambda_n b_{n-1}$. If we show that $y'$, $R$, and $(\Delta_n^y, F_n)_{n \in \N}$ satisfy the conclusions of Theorem \ref{FM}, then $y'$ will be canonical with respect to $(\Delta_n^y, F_n)_{n \in \N}$ and hence $y$ will be fundamental with respect to $(\Delta_n^y, F_n)_{n \in \N}$. So we proceed to check the numbered conclusions of Theorem \ref{FM}.

(i). Let $\gamma \in \Delta_1^y$. Let $r \in \N$ be such that $\gamma \gamma_1 F_0 \subseteq \{g_0, g_1, \ldots, g_r\}$, and let $m \in \N$ be such that $\gamma \in h_m \Delta_1^x$ and $d(h_m \cdot x, y) < 2^{-r}$. Then
$$y(\gamma \gamma_1 f) = (h_m \cdot x)(\gamma \gamma_1 f) = x(h_m^{-1} \gamma \gamma_1 f) = R(f)$$
for all $f \in F_0$ (the last equality follows since $x$ is fundamental). We notice that $\Delta_1^y \gamma_1 F_0$ is disjoint from $\bigcup_{n \geq 1} \Delta_n^y \Lambda_n b_{n-1}$ since these sets are disjoint for any pre-blueprint, as shown in the proof of Theorem \ref{FM}. Thus $y'(\gamma \gamma_1 f) = R(f)$ for all $\gamma \in \Delta_1^y$ and $f \in F_0$.

(ii). This is clear by the definition of $\Delta_n^y$ for $n \geq 1$.

(iii). By definition, $G - \dom(y') = \bigcup_{n \geq 1} \Delta_n^y \Lambda_n b_{n-1}$. In the proof of Theorem \ref{FM}, it was shown that this union is disjoint for all pre-blueprints.

(iv). Notice that $\Delta_n^x \Lambda_n b_{n-1} \subseteq \Delta_1^x b_1 \cup \Delta_1^x \Lambda_1 \subseteq \Delta_1^x D_0^1$ for $n \geq 1$. If $g \in \dom(x')$ and $x'(g) = R(1_G)$ then $g \in \Delta_1^x (\gamma_1 F_0 \cup D_0^1)$ since $x'$ is canonical. If $g \in G - \dom(x')$ then $g \in \Delta_n^x \Lambda_n b_{n-1}$ for some $n \geq 1$ and hence $g \in \Delta_1^x D_0^1$. Thus for all $g \in G$, $x(g) = R(1_G)$ implies $g \in \Delta_1^x (\gamma_1 F_0 \cup D_0^1)$. Suppose $g \in G$ satisfies $y(g) = R(1_G)$. Let $r \in \N$ be such that $g F_1^{-1} F_1 \subseteq \{g_0, g_1, \ldots, g_r\}$ and let $m \in \N$ be such that $d(h_m \cdot x, y) < 2^{-r}$. Then $x(h_m^{-1} g) = (h_m \cdot x)(g) = y(g) = R(1_G)$ so $h_m^{-1} g \in \Delta_1^x (\gamma_1 F_0 \cup D_0^1)$. Let $f \in \gamma_1 F_0 \cup D_0^1$ be such that $h_m^{-1} g f^{-1} \in \Delta_1^x$. Then for all $v \in V_1$
$$y(g f^{-1} v) = (h_m \cdot x)(g f^{-1} v) = x(h_m^{-1} g f^{-1} v) = P_1(v).$$
So $g f^{-1} \in \Delta_1^y$ and $g \in \Delta_1^y f \subseteq \Delta_1^y (\gamma_1 F_0 \cup D_0^1)$.

(v). This follows abstractly from (iii) for any pre-blueprint, as shown in the proof of Theorem \ref{FM}.

(vi). Fix $n \geq 1$, $\gamma, \sigma \in \Delta_n^y$, and
$$f \in F_n - \{a_n, b_n\} - \bigcup_{1 \leq k \leq n} D^n_k \Lambda_k b_{k-1}.$$
Let $m \in \N$ be such that $\gamma, \sigma \in h_m \Delta_n^x$, $y(\gamma f) = (h_m \cdot x)(\gamma f)$, and $y(\sigma f) = (h_m \cdot x)(\sigma f)$. Then $h_m^{-1} \gamma, h_m^{-1} \sigma \in \Delta_n^x$ so $x(h_m^{-1} \gamma f) = x(h_m^{-1} \sigma f)$ since $x$ is fundamental. It follows that
$$y(\gamma f) = (h_m \cdot x)(\gamma f) = x(h_m^{-1} \gamma f) = x(h_m^{-1} \sigma f) = (h_m \cdot x)(\sigma f) = y(\sigma f).$$

(vii). This follows abstractly from (v) and (vi) for any pre-blueprint, as shown in the proof of Theorem \ref{FM}.
\end{proof}

Conclusion (vi) of Lemma \ref{SUBFLOW BP LIST} motivates the following definition.

\begin{definition} \index{$x$-regular} \index{$x$-centered}
Let $G$ be a countably infinite group, let $(\Delta_n, F_n)_{n \in \N}$ be a centered, directed, and maximally disjoint blueprint, and let $x \in 2^G$ be fundamental with respect to this blueprint. A function $y \in \overline{[x]}$ is called \emph{$x$-regular} if for some $g \in G$
$$g \in \bigcap_{n \in \N} \Delta_n^y.$$
If $g = 1_G$, then $y$ is called \emph{$x$-centered}.
\end{definition}

Notice that the group element $g$ in the previous definition must be unique by clause (iv) of Lemma \ref{SUBFLOW BP LIST}. Also, if $y \in \overline{[x]}$ is $x$-regular then every element in the orbit of $y$ is $x$-regular, and if $g \in \bigcap_{n \in \N} \Delta_n^y$ then $g^{-1} \cdot y$ is the unique $x$-centered element in the orbit of $y$.

The next lemma presents a nontrivial way of testing when $y \in \overline{[x]}$ is $x$-regular, at least in the case of blueprints which are centered and guided by a growth sequence. The precise numbers and combinations of $F_n$'s appearing in this lemma are ad-hoc; this lemma will be used for a specific purpose in the final section of this chapter. If needed, one could find similar tests to the one below.

\begin{lem} \label{LEM REGT}
Let $G$ be a countably infinite group, let $(\Delta_n, F_n)_{n \in \N}$ be a centered blueprint guided by a growth sequence $(H_n)_{n \in \N}$, let $x \in 2^G$ be fundamental with respect to this blueprint, and let $y \in \overline{[x]}$. Then $y$ is $x$-regular if and only if for all but finitely many $n \equiv 1 \mod 10$ the set $F_n F_n^{-1} F_{n+4}^{-1} \cap \Delta_{n+17}^y$ is nonempty.
\end{lem}

\begin{proof}
First suppose that $y$ is $x$-regular with $\gamma \in \Delta_n^y$ for all $n \in \N$. Recall that by the definition of a growth sequence we have $G = \bigcup_{n \in \N} H_n$. Thus there is $N \in \N$ with $\gamma \in H_N$. By clause (ii) of Lemma \ref{LEM GUIDED BP}, $\gamma \in H_N \subseteq F_{N+2} F_0^{-1}$. Since $(\Delta_n^x, F_n)_{n \in \N}$ is centered, $(F_n)_{n \in \N}$ is increasing. So it follows that $\gamma \in F_n F_n^{-1} F_{n+4}^{-1} \cap \Delta_{n+17}^y$ for every $n \geq N+2$.

Now suppose $y \in \overline{[x]}$ and $m \equiv 1 \mod 10$ satisfy $F_n F_n^{-1} F_{n+4}^{-1} \cap \Delta_{n+17}^y \neq \varnothing$ for all $n \geq m$ congruent to $1$ modulo $10$. Let $n \geq m$ be congruent to $1$ modulo $10$ and fix $\gamma \in F_n F_n^{-1} F_{n+4}^{-1} \cap \Delta_{n+17}^y$. We will show $\gamma \in \Delta_{(n+10)+17}^y$. Since $n + 10 > m$, our assumption on $y$ gives
$$1_G \in \Delta_{n+27}^y F_{n+14} F_{n+10} F_{n+10}^{-1}.$$
By making repeated uses of clause (iii) of Definition \ref{DEFN GROWTH}, clause (ii) of Definition \ref{DEFN BP GUIDE},  and clause (ii) of Lemma \ref{LEM GUIDED BP}, we have
$$\gamma \in F_n F_n^{-1} F_{n+4}^{-1} \subseteq H_n H_{n+1} H_{n+4}^{-1} \subseteq H_{n+5}$$
$$\subseteq F_{n+7} F_{n+7}^{-1} \subseteq \Delta_{n+27}^y F_{n+14} F_{n+10} F_{n+10}^{-1} F_{n+7} F_{n+7}^{-1}$$
$$\subseteq \Delta_{n+27}^y H_{n+14} H_{n+10} H_{n+11} H_{n+12} H_{n+13} \subseteq \Delta_{n+27}^y H_{n+15}$$
By clause (i) of Lemma \ref{LEM GUIDED BP}, clause (i) of Lemma \ref{STRONG BP LIST}, and clause (i) of Lemma \ref{SUBFLOW BP LIST} we have that $(\Delta_n^y)_{n \in \N}$ is a decreasing sequence. By clause (iii) of Lemma \ref{SUBFLOW BP LIST} $(\Delta_n^y, F_n)_{n \in \N}$ is also guided by $(H_n)_{n \in \N}$. Since $\gamma \in \Delta_{n+17}^y \subseteq \Delta_{n+15}^y$ we have
$$\gamma \in \gamma F_{n+15} \cap \Delta_{n+27}^y H_{n+15}$$
and therefore by clause (iii) of Lemma \ref{LEM GUIDED BP} $\gamma \in \gamma F_{n+15} \subseteq \Delta_{n+27}^y F_{n+17}$. However, $\gamma \in \Delta_{n+17}^y$, the $\Delta_{n+17}^y$-translates of $F_{n+17}$ are disjoint, and $\Delta_{n+27}^y \subseteq \Delta_{n+17}^y$. So we must have that $\gamma \in \Delta_{n+27}^y$. In particular,
$$\gamma \in F_{n+10} F_{n+10}^{-1} F_{n+14}^{-1} \cap \Delta^y_{n +27}.$$
We can repeat the above argument and apply induction to conclude that $\gamma \in \Delta_{k+17}^y$ for all $k \geq n$ congruent to $1$ modulo $10$. Therefore $\gamma \in \Delta_n^y$ for all $n \in \N$ since $(\Delta_n^y)_{n \in \N}$ is decreasing. We conclude that $y$ is $x$-regular.
\end{proof}

\section{Pre-minimality}

This section is devoted to studying the significance of following property in the context of fundamental functions.

\begin{definition} \index{pre-minimal} \index{minimal!pre-minimal}
Let $G$ be a countably infinite group, and let $c \in 2^{\subseteq G}$. The function $c$ is called \emph{pre-minimal} if there is a minimal $c' \in 2^G$ extending $c$.
\end{definition}

Our goal in this section is to provide several ways of testing when fundamental functions are pre-minimal. It would be nice if pre-minimality was highly reliant on the structure of the blueprint used, however this turns out not to be the case. Nevertheless, with increased restrictions on the blueprint this comes closer to being the case. This section consists of several characterizations of pre-minimal functions. We begin by assuming as little as possible about the fundamental function and its blueprint, and as we proceed through the section we place more and more restrictions on the blueprint in order to arrive at nicer and nicer characterizations of pre-minimality.

Recall from clause (vii) of Lemma \ref{STRONG BP LIST} that if a blueprint $(\Delta_n, F_n)_{n \in \N}$ is directed, then $\bigcap_{n \in \N} \Delta_n b_n$ is either empty or a singleton.

\begin{lem}
Let $G$ be a countably infinite group, let $(\Delta_n, F_n)_{n \in \N}$ be a blueprint, and let $c \in 2^{\subseteq G}$ be fundamental with respect to this blueprint and have the property that $c$ is constant on $\bigcap_{n \in \N} \Delta_n b_n$. Then $c$ is pre-minimal if and only if either $c_0$ or $c_1$ is minimal, where $c_i \in 2^G$ is the function which extends $c$ and satisfies $c_i(g) = i$ for $g \not\in \dom(c)$.
\end{lem}

\begin{proof}
If $\bigcap_{n \in \N} \Delta_n b_n \neq \varnothing$, then let $i$ be the constant value that $c$ takes on this set. If $\bigcap_{n \in \N} \Delta_n b_n = \varnothing$, then let $i$ be either $0$ or $1$.

Clearly if either $c_0$ or $c_1$ is minimal then $c$ is pre-minimal. So assume that $c$ is pre-minimal. We will show that $c_i$ is minimal. Since $c$ is pre-minimal, there is a minimal $d \in 2^G$ extending $c$. We will use Lemma \ref{lem:minimallemma} to show that $c_i$ is minimal. Let $A \subseteq G$ be finite. Recall that $(\Delta_n b_n)_{n \in \N}$ is a decreasing sequence. Since $A$ is finite we may fix $k \in \N$ so that
$$\forall a \in A \ a \in \Delta_k b_k \Longrightarrow a \in \bigcap_{n \in \N} \Delta_n b_n.$$
Set
$$B = \bigcup_{0 \leq n \leq k} A F_n^{-1} F_n.$$
Since $d$ is minimal, there is a finite $T \subseteq G$ so that for all $g \in G$ there is $t \in T$ with $d(gtb) = d(b)$ for all $b \in B$.

Fix an arbitrary $g \in G$, and let $t \in T$ be such that $d(gtb) = d(b)$ for all $b \in B$. We will show that $c_i(gta) = c_i(a)$ for all $a \in A$.

Fix $a \in A$, $n \leq k$, and $f \in F_n$. We claim that $a \in \Delta_n f$ if and only if $gta \in \Delta_n f$. First suppose $a \in \Delta_n f$. Say $a = \gamma f$ with $\gamma \in \Delta_n$. Then
$$\gamma F_n \subseteq a F_n^{-1} F_n \subseteq A F_k^{-1} F_k = B.$$
Since $c$ admits a $\Delta_n$ membership test with test region a subset of $F_n$ and $d(gtb) = d(b)$ for all $b \in B$, it follows that $gtaf^{-1} \in \Delta_n$ and thus $gta \in \Delta_n f$. The argument that $gta \in \Delta_n f$ implies $a \in \Delta_n f$ is identical.

A consequence of the previous paragraph is that for $a \in A$ and $n \leq k$,
$$a \in \Delta_n \Theta_n b_{n-1} \Longleftrightarrow gta \in \Delta_n \Theta_n b_{n-1}.$$

Set $A' = A - \Delta_k b_k$. For $n > k$, $\Delta_n \Theta_n b_{n-1}$ is contained in $\Delta_k b_k$. Since $A' \cap \Delta_k b_k = \varnothing$, we have $gtA' \cap \Delta_k b_k = \varnothing$ as well. As
$$\dom(c) = G - \bigcup_{n \geq 1} \Delta_n \Theta_n b_{n-1}$$
it follows that
$$(gtA') \cap \dom(c) = gt (A' \cap \dom(c)).$$
Therefore $c_i(gta) = i = c_i(a)$ for all $a \in A' - \dom(c)$. Since $d$ extends $c$ and $d(gta) = d(a)$ for all $a \in A$, we have $c(gta) = c(a)$ for all $a \in A' \cap \dom(c)$. Putting these together we have $c_i(gta) = c_i(a)$ for all $a \in A'$. If $a \in A - A'$ then $\bigcap_{n \in \N} \Delta_n b_n \neq \varnothing$ and we have $d(gta) = d(a) = c(a) = i$. In general for $h \in G$, $d(h) = i$ implies $c_i(h)=i$. So $c_i(gta) = c_i(a) = i$ for all $a \in A-A'$. We conclude that $c_i(gta) = c_i(a)$ for all $a \in A$. Thus $c_i$ is minimal as claimed.
\end{proof}

Recall from the proof of clause (viii) of Lemma \ref{STRONG BP LIST} that if $(\Delta_n, F_n)_{n \in \N}$ is a centered and directed blueprint and $\beta_n \neq 1_G$ for infinitely many $n \in \N$, then $\bigcap_{n \in \N} \Delta_n b_n = \varnothing$.

\begin{cor} \label{PREMIN EMPTY}
Let $G$ be a countably infinite group, let $(\Delta_n, F_n)_{n \in \N}$ be a blueprint with $\bigcap_{n \in \N} \Delta_n b_n = \varnothing$, and let $c \in 2^{\subseteq G}$ be fundamental with respect to this blueprint. Then $c$ is pre-minimal if and only if for every finite $A \subseteq G$ there is a finite $T \subseteq G$ with the property that for every $g \in G$ there is $t \in T$ satisfying
$$(gtA) \cap \dom(c) = gt(A \cap \dom(c)), \text{ and}$$
$$\forall a \in A \cap \dom(c) \ c(gta) = c(a).$$
\end{cor}

\begin{proof}
If $c$ has the stated property then the function $c_0$ (as well as $c_1$) defined in the previous lemma is clearly minimal and thus $c$ is pre-minimal. Conversely, suppose $c$ is pre-minimal. If we let $A \subseteq G$ be finite and follow the argument in the proof of the previous lemma, then $A' = A$ and we find that $c$ satisfies the condition stated above.
\end{proof}

\begin{lem}
Let $G$ be a countably infinite group, and let $(\Delta_n, F_n)_{n \in \N}$ be a directed blueprint with $\bigcap_{n \in \N} \Delta_n b_n = \varnothing$ and with the property that for some $k \in \N$ the $\Delta_k$-translates of $F_k$ are maximally disjoint. Let $c \in 2^{\subseteq G}$ be fundamental with respect to this blueprint. Then $c$ is pre-minimal if and only if for every finite $A \subseteq G$ there is $n > k \in \N$ and $g \in G$ so that for all $\gamma \in \Delta_n$ there is $\lambda \in D^n_k$ satisfying:
$$\gamma \lambda g [A \cap \dom(c)] = (\gamma \lambda g A) \cap \dom(c)$$
and
$$\forall a \in A \cap \dom(c) \ c(\gamma \lambda g a) = c(a).$$
\end{lem}

\begin{proof}
First suppose that $c$ has the property stated above. Define $c': G \rightarrow 2$ to be the function which extends $c$ and satisfies $c'(g) = 0$ for all $g \not\in \dom(c)$. We will use Lemma \ref{lem:minimallemma} to show that $c'$ is minimal. So fix a finite $A \subseteq G$, and let $n > k$ and $g \in G$ be so as to satisfy the above stated condition satisfied by $c$. Let $B \subseteq G$ be finite such that $\Delta_n B = G$. Set $T = B^{-1} D^n_k g$. Now let $h \in G$ be arbitrary. Since $\Delta_n B = G$, there is $b \in B^{-1}$ with $h b \in \Delta_n$. It follows there is $\lambda \in D^n_k$ with
$$h b \lambda g [A \cap \dom(c)] = (h b \lambda g A) \cap \dom(c)$$
and
$$\forall a \in A \cap \dom(c) \ c(h b \lambda g a) = c(a).$$
It immediately follows that $c'(h b \lambda g a) = c'(a)$ for all $a \in A$. Also we have $b \lambda g \in T$. We conclude $c'$ is minimal and hence $c$ is pre-minimal.

Now assume that $c$ is pre-minimal. Let $A \subseteq G$ be finite. By enlarging $A$ if necessary, we may assume that $\psi F_k \subseteq A$ for some $\psi \in \Delta_k$. Set $g = \psi^{-1}$. By Corollary \ref{PREMIN EMPTY} there is a finite $T \subseteq G$ so that for all $h \in G$ there is $t \in T$ with
$$h t (A \cap \dom(c)) = (h t A) \cap \dom(c), \text{ and}$$
$$\forall a \in A \cap \dom(c) \ c(h t a) = c(a).$$
By clause (iii) of Lemma \ref{STRONG BP LIST}, there is $n > k$ and $\sigma \in \Delta_n$ with
$$T A F_k F_k F_k F_k^{-1} \cap \Delta_k \subseteq \sigma F_n.$$

Now let $\gamma \in \Delta_n$ be arbitrary. Let $t \in T$ be such that
$$\gamma \sigma^{-1} t (A \cap \dom(c)) = (\gamma \sigma^{-1} t A) \cap \dom(c), \text{ and}$$
$$\forall a \in A \cap \dom(c) \ c(\gamma \sigma^{-1} t a) = c(a).$$
Since $\psi F_k \subseteq A$ and $c$ has a $\Delta_k$ membership test, it must be that $\gamma \sigma^{-1} t \psi \in \Delta_k$. By clause (v) of Lemma \ref{STRONG BP LIST} we have that
$$\Delta_k \cap \gamma \sigma^{-1} T A = \gamma \sigma^{-1} (\Delta_k \cap T A) \subseteq \gamma D^n_k.$$
Therefore $\gamma \sigma^{-1} t \psi \in \gamma D^n_k$, so $\lambda = \sigma^{-1} t \psi \in D^n_k$. We have $\gamma \lambda g = \gamma \sigma^{-1} t$ so
$$\gamma \lambda g (A \cap \dom(c)) = (\gamma \lambda g A) \cap \dom(c), \text{ and}$$
$$\forall a \in A \cap \dom(c) \ c(\gamma \lambda g a) = c(a).$$
Thus $c$ has the claimed property.
\end{proof}

In the previous lemma, we only assume that the blueprint is directed and that the $\Delta_k$-translates of $F_k$ are maximally disjoint in order to apply the conclusion of clause (v) of Lemma \ref{STRONG BP LIST}. In all of the following results in this section we assume that the blueprint is directed and that the $\Delta_i$-translates of $F_i$ are maximally disjoint for either $i = 0$ or $i = 1$, but it is still the case that all we really need is to be able to apply clause (v) of Lemma \ref{STRONG BP LIST}. Clause (v) of Lemma \ref{STRONG BP LIST} therefore plays a special role in this section and the next. Really we are using clause (v) to get a more descriptive version of clause (vi) of Lemma \ref{STRONG BP LIST} where $0$ is replaced by $i$ (where $i \in \{0,1\}$ is such that the $\Delta_i$-translates of $F_i$ are maximally disjoint). Thus, the minimal property mentioned in (vi) seems to be related to the pre-minimality of fundamental functions.

\begin{lem} \label{LEM PREMF}
Let $G$ be a countably infinite group, and let $(\Delta_n, F_n)_{n \in \N}$ be a directed blueprint with $\bigcap_{n \in \N} \Delta_n b_n = \varnothing$ and with the property that the $\Delta_i$-translates of $F_i$ are maximally disjoint for either $i = 0$ or $i = 1$. Let $c \in 2^{\subseteq G}$ be fundamental with respect to this blueprint. Then $c$ is pre-minimal if and only if for every $k \geq 1$ and $\psi \in \Delta_k$ there is $n > k$ so that for all $\gamma \in \Delta_n$ there is $\lambda \in D^n_k$ satisfying:
$$(\gamma \lambda)^{-1} [(\gamma \lambda F_k) \cap \dom(c)] = \psi^{-1} [(\psi F_k) \cap \dom(c)]$$
and
$$\forall f \in \psi^{-1} [(\psi F_k) \cap \dom(c)] \ c(\gamma \lambda f) = c(\psi f).$$
\end{lem}

\begin{proof}
Notice that the last two expressions are equivalent to $[(\gamma \lambda)^{-1} \cdot c] \res F_k = [ \psi^{-1} \cdot c] \res F_k$. Fix $i \in \{0, 1\}$ so that the $\Delta_i$-translates of $F_i$ are maximally disjoint. First assume that $c$ has the stated property. We will show that $c$ is pre-minimal by applying Corollary \ref{PREMIN EMPTY}. Let $A \subseteq G$ be finite. By directedness, there is $k \geq 1$ and $\psi \in \Delta_k$ with
$$\Delta_i \cap A F_1^{-1} F_1 F_i F_i F_i^{-1} \subseteq \psi F_k.$$
Notice that the set on the left is necessarily contained in $\psi D^k_i$. Also notice that $(\Delta_1 \cap A F_1^{-1}) a_1 a_i^{-1}$ is contained in the set on the left, so by the coherent property of blueprints
$$(\Delta_1 \cap A F_1^{-1}) F_1 \subseteq \psi F_k.$$
By assumption, there is $n > k$ so that for all $\gamma \in \Delta_n$ there is $\lambda \in D^n_k$ satisfying $[(\gamma \lambda)^{-1} \cdot c] \res F_k = [ \psi^{-1} \cdot c] \res F_k$. Let $B \subseteq G$ be finite with $\Delta_n B = G$. Set $T = B^{-1} D^n_k \psi^{-1}$ and let $g \in G$ be arbitrary. Then there is $b \in B^{-1}$ with $g b = \gamma \in \Delta_n$. Let $\lambda \in D^n_k$ be such that $[(\gamma \lambda)^{-1} \cdot c] \res F_k = [\psi^{-1} \cdot c] \res F_k$. Set $t = b \lambda \psi^{-1} \in T$ and notice $\gamma \lambda = g t \psi$. So $[(g t \psi)^{-1} \cdot c] \res F_k = [ \psi^{-1} \cdot c] \res F_k$ and therefore
$$[(g t)^{-1} \cdot c] \res \psi F_k = c \res \psi F_k.$$
We will be done if we can show that $[(gt)^{-1} \cdot c] \res A = c \res A$. For this it suffices to fix $a \in A - \psi F_k$ and show that $a, g t a \in \dom(c)$ and $c(g t a) = c(a)$. Since $(\Delta_1 \cap A F_1^{-1}) F_1 \subseteq \psi F_k$, we must have that $a \not\in \Delta_1 F_1$. By our choice of $k$ and clause (v) of Lemma \ref{STRONG BP LIST} we have that
$$\Delta_1 \cap \gamma \lambda \psi^{-1} a F_1^{-1} = \gamma \lambda \psi^{-1} (\Delta_1 \cap a F_1^{-1}) = \varnothing.$$
So $\gamma \lambda \psi^{-1} a = g t a \not\in \Delta_1 F_1$. It follows from Definition \ref{DEF FUNDF} and conclusions (iii) and (iv) of Theorem \ref{FM} that $a, gta \in \dom(c)$ and $c(gta) = c(a)$.

Now assume that $c$ is pre-minimal. Fix $k \geq 1$ and $\psi \in \Delta_k$. By Corollary \ref{PREMIN EMPTY}, there is a finite $T \subseteq G$ so that for all $g \in G$ there is $t \in T$ with $[(g t)^{-1} \cdot c] \res \psi F_k = c \res \psi F_k$. By directedness, there is $n \geq k$ and $\sigma \in \Delta_n$ with
$$\Delta_i \cap T \psi F_k F_i F_i F_i^{-1} \subseteq \sigma F_n.$$
Notice that the set on the left is necessarily contained in $\sigma D^n_i$. Also notice that $(\Delta_k \cap T \psi) a_k a_i^{-1}$ is contained in the set on the left, so by the coherent property of blueprints we have
$$\Delta_k \cap T \psi \subseteq \sigma D^n_k.$$
Now let $\gamma \in \Delta_n$ be arbitrary. Let $t \in T$ be such that $[(\gamma \sigma^{-1} t)^{-1} \cdot c] \res \psi F_k = c \res \psi F_k$. Then $[(\gamma \sigma^{-1} t \psi)^{-1} c] \res F_k = [\psi^{-1} \cdot c] \res F_k$ so it suffices to show that $\sigma^{-1} t \psi \in D^n_k$. Since $c$ has a $\Delta_k$ membership test and $\psi \in \Delta_k$, it must be that $\gamma \sigma^{-1} t \psi \in \Delta_k$. By clause (v) of Lemma \ref{STRONG BP LIST} we have that
$$\Delta_k \cap \gamma \sigma^{-1} T \psi = \gamma \sigma^{-1} (\Delta_k \cap T \psi) \subseteq \gamma \sigma^{-1} \sigma D^n_k = \gamma D^n_k.$$
Therefore $\gamma \sigma^{-1} t \psi \in \gamma D^n_k$, so $\lambda = \sigma^{-1} t \psi \in D^n_k$. Thus $c$ has the stated property.
\end{proof}

\begin{cor} \label{STRONG PREMIN}
Let $G$ be a countably infinite group, and let $(\Delta_n, F_n)_{n \in \N}$ be a directed and centered blueprint with $\bigcap_{n \in \N} \Delta_n b_n = \varnothing$ and with the property that the $\Delta_i$-translates of $F_i$ are maximally disjoint for either $i = 0$ or $i = 1$. Let $c \in 2^{\subseteq G}$ be fundamental with respect to this blueprint. Then $c$ is pre-minimal if and only if for every $k \geq 1$ there is $n > k$ so that for all $\gamma \in \Delta_n$ there is $\lambda \in D^n_k$ satisfying:
$$\gamma \lambda [F_k \cap \dom(c)] =(\gamma \lambda F_k) \cap \dom(c)$$
and
$$\forall f \in F_k \cap \dom(c) \ c(\gamma \lambda f) = c(f).$$
\end{cor}

\begin{proof}
If $c$ is pre-minimal, then by picking any $k \geq 1$, setting $\psi = 1_G \in \Delta_k$, and applying the previous lemma we see that $c$ has the stated property. Now assume that $c$ has the stated property. We will apply the previous lemma to show that $c$ is pre-minimal. Pick $k \in \N$ and $\psi \in \Delta_k$. By clause (iv) of Lemma \ref{STRONG BP LIST}, there is $m \geq k$ with $\psi F_k \subseteq F_m$. Let $n > m$ be such that for all $\gamma \in \Delta_n$ there is $\lambda \in D^n_m$ with $[(\gamma \lambda)^{-1} \cdot c] \res F_m = c \res F_m$. Now let $\gamma \in \Delta_n$ be arbitrary, and let $\lambda \in D^n_m$ be such that $[(\gamma \lambda)^{-1} \cdot c] \res F_m = c \res F_m$. In particular, $[(\gamma \lambda)^{-1} \cdot c] \res F_k = c \res F_k$ (clause (i) of Lemma \ref{STRONG BP LIST}) and hence $[(\gamma \lambda \psi)^{-1} \cdot c] \res \psi F_k = [\psi^{-1} \cdot c] \res \psi F_k$. Notice that $\psi \in \Delta_k \cap F_m = D^m_k$ and therefore $\lambda \psi \in D^n_k$. So we have shown that $c$ satisfies the condition stated in the previous lemma. We conclude that $c$ is pre-minimal.
\end{proof}

The following proposition is especially useful.

\begin{prop} \label{CANONICAL PREMIN}
Let $G$ be a countably infinite group,  and let $(\Delta_n, F_n)_{n \in \N}$ be a directed blueprint with $\bigcap_{n \in \N} \Delta_n a_n = \bigcap_{n \in \N} \Delta_n b_n = \varnothing$ and with the property that the $\Delta_i$-translates of $F_i$ are maximally disjoint for either $i = 0$ or $i = 1$. Then any $c \in 2^{\subseteq G}$ which is canonical with respect to this blueprint is pre-minimal.
\end{prop}

\begin{proof}
Let $(\Delta_n, F_n)_{n \in \N}$, $i \in \{0, 1\}$, and $c \in 2^{\subseteq G}$ be as stated. We will apply Lemma \ref{LEM PREMF} to show that $c$ is pre-minimal. So fix $k \geq 1$ and $\psi \in \Delta_k$. Since $\bigcap_{n \in \N} \Delta_n a_n = \bigcap_{n \in \N} \Delta_n b_n = \varnothing$ and the blueprint is directed, there is $n \geq k$ and $\sigma \in \Delta_n$ with $\psi F_k \subseteq \sigma F_n$ and $\psi F_k \cap \Delta_n \{a_n, b_n\} = \varnothing$. So $\psi F_k \subseteq \sigma (F_n - \{a_n, b_n\})$. Notice that $\sigma^{-1} \psi \in D^n_k$. Now let $\gamma \in \Delta_n$ be arbitrary and set $\lambda = \sigma^{-1} \psi$. By conclusion (vii) of Theorem \ref{FM}
$$\gamma^{-1} [(\gamma F_n - \{\gamma a_n, \gamma b_n\}) \cap \dom(c) ] = \sigma^{-1} [(\sigma F_n - \{\sigma a_n, \sigma b_n\}) \cap \dom(c) ]$$
and
$$\forall f \in \sigma^{-1} [ (\sigma F_n - \{\sigma a_n, \sigma b_n\}) \cap \dom(c)] \ c(\gamma f) = c(\sigma f).$$
Since $\sigma^{-1} \psi F_k \subseteq F_n - \{a_n, b_n\}$ and $\psi = \sigma \lambda$, it follows that
$$(\gamma \lambda)^{-1} [(\gamma \lambda F_k) \cap \dom(c) ] = \psi^{-1} [(\psi F_k) \cap \dom(c) ]$$
and
$$\forall f \in \psi^{-1} [ (\psi F_k) \cap \dom(c)] \ c(\gamma \lambda f) = c(\psi f).$$
By Lemma \ref{LEM PREMF} we conclude that $c$ is pre-minimal.
\end{proof}

\section{$\Delta$-minimality}

In the previous section we saw that even for the most special blueprints the best characterization for pre-minimality (best in terms of relating pre-minimality with blueprint structure) we could get was Corollary \ref{STRONG PREMIN}. Since blueprints are such an important feature of fundamental functions, we want to study a property which is more closely related to blueprint structure. We also want this property to be stronger than pre-minimality.

\begin{definition} \label{DEF DMIN}\index{$\Delta$-minimal}\index{minimal!$\Delta$-minimal}
Let $G$ be a countably infinite group and let $(\Delta_n, F_n)_{n \in \N}$ be a blueprint. A function $f \in \N^{\subseteq G}$ is \emph{$\Delta$-minimal} (relative to $(\Delta_n, F_n)_{n \in \N}$) if for every finite $A \subseteq G$ there is $n \in \N$ and $\sigma \in \Delta_n$ so that for all $\gamma \in \Delta_n$
$$\gamma \sigma^{-1} A \cap \dom(f) = \gamma \sigma^{-1} ( A \cap \dom(f))$$
and
$$\forall a \in A \cap \dom(f) \ f( \gamma \sigma^{-1} a) = f(a).$$
\end{definition}

We first observe that $\Delta$-minimality is indeed stronger than pre-minimality.

\begin{lem}
Let $G$ be a countably infinite group, let $(\Delta_n, F_n)_{n \in \N}$ be a blueprint, and let $c \in 2^{\subseteq G}$ be any function. If $c$ is $\Delta$-minimal, then it is pre-minimal.
\end{lem}

\begin{proof}
Assume $c$ is $\Delta$-minimal. Let $c' \in 2^G$ be the function which extends $c$ and satisfies $c'(g) = 0$ for all $g \in G - \dom(c)$. We will show that $c'$ is minimal by applying Lemma \ref{lem:minimallemma}. So let $A \subseteq G$ be finite. Let $n \in \N$ and $\sigma \in \Delta_n$ be such that for all $\gamma \in \Delta_n$
$$\gamma \sigma^{-1} A \cap \dom(c) = \gamma \sigma^{-1} ( A \cap \dom(c))$$
and
$$\forall a \in A \cap \dom(c) \ c( \gamma \sigma^{-1} a) = c(a).$$
It follows that $c'(\gamma \sigma^{-1} a) = c'(a)$ for all $\gamma \in \Delta_n$ and $a \in A$. Let $B \subseteq G$ be finite with $\Delta_n B = G$. Set $T = B^{-1} \sigma^{-1}$ and let $g \in G$ be arbitrary. Then there is $b \in B^{-1}$ with $g b \in \Delta_n$ and hence $c'(g b \sigma^{-1} a) = c'(a)$ for all $a \in A$. We conclude that $c'$ is minimal and thus $c$ is pre-minimal.
\end{proof}

Usually we will be interested in $\Delta$-minimal functions in the context of centered and directed blueprints. As the next lemma shows, when the blueprint is centered and directed the definition for $\Delta$-minimal takes a very nice form.

\begin{lem} \label{CENTER DIRECT DMIN}
Let $G$ be a countably infinite group, let $(\Delta_n, F_n)_{n \in \N}$ be a centered and directed blueprint, and let $f \in \N^{\subseteq G}$ be any function. Then $f$ is $\Delta$-minimal if and only if for every finite $A \subseteq G$ there is $n \in \N$ so that for all $\gamma \in \Delta_n$
$$\gamma A \cap \dom(f) = \gamma (A \cap \dom(f))$$
and
$$\forall a \in A \cap \dom(f) \ f(\gamma a) = f(a).$$
\end{lem}

\begin{proof}
Clearly $f$ is $\Delta$-minimal if it satisfies the above condition since one can choose $\sigma = 1_G \in \Delta_n$ in Definition \ref{DEF DMIN}. So suppose that $f$ is $\Delta$-minimal. Let $A \subseteq G$ be finite and let $k \in \N$ and $\sigma \in \Delta_k$ be such that for all $\gamma \in \Delta_k$
$$\gamma \sigma^{-1} A \cap \dom(f) = \gamma \sigma^{-1} (A \cap \dom(f))$$
and
$$\forall a \in A \cap \dom(f) \ f(\gamma \sigma^{-1} a) = f(a).$$
By clause (iv) of Lemma \ref{STRONG BP LIST}, there is $n \geq k$ with $\sigma F_k \subseteq F_n$. In particular, $\sigma \in D^n_k$. So if $\gamma \in \Delta_n$ then $\gamma \sigma \in \Delta_k$ and since $\gamma \sigma \sigma^{-1} = \gamma$ we have
$$\gamma A \cap \dom(f) = \gamma (A \cap \dom(f))$$
and
$$\forall a \in A \cap \dom(f) \ f(\gamma a) = f(a).$$
\end{proof}

By comparing Corollary \ref{STRONG PREMIN} to the lemma below, one can see that among centered and directed blueprints the advantage to $\Delta$-minimality is that instead of having uncertainty as to which $\lambda \in D^n_k$ ``works,'' we know that specifically $1_G \in D^n_k$ ``works.''

\begin{lem} \label{LEM FDMIN}
Let $G$ be a countably infinite group, and let $(\Delta_n, F_n)_{n \in \N}$ be a centered and directed blueprint with $\bigcap_{n \in \N} \Delta_n b_n = \varnothing$ and with the property that the $\Delta_i$-translates of $F_i$ are maximally disjoint for either $i = 0$ or $i = 1$. Let $c \in 2^{\subseteq G}$ be fundamental with respect to this blueprint. Then $c$ is $\Delta$-minimal if and only if for every $k \geq 1$ there is $n > k$ with the property that for all $\gamma \in \Delta_n$
$$\gamma F_k \cap \dom(c) = \gamma (F_k \cap \dom(c))$$
and
$$\forall f \in F_k \cap \dom(c) \ c(\gamma f) = c(f).$$
\end{lem}

\begin{proof}
By taking $A = F_k$ in Lemma \ref{CENTER DIRECT DMIN}, we see that if $c$ is $\Delta$-minimal then it has the property stated above. Now assume that $c$ has the property stated above. Fix $i \in \{0, 1\}$ so that the $\Delta_i$-translates of $F_i$ are maximally disjoint. Let $A \subseteq G$ be finite. By directedness, there is $k \geq 1$ and $\psi \in \Delta_k$ with
$$\Delta_i \cap A F_1^{-1} F_1 F_i F_i F_i^{-1} \subseteq \psi F_k.$$
By assumption, there is $n > k$ so that $(\gamma^{-1} \cdot c) \res F_k = c \res F_k$ for all $\gamma \in \Delta_n$. Fix $\gamma \in \Delta_n$. We will show that $(\gamma^{-1} \cdot c) \res A = c \res A$. For this it suffices to fix $a \in A - F_k$ and show that $a, \gamma a \in \dom(c)$ and $c(\gamma a) = c(a)$. Since $\Delta_1 \cap A F_1^{-1} \subseteq F_k$, by the coherent property of blueprints we must have that $a \not\in \Delta_1 F_1$. By our choice of $k$ and clause (v) of Lemma \ref{STRONG BP LIST} we have that
$$\Delta_1 \cap \gamma a F_1^{-1} = \gamma (\Delta_1 \cap a F_1^{-1}) = \varnothing.$$
So $\gamma a \not\in \Delta_1 F_1$. It follows from Definition \ref{DEF FUNDF} and clauses (iii) and (iv) of Theorem \ref{FM} that $a, \gamma a \in \dom(c)$ and $c(\gamma a) = c(a)$.
\end{proof}

The previous lemma allows us to extend Proposition \ref{CANONICAL PREMIN} to the proposition below. The proof of this proposition is nearly identical to that of Proposition \ref{CANONICAL PREMIN} except Lemma \ref{LEM PREMF} is replaced by the lemma above. Nevertheless, we include the proof for completeness.

\begin{prop} \label{CANONICAL DMIN}
Let $G$ be a countably infinite group,  and let $(\Delta_n, F_n)_{n \in \N}$ be a centered and directed blueprint with $\bigcap_{n \in \N} \Delta_n a_n = \bigcap_{n \in \N} \Delta_n b_n = \varnothing$ and with the property that the $\Delta_i$-translates of $F_i$ are maximally disjoint for either $i = 0$ or $i = 1$. Then any $c \in 2^{\subseteq G}$ which is canonical with respect to this blueprint is $\Delta$-minimal.
\end{prop}

\begin{proof}
Let $(\Delta_n, F_n)_{n \in \N}$, $i \in \{0, 1\}$, and $c \in 2^{\subseteq G}$ be as stated. We will apply Lemma \ref{LEM FDMIN} to show that $c$ is $\Delta$-minimal. So fix $k \geq 1$. Since $\bigcap_{n \in \N} \Delta_n a_n = \bigcap_{n \in \N} \Delta_n b_n = \varnothing$, there is $n \geq k$ with $F_k \cap \Delta_n \{a_n, b_n\} = \varnothing$. Notice that $F_k \subseteq F_n$ by clause (i) of Lemma \ref{STRONG BP LIST}. So $F_k \subseteq F_n - \{a_n, b_n\}$. Now let $\gamma \in \Delta_n$ be arbitrary. By conclusion (vii) of Theorem \ref{FM}
$$\gamma^{-1} [(\gamma F_n - \{\gamma a_n, \gamma b_n\}) \cap \dom(c) ] =  (F_n - \{a_n, b_n\}) \cap \dom(c)$$
and
$$\forall f \in (F_n - \{a_n, b_n\}) \cap \dom(c) \ c(\gamma f) = c(f).$$
Since $F_k \subseteq F_n - \{a_n, b_n\}$, it follows that
$$\gamma^{-1} [(\gamma F_k) \cap \dom(c) ] = F_k \cap \dom(c)$$
and
$$\forall f \in F_k \cap \dom(c) \ c(\gamma f) = c(f).$$
By Lemma \ref{LEM FDMIN} we conclude that $c$ is $\Delta$-minimal.
\end{proof}

For the rest of this section we study properties of $\Delta$-minimal functions. We will see that $\Delta$-minimal functions behave well under unions and that the subflows they generate have nice properties. These nice properties are the reason why we define and study $\Delta$-minimality. We begin by looking at unions.

\begin{lem} \label{LEM UNION}
Let $G$ be a countably infinite group, and let $(\Delta_n, F_n)_{n \in \N}$ be a centered and directed blueprint. If $c_1$ and $c_2$ are both $\Delta$-minimal and agree on $\dom(c_1) \cap \dom(c_2)$, then $c_1 \cup c_2$ is $\Delta$-minimal as well.
\end{lem}

\begin{proof}
Let $A \subseteq G$ be finite. Let $n_1$ and $n_2$ be as in Lemma \ref{CENTER DIRECT DMIN} relative to $c_1$ and $c_2$ respectively. Since the blueprint is centered, $(\Delta_n)_{n \in \N}$ is decreasing. It is then easy to see that for $c = c_1 \cup c_2$, $n = \max(n_1, n_2)$, and $\gamma \in \Delta_n$ we have
$$\gamma A \cap \dom(c) = \gamma ( A \cap \dom(c))$$
and
$$\forall a \in A \cap \dom(c) \ c( \gamma a) = c(a).$$
Since $A$ was arbitrary, we conclude that $c = c_1 \cup c_2$ is $\Delta$-minimal.
\end{proof}

\begin{lem}
Let $G$ be a countably infinite group, and let $(\Delta_n, F_n)_{n \in \N}$ be a blueprint. If $c_1$ is fundamental with respect to this blueprint and pre-minimal, $c_2$ is $\Delta$-minimal, and $c_1$ and $c_2$ agree on $\dom(c_1) \cap \dom(c_2)$, then $c_1 \cup c_2$ is pre-minimal.
\end{lem}

\begin{proof}
Let $d \in 2^G$ be a minimal extension of $c_1$, and define $x \in 2^G$ by $x(g) = (c_1 \cup c_2)(g)$ if $g \in \dom(c_1) \cup \dom(c_2)$ and $x(g) = d(g)$ otherwise. It suffices to show that $x$ is minimal. Let $A \subseteq G$ be finite. Let $n \in \N$ and $\sigma \in \Delta_n$ be such that for all $\gamma \in \Delta_n$
$$\gamma \sigma^{-1} A \cap \dom(c_2) = \gamma \sigma^{-1} (A \cap \dom(c_2))$$
and
$$\forall a \in A \cap \dom(c_2) \ c_2(\gamma \sigma^{-1} a) = c_2(a).$$
Since $d$ is minimal, there is a finite set $T \subseteq G$ so that for all $g \in G$ there is $t \in T$ with
$$\forall h \in A \cup \sigma F_n \ d(g t h) = d(h).$$
Fix $g \in G$ and let $t \in T$ be such that the expression above is satisfied. Since $c_1 \subseteq d$ is fundamental, $d$ admits a simple $\Delta_n$-membership test with test region a subset of $F_n$. So we must have $g t \sigma \in \Delta_n$. Therefore for all $a \in A \cap \dom(c_2)$ we have $gta \in \dom(c_2)$ and
$$x(g t a) = c_2(g t a) = c_2(a) = x(a).$$
On the other hand, if $a \in A - \dom(c_2)$ then $g t a \not\in \dom(c_2)$ and hence
$$x(g t a) = d(g t a) = d(a) = x(a).$$
We concluded that $x$ is minimal and hence $c_1 \cup c_2$ is pre-minimal.
\end{proof}

\begin{lem} \label{MIN UNION}
Let $G$ be a countably infinite group, let $(\Delta_n, F_n)_{n \in \N}$ be a blueprint, and let $(c_n)_{n \in \N}$ be an increasing sequence of elements of $2^{\subseteq G}$ which are all fundamental with respect to this blueprint. If each $c_n$ is $\Delta$-minimal then $c = \bigcup_{n \in \N} c_n$ is $\Delta$-minimal as well.
\end{lem}

\begin{proof}
Let $A \subseteq G$ be finite. Let $k \in \N$ be such that
$$\forall a \in A \ a \in \Delta_k b_k \Longrightarrow a \in \bigcap_{n \in \N} \Delta_n b_n.$$
Let $m \in \N$ be such that for all $1 \leq i \leq k$
$$\Theta_i(c_m) = \Theta_i(c)$$
and
$$A \cap \dom(c_m) = A \cap \dom(c).$$
Such an $m$ exists since $(\Theta_i(c_n))_{n \in \N}$ is a decreasing sequence of subsets of the finite set $\Lambda_i$. Since $c_m$ is $\Delta$-minimal, there is $n \in \N$ and $\sigma \in \Delta_n$ so that for all $\gamma \in \Delta_n$
$$\gamma \sigma^{-1} A b_k^{-1} F_k \cap \dom(c_m) = \gamma \sigma^{-1} ( A b_k^{-1} F_k \cap \dom(c_m))$$
and
$$\forall a \in A b_k^{-1} F_k \cap \dom(c_m) \ c_m(\gamma \sigma^{-1} a) = c_m(a).$$
Notice that $b_k \in F_k$ so $A \subseteq A b_k^{-1} F_k$ and hence $(\sigma \gamma^{-1} \cdot c_m) \res A = c_m \res A$ for all $\gamma \in \Delta_n$.

Fix $\gamma \in \Delta_n$. To finish the proof it suffices to show that
$$\gamma \sigma^{-1} A \cap \dom(c) = \gamma \sigma^{-1} A \cap \dom(c_m).$$
Clearly $\gamma \sigma^{-1} A \cap \dom(c_m) \subseteq \gamma \sigma^{-1} A \cap \dom(c)$ since $c$ extends $c_m$. Pick $a \in A$ with $\gamma \sigma^{-1} a \in \dom(c)$ and towards a contradiction suppose $\gamma \sigma^{-1} a \not\in \dom(c_m)$. Since $\Theta_i(c) = \Theta_i(c_m)$ for all $1 \leq i \leq k$, we have by Definition \ref{DEF FUNDF} that
$$\dom(c) - \dom(c_m) = \bigcup_{i = k+1}^\infty \Delta_i (\Theta_i(c_m) - \Theta_i(c)) b_{i-1} \subseteq \Delta_k b_k.$$
So $\gamma \sigma^{-1} a \in \Delta_k b_k$. Let $\psi \in \Delta_k$ be such that $\gamma \sigma^{-1} a = \psi b_k$. Then $\psi$ passes the $\Delta_k$ membership test and $\psi F_k \subseteq \gamma \sigma^{-1} A b_k^{-1} F_k$. It follows that $ \sigma \gamma^{-1} \psi$ must also pass the $\Delta_k$ membership test and thus $\sigma \gamma^{-1} \psi \in \Delta_k$. Then $a = \sigma \gamma^{-1} \psi b_k \in \Delta_k b_k$ so by our choice of $k$ we have $a \in \bigcap_{i \in \N} \Delta_i b_i$. By Definition \ref{DEF FUNDF} and clause (iii) of Theorem \ref{FM} we have $a \in \dom(c_m)$. But then $\gamma \sigma^{-1} a \in \dom(c_m)$, a contradiction.
\end{proof}

Now we look at subflows generated by $\Delta$-minimal functions.

\begin{lem}\label{MIN REG1}
Let $G$ be a countably infinite group, let $(\Delta_n, F_n)_{n \in \N}$ be a centered, directed, and maximally disjoint blueprint, and let $x \in 2^G$ be fundamental with respect to this blueprint. If $c \subseteq x$ is $\Delta$-minimal and $y \in \overline{[x]}$ is $x$-centered, then
$$\forall g \in \dom(c) \ y(g) = x(g).$$
\end{lem}

\begin{proof}
Let $A$ be an arbitrary finite subset of $G$. We will show $y(a) = x(a)$ for all $a \in A \cap \dom(c)$. Since $c$ is $\Delta$-minimal and the blueprint is centered and directed, there is $n \in \N$ so that for all $\gamma \in \Delta_n$,
$$\gamma A \cap \dom(c) = \gamma (A \cap \dom(c) )$$
and
$$\forall a \in A \cap \dom(c) \ c(\gamma a) = c(a).$$
Since $y \in \overline{[x]}$, there is $g \in G$ so that $y(h) = (g^{-1} \cdot x)(h) = x(g h)$ for all $h \in F_n \cup A$. Since $y$ is $x$-centered, $1_G \in \Delta_n^y$. Consequently, we must have $g \in \Delta_n^x$. It follows that for all $a \in A \cap \dom(c)$
$$y(a) = x(g a) = c(g a) = c(a) = x(a).$$
As $A$ was an arbitrary finite subset of $G$, we conclude $y(g) = x(g)$ for all $g \in \dom(c)$.
\end{proof}

\begin{lem}\label{MIN REG2}
Let $G$ be a countably infinite group, let $(\Delta_n, F_n)_{n \in \N}$ be a centered, directed, and maximally disjoint blueprint, and let $x \in 2^G$ be fundamental with respect to this blueprint. The following are equivalent:
\begin{enumerate}
\item[\rm (i)] $x$ is $\Delta$-minimal;
\item[\rm (ii)] $\{y \in \overline{[x]} \: y \text{ is } x\text{-regular} \} = [x]$.
\end{enumerate}
\end{lem}

\begin{proof}
First assume that $x$ is $\Delta$-minimal. If $z \in \overline{[x]}$ is $x$-regular, then there is an $x$-centered $y \in [z]$. By the previous lemma, $x = y \in [z]$ so $z \in [x]$. On the other hand, it is clear that every element of $[x]$ is $x$-regular.

Now assume that $x$ is not $\Delta$-minimal. Then there is a finite $A \subseteq G$ such that for all $n \in \N$ we can find $\gamma_n \in \Delta_n$ with $x(\gamma_n a) \neq x(a)$ for some $a \in A$. Let $y$ be a limit point of the sequence $(\gamma_n^{-1} \cdot x)_{n \in \N}$. Since $(\Delta_n)_{n \in \N}$ is a decreasing sequence, $1_G \in \gamma_n^{-1} \Delta_k$ for all $n \geq k$. So by clause (i) of Proposition \ref{SUBFLOW BP}, $1_G \in \Delta_k^y$ for all $k \in \N$. Thus $y$ is $x$-centered and in particular $x$-regular. If $y \in [x]$ then there is $z \in [y]$ with $z = x$. Such a $z$ would be have to be $x$-centered, but $y$ is $x$-centered and every orbit in $\overline{[x]}$ contains at most one $x$-centered element (clause (iv) of Lemma \ref{SUBFLOW BP LIST}). So $y \in [x]$ if and only if $y = x$. However, since $A$ is finite there is $a \in A$ and $n \in \N$ with $y(a) = (\gamma_n^{-1} \cdot x)(a) = x(\gamma_n a) \neq x(a)$. Thus $y \neq x$ and $y \not\in [x]$. We conclude that the set of $x$-regular elements of $\overline{[x]}$ does not coincide with $[x]$.
\end{proof}

\section{Minimality constructions} \label{SECT MIN CONSTR}

In showing that every group has a $2$-coloring (Theorem \ref{GEN COL}) and in characterizing groups with the ACP (Theorem \ref{thm:ACP}), we have seen that putting a graph structure on the sets $\Delta_n$ can be very useful. Indeed, the undefined points of a fundamental function $c$ are very useful for controlling which elements of $[c]$ are close, or look similar, and which are far apart. Usually, when one has an application in mind one has a specific idea of which points of $[c]$ should look similar and which ones should not. Since the $\Delta_n$'s organize the positioning of the undefined points of $c$, such requirements are best studied in the form of a graph on $\Delta_n$. We will see that this approach is relied upon heavily in the next section. This section develops methods for extending $\Delta$-minimal functions to new $\Delta$-minimal functions through the analysis of graphs on $\Delta_n$.

The next definition is very much in the spirit of Definition \ref{DEF DMIN}. However, note that in the definition below $\Delta$ appears with a subscript while it does not in Definition \ref{DEF DMIN}.

\begin{definition} \label{DEF GMIN}\index{$\Delta_n$-minimal}\index{minimal!$\Delta_n$-minimal}
Let $G$ be a countably infinite group, let $(\Delta_n, F_n)_{n \in \N}$ be a blueprint, and let $n \geq 1$. A symmetric relation $E \subseteq \Delta_n \times \Delta_n$ is called \emph{$\Delta_n$-minimal} if
\begin{enumerate}
\item[ (i)] there is a finite $A \subseteq G$ such that $(\gamma, \psi) \in E$ implies $\psi \in \gamma A$;
\item[ (ii)] for every $\psi_1, \psi_2 \in \Delta_n$, there is $m \in \N$ and $\sigma \in \Delta_m$ with $\Delta_m \sigma^{-1} \{\psi_1, \psi_2\} \subseteq \Delta_n$ and such that for every $\gamma \in \Delta_m$ $(\gamma \sigma^{-1} \psi_1, \gamma \sigma^{-1} \psi_2) \in E$ if and only if $(\psi_1, \psi_2) \in E$.
\end{enumerate}
\end{definition}

As before, the definition takes a simpler form if we place more assumptions on the blueprint.

\begin{lem} \label{GMIN CENTER DIRECT}
Let $G$ be a countably infinite group, let $(\Delta_n, F_n)_{n \in \N}$ be a centered and directed blueprint, let $n \geq 1$, and let $E \subseteq \Delta_n \times \Delta_n$ be a symmetric relation. Then $E$ is $\Delta_n$-minimal if and only if the following hold:
\begin{enumerate}
\item[\rm (i)] there is a finite $A \subseteq G$ such that $(\gamma, \psi) \in E$ implies $\psi \in \gamma A$;
\item[\rm (ii)] for every $\psi_1, \psi_2 \in \Delta_n$, there is $m \geq n$ with $\psi_1, \psi_2 \in D_n^m$ and for every $\gamma \in \Delta_m$ $(\gamma \psi_1, \gamma \psi_2) \in E$ if and only if $(\psi_1, \psi_2) \in E$.
\end{enumerate}
\end{lem}

\begin{proof}
First suppose that $E$ is $\Delta_n$-minimal. Clearly property (i) above is satisfied. Fix $\psi_1, \psi_2 \in \Delta_n$ and let $m \in \N$ and $\sigma \in \Delta_m$ be such that $\Delta_m \sigma^{-1} \{\psi_1, \psi_2\} \subseteq \Delta_n$ and for every $\gamma \in \Delta_m$, $(\gamma \sigma^{-1} \psi_1, \gamma \sigma^{-1} \psi_2) \in E$ if and only if $(\psi_1, \psi_2) \in E$. By clause (iv) of Lemma \ref{STRONG BP LIST} there is $k \geq n$ with $\psi_1, \psi_2, \sigma \in F_k$ and hence $\psi_1, \psi_2 \in D^k_n$ and $\sigma \in D^k_m$. If $\gamma \in \Delta_k$, then $\gamma \sigma \in \Delta_k D^k_m \subseteq \Delta_m$. Since $\gamma \sigma \sigma^{-1} = \gamma$, we have that $(\gamma \psi_1, \gamma \psi_2) \in E$ if and only if $(\psi_1, \psi_2) \in E$. Thus $E$ has the property above.

Now suppose that $E$ has the property above. The property $\psi_1, \psi_2 \in D^m_n$ implies $\Delta_m \{\psi_1, \psi_2\} \subseteq \Delta_n$. Since the blueprint is centered, we can pick $\sigma = 1_G \in \Delta_m$ to see that $E$ is $\Delta_n$-minimal.
\end{proof}

\begin{lem}\label{SIMP MIN}
Let $G$ be a countably infinite group, and let $(\Delta_n, F_n)_{n \in \N}$ be a centered and directed blueprint. Let $m \geq n \geq 1$, let $\lambda \in D^m_n$, and let $A \subseteq G$ be finite. Define $E \subseteq \Delta_n \times \Delta_n$ by
$$(\psi_1, \psi_2) \in E \Longleftrightarrow (\psi_1 \in \Delta_m \lambda \wedge \psi_2 \in \psi_1 A) \vee (\psi_2 \in \Delta_m \lambda \wedge \psi_1 \in \psi_2 A).$$
Then $E$ is $\Delta_n$-minimal.
\end{lem}

\begin{proof}
We will apply Lemma \ref{GMIN CENTER DIRECT}. Clearly clause (i) is satisfied. We only need to check clause (ii). Fix $\psi_1, \psi_2 \in \Delta_n$. By clause (iv) of Lemma \ref{STRONG BP LIST} there is $k > m$ with $\psi_1, \psi_2 \in F_k$. Fix $\gamma \in \Delta_k$. Clearly,
$$\psi_1^{-1} \psi_2 \in A \Longleftrightarrow (\gamma \psi_1)^{-1} (\gamma \psi_2) \in A \text{ and } \psi_2^{-1} \psi_1 \in A \Longleftrightarrow (\gamma \psi_2)^{-1} (\gamma \psi) \in A,$$
and by conclusion (vii) of Lemma \ref{BP LIST}
$$\psi_i \in \Delta_m \lambda \Longleftrightarrow \gamma \psi_i \in \Delta_m \lambda.$$
Therefore $(\gamma \psi_1, \gamma \psi_2) \in E$ if and only if $(\psi_1, \psi_2) \in E$. We conclude $E$ is $\Delta_n$-minimal.
\end{proof}

\begin{lem}\label{DISJOINT}
Let $G$ be a countably infinite group, let $(\Delta_n, F_n)_{n \in \N}$ be a centered and directed blueprint, let $n \geq 1$, and let $E_1, E_2 \subseteq \Delta_n \times \Delta_n$ be symmetric relations. If $E_1$ is $\Delta_n$-minimal and $E_1 \cap E_2 = \varnothing$, then $E_1 \cup E_2$ is $\Delta_n$-minimal if and only if $E_2$ is $\Delta_n$-minimal.
\end{lem}

\begin{proof}
Since the blueprint is centered, $(\Delta_n)_{n \in \N}$ is decreasing and $(D^m_n)_{m \geq n}$ is increasing ($n$ is fixed). Suppose $E_2$ is $\Delta_n$-minimal. Pick $\psi_1, \psi_2 \in \Delta_n$ and apply Lemma \ref{GMIN CENTER DIRECT} to $E_1$ and $E_2$ to get numbers $m_1$ and $m_2$ in clause (ii). Clearly then taking $m = \max(m_1, m_2)$ shows that $E_1 \cup E_2$ satisfies (ii) for $\psi_1$ and $\psi_2$. We conclude that $E_1 \cup E_2$ is $\Delta_n$-minimal.

Now assume $E_1 \cup E_2$ is $\Delta_n$-minimal. Then $E_2$ satisfies clause (i) of Lemma \ref{GMIN CENTER DIRECT}. So we only need to check clause (ii). For any $\psi_1, \psi_2 \in \Delta_n$ we have
$$(\psi_1, \psi_2) \in E_2 \Longleftrightarrow (\psi_1, \psi_2) \in E_1 \cup E_2 \text{ and } (\psi_1, \psi_2) \not\in E_1.$$
So if $\psi_1, \psi_2 \in \Delta_n$ and $m_1$ and $m_2$ are as in clause (ii) for $E_1$ and $E_1 \cup E_2$, respectively, then $m = \max(m_1, m_2)$ shows that $E_2$ satisfies (ii) for $\psi_1$ and $\psi_2$. We conclude that $E_2$ is $\Delta_n$-minimal.
\end{proof}

\begin{lem}\label{MIN GRAPH}
Let $G$ be a countably infinite group, let $(\Delta_n, F_n)_{n \in \N}$ be a blueprint, let $c \in 2^{\subseteq G}$ be fundamental with respect to this blueprint, and let $n, t \geq 1$ satisfy $|\Theta_n| \geq t$. If $\mu : \Delta_n \rightarrow \{0, 1, \ldots, 2^t - 1\}$ then there is $c' \supseteq c$ such that
\begin{enumerate}
\item[\rm (i)] $c'$ is fundamental with respect to $(\Delta_n, F_n)_{n \in \N}$;
\item[\rm (ii)] $c' \res (\dom(c') - \dom(c))$ is $\Delta$-minimal if $\mu$ is $\Delta$-minimal;
\item[\rm (iii)] $|\Theta_n(c')| = |\Theta_n(c)| - t$;
\item[\rm (iv)] for all $\gamma, \psi \in \Delta_n$
$$\mu(\gamma) = \mu(\psi) \Longleftrightarrow \forall f \in F_n \cap (\dom(c') - \dom(c)) \ c'(\gamma f) = c'(\psi f).$$
\end{enumerate}
\end{lem}

\begin{proof}
For each $i \geq 1$, define $\B_i : \N \rightarrow \{0,1\}$ so that $\B_i(k)$ is the $i^\text{th}$ digit from least to most significant in the binary representation of $k$ when $k \geq 2^{i-1}$ and $\B_i(k) = 0$ when $k < 2^{i-1}$. Fix distinct $\theta_1, \theta_2, \ldots, \theta_t \in \Theta_n$. Define $c' \supseteq c$ by setting
$$c'(\gamma \theta_i b_{n-1}) = \B_i(\mu(\gamma))$$
for each $\gamma \in \Delta_n$ and $1 \leq i \leq t$. Clearly $c'$ satisfies (i) and (iii). Also, since $t$ binary digits are sufficient to encode the values of $\mu$, property (iv) is also satisfied. We proceed to check (ii).

Assume $\mu$ is $\Delta$-minimal. Let $A \subseteq G$ be finite. Since $\mu$ is $\Delta$-minimal, there is $m \in \N$ and $\sigma \in \Delta_m$ such that for all $\gamma \in \Delta_m$
$$\gamma \sigma^{-1} A \Theta_n(c)^{-1} \cap \Delta_n = \gamma \sigma^{-1} (A \Theta_n(c) \cap \Delta_n)$$
and
$$\forall a \in A \Theta_n(c)^{-1} \cap \Delta_n \ \mu(\gamma \sigma^{-1} a) = \mu(a).$$
Fix $\gamma \in \Delta_n$, $a \in A$, and $1 \leq i \leq t$. By the above expressions we have that $\gamma \sigma^{-1} a \in \Delta_n \theta_i$ if and only if $a \in \Delta_n \theta_i$. By varying $i$ between $1$ and $t$ we see that $\gamma \sigma^{-1} a \in \dom(c') - \dom(c)$ if and only if $a \in \dom(c') - \dom(c)$. If $a \in \dom(c') - \dom(c)$, then there are $\psi_1, \psi_2 \in \Delta_n$ and $1 \leq i \leq t$ with $a = \psi_1 \theta_i$ and $\gamma \sigma^{-1} a = \psi_2 \theta_i$. By the $\Delta$-minimality of $\mu$ we have that $\mu(\psi_1) = \mu(\psi_2)$ and therefore
$$c'(\gamma \sigma^{-1} a) = \B_i(\mu(\psi_2)) = \B_i(\mu(\psi_1)) = c'(a).$$
We conclude that $c' \res (\dom(c') - \dom(c))$ is $\Delta$-minimal.
\end{proof}

Let $G$ be a countably infinite group, let $(\Delta_n, F_n)_{n \in \N}$ be a blueprint, let $n \geq 1$, and let $\Gamma$ a graph with vertex set $\Delta_n$. If $F \subseteq E(\Gamma)$ is symmetric, define \index{$[F]$}
$$[F] = \{(\gamma, \psi) \: \exists \sigma_1, \sigma_2, \ldots, \sigma_m \in \Delta_n \ \gamma = \sigma_1 \wedge \psi = \sigma_m \wedge \forall 1 \leq i < m \ (\sigma_i, \sigma_{i+1}) \in F\}.$$
Note that $[F]$ is an equivalence relation on $\Gamma$, even in the case $F$ is empty (the reflexive property follows by taking $m = 1$). For $\gamma \in \Delta_n$, let $[\gamma]$ denote the $[F]$ equivalence class of $\gamma$. We define \index{$\Gamma / [F]$} $\Gamma / [F]$ to be the graph with vertices the equivalence classes of $[F]$ and with edge relation
$$([\gamma], [\psi]) \in E(\Gamma / [F]) \Longleftrightarrow ([\gamma] \neq [\psi]) \wedge (\exists \sigma \in [\gamma], \lambda \in [\psi] \ (\sigma, \lambda) \in E(\Gamma)).$$

The next theorem is quite useful in constructing minimal elements of $2^G$ with various properties. Here we use the term ``partition" loosely and allow the empty set to be a member of a partition.

\begin{theorem}\label{MIN EXT}
Let $G$ be a countably infinite group, let $(\Delta_n, F_n)_{n \in \N}$ be a centered blueprint guided by a growth sequence $(H_n)_{n \in \N}$, let $c \in 2^{\subseteq G}$ be fundamental with respect to this blueprint, and let $n, t \geq 1$ satisfy $|\Theta_n| \geq t$. Let $\Gamma$ be a graph with vertex set $\Delta_n$, and let $\{E_1, E_2\}$ be a partition of $E(\Gamma)$ with $E_1 \cap [E_2] = \varnothing$. If the degree of every vertex in $\Gamma / [E_2]$ is less than $2^t$ then there is $c' \supseteq c$ such that
\begin{enumerate}
\item[\rm (i)] $c'$ is fundamental with respect to $(\Delta_n, F_n)_{n \in \N}$;
\item[\rm (ii)] $c' \res (\dom(c') - \dom(c))$ is $\Delta$-minimal if $E_1$ and $[E_2]$ are $\Delta_n$-minimal;
\item[\rm (iii)] $|\Theta_n(c')| = |\Theta_n(c)| - t$;
\item[\rm (iv)] $(\gamma, \psi) \in E_1$ implies $c'(\gamma f) \neq c'(\psi f)$ for some $f \in F_n \cap (\dom(c') - \dom(c))$;
\item[\rm (v)] $(\gamma, \psi) \in E_2$ implies $c'(\gamma f) = c'(\psi f)$ for all $f \in F_n \cap (\dom(c') - \dom(c))$.
\end{enumerate}
\end{theorem}

\begin{proof}
We remind the reader that our blueprint is centered, directed, and maximally disjoint by clause (i) of Lemma \ref{LEM GUIDED BP} and that $(\Delta_k)_{k \in \N}$ is decreasing by clause (i) of Lemma \ref{STRONG BP LIST}. The plan is to construct $\mu: \Delta_n \rightarrow \{0, 1, \ldots, 2^t-1\}$ which is $\Delta$-minimal if $E_1$ and $[E_2]$ are $\Delta_n$-minimal and which satisfies
$$(\gamma, \psi) \in E_1 \Longrightarrow \mu(\gamma) \neq \mu(\psi), \text{ and}$$
$$(\gamma, \psi) \in E_2 \Longrightarrow \mu(\gamma) = \mu(\psi).$$
After doing this, applying Lemma \ref{MIN GRAPH} will complete the proof.

First suppose $E_1$ or $[E_2]$ is not $\Delta_n$-minimal. Since each vertex of $\Gamma / [E_2]$ has degree less than $2^t$, there is a graph-theoretic $(2^t)$-coloring of $\Gamma / [E_2]$, say $\overline{\mu}: V(\Gamma / [E_2]) \rightarrow \{0, 1, \ldots, 2^t - 1\}$. Define $\mu: V(\Gamma) \rightarrow \{0, 1, \ldots, 2^t - 1\}$ by $\mu(\gamma) = \overline{\mu}([\gamma])$ for $\gamma \in \Delta_n$, where $[\gamma]$ is the $[E_2]$-equivalence class containing $\gamma$. Clearly, if $(\gamma, \psi) \in E_2$ then $[\gamma] = [\psi]$ and $\mu(\gamma) = \mu(\psi)$. On the other hand, if $(\gamma, \psi) \in E_1$ then $[\gamma] \neq [\psi]$ since $E_1 \cap [E_2] = \varnothing$, and hence $([\gamma], [\psi]) \in E(\Gamma / [E_2])$ so $\mu(\gamma) = \overline{\mu}([\gamma]) \neq \overline{\mu}([\psi]) = \mu(\psi)$.

Now suppose $E_1$ and $[E_2]$ are $\Delta_n$-minimal. Then there is a finite $A \subseteq G$ such that $(\gamma, \psi) \in E_1 \cup [E_2]$ implies $\psi \in \gamma A$ (and by symmetry $\gamma \in \psi A$). Let $m(0) \geq n$ be such that $A \subseteq H_{m(0)}$. Let $i \in \N$, $\sigma \in \Delta_{m(0)+3+i}$, and $\gamma \in \sigma F_{m(0)+i} \cap \Delta_n$. If $\psi \in \Delta_n$ and $(\psi, \gamma) \in E_1 \cup [E_2]$ then
$$\psi \in \gamma A \subseteq \sigma F_{m(0)+i} H_{m(0)} \subseteq \sigma H_{m(0)+1+i}$$
by clause (ii) of Definition \ref{DEFN BP GUIDE} and clause (iii) of Definition \ref{DEFN GROWTH}. So $\psi \in \sigma F_{m(0)+3+i}$ by clause (iii) of Lemma \ref{LEM GUIDED BP}. This is an important property one should remember for this proof.

In order to define $\mu$, we construct a sequence of functions $(\mu_k)_{k \geq 1}$ mapping into $\{0, 1, \ldots, 2^t -1\}$, and a sequence $(m(k))_{k \in \N}$ satisfying for each $k \geq 1$:
\begin{enumerate}
\item[\rm (1)] $\mu_{k+1} \supseteq \mu_k$;
\item[\rm (2)] $m(k+1) \geq m(k) + 9$;
\item[\rm (3)] $\dom(\mu_k) = \bigcup_{0 \leq i < k} \Delta_{m(i+1)} D^{m(i)}_n$;
\item[\rm (4)] $([\gamma], [\psi]) \in E(\Gamma / [E_2]) \Rightarrow \mu_k(\gamma) \neq \mu_k(\psi)$ and $(\gamma, \psi) \in [E_2] \Rightarrow \mu_k(\gamma) = \mu_k(\psi)$ whenever $\gamma, \psi \in \dom(\mu_k)$;
\item[\rm (5)] For all $\sigma \in \Delta_{m(k)}$ and all $\gamma \in D^{m(k-1)}_n$, $\mu_k (\gamma) = \mu_k (\sigma \gamma)$;
\end{enumerate}

We begin by constructing $\mu_1$. We can of course fix a labeling of $D^{m(0)}_n = F_{m(0)} \cap \Delta_n$ such that two members are labeled differently if their $[E_2]$ classes are adjacent and are labeled the same if they are in the same $[E_2]$ class. We can choose our labels from the set $\{0, 1, \ldots, 2^t -1\}$ since each vertex of $\Gamma / [E_2]$ has degree less than $2^t$. Since our blueprint is centered and directed, since $(\Delta_k)_{k \in \N}$ is decreasing, since $D^{m(0)}_n$ is finite, and since $E_1$ and $[E_2]$ are $\Delta_n$-minimal, there is $m(1) \geq m(0)+9$ such that for all $\gamma, \psi \in D^{m(0)+9}_n$ and all $\sigma \in \Delta_{m(1)}$,
$$(\gamma, \psi) \in E_1 \Longleftrightarrow (\sigma \gamma, \sigma \psi) \in E_1,$$
and
$$(\gamma, \psi) \in [E_2] \Longleftrightarrow (\sigma \gamma, \sigma \psi) \in [E_2].$$
Since $E_1$ and $[E_2]$ look identical on each $\Delta_{m(1)}$-translate of $D^{m(0)+9}_n$, we can copy this labeling on every $\Delta_{m(1)}$-translate of $D^{m(0)}_n$ to get the function $\mu_1$. Clearly (2), (3), and (5) are then satisfied. Let $\gamma, \psi \in \dom(\mu_1)$, say $\gamma \in \sigma_1 D^{m(0)}_n$ and $\psi \in \sigma_2 D^{m(0)}_n$ with $\sigma_1, \sigma_2 \in \Delta_{m(1)}$. First suppose that $(\gamma, \psi) \in [E_2]$. Then $\psi \in \sigma_1 F_{m(0)+3} \cap \sigma_2 F_{m(0)}$, so we must have $\sigma_1 = \sigma_2$ since $m(1) > m(0)+3$. So $\gamma, \psi \in \sigma_1 D^{m(0)}_n$, so $(\sigma_1^{-1} \gamma, \sigma_1^{-1} \psi) \in [E_2]$ by the definition of $m(1)$. Therefore
$$\mu_1(\gamma) = \mu_1(\sigma_1^{-1} \gamma) = \mu_1(\sigma_1^{-1} \psi) = \mu_1(\psi).$$
Now suppose that $([\gamma], [\psi]) \in E(\Gamma / [E_2])$. Then there are $\gamma' \in [\gamma]$ and $\psi' \in [\psi]$ with $(\gamma', \psi') \in E_1$. Since $m(1) \geq m(0)+9$, as before we find that $\gamma' \in \sigma_1 F_{m(0)+3}$, $\psi' \in \sigma_1 F_{m(0)+6}$, and $\psi \in \sigma_1 F_{m(0)+9}$. So $\sigma_2 = \sigma_1$, $\psi \in \sigma_1 F_{m(0)}$ and $\psi' \in \sigma_1 F_{m(0)+3}$. Again by the definition of $m(1)$ we have that $\sigma_1^{-1} \gamma' \in [\sigma_1^{-1} \gamma]$, $\sigma_1^{-1} \psi' \in [\sigma_1^{-1} \psi]$, and $(\sigma_1^{-1} \gamma', \sigma_1^{-1} \psi') \in E_1$. Therefore $([\sigma_1^{-1} \gamma], [\sigma_1^{-1} \psi]) \in E(\Gamma / [E_2])$ and
$$\mu_1(\gamma) = \mu_1(\sigma_1^{-1} \gamma) \neq \mu_1(\sigma^{-1} \psi) = \mu_1(\psi).$$
Therefore $\mu_1$ satisfies (4).

Now suppose $\mu_k$ has been constructed and $m(k)$ has been defined. Again, since $D^{m(k)}_n = F_{m(k)} \cap \Delta_n$ is finite, there is $m(k+1) \geq m(k) + 9$ such that for all $\gamma, \psi \in D^{m(k)+9}_n$ and all $\sigma \in \Delta_{m(k+1)}$
$$(\gamma, \psi) \in E_1 \Longleftrightarrow (\sigma \gamma, \sigma \psi) \in E_1,$$
and
$$(\gamma, \psi) \in [E_2] \Longleftrightarrow (\sigma \gamma, \sigma \psi) \in [E_2].$$
Let $\sigma \in \Delta_{m(k+1)}$, and let $\gamma \in \dom(\mu_k) \cap F_{m(k)+9}$. Then by (3) there is $0 \leq i < k$ with $\gamma \in \Delta_{m(i+1)} F_{m(i)}$. Let $\eta \in \Delta_{m(i+1)}$ and $\lambda \in D^{m(i)}_n$ be such that $\gamma = \eta \lambda$. By the coherent property of blueprints, $\eta \in F_{m(k)+9}$. So $\eta \in D^{m(k)+9}_{m(i+1)} \subseteq D^{m(k+1)}_{m(i+1)}$ and $\sigma \eta \in \Delta_{m(i+1)}$. We apply (5) twice to get
$$\mu_k(\sigma \gamma) = \mu_k((\sigma \eta) \lambda) = \mu_{i+1}((\sigma \eta) \lambda)$$
$$= \mu_{i+1} (\lambda) = \mu_{i+1} (\eta \lambda) = \mu_k (\eta \lambda) = \mu_k (\gamma).$$
We therefore have four important facts:
$$\forall \gamma, \psi \in D^{m(k)+9}_n \ [(\gamma, \psi) \in E_1 \Longleftrightarrow \forall \sigma \in \Delta_{m(k+1)} \ (\sigma \gamma, \sigma \psi) \in E_1];$$
$$\forall \gamma, \psi \in D^{m(k)+9}_n \ [(\gamma, \psi) \in [E_2] \Longleftrightarrow \forall \sigma \in \Delta_{m(k+1)} \ (\sigma \gamma, \sigma \psi) \in [E_2]];$$
$$\forall \gamma \in D^{m(k)+9}_n \ [\gamma \in \dom(\mu_k) \Longleftrightarrow \forall \sigma \in \Delta_{m(k+1)} \ \sigma \gamma \in \dom(\mu_k)];$$
$$\forall \gamma \in D^{m(k)+9}_n \cap \dom(\mu_k) \ \forall \sigma \in \Delta_{m(k+1)} \ \mu_k(\gamma) = \mu_k(\sigma \gamma).$$
The first two follow from our choice of $m(k+1)$, the third from (3) and conclusion (vii) of Lemma \ref{BP LIST}, and the fourth was just verified. These four statements say that $\mu_k$, $E_1$, and $[E_2]$ each look the same on every $\Delta_{m(k+1)}$-translate of $F_{m(k)+9}$.

We choose an extension, $\tilde{\mu_k}$, of $\mu_k$ to $\dom(\mu_k) \cup D^{m(k)}_n$ with the property that if $([\gamma], [\psi]) \in E(\Gamma / [E_2])$ then $\gamma$ and $\psi$ are labeled differently, and if $(\gamma, \psi) \in [E_2]$ then they are labeled the same. The function $\tilde{\mu_k}$ exists since $\mu_k$ satisfies (4), and $\tilde{\mu_k}$ can be assumed to take values in the set $\{0, 1, \ldots, 2^t -1\}$ since the degree of every vertex of $\Gamma / [E_2]$ is less than $2^t$. We copy $\tilde{\mu_k} \res D^{m(k)}_n$ to all $\Delta_{m(k+1)}$-translates of $D^{m(k)}_n$ and union with $\mu_k$ to get $\mu_{k+1}$. Clearly $\mu_{k+1}$ satisfies (1), (2), (3), and (5). To verify (4), we essentially repeat the argument used to show that $\mu_1$ satisfies (4). Let $\gamma, \psi \in \dom(\mu_{k+1})$ with $(\gamma, \psi) \in [E_2]$ or $([\gamma], [\psi]) \in E(\Gamma / [E_2])$. If $\gamma, \psi \in \dom(\mu_k)$, then there is nothing to show. So we may suppose that $\gamma \not\in \dom(\mu_k)$ and hence $\gamma \in \sigma D^{m(k)}_n$ for some $\sigma \in \Delta_{m(k+1)}$. Arguing as we did for $\mu_1$, we find that $\psi \in \sigma F_{m(k)+9}$. By the definition of $m(k+1)$, we have that $(\sigma^{-1} \gamma, \sigma^{-1} \psi) \in E_1$ if and only if $(\gamma, \psi) \in E_1$. Similarly, by again using the definition of $m(k+1)$ and picking $\gamma' \in [\gamma]$ and $\psi' \in [\psi]$ with $(\gamma', \psi') \in E_1$ if necessary, we have $([\sigma^{-1} \gamma], [\sigma^{-1} \psi]) \in E(\Gamma / [E_2])$ if and only if $([\gamma], [\psi]) \in E(\Gamma / [E_2])$. Now we have
$$\mu_{k+1}(\gamma) = \tilde{\mu_k}(\sigma^{-1} \gamma) \neq \tilde{\mu_k}(\sigma^{-1} \psi) = \mu_{k+1}(\psi)$$
if $([\gamma], [\psi]) \in E(\Gamma / [E_2])$ and
$$\mu_{k+1}(\gamma) = \tilde{\mu_k}(\sigma^{-1} \gamma) = \tilde{\mu_k}(\sigma^{-1} \psi) = \mu_{k+1}(\psi)$$
if $(\gamma, \psi) \in [E_2]$. Thus $\mu_{k+1}$ satisfies (4).

We finish by setting $\mu = \bigcup_{k \geq 1} \mu_k$. By clause (iv) of Lemma \ref{STRONG BP LIST} we see that $\dom(\mu) = \Delta_n$. Property (5) together with clauses (iv) and (v) of Lemma \ref{STRONG BP LIST} imply that $\mu$ is $\Delta$-minimal. Finally, property (4) gives
$$(\gamma, \psi) \in E_1 \Longrightarrow \mu(\gamma) \neq \mu(\psi), \text{ and}$$
$$(\gamma, \psi) \in E_2 \Longrightarrow \mu(\gamma) = \mu(\psi).$$
Applying Lemma \ref{MIN GRAPH} completes the proof.
\end{proof}

We can give a few immediate consequences.

\begin{cor}\label{GEN MINCOL}
Let $G$ be a countably infinite group, let $(\Delta_n, F_n)_{n \in \N}$ be a blueprint, and for each $n \geq 1$ let $B_n$ be finite with $\Delta_n B_n B_n^{-1} = G$. If $c \in 2^{\subseteq G}$ is $\Delta$-minimal and fundamental with respect to this blueprint and $|\Theta_n| \geq \log_2 \ (2 |B_n|^4+1)$ for each $n \geq 1$, then there exists a $\Delta$-minimal and fundamental $c' \supseteq c$ with $|\Theta_n(c')| > |\Theta_n(c)| - \log_2 \ (2 |B_n|^4+1) - 1$ and with the property that any $x \in 2^G$ extending $c'$ is a 2-coloring. In particular, every countable group has a minimal $2$-coloring.
\end{cor}

\begin{proof}
Recall that in the proof of Theorem \ref{GEN COL} we defined the graph $\Gamma_n$ on the vertex set $\Delta_n$ with edge relation given by
$$(\gamma, \psi) \in E(\Gamma_n) \Longleftrightarrow \gamma^{-1} \psi \in B_n B_n^{-1} s_n B_n B_n^{-1} \text{ or } \psi^{-1} \gamma \in B_n B_n^{-1} s_n B_n B_n^{-1}.$$
We seek to apply Theorem \ref{MIN EXT} with respect to the ``partition" $\{E(\Gamma_n), \varnothing\}$ of $E(\Gamma_n)$. To establish our claims we need only check that $E(\Gamma_n)$ is $\Delta_n$-minimal for every $n \geq 1$. However, this follows from Lemma \ref{SIMP MIN}.
\end{proof}

\begin{theorem} \label{thm:minimaldensity}
If $G$ is a countably infinite group, then the collection of minimal 2-colorings is dense in $2^G$. 
\end{theorem}

We omit this proof since the following theorem is stronger.

\begin{theorem} \label{thm:minimalperfectdensity}
Let $G$ be a countably infinite group, $x \in 2^G$, and $\epsilon > 0$. Then there is a perfect set of pairwise orthogonal minimal $2$-colorings in the ball of radius $\epsilon$ about $x$.
\end{theorem}

\begin{proof}
Let $x \in 2^G$ and $\epsilon > 0$. Apply Corollary \ref{GEN DENSE} to get a nontrivial locally recognizable $R: A \rightarrow 2$. Now apply Corollary \ref{GEN SUBEXP FREE} to get a function $c$ which is canonical with respect to a centered blueprint $(\Delta_n, F_n)_{n \in \N}$ guided by a growth sequence, compatible with $R$, and satisfies
$$|\Lambda_n| \geq \log_2(4|F_n|^4+2) = 1 + \log_2 \ (2|F_n|^4+1)$$
for all $n \geq 1$. By clause (i) of Lemma \ref{LEM GUIDED BP} and clause (viii) of Lemma \ref{STRONG BP LIST} we can require the blueprint satisfy
$$\bigcap_{n \in \N} \Delta_n a_n = \bigcap_{n \in \N} \Delta_n b_n = \varnothing.$$
By Proposition \ref{CANONICAL DMIN} $c$ is $\Delta$-minimal. Now apply Corollary \ref{GEN MINCOL} to get a fundamental and $\Delta$-minimal $c'$ for which every extension is a $2$-coloring and with $|\Theta_n(c')| > 0$ for all $n \geq 1$. We apply Proposition \ref{GEN ORTH} to get a perfect set $\{x_\tau \: \tau \in 2^{\N}\}$ of pairwise orthogonal functions extending $c'$. If $n \geq 1$, $\theta \in \Theta_n(c')$, and $i \in \{0,1\}$, then the constant function $d: \Delta_n \theta b_{n-1} \rightarrow \{i\}$ is $\Delta$-minimal by clause (vii) of Lemma \ref{BP LIST}. By the proof of Proposition \ref{GEN ORTH}, we see that each $x_\tau$ can be constructed by taking the union of $c'$ with a sequence of such functions $d$ with disjoint domains. By Lemma \ref{LEM UNION}, the union of $c'$ with finitely many of these $d$'s is $\Delta$-minimal, and so by Lemma \ref{MIN UNION} each $x_\tau$ is $\Delta$-minimal. So each $x_\tau$ is minimal and extends $c$ and $c'$. So each $x_\tau$ is a $2$-coloring and if $\gamma \in \Delta_1$ then $d((\gamma \gamma_1)^{-1} \cdot x_\tau, x) < \epsilon$. So $\{(\gamma \gamma_1)^{-1} \cdot x_\tau \: \tau \in 2^{\N}\}$ is a perfect set of pairwise orthogonal minimal $2$-colorings contained in the ball of radius $\epsilon$ about $x$.
\end{proof}

The next two corollaries draw interest from descriptive set theory. The first relates to subflows of $(2^\N)^G$ (as discussed in Section \ref{SEC 2NG}), and the second is a strengthening of one of the main results of the authors' previous paper \cite{GJS}.

\begin{cor}
Let $G$ be a countably infinite group, $x \in (2^\N)^G$, and $\epsilon > 0$. Then there is a perfect set of pairwise orthogonal minimal hyper aperiodic points in the ball of radius $\epsilon$ about $x$.
\end{cor}

\begin{proof}
Let $C \subseteq G$ and $L \subseteq \N$ be finite sets with the property that for every $y \in (2^\N)^G$
$$\forall c \in C \ \forall n \in L \ y(c)(n) = x(c)(n) \Longrightarrow d(x, y) < \epsilon.$$
Fix a bijection $\phi$ between $2^L$ and $2^{|L|} = \{0, 1, \ldots, 2^{|L|}-1\}$. Set $k = 2^{|L|}$, and define $\tilde{x} \in k^G$ by $\tilde{x}(g) = \phi(x(g) \res L)$. As we have mentioned before, in this paper we work in the most restrictive case of Bernoulli flows of the form $2^G$. However, all of our proofs and results trivially generalize to arbitrary Bernoulli flows $m^G$. So by the previous theorem, there exists a perfect set $\tilde{P}$ consisting of pairwise orthogonal minimal $k$-colorings with the property that $\tilde{y} \res C = \tilde{x} \res C$ for all $\tilde{y} \in \tilde{P}$. Now define $\theta: k^G \rightarrow (2^\N)^G$ by
$$\theta(\tilde{y})(g)(n) = \begin{cases}
\phi^{-1}(\tilde{y}(g))(n) & \text{if } n \in L \\
0 & \text{otherwise.}
\end{cases}$$
Then $\theta$ is a homeomorphic embedding which commutes with the action of $G$. Thus $P = \{\theta(\tilde{y}) \: \tilde{y} \in \tilde{P}\}$ is a perfect set of pairwise orthogonal minimal hyper aperiodic points. Moreover, if $y = \theta(\tilde{y}) \in P$ then $\tilde{y} \res C = \tilde{x} \res C$ and therefore $y(g)(n) = x(g)(n)$ for all $g \in C$ and $n \in L$. Thus $P$ is contained in the ball of radius $\epsilon$ about $x$.
\end{proof}

Recall that a set $A \subseteq F(G)$ is a complete section if and only if $A \cap [x] \neq \varnothing$ for every $x \in F(G)$. We refer the reader to \cite{GJS} for the descriptive set theoretic motivation to the following result.

\begin{cor}
If $G$ is a countably infinite group and $(A_n)_{n \in \N}$ is a decreasing sequence of closed complete sections of $F(G)$, then
$$G \cdot \left( \bigcap_{n \in \N} A_n \right)$$
is dense in $2^G$.
\end{cor}

\begin{proof}
Let $x \in 2^G$ and let $\epsilon > 0$. By the previous theorem, there is a minimal $2$-coloring $y$ with $d(x, y) < \epsilon$. Now $\overline{[y]}$ is a compact set contained in $F(G)$, so there is $z \in \overline{[y]} \cap (\bigcap_{n \in \N} A_n)$. Since $y$ is minimal, $\overline{[z]} = \overline{[y]}$. So $\overline{[z]} \cap B(x; \epsilon) = \overline{[y]} \cap B(x; \epsilon) \neq \varnothing$. Since $B(x; \epsilon)$ is open, it follows that $[z] \cap B(x; \epsilon) \neq \varnothing$.
\end{proof}

\section{Rigidity constructions for topological conjugacy} \label{SEC RIGID}

In this section we develop tools to control when $\overline{[x]}$ and $\overline{[y]}$ are topologically conjugate for functions $x, y \in 2^G$ which are fundamental with respect to a centered blueprint guided by a growth sequence. Section \ref{SECT ISO} is highly dependent on our work in this section. We remind the reader the definition of topologically conjugate.

\begin{definition} \index{topological conjugacy relation} \index{topologically conjugate} \index{topological conjugate} \index{conjugacy}
Let $G$ be a countable group and let $S_1, S_2 \subseteq 2^G$ be subflows. $S_1$ is \emph{topologically conjugate} to $S_2$ or is a \emph{topological conjugate} of $S_2$ if there is a homeomorphism $\phi: S_1 \rightarrow S_2$ satisfying $\phi(g \cdot x) = g \cdot \phi(x)$ for all $x \in S_1$ and $g \in G$. Such a function $\phi$ is called a \emph{conjugacy} between $S_1$ and $S_2$. If $\phi: S_1 \rightarrow S_2$ is a function satisfying $\phi(g \cdot x) = g \cdot \phi(x)$ for all $g \in G$ and $x \in S_1$, then we say that $\phi$ \emph{commutes with the action of $G$}. The property of being topologically conjugate induces an equivalence relation on the set of all subflows of $2^G$. We call this equivalence relation the \emph{topological conjugacy relation}.
\end{definition} 

The complexity of the topological conjugacy relation will be studied in Chapter \ref{CHAP ISO}. Here we will be interested in developing a type of rigidity for the topological conjugacy relation. In other words, we seek a method for constructing a collection of functions with the property that any two of these functions generate topologically conjugate subflows if and only if there is a particularly nice conjugacy between the subflows they generate. In order to exert control over the topological conjugacy class of $\overline{[x]}$, one must be able to exert control over the behavior of elements of $\overline{[x]}$. For fundamental $x \in 2^G$ this is a difficult task. Especially bothersome are the elements of $\overline{[x]}$ which are not $x$-regular. Such elements are difficult to work with and difficult to control. The $x$-regular elements are better understood, and $x$-centered elements are best understood of all, partially in view of Lemma \ref{MIN REG1}. The usefulness of the main result of this section should therefore be appreciated. Starting from a fundamental $c \in 2^{\subseteq G}$ (with a few assumptions), we present a method to extend $c$ to a fundamental $c'$ with the property that for $x, y \in 2^G$ extending $c'$ (and satisfying a few conditions), $\overline{[x]}$ and $\overline{[y]}$ are topologically conjugate if and only if there is a conjugacy from $\overline{[x]}$ onto $\overline{[y]}$ sending $x$ to a $y$-centered element of $\overline{[y]}$. We prove this with two independent theorems.

We begin with two lemmas. In this first lemma, $\forall^\infty$ denotes ``for all but finitely many," $\exists^\infty$ denotes ``there exists infinitely many,'' and $g^A$ denotes $\{a g a^{-1} \: a \in A\}$ for $A \subseteq G$ and $g \in G$.

\begin{lem} \label{CHOICE}
Let $G$ be an infinite group, let $A \subseteq G$ be finite, and let $C \subseteq G$ be infinite. Then
$$\forall^\infty \psi \in C \ \exists^\infty \gamma \in C \ \psi \not\in (\gamma^{-1} \psi)^A \cup (\gamma^{-1})^A \cup \gamma (\gamma^{-1})^A.$$
\end{lem}

\begin{proof}
For each $a \in A$, the map $g \mapsto a g a^{-1}$ is an automorphism of $G$, so since $A$ is finite we have that for every $\psi \in C$ there can only be finitely many $\gamma \in C$ with $\psi \in (\gamma^{-1} \psi)^A \cup (\gamma^{-1})^A$. It will therefore suffice to show that for all but finitely many $\psi \in C$ there are infinitely many $\gamma \in C$ with $\psi \not\in \gamma (\gamma^{-1})^A$.

Let $D = \{\psi \in C \: \exists^\infty \gamma \in C \ \psi \not\in \gamma (\gamma^{-1})^A\}$. We will show $|C - D| \leq n = |A|$. Suppose $\psi_1, \psi_2, \ldots, \psi_{n+1} \in C - D$. Then
$$\forall 1 \leq i \leq n+1 \ \forall^\infty \gamma \in C \ \psi_i \in \gamma (\gamma^{-1})^A$$
$$\Longrightarrow \forall^\infty \gamma \in C \ \psi_1, \psi_2, \ldots, \psi_{n+1} \in \gamma (\gamma^{-1})^A.$$
However, $|\gamma (\gamma^{-1})^A| \leq |A| = n$. So for some $i \neq j$ we have $\psi_i = \psi_j$. We conclude $|C - D| \leq n$.
\end{proof}

\begin{lem} \label{F TEST}
Let $G$ be a countably infinite group, let $(\Delta_n, F_n)_{n \in \N}$ be a blueprint guided by a growth sequence $(H_n)_{n \in \N}$, and let $c \in 2^{\subseteq G}$ be fundamental with respect to this blueprint. If $n \geq 1$ and $\gamma, \psi \in \Delta_{n+2}$ satisfy
$$\gamma^{-1} [ \gamma F_{n+2} \cap \dom(c)] = \psi^{-1} [ \psi F_{n+2} \cap \dom(c)]$$
and
$$\forall f \in \gamma^{-1} [ \gamma F_{n+2} \cap \dom(c)] \ c(\gamma f) = c(\psi f),$$
then
$$\forall h \in \gamma^{-1} [ \gamma H_n \cap \dom(c)] \ c(\gamma h) = c(\psi h).$$
\end{lem}

\begin{proof}
Let $h \in H_n - F_{n+2}$. If $\lambda \in \Delta_1$ and $\gamma h \in \lambda F_1$, then $\gamma h \in \lambda F_1 \subseteq \gamma F_{n+2}$ by clause (iii) of Lemma \ref{LEM GUIDED BP} which contradicts $h \not\in F_{n+2}$. So $\gamma h \not\in \Delta_1 F_1$ and by an identical argument $\psi h \not\in \Delta_1 F_1$ as well. It follows from clause (iv) of Theorem \ref{FM} and the fact that $\bigcup_{n \geq 1} \Delta_n \Lambda_n b_{n-1} \subseteq \Delta_1 F_1$ that $\gamma h, \psi h \in \dom(c)$ and $c(\gamma h) = c(\psi h)$. Therefore $c(\gamma h) = c(\psi h)$ for all
$$h \in (H_n - F_{n+2}) \cup \gamma^{-1} [ \gamma F_{n+2} \cap \dom(c)] = \gamma^{-1} [ \gamma H_n \cap \dom(c)].$$
\end{proof}

We introduce a useful tool from symbolic dynamics known as block codes.

\begin{definition} \index{block code}\index{block code!induced by}
Let $G$ be a countable group. A \emph{block code} is any function $\hat{f} : 2^H \rightarrow \{0,1\}$ where $H$ is a finite subset of $G$. A function $f: A \rightarrow B$, where $A, B \in \Su(G)$, is \emph{induced by a block code $\hat{f}$} if for all $x \in A$ and $g \in G$
$$f(x)(g) = \hat{f}((g^{-1} \cdot x)\res H)$$
where $H = \dom(\hat{f})$.
\end{definition}

The following theorem is usually stated for $G = \ZZ$. The generalization to other groups is immediate and well known. We include a proof for completeness.

\begin{theorem}[\cite{DLBM}, Proposition 1.5.8] \label{THM BC}
Let $G$ be a countable group. A function $f : A \rightarrow B$, where $A$ and $B$ are subflows of $2^G$, is continuous and commutes with the action of $G$ if and only if $f$ is induced by a block code.
\end{theorem}

\begin{proof}
First suppose $f$ is induced by a block code $\hat{f}$. Let $H = \dom(\hat{f})$. If we fix $g \in G$, then for all $x \in A$ we have $f(x)(g) = \hat{f}((g^{-1} \cdot x)\res H)$. Thus if $y \in 2^G$ agrees with $x$ on the finite set $g H$, then $f(y)(g) = f(x)(g)$. So $f$ is continuous. Also, for $x \in A$ and $g, h \in G$
$$f(g \cdot x)(h) = \hat{f}((h^{-1} g \cdot x)\res H) = f(x) (g^{-1} h) = (g \cdot f(x))(h).$$
So $f(g \cdot x) = g \cdot f(x)$.

Now suppose $f$ is continuous and commutes with the shift action of $G$. Since $A$ is a compact metric space, $f$ is uniformly continuous. So there is a finite $H \subseteq G$ such that for all $x, y \in A$
$$\forall h \in H \ x(h) = y(h) \Longrightarrow f(x)(1_G) = f(y)(1_G).$$
So we may define $\hat{f}: 2^H \rightarrow \{0, 1\}$ by
$$\hat{f}(z) = \begin{cases}
f(x)(1_G) & \text{if there is } x \in A \text{ with } x\res H = z \\
0 & \text{otherwise}
\end{cases}$$
Then we have for all $x \in A$ and all $g \in G$
$$f(x)(g) = (g^{-1} \cdot f(x))(1_G) = f(g^{-1} \cdot x)(1_G) = \hat{f}((g^{-1} \cdot x)\res H).$$
So $f$ is induced by the block code $\hat{f}$.
\end{proof}

Now we are ready to present the first rigidity construction.

\begin{theorem} \label{RIGID1}
Let $G$ be a countably infinite group, let $(\Delta_n, F_n)_{n \in \N}$ be a centered blueprint guided by a growth sequence $(H_n)_{n \in \N}$ with $\gamma_n = 1_G$ for all $n \geq 1$, and let $c \in 2^{\subseteq G}$ be fundamental with respect to this blueprint. Suppose that for every $n \geq 1$ there are infinitely many $\gamma \in \Delta_n$ with $c \res {F_n} = (\gamma^{-1} \cdot c) \res {F_n}$ ($c$ being $\Delta$-minimal would be sufficient for this) and that $|\Theta_n| \geq \log_2 (12 |F_n|^2+1)$ for each $n \equiv 1 \mod 10$. Then there are $\nu_1^n, \nu_2^n \in \Delta_{n+5}$ for each $n \equiv 1 \mod 10$ and $c' \supseteq c$ with the following properties:
\begin{enumerate}
\item[\rm (i)] $c'$ is fundamental with respect to $(\Delta_n, F_n)_{n \in \N}$;
\item[\rm (ii)] $c' \res (\dom(c') - \dom(c))$ is $\Delta$-minimal;
\item[\rm (iii)] $|\Theta_n(c')| > |\Theta_n(c)| - \log_2 \ (12 |F_n|^2+1) - 1$ for $n \equiv 1 \mod 10$, and $\Theta_n(c') = \Theta_n(c)$ otherwise;
\item[\rm (iv)] $c'(f) = c'(\nu_1^n f) = c'(\nu_2^n f)$ for all $f \in F_{n+4} \cap \dom(c')$ and all $n \equiv 1 \mod 10$;
\item[\rm (v)] if $x, y \in 2^G$ extend $c'$ and $x(f) = x(\nu_1^n f) = x(\nu_2^n f)$ for all $f \in F_{n+4}$ and all $n \equiv 1 \mod 10$, then any conjugacy between $\overline{[x]}$ and $\overline{[y]}$ must map $x$ to a $y$-regular element of $\overline{[y]}$.
\end{enumerate}
\end{theorem}

Note that $F_{n+4} \subseteq F_{n+5} - \{b_{n+5}\}$ since $\beta_{n+5} \neq \gamma_{n+5} = 1_G$. Therefore in (iv) $\nu_1^n f, \nu_2^n f \in \dom(c')$ since $1_G, \nu_1^n, \nu_2^n \in \Delta_{n+5}$ (see Lemma \ref{FUND DOM}).

\begin{proof}
We will actually prove something slightly more general which would have been too cumbersome to include in the statement of the theorem. The overall approach will be to make use of Lemma \ref{LEM REGT}. We wish to construct a sequence of functions $(c_n)_{n \geq -1}$ and a sequence $(\nu_1^{10n+1}, \nu_2^{10n+1})_{n \geq 1}$ satisfying for each $n \in \N$:
\begin{enumerate}
\item[\rm (1)] $1_G$, $\nu_1^{10 n+1}$, and $\nu_2^{10 n+1}$ are distinct elements of $\Delta_{10 n + 6}$;
\item[\rm (2)] $c_{-1} = c$;
\item[\rm (3)] $c_n \supseteq c_{n-1}$;
\item[\rm (4)] $c_n$ is fundamental with respect to $(\Delta_n, F_n)_{n \in \N}$ and $c_n \res (\dom(c_n) - \dom(c))$ is $\Delta$-minimal;
\item[\rm (5)] $|\Theta_{10n +1} (c_n)| > |\Theta_{10n +1} (c_{n-1})| - \log_2 \ (12 |F_{10n +1}|^2+1) - 1$, and $\Theta_m(c_n) = \Theta_m(c_{n-1})$ for all $m \neq 10n +1$;
\item[\rm (6)] $c_n(f) = c_n(\nu_1^{10n+1} f) = c_n(\nu_2^{10 n+1} f)$ for all $f \in F_{10n+5} \cap \dom(c_n)$;
\item[\rm (7)] If $\gamma \in \Delta_{10n+1}$, $a \in F_{10n+1} F_{10n+1}^{-1}$, $\gamma a^{-1} \nu_1^{10n+1} a, \gamma a^{-1} \nu_2^{10n+1} a \in \Delta_{10n+1}$, and
$$c_n(\gamma f) = c_n(\gamma a^{-1} \nu_1^{10n+1} a f) = c_n(\gamma a^{-1} \nu_2^{10n+1} a f)$$
for all $f \in (F_{10n+1} - \{b_{10n+1}\}) \cap \dom(c_n)$, then $\gamma \in \Delta_{10n + 18} F_{10n + 5}$.
\end{enumerate}

Set $c_{-1} = c$, and suppose $c_{-1}$ through $c_{n-1}$ have been constructed. Here is where we introduce the extra bit of generality not mentioned in the statement of the theorem. If desired, one could extend $c_{n-1}$ to any $c_{n-1}'$ satisfying: $c_{n-1}'$ is fundamental with respect to $(\Delta_n, F_n)_{n \in \N}$; $c_{n-1}' \res (\dom(c_{n-1}') - \dom(c))$ is $\Delta$-minimal; $\Theta_m(c_{n-1}') = \Theta_m(c_{n-1})$ for all $m \leq 10(n-1)+5$; $|\Theta_m(c_{n-1}')| \geq \log_2 \ (12 |F_m|^2+1)$ for all $m > 10(n-1)+5$ congruent to $1$ modulo $10$. We will construct $c_n$ to extend $c_{n-1}'$. To arrive at the exact stated conclusions of the theorem, one must chose $c_{n-1}' = c_{n-1}$ at every stage in the construction. If at some stage one chooses to have $c_{n-1}' \neq c_{n-1}$, then (iii) and (5) will no longer be true, but the remaining properties from (i) through (v) and (1) through (7) will still hold.

Set $m = 10n +1$ and let
$$C = \{\nu \in \Delta_{m+5} \: \nu \neq 1_G \wedge \forall f \in F_{m+4} \cap \dom(c_{n-1}') \ c_{n-1}'(\nu f) = c_{n-1}'(f) \}.$$
Since $c_{n-1}' \res (\dom(c_{n-1}') - \dom(c))$ is $\Delta$-minimal, our assumption on $c$ gives that $C$ is infinite. By applying Lemma \ref{CHOICE}, we see that there are distinct $\nu_1, \nu_2 \in C$ satisfying for every $a \in F_m F_m^{-1}$ and every $\lambda \in D^{m+4}_m$:
$$\nu_2 \neq a \lambda^{-1} \nu_1^{-1} \lambda a^{-1};$$
$$\nu_2 \neq a \lambda^{-1} \nu_1^{-1} \nu_2 \lambda a^{-1};$$
$$\nu_2 \neq \lambda a^{-1} \nu_1^{-1} a \lambda^{-1};$$
$$\nu_2 \neq \nu_1 \lambda a^{-1} \nu_1^{-1} a \lambda^{-1};$$
For notational convenience, let $\nu_0 = 1_G$.

Set $A = \{a^{-1} \nu_i a \: i = 1,2 \wedge a \in F_m F_m^{-1} \}$ and let $k \geq m+17$ be such that $\{1_G, \nu_1, \nu_2 \} F_{m+4} \subseteq F_k$. Define $E_1, E_2 \subseteq \Delta_m \times \Delta_m$ by:
$$E_2 = \{(\sigma \nu_i \lambda, \sigma \nu_j \lambda) \: i \neq j \in \{0, 1, 2\}, \ \sigma \in \Delta_k, \text{ and } \lambda \in D^{m+4}_m \};$$
$$(\gamma, \psi) \in E_1 \Longleftrightarrow [(\gamma, \psi) \not\in E_2] \wedge [\psi \in \gamma A \vee \gamma \in \psi A].$$
Let $\Gamma$ be the graph with vertex set $\Delta_m$ and edge relation $E(\Gamma) = E_1 \cup E_2$. Our immediate goal is to apply Theorem \ref{MIN EXT}.

Since the $\Delta_k$-translates of $F_k$ are disjoint, our choice of $k$ implies $[E_2] = E_2$, and every $[E_2]$ equivalence class consists of three members. So $E_1 \cap [E_2] = \varnothing$. By applying Lemma \ref{SIMP MIN}
for each $i \neq j \in \{0, 1, 2\}$ and each $\lambda \in D^{m+4}_m$ and taking unions, we see that $E_2$ is $\Delta_m$-minimal (for the set $A$ in Lemma \ref{SIMP MIN}, use a singleton of the form $\{\lambda^{-1} \nu_i^{-1} \nu_j \lambda\}$). It is also clear from Lemma \ref{SIMP MIN} that $E_1 \cup E_2$ is $\Delta_m$-minimal as well (in the lemma use $m = n$ and $\lambda = 1_G$). Since $E_1 \cap E_2 = \varnothing$, $E_1$ is $\Delta_m$-minimal by Lemma \ref{DISJOINT}. Since $|A| \leq 2 |F_m|^2$, each vertex of $\Gamma$ has at most $4 |F_m|^2$ $E_1$-neighbors. Therefore each vertex of $\Gamma / [E_2]$ has degree at most $12 |F_m|^2$ (since every $[E_2]$ class consists of three vertices of $\Gamma$). Let $t$ be the least integer greater than or equal to $\log_2 \ (12|F_m|^2+1)$ and apply Theorem \ref{MIN EXT} to get $c_n$ from $c_{n-1}'$. Define $\nu_1^m = \nu_1$ and $\nu_2^m = \nu_2$. Properties (1) through (6), with the exception of (5) if $c_{n-1}' \neq c_{n-1}$, are clearly satisfied. We proceed to verify (7).

Suppose $\gamma \in \Delta_m$ and $a \in F_m F_m^{-1}$ are such that $\gamma a^{-1} \nu_1 a, \gamma a^{-1} \nu_2 a \in \Delta_m$ and for all $f \in (F_m - \{b_m\}) \cap \dom(c_n)$
$$c_n(\gamma f) = c_n(\gamma a^{-1} \nu_1 a f) = c_n(\gamma a^{-1} \nu_2 a f).$$
We cannot have $(\gamma, \gamma a^{-1} \nu_1 a) \in E_1$ nor $(\gamma, \gamma a^{-1} \nu_2 a) \in E_1$. However, we have $\gamma a^{-1} \nu_1 a, \gamma a^{-1} \nu_2 a \in \gamma A$, so it must be that $(\gamma, \gamma a^{-1} \nu_i a) \in E_2$ for $i = 1, 2$. Clearly $\gamma, \gamma a^{-1} \nu_1 a, \gamma a^{-1} \nu_2 a$ are all distinct since $1_G, \nu_1, \nu_2$ are distinct. It follows that there is $\sigma \in \Delta_k$ and $\lambda \in D^{m+4}_m$ such that
$$\{\gamma, \gamma a^{-1} \nu_1 a, \gamma a^{-1} \nu_2 a\} = \{\sigma \lambda, \sigma \nu_1 \lambda, \sigma \nu_2 \lambda\}.$$
If $\gamma = \sigma \lambda$ then we are done since $\sigma \in \Delta_{m+17}$. Towards a contradiction, suppose $\gamma \neq \sigma \lambda$. We have two cases to consider.

\underline{Case 1}: $\gamma = \sigma \nu_1 \lambda$. Observe that
$$\gamma a^{-1} \nu_2 a = \sigma \lambda \Longrightarrow \sigma \nu_1 \lambda a^{-1} \nu_2 a = \sigma \lambda \Longrightarrow \nu_2 = a \lambda^{-1} \nu_1^{-1} \lambda a^{-1},$$
and
$$\gamma a^{-1} \nu_2 a = \sigma \nu_2 \lambda \Longrightarrow \sigma \nu_1 \lambda a^{-1} \nu_2 a = \sigma \nu_2 \lambda \Longrightarrow \nu_2 = a \lambda^{-1} \nu_1^{-1} \nu_2 \lambda a^{-1}.$$
Now by our previous remarks one of the two rightmost statements must be true, but both are in contradiction of our choice of $\nu_1$ and $\nu_2$.

\underline{Case 2}: $\gamma = \sigma \nu_2 \lambda$. We have
$$\gamma a^{-1} \nu_1 a = \sigma \lambda \Longrightarrow \sigma \nu_2 \lambda a^{-1} \nu_1 a = \sigma \lambda \Longrightarrow \nu_2 = \lambda a^{-1} \nu_1^{-1} a \lambda^{-1},$$
and
$$\gamma a^{-1} \nu_1 a = \sigma \nu_1 \lambda \Longrightarrow \sigma \nu_2 \lambda a^{-1} \nu_1 a = \sigma \nu_1 \lambda \Longrightarrow \nu_2 = \nu_1 \lambda a^{-1} \nu_1^{-1} a \lambda^{-1}.$$
Again, one of the two rightmost statements must be true, however both contradict our choice of $\nu_1$ and $\nu_2$. We conclude (7) is satisfied.

Let $c' = \bigcup_{n \in \N} c_n$. Then $c'$ satisfies (i) and (ii). If $c_{n-1}' = c_{n-1}$ for all $n \in \N$ then $c'$ satisfies (iii) as well. For (iv), just note that $1_G, \nu_1^{10n+1}, \nu_2^{10n+1} \in \Delta_{10n+6}$ and $\Delta_{10n+6} F_{10n+5} \cap \dom(c') = \Delta_{10n+6} F_{10n+5} \cap \dom(c_n)$ since $\Delta_{10n+6} F_{10n+5} \cap \Delta_m \Lambda_m b_{m-1} = \varnothing$ for $m \geq 10n+6$ (since $1_G = \gamma_m \not\in \Lambda_m \cup \{\beta_m\}$).

Now let $x, y \in 2^G$ extend $c'$ with $x(f) = x(\nu_1^n f) = x(\nu_2^n f)$ for all $f \in F_{n+4}$ and all $n \equiv 1 \mod 10$. If $\overline{[x]}$ and $\overline{[y]}$ are not topologically conjugate, then there is nothing to show. So assume $\overline{[x]}$ is topologically conjugate to $\overline{[y]}$ and let $\phi: \overline{[x]} \rightarrow \overline{[y]}$ be a conjugacy. By Lemma \ref{LEM REGT} it suffices to show that $F_n F_n^{-1} F_{n+4}^{-1} \cap \Delta_{n+17}^z \neq \varnothing$ for each $n \equiv 1 \mod 10$, where $z = \phi(x) \in \overline{[y]}$.

Since $\phi$ is induced by a block code, there is a finite $K \subseteq G$ such that for all $g, h \in G$
$$\forall k \in K \ x(gk) = x(hk) \Longrightarrow (g^{-1} \cdot x) \res K = (h^{-1} \cdot x) \res K \Longrightarrow z(g) = z(h).$$
Let $n \equiv 1 \mod 10$ satisfy $K \subseteq H_n$. By clause (ii) of Lemma \ref{SUBFLOW BP LIST}, the $\Delta_n^z$-translates of $F_n$ are maximally disjoint within $G$. So there is $a \in F_n F_n^{-1}$ with $a \in \Delta_n^z$. Note that
$$a F_n K \subseteq F_n F_n^{-1} F_n H_n \subseteq H_{n+1} H_n \subseteq H_{n+2},$$
so
$$\forall g, h \in G \ (\forall k \in H_{n+2} \ x(gk) = x(hk) \Longrightarrow \forall f \in F_n \ z(gaf) = z(haf)).$$
Since $1_G, \nu_1^n, \nu_2^n \in \Delta_{n+5}$ and $x(f) = x(\nu_1^n f) = x(\nu_2^n f)$ for all $f \in F_{n+4}$, it follows from Lemma \ref{F TEST} that $x(h) = x(\nu_1^n h) = x(\nu_2^n h)$ for all $h \in H_{n+2}$. Therefore $z(a f) = z(\nu_1^n a f) = z(\nu_2^n a f)$ for all $f \in F_n$.

Let $r \in \N$ be such that
$$\{1_G, \nu_1^n, \nu_2^n\} a F_n \cup a F_{n+4}^{-1} F_{n+17} \subseteq \{g_0, g_1, \ldots, g_r\}$$
(where $g_0, g_1, \ldots$ is the enumeration of $G$ used for defining the metric on $2^G$). Let $p \in G$ be such that $d(p^{-1} \cdot y, z) < 2^{-r}$. Then
$$\forall g \in \{1_G, \nu_1^n, \nu_2^n\} a F_n \cup a F_{n+4}^{-1} F_{n+17} \ y(p g) = (p^{-1} \cdot y)(g) = z(g).$$
As $a \in \Delta_n^z$, it must be that $p a \in \Delta_n^y$. Let $\gamma = p a \in \Delta_n^y$. Then for all $f \in F_n$ and $i = 1, 2$
$$y(\gamma f) = y(p a f) = z(a f) = z(\nu_i^n a f) = y(p \nu_i^n a f) = y(\gamma a^{-1} \nu_i^n a f).$$
So by considering the $\Delta_n^y$ membership test we have that $\gamma a^{-1} \nu_i^n a \in \Delta_n^y$. It follows from (7) that $\gamma \in \Delta_{n+17}^y F_{n+4}$. In particular, there is $s \in F_{n+4}^{-1}$ with $\gamma s \in \Delta_{n+17}^y$. This gives that for all $f \in F_{n+17}$
$$z(a s f) = y(p a s f) = y(\gamma s f),$$
so $a s \in \Delta_{n+17}^z$. In particular, $F_n F_n^{-1} F_{n+4}^{-1} \cap \Delta_{n+17}^z \neq \varnothing$.
\end{proof}

\begin{theorem} \label{RIGID2}
Let $G$ be a countably infinite group, let $(\Delta_n, F_n)_{n \in \N}$ be one of the blueprints referred to in Proposition \ref{SCHEME 1} with $\gamma_n = 1_G$ for all $n \geq 1$, let $(H_n)_{n \in \N}$ be the corresponding growth sequence, and let $c \in 2^{\subseteq G}$ be fundamental with respect to this blueprint. Suppose that $c$ is $\Delta$-minimal and that $|\Theta_n| \geq \log_2 (2 |F_n|^4+1)$ for each $n \equiv 6 \mod 10$. Then for each $n \equiv 6 \mod 10$ there are $\nu_1^n, \nu_2^n, \ldots, \nu_{s(n)}^n \in \Delta_{n+5}$, where $s(n) = |F_n F_n^{-1} - \Z(G)|$, and $c' \supseteq c$ with the following properties:
\begin{enumerate}
\item[\rm (i)] $c'$ is fundamental with respect to $(\Delta_n, F_n)_{n \in \N}$;
\item[\rm (ii)] $c'$ is $\Delta$-minimal;
\item[\rm (iii)] $|\Theta_n(c')| > |\Theta_n(c)| - \log_2 \ (2 |F_n|^4+1) - 1$ for $n \equiv 6 \mod 10$, and $\Theta_n(c') = \Theta_n(c)$ otherwise;
\item[\rm (iv)] $c'(f) = c'(\nu_i^n f)$ for all $1 \leq i \leq s(n)$, $f \in F_{n+4} \cap \dom(c')$, and $n \equiv 6 \mod 10$;
\item[\rm (v)] if $x, y \in 2^G$ extend $c'$ and $x(f) = x(\nu_i^n f)$ for all $1 \leq i \leq s(n)$, $f \in F_{n+4}$, and $n \equiv 6 \mod 10$, then for any $y$-regular $z \in \overline{[y]}$
\begin{quote}
there is a conjugacy between $\overline{[x]}$ and $\overline{[y]}$ sending $x$ into $[z]$ if and only if there is a conjugacy between $\overline{[x]}$ and $\overline{[y]}$ sending $x$ to the unique $y$-centered element of $[z]$.
\end{quote}
\end{enumerate}
\end{theorem}

\begin{proof}
As in the previous theorem, we will actually prove something a little more general. We wish to construct a sequence of functions $(c_n)_{n \geq -1}$ and a collection $\{\nu_i^n \: n \equiv 6 \mod 10, \ 1 \leq i \leq s(n)\}$ satisfying for each $n \in \N$:
\begin{enumerate}
\item[\rm (1)] $1_G, \nu_1^{10n+6}, \nu_2^{10n+6}, \ldots, \nu_{s(10n+6)}^{10n+6}$ are all distinct elements of $\Delta_{10n+11}$;
\item[\rm (2)] $c_{-1} = c$;
\item[\rm (3)] $c_n \supseteq c_{n-1}$
\item[\rm (4)] $c_n$ is fundamental with respect to $(\Delta_n, F_n)_{n \in \N}$ and is $\Delta$-minimal;
\item[\rm (5)] $|\Theta_{10n+6}(c_n)| > |\Theta_{10n+6}(c_{n-1})| - \log_2 \ (2 |F_{10n+6}|^4+1) - 1$, and $\Theta_m(c_n) = \Theta_m(c_{n-1})$ for $m \neq 10n+6$;
\item[\rm (6)] $c_n(f) = c_n(\nu_i^{10n+6} f)$ for all $1 \leq i \leq s(10n+6)$, and $f \in F_{10n+10} \cap \dom(c_n)$;
\item[\rm (7)] for all $a \in F_{10n+6} F_{10n+6}^{-1} - \Z(G)$ there is $1 \leq i \leq s(10n+6)$ and $f \in F_{10n+6} \cap \dom(c_n)$ with $c_n(f) \neq c_n(a^{-1} \nu_i^{10n+6} a f)$.
\end{enumerate}

Set $c_{-1} = c$, and suppose $c_{-1}$ through $c_{n-1}$ have been constructed. As in the previous theorem, if desired, one could extend $c_{n-1}$ to any $c_{n-1}'$ satisfying: $c_{n-1}'$ is fundamental with respect to $(\Delta_n, F_n)_{n \in \N}$; $c_{n-1}'$ is $\Delta$-minimal; $\Theta_m(c_{n-1}') = \Theta_m(c_{n-1})$ for all $m \leq 10n$; $|\Theta_m(c_{n-1}')| \geq \log_2 \ (2 |F_m|^4+1)$ for all $m > 10n$ congruent to $6$ modulo $10$. We will construct $c_n$ to extend $c_{n-1}'$. To arrive at the exact stated conclusions of the theorem, one must chose $c_{n-1}' = c_{n-1}$ at every stage in the construction. If at some stage one chooses to have $c_{n-1}' \neq c_{n-1}$, then (iii) and (5) will no longer be true, but the remaining properties from (i) through (v) and (1) through (7) will still hold.

Set $m = 10n+6$, and enumerate $F_m F_m^{-1} - \Z(G)$ as $a_1, a_2, \ldots, a_{s(m)}$. Using Proposition \ref{SCHEME 1} and the fact that $c_{n-1}'$ is $\Delta$-minimal, we can pick distinct, non-identity $\nu_1, \nu_2, \ldots, \nu_{s(m)} \in \Delta_{m+5}$ one at a time so that they satisfy:
$$\forall 1 \leq i \leq s(m) \ \forall f \in F_{m+4} \cap \dom(c_{n-1}') \ c_{n-1}'(\nu_i f) = c_{n-1}'(f);$$
$$\forall 1 \leq i \leq s(m) \ a_i^{-1} \nu_i a_i \neq \nu_i;$$
$$\forall 1 \leq i \neq j \leq s(m) \ \{\nu_i, a_i^{-1} \nu_i a_i\} \cap \{\nu_j, a_j^{-1} \nu_j a_j\} = \varnothing.$$
For notational convenience, let $\nu_0 = 1_G$.

Let $k$ be such that $\nu_i F_{m+4} \subseteq F_k$ for each $1 \leq i \leq s(m)$. Define $F_1, F_2, E_1, E_2 \subseteq \Delta_m \times \Delta_m$ by:
$$F_1 = \{(\sigma, \sigma a_i^{-1} \nu_i a_i) \: 1 \leq i \leq s(m), \ \sigma \in \Delta_k, \ \sigma a_i^{-1} \nu_i a_i \in \Delta_m\};$$
$$(\gamma, \psi) \in E_1 \Longleftrightarrow [(\gamma, \psi) \in F_1] \vee [(\psi, \gamma) \in F_1];$$
$$F_2 = \{(\sigma \nu_i \lambda, \sigma \nu_j \lambda) \: 0 \leq i \neq j \leq s(m), \ \sigma \in \Delta_k, \ \lambda \in D^{m+4}_m\};$$
$$(\gamma, \psi) \in E_2 \Longleftrightarrow [(\gamma, \psi) \in F_2] \vee [(\psi, \gamma) \in F_2].$$
Let $\Gamma$ be the graph with vertex set $\Delta_m$ and edge relation $E(\Gamma) = E_1 \cup E_2$. Once again, our plan is to apply Theorem \ref{MIN EXT} with respect to $\Gamma$ and the partition $\{E_1, E_2\}$ of $E(\Gamma)$.

Note that, by our choice of $k$, $[E_2] = E_2$ and each $[E_2]$ equivalence class has exactly $s(m)+1$ members. Also, note that $E_1 \cap [E_2] = \varnothing$ (use the fact that if $(\gamma, \psi) \in E_1$ then either $\gamma \in \Delta_k$ or $\psi \in \Delta_k$). By fixing $\lambda \in D^{m+4}_m$ and $0 \leq i \neq j \leq s(m)$, applying Lemma \ref{SIMP MIN}, and then taking unions over all $\lambda$, $i$, and $j$, we see that $E_1$ and $[E_2]$ are both $\Delta_m$-minimal. Each vertex in $\Gamma$ has at most $s(m)$ many $E_1$-neighbors, so each vertex of $\Gamma / [E_2]$ has degree at most $s(m)(s(m)+1) \leq 2s(m)^2 \leq 2 |F_m|^4$. Let $t$ be the least integer greater than or equal to $\log_2 \ (2 |F_m|^4+1)$, and apply Theorem \ref{MIN EXT} to get $c_n$ from $c_{n-1}'$. Set $\nu_i^m = \nu_i$ for each $1 \leq i \leq s(m)$. Properties (1) through (6), with the exception of (5) if $c_{n-1}' \neq c_{n-1}$, are then clearly satisfied. To verify (7), just notice that $1_G \in \Delta_k$ and either $a_i^{-1} \nu_i a_i \not\in \Delta_m$, in which case it fails the $\Delta_m$ membership test while $1_G$ passes, or $(1_G, a_i^{-1} \nu_i a_i) \in E_1$, in which case the claim follows by the definition of $c_n$.

Let $c' = \bigcup_{n \in \N} c_n$. Properties (i) and (ii) clearly hold (for (ii) use Lemma \ref{MIN UNION}). If $c_{n-1}' = c_{n-1}$ for each $n \in \N$, then (iii) holds as well. To see property (iv), just note that $1_G, \nu_i^n \in \Delta_{10n+11}$, and $\Delta_{10n+11} F_{10n+10} \cap \dom(c') = \Delta_{10n+11} F_{10n+10} \cap \dom(c_n)$ since $\Delta_{10n+11} F_{10n+10} \cap \Delta_m \Lambda_m b_{m-1} = \varnothing$ for $m \geq 10n+11$ (since $1_G = \gamma_m \not\in \Lambda_m \cup \{\beta_m\}$).

Now let $x, y \in 2^G$ extend $c'$ with $x(f) = x(\nu_i^n f)$ for all $n \equiv 6 \mod 10$, $1 \leq i \leq s(n)$, and $f \in F_{n+4}$. Let $z \in \overline{[y]}$ be $y$-regular. If there is no conjugacy between $\overline{[x]}$ and $\overline{[y]}$ mapping $x$ into $[z]$ then there is nothing to show. So suppose $\phi: \overline{[x]} \rightarrow \overline{[y]}$ is a conjugacy with $\phi(x) \in [z]$. Without loss of generality, we may assume that $z$ itself is $y$-centered.

Let $a \in G$ be such that $\phi(x) = a \cdot z$. Now if $a \in \Z(G)$, then there is an automorphism (self conjugacy) of $\overline{[z]} = \phi(\overline{[x]}) = \overline{[y]}$ sending $a \cdot z$ to $z$ (the automorphism being $w \mapsto a^{-1} \cdot w$ for $w \in \overline{[y]}$; this is clearly continuous and it commutes with the action of $G$ since $a \in \Z(G)$). Since $z$ is $y$-centered, our claim would follow by composing this automorphism with $\phi$. Towards a contradiction, suppose $a \not\in \Z(G)$. Since $\phi$ is induced by a block code, there is a finite $K \subseteq G$ such that for all $g, h \in G$
$$\forall k \in K \ x(gk) = x(hk) \Longrightarrow (g^{-1} \cdot x) \res K = (h^{-1} \cdot x) \res K \Longrightarrow (a \cdot z)(g) = (a \cdot z)(h).$$
Let $n \equiv 6 \mod 10$ be such that $a \in H_{n-2} \subseteq F_n F_n^{-1}$ and $K \subseteq H_n$. Note that
$$a F_n K \subseteq F_n F_n^{-1} F_n H_n \subseteq H_{n+1} H_n \subseteq H_{n+2}.$$
Therefore
$$\forall g, h \in G \ (\forall k \in H_{n+2} \ x(gk) = x(hk) \Longrightarrow \forall f \in F_n \ (a \cdot z)(g a f) = (a \cdot z)(h a f)).$$
Since $1_G, \nu_i^n \in \Delta_{n+5}$ and $x(f) = x(\nu_i^n f)$ for all $f \in F_{n+4}$ and all $1 \leq i \leq s(n)$, it follows from Lemma \ref{F TEST} that $x(h) = x(\nu_i^n h)$ for all $h \in H_{n+2}$ and all $1 \leq i \leq s(n)$. Therefore, for all $f \in F_n$ and all $1 \leq i \leq s(n)$.
$$z(f) = (a \cdot z) (a f) = (a \cdot z) (\nu_i^n a f) = z(a^{-1} \nu_i^n a f).$$
However, $y \supseteq c'$ and $c'$ is $\Delta$-minimal, so this is in contradiction with property (7) and Lemma \ref{MIN REG1}. Thus it must be that $a \in \Z(G)$.
\end{proof}

These two proofs taken together, give us the following.

\begin{cor} \label{RIGID}
Let $G$ be a countably infinite group, let $(\Delta_n, F_n)_{n \in \N}$ be one of the blueprints referred to in Proposition \ref{SCHEME 1} with $\gamma_n = 1_G$ for all $n \geq 1$, let $(H_n)_{n \in \N}$ be the corresponding growth sequence, and let $c \in 2^{\subseteq G}$ be fundamental with respect to this blueprint. Suppose that $c$ is $\Delta$-minimal and that $|\Theta_n| \geq \log_2 (12 |F_n|^4+1)$ for each $n \equiv 1 \mod 5$. Then for every $n \in \N$ there are $\nu_1^{10n+1}, \nu_2^{10n+1} \in \Delta_{10n+6}$ and $\nu_1^{10n+6}, \nu_2^{10n+6}, \ldots, \nu_{s(10n+6)}^{10n+6} \in \Delta_{10n+11}$ (where $s(n)$ is defined as in the previous theorem) and $c' \supseteq c$ which simultaneously satisfy both the conclusions of Theorem \ref{RIGID1} and the conclusions of Theorem \ref{RIGID2}.
\end{cor}

\begin{proof}
The added generality of the two previous theorems allows the inductive constructions appearing in the proofs to be interwoven.
\end{proof}

We point out that our results in this section are not only applicable in the context of conjugacies, but also in the more general context of continuous functions commuting with the action of $G$. For example, if $x$ and $y$ are as in clause (v) of Theorem \ref{RIGID1} then any continuouos $\phi: \overline{[x]} \rightarrow \overline{[y]}$ which commutes with the action of $G$ must map $x$ to a $y$-regular element of $\overline{[y]}$. Similarly, if $x$ and $y$ are as in clause (v) of Theorem \ref{RIGID2}, then any $\phi$ as in the previous sentence must satisfy $\phi(x) \in \Z(G) \cdot z$ where $z$ is some $y$-centered element of $\overline{[y]}$. These statements can be easily verified by looking back at the proofs of Theorems \ref{RIGID1} and \ref{RIGID2}. We will not make use of these more general facts. 
\chapter{The Descriptive Complexity of Sets of $2$-Colorings} \label{CHAP DESC CMPLXTY}

In the following three chapters we study some further problems involving $2$-colorings on countably infinite groups. In doing this we make use of the fundamental method and its variations as well as employing additional results and methods in  descriptive set theory, combinatorial group theory, and topological dynamics.

In this chapter we study the descriptive complexity of sets of $2$-colorings, minimal elements,
 and minimal $2$-colorings on any countably infinite group $G$.
 It is obvious that all these sets are
 ${\bf\Pi}^0_3$ (i.e. $F_{\sigma\delta}$) subsets of $2^G$. We characterize the exact descriptive complexity of these sets.

\section{Smallness in measure and category} \label{SECT SMALLNESS}

We have shown in the preceding chapter that for any countably infinite group $G$, the set of minimal $2$-colorings on $G$ is dense (Theorem~\ref{thm:minimaldensity}). In addition, within any given open set in $2^G$ there are perfectly many minimal
$2$-colorings on $G$ (Theorem~\ref{thm:minimalperfectdensity}). These manifest the largeness of the sets in some sense. In this section we show
that in the sense of measure or category, both the sets of $2$-colorings and of minimal elements are
small.

The space $2^G$ carries the Bernoulli product measure $\mu$. For any $p\in 2^{<G}$ (with $\dom(p)$ finite),
$$ \mu(N_p)=\mu(\{x\in 2^G\,:\, p\subseteq x\})=2^{-|\dom(p)|}. $$

\begin{lem}\label{lem:smallblocking}
Let $G$ be a countably infinite group and $s\in G$ with $s\neq 1_G$.
Then the set of all elements $x\in 2^G$ blocking $s$ is meager and
has $\mu$ measure $0$ in $2^G$.
\end{lem}

\begin{proof}
For any finite $T\subseteq G$, let
$$ B_T=\{x\in 2^G\,:\, \forall g\in G\ \exists t\in T\ x(gt)\neq x(gst)\}. $$
It is clear that $B_T$ is closed. We first show that each $B_T$ is nowhere dense.

Otherwise, there is a nonempty open set in which $B_T$ is dense. It follows that there
is $p\in 2^{<G}$ such that $N_p\subseteq B_T$. Define $x\in 2^G$ by
$$ x(g)=\left\{\begin{array}{ll} p(g), & \mbox{ if $g\in\dom(p)$,} \\ 0, & \mbox{ otherwise.} \end{array}\right.
$$
Then $x\in N_p$ and so $x\in B_T$. Let $F=\dom(p)T^{-1}\cup \dom(p)T^{-1}s^{-1}$. Since $F$ is finite, there is $g_0\not\in F$. Then
for all $t\in T$, $g_0t\not\in \dom(p)$ and $g_0st\not\in\dom(p)$. This implies that for all $t\in T$,
$x(g_0t)=0=x(g_0st)$, and so $x\not\in B_T$, a contradiction.

The above argument shows that $B_T$ is meager, and in particular it cannot have a nonempty interior. To see that $B_T$ has $\mu$ measure $0$, we use the fact that $\mu$ is actually ergodic, i.e., any invariant Borel subset of $2^G$ must have $\mu$ measure $0$ or $1$. It is clear that $B_T$ is invariant. Toward a contradiction, assume $\mu(B_T)\neq 0$. Then it has $\mu$ measure $1$, and its complement, being open and of $\mu$ measure $0$, must be empty. This shows that $B_T=2^G$, a contradiction.
\end{proof}

\begin{lem}\label{lem:smallminimal}
Let $G$ be a countably infinite group and $A\subseteq G$ finite and nonempty.
Then the set of all elements $x\in 2^G$ with a finite $T\subseteq G$ such that
$$ \forall g\in G\ \exists t\in T\ \forall a\in A\ x(gta)=x(a) $$
is meager and has $\mu$ measure $0$.
\end{lem}

\begin{proof}
For any finite $T\subseteq G$, let
$$ R_T=\{ x\in 2^G\,:\, \forall g\in G\ \exists t\in T\ \forall a\in A\ x(gta)=x(a)\}. $$
Similar to the preceding proof, it suffices to show that each $R_T$ is meager and has $\mu$ measure $0$. Assume not. Since $R_T$ is closed, it has a nonempty interior. Let $p\in 2^{<G}$ be such that $N_p\subseteq R_T$. We consider two cases. Case 1: $A\subseteq \dom(p)$. Fix any $a_0\in A$.
Define $y\in 2^G$ by
$$ y(g)=\left\{\begin{array}{ll} p(g), & \mbox{ if $g\in\dom(p)$,} \\
1-p(a_0), & \mbox{ otherwise.}
\end{array}\right.
$$
Then $y\in N_p$ and so $y\in R_T$. Let $F=\dom(p)A^{-1}T^{-1}$. Since $F$ is finite, there is $g_0\not\in F$. Then for any $t\in T$ and $a\in A$, $g_0ta\not\in \dom(p)$ and $y(g_0ta)=1-p(a_0)$. In particular for any $t\in T$,
$y(g_0ta_0)\neq y(a_0)$. This shows that $y\not\in R_T$, a contradiction. Case 2: $A\not\subseteq \dom(p)$. In this case let $b_0\in A-\dom(p)$. Define $z\in 2^G$ by
$$ z(g)=\left\{\begin{array}{ll} p(g), & \mbox{ if $g\in\dom(p)$,} \\
1, & \mbox{ if $g=b_0$,} \\ 0, & \mbox{ otherwise.}
\end{array}\right.
$$
Then $z\in N_p$ and so $z\in R_T$. Let $K=(\dom(p)\cup\{b_0\})A^{-1}T^{-1}$. Since $K$ is finite, there
is $h_0\not\in K$. Then for any $t\in T$ and $a\in A$, $h_0ta\not\in \dom(p)\cup\{b_0\}$ and so
$z(h_0ta)=0$. In particular for any $t\in T$, $z(h_0tb_0)\neq z(b_0)$. This shows again that $z\not\in R_T$, a contradiction.
\end{proof}

\begin{theorem}\label{thm:coloringmeager} For any countably infinite group $G$ the set of all $2$-colorings on $G$ and the set of all minimal elements are meager and have $\mu$ measure $0$.
\end{theorem}

\begin{proof}
This follows immediately from the above lemmas.
\end{proof}

\section{${\bf\Sigma}^0_2$-hardness and ${\bf\Pi}^0_3$-completeness} \label{SECT CMPLXTY 1}

In this section we show that for any countably infinite group $G$, the set of all $2$-colorings on $G$
is ${\bf\Sigma}^0_2$-hard and the set of all minimal elements in $2^G$ is ${\bf\Pi}^0_3$-complete. This completely characterizes the descriptive complexity for the set of all minimal elements. For the set of
all $2$-colorings the complete characterization for its descriptive complexity will be given in the next two sections.

We first briefly review the related descriptive set theory. For further results and unexplained terminology the reader can consult \cite{KechrisBook} (especially \cite{KechrisBook} Sections 22 and 23).

Let $X$ be an uncountable Polish space, i.e., separable completely metrizable topological space. A subset $Y\subseteq X$ is {\it ${\bf\Sigma}^0_2$-hard} if for any ${\bf\Sigma}^0_2$ (i.e. $F_\sigma$) subset $Z$ of the Baire space $\N^\N$ there is a continuous function $f: \N^\N\to X$
such that $Z=f^{-1}[Y]$, i.e., for all $z\in \N^\N$,
$$ z\in Z\iff f(z)\in Y. $$
A set is {\it ${\bf\Sigma}^0_2$-complete} if it is ${\bf\Sigma}^0_2$ and is also ${\bf\Sigma}^0_2$-hard. Intuitively, a set is ${\bf\Sigma}^0_2$-complete if it is the most complex ${\bf\Sigma}^0_2$ set in a Polish space, and a set is ${\bf\Sigma}^0_2$-hard if it is at least as complex as any ${\bf\Sigma}^0_2$ set. Define
$$ S=\{\alpha\in 2^\N\,:\, \exists n\ \forall m>n\ \alpha(m)=0\}. $$
Then $S$ is known to be ${\bf\Sigma}^0_2$-complete. For any subset $Y$ of an uncountable Polish space $X$, $Y$ is ${\bf\Sigma}^0_2$-hard iff there is a continuous $f: 2^\N\to X$ such that $S=f^{-1}[Y]$. \index{${\bf\Sigma}^0_2$-hard}\index{${\bf\Sigma}^0_2$-complete}

Similar definitions can be given for ${\bf\Pi}^0_3$-hardness and ${\bf\Pi}^0_3$-completeness. Define
$$ P=\{ \alpha\in 2^{\N\times\N}\,:\, \forall k\geq 1\ \exists n>k\ \forall m\geq n\ \alpha(k,m)=0\}. $$
Then $P$ is known to be ${\bf\Pi}^0_3$-complete. Using this fact we can give the following alternative definitions for ${\bf\Pi}^0_3$-hardness and ${\bf\Pi}^0_3$-completeness. For any subset $Y$ of an uncountable Polish space
$X$, $Y$ is {\it ${\bf\Pi}^0_3$-hard} iff there is a continuous function $f: 2^{\N\times\N}\to X$
such that $P=f^{-1}[Y]$; $Y$ is {\it ${\bf\Pi}^0_3$-complete} if $Y$ is ${\bf\Pi}^0_3$ and ${\bf\Pi}^0_3$-hard.\index{${\bf\Pi}^0_3$-hard}\index{${\bf\Pi}^0_3$-complete}

By the definition of $2$-coloring it is obvious that for any countable group $G$ the set of all $2$-colorings on $G$ is ${\bf\Pi}^0_3$. When $G$ is finite then there are only finitely many orbits in
$2^G$, and every orbit is closed. In this case the set of all $2$-colorings on $G$ coincides with the
free part and is also closed. Below we show that for any countably infinite group $G$, the set of all $2$-colorings is ${\bf\Sigma}^0_2$-hard.

\begin{theorem}\label{thm:Sigma02hard} For any countably infinite group $G$, the set of all $2$-colorings on $G$ is ${\bf\Sigma}^0_2$-hard.
\end{theorem}

\begin{proof} We give two short proofs for this theorem. The first proof uses Wadge determinacy and shows a general claim in descriptive set theory: any dense meager set is ${\bf\Sigma}^0_2$-hard. Let $X$ be a Polish space and $C\subseteq X$ dense and meager. If $C$ is not ${\bf\Sigma}^0_2$-hard then for some ${\bf\Sigma}^0_2$ set $Y\subseteq \N^\N$, $Y\not\leq_W C$. Then by Wadge determinacy (c.f. \cite{KechrisBook} Theorem 21.14) $C\leq_W (\N^\N-Y)$.
This shows that $C$ is ${\bf\Pi}^0_2$, or $G_\delta$. Thus $C$ is a dense $G_\delta$ in $2^G$, and therefore comeager, a contradiction.

For readers not familiar with descriptive set theory we offer the following direct proof. Fix a strong $2$-coloring $x \in 2^G$ (given by Theorem~\ref{thm:strongcoloring}). Fix an increasing sequence $(A_n)_{n \in \N}$ of finite subsets of $G$ with $\bigcup_n A_n = G$, and a sequence $(h_n)_{n \in \N}$ with $h_m A_m \cap h_k A_k = \emptyset$ for $m \neq k$. For $\alpha \in 2^\N$, define $f(\alpha) \in 2^G$ by
$$f(\alpha)(g) =\left\{\begin{array}{ll}
x(g), & \mbox{if $g \not\in \bigcup_k h_k A_k$;} \\
x(g), & \mbox{if for some $k$, $g \in h_k A_k$ and $\alpha(k) = 0$;} \\
1, & \mbox{if for some $k$, $g \in h_k A_k$ and $\alpha(k) = 1$.}
\end{array}\right.
$$
Clearly $f$ is continuous. If $\alpha \in S$, then $f(\alpha)$ and $x$ differ at finitely many coordinates. Since $x$ is a strong $2$-coloring, this implies $f(\alpha)$ is a $2$-coloring. On the other hand, if $\alpha \not\in S$ then the set $N=\{k\in\N\,:\, \alpha(k)=1\}$ is infinite, and so $\lim_{k\in N}h_k^{-1}\cdot f(\alpha)$ is the constant 1 function in $2^G$; hence $f(\alpha)$ is not a $2$-coloring.
\end{proof}

The following theorem completely characterize the descriptive complexity for the set of all minimal elements in $2^G$ for any countably infinite group $G$.

\begin{theorem}\label{thm:minimalcomplexity} For any countably infinite group $G$, the set of all minimal elements of $2^G$ is
${\bf\Pi}^0_3$-complete.
\end{theorem}

\begin{proof}
Let $G$ be a countably infinite group and $(\Delta_n, F_n)_{n\in\mathbb{N}}$ a centered blueprint guided by a growth sequence $(H_n)_{n\in\mathbb{N}}$, with $\alpha_n\neq \gamma_n=1_G\neq \beta_n$ for all $n\geq 1$. Then $(\Delta_n,F_n)_{n\in\mathbb{N}}$ is in fact directed and maximally disjoint by clause (i) of Lemma~\ref{LEM GUIDED BP}, and $\bigcap_{n\in\mathbb{N}}\Delta_na_n=\bigcap_{n\in\mathbb{N}}\Delta_nb_n=\emptyset$ by clause (viii) of Lemma~\ref{STRONG BP LIST}. For such blueprints Corollary~\ref{STRONG PREMIN} applies.
If $x\in 2^G$ is fundamental with respect to this blueprint, then $x$ is minimal iff $x$ is pre-minimal iff for every $k\geq 1$ there is $n>k$ such that
$$ \forall \gamma\in \Delta_n\ \exists \lambda\in D^n_k\ \forall f\in F_k\ x(\gamma\lambda f)=x(f). $$
Let $c\in 2^{\subseteq G}$ be $\Delta$-minimal and fundamental with respect to this blueprint, with $\Theta_n=\Theta_n(c)$ nonempty for all $n\geq 1$. By Proposition~\ref{CANONICAL DMIN} it suffices to take $c$ to be canonical with respect to $(\Delta_n,F_n)_{n\in\mathbb{N}}$. Recall
from Definition~\ref{DEF FUNDF} that
$$ G-\dom(c)=\bigcup_{n\geq 1} \Delta_n\Theta_nb_{n-1} $$
and for distinct $n, m\geq 1$, $\Delta_n\Theta_nb_{n-1}$ and $\Delta_m\Theta_mb_{m-1}$ are disjoint (clause (iii) of Theorem~\ref{FM}).
We define a continuous function $\Phi:2^{\N\times\N}\to 2^G$ so that
for all $\alpha\in 2^{\N\times\N}$, $\alpha\in P$ iff $\Phi(\alpha)$ is a minimal element in $2^G$.

For $k \geq 1$, and $\gamma \in \Delta_k$ define
$$R_k(\gamma) = \max \{ n \geq k \: n = k \mbox{ or else } \gamma \in (\Delta_n - \{1_G\}) F_n \}$$
and
$$ S_k(\gamma) = \mbox{the unique $\lambda \in \Delta_{R_k(\gamma)}$ with $\gamma \in \lambda F_{R_k(\gamma)}$}. $$
Note the following basic properties of these functions. If $n>k$ and $\gamma\in F_n$ then $R_k(\gamma)<n$. If
$R_k(\gamma)=k$ then $S_k(\gamma)=\gamma$. If $R_k(\gamma)>k$ then $R_k(S_k(\gamma)^{-1}\gamma)<R_k(\gamma)$ since $S_k(\gamma)^{-1}\gamma\in F_{R_k(\gamma)}$. Intuitively the function $R_k$ is a rank function for elements of $\Delta_k$. These properties make it possible to use the following kind of induction on $\gamma$. The base case of the induction is $\gamma=1_G$. In general, the case for $\gamma$ makes use of the inductive case for $S_k(\gamma)^{-1}\gamma$.

Given $\alpha \in 2^{\N \times \N}$, we define $\Phi(\alpha) \in 2^G$ to extend $c$ as follows. For $k \geq 1$, $\gamma \in \Delta_k$, and $\theta \in \Theta_k$, we inductively define $\Phi(\alpha)(1_G \theta b_{k-1}) = 0$ and
$$\Phi(\alpha)(\gamma \theta b_{k-1}) = \max\{\alpha(k, R_k(\gamma)), \Phi(\alpha)(S_k(\gamma)^{-1}\gamma \theta b_{k-1})\}.$$
Then $\Phi$ is continuous.

Suppose $\alpha \in P$. We will apply Lemma~\ref{LEM FDMIN} to verify that $f(\alpha)$ is indeed $\Delta$-minimal. First note that our blueprint satisfies all the requirements of Lemma~\ref{LEM FDMIN}. To
check $\Delta$-minimality, fix $k \geq 1$. Since $c$ is $\Delta$-minimal, there is $K> k$ so that for all $\gamma \in \Delta_K$ we have $(\gamma^{-1} \cdot c) \res (F_k \cap\dom(c)) = c \res (F_k\cap\dom(c))$. Let $n > K$ be such that $\alpha(t, m) = 0$ for all $t \leq k$ and $m \geq n$. It suffices to verify that for all $\gamma\in \Delta_n$  and $f\in F_k-\dom(c)$, $\Phi(\alpha)(\gamma f)=\Phi(\alpha)(f)$. Consider a fixed $f \in F_k - \dom(c)$. Since $f\not\in\dom(c)$, there is $t \geq 1$
with $f \in \Delta_t \Theta_t b_{t-1}$. For $t > k$, we have $1_G = \gamma_t \not\in \Theta_t$ and therefore $\Delta_t F_{t-1}$ and $\Delta_t \Theta_t F_{t-1}$ are disjoint. So for $t > k$, $F_k \subseteq \Delta_t F_{t-1}$ is disjoint from $\Delta_t \Theta_t b_{t-1} \subseteq \Delta_t \Theta_t F_{t-1}$. Thus we must have $f \in \Delta_t \Theta_t b_{t-1}$ for some $t \leq k$. Thus there are $\lambda \in D^k_t$ and $\theta\in \Theta_t$ such that $f = \lambda \theta b_{t-1}$. Now we need to verify that for all $\gamma\in\Delta_n$,
$$ \Phi(\alpha)(\gamma\lambda\theta b_{t-1})=\Phi(\alpha)(\lambda\theta b_{t-1}). $$
We do this by induction on $R_t(\gamma\lambda)$. For the base case of the induction $\gamma=1_G$, and the identity holds trivially. For the general inductive case, since $\gamma\in \Delta_n-\{1_G\}$, we have $R_t(\gamma\lambda)\geq n$. Thus $\alpha(t, R_t(\gamma\lambda))=0$, and
$$\Phi(\alpha)(\gamma \lambda\theta b_{t-1}) = \Phi(\alpha)(S_t(\gamma\lambda)^{-1} \gamma \lambda\theta b_{t-1}).$$
Since $R_t(S_t(\gamma\lambda)^{-1}\gamma\lambda)<R_t(\gamma\lambda)$, we have by the inductive hypothesis
that
$$\Phi(\alpha)(S_t(\gamma\lambda)^{-1}\gamma\lambda\theta b_{t-1})=\Phi(\alpha)(\lambda\theta b_{t-1}). $$
This shows that $\Phi(\alpha)(\gamma\lambda \theta b_{t-1})=\Phi(\alpha)(\lambda\theta b_{t-1})$ as needed, and so $\Phi(\alpha)$ is $\Delta$-minimal by Lemma~\ref{LEM FDMIN}.

For the converse, suppose $\alpha \not\in P$. We will apply Corollary~\ref{STRONG PREMIN} to show that $\Phi(\alpha)$ is not pre-minimal. Let $k \in \N$ be such that $\alpha(k, n) =1$ for infinitely many $n \in \N$. Let $N=\{ n\in\N\,:\, \alpha(k,n)=1\}$. For any $n\in N$ and $\lambda\in D^n_k$, $R_k(\alpha_{n+1}\lambda)=n$ since $\alpha_{n+1} \in \Delta_n - \{1_G\}$. Thus for all $n\in N$, $\lambda \in D^n_k$, and $\theta \in \Theta_k$,
$$\Phi(\alpha)(\alpha_{n+1} \lambda \theta b_{k-1}) = \alpha(k,R_k(\alpha_{n+1}\lambda))= 1 \neq 0 = \Phi(\alpha)(1_G\theta b_{k-1})=\Phi(\alpha)(\theta b_{k-1}).$$
This shows that $\Phi(\alpha)$ is not pre-minimal, hence is not minimal.
\end{proof}

The above proof has the following immediate corollary.

\begin{cor}\label{cor:minimalcoloringcomplexity}
For any countably infinite group $G$, the set of all minimal $2$-colorings on $G$ is ${\bf\Pi}^0_3$-complete.
\end{cor}

\begin{proof}
In the above proof we may suppose $c$ is a $\Delta$-minimal fundamental function with the property that any $x\in 2^G$ extending $c$ is a $2$-coloring. Such elements exist by Proposition~\ref{CANONICAL DMIN} and Corollary~\ref{GEN MINCOL}. Then for any $\alpha\in 2^{\N\times\N}$, $\Phi(\alpha)$ is a $2$-coloring, and $\alpha\in P$ iff $\Phi(\alpha)$ is a minimal $2$-coloring on $G$.
\end{proof}

\section{Flecc groups} \label{sec:flecc}

In the rest of this chapter we characterize the exact descriptive
complexity for the set of all $2$-colorings on a countably infinite
group. In contrast to Theorem~\ref{thm:minimalcomplexity} and
Corollary~\ref{cor:minimalcoloringcomplexity}, the set of all
$2$-colorings is not always ${\bf\Pi}^0_3$-complete. In this section
we isolate a group theoretic concept implying that the complexity is
simpler than ${\bf\Pi}^0_3$.

\begin{definition} \label{def:flecc}
A countable group $G$ is called {\it flecc} if there exists a finite
set $A\subseteq G-\{1_G\}$ such that for all $g\in G-\{1_G\}$ there is
$i\in {\mathbb Z}$ and $h\in G$ such that
$$ h^{-1}g^ih\in A. $$
\end{definition}\index{flecc group}

We first justify the terminology by giving a characterization for flecc groups.

\begin{definition} Let $G$ be a countable group and $g\in G$.
\begin{enumerate}
\item[(1)]
The {\it extended conjugacy class (ecc)} of $g$ is defined as the
set
$$ \ecc(g)=\{ h^{-1}g^ih\,:\, i\in {\mathbb Z},\, h\in G\}. $$
\item[(2)]
For $g$ of infinite order, we call the set $ \bigcap_{n\in {\mathbb N}} \ecc(g^n) $
 the {\it limit extended conjugacy class (lecc)} of $g$, and denote it by $\lecc(g)$.
\item[(3)]
If $g\neq 1_G$ is of finite order, we call any $\ecc(g^k)$, where $\mbox{order}(g)/k$ is prime,
a {\it lecc} of $g$.
\end{enumerate}
\end{definition}
\index{ecc}\index{lecc}\index{extended conjugacy class}\index{limit extended conjugacy class}

We need the following basic property of lecc's.

\begin{lem}\label{lem:leccpartition} Two lecc's coincide if they have a nontrivial intersection.
\end{lem}

\begin{proof}
We first claim that for any $g\in G$ of infinite order and $1_G\neq g'\in \lecc(g)$, $\lecc(g)=\lecc(g')$. On the one hand, it is obvious that $g'\in \ecc(g)$ and hence $\ecc(g')\subseteq \ecc(g)$
and $\ecc(g'^{n})\subseteq \ecc(g^n)$ for any $n\in\N$. Hence $\lecc(g')\subseteq \lecc(g)$. On the other hand, let $i\in\N_+$ be the least such that for some $h\in G$, $h^{-1}g^ih=g'$. Then $\ecc(g^i)=\ecc(g')$ and
so for any $n\in\N$, $\ecc(g^{in})\subseteq \ecc(g'^n)$. This implies that
$$\lecc(g)=\bigcap_{n\in\N}\ecc(g^n)\subseteq
\bigcap_{n\in\N} \ecc(g^{in})\subseteq \bigcap_{n\in\N} \ecc(g'^n)=\lecc(g'). $$
Thus $\lecc(g)=\lecc(g')$.

By a similar argument we can show that the same holds for $g\neq 1_G$ of finite order. If $\mbox{order}(g)/k=p$ is prime, then
$\langle g^k\rangle$ is a cyclic group of order $p$. Thus any nonidentity element in $\langle g^k\rangle$ is a
generator. The rest of the proof is similar as above.

Now suppose $\lecc(g)\,\cap\,\lecc(g')\neq\{1_G\}$. Let $k\in \lecc(g)\,\cap\, \lecc(g')$ so that $k\neq 1_G$.
Then $\lecc(g)=\lecc(k)=\lecc(g')$.
\end{proof}

\begin{prop}\label{prop:fleccchar} Let $G$ be a countable group. Then $G$ is flecc iff
\begin{itemize}
\item[(a)] for any $g\in G$ of infinite order, the lecc of $g$ is not $\{1_G\}$, and
\item[(b)] there are only finitely many distinct lecc's in $G$.
\end{itemize}
\end{prop}

\begin{proof}
$(\Rightarrow)$ Suppose $G$ is flecc. Let $A\subseteq G-\{1_G\}$ be finite such that for all $g\in G-\{1_G\}$ there is $i\in{\mathbb Z}$ and $h\in G$ such that $h^{-1}g^ih\in A$.
Fix $g\in G$ of infinite order. For any
$n\in\N$, $\ecc(g^{n!})\subseteq \ecc(g^n)$. Hence the lecc of $g$ can also be written as $\bigcap_{n\in\mathbb N} \ecc(g^{n!})$. Note that $\ecc(g^{n!})\supseteq \ecc(g^{(n+1)!})$ for all $n$. If $\bigcap_{n\in\N} \ecc(g^{n!})=\{1_G\}$, then for any $a\in A$,
there is $n_a\in \N_+$ such that $a\not\in \ecc(g^{n_a!})$. Let $n\geq n_a$ for all $a\in A$. Then $\ecc(g^{n!})\cap A=\emptyset$, contradicting the definition of flecc group. We thus have shown that (a) holds. It is clear that if $g\neq 1_G$ is of finite order, then any lecc of $g$ is also nontrivial.

To prove (b) we assume there are infinitely many distinct lecc's in $G$. Then by Lemma~\ref{lem:leccpartition} the pairwise intersections of these lecc's are trivial. Thus
for any finite subset $A\subseteq G-\{1_G\}$ there is $g\in G$ such that $\lecc(g)\neq\{1_G\}$ but $\lecc(g)\cap A=\emptyset$. By an argument similar to the one in the preceding paragraph, we get that for some $n$,
$\ecc(g^{n!})\neq\{1_G\}$ and $\ecc(g^{n!})\cap A=\emptyset$, contradicting the definition of flecc group.

$(\Leftarrow)$ Suppose (a) and (b) both hold. Then let $A\subseteq G$ be finite so that for any $g\in G$, $A\cap \lecc(g)\neq\emptyset$. Then in fact for any $g\in G$, $A\cap \ecc(g)\neq\emptyset$. This shows
that $G$ is flecc.
\end{proof}

Thus the terminology {\it flecc} represents the phrase that $G$ has only finitely many distinct limit
extended conjugacy classes. It does not appear that this concept has been  studied before. The rest of this section is devoted to a study of this concept.

Next we give some further characterizations of flecc groups. For simplicity we use $\mathbb{Z}^*$ to denote the set ${\mathbb Z}-\{0\}$, the set of all nonzero integers. \index{$\mathbb{Z}^*$}We also use $\sim$ to denote the conjugacy equivalence relation in the group $G$, i.e., $g\sim g'$ iff there is $h\in G$ such that $g'=h^{-1}gh$. Using this notation, we can express the fleccness of $G$ as there being a finite set $A\subseteq G$ of nonidentity elements such that for any nonidentity $g\in G$ there is $i\in\mathbb{Z}$ with $g^i\sim a$ for some $a\in A$, i.e., any nonidentity element of the group has a power which is conjugate to some element of $A$.

We also note the following characterization of nontriviality of lecc.

\begin{lem} \label{lem:leccchar}
Let $G$ be a group and $g \in G$ of infinite order. Then $\lecc(g)\neq \{1_G\}$
iff $\exists k \in \mathbb{Z}^*\
\forall n \in \mathbb{Z}^*\ \exists i \in \mathbb{Z}^*\ (g^{in} \sim g^k)$.
\end{lem}

\begin{proof}
Let $g \in G$ have infinite order, and suppose $\lecc(g)\neq \{1_G\}$.
Let $h \in \lecc(g)-\{1_G\}$. Since $h\in \bigcap_{n\in\N} \ecc(g^n)$, we have $\forall n \in \mathbb{Z}^*\
\exists i \in \mathbb{Z}^*\ (g^{in} \sim h)$. In particular, $h$ is conjugate to a
power of $g$. Let $g^k$ be such a power. Then
$\forall n \in \mathbb{Z}^*\
\exists i \in \mathbb{Z}^*\ (g^{in} \sim h \sim g^k)$.
Conversely, suppose $k\in\mathbb{Z}^*$ is such that $\forall n \in \mathbb{Z}^*\
\exists i \in \mathbb{Z}^*\ (g^{in} \sim g^k)$. Then $g^k\neq 1_G$ and $g^k\in\bigcap_{n\in\N}\ecc(g^n)$. Thus
$\lecc(g)\neq\{1_G\}$.
\end{proof}

When we try to determine whether a given group is flecc, it is easier to
consider the elements of finite order
and those of infinite order separately. Note that the $\lecc$ classes of the elements of finite order are just the conjugacy classes
of elements of prime order. Thus in a flecc group there are only finitely many conjugacy classes among the elements of prime order. We have the following alternative characterization.

\begin{prop} \label{prop:flecccharc}
A group $G$ is flecc iff all the following hold:
\begin{enumerate}
\item \label{prop:flecccharc_1}
There are only finitely many conjugacy classes among the elements of prime order.
\item \label{prop:flecccharc_2}
For every $g \in G$ of infinite order we have:
$$
\exists k \in \mathbb{Z}^*\ \forall n \in \mathbb{Z}^*\ \exists i \in \mathbb{Z}^*\ (g^{in} \sim g^k)$$
\item \label{prop:flecccharc_3}
There is a finite set $A$ of elements of infinite order
such that for any $g \in G$ of infinite order there is an $a \in A$ and
$i, m \in \mathbb{Z}^*$ such that $g^i \sim a^m$.
\end{enumerate}
\end{prop}

\begin{proof} $(\Rightarrow)$ This is immediate from Proposition~\ref{prop:fleccchar} and Lemma~\ref{lem:leccchar}. Only note that in (\ref{prop:flecccharc_3}) we may take $m=1$.

$(\Leftarrow)$ Assume (\ref{prop:flecccharc_2}) and  (\ref{prop:flecccharc_3}). It suffices to verify that
there is a finite set $B\subseteq G$ of elements of infinite order such that for all $g\in G$ of infinite order there is $i\in \mathbb{Z}^*$ such that $g^i\sim b$ for some $b\in B$. For this, let $A=\{ a_1,\dots, a_N\}$
be given as in  (\ref{prop:flecccharc_3}). For each $a_j \in A$, let $k_j \in \mathbb{Z}^*$
be given as in  (\ref{prop:flecccharc_2}).
Set $B=\{ a_1^{k_1},\dots, a_N^{k_N}\}$. We check that $B$ is as required. Suppose
$g \in G$ has infinite order. By
 (\ref{prop:flecccharc_3}) there is an $a_j \in A$ and $p, m\in\mathbb{Z}^*$ such that
$g^{p} \sim a_j^m$. By (2) and our choice of $k_j$, there is $i\in\mathbb{Z}^*$ such that $a_j^{im}\sim a_j^{k_j}$. Thus $g^{ip}\sim a_j^{im}\sim a_j^{k_j}\in B$.
\end{proof}

We note some basic properties of flecc groups and give some examples of flecc and nonflecc groups below.

Any finite group is obviously a flecc group. The following lemma gives some simple closure properties
of the flecc groups.

\begin{lem}\label{lem:fleccproduct} Let $G, H$ be countable groups. If $G \times H$ is flecc then
$G$ and $H$ are flecc. If $G$, $H$ are flecc and one of them is finite, then $G \times H$ is flecc.
Also, if $G_n$, $n\in\N$, are nontrivial, then $\oplus_n G_n$ is not flecc.
\end{lem}

\begin{proof}
Suppose first that $G\times H$ is flecc.
If $g\in G$, then $\lecc(g)\times\{1_H\}=\lecc(g,1_H)$. It
follows from Proposition \ref{prop:fleccchar} immediately that $G$ is
flecc if $G\times H$ is.

Suppose next that $G$, $H$ are flecc and $H$ is finite.
Let $A_1 \subseteq G-\{ 1_G\}$ witness that $G$ is flecc.
Let $A= (A_1 \times H) \cup (\{ 1_G\} \times H -\{ (1_G,1_H)\}$,
so $A$ is a finite subset of $G \times H -\{ (1_G, 1_H)\}$.
To see $A$ works, let $(g,h) \in G\times H-\{ (1_G, 1_H)\}$. If $g=1_G$, then
$h \neq 1_H$ and so $(1,h)$ is an element of $A$. If $g \neq 1_G$,
then for some $i$ and $g' \in G$ we have ${g'}^{-1} g^i g' =a_1 \in A_1$.
But then $(g', 1_H)^{-1} (g,h)^i (g', 1_H) \in A$.

To see the last statement, suppose that $G=\oplus_n G_n$, where each $G_n$ is nontrivial. Then for any $n$ and $1_{G_n} \neq g_n\in
G_n$, $\lecc(g_n)\times\prod_{m\neq n} \{1_{G_m}\}$ is an lecc in
$\oplus_n G_n$. If any of them is trivial, $\oplus_n G_n$ is not flecc. Assuming all of them are nontrivial, then there are
infinitely many distinct lecc's in $\oplus_n G_n$, so again $\oplus_n G_n$ is not flecc.
\end{proof}

Among countably infinite groups, the simplest example of a flecc group is the quasicyclic group ${\mathbb Z}(p^\infty)$. The fact that it is flecc is straightforward to check. The group ${\mathbb Z}$,
however, is not flecc. The following proposition completely characterize abelian flecc groups. Note
that this class of groups coincide with the class of all abelian groups with the minimal condition
(c.f. \cite{RobinsonBook} Theorem 4.2.11).

\begin{prop}\label{prop:abelianflecc} An abelian group is flecc iff it is a direct sum of finitely many quasicyclic groups and cyclic groups of prime-power order.
\end{prop}

\begin{proof}
$(\Leftarrow)$ By Lemma~\ref{lem:fleccproduct}
it suffices to show that a finite product of quasicyclic groups is flecc.
Consider
$$G=\mathbb{Z}(p_1^\infty) \times \cdots \times \mathbb{Z}(p_n^\infty).$$
Here we regard $\mathbb{Z}(p^j)$ as the mod $1$ additive group of fractions of the form $\frac{a}{p^j}$, where
$0\leq a<p^j$. Then $\mathbb{Z}(p^\infty)=\bigcup_{j\in\N}\mathbb{Z}(p^j)$. Since every element of $G$ has finite order, and the group is abelian, we only need to verify that there are only finitely many elements of prime order in $G$. For this note that given $g\in G$, i.e.,
$$ g=\left( \frac{a_1}{p_1^{j_1}},\,\cdots,\frac{a_n}{p_n^{j_n}}\right), $$
$g$ is of prime order iff $g$ is of order $p_k$ for some $1\leq k\leq n$. Moreover, assuming $a_k\neq 0\rightarrow (a_k,p_k)=1$ for all $1\leq k\leq n$, we have that $p_kg=0$ iff
\begin{enumerate}
\item[(i)] for all $l$ with $1\leq l\leq n$ and $p_l\neq p_k$, $a_l=0$, and
\item[(ii)] for all $l$ with $1\leq l\leq n$ and $p_l=p_k$, if $a_l\neq 0$ then $j_l=1$.
\end{enumerate}
Obviously, there are only finitely many elements in $G$ of this form.

$(\Rightarrow)$ Assume $G$ is an abelian flecc
group. Then $G$ can be written as the direct sum of a divisible
subgroup $D$ and a reduced group $R$.  Recall that a divisible abelian
group is a (possibly infinite) sum of quasicyclic groups and copies of
$\Q$.  From Lemma~\ref{lem:fleccproduct} the sum is actually a finite
sum in this case.  Also by Lemma~\ref{lem:fleccproduct}, both $D$ and
$R$ are also flecc.  Since an abelian flecc group must be a torsion
group (by Proposition~\ref{prop:fleccchar} (a) an abelian flecc group
cannot contain an element of infinite order), it follows that the
divisible group $D$ is a direct sum of finitely many quasicyclic
groups.  It remains to show that the reduced group $R$ is finite.
Again by Proposition~\ref{prop:fleccchar} $R$ is a torsion group.
Also the primary decomposition of $R$ cannot contain infinitely many
summands by Lemma~\ref{lem:fleccproduct}.  Thus the proof reduces to
the case $R$ being a reduced $p$-group. Now the definition of flecc in the
case of abelian $p$-groups is equivalent to there being only finitely
many distinct subgroups of order $p$. This implies that $R$ is finite,
as follows. Define a relation $\leq$ in $R$ by letting $g\leq h$ iff
$pg=h$. Then the transitive closure of $\leq$, still denoted by
$\leq$, is a partial order. $0=1_R$ is the largest element, and by the
fleccness $0$ has only finitely many immediate predecessors. This
implies that every element has finitely many immediate predecessors,
since if $pg_1=pg_2$ then $p(g_1-g_2)=0$. Thus $\leq$ is a finite
splitting tree. If $R$ is infinite then by K\"{o}nig's lemma there is
an infinite branch, which implies that there is a divisible subgroup
of $R$. But $R$ is reduced, contradiction. Thus we have shown that an
abelian flecc group is a direct sum of finitely many quasicyclic
groups and a finite group. A finite abelian group is a direct sum of
finitely many cyclic groups of prime-power order.
\end{proof}

If a countable group has finitely many conjugacy classes then it is
obviously flecc. By a well known theorem (c.f.\ \cite{RobinsonBook}
Theorem 6.4.6) of Higman, Neumann, and Neumann using HNN extensions,
every countable torsion-free group is the subgroup of a countable
group with only two conjugacy classes. It follows that every countable
torsion-free group is the subgroup of a countable flecc group. In fact, we have the following.

\begin{prop}
A group $G$ is a subgroup of a flecc group iff there are only finitely many primes
$p$ such that $G$ has an element of order $p$.
\end{prop}

\begin{proof}
If $p$, $q$ are distinct primes, then any elements $x$, $y$ of order $p$ and $q$
respectively cannot be conjugate in any group $H$ containing $G$. So, if $G$ is contained in
a flecc group then there can be
only finitely many primes $p$ such that there is an element of order $p$ in $G$.
Conversely, suppose that there are only finitely many such primes.
Call this set $P$.
By Higman, Neumann and Neumann, there is a group $H$ containing $G$
such that any two elements of $H$ of the same order are conjugate in $H$.
So, for each of the finitely many primes $p \in P$ there is only one
conjugacy class of elements of order $p$ in $H$. Also, the
HNN extension $H$ has the property that if $H$ has an element of
finite order $n$, then so does $G$.
Thus there are only finitely many
conjugacy classes of elements of prime-power order in $H$.
This shows that there are only finitely many flecc classes in $H$ for elements of finite order.
For elements of infinite order in $H$, just note that any
two elements of $H$ of infinite order are conjugate.
\end{proof}

Although the flecc groups are not closed under subgroups, we have the following
simple fact.

\begin{prop}
If $G$ is flecc and $H \unlhd G$, then $H$ is flecc.
\end{prop}

\begin{proof}
Let $A \subseteq G-\{ 1_G\}$ be finite satisfying
Definition~\ref{def:flecc}. Let $A'=A \cap H$. Let
$1_H \neq h \in H$. Then for some $i \in \mathbb{Z}^*$ and $a\in A$,
$h^i\sim a \in A$. As $H$ is normal in $G$,
$a \in H$, so $a \in A'$. Therefore, $A'$ witnesses that $H$ is flecc.
\end{proof}

We do not know if the class of flecc groups is closed under products
or quotients. The following is the best partial result we know of.

\begin{prop}
If $G$ is a flecc group and $H$ is a torsion flecc group, then $G \times H$ is flecc.
If $T$ is the torsion subgroup of the flecc group $G$, then $G/T$ is flecc.
\end{prop}

\begin{proof} Suppose $G$, $H$ are flecc and $H$ is torsion. To show $G\times H$ is flecc, we consider its elements of prime order and those of infinite order separately. For any prime $p$, any element of
$G \times H$ of order $p$ is of the form $(g,h)$ where $g$, $h$ have order $1$ or $p$.
But if $g \sim g'$ in $G$ and
$h \sim h'$ in $H$, then $(g,h) \sim (g',h')$ in $G \times H$.
It follows that there are only finitely many primes $p$
such that some element of $G$ or $H$ has order $p$.
Since  $G$, $H$ are flecc, there are only finitely many possibilities for the conjugacy
classes of $g$ and $h$ in $G$ and $H$ respectively.
Thus, there are only finitely many conjugacy classes of elements of prime order
in $G \times H$.

Turning to elements of infinite order, let $A \subseteq G$ be finite and consist of elements of infinite order
such that for all nonidentity $g\in G$ of infinite order there is $i\in\mathbb{Z}$ and $a\in A$ with $g^i\sim a$. Suppose $(g,h)$ has infinite order in $G \times H$. Since $H$ is torsion, $g$ must have infinite order in $G$. Let $i_0\in\mathbb{Z}^*$ be such that $h^{i_0}=1_H$, so
$(g,h)^{i_0}=(g^{i_0},1_H)$. Since $g^{i_0}$ still has infinite order in $G$,
there is an $i_1 \in \mathbb{Z}^*$ such that $(g^{i_0})^{i_1} \sim a$ for some $a \in A$.
Then $(g,h)^{i_0i_1} \sim (a,1_H)$. We have shown that $A \times \{ 1_H\}$ witnesses
 the fleccness for elements of infinite order in
$G \times H$.

Suppose next that $T$ is the  torsion subgroup of the flecc group $G$.
Since $G/T$ is torsion-free, we need only consider elements of infinite order in $G/T$.
Let $A \subseteq G$ be a finite set of elements of infinite order
such that for all $g\in G$ of infinite order there is $i\in \mathbb{Z}^*$ and $a\in A$ with $g^i\sim a$.
Then $\bar{A}=\{ aT \colon a \in A\} \subseteq G/T$
is a finite set of elements of infinite order in $G/T$. If $\bar{g}=gT$
has infinite order in $G/T$, then $g$ has infinite order in $G$ and so
for some $i \in \mathbb{Z}^*$ we have $g^i \sim a \in A$. But then
$\bar{g}^i \sim \bar{a}=Ta$ in $G/T$.
\end{proof}

Flecc groups are relevant to our study because of the following
observation.

\begin{lem}\label{lem:fleccblocking} Let $G$ be a countable flecc group and $x\in 2^G$. Let
$A\subseteq G-\{1_G\}$ be finite such that for all $g\in G-\{1_G\}$
  there is $i\in{\mathbb Z}$ and $h\in G$ such that $h^{-1}g^ih\in
  A$. Then $x$ is a $2$-coloring on $G$ iff $x$ blocks all $s\in A$,
  i.e., for all $s\in A$ there is a finite $T\subseteq G$ such that
$$ \forall g\in G\ \exists t\in T\ x(gt)\neq x(gst). $$
\end{lem}

\begin{proof} The $(\Rightarrow)$ direction is trivial. We only show $(\Leftarrow)$.
Assume $x$ is not a $2$-coloring on $G$. Then there is a periodic element $y\in \overline{[x]}$ with period $g$, i.e., $g\cdot y=y$. By fleccness there is $i\in {\mathbb Z}$ and $h\in G$ with $h^{-1}g^ih\in A$, and we have $(h^{-1}g^ih)\cdot (h^{-1}\cdot y)=h^{-1}\cdot y$. This means that there is
$s=h^{-1}g^ih\in A$ and $z=h^{-1}\cdot y\in \overline{[x]}$ such that $s\cdot z=z$.
By Corollary~\ref{cor:blockinglemma} $x$ does not block $s\in A$.
\end{proof}

\begin{theorem} If $G$ is a countably infinite flecc group, then the set of all $2$-colorings on $G$ is ${\bf\Sigma}^0_2$-complete.
\end{theorem}

\begin{proof}
If $G$ is a countably infinite flecc group, then the characterization in Lemma~\ref{lem:fleccblocking} for the set of all $2$-colorings on $G$ is ${\bf\Sigma}^0_2$. By Theorem~\ref{thm:Sigma02hard} this set is ${\bf\Sigma}^0_2$-hard, hence it is ${\bf\Sigma}^0_2$-complete.
\end{proof}

This completely characterizes the descriptive complexity of the set of all $2$-colorings for any
countably infinite flecc group.

\section{Nonflecc groups} \label{SECT NONFLECC}

In this section we show that the set of all $2$-colorings
on any countably infinite nonflecc group is ${\bf\Pi}^0_3$-complete. Since the proof is rather involved, we first illustrate the main ideas of the proof.

We again consider the ${\bf\Pi}^0_3$-complete set
$$ P=\{ \alpha\in 2^{\N\times\N}\,:\, \forall k\geq 1\ \exists n>k\ \forall m\geq n\ \alpha(k,m)=0\} $$
and define a continuous function $f:2^{\N\times\N}\to 2^G$ so that
for any $\alpha\in 2^{\N\times\N}$, $f(\alpha)$ is a $2$-coloring on $G$ iff $\alpha\in P$. To define $f(\alpha)$ we start with a fixed
$2$-coloring $x$ on $G$, identify infinitely many pairwise disjoint fixed finite sets $K_n$, and modify
the detail of $x$ on $K_n$ according to $\alpha$. When $\alpha\not\in P$, i.e., when there exists $k\geq 1$ such that $\alpha(k,n)=1$ for infinitely many $n>k$, the definition of $f(\alpha)\upharpoonright K_n$ for these infinitely many $n>k$ will give rise to a periodic element in $\overline{[f(\alpha)]}$. Denote
the period for this element by $s_k$. We will make sure that for each $n$ with $\alpha(k,n)=1$,
some left translate of $f(\alpha)\!\upharpoonright\! K_n$ already has period $s_k$. On the other hand, when $\alpha\in P$, we need to make sure that $f(\alpha)$ blocks all nonidentity $s\in G$. Thus in the
situation $\alpha(k,n)=1$ but $\alpha(1,n)=\dots=\alpha(k-1,n)=0$, we will make sure that $f(\alpha)$ blocks
enough $s$, e.g. all $s\in H_{k-1}$. In particular $f(\alpha)\upharpoonright K_n$ blocks all $s\in H_{k-1}$. Putting the two requirements together, when $\alpha(k,n)=1$ and $\alpha(1,n)=\dots=\alpha(k-1,n)=0$, we need some left translate of $f(\alpha)\upharpoonright K_n$ to have a period $s_k$ (with $s_k$ not depending on $n$) and at the same time to block all $s\in H_{k-1}$. In the following we first focus on showing that it is possible to construct such partial functions for nonflecc groups.

The next two theorems produce, for any countable nonflecc group $G$, a periodic element of $2^G$ that blocks a specific finite subset of elements in $G$. The method of proof is a variation of the fundamental method: we first create some layered regular marker structures, then use a membership test and an orthogonality scheme similar to the proof of Theorem~\ref{GEN COL}. The following theorem constructs the marker structures. These marker structures use objects $\Gamma_i$ which play roles similar to the sets $\Delta_i$ used in blueprints. For clarity, in the rest of this section we do not use the abbreviated notation $D^n_k$ but instead write out $\Delta_k \cap F_n$.

\begin{theorem} \label{thm:nonflecc1}
Let $G$ be a countable nonflecc group, $(\Delta_n,F_n)_{n\in\N}$ a centered blueprint, $k\geq 1$, and $A\subseteq G$ a finite set with $F_kF_k^{-1}\subseteq A$.
Suppose for any $i<j\leq k$ and $g\in G$, if $gF_i\cap F_j\neq\emptyset$ then  $gF_i\cap (\Delta_i \cap F_j) F_i\neq\emptyset$.
Then there is $s_k\in G$ and a sequence $(\Gamma_i)_{i\leq k}$ of subsets of $G$ such that
\begin{enumerate}
\item[(i)] $1_G\in \Gamma_i$ for all $i\leq k$;
\item[(ii)] for all $i<k$, $\Gamma_{i+1}\subseteq \Gamma_i$;
\item[(iii)] for all $i\leq k$, the $\Gamma_i$-translates of $F_i$ are maximally disjoint within $G$;
\item[(iv)] for all $g\in G$ and $l\in {\mathbb Z}^*$, $g^{-1}s_k^lg\not\in A-\{1_G\}$;
\item[(v)] for all $i\leq k$ and $g\in G$, $g\in \Gamma_i$ iff $s_k g\in \Gamma_i$;
\item[(vi)] for all $i\leq j\leq k$ and $\delta\in \Gamma_j$, $\Gamma_i\cap \delta F_j=\delta(\Delta_i\cap F_j)$;
\item[(vii)] for all $i\leq j\leq k$, $\gamma\in\Gamma_i$, and $\delta\in\Gamma_j$, if $\gamma F_i\cap \delta F_j\neq\emptyset$, then $\gamma F_i\subseteq \delta F_j$.
\end{enumerate}
\end{theorem}

\begin{proof} Since $(\Delta_n, F_n)_{n\in\N}$ is centered, we have $(\Delta_n)_{n\in\N}$ is decreasing, with $1_G\in \Delta_n$ for all $n\in\N$.

Let $A_0=A-\{1_G\}$. Since $G$ is nonflecc, there is $s_k\in G-\{1_G\}$ such
that for all $g\in G$ and $l\in{\mathbb Z}$, $g^{-1}s_k^lg\not\in A_0$. Fix such an $s_k$. Then (iv) is
satisfied. Next we define $\Gamma_{k-j}$ by induction on $0\leq j\leq k$. Fix an enumeration $1_G=g_0,g_1,\dots$ of all elements of $G$.

We first define $\Gamma_k$ in infinitely many stages. At each stage $m\in\N$ we define a set $\Gamma_{k,m}$, and eventually let $\Gamma_k=\bigcup_m \Gamma_{k,m}$. The sets $\Gamma_{k,m}$ are defined by induction on $m$. For $m=0$ let $\Gamma_{k,0}=\langle s_k\rangle$. In general suppose $\Gamma_{k,m}$ is already defined. If there is an $n$ such that $g_nF_k\cap \Gamma_{k,m}F_k=\emptyset$, let $n_m$ be the least such $n$, and let $\Gamma_{k,m+1}=\Gamma_{k,m}\cup\langle s_k\rangle g_{n_m}$. Otherwise let $\Gamma_{k,m+1}=\Gamma_{k,m}$. This finishes the definition of $\Gamma_{k,m}$ for all $m$, and also of $\Gamma_k$.

It is obvious that $1_G\in \Gamma_k$. Also obvious is that $\langle s_k\rangle \Gamma_{k,m}=\Gamma_{k,m}$ for all $m$, and therefore $\langle s_k\rangle \Gamma_k=\Gamma_k$, and (v) holds for $\Gamma_k$. Before defining $\Gamma_i$, $i<k$, we verify that (iii) holds for $\Gamma_k$.

First we show by induction on $m$ that the $\Gamma_{k,m}$-translates of $F_k$ are pairwise disjoint. For $m=0$ this reduces to the statement that for all $l\neq r\in {\mathbb Z}$, $s_k^lF_k\cap s_k^rF_k=\emptyset$. Since $F_kF_k^{-1}\subseteq A$, we have this required property by (iv). In general
suppose all $\Gamma_{k,m}$-translates of $F_k$ are pairwise disjoint. Suppose also $\Gamma_{k,m+1}=\Gamma_{k,m}\cup \langle s_k\rangle g_{n_m}$. Then by (iv) we have that
$s_k^lg_{n_m}F_k\cap s_k^rg_{n_m}F_k=\emptyset$ for all $l\neq r\in {\mathbb Z}$. Also, since $g_{n_m}F_k\cap \Gamma_{k,m}F_k=\emptyset$ and $\langle s_k\rangle\Gamma_{k,m}=\Gamma_{k,m}$, we
have that $\langle s_k\rangle g_{n_m}F_k\cap \Gamma_{k,m}F_k=\emptyset$. This shows that the $\Gamma_{k,m+1}$-translates of $F_k$ are pairwise disjoint. It follows that the $\Gamma_k$-translates
of $F_k$ are all pairwise disjoint. To see that the $\Gamma_k$-translates of $F_k$ form a maximally disjoint collection, simply note that if $g_mF_k\cap \Gamma_kF_k=\emptyset$, then $m<n_m$, contradicting the definition of $n_m$. We have thus defined $\Gamma_k$ to satisfy all requirements (i) through (v).

The version of (vi) that makes sense so far states that for all $\delta\in \Gamma_k$, $\Gamma_k\cap \delta F_k=\delta (\Delta_k\cap F_k)$, which is trivially true since both sides of the equation are the singleton $\{\delta\}$. The version of (vii) that makes sense so far states that if $\gamma,\delta\in \Gamma_k$ and $\gamma F_k\cap \delta F_k\neq\emptyset$ then $\gamma F_k\subseteq \delta F_k$. This follows immediately from the disjointness of $\Gamma_k$-translates of $F_k$. In fact, under the assumption we have $\gamma=\delta$ and therefore $\gamma F_k=\delta F_k$.

In general, suppose $\Gamma_{i+1}, \dots, \Gamma_k$ have been defined to satisfy (i) through (vii), we
define $\Gamma_i\supseteq \Gamma_{i+1}$ as follows. By induction on $m\in\N$ we define two increasing sequences $S_{i,m}$ and $\Gamma_{i,m}$,  and then take $\Gamma_i=\bigcup_m \Gamma_{i,m}$. For $m=0$ let
$$ S_{i,0}=\Gamma_{i+1}F_{i+1}\cup\dots\cup\Gamma_kF_k  $$
and
$$ \Gamma_{i,0}=\Gamma_{i+1}(\Delta_i\cap F_{i+1})\cup\dots\cup\Gamma_k(\Delta_i\cap F_k). $$
In general suppose $S_{i,m}$ and $\Gamma_{i,m}$ have been defined. If there is $n\in\N$ such that
$g_nF_i\cap S_{i,m}=\emptyset$, then let $n_m$ be the least such $n$, and let
$$ S_{i,m+1}=S_{i,m}\cup \langle s_k\rangle g_{n_m}F_i $$
and
$$ \Gamma_{i,m+1}=\Gamma_{i,m}\cup \langle s_k\rangle g_{n_m}. $$
This finishes the definition of $S_{i,m}$ and $\Gamma_{i,m}$ for all $m$, and hence that of $\Gamma_i$.

We verify that (i) through (vii) hold with this inductive construction. Since $1_G\in \Delta_i\cap F_{i+1}$, we have that $\Gamma_i\supseteq \Gamma_{i,0}\supseteq \Gamma_{i+1}$. It follows that (i) and (ii) hold. Also obvious is that $\langle s_k\rangle S_{i,m}=S_{i,m}$ and $\langle s_k\rangle \Gamma_{i,m}=\Gamma_{i,m}$, and (v) follows for $\Gamma_i$. Next we proceed to verify (iii), (vi) and (vii).

To show that all $\Gamma_i$-translates of $F_i$ are pairwise disjoint, it suffices to
show that all $\Gamma_{i,0}$-translates of $F_i$ are pairwise disjoint, since then the argument
as above will show inductively that the $\Gamma_{i,m}$-translates of $F_i$ are pairwise disjoint
for all $m>0$. Note that for
all $i<j\leq k$, the $(\Delta_i\cap F_j)$-translates of $F_i$ are pairwise disjoint and are contained in $F_j$ by the disjoint and coherent properties of a blueprint. Since all $\Gamma_j$-translates of $F_j$ are pairwise disjoint, it follows that all
$\Gamma_j(\Delta_i\cap F_j)$-translates of $F_i$ are pairwise disjoint. Next suppose $i<j< j'\leq k$,
$$ \gamma\in \Gamma_j,\ \delta\in \Delta_i\cap F_j, \ \gamma'\in \Gamma_{j'}, \ \delta'\in \Delta_i\cap F_{j'}, $$
and $\gamma\delta F_i\cap\gamma'\delta' F_i\neq\emptyset$. We need to show that $\gamma\delta =\gamma'\delta'$.

Note that $\delta F_i\subseteq F_j$ and $\delta'F_i\subseteq F_{j'}$, so we have
that $\gamma F_j\cap \gamma' F_{j'}\neq\emptyset$. By the inductive hypothesis (vii), $\gamma F_j\subseteq \gamma' F_{j'}$. Thus $\gamma\in \Gamma_j\cap\gamma' F_j=
\gamma' (\Delta_j\cap F_{j'})$ by the inductive hypothesis (vi). This means that there is $\eta\in \Delta_j\cap F_{j'}$ such that $\gamma=\gamma'\eta$.

Now $\gamma\delta=\gamma'\eta\delta\in \Gamma_{j'}(\Delta_i\cap F_{j'})$. This is because $\eta\delta\in \Delta_j(\Delta_i \cap F_j) \subseteq \Delta_i$ and $\eta\delta\in \eta F_j\subseteq F_{j'}$ by the coherent property of the blueprint $(\Delta_n, F_n)_{n\in\N}$. Since the $\Gamma_{j'}(\Delta_i\cap F_{j'})$-translates of $F_i$ are pairwise disjoint, we have  $\gamma\delta=\gamma'\delta'$ as needed.

Next we verify that the $\Gamma_i$-translates of $F_i$ are maximally disjoint within $G$. For this assume $g\in G$ is such that $gF_i\cap \Gamma_iF_i=\emptyset$. We claim first that $gF_i\cap S_{i,0}\neq\emptyset$. Otherwise $gF_i\cap S_{i,0}=\emptyset$. Let $g=g_m$. Then $m<n_m$, contradicting the definition of $n_m$. Now suppose $j>i$ and $gF_i\cap \Gamma_jF_j\neq\emptyset$. By the assumption
of the theorem $gF_i\cap \Gamma_j(\Delta_i\cap F_j)F_i\neq\emptyset$. Thus there is $\gamma\in \Gamma_j(\Delta_i\cap F_j)\subseteq \Gamma_{i,0}\subseteq \Gamma_i$ such that $gF_i\cap \gamma F_i\neq\emptyset$, a contradiction.

Now the inductive version of (vi) to be verified states that if $j$ is such that $i\leq j\leq k$ and $\delta\in \Gamma_j$,
then $\Gamma_i\cap \delta F_j=\delta(\Delta_i\cap F_j)$. If $j=i$ then this is trivially true since both sides of the equality are the singleton $\{\delta\}$. We assume $j>i$. By our definition $\delta(\Delta_i\cap F_j)\subseteq \Gamma_{i,0}$, and thus $\delta(\Delta_i\cap F_j)\subseteq \Gamma_i\cap \delta F_j$. Conversely, suppose $\gamma\in \Gamma_i\cap\delta F_j$. Then $\delta^{-1} \gamma \in F_j$ and in particular $\delta^{-1} \gamma F_i \cap F_j \neq \varnothing$. So by assumption on the blueprint, $\delta^{-1} \gamma F_i \cap (\Delta_i \cap F_j) F_i \neq \varnothing$ and hence $\gamma F_i \cap \delta (\Delta_i \cap F_j) F_i \neq \varnothing$. However, $\gamma_i \in \Gamma_i$ and $\delta (\Delta_i \cap F_j) \in \Gamma_i$, so $\gamma \in \delta (\Delta_i \cap F_j)$ since the $\Gamma_i$-translates of $F_i$ are disjoint. Thus $\Gamma_i \cap \delta F_j \subseteq \delta (\Delta_i \cap F_j)$ and $\Gamma_i \cap \delta F_j = \delta (\Delta_i \cap F_j)$, establishing (vi).

Finally we verify the inductive version of (vii) which states that if $j$ is such that $i\leq j\leq k$ and $\gamma\in \Gamma_i$, $\delta\in \Gamma_j$, and $\gamma F_i\cap \delta F_j\neq\emptyset$, then $\gamma F_i\subseteq \delta F_j$. For $j=i$ this follows immediately from the pairwise disjointness of $\Gamma_i$-translates of $F_i$. We assume $j>i$. Since $\delta F_j\subseteq S_{i,0}$, we have that $\gamma\in \Gamma_{i,0}$. Thus there is $j'$ with $i<j'\leq k$ and $\delta'\in \Gamma_{j'}$ such that $\gamma\in \delta'(\Delta_i\cap F_{j'})$. Let $j'$ be the smallest such index. By the coherent property of the blueprint $(\Delta_n, F_n)_{n\in\N}$ we have $\gamma F_i\subseteq \delta' F_{j'}$. If $j'\leq j$ then by the inductive hypothesis we have $\delta' F_{j'}\subseteq \delta F_j$ since
$\delta' F_{j'}\cap \delta F_j\supseteq \gamma F_i\cap \delta F_j\neq\emptyset$. It follows that $\gamma F_i\subseteq \delta F_j$ as we needed. If $j'>j$ we have $\delta F_j\subseteq \delta' F_{j'}$ from the inductive hypothesis again since $\delta F_j\cap\delta' F_{j'}\supseteq \delta F_j\cap \gamma F_i\neq\emptyset$. In particular $\delta\in \delta' F_{j'}$. By (vi) $\delta\in \delta'(\Delta_j\cap F_{j'})$. Now $\delta'^{-1}\delta\in \Delta_j\cap F_{j'}$ and $\delta'^{-1}\gamma\in \Delta_i\cap F_{j'}$, and $\delta'^{-1}\delta F_j\cap \delta'^{-1}\gamma F_i\neq\emptyset$. By the coherent property of the blueprint $(\Delta_n, F_n)_{n\in\N}$, we conclude that $\delta'^{-1}\gamma F_i\subseteq \delta'^{-1}\delta F_j$, and therefore $\gamma F_i\subseteq \delta F_j$ as we needed.
\end{proof}

Note that the assumption in the above theorem is true for the blueprint constructed in Theorem~\ref{EXIST STRONG BP} (clause (2) in that proof).
Thus it does not lose generality to assume that all blueprints we are working with have this property. In fact the sequence $(\Gamma_i, F_i)_{i\leq k}$ satisfies all axioms for a centered and maximally disjoint blueprint except that the length of the sequence is finite.

The next theorem gives the promised periodic element blocking a finite set of elements. The proof uses $\Gamma_i$ membership tests that are taken directly from the proof of Theorem~\ref{FM}. It also employs a coding technique similar to the
proof of Theorem~\ref{GEN COL}.

\begin{theorem} \label{thm:nonflecc2}
Let $G$ be a countable nonflecc group, $(\Delta_n, F_n)_{n\in\N}$ a centered blueprint guided by a growth sequence $(H_n)_{n\in\N}$ with $\gamma_1=1_G$, $R:H_0\to 2$ a nontrivial locally recognizable function, $k\geq 1$, and $A\subseteq G$ a finite set. Let $M_n=H_n\cup H_n^{-1}$ for all $n\in\N$.
Assume that
\begin{enumerate}
\item[(a)] for all $i<k$, $M_i^4\subseteq H_{i+1}$;
\item[(b)] for all $i\leq k$, $|\Lambda_i|>23\log_2|M_i|$;
\item[(c)] $M_k^{23}\subseteq A$;
\item[(d)] for any $i<j\leq k$ and $g\in G$, if $gF_i\cap F_j\neq\emptyset$ then $gF_i\cap (\Delta_i\cap F_j)F_i\neq\emptyset$.
\end{enumerate}
Let $s_k\in G$ and $\Gamma_i$, $i\leq k$, satisfy the conclusions of Theorem~\ref{thm:nonflecc1}.
Then there is $x_k\in 2^G$ such that
\begin{enumerate}
\item[(i)] $s_k^{-1}\cdot x_k=x_k$, i.e., for all $g\in G$, $x_k(s_kg)=x_k(g)$;
\item[(ii)] for all $1\leq i\leq k$, $x_k$ admits a simple $\Gamma_i$ membership test with test region a subset of $F_i$;
\item[(iii)] $x_k$ is compatible with $R$, and $x_k(g)=1-R(1_G)$ for all $g\in G-\Gamma_1(F_0\cup D^1_0)$;
\item[(iv)] for all $i\leq k$, if $\gamma, \gamma'\in \Gamma_i$ and $\gamma'\in \gamma M_i^{23}$,
then there is $t\in F_i$ such that $x_k(\gamma t)\neq x_k(\gamma't)$.
\end{enumerate}
\end{theorem}

\begin{proof}
Let $D=G-\bigcup_{1\leq i\leq k} \Gamma_i\Lambda_ib_{i-1}$. The displayed union in the equality is disjoint. Note that $\langle s_k\rangle D=D$.
Since $(\Gamma_i,F_i)_{i\leq k}$ satisfies all conditions required of a pre-blueprint for $i\leq k$, a definition of $x_k$ on $D$ can be made in the same fashion as in the proof of Theorem~\ref{FM}. This ensures
(iii) by clause (iv) of Theorem~\ref{FM} and the assumption $\gamma_1=1_G$. Also
$x_k\upharpoonright D$ admits a simple $\Gamma_i$ membership test for $i\geq 1$, with test
region a subset of $F_i$, by clause (ii) of Theorem~\ref{FM}. In fact, all conclusions of Theorem~\ref{FM} are true for $n\leq k$.
Since $\langle s_k\rangle \Gamma_i=\Gamma_i$ for all $i\leq k$, by Theorem~\ref{FM} (vi) and (vii), $x_k(s_kg)=x_k(g)$ for all $g\in D$.

To define $x_k$ on $G-D$ we use the technique presented in the proof of Theorem~\ref{GEN COL}. For each
$1\leq i\leq k$ let $R_i$ be the graph with vertex set $\Gamma_i$ and edge relation given by
$$ (\gamma,\gamma')\in E(R_i)\iff \gamma'\in \gamma M_i^{23}. $$
Since $M_i^{-1}=M_i$, $E(R_i)$ is a symmetric relation. A usual greedy algorithm gives a graph theoretic
$(|M_i|^{23}+1)$-coloring of $R_i$, $\mu_i: \Gamma_i\to \{0,1,\dots, |M_i|^{23}\}$. We claim that $\mu_i$
can be obtained so that $\mu_i(s_k\gamma)=\mu_i(\gamma)$ for all $\gamma\in \Gamma_i$. To see this,
apply the greedy algorithm in infinitely many stages as follows. Enumerate all elements of $\Gamma_i$ as $1_G=\gamma_0, \gamma_1,\dots$. At stage $m=0$, let $S_0=\langle s_k\rangle \gamma_0$, assign $\mu_i(\gamma_0)$ arbitrarily, and then let $\mu_i(s_k^l\gamma_0)=\mu_i(\gamma_0)$ for all $l\in{\mathbb Z}$. Since for any $g\in G$ and $l\in{\mathbb Z}^*$, $g^{-1}s_k^lg\not\in A-\{1_G\}$ and $M_i^{23}\subseteq M_k^{23}\subseteq A$, we have $g^{-1}s_k^lg\not\in M_i^{23}-\{1_G\}$. It follows that
for any $l\neq r\in {\mathbb Z}$, $(s_k^l\gamma_0, s_k^r\gamma_0)\not\in E(R_i)$. In general suppose $S_m$
 and $\mu_i\upharpoonright S_m$ have been defined. We define $S_{m+1}$ by induction. If there is $n$ such
 that $\gamma_n\not\in S_m$, let $n_m$ be the least such $n$, and let $S_{m+1}=S_m\cup\langle s_k\rangle \gamma_{n_m}$. Define $\mu_i(\gamma_{n_m})$ arbitrarily using the greedy algorithm: since $\gamma_{n_m}$ has at most $|M_i|^{23}$ many adjacent vertices there is some $v\in \{0,1,\dots, |M_i|^{23}\}$ such that by assigning $\mu_i(\gamma_{n_m})=v$ the resulting function is a partial graph-theoretic coloring. Arbitrarily choose such a $v$, and let $\mu_i(s_k^l\gamma_{n_m})=v$ for all $l\in{\mathbb Z}$. By a similar argument as above, $(s_k^l\gamma_{n_m}, s_k^r\gamma_{n_m})\not\in E(R_i)$ for any $l\neq r\in {\mathbb Z}$. Suppose $(s_k^l\gamma_{n_m},\psi)\in E(R_i)$ for some $l\in\mathbb{Z}$ and $\psi\in S_m$. Then $(\gamma_{n_m},s_k^{-l}\psi)\in E(R_i)$, where $s_k^{-l}\psi\in S_m$ since $\langle s_k\rangle S_m=S_m$. By induction $\mu_i(\psi)=\mu_i(s_k^{-l}\psi)$. Thus
 $$\mu_i(s_k^l\gamma_{n_m})=v=\mu_i(\gamma_{n_m})\neq \mu_i(s_k^{-1}\psi)=\mu_i(\psi). $$
 We conclude that $\mu_i\upharpoonright S_{m+1}$ is a partial
 graph-theoretic coloring. To summarize, we have defined $\mu_i$ on all of $\Gamma_i$ so that $\mu_i(s_k\gamma)=\mu_i(\gamma)$ for all $\gamma\in \Gamma_i$.

 The rest of the proof follows that of Theorem~\ref{GEN COL}. We thus constructed $x_k$ to satisfy (i) through (iv).
\end{proof}

Again, the assumptions of the theorem are easy to arrange. Growth sequences satisfying (a) and (b) can be obtained by Corollary~\ref{GROW BP}, since $|M_n|\leq 2|H_n|$; condition (c) is simply a largeness condition for the finite set $A$; blueprints satisfying (d), as we remarked before, arise from the proof of Theorem~\ref{EXIST STRONG BP}. One might have noticed that some of these hypotheses are not needed in the theorem. They will be only needed in some of our later proofs, but we state them here just to keep our assumptions about the blueprint explicitly isolated.

We note the following corollary of the proof of Theorem~\ref{thm:nonflecc2}.

\begin{cor}\label{cor:nonflecc3} Let $G$ be a countable nonflecc group and $A\subseteq G-\{1_G\}$ be finite. Then there is a periodic $x\in 2^G$ which blocks all elements of $A$.
\end{cor}

For the convenience of our later arguments we also note some finer details of the $\Gamma_i$ membership test constructed in the above theorem. By the proof of Theorem~\ref{FM}, we have the following ``standard form" of the membership tests. For the $\Gamma_1$ membership test, we have
$$ g\in \Gamma_1\mbox{ iff } \forall f\in F_0\ x_k(gf)=R(f). $$
Note that here we used the hypothesis that $\gamma_1=1_G$. The general $\Gamma_i$ membership test for $i>1$ is defined by induction on $i$ with test region $V_i=\gamma_i(V_{i-1}\cup\{a_{i-1},b_{i-1}\})$. Specifically, for $i>1$,
$$ g\in \Gamma_i\mbox{ iff } g\gamma_i\in \Gamma_{i-1} \wedge x_k(g\gamma_ia_{i-1})=x_k(g\gamma_ib_{i-1})=1. $$
An important feature of these membership tests is that they only depend on the locally recognizable function $R$ and the blueprint $(\Delta_n,F_n)_{n\in\N}$. In particular, they do not depend on the number $k\geq 1$. In other words, if $k\neq k'\geq 1$ and if the above theorem is applied to $k$ and $k'$ respectively, with the same locally recognizable function $R$ and the same blueprint $(\Delta_n, F_n)_{n\in\N}$, then the resulting membership tests take the same form. This point will come up in the proof of our main theorem below.

 In the proof of our main theorem only a finite part of $x_k$ will be used at each stage. However, we need such finite parts to maintain their integrity when it comes to $\Gamma_i$ membership tests for $i\leq k$. For this purpose we define the
following saturation operation for finite sets. Given a finite set $B\subseteq G$, let
$$ \mathrm{sat}_0(B)=\bigcup\{\gamma F_0\,:\, \gamma\in \Gamma_0, \ \gamma F_0\cap B\neq \emptyset\} $$
and
$$\mathrm{sat}_k(B)=\mathrm{sat}_0(B)\cup \bigcup_{i\leq k}(\Gamma_i\cap \mathrm{sat}_0(B))F_i. $$
It is important to note that $B$ is not necessarily contained in either $\mathrm{sat}_0(B)$ or $\mathrm{sat}_k(B)$, but $B\cap \Gamma_0F_0\subseteq \mathrm{sat}_0(B)$ by definition. Moreover, $\mathrm{sat}_k(B)$ has the obvious property that for all $\gamma\in \Gamma_0$, if $\gamma F_0\cap B\neq \emptyset$ then $\gamma F_0\subseteq \mathrm{sat}_k(B)$. We also have the following strengthened property.

\begin{lem}\label{lem:satk} For all $i\leq k$ and $\gamma\in \Gamma_i\cap\mathrm{sat}_k(B)$, $\gamma F_i\subseteq \mathrm{sat}_k(B)$.
\end{lem}

\begin{proof}
Fix $i\leq k$ and $\gamma\in \Gamma_i\cap \mathrm{sat}_k(B)$. If $\gamma\in \mathrm{sat}_0(B)$ then $\gamma\in \Gamma_i\cap \mathrm{sat}_0(B)$ and therefore $\gamma F_i\subseteq \mathrm{sat}_k(B)$ by definition. Suppose instead $\gamma\in (\Gamma_j\cap \mathrm{sat}_0(B))F_j$ for some $j\leq k$. Then there is some $\delta\in \Gamma_j\cap \mathrm{sat}_0(B)$ such that $\gamma\in \delta F_j$. If $j\geq i$, then by clause (vii) of Theorem~\ref{thm:nonflecc1} $\gamma F_i\subseteq \delta F_j$, and therefore $\gamma F_i\subseteq \delta F_j\subseteq \mathrm{sat}_k(B)$ by definition. If $j<i$, then by clause (ii) of Theorem~\ref{thm:nonflecc1} $\gamma\in \Gamma_j$. Since the $\Gamma_j$-translates of $F_j$ are pairwise disjoint, we have $\gamma=\delta$. This means that $\delta\in \Gamma_i$, and so $\gamma F_i=\delta F_i\subseteq (\Gamma_i\cap \mathrm{sat}_0(B))F_i\subseteq\mathrm{sat}_k(B)$.
\end{proof}

For any $n>k$, we also define
$$ K_{n,k}=\mathrm{sat}_k(M_nM_k^3M_{k-1}^3\dots M_0^3), $$
where $M_i=H_i\cup H_i^{-1}$,
and define $x_k^n=x_k\upharpoonright K_{n,k}$. $x^n_k$ will be the finite part of $x_k$ used in our main construction.

In our main construction the background will be colored differently from some translates of the regions $K_{n,k}$. However, in the coloring of both parts we use some blueprint and some membership test for center points. To make the membership tests distinct we need the following lemma similar to Proposition~\ref{LR EXT}.

\begin{lem}\label{lem:locallyorthogonal} Let $G$ be a countably infinite group, $B\subseteq G$ a finite set, and $Q: B\to 2$ any function. Then there exist a finite set
$A\supseteq B$ and two nontrivial locally recognizable functions $R, R': A\to 2$ both extending $Q$ such that
for all $x, x'\in 2^G$ with $x\upharpoonright A=R$ and $x'\upharpoonright A=R'$,
$$ \forall g\in A\ \exists h\in A\ x(gh)\neq x'(h) $$
and
$$ \forall g\in A\ \exists h\in A\ x'(gh)\neq x(h). $$
\end{lem}

\begin{proof} The proof is also similar to that of Proposition~\ref{LR EXT}. For clarity we give a self-contained argument below.
By defining $Q(1_G) = 0$ if necessary, we may assume $1_G \in B$. Set $B_1 = B$. Choose distinct elements $a, b, a', b' \in G - B_1$ and set $B_2 = B_1 \cup \{a, b, a', b'\}$. Next chose any $c \in G - (B_2 B_2 \cup B_2 B_2^{-1})$ and set $B_3 = B_2 \cup \{c\} = B_1 \cup \{a, b, a', b', c\}$. Let $A = B_3 B_3$. Define $R: A \rightarrow 2$ by
$$R(g) = \begin{cases}
Q(g) & \text{if } g \in B_1 \\
Q(1_G) & \text{if } g \in \{a, b, c\} \\
1- Q(1_G) & \text{if } g\in \{a', b'\} \\
1 - Q(1_G) & \text{if } g \in A - B_3.
\end{cases}$$
$R'$ is similarly defined, with the role of $a, b$ and respectively of $a', b'$ interchanged. Thus for all $g\in A-\{a, b, a', b'\}$, $R(g)=R'(g)$, and for
$g\in \{a, b, a', b'\}$, $R(g)=1-R'(g)$. It is obvious that both $R$ and $R'$ are nontrivial (Definition \ref{DEFN LR}).

We claim that for any nonidentity $g\in B_3$, at least one of $a$, $b$, or $c$ is not an element of $gB_3$. To see this, consider the following cases. \underline{Case 1}: $g \in B_2$. Then $c \not\in g B_2 \subseteq B_2 B_2$ and $c \neq g c$ since $g \neq 1_G$. Thus $c \not\in g B_3$. \underline{Case 2}: $g \in B_3 - B_2 = \{c\}$. Then $g = c$. Since $c \not\in B_2 B_2^{-1}$, we have $a, b \not\in c B_2$. If $a, b \in c B_3$ then we must have $a = c^2 = b$, contradicting $a \neq b$. We conclude $\{a, b\} \not\subseteq c B_3 = g B_3$.

By symmetry we also have that for any $g\in B_3$, at least one of $a'$, $b'$, or $c$ is not an element of $gB_3$.

We verify that $R$ is locally recognizable. Towards a contradiction, suppose there is $y \in 2^G$ extending $R$ such that for some $1_G \neq g \in A$, $y(g h) = y(h)$ for all $h \in A$. In particular, $y(g) = y(1_G) = R(1_G)$ so $g \in B_3$. By the above claim we have $\{a, b, c\} \not\subseteq g B_3 \subseteq A$, but
$\{a, b, c\} \subseteq \{h \in A \ | \ y(h) = y(1_G)\} \subseteq B_3$. Therefore
$$|\{h \in B_3 \ | \ y(g h) = y(1_G)\}| < |\{h \in A \ | \ y(h) = y(1_G)\}|$$
$$= |\{h \in B_3 \ | \ y(h) = y(1_G) \}| = |\{h \in B_3 \ | \ y(gh) = y(1_G)\}|,$$
a contradiction.

A symmetric argument with $a$, $b$ replaced by $a'$, $b'$, respectively proves that $R'$ is locally recognizable.

To complete the proof of the lemma, we let $x, x'\in 2^G$
extending $R$, $R'$, respectively, and toward a contradiction assume that for some $g\in A$, $x(gh)=x'(h)$ for all $h\in A$. Since $R\neq R'$ we know that $g\neq 1_G$. Since
$R(g)=x(g)=x'(1_G)=Q(1_G)$, we also know that $g\in B_3$. Again we have $\{a, b, c\}\not\subseteq gB_3\subseteq A$, and therefore
$$\begin{array}{rl}&|\{h \in B_3 \ | \ x(g h) = Q(1_G)\}| < |\{h\in A\ |\ x(h)=Q(1_G)\}| \\
=& |\{h \in A \ | \ x'(h) = Q(1_G)\}| = |\{h \in B_3 \ | \ x'(h) = Q(1_G) \}| \\
=& |\{h \in B_3 \ | \ x(gh) = Q(1_G)\}|,
\end{array}
$$
a contradiction. Finally a symmetric argument gives that for all $g\in A$ there is $h\in A$ such that $x'(gh)\neq x(h)$. This finishes the proof of the lemma.
\end{proof}

We are now ready to prove our main theorem of this section.

\begin{theorem} Let $G$ be a countable nonflecc group. Then the set of all $2$-colorings on $G$ is
${\bf\Pi}^0_3$-complete.
\end{theorem}

\begin{proof}
Let $H_0\subseteq G$ be a finite set with $1_G\in H_0$, and $R, R': H_0\to 2$ be two distinct nontrivial locally recognizable functions
given by Lemma~\ref{lem:locallyorthogonal}, with $R(1_G)=R'(1_G)=1$. Let $(H_n)_{n\in\N}$ be a growth sequence with the additional properties that, for all $n\in\N$, letting $M_n=H_n\cup H_n^{-1}$,
\begin{enumerate}
\item[(1)] $M_n^4\subseteq H_{n+1}$;
\item[(2)] $M_n^3M_{n-1}^3\dots M_0^3\subseteq H_{n+1}$.
\end{enumerate}
Such sequences are easy to construct. Let $(\Delta_n, F_n)_{n\in\N}$ be a centered blueprint guided by $(H_n)_{n\in\N}$, with $\gamma_1=1_G$, such that
\begin{enumerate}
\item[(3)] for all $n\in\N$, $|\Lambda_n|>23\log_2|M_n|$+1.
\end{enumerate}
The existence of such blueprints follows from Corollary~\ref{GROW BP}, since $|M_n|\leq 2|H_n|$ for all $n\in\N$. The construction of the blueprint $(\Delta_n, F_n)_{n\in\N}$ follows the proof of Theorem~\ref{EXIST STRONG BP}, and therefore we have
\begin{enumerate}
\item[(4)] for all $i<j$ and $g\in G$, if $gF_i\cap F_j\neq\emptyset$, then $gF_i\cap (\Delta_i\cap F_j)F_i\neq\emptyset$.
\end{enumerate}
Now apply Theorem~\ref{thm:strongcoloring} to obtain a strong $2$-coloring $z\in 2^G$ fundamental with respect to $(\Delta_n, F_n)_{n\in\N}$ and compatible with $R'$.

For each $n\geq 1$, let $K_n=M_n^4$. Fix an enumeration of $G-\{1_G\}$ as $\sigma_1,\sigma_2,\dots$ so that each $s\in G-\{1_G\}$ is enumerated infinitely many times. We inductively define two sequences $(\pi_n)_{n\geq 1}$ and $(w_n)_{n\geq 1}$ of elements of $G$ so that
\begin{enumerate}
\item[(5)] for all $n\neq m\in\N$, the sets $\pi_nK_n$, $\pi_mK_m$, $\{w_n,\sigma_nw_n\}$, and $\{w_m, \sigma_mw_m\}$ are pairwise disjoint;
\item[(6)] for all $n\geq 1$, $z(w_n)\neq z(\sigma_nw_n)$; and
\item[(7)] for all $1\leq n<n'$, $\pi_nK_nM_n^{23}\cap \pi_{n'}K_{n'}=\emptyset$.
\end{enumerate}
For $n=1$ let $\pi_1=1_G$. Since $z$ is a strong $2$-coloring, there are infinitely many $w\in G$ such that $z(\sigma_1w)\neq z(w)$. Let $w_1\not\in \pi_1K_1\cup\sigma_1^{-1}\pi_1K_1$ be such a $w$. Then $w_1, \sigma_1w_1\not\in \pi_1K_1$. So
$\{w_1,\sigma_1w_1\}\cap \pi_1K_1=\emptyset$. In general suppose $\pi_m$ and $w_m$ have been
defined for all $1\leq m\leq n$ to satisfy (5) through (7). Let
$$ B=\bigcup_{1\leq m\leq n}\left(\{w_m,\sigma_mw_m\}\cup \pi_mK_mM_m^{23}\right). $$
Then $B$ is finite and we may find $\pi_{n+1}\not\in BK_{n+1}^{-1}$ and $w_{n+1}\not\in \{1_G, \sigma^{-1}\}(B\cup \pi_{n+1}K_{n+1})$. Then we have that $\pi_{n+1}K_{n+1}\cap B=\emptyset$
and $\{w_{n+1},\sigma_{n+1}w_{n+1}\}\cap (B\cup \pi_{n+1}K_{n+1})=\emptyset$. Therefore the sequences
are as required.

For each $k\geq 1$ we apply Theorem~\ref{thm:nonflecc1} with $A=M_k^{23}$ to obtain an $s_k\in G$ and a sequence $(\Gamma_i)_{i\leq k}$. Note that $s_k$ only depends on $k$. The sequence $(\Gamma_i)_{i\leq k}$ also depends on $k$. However, for simplicity we will not introduce $k$ into the notation. The reader should be aware that at various places we might be referring to different $\Gamma_i$'s coming from different $k$ values. With $k$ fixed we continue to apply Theorem~\ref{thm:nonflecc2} to obtain $x_k$ compatible with $R$. For $n>k$, we also defined the set $K_{n,k}$ and $x^n_k=x_k\upharpoonright K_{n,k}$. Note that by (1) and (2) above,
$$\begin{array}{rl} K_{n,k}=\mathrm{sat}_k(M_nM_k^3\dots M_0^3)&\subseteq (M_nM_k^3\dots M_0^3)F_0^{-1}F_0F_k \\ \\
&\subseteq M_nM_{k+1}M_k^3\subseteq M_n^3\subseteq K_n.
\end{array}$$

We are finally ready to define a continuous function $f: 2^{\N\times\N}\to 2^G$ so that $f(\alpha)$ is a $2$-coloring on $G$ iff $\alpha\in P$, where
$$ P=\{ \alpha\in 2^{\N\times\N}\,:\, \forall k\geq 1\ \exists n>k\ \forall m\geq n\ \alpha(k,m)=0\}. $$
Given $\alpha\in 2^{\N\times\N}$, we let $f(\alpha)(g)=z(g)$ if $g\not\in \bigcup_{n\geq 1}\pi_nK_n$. Then for each $n\geq 1$, $f(\alpha)\upharpoonright \pi_nK_n$ will be defined according to
the values $\alpha(1,n),\dots,\alpha(n,n)$, as follows. Let $1\leq k<n$ be the least such that $\alpha(k,n)=1$. If $k$ is undefined then we define $f(\alpha)(g)=z(g)$ for all $g\in \pi_nK_n$.
Suppose $k$ is defined. Then we let $f(\alpha)(\pi_ng)=x^n_k(g)$ for all $g\in K_{n,k}$. For $g\in \pi_n(K_n-K_{n,k})$, if $g\not\in \Delta_0F_0$ then note that $z(g)=0$ and we let $f(\alpha)(g)=0$ as well.
If $g\in \pi_n(K_n-K_{n,k})$ but $g\in\Delta_0F_0$, let $\gamma\in \Delta_0$ be the unique element such that
$g\in\gamma F_0$. If $\gamma F_0\cap \pi_nK_{n,k}=\emptyset$ we let $f(\alpha)(g)=z(g)$, otherwise let $f(\alpha)(g)=0$. Note that in case $\gamma F_0\cap \pi_nK_{n,k}\neq\emptyset$, we have
$$ \gamma F_0\subseteq \pi_nK_{n,k}F_0^{-1}F_0\subseteq \pi_n M_n^3M_1\subseteq \pi_n K_n, $$
and therefore $f(\alpha)\upharpoonright \gamma F_0$ is well defined. This finishes the definition of $f(\alpha)$.

It is obvious that each $f(\alpha)\in 2^G$ and it is routine to check that $f$ is a continuous function. We argue first that if $\alpha\not\in P$ then $f(\alpha)$ is not a $2$-coloring. Assume $\alpha\not\in P$. Let $k\geq 1$ be the least such that for infinitely many $n>k$, $\alpha(k,n)=1$. Let $n_0$ be large enough such that for all $1\leq i<k$ and $n\geq n_0$, $\alpha(i,n)=0$. Then for infinitely many $n>n_0$, $k$ is the least $m$ such that $\alpha(m,n)=1$. By the definition of $f(\alpha)$, for infinitely many such $n>k$, $f(\alpha)(\pi_ng)=x^n_k(g)$ for $g\in K_{n,k}$. In other words, for infinitely many $n>k$,
$(\pi_n^{-1}\cdot f(\alpha))\upharpoonright K_{n,k}=x^n_k$. We claim that $(\pi_n^{-1}\cdot f(\alpha))\upharpoonright H_n=x_k\upharpoonright H_n$. This implies that
$x_k\in \overline{[f(\alpha)]}$. Since $x_k$ is periodic, $f(\alpha)$ is not a $2$-coloring. To prove the claim, fix such an $n>k$ and let $g\in H_n$.
If $g\in K_{n,k}$ then there is nothing to prove. Assume $g\not\in K_{n,k}$. Since the $\Gamma_0$-translates of $F_0$ are maximally disjoint within $G$, there is $\gamma\in \Gamma_0$ such that $gF_0\cap \gamma F_0\neq\emptyset$. For any such $\gamma$, we have
$\gamma\in gF_0F_0^{-1}\subseteq H_nF_0F_0^{-1}\subseteq K_n$, and therefore $\gamma F_0\subseteq K_{n,k}$. This implies that $g\not\in \Gamma_0F_0$. It follows from our definition of $f(\alpha)$ that $f(\alpha)(\pi_ng)=0$. Also by Theorem~\ref{thm:nonflecc2} (iii) we have $x^n_k(g)=0$. This completes the proof of the claim, and hence we have shown that if $\alpha\not\in P$ then $f(\alpha)$ is not a $2$-coloring.

The rest of the proof is devoted to showing that if $\alpha\in P$, then $f(\alpha)$ is a $2$-coloring. We first note that for any $n\geq 1$, since $w_n, \sigma_nw_n\not\in \bigcup_{m\geq 1} \pi_mK_m$, we
have $f(\alpha)(w_n)=z(w_n)\neq z(\sigma_nw_n)=f(\alpha)(\sigma_nw_n)$. Thus in particular $f(\alpha)$ is aperiodic. By Lemma~\ref{lem:almostcoloringlemma}, to show that $f(\alpha)$ is a $2$-coloring it suffices to show that it is a near $2$-coloring. Fix $\alpha\in P$ and any $i\geq 1$.  Fix any $s\in H_i$ with $s\neq 1_G$. Let $n_0>i$ be large enough such that for all $n\geq n_0$ and $1\leq j\leq i$, $\alpha(j,n)=0$. Let $S=\bigcup_{1\leq m\leq n_0} \pi_mK_m$ and $T=M_i^{11}F_i$. We verify that $f(\alpha)$ nearly blocks $s$
by showing that for all $g\not\in SF_i^{-1}F_iF_i^{-1}\{1_G, s^{-1}\}$ there is $t\in T$ with
$f(\alpha)(gst)\neq f(\alpha)(gt)$.

To simplify notation denote $f(\alpha)$ by $y$. Also denote, for $k\geq 1$,
$$ N_k=\{ n\geq 1\,:\, n>k \mbox{ and $k$ is the least $m$ with } \alpha(m,n)=1\} $$
and
$$ X_k=\bigcup\{ \pi_nK_{n,k}\,:\, n\in N_k \}. $$
Then $X_k$ is the set on which the definition of $y$ is given by the periodic element $x_k$. Let
$$ N=\bigcup_{k\geq 1}N_k \ \ \mbox{ and }\ \ X=\bigcup_{k\geq 1} X_k. $$
Then $N$ contains (and is most likely equal to) the set of all $n\geq 1$ such that $y\upharpoonright \pi_nK_n\neq z\upharpoonright \pi_nK_n$. Note that $N_k\cap N_{k'}=\emptyset$ and $X_k\cap X_{k'}=\emptyset$ for $k\neq k'\geq 1$.

For a fixed $k\geq 1$, if $j\leq k$, we define
$$ \Gamma_j^*=\bigcup\{ \pi_n\Gamma_j\cap X_k\,:\, n\in N_k\}. $$
Then $\Gamma_j^*F_j\subseteq X_k$ by the definition of $K_{n,k}$ and Lemma~\ref{lem:satk}. By Theorem~\ref{thm:nonflecc2} elements of $\Gamma_j^*$, $j\geq 1$, satisfy a simple membership test induced by $R$ with test region a subset of $F_j$. Moreover, by the remark following Theorem~\ref{thm:nonflecc2}, these membership tests do not depend on $k$ since they are obtained from applying Theorem~\ref{FM} by using the same locally recognizable function $R$ and the same blueprint. In other words, if $k\neq k'\geq 1$ and $j\leq k, k'$, then the membership tests for elements of $\Gamma_j$ in the constructions of $x_k$ and $x_{k'}$ take the same form. For this reason, and for simplicity of notation, we refrained from mentioning $k$ in the notation $\Gamma_j^*$. We will refer to the membership test involved as simply the {\it $\Gamma_j$ membership test}.

For any $j\in\N$ define
$$ \Delta_j^*=\{ \gamma\in \Delta_j\,:\, \gamma F_j\cap X=\emptyset\}. $$
In particular, if $k\geq 1$ and $j\leq k$, then $\Delta_j^*F_j\cap X_k=\emptyset$.
Since the construction of the strong $2$-coloring $z$ also comes from the proof of Theorem~\ref{FM} by using the same blueprint, elements of $\Delta_j$, $j\geq 1$, satisfy a similar simple membership test (for $z$), except it is induced
by $R'$ instead, also with test region a subset of $F_j$. We refer to it as the {\it $\Delta_j$ membership test}.

We remark that the $\Delta_j$ membership test (for $z$) takes exactly the same form as the $\Gamma_j$ membership test (for any $x_k$ with $k \geq j$), except that instead of the locally recognizable function $R$ we use $R'$.
Since $R$ and $R'$ are distinct, the $\Gamma_1$ membership test (for any $x_k$ with $k \geq 1$) and the $\Delta_1$ membership test (for $z$) are different. For $j>1$, the $\Gamma_j$ ($\Delta_j$) membership test is constructed by the same induction using $\Gamma_{j-1}$ ($\Delta_{j-1}$) membership tests. Hence the $\Gamma_j$ membership test and the $\Delta_j$ membership tests are also different.

We note that elements of $\Delta_j^*$ satisfy the $\Delta_j$ membership test on $y$. Conversely, we do not necessarily have that elements satisfying the $\Delta_j$ membership test on $y$ must be in $\Delta_j^*$. Instead, we note that if $g\in G$ is such that
$gF_j\cap X=\emptyset$ and $g$ satisfies the $\Delta_j$ membership test, then $g\in \Delta_j^*$. This is easily seen by induction on $j\geq 1$. When $j=1$ we assume $gF_1\cap X=\emptyset$ and $g$ satisfies the $\Delta_1$ membership test, which is $y(ga)=R'(a)$ for all $a\in F_0$. By the properties of $R'$ and our definition of $y$, $g\in \Delta_1$. Since $gF_1\cap X=\emptyset$, we have $g\in \Delta_1^*$. The proof of the inductive step follows routinely from the definition of the $\Gamma_j$ membership test.

We now claim that for any $j\geq 1$ and $g\in G$, $g\in \Gamma_j^*$
iff $g$ satisfies the $\Gamma_j$ membership test in $y$. In other
words, the $\Gamma_j^*$ membership test on $y$ takes exactly the same
form as the $\Gamma_j$ membership test on $x_k$ for $k \geq j$. We first verify this claim for $j=1$.
Thus we are to show that $g\in \Gamma_1^*$ iff $g$ satisfies the $\Gamma_1$ membership test in $y$.
The nontrivial direction is to show that if $g$ satisfies the $\Gamma_1$ membership test, 
then $g\in \Gamma_1^*$.  Since  $y(gf)=R(f)$ for all $f\in F_0$ we in particular have $y(g)=R(1_G)=1$.
If $g \notin X$, then since $y(g)=R(1_G)=1$ we must have that 
$g \in \gamma F_0$ for some unique $\gamma \in \Delta_0$. 
From the definition of $y$ we cannot have that 
$\gamma F_0\cap X \neq \emptyset$ as otherwise $y(g)=0$. So, $y \restriction \gamma F_0=
z \restriction \gamma F_0$. But then if follows from Lemma~\ref{lem:locallyorthogonal} that $g$ 
cannot satisfy the $\Gamma_1$ membership test in $y$. So, we may assume $g \in X$. 
We may therefore assume $g\in X_k$ for some $k\geq 1$. Fix $n>k$ such that
$g\in \pi_nK_{n,k}$. Since $y(g)=1$ we have that for some $\gamma \in \Gamma_0$ 
that $g \in \pi_n \gamma F_0$. By the $0$-saturation of $K_{n,k}$ we have that 
$\pi_n\gamma F_0 \subseteq \pi_n K_{n,k}$. We must have that $g=\pi_n \gamma$ and $\gamma \in \Gamma_1$
 as otherwise $g$ would not pass
the $\Gamma_1$ membership test in $y$ (c.f.\ Lemma~\ref{LEM LR EQUIV}). To see this, 
note that if $\gamma \notin \Gamma_1$, then since $y(g)=1$ we have 
$\gamma= \gamma' \delta$ for some $\gamma' \in \Gamma_1$
and $\delta \in D^1_0-\{ 1_G\}$. Also in this case we must have $g=\pi_n \gamma' \delta$ as $y(g)=1$. 
However, on the one hand, $g$ satisfies the $\Gamma_1$ membership test, which means that for all $f\in F_0$,
 $$ y(gf)=x_k(\pi_n^{-1}gf)=R(f). $$
On the other hand, for any $\gamma' \in \Gamma_1$ and $\delta\in D^1_0-\{1_G\}=D^1_0-\{\gamma_1\}$, 
the construction of the $\Gamma_1$ membership test using Theorem~\ref{FM} gives that
 $$ |\{f\in F_0\,:\, x_k(\gamma' \delta f)=R(1_G)=1\}|\leq 1. $$
Since $R$ is nontrivial, it follows that $\pi_n^{-1}g\not\in \gamma'(D^1_0-\{1_G\})$. 
Thus we must have $\pi_n^{-1}g\in \gamma F_0$ for $\gamma \in \Gamma_1$. Since 
$\pi_n^{-1}g \in K_{n,k}$, $\gamma \in \Gamma_1^*$. As $\pi_n \gamma F_0 \subseteq \pi_n K_{n,k}$
by $0$-stauration, we have that $y \res \pi_n \gamma F_0= x_k \res \pi_n \Gamma F_0$. 
Since $R$ is locally recognizable, we have $g=\pi_n\gamma\in
\Gamma_1^*$, as required.

Suppose next that $j>1$ and $g$ passes the $\Gamma_j$ membership test in $y$.
Since $g$ also passes the $\Gamma_1$ membership test we have that $g\in \Gamma_1^*$, say 
$g=\pi_n \gamma$ where $\gamma \in \Gamma_1 \cap K_{n,k}$. Suppose first that $j \leq k$ 
and assume inductively that $\gamma \in \Gamma_{j-1}$. 
As $g$ satisfies the $\Gamma_j$ membership test in $y$ we have that $\pi_n \gamma \gamma_j$ 
satisfies the $\Gamma_{j-1}$
membership test in $y$. From the $j=1$ case we have that $\pi_n \gamma\gamma_j \in X$
and hence $\pi_n \gamma \gamma_j \in \pi_n K_{n,k}$ as $\gamma_j \in F_k$ and 
$\pi_n K_n F_k$ is disjoint from all $K_m$ for $m \neq n$. Since $j \leq k$ we have by saturation that 
$\pi_n \gamma\gamma_j F_j \subseteq \pi_n K_{n,k}$. It follows that 
$\gamma$ passes the $\Gamma_j$ membership test in $x_k$. Thus, $\gamma \in \Gamma_j$ and so 
$g=\pi_n \gamma \in \pi_n \Gamma_j \cap \pi_n K_{n,k} \subseteq \Gamma^*_j$. Suppose 
now $j>k$, and  $g$ passes the $k+1$ membership test in $y$. 
As above we get that $g=\pi_n \gamma$ where $\gamma \in \Gamma_k$ and $\gamma \gamma_k \in \Gamma_k$
where $\pi_n \gamma \gamma_k F_k \subseteq \pi_n K_{n,k}$. Since $g$ passes the $k+1$ menbership
test in $y$ we see that $y(\pi_n \gamma \gamma_k a_{k})= y(\pi_n \gamma \gamma_k b_{k})=1$ and so 
$x_k(\gamma \gamma_k a_{k})=x_k(\gamma \gamma_k b_{k})=1$. This is a contradiction as 
for any $\gamma'\in \Gamma_k$ we have that $x_k(\gamma' a_{k})$ and $x_k(\gamma' b_{k})$ 
are not both $1$ from the definition of
$x_k$ (we may assume $x_k$ has this property without loss of generality).
So, there is no $g \in \pi_nK_{n,k}$ which passes the $\gamma_{k+1}$ membership test. 
This establishes the claim.

We fix $g\not\in SF_i^{-1}F_iF_i^{-1}\{1_G,s^{-1}\}$. Consider the following cases below.

Case 1a: $gF_iF_i^{-1}F_i\cap X\neq\emptyset$. Thus there is $\delta_0\in gF_iF_i^{-1}$ such that $\delta_0 F_i\cap X\neq\emptyset$. Then for some $k\geq 1$ and $n>k$, $\delta_0 F_i\cap \pi_nK_{n,k}\neq\emptyset$. Fix such $\delta_0$, $k\geq 1$ and $n>k$.

Since $g\not\in SF_i^{-1}F_iF_i^{-1}$ but $g\in \delta_0 F_iF_i^{-1}$, we have $\delta_0\not\in SF_i^{-1}$ and $\delta_0 F_i\cap S=\emptyset$. Thus $n>n_0$, where $n_0$ is defined in the definition of $S$. Recall that
$k$ is the least integer with $1\leq k<n$ such that $\alpha(k,n)=1$. Since $\alpha(j,n)=0$ for all $1\leq j\leq i$, we know that $i<k$.

Let $C=M_nM_k^3\dots M_i^3\dots M_0^3$. Recall that $K_{n,k}=\mathrm{sat}_k(C)$.
It follows that there is $1\leq j\leq k$ and $\delta_1\in \pi_n\mbox{sat}_0(C)\cap\Gamma_j^*$ such that
$\delta_0 F_i\cap \delta_1 F_j\neq\emptyset$.  If $j\geq i$ then by (4) we may assume $j=i$, and thus we
have found $\delta_1\in \Gamma_i^*$ with $\delta_0F_i \cap\delta_1 F_i\neq\emptyset$. Noting that
$g^{-1}\delta_0, \delta_0^{-1}\delta_1\in F_iF_i^{-1}$, we have that $g^{-1}\delta_1\in M_i^4\subseteq M_i^{11}$. Alternatively, assume $j<i$. Then $\delta_0^{-1}\delta_1\in F_iF_j^{-1}$. Since $\pi_n^{-1}\delta_1\in \mbox{sat}_0(C)\subseteq M_nM_k^3\dots M_0^3F_0F_0^{-1}$, then there is $\delta_2\in \pi_nM_nM_k^3\dots M_{i+1}^3$ such that $\delta_2^{-1}\delta_1\in M_i^3M_{i-1}^3\dots M_0^3F_0F_0^{-1}$. By (2)
$\delta_2^{-1}\delta_1\in M_i^3M_iF_0F_0^{-1}=M_i^4F_0F_0^{-1}$. Since $M_i^{-1}=M_i$ we have $\delta_1^{-1}\delta_2\in F_0F_0^{-1}M_i^4$. Now $\pi_n^{-1}\delta_2 F_i\subseteq M_n M_k^3\dots M_{i+1}^3F_i\subseteq K_{n,k}$, and thus there is $\delta_3\in \Gamma_i$ such that
$\pi_n^{-1}\delta_2F_i\cap \delta_3 F_i\neq\emptyset$ by the maximal disjointness of $\Gamma_i$-translates of $F_i$ by Theorem~\ref{thm:nonflecc1} (iii). Since $\delta_2^{-1}\pi_n\delta_3\in F_iF_i^{-1}$, we have $$\delta_3=\pi_n^{-1}\delta_2(\delta_2^{-1}\pi_n\delta_3)\in M_n M_k^3\dots M_{i+1}^3F_iF_i^{-1}\subseteq M_n M_k^3\dots M_{i+1}^3M_i^2\subseteq C, $$
and $\delta_3F_i\subseteq M_n M_k^3\dots M_i^3\subseteq C$. This shows that $\pi_n\delta_3\in \Gamma_i^*$. Thus we have found $\pi_n\delta_3\in \Gamma_i^*$ such that
$$\begin{array}{rcl}
g^{-1}\pi_n\delta_3&=&(g^{-1}\delta_0)(\delta_0^{-1}\delta_1)(\delta_1^{-1}\delta_2)(\delta_2^{-1}\pi_n\delta_3)\\
&\in&
(F_iF_i^{-1})(F_iF_j^{-1})(F_0F_0^{-1}M_i^4)(F_iF_i^{-1})\subseteq M_i^{11}.
\end{array} $$
In either case of $j<i$ or $j\geq i$, we have found $\gamma\in\Gamma_i^*$ such that $g^{-1}\gamma\in M_i^{11}$.

Let $t_0=g^{-1}\gamma$. Then $gt_0$ satisfies the $\Gamma_i$ membership test. If $gst_0$ does not satisfy the $\Gamma_i$ membership test, then there is $t_1\in F_i$ such that $y(gt_0t_1)\neq y(gst_0t_1)$ by Theorem~\ref{thm:nonflecc2} (ii). Since $t_0t_1\in M_i^{11}F_i=T$, we are done. Otherwise, assume $gst_0$ satisfies the $\Gamma_i$ membership test. We have
$$ (gt_0)^{-1}(gst_0)=t_0^{-1}st_0\in M_i^{11}H_iM_i^{11}=M_i^{23}. $$
By (7) $gst_0\in \Gamma_i^*\cap \pi_nK_n$. By Theorem~\ref{thm:nonflecc2} (iv) there is $t_1\in F_i$ such that $y(gt_0t_1)\neq y(gst_0t_1)$. Again $t_0t_1\in M_i^{11}F_i$. 

Case 1b: $gsF_iF_i^{-1}F_i\cap X\neq\emptyset$. The argument is similar to the above argument in Case 1a, with $gs$ now playing the role of $g$ in that argument.

Case 2: Otherwise, $gF_iF_i^{-1}F_i\cap X=\emptyset$ and $gsF_iF_i^{-1}F_i\cap X=\emptyset$. In particular, for every $\delta\in \Delta_i$ with $gF_i\cap \delta F_i\neq\emptyset$, we have that $\delta F_i\cap X=\emptyset$. Thus for every $\delta\in \Delta_i$ with $g F_i \cap \delta F_i\neq\emptyset$, $\delta\in \Delta_i^*$. Similarly for every $\delta\in \Delta_i$ with $gs F_i\cap \delta F_i\neq\emptyset$, we also have $\delta\in \Delta_i^*$. As usual there is $t_0\in F_iF_i^{-1}$ such that $gt_0\in \Delta_i^*$. If $gst_0\in \Delta_i^*$ then we may find $t_1\in F_i$ so that $y(gt_0t_1)\neq y(gst_0t_1)$, since $(gt_0)^{-1}(gst_0)\in M_i^{23}$, and we are done. Assume $gst_0\not\in \Delta_i^*$. Since $gst_0F_i\cap X=\emptyset$, we have that $gst_0$ fails the $\Delta_i$ membership test on $y$. Now that $gt_0$ does satisfy the $\Delta_i$ membership test, we routinely find $t_1\in F_i$ with $y(gt_0t_1)\neq y(gst_0t_1)$.

This shows that $y=f(\alpha)$ is a $2$-coloring, and our proof is complete.
\end{proof}

We also draw the following corollaries from the proof.

\begin{theorem} For any countable nonflecc group $G$ the set of all strong $2$-colorings on $G$ is
${\bf\Pi}^0_3$-complete.
\end{theorem}

\begin{proof}
It suffices to note that, in the above proof if $\alpha\in P$ then $f(\alpha)$ is in fact a strong $2$-coloring on $G$. This is because $y(w_n)=z(w_n)$ and $y(\sigma_nw_n)=z(\sigma_nw_n)$ for all $n\geq 1$ by (5) and the definition of $y$. By (6), $y(w_n)\neq y(\sigma_nw_n)$ for all $n\geq 1$. Thus for each $s\neq 1_G$ there are infinitely many $t\in G$ such that $y(t)\neq y(st)$.
\end{proof}

The following corollary summarizes our findings.

\begin{cor} Let $G$ be a countable group. Then the following hold:
\begin{enumerate}
\item[(1)] If $G$ is finite, then the set of all $2$-colorings on $G$ is closed.
\item[(2)] If $G$ is an infinite flecc group, then the set of all $2$-colorings on $G$ is ${\bf\Sigma}^0_2$-complete;
\item[(3)] If $G$ is not flecc, then the set of all $2$-colorings on $G$ is ${\bf\Pi}^0_3$-complete.
\end{enumerate}
\end{cor}

\chapter{The Complexity of the Topological Conjugacy Relation} \label{CHAP ISO}

In this chapter we study the complexity of the topological conjugacy relation among subflows of $2^G$. We remind the reader the definition of topological conjugacy.

\begin{definition}
Let $G$ be a countable group and let $S_1, S_2 \subseteq 2^G$ be subflows. $S_1$ is \emph{topologically conjugate} to $S_2$ or is a \emph{topological conjugate} of $S_2$ if there is a homeomorphism $\phi: S_1 \rightarrow S_2$ satisfying $\phi(g \cdot x) = g \cdot \phi(x)$ for all $x \in S_1$ and $g \in G$. Such a function $\phi$ is called a \emph{conjugacy} between $S_1$ and $S_2$. The property of being topologically conjugate induces an equivalence relation on the set of all subflows of $2^G$. We call this equivalence relation the \emph{topological conjugacy relation}.
\end{definition} 

The purpose of this chapter is to study the complexity of the topological conjugacy relation, meaning, in some sense, how difficult it is to determine when two subflows are topologically conjugate. The precise mathematical way of discussing the complexity of equivalence relations is via the theory of Borel equivalence relations. In the first section, we present a basic introduction to the aspects of the theory of countable Borel equivalence relations which will be needed in this chapter. The second section consists mostly of preparatory work and basic lemmas. In the third section, we show that for every countably infinite group the equivalence relation $E_0$, which we will define in section one, is a lower bound to the complexity of the topological conjugacy relation restricted to free minimal subflows. In the fourth section, we give a complete classification of the complexity of both the topological conjugacy relation and the restriction of the topological conjugacy relation to free subflows.

\section{Introduction to countable Borel equivalence relations} \label{INTRO RELATNS}

In this section we review common notation and terminology and basic facts related to the theory of countable Borel equivalence relations. Some references for this material include \cite{JKL} and \cite{Gao}. Throughout this chapter we will also work with descriptive set theory, and we refer the reader to \cite{KechrisBook} for any missing details.

\index{Polish space}\index{Borel sets}\index{Borel!sets}\index{Borel equivalence relation}\index{Borel!equivalence relation}\index{Borel function}\index{Borel!function} Let $X$ be a \emph{Polish space}, that is, a topological space which is separable and which admits a complete metric compatible with its topology. Recall that the \emph{Borel sets} of $X$ are the members of the $\sigma$-algebra generated by the open sets. Informally, the Borel subsets of $X$ are considered to be the \emph{definable} subsets of $X$. A \emph{Borel equivalence relation} on $X$ is an equivalence relation on $X$ which is a Borel subset of $X \times X$, where $X \times X$ has the product topology. Given two Polish spaces $X$ and $Y$, a function $f : X \rightarrow Y$ is \emph{Borel} if the pre-image of every Borel set in $Y$ is Borel in $X$. As with Borel sets, Borel functions are viewed informally as being \emph{definable}.

\index{Borel reducible}\index{Borel!reducible}\index{Borel embeddable}\index{Borel!embeddable}\index{continuously reducible}\index{continuously embeddable}\index{$\leq_B$}\index{$\leq_c$}\index{$\sqsubseteq_B$}\index{$\sqsubseteq_c$} We compare Borel equivalence relations and discuss their complexity relative to one another via the notion of Borel reducibility. If $E$ is a Borel equivalence relation on $X$ and $F$ is a Borel equivalence relation on $Y$, then $E$ is \emph{Borel reducible} to $F$, written $E \leq_B F$, if there is a Borel function $f: X \rightarrow Y$ such that $f(x_1) \ F \ f(x_2) \Leftrightarrow x_1 \ E \ x_2$ for all $x_1, x_2 \in X$. Such a function $f$ is called a \emph{reduction}, and if $f$ is injective then we say $E$ is \emph{Borel embeddable} into $F$, written $E \sqsubseteq_B F$. Furthermore, $E$ is \emph{continuously reducible} (or \emph{continuously embeddable}) to $F$, written $E \leq_c F$ (respectively $E \sqsubseteq_c F$), if the reduction (respectively embedding) $f$ is continuous.  

Intuitively, if $E$ is Borel reducible to $F$, then $E$ is considered to be no more complicated than $F$, and $F$ is considered to be at least as complicated as $E$. To illustrate, suppose $E$ is Borel reducible to $F$ and $f : X \rightarrow Y$ is a Borel reduction. If there were a definable (Borel) way to determine when two elements of $Y$ are $F$-equivalent, then by using the Borel function $f$ there would be a definable (Borel) way to determine when two elements of $X$ are $E$-equivalent. The theory of Borel equivalence relations therefore allows us to compare the relative complexity of classification problems.

\index{finite equivalence relation}\index{countable equivalence relation}\index{universal countable Borel equivalence relation}\index{equivalence relation!finite}\index{equivalence relation!Borel}\index{equivalence relation!countable}\index{Borel!universal countable equivalence relation}\index{equivalence relation!universal countable Borel}An equivalence relation $E$ is \emph{finite} if every $E$-equivalence class is finite, and $E$ is \emph{countable} if every $E$-equivalence class is countable. A \emph{universal countable} Borel equivalence relation $F$ is a countable Borel equivalence relation with the property that if $E$ is any other countable Borel equivalence relation then $E$ is Borel reducible to $F$. Thus, the universal countable Borel equivalence relations are the most complicated among all countable Borel equivalence relations. Let $\F$ be the nonabelian free group on two generators. Then the equivalence relation $E_\infty$ on $2^\F$ defined by $x \ E_\infty \ y \Leftrightarrow \exists f \in \F \ f \cdot x = y$ is a universal countable Borel equivalence relation.

\index{smooth equivalence relation}\index{equivalence relation!smooth} On the other hand, one of the least complicated classes of Borel equivalence relations are the smooth equivalence relations. A Borel equivalence relation $E$ is \emph{smooth} if there is a Polish space $Y$ and a Borel $f : X \rightarrow Y$ such that $x_1 \ E \ x_2 \Leftrightarrow f(x_1) = f(x_2)$ for all $x_1, x_2 \in X$. This condition is equivalent to $E$ being Borel reducible to the equality equivalence relation on $Y$. Smooth equivalence relations are considered to be the simplest Borel equivalence relations because there is a definable way to determine when two elements are equivalent. The universal countable Borel equivalence relations are not smooth. Note that if $E$ is not smooth and is Borel reducible to $F$ then $F$ is not smooth.

\index{hyperfinite}\index{equivalence relation!hyperfinite} A Borel equivalence relation $E$ is \emph{hyperfinite} if $E = \bigcup_{n \in \N} E_n$, where $(E_n)_{n \in \N}$ is an increasing sequence of finite Borel equivalence relations. The canonical example of a hyperfinite equivalence relation is $E_0$, which is the equivalence relation on $2^\N$ defined by
$$x \ E_0 \ y \Longleftrightarrow \exists m \ \forall n \geq m \ x(n) = y(n).$$
$E_0$ is not smooth.

\section{Basic properties of topological conjugacy} \label{SECT BASIC TC}

The purpose of this section is to develop some of the basic facts regarding the topological conjugacy relation which will be needed in later sections. Of particular importance is to prove that the topological conjugacy relation is a countable Borel equivalence relation.

When discussing the topological conjugacy relation, we will employ the following notation: \index{$\Su(G)$}\index{$\SuM(G)$}\index{$\SuF(G)$}\index{$\SuMF(G)$}\index{$\TC(G)$}\index{$\TCM(G)$}\index{$\TCF(G)$}\index{$\TCMF(G)$}\index{$\TCP(G)$}
$$\Su(G) = \{A \subseteq 2^G \: A \text{ is a subflow of } 2^G\};$$
$$\SuM(G) = \{A \in \Su(G) \: A \text{ is minimal}\};$$
$$\SuF(G) = \{A \in \Su(G) \: A \text{ is free}\};$$
$$\SuMF(G) = \SuM(G) \cap \SuF(G);$$
$$\TC(G) = \text{the topological conjugacy relation on } \Su(G);$$
$$\TCM(G) = \TC(G)\res (\SuM(G) \times \SuM(G));$$
$$\TCF(G) = \TC(G)\res (\SuF(G) \times \SuF(G));$$
$$\TCMF(G) = \TCM(G) \cap \TCF(G);$$
$$\TCP(G) = \{(x, y) \in 2^G \times 2^G \: \exists \text{ conjugacy } \phi: \overline{[x]} \rightarrow \overline{[y]} \text{ with } \phi(x) = y\};$$

We easily have the following.

\begin{lem}
For countable groups $G$, the equivalence relations $\TC(G)$, $\TCM(G)$, $\TCF(G)$, $\TCMF(G)$, and $\TCP(G)$ are all countable equivalence relations.
\end{lem}

\begin{proof}
Let $A \in \Su(G)$. Then for every $B \in \Su(G)$ topologically conjugate to $A$ there is a conjugacy $\phi_B: A \rightarrow B$ which is induced by a block code $\hat{\phi}_B$ (Theorem \ref{THM BC}). Clearly if $\hat{\phi}_B = \hat{\phi}_C$ then $\phi_B = \phi_C$ and $B = C$. Since there are only countably many block codes, the $\TC(G)$-equivalence class of $A$ must be countable. Similar arguments work for the other equivalence relations.
\end{proof}

\index{Vietoris topology}\index{Hausdorff metric} Let $X$ be a Polish space, and let $K(X) = \{K \subseteq X \: K \text{ compact}\}$. The \emph{Vietoris topology} on $K(X)$ is the topology generated by subbasic open sets of the form
$$\{K \in K(X) \: K \subseteq U\} \text{ and } \{K \in K(X) \: K \cap U \neq \varnothing\}$$
where $U$ varies over open subsets of $X$. It is well known that $K(X)$ with the Vietoris topology is a Polish space (see for example \cite{KechrisBook}). In fact, a compatible complete metric on $K(X)$ is the \emph{Hausdorff metric}. The Hausdorff metric, $d_H$, is defined by
$$d_H(A, B) = \max \left( \sup_{a \in A} \inf_{b \in G} d(a, b), \ \sup_{b \in B} \inf_{a \in A} d(a, b) \right),$$
where $A, B \in K(X)$ and $d$ is a complete metric on $X$ compatible with its topology.

\begin{lem}
For every countable group $G$, $\Su(G)$ and $\SuF(G)$ are Polish spaces with the subspace topology inherited from $K(2^G)$.
\end{lem}

\begin{proof}
Since $G_\delta$ subsets of Polish spaces are Polish (\cite{KechrisBook}), it will suffice to show that $\Su(G)$ and $\SuF(G)$ are $G_\delta$ in $K(2^G)$. Let $\{U_n \: n \in \N\}$ be a countable base for the topology on $2^G$ consisting of clopen sets. For $n \in \N$ and $g \in G$ define
$$V_{n,g} = \{K \in K(2^G) \: (K \cap U_n \neq \varnothing \wedge K \cap g \cdot U_n \neq \varnothing) \vee K \subseteq (2^G - (U_n \cup g \cdot U_n))\}$$
Notice that $V_{n,g}$ is open in $K(2^G)$ since $U_n$ is clopen and $G$ acts on $2^G$ by homeomorphisms. For $A \in K(2^G)$ we have
$$A \in \Su(G) \Longleftrightarrow \forall g \in G \ g \cdot A = A$$
$$\Longleftrightarrow \forall n \in \N \ \forall g \in G \ (A \cap U_n \neq \varnothing \Leftrightarrow A \cap g \cdot U_n \neq \varnothing)$$
$$\Longleftrightarrow A \in \bigcap_{n \in \N} \ \bigcap_{g \in G} V_{n,g}$$
So $\Su(G)$ is $G_\delta$ in $K(2^G)$ and hence Polish.

It now suffices to show $\SuF(G)$ is $G_{\delta}$ in $\Su(G)$. A modification of the proof of Lemma \ref{lem:basiccoloringlemma} shows that for $A \in \Su(G)$
$$A \in \SuF(G) \Longleftrightarrow \forall s \in G - \{1_G\} \ \exists \text{ finite } T \subseteq G \ \forall x \in A \ \exists t \in T \ x(st) \neq x(t)$$
$$\Longleftrightarrow \forall s \in G - \{1_G\} \ \exists \text{ finite } T \subseteq G \ A \subseteq \{y \in 2^G \: \exists t \in T \ y(st) \neq y(t)\}$$
$$\Longleftrightarrow A \in \bigcap_{s \in G - \{1_G\}} \ \bigcup_{\text{finite } T \subseteq G} \{B \in \Su(G) \: B \subseteq \{y \in 2^G \: \exists t \in T \ y(st) \neq y(t)\}\}.$$
The set $\{y \in 2^G \: \exists t \in T \ y(st) \neq y(t)\}$ is open, and therefore $\SuF(G)$ is $G_{\delta}$ in $\Su(G)$.
\end{proof}

\begin{lem} \label{LEM SUBMAP}
For countable groups $G$, the map $x \in 2^G \mapsto \overline{[x]}$ is Borel.
\end{lem}

\begin{proof}
For $x \in 2^G$, define $f(x) = \overline{[x]}$. We check that the inverse images of the subbasic open sets in $K(2^G)$ are Borel. If $U$ is open in $2^G$, then $f(x) \cap U \neq \varnothing$ if and only if $[x] \cap U \neq \varnothing$. Therefore
$$f^{-1}(\{A \in \Su(G) \: A \cap U \neq \varnothing\}) = \bigcup_{g \in G} g \cdot U$$
which is Borel (in fact open). For an open $U \subseteq 2^G$, define $U_n = \{y \in U \: d(y, 2^G - U) \geq 1/n\}$. Then each $U_n$ is closed and $f(x) \subseteq U$ if and only if $\overline{[x]} \subseteq U_n$ for some $n$ (by compactness). This is equivalent to the condition that $x \in \bigcap_{g \in G} g \cdot U_n$ for some $n$. Thus
$$f^{-1}(\{A \in \Su(G) \: A \subseteq U\}) = \bigcup_{n \in \N} \bigcap_{g \in G} g \cdot U_n.$$
We conclude $f$ is Borel.
\end{proof}

For the next lemma we need to review some terminology. A measurable space $(X, S)$ (a set $X$ and a $\sigma$-algebra $S$ on $X$) is said to be a \emph{standard Borel space} if there is a Polish topology on $X$ for which $S$ coincides with the collection of Borel sets in this topology. Thus standard Borel spaces are essentially Polish spaces, but the topology is not emphasized. If $(X, S)$ is a measurable space and $Y \subseteq X$, then the \emph{relative} $\sigma$-algebra on $Y$ inherited from $X$ is the $\sigma$-algebra consisting of sets of the form $Y \cap A$, where $A$ ranges over all elements of $S$. We will view every Polish space as a measurable space with the $\sigma$-algebra of Borel sets. A well known result is that if $X$ is a Polish space and $Y \subseteq X$ is Borel, then $Y$ is a standard Borel space with the relative $\sigma$-algebra inherited from $X$ (see \cite{KechrisBook}).

\begin{lem}
For every countable group $G$, $\SuM(G)$ and $\SuMF(G)$ are standard Borel spaces with the relative $\sigma$-algebra inherited from $\Su(G)$.
\end{lem}

\begin{proof}
It suffices to show that $\SuM(G)$ is a Borel subset of $\Su(G)$. We will need a Borel function $f: \Su(G) \rightarrow 2^G$ for which $f(A) \in A$ for every $A \in \Su(G)$. Such a function is called a \emph{Borel selector}, and by standard results in descriptive set theory they are known to exist within this context (see \cite{KechrisBook}). For clarity and to minimize pre-requisites, we construct a Borel selector $f$ explicitly. Fix an enumeration $g_0, g_1, \ldots$ of $G$, and define a partial order, $\prec$, on $2^G$ by
$$x \prec y \Longleftrightarrow (x = y) \vee (\exists n \in \N \ \forall k < n \ x(g_k) = y(g_k) \wedge x(g_n) < y(g_n)).$$
By compactness, if $A \in \Su(G)$ then $A$ contains a $\prec$-least element. Define $f: \Su(G) \rightarrow 2^G$ by letting $f(A)$ be the $\prec$-least element of $A$. One can show that $f$ is continuous and hence Borel.

By Lemma \ref{lem:minimallemma} we have that for $A \in \Su(G)$
$$A \in \SuM(G) \Longleftrightarrow A = \overline{[f(A)]} \wedge (\forall \text{ finite } H \subseteq G \ \exists \text{ finite } T \subseteq G$$
$$\forall g \in G \ \exists t \in T \ \forall h \in H \ f(A)(gth) = f(A)(h))$$
$$\Longleftrightarrow A = \overline{[f(A)]} \wedge f(A) \in \bigcap_{\text{finite } H \subseteq G} \ \bigcup_{\text{finite } T \subseteq G} \ \bigcap_{g \in G} \ \bigcup_{t \in T} \ \bigcap_{h \in H} \{y \in 2^G \: y(gth) = y(h)\}.$$
Note that this last set on the right is Borel. Finally, if we define $g: \Su(G) \rightarrow \Su(G) \times \Su(G)$ by $g(A) = (A, \overline{[f(A)]})$, then $g$ is Borel and
$$A = \overline{[f(A)]} \Longleftrightarrow A \in g^{-1}(\{(B, B) \: B \in \Su(G)\}).$$
We conclude $\SuM(G)$ is a Borel subset of $\Su(G)$, in fact $\SuM(G)$ is $\mathbf{\Pi}^0_3$ in $\Su(G)$. Clearly $\SuMF(G) = \SuM(G) \cap \SuF(G)$ is a Borel ($\mathbf{\Pi}^0_3$) subset of $\Su(G)$ as well.
\end{proof}

We now prove that all of the equivalence relations we are working with are countable Borel equivalence relations.

\begin{prop} \label{TC BOREL EQREL}
For countable groups $G$, the equivalence relations $\TC(G)$, $\TCM(G)$, $\TCF(G)$, $\TCMF(G)$, and $\TCP(G)$ are all countable Borel equivalence relations.
\end{prop}

\begin{proof}
We saw at the beginning of this section that they are all countable equivalence relations. So we only need to check that they are all Borel. Since $\SuM(G)$, $\SuF(G)$, and $\SuMF(G)$ are Borel subsets of $\Su(G)$, we only need to check that $\TC(G)$ and $\TCP(G)$ are Borel.

For a block code $\hat{f}$, we will let $f : 2^G \rightarrow 2^G$ be the function induced by $\hat{f}$. Let $\{U_n \: n \in \N\}$ be a countable base for the topology on $2^G$. For $A, B \in \Su(G)$ we have
$$(A, B) \in \TC(G) \Longleftrightarrow \exists \text{ block codes } \hat{f_1}, \hat{f_2}$$
$$f_1(A) = B \wedge f_2(B) = A \wedge (f_2 \circ f_1)\res A = \textrm{id}_A \wedge (f_1 \circ f_2)\res B = \textrm{id}_B.$$
For a fixed block code $\hat{f_1}$ the set $\{(A, B) \in \Su(G)^2 \: f_1 (A) = B\}$ is Borel (in fact $G_\delta$) since
$$f_1(A) = B \Longleftrightarrow (\forall n \in \N \ B \cap U_n \neq \varnothing \Leftrightarrow A \cap f_1^{-1}(U_n) \neq \varnothing)$$
$$\Longleftrightarrow (A, B) \in \bigcap_{n \in \N} (\{(K_1, K_2) \in \Su(G)^2 \: K_1 \cap f_1^{-1}(U_n) \neq \varnothing \wedge K_2 \cap U_n \neq \varnothing\}$$
$$\cup \{(K_1, K_2) \in \Su(G)^2 \: K_1 \subseteq 2^G - f_1^{-1}(U_n) \wedge K_2 \subseteq 2^G - U_n\}).$$
Also, for fixed block codes $\hat{f_1}$ and $\hat{f_2}$, the set of $A \in \Su(G)$ with $(f_2 \circ f_1)\res A = \textrm{id}_A$ is Borel (in fact closed) since $\{x \in 2^G \: f_2 \circ f_1(x) \neq x\}$ is open and
$$(f_2 \circ f_1)\res A = \textrm{id}_A \Longleftrightarrow A \cap \{x \in 2^G \: f_2 \circ f_1(x) \neq x\} = \varnothing.$$
So we conclude that $\TC(G)$ is a Borel equivalence relation (in fact it is $\mathbf{\Sigma}^0_3$).

Now we consider $\TCP(G)$. Note that if $\hat{f_1}$ is a block code and $f_1(x) = y$, then $f_1(\overline{[x]}) = \overline{[y]}$ since $f_1$ is continuous and $\overline{[x]}$ is compact. Also, if $\hat{f_2}$ is another block code, then $\{z \in 2^G \: f_2 \circ f_1 (z) = z\}$ is closed and $G$-invariant. So $(f_2 \circ f_1)\res \overline{[x]} = \textrm{id}_{\overline{[x]}}$ if and only if $f_2 \circ f_1 (x) = x$. Therefore
$$(x, y) \in \TCP(G) \Longleftrightarrow \exists \text{ block codes } \hat{f_1}, \hat{f_2} \ f_1(x) = y \wedge f_2(y) = x$$
$$\Longleftrightarrow (x, y) \in \bigcup_{\text{block codes }\hat{f_1}, \hat{f_2}} (\{(z, f_1(z)) \: z \in 2^G\} \cap \{(f_2(z), z) \: z \in 2^G\}).$$
We conclude that $\TCP(G)$ is a Borel equivalence relation (in fact it is $F_\sigma$).
\end{proof}

\begin{cor} \label{COR ISOTEST}
For countable groups $G$ and $x, y \in 2^G$, $x \ \TCP(G) \ y$ if and only if there is a finite set $H$ such that
$$\forall g_1, g_2 \in G \ (\forall h \in H \ x(g_1 h) = x(g_2 h) \Longrightarrow y(g_1) = y(g_2)), \text{ and}$$
$$\forall g_1, g_2 \in G \ (\forall h \in H \ y(g_1 h) = y(g_2 h) \Longrightarrow x(g_1) = x(g_2)).$$
\end{cor}

\begin{proof}
We showed in the proof of the previous proposition that $x$ and $y$ are $\TCP(G)$-equivalent if and only if there are block codes $\hat{f_1}$ and $\hat{f_2}$ such that $f_1(x) = y$ and $f_2(y) = x$, where $f_1$ and $f_2$ are the functions induced by $\hat{f_1}$ and $\hat{f_2}$ respectively. The existence of such block codes is equivalent to the condition in the statement of this corollary.
\end{proof}

\begin{lem} \label{LEM EQ}
Let $G$ be a countable group. Then
\begin{enumerate}
\item[\rm (i)] $\TCMF(G) \sqsubseteq \TCM(G)$;
\item[\rm (ii)] $\TCMF(G) \sqsubseteq \TCF(G)$;
\item[\rm (iii)] $\TCM(G) \sqsubseteq \TC(G)$;
\item[\rm (iv)] $\TCF(G) \sqsubseteq_c \TC(G)$;
\end{enumerate}
\end{lem}

\begin{proof}
Use the inclusion map for each embedding. The first three embeddings are only Borel because we never formally fixed Polish topologies on $\SuM(G)$ and $\SuMF(G)$.
\end{proof}

In section four, after presenting a complete classification of $\TC(G)$ and $\TCF(G)$ it will be a corollary that $\TCF(G)$ and $\TC(G)$ are Borel bi-reducible. In the remainder of this section, we present a relationship between the topological conjugacy relations on $2^H$, $2^K$, and $2^{H \times K}$ for countable groups $H$ and $K$.

Let $2^H \times 2^K$ have the product topology, and let $H \times K$ act on $2^H \times 2^K$ in the obvious way. For $(x,y) \in 2^H \times 2^K$ we let $[(x,y)]$ denote the orbit of $(x,y)$. We call a closed subset of $2^H \times 2^K$ which is invariant under the action of $H \times K$ a \emph{subflow}. A subflow of $2^H \times 2^K$ is \emph{free} if every point in the subflow has trivial stabilizer, and it is \emph{minimal} if every orbit in the subflow is dense. Two subflows of $2^H \times 2^K$ are \emph{topologically conjugate} if there is a homeomorphism between them which commutes with the action of $H \times K$. Such a homeomorphism is called a \emph{conjugacy}. Finally, $\TC(2^H \times 2^K)$, $\TCF(2^H \times 2^K)$, $\TCMF(2^H \times 2^K)$, $\TCM(2^H \times 2^K)$, and $\TCP(2^H \times 2^K)$ denote the obvious equivalence relations.

\begin{lem}
Let $H$ and $K$ be countable groups, let $A_1, A_2 \in \Su(H)$ and $B_1, B_2 \in \Su(K)$, and let $x_1, x_2 \in 2^H$ and $y_1, y_2 \in 2^K$. Then
\begin{enumerate}
\item[\rm (i)] $A_1 \times B_1$ is free if and only if $A_1$ and $B_1$ are free;
\item[\rm (ii)] $A_1 \times B_1$ is minimal if and only if $A_1$ and $B_1$ are minimal;
\item[\rm (iii)] if $A_1$, $A_2$, $B_1$, and $B_2$ are minimal then $A_1 \times B_1$ is topologically conjugate to $A_2 \times B_2$ if and only if $A_1 \ \TC(H) \ A_2$ and $B_1 \ \TC(K) \ B_2$;
\item[\rm (iv)] $(x_1,y_1) \ \TCP(2^H \times 2^K) \ (x_2,y_2)$ if and only if both $x_1 \ \TCP(H) \ x_2$ and $y_1 \ \TCP(K) \ y_2$.
\end{enumerate}
\end{lem}

\begin{proof}
The proofs of clauses (i) and (ii) are trivial. For (iii) it is clear that if $\phi: A_1 \rightarrow A_2$ and $\psi: B_1 \rightarrow B_2$ are conjugacies then $\phi \times \psi$ is a conjugacy between $A_1 \times B_1$ and $A_2 \times B_2$. Now suppose that $A_1$, $A_2$, $B_1$, and $B_2$ are minimal and that $\theta$ is a conjugacy between $A_1 \times B_1$ and $A_2 \times B_2$. Fix $y_1 \in B_1$ and $x_1 \in A_1$, let $p_H: 2^H \times 2^K \rightarrow 2^H$ and $p_K: 2^H \times 2^K \rightarrow 2^K$ be the projection maps, and set $y_2 = p_K(\theta(x_1,y_1))$. Then for every $h \in H$ we have
$$y_2 = p_K(\theta(x_1,y_1)) = p_K(h \cdot \theta(x_1, y_1)) = p_K(\theta(h \cdot x_1, y_1)).$$
Since $A_1$ is minimal, $A_1 = \overline{[x_1]}$ and therefore for every $x \in A_1$ we have $y_2 = p_K(\theta(x,y_1))$. The same reasoning shows that for every $x \in A_2$ we have $y_1 = p_K(\theta^{-1}(x, y_2))$. Define $\phi: A_1 \rightarrow A_2$ by $\phi(x) = p_H(\theta(x,y_1))$. This is clearly continuous and commutes with the action of $H$. It is injective and surjective since $\phi^{-1}(x) = p_H(\theta^{-1}(x,y_2))$. Thus $A_1 \ \TC(H) \ A_2$. A similar argument shows that $B_1 \ \TC(K) \ B_2$. The proof of (iv) is essentially the same.
\end{proof}

\begin{cor} \label{COR TC PROD}
For countable groups $H$ and $K$ we have
\begin{enumerate}
\item[\rm (i)] $\TCMF(H) \times \TCMF(K) \sqsubseteq_c \TCMF(2^H \times 2^K)$;
\item[\rm (ii)] $\TCM(H) \times \TCM(K) \sqsubseteq_c \TCM(2^H \times 2^K)$;
\item[\rm (iii)] $\TCP(H) \times \TCP(K) \sqsubseteq_c \TCP(2^H \times 2^K)$.
\end{enumerate}
\end{cor}

Now we want to relate topological conjugacy in $2^H \times 2^K$ to topological conjugacy in $2^{H \times K}$. The following lemma makes this easy. In the rest of this section we let $0$ denote the element of $2^H$ or $2^K$ which is identically zero.

\begin{lem}
Let $H$ and $K$ be countable groups. There exists a function $f : 2^H \times 2^K \rightarrow 2^{H \times K}$ with the following properties:
\begin{enumerate}
\item[\rm (i)] $f$ restricted to $(2^H - \{0\}) \times (2^K - \{0\})$ is a homeomorphic embedding;
\item[\rm (ii)] $f$ commutes with the action of $H \times K$;
\item[\rm (iii)] if $A \in \Su(H)$ and $B \in \Su(K)$ then $f(A \times B) \in \Su(H \times K)$;
\item[\rm (iv)] if $A \in \SuM(H)$ and $B \in \SuM(K)$, then $f(A \times B) \in \SuM(H \times K)$;
\item[\rm (v)] if $A \in \SuF(H)$ and $B \in \SuF(K)$, then $f(A \times B) \in \SuF(H \times K)$.
\end{enumerate}
\end{lem}

\begin{proof}
For notational convenience, we denote $(h, k) \in H \times K$ by $h k$. For $x \in 2^H$, $y \in 2^K$, $h \in H$, and $k \in K$ define
$$f(x,y)(hk) = \min(x(h), y(k)) = x(h) \cdot y(k).$$
So $f(x,y) = xy$ is the product of $x$ and $y$, as defined at the beginning of Section \ref{SECT COLOR GRP EXT}.

(i). Clearly $f$ is continuous. Suppose $x_0, x_1 \in 2^H - \{0\}$ and $y_0, y_1 \in 2^K - \{0\}$ satisfy $f(x_0, y_0) = f(x_1, y_1)$. We claim $x_0 = x_1$ and $y_0 = y_1$. Towards a contradiction, suppose $x_0 \neq x_1$ (the case $y_0 \neq y_1$ is similar). Let $h \in H$ be such that $x_0(h) \neq x_1(h)$. Then for some $i = 0, 1$ we have $x_i(h) = 0$. Therefore for all $k \in K$
$$y_{1-i}(k) = x_{1-i}(h) \cdot y_{1-i}(k) = f(x_{1-i}, y_{1-i})(hk)$$
$$= f(x_i, y_i)(h k) = x_i(h) \cdot y_i(k) = 0,$$
contradicting $y_{1-i} \neq 0$. We conclude $f$ is one-to-one on $(2^H - \{0\}) \times (2^K - \{0\})$. Now let $U \subseteq (2^H - \{0\}) \times (2^K - \{0\})$ be open. Then $U$ is open in $2^H \times 2^K$, and since $f(\{0\} \times 2^K \cup 2^H \times \{0\}) = 0 \not\in f(U)$, we have
$$f(U) = f(2^H \times 2^K) - f(2^H \times 2^K - U)$$
is open in $f(2^H \times 2^K)$ since $f(2^H \times 2^K - U)$ is compact. We have verified (i).

(ii). This is easily checked.

(iii). By (ii) $f(A \times B)$ is invariant under the action of $H \times K$. Since $A \times B$ is compact, $f(A \times B)$ is also compact and hence closed.

(iv). Fix $x \in A$ and $y \in B$. Since $A$ and $B$ are minimal, $x$ and $y$ are minimal and $\overline{[x]} = A$ and $\overline{[y]} = B$. So $A \times B = \overline{[(x,y)]}$ and $f(A \times B) = \overline{[f(x,y)]}$ since $f(A \times B)$ is compact. So it suffices to show that $f(x,y) \in 2^{H \times K}$ is minimal. If $f(x,y)$ is identically $0$, then it is trivially minimal. Otherwise, $f(x,y)$ is minimal by clause (iv) of Proposition \ref{prop:productcoloring}.

(v). It suffices to show that $f(x,y)$ is a $2$-coloring for every $x \in A$ and $y \in B$. If $x \in A$ and $y \in B$, then $x$ and $y$ are $2$-colorings since $A$ and $B$ are free. So $f(x,y)$ is a $2$-coloring by clause (i) of Proposition \ref{prop:productcoloring}.
\end{proof}

Define $P(H,K) = (2^H - \{0\}) \times (2^K - \{0\})$. Clauses (i) and (ii) of the previous lemma say that $P(H,K)$ and $f(P(H,K))$ are topologically conjugate. In other words, they have identical topology and dynamics arising from the action of $H \times K$. If we define $A(H) = \{x \in 2^H \: 0 \not\in \overline{[x]}\}$ and $A(K) = \{y \in 2^K \: 0 \not\in \overline{[y]}\}$ then we have the following.

\begin{theorem}
If $H$ and $K$ are countable groups then
\begin{enumerate}
\item[\rm (i)] $(\TCP(H) \res A(H)) \times (\TCP(K) \res A(K)) \sqsubseteq_c \TCP(H \times K)$;
\item[\rm (ii)] $\TCM(H) \times \TCM(K) \leq_B \TCM(H \times K)$;
\item[\rm (iii)] $\TCMF(H) \times \TCMF(K) \sqsubseteq_c \TCMF(H \times K)$.
\end{enumerate}
\end{theorem}

\begin{proof}
(i) and (iii) follow immediately from Corollary \ref{COR TC PROD} and the previous lemma. Let $1$ denote the element of $2^H$ which has constant value $1$. Define $Q_H: \SuM(H) \rightarrow \SuM(H)$ by
$$Q_H(A) = \begin{cases}
A & \text{if } 0 \not\in A \\
1 & \text{if } 0 \in A.
\end{cases}$$
Notice that if $0 \in A$ then $A = \{0\}$ by minimality of $A$. Also, note $Q_H(A) \ \TCM(H) \ A$, so that $Q_H(A_1) \ \TCM(H) \ Q_H(A_2)$ if and only if $A_1 \ \TCM(H) \ A_2$. Define $Q_K$ similarly. Then $Q_H$ and $Q_K$ are Borel. Now the map $f \circ (Q_H \times Q_K): \SuM(H) \times \SuM(K) \rightarrow \SuM(H \times K)$ is a Borel reduction of $\TCM(H) \times \TCM(K)$ to $\TCM(H \times K)$.
\end{proof}

\section{Topological conjugacy of minimal free subflows} \label{SECT ISO}

In this section, we show that $E_0$ continuously embeds into $\TCP(G)$ and Borel embeds into $\TCMF(G)$ for every countably infinite group $G$. Something which makes $E_0$ easy to work with is that it deals with one-sided infinite sequences, as do our blueprints. The basic idea will be the following. We will fix a fundamental function $c \in 2^{\subseteq G}$, and for each $x \in 2^\N$ we will build a function $e(x) \in 2^G$ extending $c$ in such a way that, for every $n \geq 1$, $e(x) \res \Delta_n \Theta_n(c)$ will depend only on $x \res \{i \in \N \: i \geq n-1\}$. If $x, y \in 2^\N$ are $E_0$-equivalent, then $e(x)$ and $e(y)$ will only be different on a small scale. So, we should be able to build a conjugacy between $\overline{[e(x)]}$ and $\overline{[e(y)]}$ using a block code with a large domain. On the other hand, if $x$ and $y$ are not $E_0$-equivalent, then on arbitrarily large subsets of $G$ $e(x)$ and $e(y)$ will have distinctly different behavior and therefore will not be topologically conjugate. With this basic outline of the proof in mind, the details should be easy to follow.

We begin with a very simple lemma. We point out that an immediate consequence of Theorem \ref{THM BC} is that every continuous function which commutes with the action of $G$ and is defined on a subflow of $2^G$ can be extended (not necessarily uniquely) to a continuous function commuting with the action of $G$ defined on all of $2^G$.

\begin{lem} \label{LEM FISO}
Let $G$ be a countably infinite group, let $(\Delta_n, F_n)_{n \in \N}$ be a blueprint guided by a growth sequence $(H_n)_{n \in \N}$, and let $c \in 2^G$ be fundamental with respect to this blueprint. If $y \in 2^G$ and if there is a continuous function from $2^G$ to itself which commutes with the action of $G$ and sends $c$ to $y$, then there is $n \geq 1$ so that
$$\forall \gamma, \psi \in \Delta_{n+3} \ ( \forall f \in F_{n+3} \ c(\gamma f) = c(\psi f) \Longrightarrow \forall f \in F_n \ y(\gamma f) = y(\psi f)).$$
\end{lem}

\begin{proof}
Let $\phi: 2^G \rightarrow 2^G$ be a function satisfying the hypothesis. Then $\phi$ is induced by a block code $\hat{\phi}$, and there is $n \in \N$ with $\dom(\hat{\phi}) \subseteq H_n$. Let $\gamma, \psi \in \Delta_{n+3}$ satisfy $c(\gamma f) = c(\psi f)$ for all $f \in F_{n+3}$. Then it follows from Lemma \ref{F TEST} that $c(\gamma h) = c(\psi h)$ for all $h \in H_{n+1}$. Thus, for $f \in F_n$ we have $f \dom(\hat{\phi}) \subseteq H_n H_n \subseteq H_{n+1}$, so
$$(f^{-1} \gamma^{-1} \cdot c)\res \dom(\hat{\phi}) = (f^{-1} \psi^{-1} \cdot c)\res \dom(\hat{\phi}).$$
It follows that $y(\gamma f) = \phi(c)(\gamma f) = \phi(c)(\psi f) = y(\psi f)$.
\end{proof}

Let $(\Delta_n, F_n)_{n \in \N}$ be a centered blueprint guided by a growth sequence $(H_n)_{n \in \N}$. Recall that such a blueprint is necessarily directed and maximally disjoint and furthermore $(F_n)_{n \in \N}$ is an increasing sequence and $(\Delta_n)_{n \in \N}$ is a decreasing sequence (see Lemma \ref{LEM GUIDED BP} and clause (i) of Lemma \ref{STRONG BP LIST}). The following two functions will be very useful in defining the function $e: 2^\N \rightarrow 2^G$. Define
$$r: \bigcup_{n \geq 1} \left( \Delta_n \times \{n\} \right) \rightarrow \N$$
by
$$r(\gamma, n) = \min\{k > n \: \gamma \in \Delta_k F_k\}$$
for $(\gamma, n) \in \dom(r)$. Additionally, define
$$L: \bigcup_{n \geq 1} \left( \Delta_n \times \{n\} \right) \rightarrow \Delta_1$$
so that for $(\gamma, n) \in \dom(L) \ L(\gamma, n) = \psi$, where $\psi$ is the unique element of $\Delta_{r(\gamma, n)}$ with $\gamma \in \psi F_{r(\gamma, n)}$. Intuitively, the functions $r$ and $L$ together allow one to ``lift'' a $\Delta_k$-translate of $F_k$, say $\gamma F_k$, to a $\Delta_m$-translate of $F_m$ containing $\gamma F_k$, where $m > k$ is least with $\gamma F_k \subseteq \Delta_m F_m$. Formally this is expressed as $\gamma F_k \subseteq L(\gamma, k) F_{r(\gamma, k)}$. For convenience we let $L^0 (\gamma, n) = \gamma, \ r^0 (\gamma, n) = n, \ L^1 = L$, and $r^1 = r$. In general, for $k > 1$ let
$$r^k (\gamma, n) = r( L^{k-1} (\gamma, n), r^{k-1} (\gamma, n) ),$$
and
$$L^k ( \gamma, n) = L( L^{k-1} (\gamma, n), r^{k-1} (\gamma, n) ).$$
These functions will only be used in this section.

\begin{lem} \label{LIFTS}
Let $G$ be a countably infinite group and let $(\Delta_n, F_n)_{n \in \N}$ be a centered blueprint guided by a growth sequence $(H_n)_{n \in \N}$. Then we have the following:
\begin{enumerate}
\item[\rm (i)] if $n \geq 1$, $\gamma \in \Delta_n$, and $1 \leq k < n$ then $r(\gamma, k) = k+1$ and $L(\gamma, k) = \gamma$;
\item[\rm (ii)] $\gamma \in L^k(\gamma, n) F_{r^k(\gamma, n)}$ for all $n \geq 1$, $\gamma \in \Delta_n$, and $k \in \N$;
\item[\rm (iii)] if $\gamma \in \Delta_n$, $m \geq n \geq 1$, $\sigma \in \Delta_m$, and $\gamma \in \sigma F_m$, then there exists $k \in \N$ with $r^k (\gamma, n) = m$ and $L^k(\gamma, n) = \sigma$;
\item[\rm (iv)] for all $n \geq 1$ and $\gamma \in \Delta_n$, there is $N \in \N$ so that for all $k \geq N \ L^k(\gamma, n) = 1_G$ and $r^k(\gamma, n) = k - N + r^N(\gamma, n)$;
\item[\rm (v)] for all $1 \leq k \leq n$, $\lambda \in D^n_k$, and $\gamma \in \Delta_n$, if $m \in \N$ satisfies either $r^m(\lambda, k) = n$ or $r^m(\gamma \lambda, k) = n$ then for all $0 \leq i \leq m$
$$r^i(\gamma \lambda, k) = r^i(\lambda, k), \text{ and}$$
$$L^i(\gamma \lambda, k) = \gamma L^i(\lambda, k);$$
\item[\rm (vi)] for all $n \geq 1$, $\gamma \in \Delta_n$, and $m \geq 1$
$$L( [L^m(\gamma, n)]^{-1} L^{m-1}(\gamma, n), r^{m-1}(\gamma, n) ) = 1_G$$
and
$$r( [L^m(\gamma, n)]^{-1} L^{m-1}(\gamma, n), r^{m-1}(\gamma, n) ) = r^m(\gamma, n).$$
\end{enumerate}
\end{lem}

\begin{proof}
(i). Clearly $\gamma \in \Delta_k$ and $\gamma \in \gamma F_k$ for all $1 \leq k \leq n$.

(ii). By definition $\gamma \in L^1(\gamma, n) F_{r^1(\gamma, n)}$. Suppose $\gamma \in L^m(\gamma, n) F_{r^m(\gamma, n)}$. Then $L^m(\gamma, n) \in \Delta_{r^m(\gamma, n)}$ and $L^m(\gamma, n) \in L^{m+1}(\gamma, n) F_{r^{m+1}(\gamma, n)}$. So by the coherent property of blueprints
$$\gamma \in L^m(\gamma, n) F_{r^m(\gamma, n)} \subseteq L^{m+1}(\gamma, n) F_{r^{m+1}(\gamma, n)}.$$

(iii). Note that in general for $(\psi, i) \in \dom(r)$, $r(\psi, i) > i$. If $n = m$ then $\gamma = \sigma$ and (iii) is satisfied by taking $k = 0$. Otherwise, let $k \in \N$ be maximal with $r^k(\gamma, n) < m$. Then by (ii) $L^k(\gamma, n) F_{r^k(\gamma, n)} \cap \sigma F_m \neq \varnothing$. By the coherent property of blueprints $L^k(\gamma, n) \in \sigma F_m$. Since $k$ is maximal with $r^k(\gamma, n) < m$, it follows from the definition of $r$ and $L$ that $r^{k+1}(\gamma, n) = m$ and $L^{k+1}(\gamma, n) = \sigma$.

(iv). This follows from clause (iv) of Lemma \ref{STRONG BP LIST} together with (iii) and (i).

(v). Notice that $m$ must exist by (iii). If $m = 0$ then $k = n$, $\lambda = 1_G$, and the claim is trivial. So assume $m > 0$. Clearly $\lambda, \gamma \lambda \in \Delta_n F_n$, so $r(\lambda, k), r(\gamma \lambda, k) \leq n$. By clause (vii) of Lemma \ref{BP LIST}, $\lambda \in \Delta_s F_s$ if and only if $\gamma \lambda \in \Delta_s F_s$ for $k < s \leq n$. It then follows from the definition of $r$ that $r(\lambda, k) = r(\gamma \lambda, k) = t \leq n$. Set $\psi = L(\lambda, k) \in \Delta_t$. We have $\lambda \in F_n \cap \psi F_t$, so by the coherent property of blueprints $\psi \in D^n_t$. Since $\gamma \psi \in \Delta_t$ and $\gamma \lambda \in \gamma \psi F_t$, we have $L(\gamma \lambda, k) = \gamma \psi$. Thus we have verified the claim for $i = 0$ and $i = 1$. The claim then follows by induction: replace $\lambda$ with $\psi$ and $k$ with $t$.

(vi). $L^{m-1} (\gamma, n) \in L^m (\gamma, n) F_{r^m (\gamma, n)}$, so there is $\lambda \in D^{r^m (\gamma, n)}_{r^{m-1} (\gamma, n)}$ such that $L^{m-1} (\gamma, n) = L^m (\gamma, n) \lambda$. Then by (v)
$$L^m (\gamma, n) = L(L^{m-1}(\gamma, n), r^{m-1}(\gamma, n)) = L^m (\gamma, n) L(\lambda, r^{m-1} (\gamma, n)),$$
which implies
$$L( [L^m(\gamma, n)]^{-1} L^{m-1}(\gamma, n), r^{m-1}(\gamma, n) ) = L(\lambda, r^{m-1} (\gamma, n) ) = 1_G.$$
Clause (v) also implies that
$$r([L^m(\gamma, n)]^{-1} L^{m-1}(\gamma, n), r^{m-1}(\gamma, n) ) = r(L^{m-1}(\gamma, n), r^{m-1}(\gamma, n)) = r^m(\gamma, n).$$
\end{proof}

We are now prepared to prove the main theorem of this section. The following theorem appears to be quite nontrivial as it relies on all of the machinery developed in Chapters \ref{CHAP FM} and \ref{CHAP STUDY}.

\begin{theorem} \label{THM MF}
For any countably infinite group $G$, $E_0$ continuously embeds into $\TCP(G)$ and embeds into $\TC(G)$, $\TCF(G)$, $\TCM(G)$, and $\TCMF(G)$.
\end{theorem}

\begin{proof}
For $n \geq 1$ and $k \in \N$ define $p_n(k) = 4 \cdot (2 k^4 +1) \cdot (12 k^4+1)$ and $q_n(k) = 2k^3$. Then $(p_n)_{n \geq 1}$ is a sequence of functions of subexponential growth. By Proposition \ref{SCHEME 1}, there is a centered blueprint $(\Delta_n, F_n)_{n \in \N}$ guided by a growth sequence $(H_n)_{n \in \N}$ with $|\Lambda_n| \geq q_n(|F_{n-1}|) + \log_2 \ p_n(|F_n|)$ for each $n \geq 1$ and such that for every $g \in G - \Z(G)$ and every $n \geq 1$ there are infinitely many $\gamma \in \Delta_n$ with $\gamma g \neq g \gamma$. We are free to pick any distinct $\alpha_n, \beta_n, \gamma_n \in D^n_{n-1}$ for each $n \geq 1$. We choose $\gamma_n = 1_G$ and let $\alpha_n$ and $\beta_n$ be arbitrary for every $n \geq 1$. By clause (i) of Lemma \ref{LEM GUIDED BP} the blueprint is directed and maximally disjoint, and by clause (viii) of Lemma \ref{STRONG BP LIST} $\bigcap_{n \in \N} \Delta_n a_n = \bigcap_{n \in \N} \Delta_n b_n = \varnothing$. Apply Theorem \ref{FM} to get a function $c \in 2^{\subseteq G}$ which is canonical with respect to this blueprint. By Proposition \ref{CANONICAL DMIN} the function $c$ is $\Delta$-minimal. Apply Corollary \ref{GEN MINCOL} to get a function $c' \supseteq c$ which is fundamental with respect to $(\Delta_n, F_n)_{n \in \N}$, is $\Delta$-minimal, has $|\Theta_n(c')| > 1 + \log_2 \ (12|F_n|^4+1)$ for each $n \geq 1$, and has the property that every extension of $c'$ to all of $G$ is a $2$-coloring. Now apply Corollary \ref{RIGID} to get a fundamental and $\Delta$-minimal $c'' \supseteq c'$ and a collection $\{\nu_i^n \in \Delta_{n+5} \: n \equiv 1 \mod 5, 1 \leq i \leq s(n)\}$ where $s(n) = 2$ if $n \equiv 1 \mod 10$ and $s(n) = |F_n F_n^{-1} - \Z(G)|$ otherwise. We have that $|\Theta_n(c'')| \geq 1$ for all $n \geq 1$, $c''(f) = c''(\nu^n_i f)$ for all $n \equiv 1 \mod 5$, $1 \leq i \leq s(n)$, and $f \in F_{n+4} \cap \dom(c'')$, and if $x, y \in 2^G$ extend $c''$ and $x(f) = x(\nu^n_i f)$ for all $n \equiv 1 \mod 5$, $1 \leq i \leq s(n)$, and $f \in F_{n+4}$ then $\overline{[x]}$ and $\overline{[y]}$ are topologically conjugate if and only if there is a conjugacy mapping $x$ to a $y$-centered element of $\overline{[y]}$. By using Lemma \ref{MIN GRAPH} (with $\mu$ identically $0$), Lemma \ref{LEM UNION}, and Lemma \ref{MIN UNION}, we may suppose without loss of generality that $|\Theta_n(c'')| = 1$ for all $n \geq 1$.

For $k \in \N$, let $r^k$ and $L^k$ be defined as in the paragraph preceding Lemma \ref{LIFTS}. We wish to find a function
$$\mu : \bigcup_{n \geq 2} \left( \{\gamma \in D^n_k \: 1 \leq k < n, \ r(\gamma, k) = n\}  \times \{n\} \right) \rightarrow \{0, 1\}$$
which satisfies:
\begin{enumerate}
\item[\rm (1)] for each $n \geq 2 \ \mu(1_G, n) = 0$;
\item[\rm (2)] for each $n \geq 2$ there is $\psi \in D^n_{n-1}$ with $\mu(\psi, n) = 1$;
\item[\rm (3)] if $(\gamma, n) \in \dom(\mu)$ and $\gamma \not\in D^n_{n-1}$ then $\mu(\gamma, n) = 0$;
\item[\rm (4)] for every $k \equiv 1 \mod 5$, $1 \leq i \leq s(k)$, and $m \in \N$
$$\mu([L^{m+1} (\nu_i^k, k+5) ]^{-1} L^m (\nu_i^k, k+5), r^{m+1}(\nu_i^k, k+5))$$
$$= \mu(1_G, r^{m+1}(\nu_i^k, k+5)) = 0.$$
\end{enumerate}
It may aid the reader to note that $\dom(\mu)$ may be expressed in a possibly more understandable form. Let $\pi_1: G \times \N \rightarrow G$ be the first component projection map. Then
$$\dom(\mu) = \bigcup_{n \geq 2} \pi_1 [ L^{-1}(1_G) \cap r^{-1}(n) ] \times \{n\}.$$
Note that $D^n_{n-1} \times \{n\} \subseteq \dom(\mu)$ for each $n \geq 2$, and that the expression in (4) lies in the domain of $\mu$ by conclusion (vi) of Lemma \ref{LIFTS}.

Clearly (1), (2), and (3) are achievable, and (4) is consistent with (1) and (3). The only difficulty is to show that (2) and (4) can be simultaneously achieved. However, we only need to observe that if $r^{m+1}(\nu_i^k, k+5) = n$ then $k+5 < n$ and
$$\sum\limits_{ \begin{array}{c} k \equiv 1 \mod 5 \\ k < n-5 \end{array} } s(k) < \sum\limits_{ \begin{array}{c} k \equiv 1 \mod 5 \\ k < n-5 \end{array} } |F_k|^2 < n \cdot |F_{n-1}|^2 < |F_{n-1}|^3 \leq \frac{1}{2} |D^n_{n-1}|.$$
The last inequality holds due to the definition of $q_n$ and $c$ in the first paragraph. Therefore (2) and (4) can be simultaneously achieved, and such a function $\mu$ exists.

For each $n \geq 1$, let $\theta_n$ be the unique element of $\Theta_n(c'')$. For $x \in 2^\N$ define $e(x) \in 2^G$ so that $e(x) \supseteq c''$ and for $n \geq 1$ and $\gamma \in \Delta_n$
$$e(x)(\gamma \theta_n b_{n-1}) = \sum\limits_{k=0}^\infty x( r^k(\gamma, n) - 1) \cdot \mu( [L^{k+1} (\gamma, n) ]^{-1} L^k (\gamma, n), r^{k+1} (\gamma, n) ) \mod 2.$$
This sum is finite by clause (iv) of Lemma \ref{LIFTS} and property (1) of $\mu$. Moreover, the number of indices for which the summand is nonzero is bounded independent of $x \in 2^\N$. It is therefore easy to see that $e: 2^\N \rightarrow 2^G$ is continuous (where $2^\N$ has the product topology). By Lemma \ref{LEM SUBMAP} the map $x \mapsto \overline{[e(x)]}$ is Borel. The expression above is well defined as clause (vi) of Lemma \ref{LIFTS} implies that
$$( [L^{k+1} (\gamma, n) ]^{-1} L^k (\gamma, n), r^{k+1} (\gamma, n) ) \in \dom(\mu)$$
for all $k \in \N$, $n \geq 1$, and $\gamma \in \Delta_n$. The function $e(x)$ has two useful properties which we list below.
\begin{enumerate}
\item[\rm (a)] Let $x \in 2^\N$, $1 \leq k \leq n$, $\gamma \in \Delta_n$, and $\lambda_1, \lambda_2 \in D^n_k$ satisfy $r(\lambda_1, k) = r(\lambda_2, k) = n$. If $x(k-1) = 1$ then
$$e(x)(\gamma \lambda_1 \theta_k b_{k-1}) = e(x)(\gamma \lambda_2 \theta_k b_{k-1}) \Longleftrightarrow \mu(\lambda_1, n) = \mu(\lambda_2, n)$$
and if $x(k-1) = 0$ then $e(x)(\gamma \lambda_1 \theta_k b_{k-1}) = e(x)(\gamma \lambda_2 \theta_k b_{k-1})$ always.
\item[\rm (b)] Let $x \in 2^\N$, $n \geq 1$, and $\gamma \in \Delta_{n+1}$. Then for every $f \in F_n - \dom(c'')$
$$e(x)(\gamma \theta_{n+1} b_n) = e(x)(\theta_{n+1} b_n) \Longleftrightarrow e(x)(\gamma f) = e(x)(f).$$
\end{enumerate}
We spend the next two paragraphs establishing the validity of (a) and (b).

(a). By clause (v) of Lemma \ref{LIFTS} $r(\gamma \lambda_1, k) = r(\gamma \lambda_2, k) = n$, and so by the definition of $L$ we must have $L(\gamma \lambda_1, k) = L(\gamma \lambda_2, k) = \gamma$. By the definition of $L^m$ and $r^m$, it follows that $r^m(\gamma \lambda_1, k) = r^m(\gamma \lambda_2, k)$ and $L^m(\gamma \lambda_1, k) = L^m(\gamma \lambda_2, k)$ for all $m \geq 1$. Thus when considering the summations defining $e(x)(\gamma \lambda_1 \theta_k b_{k-1})$ and $e(x)(\gamma \lambda_2 \theta_k b_{k-1})$, we see that all the summands are equal except possibly the first. If $x(k-1) = 1$, then the first summands are equal if and only if $\mu(\lambda_1, n) = \mu(\lambda_2, n)$. If $x(k-1) = 0$, then the first summands are always equal. Property (a) now clearly follows.

(b). Fix $f \in F_n - \dom(c'')$. Since $G - \dom(c'') = \bigcup_{k \geq 1} \Delta_k \theta_k b_{k-1}$, there is $k \geq 1$ with $f \in \Delta_k \theta_k b_{k-1}$. Since $\gamma_{n+1} = 1_G$, $1_G \in \Delta_{n+1}$, and $\beta_{n+1} \neq \gamma_{n+1} \in D^{n+1}_n$, we have that if $k > n+1$ then
$$\varnothing \neq F_n \cap \Delta_k \theta_k b_{k-1} \subseteq \Delta_{n+1} \gamma_{n+1} F_n \cap \Delta_{n+1} \beta_{n+1} F_n = \varnothing,$$
a contradiction. So $k \leq n+1$ (one can further show that $k \leq n$, but we do not need this). Since $f \in F_{n+1} \cap \Delta_k F_k$ (recall $F_n \subseteq F_{n+1}$ since our blueprint is centered), it follows by the coherent property of blueprints that there is $\lambda \in D^{n+1}_k$ with $f = \lambda \theta_k b_{k-1}$. Let $m \in \N$ be such that $r^m(\lambda, k) = n+1$. Then $L^m(\lambda, k) = 1_G$. By clause (v) of Lemma \ref{LIFTS} we have that $r^i(\gamma \lambda, k) = r^i(\lambda, k)$ and $L^i(\gamma \lambda, k) = \gamma L^i(\lambda, k)$ for all $0 \leq i \leq m$. It follows that in the summations defining $e(x)(\gamma \lambda \theta_k b_{k-1})$ and $e(x)(\lambda \theta_k b_{k-1})$, the first $m$ terms (the terms where the index of the sum is between $0$ and $m-1$, inclusive) are respectively equal. Let $S$ denote the common value of the sums of the first $m$ terms. For $i > m$ we have that $L^i(\lambda, k) = L^{i-m}(1_G, n+1)$, $r^i(\lambda, k) = r^{i-m}(1_G, n+1)$, $r^i(\gamma \lambda, k) = r^{i-m}(\gamma, n+1)$, and $L^i(\gamma \lambda, k) = L^{i-m}(\gamma, n+1)$. Therefore we see that
$$e(x)(\gamma \lambda \theta_k b_{k-1}) = S + e(x)(\gamma \theta_{n+1} b_n)$$
and
$$e(x)(\lambda \theta_k b_{k-1}) = S + e(x)(\theta_{n+1} b_n).$$
Recalling that $f = \lambda \theta_k b_{k-1}$, we conclude
$$e(x)(\gamma \theta_{n+1} b_n) = e(x)(\theta_{n+1} b_n) \Longleftrightarrow e(x)(\gamma f) = e(x)(f).$$

For $x \in 2^\N$, define $\overline{e}(x) = \overline{[e(x)]}$. We will show that $e: 2^\N \rightarrow 2^G$ is a continuous embedding of $E_0$ into $\TCP(G)$ and that $\overline{e}: 2^\N \rightarrow \SuMF(G)$ is a Borel embedding of $E_0$ into $\TCMF(G)$. From this the validity of the theorem will follow by Lemma \ref{LEM EQ}. As mentioned immediately after the definition of $e$, the function $e$ is indeed continuous and the function $\overline{e}$ is indeed Borel. In the next two paragraphs we prove that the image of $\overline{e}$ is contained in $\SuMF(G)$ and that $\overline{e}$ is injective. An immediate consequence of this is that $e$ is also injective.

We check that $e(x)$ is a minimal $2$-coloring for all $x \in 2^\N$. From this it will follow that the image of $\overline{e}$ is contained in $\SuMF(G)$. The fact that $e(x)$ is a $2$-coloring is immediate since $e(x)$ extends $c''$. So we check that $e(x)$ is minimal. The function $e(x)$ is defined on all of $G$ and extends the fundamental function $c''$. Thus $e(x)$ is fundamental. Since our blueprint is centered, directed, and $\alpha_n \neq \gamma_n = 1_G \neq \beta_n$ for all $n \geq 1$, we have that $\bigcap_{n \in \N} \Delta_n b_n = \varnothing$ by clause (viii) of Lemma \ref{STRONG BP LIST}. So by Corollary \ref{STRONG PREMIN}, it suffices to show that for every $k \geq 1$ there is $n > k$ so that for all $\gamma \in \Delta_n$ there is $\lambda \in D^n_k$ with $e(x)(\gamma \lambda f) = e(x)(f)$ for all $f \in F_k$. So fix $k \geq 1$. Since $c''$ is $\Delta$-minimal, there is $m > k$ so that for all $\gamma \in \Delta_m$ we have $c''(\gamma f) = c''(f)$ for all $f \in F_k \cap \dom(c'')$. We now proceed by cases. \underline{Case 1}: $x(i-2) = 0$ for all $i > m$. Set $n = m+1$ and let $\gamma \in \Delta_n$. Set $\lambda = 1_G \in D^n_{n-1}$. After inspecting the summation defining $e(x)$ we see that
$$e(x)(\gamma \lambda \theta_{n-1} b_{n-2}) = e(x)(\gamma \theta_{n-1} b_{n-2}) = 0 = e(x)(\theta_{n-1} b_{n-2}).$$
We have that $\gamma \lambda \in \Delta_{n-1} = \Delta_m$, so by our choice of $m$ and property (b) we have $e(x)(\gamma \lambda f) = e(x)(f)$ for all $f \in F_k$. Thus $e(x)$ is minimal. \underline{Case 2}: There is $n > m$ with $x(n-2) = 1$. Fix $\gamma \in \Delta_n$. By (a) and properties (1) and (2) of $\mu$ there must be $\lambda \in D^n_{n-1}$ with $e(x)(\gamma \lambda \theta_{n-1} b_{n-2}) = e(x)(\theta_{n-1} b_{n-2})$. Again by (b) and our choice of $m$ we have $e(x)(\gamma \lambda f) = e(x)(f)$ for all $f \in F_k$. We conclude that $e(x)$ is minimal.

Now we check that $\overline{[e(x)]} \neq \overline{[e(y)]}$ for $x \neq y \in 2^\N$. It will suffice to show that $e(x)$ and $e(y)$ are orthogonal. Fix $x \neq y \in 2^\N$, and let $n \geq 1$ be such that $x(n-1) \neq y(n-1)$. Set $T = F_{n+1} F_{n+1}^{-1} F_{n+1}$, and let $g_1, g_2 \in G$ be arbitrary. Let $t \in F_{n+1} F_{n+1}^{-1}$ be such that $g_1 t \in \Delta_{n+1}$. If $g_2 t \not\in \Delta_{n+1}$ or if $e(x)(g_1 t \theta_n b_{n-1}) \neq e(y)(g_2 t \theta_n b_{n-1})$ then we are done. So suppose $g_2 t \in \Delta_{n+1}$ and $e(x)(g_1 t \theta_n b_{n-1}) = e(y)(g_2 t \theta_n b_{n-1})$. By property (2) of $\mu$, let $\psi \in D^{n+1}_n$ be such that $\mu(\psi, n+1) = 1 \neq 0 = \mu(1_G, n+1)$. Then it follows from (a) that $e(x)(g_1 t \psi \theta_n b_{n-1}) \neq e(y)(g_2 t \psi \theta_n b_{n-1})$. As $t \psi \theta_n b_{n-1} \in T$, we conclude $e(x)$ is orthogonal to $e(y)$ and $\overline{[e(x)]} \neq \overline{[e(y)]}$.

Now we move into the final stage of the proof. The remaining task is to show that for $x, y \in 2^\N$
$$x \ E_0 \ y \Longleftrightarrow e(x) \ \TCP(G) \ e(y) \Longleftrightarrow \overline{e}(x) \ \TCMF(G) \ \overline{e}(y).$$
In order to achieve this task, we rely on the rigidity constructions from Section \ref{SEC RIGID}. In particular, we invoke the fact that $c''$ originated from Corollary \ref{RIGID}. In the next paragraph, we will show that $e(x)(\nu_i^n f) = e(x)(f)$ for all $n \equiv 1 \mod 5$, $1 \leq i \leq s(n)$, $f \in F_{n+4}$, and $x \in 2^\N$. So by Corollary \ref{RIGID}, we have that $\overline{e}(x) \ \TCMF(G) \ \overline{e}(y)$ if and only if there is a conjugacy sending $e(x)$ to an $e(y)$-centered element of $\overline{[e(y)]}$. Therefore, in the second paragraph we briefly study $e(x)$-centered elements of $\overline{[e(x)]}$ for $x \in 2^\N$. Then in the final two paragraphs we prove the validity of the displayed expression above, completing the proof of the theorem.

For $x \in 2^\N$, $n \equiv 1 \mod 5$, and $1 \leq i \leq s(n)$, by considering property (4) of $\mu$ and the summation defining $e(x)$, we see that
$$e(x)(\nu_i^n \theta_{n+5} b_{n+4}) = 0 = e(x)(\theta_{n+5} b_{n+4}).$$
By the definition of $c''$ we have that $e(x)(\nu^n_i f) = e(x)(f)$ for all $f \in F_{n+4} \cap \dom(c'')$. So by (b) we have
$$e(x)(\nu_i^n f) = e(x)(f)$$
for all $n \equiv 1 \mod 5$, $1 \leq i \leq s(n)$, $f \in F_{n+4}$, and $x \in 2^\N$.

We now briefly discuss $e(x)$-centered elements of $\overline{[e(x)]}$. Let $w \in \overline{[e(x)]}$ be $e(x)$-centered. By Lemma \ref{MIN REG1} we have that $w(g) = e(x)(g)$ for all $g \in \dom(c'')$. Fix $g \not\in \dom(c'')$. Then there is $k \geq 1$ and $\psi \in \Delta_k$ with $g = \psi \theta_k b_{k-1}$. Since our blueprint is centered and directed, there is $n \geq k$ with $\psi F_k \subseteq F_n$ (clause (iv) of Lemma \ref{STRONG BP LIST}). Let $m > n$ and note $\psi \in D^n_k \subseteq D^m_k$. Since $w$ is an $e(x)$-centered element of $\overline{[e(x)]}$, by clause (i) of Proposition \ref{SUBFLOW BP} there is $\gamma \in \Delta_m$ with $w(f) = (\gamma^{-1} \cdot e(x))(f) = e(x)(\gamma f)$ for all $f \in F_m$. We have the following equivalences:
$$w(g) = e(x)(g) \Longleftrightarrow w(\psi \theta_k b_{k-1}) = e(x)(\psi \theta_k b_{k-1}) \Longleftrightarrow$$
$$e(x)(\gamma \psi \theta_k b_{k-1}) = e(x)(\psi \theta_k b_{k-1}) \Longleftrightarrow e(x)(\gamma \theta_m b_{m-1}) = e(x)(\theta_m b_{m-1})$$
$$\Longleftrightarrow w(\theta_m b_{m-1}) = e(x)(\theta_m b_{m-1})$$
(the second line is due to (b)). Therefore, we have that $w(g) = e(x)(g)$ if and only if $w(\theta_m b_{m-1}) = e(x)(\theta_m b_{m-1})$ for all sufficiently large $m$. Since the second condition does not depend on $g \in G - \dom(c'')$, we have that either $w = e(x)$ or else
$$w(g) = \begin{cases}
e(x)(g) & \text{if } g \in \dom(c'') \\
1 - e(x)(g) & \text{otherwise.}
\end{cases}$$
So $\overline{[e(x)]}$ contains at most two $e(x)$-centered elements (counting $e(x)$ itself). Moreover, an important observation is that if $w \neq e(x)$ is an $e(x)$-centered element of $\overline{[e(x)]}$, then $w$ satisfies properties (a) and (b) with respect to the sequence $x \in 2^\N$. Although we will not make any use of this fact whatsoever, we do mention that if $\overline{[e(x)]}$ does contain two $e(x)$-centered elements, then they are never $\TCP(G)$-equivalent, despite their strong similarities.

In the remaining two paragraphs we prove that for $x, y \in 2^\N$
$$x \ E_0 \ y \Longleftrightarrow e(x) \ \TCP(G) \ e(y) \Longleftrightarrow \overline{e}(x) \ \TCMF(G) \ \overline{e}(y).$$
Here we prove that the negation of the leftmost expression implies the negations of the other two. In the next and final paragraph we prove that the leftmost expression implies the other two expressions. Let $x, y \in 2^\N$ be such that $\neg(x \ E_0 \ y)$, or in other words $x(n) \neq y(n)$ for infinitely many $n \in \N$. Without loss of generality we may assume $x(n) = 0 \neq y(n)$ for infinitely many $n \in \N$. Towards a contradiction, suppose $e(x)$ is $\TCP(G)$-equivalent to some $e(y)$-centered $z \in \overline{[e(y)]}$ (potentially $z = e(y)$). By Lemma \ref{LEM FISO}, there is $n \geq 1$ such that if $\gamma, \psi \in \Delta_{n+3}$ and $e(x)(\gamma f) = e(x)(\psi f)$ for all $f \in F_{n+3}$ then $z(\gamma f) = z(\psi f)$ for all $f \in F_n$. Since $c''$ is $\Delta$-minimal, there is $m > n + 3$ so that $e(x)(\gamma f) = e(x)(f)$ for all $\gamma \in \Delta_m$ and $f \in F_{n+3} \cap \dom(c'')$. Pick $k > m$ with $x(k-1) = 0 \neq y(k-1)$, and pick $\gamma \in D^{k+1}_k \subseteq \Delta_k$ with $\mu(\gamma, k+1) = 1$. Then $L(\gamma, k) = 1_G = L(1_G, k)$ and $r(\gamma, k) = k+1 = r(1_G, k)$, so by (a) $e(x)(\gamma \theta_k b_{k-1}) = e(x)(\theta_k b_{k-1})$ and $z(\gamma \theta_k b_{k-1}) \neq z(\theta_k b_{k-1})$ (since $z$ satisfies properties (a) and (b) with respect to the sequence $y \in 2^\N$). Since $\gamma \in \Delta_k \subseteq \Delta_m$ and $k - 1 \geq n+3$, we have by (b) that $e(x)(\gamma f) = e(x)(f)$ for all $f \in F_{n+3}$. It follows that $z(\gamma f) = z(f)$ for all $f \in F_n$. In particular, $z(\gamma \theta_n b_{n-1}) = z(\theta_n b_{n-1})$ which, by (b), is in contradiction with $z(\gamma \theta_k b_{k-1}) \neq z(\theta_k b_{k-1})$.

Now let $x, y \in 2^\N$ be such that $x \ E_0 \ y$. To complete the proof of the theorem, it suffices to show that $e(x) \ \TCP \ e(y)$. Let $N \in \N$ be such that $x(n-1) = y(n-1)$ for all $n \geq N$. Set $K = F_N^{-1} F_N$. To show $e(x) \ \TCP(G) \ e(y)$, it is sufficient, by Corollary \ref{COR ISOTEST}, to show that
$$\forall g, h \in G \ (\forall k \in K \ e(x)(gk) = e(x)(hk) \Longrightarrow e(y)(g) = e(y)(h)), \text{ and}$$
$$\forall g, h \in G \ (\forall k \in K \ e(y)(gk) = e(y)(hk) \Longrightarrow e(x)(g) = e(x)(h)).$$
By symmetry of information regarding $x$ and $y$, it will be enough to verify the first property above. Let $g, h \in G$ be such that $e(x)(gk) = e(x)(hk)$ for all $k \in K$. We will show $e(y)(g) = e(y)(h)$. Note that $e(x)$ and $e(y)$ agree on $\dom(c'')$ and on $\bigcup_{n \geq N} \Delta_n \theta_n b_{n-1}$. So we may suppose at least one of $g, h$ is in
$$G - \left(\dom(c'') \cup \bigcup_{m \geq N} \Delta_m \theta_m b_{m-1} \right) = \bigcup_{1 \leq m < N} \Delta_m \theta_m b_{m-1}.$$
However, by our choice of $K$, for $m < N$ one of $g$ or $h$ is in $\Delta_m \theta_m b_{m-1}$ if and only if both are (since $e(x)$ has a $\Delta_m$ membership test with test region a subset of $F_m \subseteq F_N$). So let $1 \leq m < N$ be such that $g, h \in \Delta_m \theta_m b_{m-1}$. Let $n \leq N$ be maximal with $g \in \Delta_n F_n$. Again, since $e(x)$ has a $\Delta_n$ membership test, this same $n$ equals the maximal $i \leq N$ with $h \in \Delta_i F_i$. It follows that there are $\gamma, \psi \in \Delta_n$ and $\lambda \in D^n_m$ with $g = \gamma \lambda \theta_m b_{m-1}$ and $h = \psi \lambda \theta_m b_{m-1}$. Let $k \geq -1$ be such that $r^{k+1}(\gamma \lambda, m) = n$. By conclusion (v) of Lemma \ref{LIFTS}, for any $w \in 2^\N$
$$\sum\limits_{i=0}^k w(r^i(\gamma \lambda, m) - 1) \cdot \mu( [ L^{i+1} (\gamma \lambda, m) ]^{-1} L^i (\gamma \lambda, m), r^{i+1} (\gamma \lambda, m) ) \mod 2$$
$$= \sum\limits_{i=0}^k w(r^i(\psi \lambda, m) - 1) \cdot \mu( [ L^{i+1} (\psi \lambda, m) ]^{-1} L^i (\psi \lambda, m), r^{i+1} (\psi \lambda, m) ) \mod 2.$$
In particular, the first $k+1$ terms of the sums defining $e(y)(g)$ and $e(y)(h)$ are respectively equal. Since $e(x)(g) = e(x)(h)$, the above equality implies
$$\sum\limits_{i=k+1}^\infty x(r^i(\gamma \lambda, m) - 1) \cdot \mu( [ L^{i+1} (\gamma \lambda, m) ]^{-1} L^i (\gamma \lambda, m), r^{i+1} (\gamma \lambda, m) ) \mod 2$$
$$= \sum\limits_{i=k+1}^\infty x(r^i(\psi \lambda, m) - 1) \cdot \mu( [ L^{i+1} (\psi \lambda, m) ]^{-1} L^i (\psi \lambda, m), r^{i+1} (\psi \lambda, m) ) \mod 2.$$
If $n = N$, then $r^{k+1} (\gamma \lambda, m) = r^{k+1} (\psi \lambda, m) = N$ so
$$\sum\limits_{i=k+1}^\infty x( r^i(\gamma \lambda, m) - 1) \cdot \mu( [ L^{i+1} (\gamma \lambda, m) ]^{-1} L^i (\gamma \lambda, m), r^{i+1} (\gamma \lambda, m) ) \mod 2$$
$$= \sum\limits_{i=k+1}^\infty y(r^i(\gamma \lambda, m) - 1) \cdot \mu( [ L^{i+1} (\gamma \lambda, m) ]^{-1} L^i (\gamma \lambda, m), r^{i+1} (\gamma \lambda, m) ) \mod 2,$$
and similarly for $\psi$ in place of $\gamma$. On the other hand, if $n < N$, then since $r^{k+2} (\gamma \lambda, m), r^{k+2} (\psi \lambda, m) > N$, property (3) of $\mu$ gives
$$\sum\limits_{i=k+1}^\infty x( r^i(\gamma \lambda, m) - 1) \cdot \mu( [ L^{i+1} (\gamma \lambda, m) ]^{-1} L^i (\gamma \lambda, m), r^{i+1} (\gamma \lambda, m) ) \mod 2$$
$$= 0 + \sum\limits_{i=k+2}^\infty x( r^i(\gamma \lambda, m) - 1) \cdot \mu( [ L^{i+1} (\gamma \lambda, m) ]^{-1} L^i (\gamma \lambda, m), r^{i+1} (\gamma \lambda, m) ) \mod 2$$
$$= 0 + \sum\limits_{i=k+2}^\infty y(r^i(\gamma \lambda, m) - 1) \cdot \mu( [ L^{i+1} (\gamma \lambda, m) ]^{-1} L^i (\gamma \lambda, m), r^{i+1} (\gamma \lambda, m) ) \mod 2$$
$$= \sum\limits_{i=k+1}^\infty y(r^i(\gamma \lambda, m) - 1) \cdot \mu( [ L^{i+1} (\gamma \lambda, m) ]^{-1} L^i (\gamma \lambda, m), r^{i+1} (\gamma \lambda, m) ) \mod 2,$$
and similarly for $\psi$ in place of $\gamma$. Therefore all terms after the $(k+1)^\text{st}$ term in the sums defining $e(y)(g)$ and $e(y)(h)$ are respectively equal. We conclude that $e(y)(g) = e(y)(h)$.
\end{proof}

The above theorem has two immediate corollaries. We point out that on the space of all subflows of $k^G$ we use the Vietoris topology (see Section \ref{SECT BASIC TC}), or equivalently the topology induced by the Hausdorff metric. In symbolic and topological dynamics there is a lot of interest in finding invariants, and in particular searching for complete invariants, for topological conjugacy, particularly for subflows of Bernoulli flows over $\ZZ$ or $\ZZ^n$. The following corollary says that, up to the use of Borel functions, there are no complete invariants for the topological conjugacy relation on any Bernoulli flow.

\begin{cor}
Let $G$ be a countably infinite group and let $k > 1$ be an integer. Then there is no Borel function defined on the space of subflows of $k^G$ which computes a complete invariant for any of the equivalence relations $\TC$, $\TCF$, $\TCM$, or $\TCMF$. Similarly, there is no Borel function on $k^G$ which computes a complete invariant for the equivalence relation $\TCP$.
\end{cor}

The above theorem and corollary imply that from the viewpoint of Borel equivalence relations, the topological conjugacy relation on subflows of a common Bernoulli flow is quite complicated as no Borel function can provide a complete invariant. However, the above results do not rule out the possibility of the existence of algorithms for computing complete invariants among subflows described by finitary data, such as subflows of finite type.

The above theorem also leads to another nice corollary. We do not know if the truth of the following corollary was previously known.

\begin{cor}
For every countably infinite group $G$, there are uncountably many pairwise non-topologically conjugate free and minimal continuous actions of $G$ on compact metric spaces.
\end{cor}

\section{Topological conjugacy of free subflows} \label{SECT TC FREE}

In this section we present a complete classification of the complexity of both $\TC(G)$ and $\TCF(G)$ for every countably infinite group $G$. We show that for a countably infinite group $G$, the equivalence relations $\TC(G)$ and $\TCF(G)$ are both Borel bi-reducible with $E_0$ if $G$ is locally finite and are both universal countable Borel equivalence relations if $G$ is not locally finite. In particular, by Lemma \ref{LEM EQ} we have that for every countably infinite locally finite group $G$, all of the equivalence relations $\TCP(G)$, $\TCMF(G)$, $\TCM(G)$, $\TCF(G)$, and $\TC(G)$ are Borel bi-reducible with $E_0$. We remind the reader the definition of a locally finite group. 

\begin{definition} \index{locally finite group}
A group $G$ is \emph{locally finite} if every finite subset of $G$ generates a finite subgroup.
\end{definition}

We first consider locally finite groups. The main theorem of the previous section allows us to quickly classify the associated equivalence relations. We need the following simple lemma.

\begin{lem}
Let $G$ be a countable group, and let $f, g: 2^G \rightarrow 2^G$ be functions induced by the block codes $\hat{f}: 2^H \rightarrow 2$ and $\hat{g}: 2^K \rightarrow 2$, respectively. Then $f \circ g$ is induced by a block code on $HK$.
\end{lem}

\begin{proof}
Clearly, $f \circ g$ is continuous and commutes with the shift action of $G$. So by Theorem \ref{THM BC}, $f \circ g$ is induced by a block code. It therefore suffices to show that if $x, y \in 2^G$ agree on $H K$ then $[f \circ g (x)](1_G) = [f \circ g (y)](1_G)$. So fix $x, y \in 2^G$ with $x \res H K = y \res H K$. Then for each $h \in H$ we have $(h^{-1} \cdot x)\res K = (h^{-1} \cdot y)\res K$. Therefore for $h \in H$
$$g(x)(h) = \hat{g}((h^{-1} \cdot x)\res K) = \hat{g}((h^{-1} \cdot y)\res K) = g(y)(h).$$
So $g(x)\res H = g(y)\res H$ and therefore $f(g(x))(1_G) = f(g(y))(1_G)$.
\end{proof}

\begin{theorem}
Let $G$ be a countably infinite, locally finite group. Then $\TC(G)$, $\TCF(G)$, $\TCM(G)$, $\TCMF(G)$, and $\TCP(G)$ are all Borel bi-reducible with $E_0$. In particular, these equivalence relations are nonsmooth and hyperfinite.
\end{theorem}

\begin{proof}
By Theorem \ref{THM MF} we have that $E_0$ Borel embeds into each of the equivalence relations. By Lemma \ref{LEM EQ}, it will suffice to show that both $\TC(G)$ and $\TCP(G)$ are hyperfinite since it is well known that hyperfinite equivalence relations Borel reduce to $E_0$ (\cite{DJK}). We remind the reader that a Borel equivalence relation is hyperfinite if it is the increasing union of finite Borel equivalence relations.

Since $G$ is locally finite, we can find an increasing sequence, $(H_n)_{n \in \N}$, of finite subgroups of $G$ whose union is $G$. For each $n \in \N$, define $E_n \subseteq \Su(G) \times \Su(G)$ by the rule: $A \ E_n \ B$ if and only if there is a conjugacy $\phi$ between $A$ and $B$ for which both $\phi$ and $\phi^{-1}$ are induced by block codes on $H_n$. Then $E_n$ is an equivalence relation as transitivity follows from the previous lemma (since $H_n$ is a subgroup of $G$). Also the proof of Proposition \ref{TC BOREL EQREL} immediately shows that each equivalence relation $E_n$ is Borel. Since $\bigcup_{n \in \N} H_n = G$, we have that $\TC(G) = \bigcup_{n \in \N} E_n$. Now we use the fact that $G$ is locally finite. Each $H_n$ is finite, so there are only finitely many block codes on $H_n$ and hence each equivalence relation $E_n$ is finite. We conclude that $\TC(G)$ is hyperfinite, and therefore Borel reducible to $E_0$. A similar argument shows that $\TCP(G)$ is hyperfinite as well.
\end{proof}

The proof of the previous theorem seems quite simple, but one should not overlook the fact that it relies on Theorem \ref{THM MF}. The authors do not know if there is a simpler proof of Theorem \ref{THM MF} in the context of locally finite groups.

Now we change our focus to nonlocally finite groups. We prove that for countably infinite nonlocally finite groups $G$ the equivalence relations $\TC(G)$ and $\TCF(G)$ are universal countable Borel equivalence relations. Unfortunately, we are unable to classify $\TCM(G)$, $\TCMF(G)$, and $\TCP(G)$ for nonlocally finite groups.

The authors' original interest in studying the complexity of the topological conjugacy relation, $\TC(G)$, stemmed from the following theorem of John Clemens.

\begin{theorem}[\cite{JC}]
$\TC(\ZZ^n)$ is a universal countable Borel equivalence relation for every $n \geq 1$.
\end{theorem}

Our proof roughly follows Clemens' proof for $\ZZ$. However, substantial additions and changes to his proof are required since we want to both extend his result to all nonlocally finite groups and extend it from $\TC$ to $\TCF$. One of the crucial components of our proof is constructing elements of $2^G$ which mimic the behavior of elements of $2^\ZZ$. The following lemma is a small step towards this construction. After this lemma are three more lemmas followed by the main theorem.

In this section, for $x \in 2^\ZZ$ we let $-x$ denote the element of $2^\ZZ$ defined by $-x(n) = x(-n)$ for all $n \in \ZZ$. Clearly $x \in 2^\ZZ$ is a $2$-coloring if and only if $-x$ is a $2$-coloring.

\begin{lem} \label{MIRRORTH}
There is a $2$-coloring $\pi \in 2^\ZZ$ for which $\pi$ and $-\pi$ are orthogonal.
\end{lem}

\begin{proof}
Let $c$ be any $2$-coloring on $\ZZ$. Define
$$\pi(n) = \begin{cases}
1 & \text{if } n \equiv 0 \mod 8 \\
1 & \text{if } n \equiv 1 \mod 8 \\
0 & \text{if } n \equiv 2 \mod 8 \\
1 & \text{if } n \equiv 3 \mod 8 \\
0 & \text{if } n \equiv 4 \mod 8 \\
c(m) & \text{if } n = 8m + 5 \\
0 & \text{if } n \equiv 6 \mod 8 \\
0 & \text{if } n \equiv 7 \mod 8
\end{cases}$$
Then $\pi$ is a $2$-coloring since it clearly blocks $8n$ for all $n \in \ZZ$ (see Corollary \ref{cor:blockinglemma}). Let $g_1, g_2 \in \ZZ$ and set $T = \{0, 1, 2, \ldots, 10\}$. Clearly there is $0 \leq n \leq 7$ with $g_1 + n \equiv 0 \mod 8$ and hence $\pi(g_1+n) = \pi(g_1+n+1) = 1$. If $-\pi(g_2+n) \neq 1$ or $-\pi(g_2+n+1) \neq 1$ then we are done since $n, n+1 \in T$. Otherwise we must have that $-g_2 - n \equiv 1 \mod 8$. It follows that $-g_2-n-3 \equiv 6 \mod 8$ and thus $\pi(g_1+n+3) = 1 \neq 0 = -\pi(g_2+n+3)$. Since $n+3 \in T$, this completes the proof that $\pi$ and $-\pi$ are orthogonal.
\end{proof}

The following is a technical lemma which will be needed briefly for a very specific purpose in the proof of the main theorem.

\begin{lem} \label{LEM SN}
Let $X$ be a compact metric space, let $\ZZ$ act continuously on $X$, let $y \in X$ be minimal, and let $d \in \N$. Let $(\xi_n)_{n \in \N}$ and $(\nu_n)_{n \in \N}$ be sequences of functions from $\N$ to $\N$. Suppose that each such function is monotone increasing and tends to infinity. Then there exists an increasing function $s: \N \rightarrow \N_+$ such that $y = \lim s(n) \cdot y$ and for $n, n', k \in \N$ we have the implication:
$$-2d + \nu_k((4d+1)(s(k+1) - s(k) - 8) - 2d)$$
$$\leq (4d+1)|s(n) - s(n')| \leq$$
$$2d + \xi_k((4d+1)(s(k+1)-s(k) + 8) + 2d)$$
implies $\max(n, n') = k+1$.
\end{lem}

\begin{proof}
Let $\rho$ be the metric on $X$. We claim that if $\epsilon > 0$ and $n \in \N$, then there is $k \geq n$ with $\rho(k \cdot y, y) < \epsilon$. To see this, let $z$ be any limit point of $\{k \cdot y \: k \geq n\}$ (a limit point must exist since $X$ is compact). By minimality of $y$, we must have $y \in \overline{[z]}$. In particular, there is $m \in \ZZ$ with $\rho(m \cdot z, y) < \frac{\epsilon}{2}$. Since $\ZZ$ acts continuously on $X$, $m \cdot z$ is a limit point of $\{(m+k) \cdot y \: k \geq n\}$. So there is $k \geq n$ with $m + k \geq n$ and $\rho((m + k) \cdot y, m \cdot z) < \frac{\epsilon}{2}$. Then $\rho((m + k) \cdot y, y) < \epsilon$ and $m + k \geq n$, completing the proof of the claim.

Fix a sequence $(\epsilon_n)_{n \in \N}$ of positive real numbers tending to $0$. We will choose a function $s: \N \rightarrow \N_+$ which will have the additional property that $\rho(s(n) \cdot y, y) < \epsilon_n$ for all $n \in \N$. Pick $s(0) > 0$ with $\rho(s(0) \cdot y, y) < \epsilon_0$. Let $t_1 \geq s(0)$ be such that
$$-2d + \nu_0((4d+1)(t_1 - 8) - 2d) > 0$$
($t_1$ exists since $\nu_0$ tends to infinity). Pick $s(1) > s(0) + t_1$ with $\rho(s(1) \cdot y, y) < \epsilon_1$. Suppose that $s(0), s(1), \cdots, s(n-1)$ have been defined and satisfy all of the required properties. Let $m \in \N$ be the maximal element of the union
$$\{2d + \xi_k((4d + 1)(s(k+1) - s(k) + 8) + 2d \: k + 1 < n\}$$
$$\bigcup \{(4d + 1)|s(k') - s(k)| \: k, k' < n\}.$$
Let $t_n \geq m$ be such that
$$-2d + \nu_{n-1}((4d+1)(t_n - 8) - 2d) > m > 0$$
($t_n$ exists since $\nu_{n-1}$ tends to infinity). Now pick $s(n) > s(n-1) + t_n$ with $\rho(s(n) \cdot y, y) < \epsilon_n$. This defines the function $s: \N \rightarrow \N_+$.

Clearly we have $y = \lim s(n) \cdot y$. Let $n, n', k \in \N$ satisfy $\max(n, n') \neq k+1$. We must show that either
$$-2d + \nu_k((4d+1)(s(k+1) - s(k) - 8) - 2d) > (4d+1)|s(n) - s(n')|$$
or
$$(4d+1)|s(n) - s(n')| > 2d + \xi_k((4d+1)(s(k+1)-s(k) + 8) + 2d).$$
By swapping $n$ and $n'$, we may suppose that $n \geq n'$. If $n = n'$ then we are done since
$$(4d+1)|s(n) - s(n')| = 0 < -2d + \nu_k((4d+1)(t_{k+1} - 8) - 2d) \leq$$
$$-2d + \nu_k((4d+1)(s(k+1) - s(k) - 8) - 2d),$$
where the second inequality follows from $\nu_k$ being monotone increasing. If $n > k+1$ then by our construction we have
$$(4d+1)|s(n) - s(n')| \geq (4d+1)(s(n) - s(n-1)) \geq s(n) - s(n-1)$$
$$> t_n \geq 2d + \xi_k((4d+1)(s(k+1) - s(k) + 8) + 2d).$$
Finally, if $n < k+1$ then $n' < k+1$ and
$$(4d+1)|s(n) - s(n')| < -2d + \nu_k ((4d+1)(t_{k+1} - 8) - 2d)$$
$$\leq -2d + \nu_k((4d+1)(s(k+1) - s(k) - 8) - 2d),$$
where again the second inequality follows from $\nu_k$ being monotone increasing. This completes the proof.
\end{proof}

Let $G$ be a finitely generated countably infinite group. Call a set $S \subseteq G$ symmetric if $S = S^{-1}$. Let $S$ be a finite symmetric set which generates $G$. The (right) Cayley graph of $G$ with respect to $S$, $\Gamma_S$, is the graph with vertex set $G$ and edge set $\{(g, gs) \: g \in G, s \in S\}$. Since $S$ generates $G$, this graph is connected. Also, $G$ acts on $\Gamma_S$ by multiplication on the left and this action is by automorphisms of $\Gamma_S$. This is the only action of $G$ on $\Gamma_S$ which we will discuss. We define a metric, $\rho_S$, on $G = \mathrm{V}(\Gamma_S)$ by setting $\rho_S(g, h)$ equal to the length of the shortest path joining $g$ and $h$ in $\Gamma_S$. This metric is left-invariant, meaning that $\rho_S(tg, th) = \rho_S(g, h)$ for all $t, g, h \in G$. In particular, $G$ acts on $\Gamma_S$ by isometries. We will call $\rho_S$ the left-invariant word length metric associated to $S$. For $g, h \in G$, we let $[g, h]_S$ denote the set of shortest paths $P: \{0, 1, \ldots, \rho_S(g,h)\} \rightarrow \mathrm{V}(\Gamma_S)$ which begin at $g$ and end at $h$. Notice that for $t, g, h \in G$ we have that $t \cdot [g, h]_S = [tg, th]_S$, where $(t \cdot P)(n) = t \cdot P(n)$ for paths $P$.

\begin{lem} \label{PATH}
Let $G$ be a countably infinite group generated by a finite symmetric set $S$. Let $\rho_S$ be the left-invariant word-length metric associated to $S$, and let $d \geq 1$. Then there is a bi-infinite sequence $P: \ZZ \rightarrow G$ such that $\rho_S(P(n), P(k)) = d|n-k|$ for all $n, k \in \ZZ$.
\end{lem}

\begin{proof}
Let $\Gamma_S$ be the (right) Cayley graph of $G$ with respect to $S$. Since $G = \bigcup_{n \in \N} S^n$ is infinite and $S$ is finite, we must have that $S^{n+1} \not\subseteq S^n$ for all $n \in \N$. For every $n \geq 1$, pick $g_n \in S^{2n} - S^{2n-1}$ and $Q_n \in [1_G, g_n]_S$. Notice that $\rho_S(1_G, g_n) = 2n$ and therefore $\dom(Q_n) = \{0, 1, \ldots, 2n\}$. Also notice that $\rho_S(Q_n(k_1), Q_n(k_2)) = |k_1 - k_2|$ whenever $k_1, k_2 \in \dom(Q_n)$. For $n \geq 1$ define $P_n: \{-n, -n+1, \ldots, n\} \rightarrow \mathrm{V}(\Gamma_S)$ by setting
$$P_n(k) = Q_n(n)^{-1} \cdot Q_n(k+n).$$
Clearly each $P_n$ is a path in $\Gamma_S$ and furthermore by the left-invariance of $\rho_S$ we have that $\rho_S(P_n(k_1), P_n(k_2)) = |k_1 - k_2|$ whenever $k_1, k_2 \in \dom(P_n)$. Clearly $P_n(0) = 1_G$, and since $P_n$ is a path in $\Gamma_S$ we must have that $P_n(k) \in S^k$ for all $k \in \dom(P_n)$. Since $S^k$ is a finite set, there is a subsequence $(P_{n(i)})_{i \in \N}$ such that for all $k \in \ZZ$ the sequence of group elements $(P_{n(i)}(k))_{i \in \N}$ is eventually constant. Define $\tilde{P}: \ZZ \rightarrow G$ by letting $\tilde{P}(k)$ be the eventual value of $(P_{n(i)}(k))_{i \in \N}$. Clearly $\tilde{P}(0) = 1_G$ and $\rho_S(\tilde{P}(k_1), \tilde{P}(k_2)) = |k_1 - k_2|$ for all $k_1, k_2 \in \ZZ$. The proof is complete after defining $P: \ZZ \rightarrow G$ by $P(n) = \tilde{P}(dn)$.
\end{proof}

We let $\F = \langle a, b \rangle$ denote the nonabelian free group on the generators $a$ and $b$. Recall that $E_\infty$ denotes the equivalence relation on $2^{\F}$ given by $x \ E_\infty \ y \Leftrightarrow [x] = [y]$. $E_\infty$ is a universal countable Borel equivalence relation, or in other words, it is the most complicated countable Borel equivalence relation.

We introduce some terminology which will be helpful in the next lemma. If $g \in \F$ is not the identity element, then the \emph{reduced word representation} of $g$ is the unique ordered tuple $(s_1, s_2, \ldots, s_n)$ where $g = s_1 \cdot s_2 \cdots s_n$, each $s_i \in \{a, a^{-1}, b, b^{-1}\}$, and $s_i \neq s_{i-1}^{-1}$ for $1 < i \leq n$. If $s \in \{a, a^{-1}, b, b^{-1}\}$ then we say that $g \in \F$ \emph{begins with $s$} if $g$ is not the identity element and the first member of the reduced word representation of $g$ is $s$. We call the nonidentity elements of $\langle a \rangle \cup \langle b \rangle$ \emph{segments}. A segment is \emph{even} if the length of its reduced word representation is even and is otherwise called \emph{odd}. For $s \in \{a, a^{-1}, b, b^{-1}\}$ we say a segment $g$ is \emph{of type $s$} if there is $n \geq 1$ with $g = s^n$. For nonidentity $g \in \F$, the \emph{segment representation of $g$} is the unique ordered tuple ($s_1, s_2, \ldots, s_n)$ where $g = s_1 \cdot s_2 \cdots s_n$, each $s_i$ is a segment, and for $1 < i \leq n$ the type of $s_i$ is neither the type of $s_{i-1}$ nor the inverse of the type of $s_{i-1}$. For example, if $g = a^3 b^{-2} a$ then the reduced word representation of $g$ is
$$(a, a, a, b^{-1}, b^{-1}, a)$$
and the segment representation of $g$ is
$$(a^3, b^{-2}, a).$$
For nonidentity $g \in \F$, the \emph{segments of $g$} are the members of the segment representation of $g$, and the \emph{$n^\text{th}$ segment of $g$} is the $n^\text{th}$ member of the segment representation of $g$.

The following lemma is due to John Clemens. We include a proof for completeness.

\begin{lem}[\cite{JC}] \label{LEM UNIVERSAL}
There is a Borel set $J \subseteq 2^{\F}$ which is invariant under the action of $\F$ and satisfies:
\begin{enumerate}
\item[\rm (i)] $E_\infty \sqsubseteq_B E_\infty \res J$, so $E_\infty \res J$ is a universal countable Borel equivalence relation;
\item[\rm (ii)] For $x,y \in J$ with $\neg (x \ E_\infty \ y)$ there are infinitely many $g \in \F$ with $x(g) \neq y(g)$;
\item[\rm (iii)] For every $x \in J$ there are infinitely many $g \in \F$ with $x(g) = 1$.
\end{enumerate}
\end{lem}

\begin{proof}
Let $H \subseteq \F$ be the subgroup generated by $a^2$ and $b^2$. Let $\phi: \F \rightarrow H$ be the isomorphism induced by $\phi(a) = a^2$ and $\phi(b) = b^2$. For $x \in 2^{\F}$ define $f(x) \in 2^{\F}$ by
$$f(x)(w) = \begin{cases}
x(u) & \text{if } w = \phi(u) \text{ for some } u \text{ or } w = \phi(u)abv \text{ for some } u \text{ and } v \\
 & \  \text{ with } v \text{ not beginning with } b^{-1}; \\
1 & \text{otherwise.}
\end{cases}$$
Then $f$ is a continuous injection, so the image of $f$ is Borel. Let $J = \bigcup_{g \in \F} g \cdot f(2^{\F})$. Clearly $J$ is Borel. Clause (iii) is immediately satisfied.

Suppose $x, y \in 2^{\F}$ satisfy $x \ E_\infty \ y$. Then $y = g \cdot x$ for some $g \in \F$ and it is easy to check that $f(y) = f(g \cdot x) = \phi(g) \cdot f(x)$. Thus $f(x) \ E_\infty \ f(y)$. Pick any $x, y \in 2^{\F}$ and $g \in \F$. To complete the proof, it suffices to show that if $f(y)$ and $g \cdot f(x)$ agree at all but finitely many coordinates then $[y] = [x]$ (and hence $[f(x)] = [f(y)]$). If $f(y)$ has value $1$ at all but finitely many coordinates, then $f(y)$ is identically $1$, as are $f(x)$, $x$, and $y$ and hence $[y] = [x]$. So we may suppose that there is $k \in \F$ with $y(k) = 0$ and hence $f(y)(\phi(k)) = f(y)(\phi(k) a b v) = 0$ for all $v \in \F$ which do not begin with $b^{-1}$. Let $t, h \in \F$ be such that $g = h \phi(t)$ and such that the reduced word representation of $h$ does not end in $aa$, $a^{-1} a^{-1}$, $b b$, or $b^{-1} b^{-1}$. Then
$$f(y) =^* g \cdot f(x) = h \cdot f(t \cdot x) = h \cdot f(x')$$
($=^*$ denotes equality at all but finitely many coordinates) where $x' = t \cdot x$. If $h = 1_{\F}$, then we are done since $f$ is injective. Towards a contradiction, suppose that $h \neq 1_{\F}$.

If $h^{-1} \phi(k) a b = b^{2m}$ for $m \in \ZZ$ (possibly $m = 0$) then set $s_1 = b$ and $s_2 = a^{-1}$. Otherwise let $s_1 \in \{a, a^{-1}\}$ be such that $s_1^{-1}$ is not the last element in the reduced word representation of $h^{-1} \phi(k) a b$, and let $s_2 = b$. We must have that for some $n \geq 1$
$$0 = f(y)(\phi(k) a b (s_1 s_2)^n) = [h \cdot f(x')](\phi(k) a b (s_1 s_2)^n) = f(x')(h^{-1} \phi(k) a b (s_1 s_2)^n).$$
By the definition of $f$, there must be $p, v \in \F$ with $v$ not beginning with $b^{-1}$ and
$$h^{-1} \phi(k) a b (s_1 s_2)^n = \phi(p) a b v \text{ or } h^{-1} \phi(k) a b (s_1 s_2)^n = \phi(p).$$
By choosing a larger value of $n$ if necessary, we can assume $h^{-1} \phi(k) a b (s_1 s_2)^n = \phi(p) a b v$. Notice that since the first segment of $h^{-1}$ is odd, the first segment of $h^{-1} \phi(k)$ must also be odd. If the initial segment of $h^{-1} \phi(k) a b$ is not odd, then $h^{-1} \phi(k)$ must have at most two segments and $h^{-1} \phi(k) a b$ must be either $a^{2m} b$ or $b^{2m}$ for some $m \in \ZZ$ (possibly $m = 0$). In $\phi(p) a b v$, the first odd segment of type $b$ is preceded by an odd segment of type $a$ or $a^{-1}$. So for $m \in \ZZ$ we have
$$a^{2m} b a^{\pm 1} s_2 (s_1 s_2)^{n-1} \neq \phi(p) a b v \neq b^{2m} b a^{-1} (s_1 s_2)^{n-1}.$$
Therefore $a^{2m} b \neq h^{-1} \phi(k) a b \neq b^{2m}$ (recall the definition of $s_1$ and $s_2$). Thus the initial segment of $h^{-1} \phi(k) a b$ must be odd. So the initial segment of $\phi(p) a b v$ must be odd and thus we must have that $\phi(p) = a^{2m}$ for some $m \in \ZZ$. Then the initial segment of $h^{-1} \phi(k) a b$ must be of type $a$ or $a^{-1}$ and hence the initial segment of $h^{-1} \phi(k)$ is of type $a$ or $a^{-1}$. We cannot have $h^{-1} \phi(k) \in \langle a \rangle$ as otherwise the initial segment of $h^{-1} \phi(k) a b$ would be even. So $h^{-1} \phi(k)$ has at least two segments and the first segment is of type $a$ or $a^{-1}$. Since
$$h^{-1} \phi(k) a b (s_1 s_2)^n = a^{2m+1} b v$$
the second segment of $h^{-1} \phi(k)$ must be of type $b$ (as opposed to being of type $b^{-1}$).

Let $t_1 \in \{b, b^{-1}\}$ be such that $t_1^{-1}$ is not the last element of the reduced word representation of $h^{-1} \phi(k)$. Set $t_2 = a$. Then there is $N \in \N$ with
$$f(y)(\phi(k) (t_1 t_2)^N) = [h \cdot f(x')](\phi(k) (t_1 t_2)^N).$$
Clearly $\phi(k) (t_1 t_2)^N$ is not in the image of $\phi$. Also, the first odd segment of $\phi(k) (t_1 t_2)^N$ is of type $b$ or $b^{-1}$, so there cannot exist $k', v' \in \F$ with $v'$ not beginning with $b^{-1}$ and $\phi(k) (t_1 t_2)^N = \phi(k') a b v'$. By the definition of $f$ we have $f(y)(\phi(k) (t_1 t_2)^N) = 1$. However, there is $v' \in \F$ not beginning with $b^{-1}$ with
$$h^{-1} \phi(k) (t_1 t_2)^N = \phi(p) a b v'$$
(where $p$ is the same as in the last paragraph). We have
$$[h \cdot f(x')](\phi(k) (t_1 t_2)^N) = f(x')(\phi(p) a b v') = f(x')(\phi(p) a b v) = 0.$$
This is a contradiction. We conclude $h = 1_{\F}$.
\end{proof}

We are now ready for the final theorem of this chapter. This theorem states that $\TC(G)$ and $\TCF(G)$ are universal countable Borel equivalence relations when $G$ is not locally finite. We mention that John Clemens claims to have an independent proof of this theorem, however as of yet he has not made his proof public.

To give a very rough outline of the proof, we will construct elements of $2^G$ which have behavior very similar to elements of $2^\ZZ$ and then we will adapt and implement Clemens' proof of this result for $\TC(\ZZ)$.

\begin{theorem}
Let $G$ be a countably infinite, nonlocally finite group. Then $\TC(G)$ and $\TCF(G)$ are universal countable Borel equivalence relations.
\end{theorem}

\begin{proof}
Since $G$ is not locally finite, there is a finite symmetric $1_G \in S \subseteq G$ with $\langle S \rangle$ infinite. Set $p_1(k) = 16 \cdot (2k^4+1)$ and for $n > 1$ set $p_n(k) = 2k^4+1$. Then $(p_n)_{n \geq 1}$ is a sequence of functions of subexponential growth. By Corollary \ref{GROW BP} there is a centered, directed, and maximally disjoint blueprint $(\Delta_n, F_n)_{n \in \N}$ with $|\Lambda_n| \geq \log_2 \ p_n(|F_n|)$ for all $n \geq 1$. It is easy to see from the proof of Corollary \ref{GROW BP} that the blueprint can be chosen to have the additional property that $F_1 \subseteq \langle S \rangle$. Recall that we are free to fix a choice of distinct $\alpha_n, \beta_n, \gamma_n \in D^n_{n-1}$ for all $n \geq 1$. For $n \geq 1$ set $\gamma_n = 1_G$ and let $\alpha_n, \beta_n \in D^n_{n-1} - \{1_G\}$ be arbitrary but distinct. By clause (viii) of Lemma \ref{STRONG BP LIST}, we have $\bigcap_{n \in \N} \Delta_n a_n = \bigcap_{n \in \N} \Delta_n b_n = \varnothing$. Apply Theorem \ref{FM} to get a function $c$ which is canonical with respect to this blueprint. By Proposition \ref{CANONICAL DMIN} $c$ is $\Delta$-minimal. Apply Corollary \ref{GEN MINCOL} to get a fundamental and $\Delta$-minimal $c'$ with $|\Theta_1(c')| \geq \log_2 (16) = 4$ and with the property that every element of $2^G$ extending $c'$ is a $2$-coloring. A trivial application of Lemma \ref{MIN GRAPH} (with $\mu$ identically $0$), Lemma \ref{LEM UNION}, and Lemma \ref{MIN UNION} gives us a fundamental and $\Delta$-minimal $x \in 2^{\subseteq G}$ extending $c'$ and with the property that $|\Theta_1(x)| = 4$ and $\Theta_n(x) = \varnothing$ for $n > 1$. Let $\theta_1, \theta_2, \theta_3, \theta_4$ be the distinct elements of $\Theta_1(x)$. From this point forward $\Theta_1$ will denote $\Theta_1(x)$.

Equip $\langle S \rangle$ with the left-invariant word length metric $\rho$ induced by the generating set $S$. We define a norm on $\langle S \rangle$ and on finite subsets of $\langle S \rangle$ by
$$\|g\| = \rho(g, 1_G), \text{ and}$$
$$\|A\| = \max\{\rho(a, 1_G) \: a \in A\}$$
for $g \in \langle S \rangle$ and finite $A \subseteq \langle S \rangle$. Notice that if $A, B \subseteq \langle S \rangle$ are finite then $\|A B\| \leq \|A\| + \|B\|$. Set
$$d = \|F_1 F_1^{-1}\|.$$
This expression is meaningful since $F_1 \subseteq \langle S \rangle$. Apply Lemma \ref{PATH} with respect to the number $4d+1$ to get $P: \ZZ \rightarrow \langle S \rangle \subseteq G$. Recall that $1_G \in S$ and therefore for $m \in \N$
$$\bigcup_{i = 0}^m S^i = S^m.$$
We claim that for every $n \in \ZZ$
$$\{k \in \ZZ \: P(k) F_1 F_1^{-1} \cap P(n) F_1 F_1^{-1} F_1 F_1^{-1} S^{4d+1} F_1 F_1^{-1} \neq \varnothing\} = \{n-1, n, n+1\}.$$
Indeed, if $k \in \ZZ$ is in the set on the left then
$$P(n)^{-1} P(k) \in F_1 F_1^{-1} F_1 F_1^{-1} S^{4d+1} F_1 F_1^{-1} F_1 F_1^{-1}$$
so
$$|k - n| (4d+1) = \|P(k)^{-1} P(n)\| \leq d + d + (4d+1) + d + d = 8d+1.$$
Thus $|k - n| \leq 1$. On the other hand, if $k \in \ZZ$ and $|k - n| \leq 1$ then $P(k) \in P(n) S^{4d+1}$ since $\|P(n)^{-1} P(k)\| = |k-n|(4d+1) \leq 4d+1$. For each $n \in \ZZ$, pick $Q(n) \in \Delta_1$ with $Q(n) \in P(n) F_1 F_1^{-1}$. Note that $Q(n) \neq Q(k)$ for $n \neq k$. Also, we have that for all $n \in \ZZ$
$$\{k \in \ZZ \: Q(k) \cap Q(n) F_1 F_1^{-1} S^{4d+1} F_1 F_1^{-1} \neq \varnothing\} = \{n-1, n, n+1\}.$$

Let $\pi \in 2^\ZZ$ be a $2$-coloring with $\pi$ and $-\pi$ orthogonal (see Lemma \ref{MIRRORTH}). Define $x' \in 2^G$ by
$$x'(g) = \begin{cases}
x(g) & \text{if } g \in \dom(x) \\
\pi(n) & \text{if } g = Q(n) \theta_2 \\
1 & \text{if } g \in Q(\ZZ) \theta_3 \\
0 & \text{otherwise} \\
\end{cases}.$$
Note that $x'$ admits a $Q(\ZZ)$ membership test. Specifically,
$$g \in Q(\ZZ) \Longleftrightarrow g \in \Delta_1 \text{ and } x'(g \theta_3) = 1.$$
Set
$$X = \{ w \in \overline{[x']} \: 1_G \in \Delta_1^w \text{ and } w(\theta_3) = 1\}.$$
Notice that $X \cap [x'] = Q(\ZZ)^{-1} \cdot x'$. It follows from the definition of $\Delta_1^w$ (see Section \ref{SECT SUBFLOW FUND FNCTN}) that $X$ is a clopen subset of $\overline{[x']}$. Therefore, if $w = \lim g_n \cdot x' \in X$ then $g_n \cdot x' \in X$ and $g_n^{-1} \in Q(\ZZ)$ for all but finitely many $n \in \N$. So approximating $w \in X$ by points in $X \cap [x'] = Q(\ZZ)^{-1} \cdot x'$ gives
$$|F_1 F_1^{-1} S^{4d+1} F_1 F_1^{-1} \cdot w \cap X| = 3.$$
Notice that $w$ is in the set on the left. We put a graph structure, $\Gamma$, on $X$ as follows: $w, w' \in X$ are adjacent in $\Gamma$ if and only if $w \neq w'$ and
$$w' \in F_1 F_1^{-1} S^{4d+1} F_1 F_1^{-1} \cdot w.$$
Then every element of $X$ has degree precisely $2$ in $\Gamma$. Note that for $k \in \ZZ$ the vertices adjacent to $Q(k)^{-1} \cdot x'$ are precisely $Q(k-1)^{-1} \cdot x'$ and $Q(k+1)^{-1} \cdot x'$. Again, approximating points in $X$ by points in $X \cap [x'] = Q(\ZZ)^{-1} \cdot x'$ shows that every connected component of $\Gamma$ is infinite and in particular is isomorphic to the standard Cayley graph of $\ZZ$.

We claim that for $w \in X$ the connected component of $\Gamma$ containing $w$ is precisely $[w] \cap X$. Clearly the connected component of $\Gamma$ containing $w$ is contained in $[w] \cap X$. We now show the opposite inclusion. So suppose that $h \in G$ and $h \cdot w \in X$. We can approximate $w$ by points in $[x'] \cap X$ to find $g \in G$ such that $g \cdot x', h g \cdot x' \in X$. Since $[x'] \cap X \subseteq \langle S \rangle \cdot x'$ and $x'$ has trivial stabilizer, we have that $g, g h \in \langle S \rangle$ and thus $h \in \langle S \rangle$. Let $n \in \N$ be such that $\|h\| \leq n (4d + 1) + 2d$. Again, approximating $w$ by points in $[x'] \cap X$, we find $g' \in G$ such that $t g' \cdot x' \in X \Leftrightarrow t \cdot w \in X$ for every $t \in S^{n (4d + 1) + 2d}$. So $g' \cdot x', h g' \cdot x' \in X$ and thus there are $-n \leq m \leq n$ and $k \in \ZZ$ with $g' = Q(k)^{-1}$ and $h = Q(k+m)^{-1} Q(k)$. So for $|i| \leq |m|$ we have $Q(k+i)^{-1} Q(k) g' \cdot x' \in X$ and thus $Q(k+i)^{-1} Q(k) \cdot w \in X$. Therefore there is a path in $\Gamma$ from $w = Q(k)^{-1} Q(k) \cdot w$ to $h \cdot w = Q(k+m)^{-1} Q(k) \cdot w$. We conclude that $[w] \cap X$ is precisely the connected component of $\Gamma$ containing $w$. We point out for future reference that this argument showed that if $w \in X$ and $h \cdot w \in X$ then $h \in Q(\ZZ)^{-1} Q(\ZZ)$.

We now define an action, $*$, of $\ZZ$ on $X$. Let $M \in \N$ be such that for all $n_1, n_2 \in \ZZ$ there is $-M \leq t \leq M$ with $\pi(n_1+t) \neq \pi(n_2-t)$ ($M$ exists since $\pi$ is orthogonal to $-\pi$). Fix $w \in X$. Set $w_0 = w$ and let $w_{-M}, w_{-M+1} \ldots, w_{M-1}, w_M \in X$ be the vertices which can be joined to $w_0$ in $\Gamma$ by a path of length at most $M$. Rearranging the indices if necessary, we may assume that $(w_i, w_{i+1}) \in \mathrm{E}(\Gamma)$ for each $-M \leq i < M$. By approximating $w_0$ by elements of $[x'] \cap X$, we see that there is $n \in \ZZ$ with either $w_i(\theta_2) = \pi(n+i)$ for all $-M \leq i \leq M$ or $w_i(\theta_2) = -\pi(-n+i) = \pi(n-i)$ for all $-M \leq i \leq M$. By the definition of $M$, one of these two possibilities must fail for all $n \in \ZZ$ so in particular the two possibilities cannot be simultaneously satisfied. By swapping $w_i$ and $w_{-i}$ for each $-M \leq i \leq M$ if necessary, we may assume that there is $n \in \ZZ$ with $w_i(\theta_2) = \pi(n+i)$ for all $-M \leq i \leq M$. We define $0*w_0 = w_0$, $1*w_0 = w_1$, and $-1*w_0 = w_{-1}$. In general, recursively define $k*w_0 = 1*((k-1)*w_0)$ and $-k*w_0 = -1*((-k+1)*w_0)$ for $k > 1$ and $w_0 \in X$.

We claim this action of $\ZZ$ on $X$ is continuous. Since $X$ is clopen in $\overline{[x']}$, there is a finite $1_G \in B \subseteq G$ and $V \subseteq 2^B$ such that for $w' \in \overline{[x']}$ we have $w' \in X$ if and only if $w' \res B \in V$. By approximating $w = w_0 \in X$ by elements of $[x'] \cap X$ and recalling the definitions of $P$ and $Q$ we see that
$$w_{-M}, w_{-M+1}, \ldots, w_M \in F_1 F_1^{-1} S^{M(4d+1)} F_1 F_1^{-1} \cdot w_0 \subseteq S^{M(4d+1)+2d} \cdot w_0.$$
Let $w' \in X$ and let $w_{-M}', \ldots, w_M'$ be defined similarly to before, with $w_0' = w'$. For $-M \leq i \leq M$ let $g_i \in G$ be such that $w_i = g_i \cdot w_0$. If 
$$w_0' \res S^{M(4d+1)+2d} B \{1_G, \theta_2\} = w_0 \res S^{M(4d+1)+2d} B \{1_G, \theta_2\}$$
then $(f \cdot w_0') \res B = (f \cdot w_0) \res B$ for all $f \in S^{M(4d+1)+2d}$ and hence $w_i' = g_i \cdot w_0'$ and $w_i'(\theta_2) = w_i(\theta_2)$ for all $-M \leq i \leq M$. So if $w'$ and $w$ satisfy the displayed expression above and if $w'$ is sufficiently close to $w$, then $1*w' = g_1 \cdot w'$ is close to $1*w = g_1 \cdot w$ since the action of $G$ on $2^G$ is continuous. Similarly $-1*w'$ is close to $-1*w$. We conclude that the action of $\ZZ$ on $X$ is continuous. In fact, since $X$ is compact, for each $k \in \ZZ$ the map $w \mapsto k*w$ is uniformly continuous.

$X$ is a clopen subset of $\overline{[x']}$ and is therefore compact and Hausdorff. Therefore, there is $y \in X$ which is minimal with respect to the action of $\ZZ$ (Lemma \ref{EXIST MIN}). For notational simplicity for the rest of the proof, we redefine $\pi \in 2^\ZZ$ by
$$\pi(n) = (n*y)(\theta_2).$$
Notice that this new $\pi$ is in the closure of the orbit of the old $\pi$, the new $\pi$ is a $2$-coloring, and it is orthogonal to its reflection $-\pi$.

Fix an increasing sequence $(C_n)_{n \in \N}$ of finite subsets of $G$ with $1_G \in C_0$, $C_n = C_n^{-1}$, and $\bigcup_{n \in \N} C_n = G$. For $n \in \N$ define $\xi_n: \N \rightarrow \N$ by
$$\xi_n(m) = \|C_n^{-1} S^m C_n \cap \langle S \rangle\|.$$
Then for each $n \in \N$ the function $\xi_n$ is monotone increasing, tends to infinity, and $m \leq \xi_n(m) \leq \xi_{n+1}(m)$ for all $m \in \N$. The functions $\xi_n$ may not map bijectively onto $\N$, but we define functions $\nu_n$ which behave similar to $\xi_n^{-1}$. For $n, k \in \N$ we define
$$\nu_n(k) = \min \{m \in \N \: \xi_n(m) \geq k\}.$$
We again have that for each $n \in \N$ the function $\nu_n$ is monotone increasing, tends to infinity, and $\nu_{n+1}(k) \leq \nu_n(k) \leq k$.

Apply Lemma \ref{LEM SN} to get an increasing function $s: \N \rightarrow \N_+$ such that $\lim s(n) * y = y$ and for all $n, n', k \in \N$ the following implication holds:
$$-2d + \nu_k((4d+1)(s(k+1) - s(k) - 8) - 2d)$$
$$\leq (4d+1)|s(n) - s(n')| \leq$$
$$2d + \xi_k((4d+1)(s(k+1)-s(k) + 8) + 2d)$$
implies $\max(n, n') = k+1$. The reader is discouraged from thinking too much about the technical condition above. The technical requirement on the function $s$ is needed briefly for a very specific purpose near the end of the proof. Aside from this, we will only make use of the fact that $\lim s(n)*y = y$ and that $s(n) > 0$ for all $n \in \N$.

For $n \in \ZZ$ let $q_n \in G$ be the unique group element with $q_n^{-1} \cdot y = n * y$. Note that $q_0 = 1_G$, $\{q_n \: n \in \ZZ\} \subseteq \Delta_1^y \subseteq \langle S \rangle$, and for $n, k \in \ZZ$
$$(4d+1)|n-k| - 2d \leq \| q_n^{-1} q_k \| \leq (4d+1)|n-k| + 2d.$$
Although it is a bit of a misnomer, we will refer to $(q_n)_{n \in \ZZ}$ as a bi-infinite path, and for $t \in \ZZ$ we will refer to $(q_n)_{n \geq t}$ as a right-infinite path and $(q_n)_{n \leq t}$ as a left-infinite path. The point $y \in 2^G$ in some sense mimics $\pi \in 2^\ZZ$. The rest of the proof will proceed by working with certain elements of $2^G$ which agree with $y$ on all coordinates not in $\Delta_1^y \Theta_1$. Find a sequence $(g_n)_{n \in \N}$ in $G$ with $y = \lim g_n \cdot x'$ and set $z = \lim g_n \cdot x$. Note that the limit exists,$1_G \in \Delta_1^z = \Delta_1^y$, and for all $g \not\in \Delta_1^y \Theta_1$ $z(g) = y(g)$. Since $z \in \overline{[x]}$, $z$ is $\Delta$-minimal, every element of $2^G$ extending $z$ is a $2$-coloring, and $G - \dom(z) = \Delta_1^z \Theta_1$. We also have the useful property that $\lim q_{s(n)}^{-1} \cdot z = z$. From this point forward we will work with $y$ and extensions of $z$ and therefore no longer need $x$ or $x'$.

Let $\F$ denote the nonabelian free group with two generators $a$ and $b$. Let $J \subseteq 2^\F$ be as referred to in Lemma \ref{LEM UNIVERSAL}. Define $c: \F \rightarrow 2^{<\N}$ by setting
$$\begin{array}{ccc}
c(1_\F) & = & 11000011, \\
c(a) & = & 11100011, \\
c(a^{-1}) & = & 11010011, \\
c(b) & = & 11001011, \\
c(b^{-1}) & = & 11000111,
\end{array}$$
and
$$c(g) = c(e_1) ^\frown c(e_2) ^\frown \cdots ^\frown c(e_n)$$
where $^\frown$ denotes concatenation, $g \in \F - \{1_\F\}$, $e_i \in \{a, a^{-1}, b, b^{-1}\}$, and $g = e_1 \cdot e_2 \cdots e_n$ is the unique reduced word representation of $g$. Note that $c(g)$ has length $8$ times as long as the length of the reduced word representation of $g$ (for $g \neq 1_\F$).

\begin{figure}[ht]\label{FIG FUIK}
\begin{center}
\setlength{\unitlength}{3mm}
\begin{picture}(30,45)(5,-.5)

\put(0,10){
\put(5,27){\makebox(0,0)[b]{$G$}}
\qbezier(0,10)(0,30)(30,30)
\qbezier(0,10)(0,0)(15,0)
\qbezier(15,0)(40,0)(40,20)
\qbezier(30,30)(40,30)(40,20)

\qbezier(15,15)(17,14)(18,16)
\qbezier(18,16)(19,18)(21,17)
\qbezier(21,17)(23,16)(24,18)
\qbezier(24,18)(25,20)(27,19)
\qbezier(27,19)(30,18)(32,20)
\qbezier(32,20)(34,22)(36,21)
\qbezier(36,21)(39,20)(40,21)

\qbezier(15,16)(17,15)(18,17)
\qbezier(18,17)(19,19)(21,18)
\qbezier(21,18)(23,17)(24,19)
\qbezier(24,19)(25,21)(27,20)

\put(15,15){\line(0,1){1}}
\put(27,19){\line(0,1){1}}

\put(15,15){\circle*{.3}}
\put(16,13.5){\makebox(0,0)[b]{$q_{s(k)}$}}
\put(22.5,16.7){\circle*{.3}}
\put(23,15.4){\makebox(0,0)[b]{$q_n$}}
\put(29,13){\makebox(0,0)[b]{$q_n\theta_1\mapsto c(h_i)(n-s(k))$}}
\put(26.3,11.5){\makebox(0,0)[b]{$q_n\theta_2\mapsto \pi(n)$}}
\put(25.4,10){\makebox(0,0)[b]{$q_n\theta_3\mapsto 1$}}
\put(5,13){\circle*{0.3}}
\put(5.5,11.5){\makebox(0,0)[b]{$1_G$}}
\put(5,14){\circle*{.3}}
\put(5.5,14.5){\makebox(0,0)[b]{$\theta_4$}}
\put(8.5,16){\makebox(0,0)[b]{$\theta_4\mapsto (h_i\cdot u)(h_k)$}}
\put(10,21){\circle*{.3}}
\put(13.5,20.5){\makebox(0,0)[b]{$g\in\mbox{dom}(z)$}}
\put(13,22){\makebox(0,0)[b]{$g\mapsto z(g)$}}

\begin{dashjoin}{0.15}
\jput(15,15){}
\jput(14.5,15.4){}
\jput(14,15.5){}
\jput(13.5,15.5){}
\jput(13,15.4){}
\jput(12.5,15.2){}
\jput(12,14.8){}
\jput(11.5,14.4){}
\jput(11,14.1){}
\jput(10.5,13.7){}
\jput(10,13.4){}
\jput(9.5,13.3){}
\jput(9,13.2){}
\jput(8.5,13,2){}
\jput(8,13,3){}
\jput(7.5,13.4){}
\jput(7,13.5){}
\jput(6.5,13.5){}
\jput(6,13.4){}
\jput(5.5,13.2){}
\jput(5,13){}
\jput(4.5,12.6){}
\jput(4,12.2){}
\jput(3.5, 11.6){}
\jput(3,11.3){}
\jput(2.5,11){}
\jput(2, 11){}
\jput(1.5, 11.1){}
\jput(1,11.2){}
\jput(.5,11.3){}
\jput(0,11.4){}

\end{dashjoin}

\put(15,15.8){\line(1,1){.1}}
\put(15,15.6){\line(1,1){.3}}
\put(15,15.4){\line(1,1){.43}}
\put(15,15.2){\line(1,1){.6}}
\put(15,15){\line(1,1){.72}}

\put(15.2,15){\line(1,1){.7}}
\put(15.4,15){\line(1,1){.66}}
\put(15.6,15){\line(1,1){.66}}
\put(15.8,15){\line(1,1){.66}}
\put(16,15){\line(1,1){.74}}

\put(15.2,15){\line(-1,-1){.1}}
\put(15.4,15){\line(-1,-1){.14}}
\put(15.6,15){\line(-1,-1){.19}}
\put(15.8,15){\line(-1,-1){.24}}
\put(16,15){\line(-1,-1){.26}}

\put(16.2,15){\line(1,1){.85}}
\put(16.4,15){\line(1,1){3.32}}
\put(16.6,15){\line(1,1){3.32}}
\put(16.8,15){\line(1,1){1.5}}
\put(17,15){\line(1,1){1.0}}
\put(17.2,15){\line(1,1){.6}}

\put(16.2,15){\line(-1,-1){.3}}
\put(16.4,15){\line(-1,-1){.3}}
\put(16.6,15){\line(-1,-1){.3}}
\put(16.8,15){\line(-1,-1){.28}}
\put(17,15){\line(-1,-1){.22}}

\put(18,17){\line(1,1){1.26}}
\put(17.8,16.6){\line(1,1){1.7}}

\put(18.8,17){\line(1,1){1.33}}
\put(19.3,17.3){\line(1,1){1.}}
\put(19.5,17.3){\line(1,1){.95}}
\put(19.73,17.33){\line(1,1){.88}}
\put(19.92,17.32){\line(1,1){.82}}

\put(20.1,17.3){\line(1,1){0.78}}
\put(20.26,17.26){\line(1,1){0.72}}
\put(20.42,17.22){\line(1,1){0.72}}
\put(20.6,17.2){\line(1,1){0.69}}
\put(20.76,17.16){\line(1,1){0.68}}

\put(20.86,17.06){\line(1,1){0.7}}
\put(21,17){\line(1,1){0.72}}

\put(6,2){
\put(15,15.8){\line(1,1){.1}}
\put(15,15.6){\line(1,1){.3}}
\put(15,15.4){\line(1,1){.43}}
\put(15,15.2){\line(1,1){.6}}
\put(15,15){\line(1,1){.72}}

\put(15.2,15){\line(1,1){.7}}
\put(15.4,15){\line(1,1){.66}}
\put(15.6,15){\line(1,1){.66}}
\put(15.8,15){\line(1,1){.66}}
\put(16,15){\line(1,1){.74}}

\put(15.2,15){\line(-1,-1){.1}}
\put(15.4,15){\line(-1,-1){.14}}
\put(15.6,15){\line(-1,-1){.19}}
\put(15.8,15){\line(-1,-1){.24}}
\put(16,15){\line(-1,-1){.26}}

\put(16.2,15){\line(1,1){.85}}
\put(16.4,15){\line(1,1){3.32}}
\put(16.6,15){\line(1,1){3.32}}
\put(16.8,15){\line(1,1){1.5}}
\put(17,15){\line(1,1){1.0}}
\put(17.2,15){\line(1,1){.6}}

\put(16.2,15){\line(-1,-1){.3}}
\put(16.4,15){\line(-1,-1){.3}}
\put(16.6,15){\line(-1,-1){.3}}
\put(16.8,15){\line(-1,-1){.28}}
\put(17,15){\line(-1,-1){.22}}

\put(18,17){\line(1,1){1.26}}
\put(17.8,16.6){\line(1,1){1.7}}

\put(18.8,17){\line(1,1){1.33}}
\put(19.3,17.3){\line(1,1){1.}}
\put(19.5,17.3){\line(1,1){.95}}
\put(19.73,17.33){\line(1,1){.88}}
\put(19.92,17.32){\line(1,1){.82}}

\put(20.1,17.3){\line(1,1){0.78}}
\put(20.26,17.26){\line(1,1){0.72}}

\put(20.42,17.22){\line(1,1){0.6}}
\put(20.6,17.2){\line(1,1){0.38}}
\put(20.76,17.16){\line(1,1){0.22}}
\put(20.86,17.06){\line(1,1){0.12}}

}
}

\put(14,6.5){
\begin{dashjoin}{0.15}
\jput(0,0){}
\jput(.5,0){}
\jput(1,0){}
\jput(1.5,0){}
\jput(2,0){}
\jput(2.5,0){}
\jput(3,0){}
\jput(3.5,0){}
\jput(4,0){}
\end{dashjoin}
\put(9,-0.5){\makebox(0,0)[b]{invisible path}}
}

\put(14,5){\line(1,0){4}}
\put(23.8,4.5){\makebox(0,0)[b]{visible path data}}

\put(14,4){\line(1,0){4}}
\put(14,4){\line(0,-1){1}}
\put(14,3){\line(1,0){4}}
\put(18,3){\line(0,1){1}}
\multiput(14,3)(.2,0){16}{\line(1,1){1}}
\put(14,3.2){\line(1,1){.8}}
\put(14,3.4){\line(1,1){.6}}
\put(14,3.6){\line(1,1){.4}}
\put(14,3.8){\line(1,1){.2}}
\put(17.2,3){\line(1,1){.8}}
\put(17.4,3){\line(1,1){.6}}
\put(17.6,3){\line(1,1){.4}}
\put(17.8,3){\line(1,1){.2}}
\put(25,3){\makebox(0,0)[b]{$c$ data and path data}}

\end{picture}
\caption{\label{fig:FF} The values of $f(u,i,k)$ for $n \geq s(k)$ and for $g \in \mbox{dom}(z)$.}
\end{center}
\end{figure}
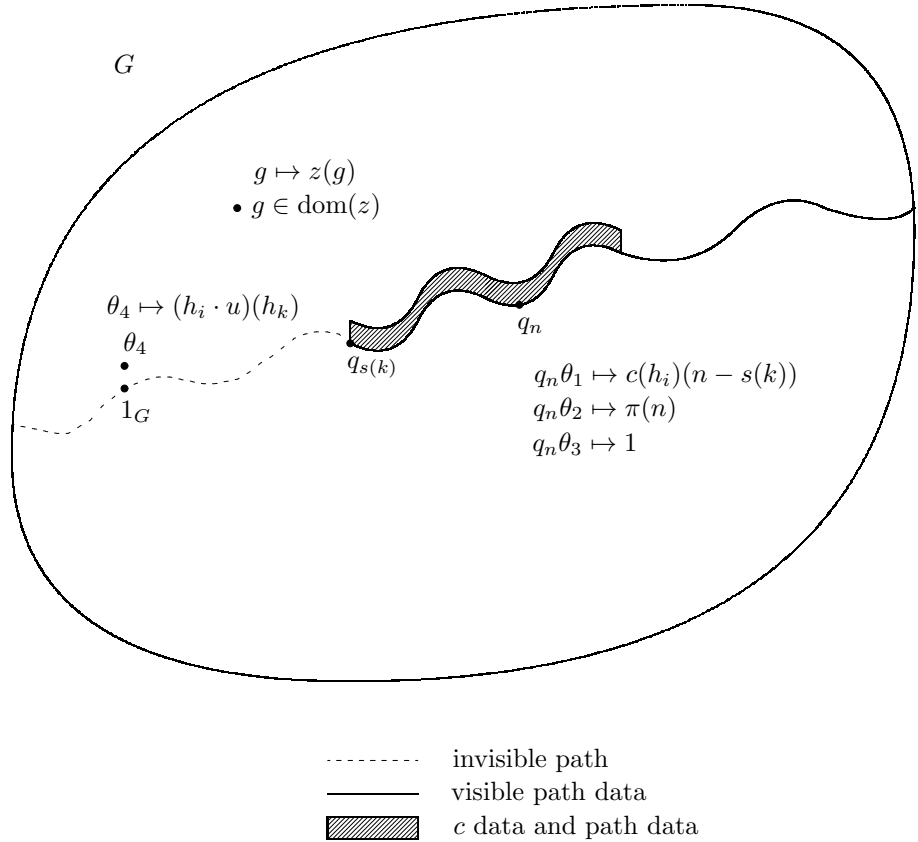

Let $(h_i)_{i \in \N}$ be an enumeration of $\F$ with $h_0 = 1_\F$. For $i, k \in \N$ and $u \in J$ define
$$f(u, i, k)(g) = \begin{cases}
z(g) & \text{if } g \in \dom(z) \\
c(h_i)(n - s(k)) & \text{if } n \geq s(k) \wedge n - s(k) \in \dom(c(h_i)) \wedge g = q_n \theta_1 \\
\pi(n) & \text{if } n \geq s(k) \wedge g = q_n \theta_2 \\
1 & \text{if } n \geq s(k) \wedge g = q_n \theta_3 \\
(h_i \cdot u)(h_k) & \text{if } g = \theta_4 \\
0 & \text{otherwise} \\
\end{cases}.$$
The values of $f(u,i,k)$ are shown in Figure \ref{FIG FUIK}. Basically, $f(u,i,k)$ extends $z$ and has special values at $\theta_4$ and along the right-infinite path $(q_n)_{n \geq s(k)}$ and is zero elsewhere. Along the path $(q_n)_{n \geq s(k)}$, $\theta_1$ is used to record $c(h_i)$, $\theta_2$ is used to record $\pi$, and $\theta_3$ provides a simple $\{q_n \: n \geq s(k)\}$ membership test. The single point $\theta_4$ is used to record $(h_i \cdot u)(h_k)$. For $u \in J$ we define the subflow $A(u)$ to be
$$A(u) = \overline{G \cdot \{f(u,i,k) \: i, k \in \N\} }.$$
By clause (ii) of Proposition \ref{EXTENSION FREE}, $A(u)$ is a free subflow of $2^G$ for every $u \in J$. To complete the proof, it suffices, by Lemma \ref{LEM EQ} and clause (i) of Lemma \ref{LEM UNIVERSAL}, to show that $A$ is a Borel reduction from $E_\infty \res J$ to $\TCF(G)$. In other words, we must show that $A: J \rightarrow \SuF(G)$ is a Borel function and
$$u \ E_\infty \ v \Longleftrightarrow \ A(u) \ \TCF(G) \ A(v)$$
for all $u, v \in J$. Recall that by definition $u \ E_\infty \ v$ if and only if $[u] = [v]$.

We first prove that $A: J \rightarrow \SuF(G)$ is Borel. Recall that the topology on $\SuF(G)$ is generated by the subbasic open sets
$$\{K \in \SuF(G) \: K \subseteq U\} \text{ and } \{K \in \SuF(G) \: K \cap U \neq \varnothing\}$$
where $U$ varies over the open subsets of $2^G$ (see Section \ref{SECT BASIC TC}). Let $U \subseteq 2^G$ be open. Temporarily define $f_{i,k}(u) = f(u,i,k)$. Then each $f_{i,k}: J \rightarrow 2^G$ is continuous. Since $U$ is open we have $A(u) \cap U \neq \varnothing$ if and only if there are $i, k \in \N$ with $[f(u,i,k)] \cap U \neq \varnothing$, or equivalently $f(u,i,k) \in \bigcup_{g \in G} g \cdot U$. Therefore
$$A^{-1}(\{K \in \SuF(G) \: K \cap U \neq \varnothing\}) = \bigcup_{i, k \in \N} f_{i,k}^{-1} \left( \bigcup_{g \in G} g \cdot U \right)$$
which is Borel since each $f_{i,k}$ is continuous. Define $U_n = \{x \in 2^G \: d(x, 2^G - U) \geq 1/n\}$. Then each $U_n$ is closed and we have that $A(u) \subseteq U$ if and only if there is $n \geq 1$ so that for all $i, k \in \N$ $[f(u,i,k)] \subseteq U_n$ (this follows from the compactness of $A(u)$). The condition $[f(u,i,k)] \subseteq U_n$ is equivalent to $f(u,i,k) \in \bigcap_{g \in G} g \cdot U_n$. Therefore
$$A^{-1}(\{K \in \SuF(G) \: K \subseteq U\}) = \bigcup_{n \geq 1} \bigcap_{i, k \in \N} f_{i,k}^{-1} \left( \bigcap_{g \in G} g \cdot U_n \right)$$
which is Borel since each $f_{i,k}$ is continuous. We conclude $A: J \rightarrow \SuF(G)$ is Borel.

For $u \in J$ we define $X(u) \subseteq A(u)$ similar to how we defined $X$ before. Specifically, $w \in A(u)$ is an element of $X(u)$ if and only if $1_G \in \Delta_1^w$ and $w(\theta_3) = 1$. We point out that $X(u)$ is a clopen subset of $A(u)$. Notice that for $i, k \in \N$ the elements of $[f(u,i,k)] \cap X(u)$ are precisely the points $q_n^{-1} \cdot f(u,i,k)$ for $n \geq s(k)$. We put a graph structure, $\Gamma(u)$, on $X(u)$ just as before. There is an edge between $w, w' \in X(u)$ if and only if $w \neq w'$ and
$$w' \in F_1 F_1^{-1} S^{4d+1} F_1 F_1^{-1} \cdot w.$$
It is clear from the definition of $A(u)$ that every element of $X(u)$ has degree either $1$ or $2$ in $\Gamma(u)$. As before, by approximating $w \in X(u)$ by elements of $X(u) \cap G \cdot \{f(u,i,k) \: i, k \in \N\}$ we see that every connected component of $\Gamma(u)$ is infinite and that for $w \in X(u)$ the connected component of $\Gamma(u)$ containing $w$ is precisely $[w] \cap X(u)$. Since some vertices in $\Gamma(u)$ have degree $1$, we cannot define an action of $\ZZ$ on $X(u)$ similar to before. We can however define an action of $\N$ on $X(u)$ which we again denote by $*$. If $w \in X(u)$, then $w$ lies in an infinite connected component of $\Gamma(u)$, so there are $w_0, w_1, \ldots, w_{2M} \in X(u)$ (where $M$ is as before) with $(w_i, w_{i+1}) \in \mathrm{E}(\Gamma(u))$ for each $0 \leq i < 2M$ and with $w = w_k$ for some $0 \leq k \leq 2M$. By re-indexing the $w_i$'s if necessary, we may assume that there is $n \in \ZZ$ with $w_i(\theta_2) = \pi(n+i)$ for all $0 \leq i \leq 2M$. If $k < 2M$ then we define $1*w = 1*w_k = w_{k+1}$. If $k = 2M$ then an approximation by elements of $X(u) \cap G \cdot \{f(u,i,k) \: i, k \in \N\}$ shows that $w$ has degree $2$ in $\Gamma(u)$ and hence there is $w_{2M+1} \neq w_{2M-1}$ which is adjacent to $w_k = w_{2M}$. In this case, we define $1*w = w_{2M+1}$. In general, define $m*w = 1*((m-1)*w)$ and $0*w = w$. This action is well defined due to the properties of $M$ and $\pi$. The same argument as for the action of $\ZZ$ on $X$ shows that the action of $\N$ on $X(u)$ is continuous, and hence for each $k \in \N$ the map $w \mapsto k*w$ is uniformly continuous. In fact, that previous argument shows something stronger which we will need. There is a finite set $B \subseteq G$ such that if $t \geq 2M$ and $w, w' \in X(u)$ satisfy
$$w \res S^{t(4d+1)+2d} B \Theta_1 = w' \res S^{t(4d+1)+2d} B \Theta_1$$
then for $g \in G$ and $0 \leq i \leq t$ we have
$$i*w = g \cdot w \Longleftrightarrow i*w' = g \cdot w',$$
$$i*(g \cdot w) = w \Longleftrightarrow i*(g \cdot w') = w',$$
and if $g \cdot w \in X(u)$ and $i*(g \cdot w) = w$ or $i*w = g \cdot w$ then $(g \cdot w) \res \Theta_1 = (g \cdot w') \res \Theta_1$.

We will now prove that if $u, v \in J$ and $u \ E_\infty \ v$ then $A(u) \ \TCF(G) \ A(v)$. So fix $u, v \in J$ with $u \ E_\infty \ v$, or in other words $[u] = [v]$. Let $g \in \F$ be such that $u = g \cdot v$. By considering the reduced word representation of $g$ in the generators $a$ and $b$ and using the fact that $\TCF(G)$, being an equivalence relation, is transitive and symmetric, we see that it suffices to consider the cases $u = a \cdot v$ and $u = b \cdot v$. We will treat the case $u = a \cdot v$ as the case $u = b \cdot v$ is nearly identical. We must show that $A(u) \ \TCF(G) \ A(v)$. Define the permutation $\sigma: \N \rightarrow \N$ by setting $\sigma(i) = j$ if $h_i a = h_j$. Let $\phi$ be the function sending $g \cdot f(u,i,k)$ to $g \cdot f(v, \sigma(i),k)$ for each $g \in G$ and $i, k \in \N$. The function $\phi$ is well defined as it is easy to check that $[f(u,i,k)] \neq [f(u,j,m)]$ for $i,j,k,m \in \N$ with $(i,k) \neq (j,m)$ (since $f(u,i,k)$ and $f(u,j,m)$ extend $z$, one can apply clause (i) of Proposition \ref{EXTENSION FREE}). Notice that
$$f(u,i,k)(\theta_4) = (h_i \cdot u)(h_k) = (h_i a \cdot v)(h_k) = (h_{\sigma(i)} \cdot v)(h_k) = f(v,\sigma(i),k)(\theta_4).$$
Therefore $\phi$ is only changing the $c(h_i)$ data in $g \cdot f(u,i,k)$ to the $c(h_{\sigma(i)})$ data in $g \cdot f(v,\sigma(i),k)$. Notice that in order to change $c(h_i)$ to $c(h_{\sigma(i)})$ one only needs to either append $c(a)$ (if $c(h_i)$ does not end with $c(a^{-1})$ and $h_i \neq 1_{\F}$), change the last $8$ digits of $c(h_i)$ (if $h_i = a^{-1}$ or $h_i = 1_{\F}$), or delete the last $8$ digits of $c(h_i)$ (if $c(h_i)$ ends with $c(a^{-1})$ and $h \neq a^{-1}$). Therefore if $\ell_i$ denotes the length of $c(h_i)$ then $g \cdot f(u,i,k)$ and $g \cdot f(v,\sigma(i),k)$ agree on
$$G - g \{q_n \: s(k) + \ell_i - 8 \leq n \leq s(k) + \ell_i + 8\} \theta_1.$$
To show that $\phi$ is induced by a block code, it suffices to show that the map $w \mapsto \phi(w)(1_G)$ is uniformly continuous for $w \in G \cdot \{f(u,i,k) \: i, k \in \N\}$. If $\phi(w)(1_G) \neq w(1_G)$, then there must be $i, k, n \in \N$ with $n \geq s(k)$ and $w = (q_n \theta_1)^{-1} \cdot f(u,i,k)$, in which case $w \in \theta_1^{-1} \cdot X(u)$ as $q_n^{-1} \cdot f(u,i,k) \in X(u)$. So for $w$ not in  $\theta_1^{-1} \cdot X(u)$ we have $\phi(w)(1_G) = w(1_G)$, and thus the map $w \mapsto \phi(w)(1_G)$ is uniformly continuous outside of $\theta_1^{-1} \cdot X(u)$. Since $\theta_1^{-1} \cdot X(u)$ is clopen, it suffices to show that the map $w \mapsto \phi(w)(1_G)$ is uniformly continuous on
$$(\theta_1^{-1} \cdot X(u)) \cap (G \cdot \{f(u,i,k) \: i, k \in \N\}).$$
But on this set the map $w \mapsto \phi(w)(1_G)$ is the composition of the maps $w \mapsto \theta_1 \cdot w$ (with domain the set above) and $g \cdot f(u,i,k) \mapsto (g \cdot f(v,\sigma(i),k))(\theta_1)$ (with domain $X(u) \cap G \cdot \{f(u,i,k) \: i, k \in \N\}$). The first map is clearly uniformly continuous, and the uniform continuity of the second map follows from our discussion on how to change $c(h_i)$ to $c(h_{\sigma(i)})$ and from the final remark of the previous paragraph. Therefore $\phi$ is induced by a block code and so extends to a continuous function $\phi: A(u) \rightarrow A(v)$ which commutes with the action of $G$. The set $\phi(A(u))$ is compact, hence closed, and contains a dense subset of $A(v)$. So $\phi(A(u)) = A(v)$. By considering the block code for $\phi$, it is easy to see that $\phi$ is injective (alternatively, we could have just as easily shown that the inverse map $\phi^{-1}: G \cdot \{f(v,i,k) \: i, k \in \N\} \rightarrow G \cdot \{f(u,i,k) \: i, k \in \N\}$ is induced by a block code and so extends to $\phi^{-1}: A(v) \rightarrow A(u)$). We conclude that $A(u) \ \TCF(G) \ A(v)$.

Before proving $A(u) \ \TCF(G) \ A(v)$ implies $u \ E_\infty \ v$, we first have to shed more light on the properties of $A(u)$. In order to better understand the subflow $A(u)$, we make a few more definitions. For $i \in \N$ define
$$f(i)(g) = \begin{cases}
z(g) & \text{if } g \in \dom(z) \\
c(h_i)(n) & \text{if } n \geq 0 \wedge n \in \dom(c(h_i)) \wedge g = q_n \theta_1 \\
\pi(n) & \text{if } n \geq 0 \wedge g = q_n \theta_2 \\
1 & \text{if } n \geq 0 \wedge g = q_n \theta_3 \\
0 & \text{otherwise} \\
\end{cases}.$$
The function $f(i)$ is very similar to $f(u,i,k)$, but differs in two ways. First, $f(i)$ has special values along the path $(q_n)_{n \geq 0}$, while $f(u,i,k)$ has special values along the path $(q_n)_{n \geq s(k)}$. Second, $f(i)(\theta_4) = 0$ while $f(u,i,k)(\theta_4) = (h_i \cdot u)(h_k)$. We define $z_0$ to be the extension of $z$ which is zero everywhere $z$ is not defined, and $z_1$ by
$$z_1(g) = \begin{cases}
z(g) & \text{if } g \in \dom(z) \\
1 & \text{if } g = \theta_4 \\
0 & \text{otherwise} \\
\end{cases}.$$
We call the points in $X(u)$ \emph{path data points}. If $w \in X(u)$ and there are $0 \leq i, j \leq 8$ and $w' \in X(u)$ with $w'(\theta_1) = 1$, $i * w' = w$, and $(j*w)(\theta_1) = 1$ then we call $w$ a \emph{$c$ data point} (here the ``c'' is in reference to the function $c:\F \rightarrow 2^{<\N}$). The final remark of the paragraph in which $X(u)$ was defined shows that the set of $c$ data points in $A(u)$ is a clopen subset of $A(u)$. By approximating $w \in X(u)$ by elements of $G \cdot \{f(u,i,k) \: i, k \in \N\}$ we see that the set of $c$ data points in $[w] \cap X(u)$ is connected in $\Gamma(u)$. The action of $\N$ on $X(u)$ gives us a well defined notion of right and left directions in $X(u)$. Namely, if $w \in X(u)$ then $1*w$ is \emph{to the right of $w$}, and if $w$ has degree $2$ in $\Gamma(u)$ then the unique $w' \in X(u)$ with $1*w' = w$ is \emph{to the left of $w$}. With this convention, calling a connected subset of $\Gamma(u)$ \emph{left-infinite}, \emph{right-infinite}, and \emph{bi-infinite} has the obvious meaning. We use left-infinite and right-infinite in the strict sense, meaning that any set which is left-infinite or right-infinite cannot be bi-infinite. We say $x \in A(u)$ \emph{has infinite, right-infinite, left-infinite, or bi-infinite path data} (\emph{$c$ data}), if $[x] \cap X(u)$ (respectively the $c$ data points in $[x] \cap X(u)$) is infinite, right-infinite, left-infinite, or bi-infinite respectively. We say $x \in A(u)$ \emph{has finite path data} (\emph{finite $c$ data}) if $[x] \cap X(u)$ (respectively the $c$ data points in $[x] \cap X(u)$) is nonempty and finite. Finally we say $x \in A(u)$ \emph{has no path data} (\emph{no $c$ data}) if $[x] \cap X(u)$ (respectively the set of $c$ data points in $[x] \cap X(u)$) is empty.

For $u \in J$, we define:
\begin{enumerate}
\item[\rm $A_1(u) =$] $\bigcup_{i,k \in \N} [f(u,i,k)]$;
\item[\rm $A_2(u) =$] $\bigcup_{i \in \N} [f(i)]$;
\item[\rm $A_3(u) =$] $\{w \in A(u) \: w \text{ has infinite } c \text{ data} \}$;
\item[\rm $A_4(u) =$] $\{w \in A(u) \: w \text{ has bi-infinite path data but no } c \text{ data}\}$;
\item[\rm $A_5(u) =$] $[z_1]$;
\item[\rm $A_6(u) =$] $\overline{[z_0]}$.
\end{enumerate}
Note that $A_1(u)$ through $A_5(u)$ may not be subflows of $2^G$. We spend the next few paragraphs proving the following:
\begin{enumerate}
\item[\rm (i)] $\{A_1(u), A_2(u), \ldots, A_6(u)\}$ is a partition of $A(u)$;
\item[\rm (ii)] $A_1(u)$ is the set of isolated points of $A(u)$;
\item[\rm (iii)] $A_6(u)$ is minimal;
\item[\rm (iv)] $A_4(u) \subseteq \overline{[w]}$ for all $w \in A_4(u)$;
\item[\rm (v)] $A_4(u) \cup A_6(u) = \bigcap_{i,k \in \N} \overline{[f(u,i,k)]}$.
\end{enumerate}

(i). We first check that the sets are pairwise disjoint. Notice that the points in $A_1(u) \cup A_2(u)$ have finite $c$ data, and the points in $A_5(u) \cup A_6(u)$ have no path data and no $c$ data. So we only need to show that $A_1(u) \cap A_2(u) = A_5(u) \cap A_6(u) = \varnothing$. Since $s(k) > 0$ for all $k \in \N$, for $i, j \in \N$ we have
$$f(u,i,k)(\theta_1) = 0 \neq 1 = f(j)(\theta_1)$$
and hence $f(u,i,k) \neq f(j)$. If $s \neq 1_G$, then by clause (i) of Proposition \ref{EXTENSION FREE} there is $t \in G$ with $t, s^{-1} t \in \dom(z)$ and
$$(s \cdot f(u,i,k))(t) = f(u,i,k)(s^{-1} t) = z(s^{-1} t) \neq z(t) = f(j)(t).$$
Therefore $f(j) \not\in [f(u,i,k)]$ and $[f(j)] \cap [f(u,i,k)] = \varnothing$. It follows that $A_1(u)$ and $A_2(u)$ are disjoint. Finally, we have $z_1(\theta_4) = 1$ but $w(\theta_4) = 0$ for all $w \in A_6(u)$ with $1_G \in \Delta_1^w$. So $A_5(u)$ and $A_6(u)$ are disjoint as well. Now we move on to showing that $A_r(u) \subseteq A(u)$ for $1 \leq r \leq 6$. This is clear for $r = 1$, $r = 3$, and $r = 4$. Since $z = \lim q_{s(n)}^{-1} \cdot z$, a comparison of $f(u,i,k)$ with $z$ shows that
$$f(i) = \lim q_{s(k)}^{-1} \cdot f(u,i,k).$$
Therefore $A_2(u) \subseteq A(u)$. Clause (iii) of Lemma \ref{LEM UNIVERSAL} implies that $z_1$ is an accumulation point of $(f(u,0,k))_{k \in \N}$, so $A_5(u) \subseteq A(u)$. Since
$$(4d+1)|n-m| - 2d \leq \|q_n^{-1} q_m\| \leq (4d+1)|n-m| + 2d,$$
some element of $A_6(u)$ is an accumulation point of $(q_n^{-1} \cdot f(u,0,0))_{n = -1}^{-\infty}$. So, temporarily assuming the validity of clause (iii), we have $A_6(u) \subseteq A(u)$. We do point out that $A_3(u)$ and $A_4(u)$ are nonempty. Any accumulation point of $(q_n^{-1} \cdot f(u,0,0))_{n \in \N}$ is an element of $A_4(u)$ and any accumulation point of $(q_{t(i)}^{-1} \cdot f(u,i,0))_{i \in \N}$ is an element of $A_3(u)$, where $t(i)$ is half the length of $c(h_i)$. Finally, we must show that if $w \in A(u)$ then $w \in A_r(u)$ for some $1 \leq r \leq 6$. If $w$ has infinite $c$ data then $w \in A_3(u)$. If $w$ has no $c$ data but has bi-infinite path data then $w \in A_4(u)$. If $w$ has no $c$ data and no path data then it is easy to see that $w \in A_5(u) \cup A_6(u)$. Finally, if $w$ has finite $c$ data and $w \not\in A_1(u)$ then $w$ must be in the orbit of a limit point of the form $\lim q_{s(k(n))}^{-1} \cdot f(u,i(n),k(n))$, where $k : \N \rightarrow \N$ tends to infinity. However, since this limit has finite $c$ data, $i(n)$ must eventually be constant, say of value $i$, and hence $\lim q_{s(k(n))}^{-1} \cdot f(u,i(n),k(n)) = f(i)$ (since $\lim q_{s(k(n))}^{-1} \cdot z = z$). Therefore $w \in [f(i)] \subseteq A_2(u)$.

(ii). From the definition of $A(u)$ it follows that every isolated point of $A(u)$ must lie in $A_1(u)$. Since $G$ acts on $A(u)$ by homeomorphisms, it suffices to show that $\{q_{s(k)}^{-1} \cdot f(u,i,k) \: i, k \in \N\}$ consists of isolated points of $A(u)$. Fix $i, k \in \N$ and towards a contradiction suppose there are $q_{s(k)}^{-1} \cdot f(u,i,k) \neq x_n \in A(u)$ with $q_{s(k)}^{-1} \cdot f(u,i,k) = \lim x_n$. Since $A_1(u)$ is dense in $A(u)$, we may suppose without loss of generality that $x_n \in A_1(u)$ for every $n \in \N$. Since $q_{s(k)}^{-1} \cdot f(u,i,k) \in X(u)$ is not in $1*X(u)$ (or equivalently, has no member of $X(u)$ to the left of it), we must have that $x_n \in X(u) - 1*X(u)$ for all sufficiently large $n \in \N$. For $j, m \in \N$ we have $[f(u,j,m)] \cap (X(u) - 1*X(u)) = q_{s(m)}^{-1} \cdot f(u,j,m)$. Therefore there are functions $j, m: \N \rightarrow \N$ such that $x_n = q_{s(m(n))}^{-1} \cdot f(u,j(n),m(n))$ for all sufficiently large $n \in \N$. By recalling $z$ and the definition of $s(n)$ and $q_n$, we see that $\lim q_{s(m(n))}^{-1} \cdot f(u,j(n),m(n)) = \lim f(j(n))$. So $q_{s(k)}^{-1} \cdot f(u,i,k) = \lim f(j(n))$. Clearly we must have that $j(n) = i$ for sufficiently large $n$. Therefore $q_{s(k)}^{-1} \cdot f(u,i,k) = f(i)$, contradicting the fact that $A_1(u)$ and $A_2(u)$ are disjoint.

(iii). This follows immediately from the definition of $z_0$ and the fact that $z$ is $\Delta$-minimal.

(iv). Let $w \in A_4(u)$ and let $h \in G$ be such that $ h \cdot w \in X(u)$. By definition of $A(u)$, there are $g_n \in G$ and functions $i, k: \N \rightarrow \N$ with $w = \lim g_n \cdot f(u,i(n),k(n))$. Since $X(u)$ is clopen, $h g_n \cdot f(u,i(n),k(n)) \in X(u)$ for sufficiently large $n \in \N$. So for large $n$ we have $h g_n = q_{m(n)}^{-1}$ where $m(n) \geq k(n)$ (since the limit $h \cdot w$ is a path data point). Since $w$ has no $c$ data and has bi-infinite path data, we have that $m(n) >> k(n) + 8 \|i(n)\|$ for large $n$ (here $\|i(n)\|$ denotes the reduced word length of $i(n) \in \F$). Since $y$ and $f(u,i(n),k(n))$ agree outside of a neighborhood of $\{q_t \: -\infty < t \leq k(n) + 8 \|i(n)\|\}$ we have that
$$h \cdot w = \lim q_{m(n)}^{-1} \cdot f(u,i(n),k(n)) = \lim m(n)*y.$$
So $h \cdot w \in \overline{\ZZ * y}$. If $x \in A_4(u)$ then we similarly have that there is $g \in G$ with $g \cdot x \in \overline{\ZZ * y}$. Since $y$ is minimal with respect to the action of $\ZZ$, we have $g \cdot x \in \overline{\ZZ * y} = \overline{\ZZ * (h \cdot w)}$ (the action of $\ZZ$ on $h \cdot w$ is defined since $h \cdot w \in \overline{\ZZ * y}$). Therefore $x \in \overline{[w]}$. We conclude $A_4(u) \subseteq \overline{[w]}$.

(v). Let $i,j,k,m \in \N$ with $(i,k) \neq (j,m)$. If $i \neq j \in \N$ then clearly $[f(u,i,k)] \neq [f(u,j,m)]$ since $c(h_i)$ and $c(h_j)$ are distinct. If $i = j$ and $k \neq m$ then $f(u,i,k) \neq f(u,i,m)$ and furthermore by clause (i) of Proposition \ref{EXTENSION FREE} we have $[f(u,i,k)] \neq [f(u,i,m)]$. So $(i,k) \neq (j,m)$ implies $[f(u,i,k)] \neq [f(u,j,m)]$. So by clause (ii) we have that $A_1(u)$ is disjoint from $\bigcap_{i,k \in \N} \overline{[f(u,i,k)]}$. Now fix $i, k \in \N$. Since $f(u,i,k)$ has finite $c$ data, every point of $\overline{[f(u,i,k)]} - [f(u,i,k)]$ must have no $c$ data. Also, if $\gamma \in \Delta_1^z$ and $f(u,i,k)(\gamma \theta_4) = 1$ then $\gamma = 1_G$. So we cannot have $z_1 \in \overline{[f(u,i,k)]}$. Therefore
$$\bigcap_{j,m \in \N} \overline{[f(u,j,m)]} \subseteq A_4(u) \cup A_6(u).$$
Since
$$(4d+1)|n-m| - 2d \leq \|q_n q_m^{-1}\| \leq (4d+1)|n-m| + 2d,$$
there are subsequences of $(q_n^{-1} \cdot f(u,i,k))_{n = -1}^{-\infty}$ and $(q_n^{-1} \cdot f(u,i,k))_{n = 1}^\infty$ converging to points in $A_6(u)$ and $A_4(u)$ respectively. So by clauses (iii) and (iv) $A_4(u) \cup A_6(u) \subseteq \overline{[f(u,i,k)]}$. We conclude
$$A_4(u) \cup A_6(u) = \bigcap_{i,k \in \N} \overline{[f(u,i,k)]}.$$

Now we begin the final phase of the proof. All that remains is to show that if $u, v \in J$ and $A(u) \ \TCF(G) \ A(v)$ then $u \ E_\infty \ v$ (or equivalently $[u] = [v]$). So suppose that $\phi: A(u) \rightarrow A(v)$ is a conjugacy. Our proof proceeds by verifying the following facts:
\begin{enumerate}
\item[\rm (vi)] $\phi(A_r(u)) = A_r(v)$ for all $1 \leq r \leq 6$;
\item[\rm (vii)] there is a fixed $g \in G$ with $\phi(\{f(u,i,k) \: i, k \in \N\}) = g \cdot \{f(v,i,k) \: i, k \in \N\}$;
\item[\rm (viii)] there is a finite $D \subseteq G$ such that if $w \in A(u)$ and $\phi(w)$ is a $c$ data point then there is at least one $c$ data point in $D \cdot w$;
\item[\rm (ix)] there is a fixed $j \in \N$ so that $\phi(f(u,0,k)) = g \cdot f(v,j,k)$ for all sufficiently large $k \in \N$;
\item[\rm (x)] $u \ E_\infty \ v$.
\end{enumerate}

(vi). By clause (ii) we immediately have that $\phi(A_1(u)) = A_1(v)$ and therefore by clause (v) we have $\phi(A_4(u) \cup A_6(u)) = A_4(v) \cup A_6(v)$. If $w \in A_6(u)$ and $\phi(w) \in A_6(v)$, then since $A_6(u)$ and $A_6(v)$ are minimal it would immediately follow that $\phi(A_6(u)) = A_6(v)$. On the other hand, if $w \in A_4(u)$ and $\phi(w) \in A_6(v)$ then by clause (iv) we would have $\phi(A_4(u)) \subseteq \phi(\overline{[w]}) \subseteq A_6(v)$. In this case we must have $\phi(A_4(u)) = A_6(v)$ as otherwise there is $w' \in A_6(u)$ with $\phi(w') \in A_6(v)$ and hence $\phi(A_6(u)) = A_6(v)$, contradicting the fact that $\phi$ is one-to-one. So either $A_6(v) = \phi(A_6(u))$ or $A_6(v) = \phi(A_4(u))$. The same applies to $\phi^{-1}$, so either $\phi(A_6(u)) = A_6(v)$ or else $\phi(A_6(u)) = A_4(v)$. Towards a contradiction, suppose $\phi(A_6(u)) = A_4(v)$. Then $\phi(z_0)$ has bi-infinite path data and no $c$ data. The functions $z_0$ and $z_1$ differ at only one point, so $\phi(z_0)$ and $\phi(z_1)$ must differ at finitely many points (since $\phi$ is induced by a block code). So $\phi(z_1)$ has infinite path data and hence $\phi(z_1) \not\in A_5(v) \cup A_6(v)$. By making only a finite number of changes, one cannot turn a bi-infinite path into a right-infinite path. Therefore $\phi(z_1) \not\in A_1(v) \cup A_2(v)$. Also, $\phi(z_1)$ can have at most finite $c$ data, so $\phi(z_1) \not\in A_3(v)$. Thus we must have $\phi(z_1) \in A_4(v) = \phi(A_6(u))$. This contradicts the fact that $\phi$ is one-to-one. Therefore we now know that $\phi(A_r(u)) = A_r(v)$ for $r = 1, 4, 6$. Since $z_1$ and $z_0$ differ at only finitely many places and $\phi(A_6(u)) = A_6(v)$, we must have that $\phi(A_5(u)) = A_5(v)$. Finally, we have $\phi(A_2(u) \cup A_3(u)) = A_2(v) \cup A_3(v)$. For $i,k \in \N$ $f(i)$ differs from $f(u,i,k)$ at only finitely many points. So every member of $A_2(u)$ differs at only finitely many points from some member of $A_1(u)$. Therefore we must have $\phi(A_2(u)) \subseteq A_2(v)$. However, the same argument applies to $\phi^{-1}$, so $\phi(A_2(u)) = A_2(v)$. We conclude $\phi(A_r(u)) = A_r(v)$ for all $1 \leq r \leq 6$.

(vii). Pick $i,k,j,m \in \N$. Let $g, h \in G$ and $i',k',j',m' \in \N$ be such that
$$\phi(f(u,i,k)) = g \cdot f(v,i',k') \text{ and } \phi(f(u,j,m)) = h \cdot f(v,j',m').$$
We will show that $g = h$. To simplify notation, set
$$x_1 = f(u,i,k); \ \ x_2 = f(u,j,m); \ \ y_1 = f(v,i',k'); \ \ y_2 = f(v,j',m').$$
Then $\phi(x_1) = g \cdot y_1$ and $\phi(x_2) = h \cdot y_2$. Recall that, as stated after the definition of $X(v)$, there is a finite set $B \subseteq G$ such that if $w, w' \in X(v)$ satisfy
$$w \res S^{6d+1} B \Theta_1 = w' \res S^{6d+1} B \Theta_1$$
then for every $g \in G$ we have
$$1 * w = g \cdot w \Longleftrightarrow 1 * w' = g \cdot w'.$$
Since $x_1$ and $x_2$ differ at only finitely many places, so do $g \cdot y_1$ and $h \cdot y_2$. Let $p \in \N$ be such that
$$\forall n \geq p \ (g \cdot y_1) \res g q_n S^{6d+1} B \Theta_1 = (h \cdot y_2) \res g q_n S^{6d+1} B \Theta_1.$$
Such $p$ exists since $g \cdot y_1$ and $h \cdot y_2$ differ at only finitely many coordinates and the elements $(q_n)_{n \in \N}$ are pairwise distinct.

Since $X(v)$ is clopen in $A(v)$ and $g \cdot y_1$ and $h \cdot y_2$ differ at only finitely many coordinates, the sets $g \cdot \{q_n \: n \geq s(k')\}$ and $h \cdot \{q_n \: n \geq s(m')\}$ differ by only finitely many elements. Let $t \geq \max(p, s(k'))$ be such that $g q_t \in h \cdot \{q_n \: n \geq s(m')\}$. Say $g q_t = h q_r$. Then by definition of $p$
$$(g \cdot y_1) \res g q_t S^{6d+1} B \Theta_1 = (h \cdot y_2) \res g q_t S^{6d+1} B \Theta_1 = (h \cdot y_2) \res h q_r S^{6d+1} B \Theta_1.$$
This implies
$$(q_t^{-1} \cdot y_1) \res S^{6d+1} B \Theta_1 = (q_r^{-1} \cdot y_2) \res S^{6d+1} B \Theta_1.$$
Let $f \in G$ be such that $1 * (q_t^{-1} \cdot y_1) = f \cdot (q_t^{-1} \cdot y_1)$. The equality above and the fact that $q_t^{-1} \cdot y_1, q_r^{-1} \cdot y_2 \in X(v)$ implies that $1 * (q_r^{-1} \cdot y_2) = f \cdot (q_r^{-1} \cdot y_2)$. However,
$$q_{t+1}^{-1} \cdot y_1 = 1 * (q_t^{-1} \cdot y_1) \text{ and } q_{r+1}^{-1} \cdot y_2 = 1 * (q_r^{-1} \cdot y_2).$$
So $f = q_{t+1}^{-1} q_t = q_{r+1}^{-1} q_r$. It follows that
$$g q_{t+1} = g q_t f^{-1} = h q_r f^{-1} = h q_{r+1}.$$
Since $t+1 > t \geq \max(p, s(k'))$, we can repeat this argument and conclude that $g q_{t+n} = h q_{r+n}$ for all $n \in \N$. This gives
$$\pi(t+n) = (g \cdot y_1)(g q_{t+n} \theta_2) = (h \cdot y_2)(h q_{r+n} \theta_2) = \pi(r+n)$$
for all $n \in \N$. Since $\pi$ is a $2$-coloring we must have $t = r$, so $g q_r = h q_r$, and hence $g = h$.

(viii). Let $(D_n)_{n \in \N}$ be an increasing sequence of finite subsets of $G$ with $\bigcup_{n \in \N} D_n = G$. Towards a contradiction, suppose that for every $n \in \N$ there is $w_n \in A(u)$ with $D_n \cdot w_n$ containing no $c$ data points and with $\phi(w_n)$ a $c$ data point. Let $w \in A(u)$ be an accumulation point of $(w_n)_{n \in \N}$. Since $\phi$ is continuous and the set of $c$ data points is a closed subset of $A(v)$, we have that $\phi(w)$ is a $c$ data point. If $g \in G$ then $g \cdot w$ is an accumulation point of $(g \cdot w_n)_{n \in \N}$ and there is $n \in \N$ with $g \in D_n$. For $k \geq n$ $D_n \subseteq D_k$, so $g \cdot w_k$ is not a $c$ data point. Since the $c$ data points in $A(u)$ is an open subset of $A(u)$ we have that $g \cdot w$ is not a $c$ data point. Since $g \in G$ was arbitrary, we have that $w \in A_4(u) \cup A_5(u) \cup A_6(u)$. However, $\phi(w)$ is a $c$ data point and therefore $\phi(w) \not\in A_4(v) \cup A_5(v) \cup A_6(v)$. This contradicts clause (vi). 

(ix). Let $g \in G$ be as in clause (vii) and let $D \subseteq G$ be as in clause (viii). Recall the increasing sequence $(C_n)_{n \in \N}$ of finite symmetric subsets of $G$ and the functions $(\xi_n)_{n \in \N}$ and $(\nu_n)_{n \in \N}$ used in defining the function $s: \N \rightarrow \N$. Let $n \in \N$ be such that $D \subseteq C_n$. Let $k \geq n + 2$ and let $j_0, j_1, j_2, m_0, m_1, m_2 \in \N$ be such that
$$\phi(f(u,0,k-t)) = g \cdot f(v,j_t,m_t)$$
for $0 \leq t \leq 2$. Then
$$\phi(q_{s(m_t)}^{-1} g^{-1} \cdot f(u,0,k-t)) = q_{s(m_t)}^{-1} \cdot f(v,j_t,m_t)$$
is a $c$ data point. As the $c$ data points of $[f(u,0,k-t)]$ are precisely
$$\{q_r^{-1} \: s(k-t) \leq r \leq s(k-t)+8\} \cdot f(u,0,k-t)$$
we have that by (viii)
$$\{q_r^{-1} \: s(k-t) \leq r \leq s(k-t)+8\} \cap C_n q_{s(m_t)}^{-1} g^{-1} \neq \varnothing.$$
Therefore
$$q_{s(m_t)} \in g^{-1} \{q_r \: s(k-t) \leq r \leq s(k-t)+8\} C_n$$
and hence for $t = 0, 1$ we have that $q_{s(m_t)}^{-1} q_{s(m_{t+1})}$ is an element of
$$C_n^{-1} \{q_r^{-1} \: s(k-t) \leq r \leq s(k-t)+8\} \{q_r \: s(k-t-1) \leq r \leq s(k-t-1)+8\} C_n.$$
Notice that for $h_1, h_2 \in \langle S \rangle$
$$h_1 \in C_n^{-1} h_2 C_n \Longrightarrow \nu_n(\|h_2\|) \leq \|h_1\| \leq \xi_n(\|h_2\|).$$
After recalling that
$$(4d+1)|r-t| - 2d \leq \|q_r^{-1} q_t\| \leq (4d+1)|r-t| + 2d$$
for all $r, t \in \N$, we have that
$$\nu_n((4d+1)(s(k-t) - s(k-t-1) - 8) - 2d)$$
$$\leq \|q_{s(m_t)}^{-1} q_{s(m_{t+1})}\| \leq$$
$$(4d+1)|s(m_{t+1})-s(m_t)|+2d$$
and
$$(4d+1)|s(m_{t+1})-s(m_t)| - 2d$$
$$\leq \|q_{s(m_t)}^{-1} q_{s(m_{t+1})}\| \leq$$
$$\xi_n((4d+1)(s(k-t) - s(k-t-1) + 8) + 2d).$$
However, for all $r \in \N$ we have $\xi_n(r) \leq \xi_{k-t}(r)$ and $\nu_{k-t}(r) \leq \nu_n(r)$. So
$$\nu_{k-t}((4d+1)(s(k-t) - s(k-t-1) - 8) - 2d)$$
$$\leq \nu_n((4d+1)(s(k-t) - s(k-t-1) - 8) - 2d)$$
and
$$\xi_n((4d+1)(s(k-t) - s(k-t-1) + 8) + 2d)$$
$$\leq \xi_{k-t}((4d+1)(s(k-t) - s(k-t-1) + 8) + 2d).$$
Therefore by the definition of $s: \N \rightarrow \N$ we must have that $\max \{m_t, m_{t+1}\} = k-t$. So $m_1 \leq \max(m_1, m_2) = k-1$ and $\max(m_0, m_1) = k$ together imply that $m_0 = k$. This gives $\phi(f(u,0,k)) = g \cdot f(v, j_0, k)$. Since $k \geq n+2$ was arbitrary, we conclude that for all $k \geq n+2$ there is $j(k) \in \N$ with
$$\phi(f(u,0,k)) = g \cdot f(v,j(k),k).$$

Since $z = \lim q_{s(r)}^{-1} \cdot z$, we see that $f(0) = \lim q_{s(k)}^{-1} \cdot f(u,0,k)$. By clause (vi) there is $j \in \N$ and $h \in G$ such that $\phi(f(0)) = h \cdot f(j)$. Then $h \cdot f(j) = \lim q_{s(k)}^{-1} g \cdot f(v,j(k),k)$. We have that $f(j) = h^{-1} \cdot (h \cdot f(j))$ is in the clopen set $X(v) - 1*X(v)$ and therefore for sufficiently large $k$ we must have $h^{-1} q_{s(k)}^{-1} g \cdot f(v,j(k),k) \in X(v) - 1*X(v)$. Since $[f(v,j(k),k)] \cap (X(v) - 1*X(v)) = q_{s(k)}^{-1} \cdot f(v,j(k),k)$, we must have $h^{-1} q_{s(k)}^{-1} g = q_{s(k)}^{-1}$ for all sufficiently large $k \in \N$. Therefore
$$f(j) = h^{-1} \cdot \lim q_{s(k)}^{-1} g \cdot f(v,j(k),k) = \lim q_{s(k)}^{-1} f(v,j(k),k)$$
so $j(k) = j$ for sufficiently large $k \in \N$.

(x). Let $g \in G$ be as in clause (vii) and let $j \in \N$ be as in clause (ix). For $i = 0, 1$ let $K_i^u = \{k \in \N \: u(h_k) = i\}$ and $K_i^v = \{k \in \N \: (h_j \cdot v)(h_k) = i\}$. Pick any $i = 0, 1$ with $K_i^u$ infinite. We have
$$\lim_{k \in K_i^u} f(u,0,k) = z_i.$$
Applying $\phi$ to both sides we get
$$\phi(z_i) = \phi \left(\lim_{k \in K_i^u} f(u,0,k) \right) = \lim_{k \in K_i^u} \phi(f(u,0,k))$$
$$= \lim_{k \in K_i^u} g \cdot f(v,j,k) = g \cdot \lim_{k \in K_i^u} f(v,j,k).$$
A priori we know that if $\lim_{k \in K_i^u} f(v,j,k)$ exists then this limit must be either $z_0$ or $z_1$. Since the limit does exist, $(h_j \cdot v)(h_k)$ must be constant for all but finitely many $k \in K_i^u$ and by clause (vi) and the equations above we have that this constant value must be $i$. Thus $K_i^u - K_i^v$ is finite.

By clause (iii) of Lemma \ref{LEM UNIVERSAL} we have that $K_1^u$ is infinite. Thus $K_1^u - K_1^v$ is finite. We now consider two cases. \underline{Case 1}: $K_0^u$ is infinite. Then $K_0^u - K_0^v$ is finite. Since $\N = K_0^u \cup K_1^u$, we have that $K = (K_0^u \cap K_0^v) \cup (K_1^u \cap K_1^v)$ is cofinite in $\N$ and $u(h_k) = (h_j \cdot v)(h_k)$ for all $k \in K$. Thus $u$ and $h_j \cdot v$ differ at only finitely many coordinates. So by clause (ii) of Lemma \ref{LEM UNIVERSAL} we have $u \ E_\infty \ v$. \underline{Case 2}: $K_0^u$ is finite. Since $\N = K_0^u \cup K_1^u$, we have that $K_1^u$ is cofinite in $\N$. Thus both $K_1^v$ and $K_1^u \cap K_1^v$ are cofinite in $\N$. So for all but finitely many $k \in \N$ we have $u(h_k) = 1 = (h_j \cdot v)(h_k)$. Thus $u$ and $h_j \cdot v$ differ at only finitely many coordinates. So by clause (ii) of Lemma \ref{LEM UNIVERSAL} we have $u \ E_\infty \ v$.
\end{proof}

\begin{cor}
For every countably infinite group $G$, $\TC(G)$ and $\TCF(G)$ are Borel bi-reducible.
\end{cor}

\begin{problem}
For a countably infinite nonlocally finite group $G$, what are the complexities of $\TCP(G)$, $\TCM(G)$, and $\TCMF(G)$?
\end{problem}

We point out that strangely we do not even know the answer to the above question in the case $G = \ZZ$. 
\chapter{Extending Partial Functions to $2$-Colorings}

In this chapter we study the problem when a partial function on a countably infinite group can be extended to a $2$-coloring on the entire group. The answer is immediate (and affirmative) if the partial function has a finite domain, since the set of all $2$-colorings is dense (Theorem~\ref{thm:density}). Results in this chapter can therefore be regarded as a strengthened form of density. A partial function with cofinite domain has only finitely many extensions. In a group with the ACP such a function can be extended to a $2$-coloring iff any extension of it is a $2$-coloring. However, in a non-ACP group $G$ we know that there are functions with domain $G-\{1_G\}$ so that one of the two extensions is a $2$-coloring and the other is periodic
(c.f., e.g., the proof of Theorem~\ref{thm:ACP}). These results suggest that it might be difficult to provide a unified solution to the above problem by stating an intrinsic condition on the partial function. Thus in this chapter we focus on the domain of the partial function. In Sections ~\ref{sec:suffext} and \ref{sec:charext} we characterize subsets of the group on which any partial function can be extended to a $2$-coloring of the full group. In Section~\ref{sec:autoext} we determine the countable group(s) $G$ for which any extension of a $2$-coloring on a nontrivial subgroup is a $2$-coloring on $G$.

\section{\label{sec:suffext}A sufficient condition for extendability}

For a countably infinite group $G$ and a subset $A\subseteq G$, we ask when any function with domain $A$ can be extended to 
a $2$-coloring on $G$. This problem will be the focus of the first two sections of this chapter. 

In this section we give a sufficient condition on $A$ for this extendability to hold. Although this result will soon be superseded by the results of the next
section, we include it here for two reasons. First, its proof is much easier and shorter than those in the next section. Second, the proof involves a new way to apply the fundamental method. In fact, the following proposition is a strengthening 
of Theorem~\ref{GEN COL} with a similar proof. It results from a careful scrutiny of  what the technique used in the proof of Theorem~\ref{GEN COL} can achieve.

\begin{prop}
Let $G$ be a countably infinite group and let $r: \N \rightarrow \N$ be any function. Then there exists a sequence $(T_n)_{n \in \N}$ of finite subsets of $G$ with the following property: if $\{u_{n,i} \in G \: n \in \N, \ 0 \leq i \leq r(n)\}$ is any collection of group elements then there is a fundamental $c \in 2^{\subseteq G}$ such that for every $n \in \N$ and $0 \leq i \leq r(n)$, $T_n$ witnesses that $c$ blocks $u_{n,i}$.
\end{prop}

\begin{proof}
For $n \geq 1$, define $p_n(k) = 2 \cdot k^4 \cdot (r(n-1)+1)$. Then $(p_n)_{n \geq 1}$ is a sequence of functions of subexponential growth. By Corollaries ~\ref{GROW BP} and \ref{GEN SUBEXP FREE} there is a blueprint $(\Delta_n, F_n)_{n \in \N}$ and a fundamental $c_0 \in 2^{\subseteq G}$ with
$$|\Theta_n| > \log_2 (2 \cdot |B_n|^4 \cdot (r(n-1)+1))$$
for each $n \geq 1$, where $B_n$ satisfies $\Delta_n B_n B_n^{-1} = G$. For $n \geq 1$, let $V_n$ be the test region for the $\Delta_n$ membership test admitted by $c_0$. For $n \geq 1$ define $T_n = B_{n+1} B_{n+1}^{-1} (V_{n+1} \cup \Theta_{n+1} b_n)$. We claim $(T_n)_{n \in \N}$ has the desired property.

Let $\{u_{n,i} \in G \: n \in \N, \ 0 \leq i \leq r(n)\}$ be any collection of elements of $G$. For $i, k \in \N$ let $\B_i(k)$ be the $i^\text{th}$ digit from least to most significant in the binary representation of $k$ when $k \geq 2^{i-1}$ and $\B_i(k) = 0$ when $k < 2^{i-1}$. For $n \geq 1$ let $s(n) = |\Theta_n|$ and let $\theta_1^n, \ldots, \theta_{s(n)}^n$ be an enumeration of $\Theta_n$. For $n \geq 1$, let $\Gamma_n$ be the graph with vertex set $\Delta_n$ and edge relation
$$\begin{array}{l}(\gamma, \psi) \in E(\Gamma_n) \Longleftrightarrow \\ 
\exists 0 \leq i \leq r(n-1) \ \gamma^{-1} \psi \in B_n B_n^{-1} u_{n,i} B_n B_n^{-1} \text{ or } \psi^{-1} \gamma \in B_n B_n^{-1} u_{n,i} B_n B_n^{-1}.
\end{array}
$$
for distinct $\gamma, \psi \in \Delta_n$. Then $\deg_{\Gamma_n} (\gamma) \leq 2 \cdot |B_n|^4 \cdot (r(n-1)+1)$ for each $\gamma \in \Delta_n$. We can therefore find, via the greedy algorithm, a graph theoretic $(2 |B_n|^4 (r(n-1)+1) + 1)$-coloring of $\Gamma_n$, say $\mu_n : \Delta_n \rightarrow \{0, 1, \ldots, 2|B_n|^4 (r(n-1)+1)\}$.

Define $c \supseteq c_0$ by setting
$$c(\gamma \theta_i^n b_{n-1}) = \mathbb{B}_i (\mu_n(\gamma))$$
for each $n \geq 1$, $\gamma \in \Delta_n$, and $1 \leq i \leq s(n)$. Since $2^{s(n)} > 2 |B_n|^4 (r(n-1)+1)$, all integers $0$ through $2 |B_n|^4 (r(n-1)+1)$ can be written in binary using $s(n)$ digits. Thus no information is lost between the $\mu_n$'s and $c$. Setting $\Theta_n(c) = \Theta_n(c_0) - \{\theta_1^n, \ldots, \theta_{s(n)}^n\}$ we clearly have that $c$ is fundamental.

Now an argument identical to that appearing in the proof of Theorem \ref{GEN COL} shows that $T_n$ witnesses that $c$ blocks $u_{n,i}$.
\end{proof}

\begin{theorem}\label{thm:FAAF} If $G$ is a countably infinite group and $A \subseteq G$ satisfies $F A^{-1}A  F \neq G$ for all finite sets $F \subseteq G$, then every partial function $c: A \rightarrow 2$ can be extended to a $2$-coloring on $G$.
\end{theorem}

\begin{proof}
For each $n \geq 1$ let $r_n : \N \rightarrow \N$ be the function which is constantly $1$. Let $(T_n)_{n \in \N}$ be as in the previous proposition. Fix an enumeration $s_1, s_2, \ldots$ of the nonidentity group elements of $G$. By assumption we have that
$$\{1_G, s_n^{-1}\} \{1_G, s_n\} T_n A^{-1} A T_n^{-1} \neq G.$$
For each $n \geq 1$ pick
$$h_n \not\in \{1_G, s_n^{-1}\} \{1_G, s_n\} T_n A^{-1} A T_n^{-1}.$$
For each $n \geq 1$ set $u_{n,0} = s_n$ and $u_{n,1} = h_n^{-1} s_n h_n$. Let $c \in 2^G$ be fundamental and such that $T_n$ witnesses that $c$ blocks $u_{n,0}$ and $u_{n,1}$ for each $n \geq 1$.

Let $x \in 2^A$ be an arbitrary function. Define $y \in 2^G$ by
$$y(g) = \begin{cases}
x(g) & \text{if } g \in A \\
c(g) & \text{otherwise.}
\end{cases}$$
So $y$ extends $x$. We will show that $y$ is a $2$-coloring of $G$. Fix $1_G \neq s \in G$. Then for some $n \geq 1$ we have $s = s_n$. Set $T = T_n \cup h_n T_n$ and let $g \in G$ be arbitrary. Notice that $A \cap (g T_n \cup g s_n T_n) \neq \varnothing$ if and only if $g \in A T_n^{-1} \cup A T_n^{-1} s_n^{-1}$. So by the definition of $h_n$ we have that
$$A \cap (g T_n \cup g s_n T_n) \neq \varnothing \Longrightarrow A \cap (g h_n T_n \cup g s_n h_n T_n) = \varnothing.$$
If $A \cap (g T_n \cup g s_n T_n) = \varnothing$ then set $k = 1_G$. Otherwise set $k = h_n$. In either case we have $A \cap (g k T_n \cup g s_n k T_n) = \varnothing$. Therefore for all $t \in T_n$
$$y(g k t) = c(g k t) \text{ and } y(g s_n k t) = c(g s_n k t).$$
Notice that $T_n$ witnesses that $c$ blocks $k^{-1} s_n k$. Therefore, there is $t \in T_n$ with
$$y(g k t) = c(g k t) \neq c(g k (k^{-1} s_n k) t) = c(g s_n k t) = y(g s_n k t).$$
This completes the proof since $k t \in T$.
\end{proof}

Note that this gives another proof for the density of $2$-colorings.

\section{\label{sec:charext}A characterization for extendability}

In this section we continue to consider the problem when any partial function with domain $A$ can be extended to a $2$-coloring on $G$.

There is an obvious obstacle if the set $A$ is too large in the following sense. We let $\chi_A$ denote the characteristic function of $A \subseteq G$:
$$\chi_A(g) = \begin{cases}
1 & \text{if } g \in A \\
0 & \text{if } g \not\in A.
\end{cases}$$
If $1\in \overline{[\chi_A]}$ then there is no $2$-coloring on $G$ extending $1\in 2^A$, the constant $1$ function with domain $A$. This is because, for any $x \in 2^G$ extending $1 \in 2^A$ and sequence $(g_n)_{n\in \N}$ of elements of $G$ so that  $1 = \lim g_n \cdot \chi_A$, we must have $\lim g_n\cdot x=1$, since $x$ can have value $0$ only when $\chi_A$ has value $0$. This shows that $x$ is not a $2$-coloring. A moment of reflection shows that $1\in \overline{[\chi_A]}$ is indeed a largeness condition since it is equivalent to saying that $A$ contains a translate of any finite subset of $G$. 

The objective of this section is to show that this is the only obstacle for the extendability. We introduce the following terminology.

\begin{definition} We say that $A\subseteq G$ is {\it slender} if $1\not\in\overline{[\chi_A]}$.
\end{definition}\index{slender set}

The following characterizations of slenderness are immediate. We state them without proof.

\begin{lem} Let $G$ be a countably infinite group. The following are equivalent for $A\subseteq G$:
\begin{itemize}
\item[(i)] $A$ is slender, i.e., $1\not\in\overline{[\chi_A]}$;
\item[(ii)] $1\,\bot\, \chi_A$;
\item[(iii)] there exists a finite $T \subseteq G$ so that for every $g \in G$ there is $t \in T$ with $gt \not\in A$.
\end{itemize}
\end{lem}

It is easy to see that any proper subgroup is slender. In fact, if $H<G$ and $T=\{1_G, a\}$ for some arbitrary $a\in G-H$, then for any $g\in G$, $gT\not\subseteq H$.

The main technical result of this section is the following theorem.

\begin{theorem}\label{thm:slenderextension}
Let $G$ be a countably infinite group, $A\subseteq G$ a slender subset, $s\in G$ a nonidentity element, and $y\in 2^A$.
There exist a slender $A'\supseteq A$ and $x_0, x_1\in 2^{A'}$ extending $y$ such that 
\begin{itemize}
\item[(a)] any extension of $x_0$ or $x_1$ blocks $s$, and 
\item[(b)] for any $z_0, z_1\in 2^G$ extending $x_0, x_1$ respectively, $z_0\,\bot\, z_1$.
\end{itemize}
\end{theorem}

Before giving the long and technical proof, we offer some remarks on its structure and main ideas. Then during the proof we give more commentaries to elaborate on the ideas. In this proof we will try to recreate, as much as is possible and needed, the machinery used in our standard construction of a $2$-coloring. However, since we only need to block a single element $s$, we do not need to construct an entire blueprint, but just a single $\Delta$ and $F$. The construction of $\Delta$ and $F$ will run parallel to extending $y$. In extending $y$, we shall create a membership test for the set $\Delta$. The membership test will rely on counting the number of $1$'s within a finite test region, but it will not be a simple membership test as used in our original fundamental method. The bulk of the work is to strategically add a lot of $1$'s at select locations but at the same time make sure they are not visible from other unwanted locations. This makes the $\Delta$-translates of $F$ look different from the background.

\begin{proof}[Proof of Theorem~\ref{thm:slenderextension}]
Let $G$, $A$, $s$ and $y$ be given. Let $T\subseteq G$ be a finite set such that $gT\not\subseteq A$ for any $g\in G$. We fix a finite $B\subseteq G$ with $1_G\in B$ 
such that for every $g \in G$, 
$$|g B \cap (G-A)| \geq 2.$$ 
Such a $B$ can be taken to be the union of two disjoint (left) translates of $T$, with one of them containing $1_G$. 

Much of this proof will rely on counting the number of $1$'s of a partial function within some left translate of $B$. 
We define a counting function as follows. For $z \in 2^{\subseteq G}$ and $g \in G$, let
$$r_z(g) = |\{b \in B \: gb \in \dom(z) \mbox{ and } z(gb) = 1\}|.$$
Note that $r_z$ is total even if $z$ is partial, and $r_z(g)\leq |B|$ in general. Set $N = |B|-2$ so that $\max\{r_y(g) \: g \in G\} \leq N$.

If we add in more $1$'s in a particular translate $g B$, then these $1$'s will be visible in at most the $g B B^{-1}$-translates of $B$. In order to control the number of $1$'s seen within translates of $B$, we will often insist that translates of $BB^{-1}$ be disjoint. Let $C$ be a finite symmetric set (meaning $C = C^{-1}$) containing $B B^{-1}$. Elements outside $gC$ will not see in their translate of $B$ newly added $1$'s in $gB$. However, since we can not expect to have a locally recognizable function on $gC$, a difficulty is how to pinpoint the precise element of the prospective $\Delta$ if we already know 
it is within $gC$. A natural solution is to look beyond $gC$ and use information in other parts of $gF$ (which requires that the prospective $F$ be sufficiently big). The following function $a$ tells us where to look for this additional information.

For every $c_1, c_2 \in C$, fix $a(c_1, c_2) \in G$ so that the following conditions hold:
$$\forall c_1, c_2 \in C \ a(c_1, c_2) = a(c_2, c_1);$$
$$\forall c_1, c_2 \in C \ c_1 a(c_1, c_2) C \cap C = \varnothing;$$
$$\forall c_1, c_2, c_3, c_4 \in C \ \{c_1, c_2\} \neq \{c_3, c_4\} \Longrightarrow c_1 a(c_1, c_2) C \cap c_3 a(c_3, c_4) C = \varnothing.$$
Note that the symmetry in the first condition implies that if $\{c_1, c_2\} \neq \{c_3, c_4\}$ then $c_2 a(c_1, c_2) C \cap c_3 a(c_3, c_4) C = \varnothing$. Such a function $a$ exists since $G$ is infinite.

Now let $F \subseteq G$ be finite with
$$F \supseteq C^3 \cup C \cdot \bigcup_{c_1, c_2 \in C} a(c_1, c_2) C$$
and satisfying
$$\rho(F; C) - |C|^6 - |C|^3 \geq \log_2 \ (2|C|^{3N+3}|F|^2+1) + \log_2 \ (|C|) + 4.$$
Such an $F$ exists by Lemma~\ref{FAST GROWTH}.

Let $D_0 \subseteq r_y^{-1}(N)=\{g\in G \: r_y(g)=N\}$ be a maximal subset of $r_y^{-1}(N)$ with the $D_0$-translates of $C F$ disjoint. Similarly, let $D_1 \subseteq r_y^{-1}(N-1)$ be a maximal subset of $r_y^{-1}(N-1)$ with the $D_1$-translates of $C F$ disjoint and $D_1 C^4 F \cap D_0 C F = \varnothing$. In general, once $D_0$ through $D_{m-1}$ have been defined ($1 < m \leq N$), let $D_m \subseteq r_y^{-1}(N-m)$ be a maximal subset of $r_y^{-1}(N-m)$ with the $D_m$-translates of $C F$ disjoint and
$$D_m C^{3m+1} F \cap \bigcup_{0 \leq i < m} D_i C F = \varnothing.$$
We set $D = \bigcup_{0 \leq i \leq N} D_i$. One or more $D_i$ might be finite or even empty (including $D_N$), but $D$ is always infinite. To provide some perspective, the set $D$ will soon be modified slightly (each element of $D$ right translated by an element of $C$) to create $\Delta$.

We point out two important properties of $D$. First, let $g \in G$ and let $0 \leq m \leq N$ be such that $r_y(g) = N-m$. Then by the definition of $D_m$ either
$$g C^{3m+1} F \cap \bigcup_{0 \leq i < m} D_i C F \neq \varnothing$$
or else $g C F \cap D_m C F \neq \varnothing$. In any case,
$$g \in \bigcup_{0 \leq i \leq m} D_i C F F^{-1} C^{3m+1} \subseteq D C F F^{-1} C^{3N+1}.$$
Second, let $0 \leq m \leq N$, let $d \in D_m$, and let $g \in d C^3$. Then for any $t < m$
$$g C^{3t+1} F \cap \bigcup_{0 \leq i \leq t} D_i C F \subseteq D_m C^{3m+1} F \cap \bigcup_{0 \leq i < m} D_i C F = \varnothing$$
and therefore $r_y(g) \neq N - t$. It follows that $r_y(g) \leq N-m$ for all $g \in D_m C^3$. Thus each point in $D$ achieves a
local maximum $r_y$-value.

In what follows we will go through a number of consecutive extensions of $y$, and in the middle of these extensions we will also define $\Delta$. 
We first extend $y$ to $y_1$ so that
$$\dom(y_1) = \dom(y) \cup D C$$
and for every $d \in D$ all elements of $dC - \dom(y)$ are assigned the value $0$ except for precisely 2 elements in $dB - \dom(y)$ which are assigned the value $1$. $y_1$ exists since the $D$ translates of $C$ are disjoint ($1_G \in C \subseteq F$) and every left translate of $B$ contains at least 2 elements not in $\dom(y)$. Notice that for $g \in G$,
$$r_y(g)\leq r_{y_1}(g) \leq r_y(g) + 2$$
and
$$r_{y_1}(g) > r_y(g) \Longrightarrow g \in DC.$$ 

Next we extend $y_1$ to $y_2$ where $y_2$ has domain
$$\dom(y_1) \cup D \{c_1 a(c_1, c_2) \: c_1 \neq c_2 \in C\} B$$
and satisfies for each $d \in D$ and each $c_1 \neq c_2 \in C$
$$\exists b \in B \ y_2(d c_1 a(c_1, c_2) b) \neq y_2(d c_2 a(c_1, c_2) b),$$
and
$$r_{y_2}(d c_1 a(c_1, c_2)) + r_{y_2}(d c_2 a(c_1, c_2)) \leq r_{y_1}(d c_1 a(c_1,c_2)) + r_{y_1}(d c_2 a(c_1, c_2)) + 1.$$
Notice that for $c_1, c_2 \in C$, $c_1 a(c_1, c_2) B \subseteq F$, so one can extend $y_1$ to $y_2$ by considering one $D$-translate of $F$ at a time. Also, by the construction of the function $a$ we have that for any $d \in D$ and $c_1, c_2 \in C$, $dc_1 a(c_1, c_2) B \cap D C = \varnothing$. Thus for $d \in D$ and $c_1, c_2 \in C$, $d c_1 a(c_1, c_2) B$ contains at least 2 elements not in $\dom(y_1)$. When $c_1 \neq c_2$, $d c_1 a(c_1, c_2) \neq d c_2 a(c_1, c_2)$ so one can achieve the last two requirements listed above. Finally, by construction $d c_1 a(c_1, c_2) B \cap d c_3 a(c_3, c_4) B = \varnothing$ for any $d \in D$ and $c_1 \neq c_2, c_3 \neq c_4 \in C$ with $\{c_1, c_2\} \neq \{c_3, c_4\}$. Therefore the function $y_2$ exists. In the above argument, we showed that various group elements had disjoint $B$ translates. However, the reader can easily check that they all have disjoint $C$ translates. We therefore notice that $y_2$ has the following properties for every $g \in G$:
$$r_y(g)\leq r_{y_2}(g) \leq r_y(g) + 2;$$
$$r_{y_2}(g) > r_y(g) \Longrightarrow g \in D(C \cup \bigcup_{c_1 \neq c_2 \in C} c_1 a(c_1, c_2) C);$$
$$r_{y_2}(g) > r_y(g) + 1 \Longrightarrow g \in DC.$$

For further extensions we set
$$V = \{a(c_1, c_2) \: c_1 \neq c_2 \in C\} C.$$
Now extend $y_2$ to $y_3$ where
$$\dom(y_3) = \dom(y_2) \cup D C V$$
and for all $g \in \dom(y_3) - \dom(y_2)$, $y_3(g) = 0$. Since $y_3$ extends $y_2$ using only the value $0$, $y_3$ also satisfies the three properties listed above for $y_2$.

Now we can modify $D$ to get $\Delta$. The construction of $y_1$ was a naive attempt of creating a membership test for $D$ by placing several $1$'s within a translate of $B$. The problem is that we may not be able to uniformly tell the difference between elements of $DC$ and elements of $D$ only using a finite set. As we will now see, the construction of $y_2$ and $y_3$ allow us to recognize a particular element of $d C$ for each $d \in D$, though this recognized element may not be $d$ itself.

By the definition of $y_2$, we have for any $d \in D$ and $c_1 \neq c_2 \in C$, there is $b \in B$ with $y_2(d c_1 a(c_1, c_2) b) \neq y_2(d c_2 a(c_1, c_2) b)$ and hence 
$$ ((d c_1)^{-1}\cdot y_3) \upharpoonright V \neq ((dc_2)^{-1}\cdot y_3) \upharpoonright V. $$
Also, by definition of $y_3$, for every $d \in D$ and $c \in C$, $((dc)^{-1} \cdot y_3) \upharpoonright V \in 2^V$. Arbitrarily pick a total ordering, $\prec$, on $2^V$. We define a function 
$$q: D=\bigcup_{0\leq m\leq N} D_m \rightarrow DC$$ 
as follows. If $d \in D_m$, then we let $q(d) = d c$, where $c$ is the unique element of
$$S = \{c' \in C \: r_{y_3}(dc') = N-m+2\}$$
so that $((dc)^{-1}\cdot y_3)\upharpoonright V$ is the $\prec$-largest among
$$\{((dc')^{-1} \cdot y_3) \upharpoonright V \: c' \in S\}.$$
We define $\Delta = q(D)$, and for each $0 \leq m \leq N$ we define $\Delta_m = q(D_m)$. Since the $D$-translates of $C F$ are disjoint, from the definition of $q$ it follows that the $\Delta$-translates of $F$ are disjoint.

Although we have defined the sets $\Delta$ and $F$, a membership test for $\Delta$ has to wait until we have further extended $y_3$. One can expect that the prospective membership test will be inductive on $\Delta_m$, and it will necessarily use the property that $\Delta_m$ elements not only have local maximum values for the counting function, but also achieve the maximum in the $\prec$-ordering among their local competitors. All these desired properties are already in place by our construction so far. However, it seems that, in order to carry out the induction successfully, we need to make use of the maximal disjointness properties imposed on elements of $D$ when it was defined. It is therefore necessary for us to create also a membership test for $D$ at the same time. And this will be implemented by trying to code information contained in the function $q$
by parities of the counting function values. This extra coding is the main idea of the following extensions.

Let $\Lambda \subseteq F$ be such that the $\Lambda$-translates of $C$ are contained and maximally disjoint within $F$. Notice that:
$$|\{ \lambda \in \Lambda \: C \lambda C \cap C \neq \varnothing\}| = |\{\lambda \in \Lambda \: \lambda \in C^3\}| \leq |C|^3;$$
$$|\{ \lambda \in \Lambda \: C \lambda C \cap C V\neq \varnothing \}| = |\{\lambda \in \Lambda \: \lambda \in C^2 V C \}| \leq |C|^6.$$
We have $|\Lambda| \geq \rho(F; C)$, and therefore there are at least
$$\log_2 \ (2|C|^{3N+3}|F|^2+1) + \log_2 \ (|C|) + 4$$
elements of $\Lambda$ which are in neither of the sets above. Let $\lambda_1, \lambda_2, \ldots, \lambda_M$ enumerate the elements of $\Lambda$ which are in neither of the above sets. Let $K$ be the least integer greater than $\log_2 \ (|C|)$. Notice that $M - K - 2 \geq \log_2 \ (2|C|^{3N+3}|F|^2+1)$. 

As $\{\lambda_1, \lambda_2, \ldots, \lambda_M\} C \subseteq F$ and the $D$-translates of $C F$ are disjoint, it follows that the $\Delta$-translates of $\{\lambda_1, \lambda_2, \ldots, \lambda_M\} C$ are disjoint. Also, by the choice of $\lambda_1$ through $\lambda_M$, we have for each $1 \leq i \leq M$
$$\Delta \lambda_i C \cap D (C \cup C V) \subseteq D C \lambda_i C \cap D (C \cup C V) = \varnothing.$$
Since
$$y_3 \upharpoonright (G - D(C \cup C V)) = y \upharpoonright (G - D(C \cup C V))$$
each $\Delta \{\lambda_1, \lambda_2, \ldots, \lambda_M\}$-translate of $B$ contains at least 2 elements not in the domain of $y_3$. Let $y_4$ extend $y_3$ by only assigning the value $0$ and have a domain with the property that every $\Delta \{\lambda_1, \ldots, \lambda_M\}$-translate of $B$ has exactly one undefined point and that no other points besides these are undefined. To be more precise, we require
$$\forall g \in \dom(y_4) - \dom(y_3) \ y_4(g) = 0,$$
$$G - \Delta \{\lambda_1, \ldots, \lambda_M\} B \subseteq \dom(y_4),$$
and
$$|(G - \dom(y_4)) \cap \gamma \lambda_i B| = 1$$
for every $\gamma \in \Delta$ and $1 \leq i \leq M$. It follows from the properties of $y_3$, the properties of the $\lambda_i$'s, and the definition of $y_4$ that for every $g \in G$:
$$r_y(g)\leq r_{y_4}(g) \leq r_y(g) + 2;$$
$$r_{y_4}(g) > r_y(g) \Longrightarrow g \in D(C \cup CV);$$
$$r_{y_4}(g) > r_y(g) + 1 \Longrightarrow g \in DC.$$
Moreover, we have that if $z \in 2^G$ is any extension of $y_4$ then for every $g \in G$:
$$r_y(g)\leq r_z(g) \leq r_y(g) + 2;$$
$$r_z(g) > r_y(g) \Longrightarrow g \in D(C \cup CV \cup C \{\lambda_1, \ldots, \lambda_M\}C);$$
$$r_z(g) > r_y(g) + 1 \Longrightarrow g \in DC.$$

For $i, k \in \N$, we let $\mathbb{B}_i(k) \in \{0, 1\}$ be the $i^\text{th}$ digit from least to most significant in the binary representation of $k$ when $k \geq 2^{i-1}$ and $\mathbb{B}_i(k) = 0$ when $k < 2^{i-1}$. Let $p: C \rightarrow 2^K$ be any injective function. Extend $y_4$ to $y_5$ so that
$$\dom(y_5) = \dom(y_4) \cup \Delta \{\lambda_1, \lambda_2, \ldots, \lambda_K\} B$$
and for every $\gamma \in \Delta$ and $1 \leq i \leq K$
$$r_{y_5}(\gamma \lambda_i) \equiv \mathbb{B}_i(p(d^{-1} \gamma)) \mod 2$$
where $d \in D$ is such that $\gamma \in dC$ (or equivalently, $d = q^{-1}(\gamma)$). Recall that $K > \log_2 \ (|C|)$, so no information is lost between the function $p$ and its ``encoding'' into $y_5$. The function $y_5$ and any extension of it still satisfy the last 3 properties listed at the end of the previous paragraph.

We claim that for any $0 \leq m \leq N$ $y_5$ admits a $\Delta_m$ membership test if and only if it admits a $D_m$ membership test. Fix $0 \leq m \leq N$ and first suppose that $y_5$ admits a $\Delta_m$ membership test. Let $z \in 2^G$ be an arbitrary extension of $y_5$. Then for $g \in G$, $g \in D_m$ iff 
\begin{quote}
there is $c \in C$ so that $gc \in \Delta_m$ and for every $1 \leq i \leq K$
$$r_z(g c \lambda_i) \equiv \mathbb{B}_i(p(c)) \mod 2.$$
\end{quote}
If $d \in D_m$, then clearly there is $c \in C$ with $dc \in \Delta_m$ and the last requirement is satisfied by definition of $y_5$. Now suppose $g \in G$ has the stated properties. Let $c \in C$ be such that $g c \in \Delta_m$. It is clear that for any $\gamma \in \Delta$ the parities of $r_z(\gamma \lambda_i)$ ($1 \leq i \leq K$) uniquely determine an element $c' \in C$ with $\gamma (c')^{-1} \in D$. Therefore, the property satisfied by $g$ clearly implies $g \in D_m$. Now suppose that $y_5$ admits a $D_m$ membership test. Let $z \in 2^G$ be an arbitrary extension of $y_5$. Since $z$ is defined on all of $G$, we can define a relation $\prec_z$ on $G$ by
$$g \prec_z h \Longleftrightarrow (g^{-1} \cdot z) \upharpoonright V \prec (h^{-1} \cdot z) \upharpoonright V.$$
The relation $\prec_z$ is a quasiordering, i.e. it is reflexive and transitive, but fails to be antisymmetric. An important property of $\prec_z$ is that the construction of $y_3$ ensures that $\prec_z$ is a total ordering whenever restricted to a single $D$-translate of $C$. We claim that for $g \in G$, $g \in \Delta_m$ iff 
\begin{quote}
$r_z(g) = N - m + 2$, there is $c \in C$ with $g c \in D_m$, and $h \prec_z g$ for all elements $h \in g c C$ with $r_z(h) = N - m +2$. 
\end{quote}
By the way $\Delta_m$ was defined, it is clear that elements of $\Delta_m$ have these properties. If $g \in G$ satisfies these properties, then let $c \in C$ be such that $g c = d \in D_m$. Let $\gamma \in \Delta_m \cap d C$. By definition, $\gamma$ is $\prec_z$-largest among all elements $h \in d C$ with $r_z(h) = N - m +2$. So $g \prec_z \gamma$, $\gamma \prec_z g$, and $\gamma, g \in d C$. Therefore $g = \gamma \in \Delta_m$. 

Now we will show that $y_5$ admits a $\Delta$ membership test. For this, it will suffice to show that each $\Delta_m$ has a membership test. We will use induction on $m$. 

We begin with $\Delta_0$. Let $z \in 2^G$ be an arbitrary extension of $y_5$. Let $\prec_z$ be defined as above. We claim that for $g \in G$, $g \in \Delta_0$ iff 
\begin{quote}
$r_z(g) = N+2$ and $h \prec_z g$ for all $h \in g C^2$ with $r_z(h) = N+2$.
\end{quote}
First, let $\gamma \in \Delta_0$. Then clearly $r_z(\gamma) = N+2$. By construction, $C^3 \subseteq F$, so the $D$-translates of $C^3$ are disjoint. Let $d \in D_0$ be such that $\gamma \in d C$. Then $\gamma C^2 \subseteq d C^3$ and therefore if $d' \in D$ and $\gamma C^2 \cap d' C \neq \varnothing$, then $d' = d$. Additionally, every element $h \in G$ with $r_z(h) = N+2$ is an element of $DC$. Therefore,
$$\{h \in \gamma C^2 \: r_z(h) = N+2\} \subseteq dC.$$
It then follows from the construction of $\Delta_0$ that $h \prec_z \gamma$ for all such $h$. Conversely, let $g \in G$ be such that $r_z(g) = N+2$ and $h \prec_z g$  for all $h \in g C^2$ with $r_z(h) = N+2$. By construction, we immediately have $g \in D_0C$. Let $d \in D_0$ be such that $g \in d C$. Let $\gamma \in d C \cap \Delta_0$. Now $\gamma \in dC \subseteq g C^2$ and $r_z(\gamma) = N + 2$, so $\gamma \prec_z g$. However, $g$ and $\gamma$ are both in $dC$ and by construction $g \prec_z \gamma$. It follows $g = \gamma \in \Delta_0$.   

For the inductive step, let $0 < m \leq N$ and suppose that $y_5$ admits a $\Delta_t$-membership test for all $t < m$. Then $y_5$ also admits a $D_t$-membership test for all $t < m$. Let $z \in 2^G$ be an arbitrary extension of $y_5$, and define $\prec_z$ as before. We claim that for $g \in G$, $g \in \Delta_m$ iff there is $c \in C$ such that
$$r_z(g) = N - m + 2;$$
$$r_z(g c) = N - m +2;$$
$$gc C^{3m+1} F \cap (\bigcup_{0 \leq t < m} D_t C F ) = \varnothing;$$
$$\text{and } h \prec_z g \text{ for all } h \in g C^2 \text{ with } r_z(h) = N - m + 2.$$
Note that one can test the truth of the third requirement by only checking the values of $z$ on a finite set. One only needs to run the $D_t$-membership test for each $0 \leq t < m$ and each element of $g c C^{3m+1} F F^{-1} C$. First suppose $\gamma \in \Delta_m$. Then clearly there is $c \in C$ with $\gamma c = d \in D_m$. By construction, $r_z(\gamma) = r_z(\gamma c) = N - m + 2$. Also, from the definition of $D_m$ we have that $d = \gamma c$ satisfies the third condition listed above. It follows from the paragraph following the definition of $D$ that $r_y(h) \leq N - m$ for all $h \in \gamma C^2 \subseteq d C^3$. Therefore every $h \in \gamma C^2$ with $r_z(h) = N - m + 2$ must lie in $D C$. Since the $D$-translates of $C^3$ are disjoint, it follows that
$$\{h \in \gamma C^2 \: r_z(h) = N - m + 2\} \subseteq d C.$$
So it follows from the definition of $\Delta_m$ that $\gamma$ is $\prec_z$ largest among all such $h$.

Conversely, suppose $g \in G$ and $c \in C$ satisfy the conditions above. Then for every $i < m$ we have
$$g c C^{3i+1} F \cap (\bigcup_{0 \leq t < i} D_t C F) = \varnothing$$
and
$$g c C F \cap D_i C F \subseteq g c C^{3m+1} F \cap D_i C F = \varnothing.$$
Therefore, it follows from the definition of $D_i$ that for every $0 \leq i < m$ $r_y(g c) \neq N-i$. So $r_y(g c) \leq N-m$. Since $r_z(g c) = N-m+2$, it must be that $r_y(g c) = N-m$ and $g c \in D_i C$ for some $m \leq i \leq N$. As mentioned in the paragraph following the definition of $D$, we have $r_y(h) \leq N - i$ for all $h \in D_i C^3$. So we must have $g c \in D_m C$. Let $d \in D_m$ be such that $g c \in d C$. Then $g \in d C^2$. We have $r_y(g) \leq N - m$. Since $r_z(g) = N - m + 2$, it must be that $r_y(g) = N - m$ and $g \in D C$. Since the $D$-translates of $C^2$ are disjoint, $g \in d C$. Then $d C \subseteq g C^2 \subseteq d C^3$. It follows that
$$\{h \in g C^2 \: r_z(h) = N - m + 2 \} \subseteq d C.$$
If $\gamma \in d C \cap \Delta_m$, then $r_z(\gamma) = N - m + 2$ so $\gamma \prec_z g$. By construction of $\Delta_m$, $g \prec_z \gamma$. Therefore $g = \gamma \in \Delta_m$. This finishes the proof that $y_5$ admits a $\Delta$ membership test. It follows also that $y_5$ admits a $D$ membership test.

The rest of the proof uses familiar techniques in the fundamental method. The remaining tasks are to extend $y_5$ to block the element $s$, to reach further orthogonal extensions, and at the same time to maintain the slenderness of the domain of the partial functions constructed. 

Recall that for any $g \in G$
$$g \in D C F F^{-1} C^{3N+1}.$$
Therefore, for any $g \in G$
$$g \in \Delta C^2 F F^{-1} C^{3N+1}.$$
Let $\Gamma$ be the graph with vertex set $\Delta$ and edge relation given by
$$(\gamma, \psi) \in E(\Gamma) \Longleftrightarrow \exists h \in C^2 F F^{-1} C^{3N+1} \ (\gamma h s h^{-1} = \psi \vee \psi h s h^{-1} = \gamma)$$
for $\gamma, \psi \in \Delta$. Then each element of $\Delta$ has degree at most $2|C|^{3N+3}|F|^2$ in $\Gamma$. By applying the greedy algorithm, we can find a graph-theoretic coloring of $\Gamma$ using $2|C|^{3N+3}|F|^2+1$ many colors, say $\mu: \Delta \rightarrow (2|C|^{3N+3}|F|^2+1)$. Now we extend $y_5$ to $x$ so that
$$\dom(x) = \dom(y_5) \cup \Delta \{\lambda_{K+1}, \lambda_{K+2}, \ldots, \lambda_{M-2}\} B$$
and
$$r_x(\gamma \lambda_{K+i}) \equiv \mathbb{B}_i(\mu(\gamma)) \mod 2$$
for all $\gamma \in \Delta$ and $1 \leq i \leq M - K - 2$. Since
$$M - K - 2 \geq \log_2 \ (2|C|^{3N+3}|F|^2+1)$$
no information is lost between $\mu$ and $x$.

We claim that any extension of $x$ blocks $s$. Let $W$ be the test region for the $\Delta$-membership test admitted by $x$. Let $T_0 = C^{3N+1} F F^{-1} C^2$, $T = T_0 (F \cup W)$, and let $g \in G$ be arbitrary. Since $\Delta T_0^{-1} = G$, there is $t \in T_0$ with $gt \in \Delta$. If $g s t \not\in \Delta$, then $gt$ passes the $\Delta$-membership test while $g s t$ fails. Therefore there is $w \in W$ with $gtw, gstw \in \dom(x)$ and $x(gtw) \neq x(gstw)$. As $tw \in T$, this completes this case. Now suppose $gst \in \Delta$. Then $gst = gt(t^{-1} s t)$ so $gt$ and $gst$ are joined by an edge in $\Gamma$. It follows $\mu(gt) \neq \mu(gst)$. So there is $1 \leq i \leq M - K - 2$ with $\mathbb{B}_i(\mu(gt)) \neq \mathbb{B}_i(\mu(gst))$. Thus
$$r_x(gt \lambda_{K+i}) \not\equiv r_x(gst \lambda_{K+i}) \mod 2$$
and therefore there is $b \in B$ with $x(gt \lambda_{K+i} b) \neq x(gst \lambda_{K+i} b)$ (recall $\Delta \lambda_{K+i} B \subseteq \dom(x)$). This completes the proof of the claim since $t \lambda_{K+i} b \in T_0 F \subseteq T$.

Now for $\iota = 0, 1$ define $x_\iota$ so that
$$\dom(x_\iota) = \dom(x) \cup \Delta \lambda_{M-1} B$$
and
$$r_x(\gamma \lambda_{M-1}) \equiv \iota \mod 2$$
for all $\gamma \in \Delta$.  Let $T_0$, $W$, and $T$ be as in the previous paragraph. We claim that for any $g_0,g_1\in G$ there is $\tau\in T$ such that $g_0\tau, g_1\tau\in \dom(x_0)=\dom(x_1)$ and $x_0(g_0\tau)\neq x_1(g_1\tau)$. Clearly it follows that $z_0\,\bot\, z_1$ for any $z_0, z_1\in 2^G$ extending $x_0,x_1$ respectively. To prove the claim, let $g_0, g_1 \in G$. Then there is $t \in T_0$ with $g_0 t \in \Delta$. If $g_1 t \not\in \Delta$ then $g_0 t$ passes the $\Delta$-membership test while $g_1 t$ fails. Thus, if $g_1 t \not\in \Delta$ then there is $w \in W$ with $x(g_0 t w) \neq x(g_1 t w)$. Clearly $tw \in T$. So suppose $g_1 t \in \Delta$. Then
$$r_{x_0}(g_0 t \lambda_{M-1}) \not\equiv r_{x_1}(g_1 t \lambda_{M-1}) \mod 2$$
so there is $b \in B$ with $x_0(g_0 t \lambda_{M-1} b) \neq x_1(g_1 t \lambda_{M-1} b)$.

Finally, we verify that $A'=\dom(x_i)$ is slender. This is immediate since if $T_0$ is as before then for every $g \in G$ there is $t \in T_0$ with $g t \in \Delta$ and
$$|g t \lambda_M B \cap (G - \dom(x_i))| = 1.$$
\end{proof}

The main result of this section now follows quickly.

\begin{theorem} \label{THM EXT SLENDER}
Let $G$ be a countably infinite group and $A \subseteq G$. Then the following are equivalent:
\begin{enumerate}
\item[\rm (i)] $A$ is slender;
\item[\rm (ii)] There exists a $2$-coloring on $G$ extending the constant function $1 \in 2^A$;
\item[\rm (iii)] For every $y \in 2^A$ there exists a $2$-coloring on $G$ extending $y$;
\item[\rm (iv)] For every $y \in 2^A$ there exists a perfect set of pairwise orthogonal $2$-colorings on $G$ extending $y$.
\end{enumerate}
\end{theorem}

\begin{proof} We have shown (ii) $\Rightarrow$ (i) at the beginning of this section by noting that if $A$ is not slender then $1\in 2^A$ cannot be extended to any $2$-coloring on $G$. The implications (iv) $\Rightarrow$ (iii) $\Rightarrow$ (ii)
are obvious.

It suffices to verify (i) $\Rightarrow$ (iv). For this we consider an additional clause:
\begin{enumerate}
\item[(v)] For any enumeration $s_1, s_2, \ldots$ of the nonidentity elements of $G$ and $y \in 2^A$, there is a collection  $\{x_u\in 2^{\subseteq G} \: u \in 2^{<\N}\}$ such that for all $u\in 2^{<\N}$,
\begin{itemize}
\item[(va)] $y \subseteq x_u \subseteq x_{u^\smallfrown 0}, x_{u^\smallfrown 1}$, 
\item[(vb)] $x_u$ blocks $s_i$ for all $i \leq |u|$, 
\item[(vc)] for $u \neq v \in 2^{<\N}$ with $|u| = |v|$, if $z_u, z_v\in 2^G$ extend $x_u, x_v$ respectively, then $z_u\,\bot\, z_v$,
\item[(vd)] $\{1_G, s_1, s_2, \ldots, s_{|u|}\} \subseteq \dom(x_u)$.
\end{itemize}
\end{enumerate}
We show (i) $\Rightarrow$ (v) $\Rightarrow$ (iv). 

(i) $\Rightarrow$ (v). Let $s_1, s_2, \ldots$ be an enumeration of the nonidentity elements of $G$ and let $y \in 2^A$. If necessary, arbitrarily define $y(1_G)$ so that $1_G \in \dom(y)$.  Applying Theorem~\ref{thm:slenderextension} to $s_1$ and $u_\emptyset=y$, we obtain $x_0$ and $x_1$. Again, if necessary, arbitrarily define $x_0(s_1)$ and $x_1(s_1)$ so that $s_1 \in \dom(x_0), \dom(x_1)$. In general, given $x_u$, apply Theorem~\ref{thm:slenderextension} to $s_{|u|+1}$ and $x_u$ to obtain $x_{u^\smallfrown 0}$ and $x_{u^\smallfrown 1}$. Again, arbitrarily define $x_{u^\smallfrown 0}$ and $x_{u^\smallfrown 1}$ on $s_{|u|+1}$ so that $s_{|u|+1} \in \dom(x_{u^\smallfrown 0}), \dom(x_{u^\smallfrown 1})$. The collection $\{x_u \: u \in 2^{<\N}\}$ has the desired properties.

(v) $\Rightarrow$ (iv). Let $y \in 2^A$, and let $\{x_u \: u \in 2^{<\N}\}$ be as in (v) with respect to any enumeration of $G - \{1_G\}$. For $\tau \in 2^\N$, let $z_\tau = \bigcup_{n \in \N} x_{\tau \upharpoonright n}$. By clause (vd), $z_\tau \in 2^G$. Each $z_\tau$ is a 2-coloring, and for $\sigma \neq \tau \in 2^\N$, $z_\sigma \bot z_\tau$. To see $\{z_\tau \: \tau \in 2^\N\}$ is a perfect set, it suffices to verify that the mapping $\sigma\mapsto z_\sigma$ is continuous, since it is obviously injective. Given $\epsilon > 0$ there corresponds a finite $B \subseteq G$ with the property that $d(w_1, w_2) < \epsilon$ if and only if $w_1$ and $w_2$ agree on $B$. By clause (vd) there is $n \in \N$ so that $B \subseteq \dom(x_{\tau \upharpoonright n})$. So $z_\tau$ and $z_\sigma$ agree on $B$ for any $\sigma \in 2^\N$ with $\sigma \upharpoonright n = \tau \upharpoonright n$. Thus the map $\sigma \mapsto z_\sigma$ is continuous and injective, so the image, $\{z_\tau \: \tau \in 2^\N\}$, is perfect.
\end{proof}

Since every proper subgroup of $G$ is slender, the following corollary is immediate.

\begin{cor} Let $G$ be a countably infinite group and $H<G$ a proper subgroup of $G$. Then every partial function with domain $H$ can be extended to a $2$-coloring on $G$.
\end{cor}

It may not be so obvious to the reader why the condition in Theorem~\ref{thm:FAAF}
implies the slenderness of $A$. Here we give a direct argument. Suppose $FA^{-1}AF\neq G$ for any finite $F\subseteq G$. In particular $A^{-1}A\neq G$. Let $a\in G-A^{-1}A$ be arbitrary. Then $a\neq 1_G$. Let $T=\{1_G,a\}$. Then $T^{-1}T\not\subseteq A^{-1}A$. We claim that $T$ witnesses the slenderness of $A$. For this let $g\in G$. Then $gT\not\subseteq A$, since otherwise $T^{-1}T=T^{-1}g^{-1}gT\subseteq A^{-1}A$, a contradiction.

By carefully analyzing the proof of Theorem \ref{thm:slenderextension}, we see that we can actually prove something stronger.

\begin{cor}
Let $G$ be a countably infinite group, and let $A \subseteq G$ be slender. Then there is a set of continuous functions $\{E_\tau \: \tau \in 2^\N\}$ such that $E_\tau(y)$ is a $2$-coloring on $G$ extending $y$ for each $y \in 2^A$ and $\tau \in 2^\N$. Moreover, if $\sigma \neq \tau \in 2^\N$, then $E_\sigma(y)$ and $E_\tau(y)$ are orthogonal.
\end{cor}

\begin{proof}
Fix an enumeration $s_1, s_2, \ldots$ of the non-identity elements of $G$. Let $Z$ denote the set of ordered triples $(y, T, i)$ where $y \in 2^{\subseteq G}$ has slender domain, $T \subseteq G$ is finite and witnesses the slenderness of $\dom(y)$ (meaning that for every $g \in G$ there is $t \in T$ with $gt \not\in \dom(y)$), and $i \in \N$ is either $0$ or else $y$ blocks $s_l$ for every $l \leq i$. We view $Z$ as a subset of $2^{\subseteq G} \times \pow_{\text{fin}}(G) \times \N$, where $\N$ has the discrete topology, $\pow_{\text{fin}}(G)$ is the collection of finite subsets of $G$ and has the discrete topology, and $2^{\subseteq G}$ has the topology coming from $3^G$ as described at the end of Section \ref{SEC BERNFLOW2}. Let $\pi_1$, $\pi_2$, and $\pi_3$ denote the projections from $Z$ to its three components. Denote by $Z_i$ the set $\pi_3^{-1}(i) \subseteq Z$. We claim that there are two continuous functions $E_0, E_1: Z \rightarrow Z$ satisfying the following:
\begin{enumerate}
\item[\rm (i)] $E_0(Z_i), E_1(Z_i) \subseteq Z_{i+1}$ for all $i \in \N$;
\item[\rm (ii)] $\{1_G, s_1, s_2, \ldots, s_{i+1}\} \subseteq \dom(\pi_1(E_0(z))), \dom(\pi_1(E_1(z)))$ for all $i \in \N$ and $z \in Z_i$;
\item[\rm (iii)] $\pi_1(E_0(z))$ and $\pi_1(E_1(z))$ both extend $\pi_1(z)$ for every $z \in Z$;
\item[\rm (iii)] any extensions of $\pi_1(E_0(z))$ and of $\pi_1(E_1(z))$ to $2^G$ are orthogonal, for every $z \in Z$.
\end{enumerate}

Before defining $E_0$ and $E_1$, we show how the statement of this corollary follows from their existence. Fix a slender set $A \subseteq G$. Let $T \subseteq G$ be finite and have the property that for every $g \in G$ there is $t \in T$ with $g t \not\in A$. Notice that the inclusion $y \in 2^A \mapsto (y, T, 0) \in Z$ is continuous. For $y \in 2^A$ and $\tau \in 2^\N$ we define
$$E_\tau(y) = \bigcup_{n \in \N} \pi_1 \circ E_{\tau(n)} \circ E_{\tau(n-1)} \circ \cdots \circ E_{\tau(0)}(y, T, 0).$$
By clauses (ii) and (iii) we have that $E_\tau(y)$ is an element of $2^G$ which extends $y$. Since $E_0$ and $E_1$ are continuous, it follows from clause (ii) that the map $y \mapsto E_\tau(y)$ is continuous for $y \in 2^A$ and fixed $\tau \in 2^\N$. Also, by clause (i) and the definition of $Z$ it follows that $E_\tau(y)$ is a $2$-coloring. Finally, by clause (iii) we have that $E_\tau(y)$ and $E_\sigma(y)$ are orthogonal for every $y \in 2^A$ and $\tau \neq \sigma \in 2^\N$.

Now we define $E_0$ and $E_1$. This is not too difficult of a task since $E_0$ and $E_1$ were in some sense implicitly defined in the proof of Theorem \ref{thm:slenderextension}. In the proof of that theorem, we began with $y \in 2^A$ with $A$ slender a non-identity $s \in G$, we picked a finite set $T \subseteq G$ witnessing the slenderness of $A$, and we described how to construct a sequence of functions, namely $y_1$, $y_2$, $y_3$, $y_4$, $y_5$, and $x$, which ultimately led to two functions $x_0$ and $x_1$. The functions $x_0$ and $x_1$ had the property that they have slender domain, they extended $y$, they block $s$, and any extensions of $x_0$ and $x_1$ to $2^G$ are orthogonal. If $s = s_{i+1}$ and $y$ blocks all $s_l$ for $l \leq i$ then we essentially want to define $E_0(y, T, i) = (x_0, T_0, i+1)$ and $E_1(y, T, i) = (x_1, T_1, i+1)$, for some finite sets $T_0, T_1 \subseteq G$ to be specified. So we must simply trace through the constructions appearing in the proof of Theorem \ref{thm:slenderextension} and show that these can be done in a continuous manner.

Let $z = (y, T, i) \in Z$. Set $A = \dom(y)$. So $A$ is slender and for every $g \in G$ there is $t \in T$ with $g t \not\in A$. Set $s = s_{i+1}$. Looking back at the proof of Theorem \ref{thm:slenderextension}, we see that many objects appearing in that proof are defined solely in terms of $T$. Namely, the following objects are defined solely in terms of $T$: $B$, $N$, $C$, $a: C \times C \rightarrow G$, $V$, $\prec$ (the total order on $2^V$), $F$, $\Lambda$, $\{\lambda_1, \lambda_2, \ldots, \lambda_M\}$, $K$, and $p: C \rightarrow 2^K$. By fixing a choice of each of these objects for every finite $T \subseteq G$, we get that these objects vary continuously with $z \in Z$ (we view these objects as points in a discrete topological space). If $x_0$ and $x_1$ are as constructed in Theorem \ref{thm:slenderextension} with respect to $z$ and $s$, then it is explicit in that proof that the finite set $C^{3N+1} F F^{-1} C^2 \lambda_M B$ witnesses the slenderness of $\dom(x_0) = \dom(x_1)$. We will define $E_0(z)$ and $E_1(z)$ so that $\pi_2(E_0(z)) = \pi_2(E_1(z)) = C^{3N+1} F F^{-1} C^2 \lambda_M B$ and $\pi_3(E_0(z)) = \pi_3(E_1(z)) = i+1$. Thus $E_0$ and $E_1$ will be continuous as long as $\pi_1 \circ E_0$ and $\pi_1 \circ E_1$ are continuous. So we only need to show that the functions $x_0$ and $x_1$ in the proof of Theorem \ref{thm:slenderextension} can be constructed continuously from $y$, $T$, and $i$.

We first show that $D_N, D_{N-1}, \ldots, D_0 \in \pow(G)$ can be defined continuously in terms of $(y, T, i)$. Here $\pow(G)$ is the power set of $G$, $\pow(G) = \{M \: M \subseteq G\}$. We identify $\pow(G)$ with $2^G$, the correspondence being given by characteristic functions, and use the corresponding topology. Fix an enumeration $g_0, g_1, \ldots$ of $G$. This enumeration of $G$ gives rise to an order on $\pow(G)$ which we will denote $\leq_{\pow(G)}$. Specifically, if $M, N \subseteq G$ then $M \leq_{\pow(G)} N$ if either $M = N$ or else $g_n \in M - N$ where $n$ is least with $g_n \in (M - N) \cup (N - M)$. Now if we define $D_N$ to be the $\leq_{\pow(G)}$-least set satisfying the conditions stated in the proof of Theorem \ref{thm:slenderextension}, then it is not difficult to see that $D_N$ depends continuously on $(y, T, i)$. Once $D_N$ has been defined in this way, we then define $D_{N-1}$ to be the $\leq_{\pow(G)}$-least set satisfying the conditions stated in the proof of Theorem \ref{thm:slenderextension}. Again, $D_{N-1}$ depends continuously on $(y, T, i, D_N)$ and hence continuously on $(y, T, i)$. Continuing in this manner, we see that $D_N, D_{N-1}, \ldots, D_0$ can be defined continuously in terms of $(y, T, i)$. In particular, the map $(y, T, i) \mapsto (y, T, i, D_N, D_{N-1}, \ldots, D_0)$ is continuous. Similar methods of choosing well orderings and picking the least element satisfying conditions given in the proof of Theorem \ref{thm:slenderextension} show that the functions $y_1$, $y_2$, and $y_3$ can be defined continuously in terms of $(y, T, i, D_N, \ldots, D_0)$.

As mentioned earlier, both $V$ and the total order $\prec$ on $2^V$ depend continuously on $T$. Thus it is easy to check that each $\Delta_r$ is a continuous function of $(y_3, T, i, D_N, \ldots, D_0)$. So $(y_3, T, i, D_N, \ldots, D_0, \Delta_N, \ldots, \Delta_0)$ depends continuously on $(y, T, i)$. Again, by picking total orderings and choosing the smallest elements relative to the appropriate conditions, we see that $y_4$ and $y_5$ depend continuously on $(y_3, T, i, D_N, \ldots, D_0, \Delta_N, \ldots, \Delta_0)$. This uses the fact that $K$ and $p: C \rightarrow 2^K$ depend continuously on $T$. Now we must consider the graph coloring of $\Gamma$, $\mu: \Delta \rightarrow (2|C|^{3N+3}|F|^2+1)$ (the vertex set of $\Gamma$ is $\Delta$). We view such functions $\mu$ as elements of $(2|C|^{3N+3}|F|^2+2)^G$ by declaring $\mu$ to be constantly $2|C|^{3N+3}|F|^2+1$ on $G - \Delta$. We define $\mu$ continuously as follows. Recall that $g_0, g_1, \ldots$ is an enumeration of $G$. We let $\mu(g_0)$ be $0$ if $g_0 \in \Delta$ and otherwise we let $\mu(g_0)$ be $2|C|^{3N+3}|F|^2+1$. In general, after $\mu$ has been defined on $g_0, g_1, \ldots, g_n$, we define $\mu(g_{n+1})$ as follows. If $g_{n+1} \not\in \Delta$ then we set $\mu(g_{n+1}) = 2|C|^{3N+3}|F|^2+1$. Otherwise we let $\mu(g_{n+1})$ be the least number among $\{0, 1, \ldots, 2|C|^{3N+2}|F|^2\}$ not appearing in the set
$$\{\mu(g_r) \: 0 \leq r \leq n, \ g_r \in \Delta, \text{ and } (g_r, g_{n+1}) \in E(\Gamma)\}.$$
This definition is valid since each vertex of $\Gamma$ has at most $2|C|^{3N+3}|F|^2$ neighbors. Recall that $(\gamma, \psi) \in E(\Gamma)$ if there is $h \in C^2 F F^{-1} C^{3N+1}$ with either $\gamma h s h^{-1} = \psi$ or $\psi h s h^{-1} = \gamma$. The set $C^2 F F^{-1} C^{3N+1}$ depends continuously on $T$, and $\Delta$ depends continuously on $(y, T, i)$. Therefore $\mu$ depends continuously on $(y, T, i)$. Now it is easy to check that $x$ and $x_0 = \pi_1(E_0(y, T, i))$ and $x_1 = \pi_1(E_1(y, T, i))$ can be defined continuously. Note that we must slightly enlarge the $x_0$ and $x_1$ appearing in the proof of Theorem \ref{thm:slenderextension} in order to have $\{1_G, s_1, s_2, \ldots, s_{i+1}\} \subseteq \dom(x_0), \dom(x_1)$. However, it is easy to see that this can be done continuously. This completes the proof.
\end{proof}

\section{Almost equality and cofinite domains} \label{SECT ALMOST EQUAL}

Recall that for $x, y \in 2^G$ we say that $x$ and $y$ are almost equal, written $x =^* y$, if $x$ and $y$ disagree on only finitely many elements of $G$. In this section we take a closer look at the behaviour of points under almost equality. We study the relationship between almost equality and periodicity, and the relationship between almost equality and $2$-colorings. From these results we will derive new findings regarding the extendability of partial functions to $2$-colorings. This section involves a substantial amount of geometric group theory - specifically the notion of the number of ends of a group and Stallings' Theorem. These geometric group theory concepts are introduced below.

We first review the notion of the number of ends of a finitely generated group. Let $G$ be a finitely generated infinite group, and let $S$ be a finite symmetric set of generators. The right Cayley graph \index{Cayley graph}\index{Cayley graph!Right Cayley Graph}\index{Cayley graph!Left Cayley Graph} of $G$ with respect to $S$, $\Gamma_{R,S}$, is the graph with vertex set $G$ and edge relation $\{(g, gs) \: g \in G, s \in S\}$. The left Cayley graph of $G$ with respect to $S$, $\Gamma_{L,S}$, is defined similarly but with edge set $\{(g, sg) \: g \in G, s \in S\}$. It is traditional to use right Cayley graphs, however $\Gamma_{L,S}$ and $\Gamma_{R,S}$ are isomorphic as graphs (the isomorphism sending $g \in V(\Gamma_{L,S})$ to $g^{-1} \in V(\Gamma_{R,S})$). The \emph{number of ends} of $G$ is defined to be the supremum of the number of infinite connected components of $\Gamma_{R,S} - A$ as $A$ ranges over all finite subsets of $G$. The number of ends of a group is in fact independent of the finite generating set used and is thus an intrinsic property of the group $G$. We recall here a well known theorem of geometric group theory.

\begin{theorem}[\cite{BH}, Theorem 8.32] \label{BH THM}
Let $G$ be a finitely generated group.
\begin{enumerate}
\item[\rm (i)] $G$ must have either $0$, $1$, $2$, or infinitely many ends.
\item[\rm (ii)] $G$ has $0$ ends if and only if $G$ is finite.
\item[\rm (iii)] $G$ has $2$ ends if and only if $G$ contains $\ZZ$ as a normal subgroup of finite index.
\item[\rm (iv)] $G$ has infinitely many ends if and only if $G$ can be expressed as an almagamated product $A *_C B$ or HNN extension $A*_C$ with $C$ finite, $[A: C] \geq 3$, and $[B: C] \geq 2$.
\end{enumerate}
\end{theorem}

As pointed out in \cite{BH}, clauses (i), (ii), and (iii) are due to Hopf while clause (iv) is a celebrated theorem of Stallings. In \cite{BH}, clause (iii) does not contain the word ``normal,'' however normality easily follows in view of Lemma \ref{lem:Poincare} of this paper.

The following theorem addresses the question of whether or not periodic points are almost equal to aperiodic points.

\begin{theorem} \label{ALMOST PERIOD}
Let $G$ be a countable group. The following are all equivalent:
\begin{enumerate}
\item[\rm (i)] there is $x \in 2^G$ such that every $y \in 2^G$ almost equal to $x$ is periodic;
\item[\rm (ii)] $G$ contains a finitely generated subgroup having infinitely many ends;
\item[\rm (iii)] $G$ contains a subgroup which is the free product of two nontrivial groups at least one of which has more than two elements;
\item[\rm (iv)] $G$ contains a nonabelian free subgroup.
\end{enumerate}
\end{theorem}

\begin{proof}
We first show that (ii), (iii), and (iv) are equivalent.

(ii) $\Rightarrow$ (iii). Let $H \leq G$ be finitely generated and have infinitely many ends. By Theorem \ref{BH THM}, $H$ can be expressed as an almagamated product $A *_C B$ or HNN extension $A*_C$ with $C$ finite, $[A: C] \geq 3$, and $[B: C] \geq 2$. First suppose that $H$ is of the form $A *_C B$. If $C$ is trivial then (iii) holds and we are done. So we may suppose that $C$ is nontrivial. We will use Lemma 6.4 of Section III.$\Gamma$.6 on page 498 of \cite{BH}. This lemma states, in particular, that if $a_1 \in A$, $a_2, a_3, \ldots, a_n \in A - C$, $b_1, \ldots, b_{n-1} \in B - C$, and $b_n \in B$, then $a_1 b_1 a_2 b_2 \cdots a_n b_n$ is not the identity element in $H = A *_C B$. Pick any $b_1, b_2 \in B - C$, and pick any $a_1 \in A - C$. Since $[A : C] \geq 3$, we can also pick $a_2 \in A - C$ with $a_2 a_1 \not\in C$. If $b_2 b_1 \in C$, then pick $a_3 \in A - C$ with $a_3 b_2 b_1 a_1 \not\in C$. If $b_2 b_1 \not\in C$, then pick any $a_3 \in A - C$. Set $u = a_1 b_1 a_2$ and $v = b_1 a_1 b_2 a_3 b_2$. Then $u$ and $v$ generate a nonabelian free group. In particular, $\ZZ * \ZZ \cong \langle u, v \rangle \leq H \leq G$ so (iii) is satisfied. Now suppose that $H$ is of the form $A*_C = \langle A, t | t^{-1} c t = \phi(c) \rangle$ where $[A : C] \geq 3$ and $\phi:C \rightarrow A$ is an injective homomorphism. Let $C' = \phi(C)$. We again use Lemma 6.4 of Section III.$\Gamma$.6 on page 498 of \cite{BH}. This lemma states, in particular, that if $a_1, a_2, \ldots, a_{n-1} \in A - (C \cup C')$ and $m_1, m_2, \ldots, m_n \in \Z - \{0\}$ then $t^{m_1} a_1 t^{m_2} a_2 \cdots t^{m_{n-1}} a_{n-1} t^{m_n}$ is not the identity element in $H = A*_C$. Since $C$ is finite, $C'$ cannot properly contain $C$ and therefore $C'$ cannot contain any coset $a C$ with $a \not\in C$. Thus there is $a \in A - C \cup C'$. Set $u = t a t$ and $v = t^2 a t^2$. Then $u$ and $v$ generate a nonabelian free group. Thus $\ZZ * \ZZ \cong \langle u, v \rangle \leq H \leq G$ so (iii) is satisfied.

(iii) $\Rightarrow$ (iv). It suffices to show that if $A$ and $B$ are nontrivial groups and $|B| > 2$, then $H = A * B$ contains a free nonabelian subgroup. By Proposition 4 on page 6 of \cite{S}, the kernel of the homomorphism $A*B \rightarrow A \times B$ is a free group with free basis $\{a^{-1} b^{-1} a b \: 1_A \neq a \in A, \ 1_B \neq b \in B\}$. So we only need to show that this free basis contains more than one element. Since $|B| > 2$, there are nonidentity $b_1 \neq b_2 \in B$. Fix any nonidentity $a \in A$. Since $H = A * B$, it is clear that $a^{-1} b_1^{-1} a b_1$ and $a^{-1} b_2^{-1} a b_2$ are distinct. Thus the free basis of the kernel contains at least two elements and therefore the kernel is a free nonabelian subgroup.

(iv) $\Rightarrow$ (ii). Clearly every nonabelian free group contains a finitely generated nonabelian free group. Also, it is easy to see that any finitely generated nonabelian free group has infinitely many ends (alternatively, one could use clause (iv) of Theorem \ref{BH THM} to see this). Thus every nonabelian free group contains a subgroup with infinitely many ends.

Now we show that the negation of (ii) implies the negation of (i). Assume that every finitely generated subgroup of $G$ has finitely many ends. Let $x \in 2^G$. Set
$$P = \{u \in G \: u \neq 1_G \text{ and } \exists y =^* x \text{ with } u \cdot y = y\}.$$
We point out that if $u \in G$ and $u \cdot x =^* x$ then we may not have that $u \in P$ (consider $G = \ZZ$ and $x \in 2^\ZZ$ given by: $x(n) = 1$ if and only if $n \geq 0$). We now point out a key observation.

\underline{Key Observation}: Consider $u_1, u_2, \ldots, u_k \in P$ and the subgroup they generate $H = \langle u_1, u_2, \ldots, u_k \rangle$. Define
$$B = \{g \in G \: \exists 1 \leq i \leq k \ x(u_i g) \neq x(g) \text{ or } x(u_i^{-1} g) \neq x(g) \}.$$
It follows from the definition of $P$ that $B$ is finite. Now put a graph structure on $G$ by using the edge set $\{(g, u_i^{\pm 1} g) \: g \in G, \ 1 \leq i \leq k\}$. Then the connected components of $G$ are precisely the right cosets of $H$ in $G$ and $x$ is constant on the connected components of $G - B$. Notice that only finitely many right cosets of $H$ intersect $B$, and therefore $x$ is constant on all but finitely many right cosets of $H$. Also notice that every connected component of $G$ (i.e. every right coset of $H$) is (graph) isomorphic to the left Cayley graph of $H$ with respect to the generating set $\{u_1^{\pm 1}, u_2^{\pm 1}, \ldots, u_k^{\pm 1}\}$. The number of ends of $H$ will therefore give us useful information about the behaviour of $x$ on the right cosets of $H$ which meet $B$.

Since $G$ does not contain any finitely generated subgroups with infinitely many ends, every finitely generated subgroup of $G$ must have either $0$, $1$, or $2$ ends. The proof proceeds in cases. The first three cases handle the scenario where $\langle P \rangle$ is finitely generated, and the last three cases handle the scenario where $\langle P \rangle$ is not finitely generated. Note that if $\langle P \rangle$ is finitely generated, then there exists a finite set $S \subseteq P$ with $\langle S \rangle = \langle P \rangle$.

\underline{Case 1:} $\langle P \rangle$ is finitely generated and has $0$ ends. Then $\langle P \rangle$ is finite and thus $P$ is finite. We define $y \in 2^G$ by
$$y(g) = \begin{cases}
1 & \text{if } g = 1_G \\
0 & \text{if } g \in P \\
x(g) & \text{otherwise}.
\end{cases}$$
Clearly $y =^* x$. If $h \in G$ and $h \cdot y = y$ then
$$y(h) = (h^{-1} \cdot y)(1_G) = y(1_G) = 1.$$
So $h \not\in P$ and therefore $h = 1_G$ and $y$ is an aperiodic point almost equal to $x$.

\underline{Case 2:} $\langle P \rangle$ is finitely generated and has $1$ end. Use $H = \langle P \rangle$ and consider the Key Observation. Using the fact that $\langle P \rangle$ is one ended, we get that $x$ is constant on a cofinite subset of $\langle P \rangle$. Let $i$ be the value of $x$ on this cofinite set and define $y \in 2^G$ by
$$y(g) = \begin{cases}
x(g) & \text{if } g \not\in \langle P \rangle \\
i & \text{if } g \in \langle P \rangle \text{ and } g \neq 1_G \\
1 - i & \text{if } g = 1_G.
\end{cases}$$
If $h \cdot y = y$ then $h \in P \cup \{1_G\}$ since $y =^* x$. Then $y(h^{-1}) = (h \cdot y)(1_G) = y(1_G) = 1-i$ so $h = 1_G$. Thus $y$ is aperiodic and is almost equal to $x$.

\underline{Case 3:} $\langle P \rangle$ is finitely generated and has two ends. If $x$ is constant on a cofinite subset of $\langle P \rangle$, then we can repeat the argument appearing in Case 2. So we may suppose that $x$ is not constant on any cofinite subset of $\langle P \rangle$. Since $\langle P \rangle$ has two ends, it follows from the Key Observation that there are disjoint infinte sets $L, R \subseteq \langle P \rangle$ (in some sense, the two ends) such that $x$ is constant on $L$ and $R$ (separately) and $L \cup R$ is cofinite in $\langle P \rangle$. By Theorem \ref{BH THM}, $\langle P \rangle$ contains $\ZZ$ as a finite index normal subgroup. Although we will not make direct use of this, we point out that intuitively the Cayley graph of $\langle P \rangle$ is ``cylinder-like'' and this accounts for the ``two ends.'' Let $v \in \langle P \rangle$ generate an infinite normal cyclic subgroup of finite index. Notice that the map $p \mapsto p v$ induces an isometry of every left Cayley graph of $\langle P \rangle$. Therefore, if $p, q \in \langle P \rangle$ then for all but finitely many $k \in \ZZ$ we have either $p v^k, q v^k \in L$ or $p v^k, q v^k \in R$. Let $a_0 = 1_G, a_1, \cdots, a_n$ be a complete set of distinct coset representatives for the cosets of $\langle v \rangle$ in $\langle P \rangle$. In other words, $a_s \langle v \rangle \cap a_t \langle v \rangle = \varnothing$ for $s \neq t$ and $\langle P \rangle = \bigcup a_s \langle v \rangle$. By swapping $L$ and $R$ if necessary, we may suppose that for all sufficiently large $k > 0$ we have
$$\{a_0, a_1, \ldots, a_n\} v^{-k} \subseteq L \text{ and } \{a_0, a_1, \ldots, a_n\} v^k \subseteq R.$$
Since $x$ is constant on $L$, we can let $i$ denote this constant value. Then $x$ is identically $1-i$ on $R$. Define $\phi, \theta \in 2^\ZZ$ by
$$\phi(k) = \begin{cases}
i & \text{if } k < 0 \\
1-i & \text{otherwise}
\end{cases} \ ; \ \ \ \ \
\theta(k) = \begin{cases}
i & \text{if } k < 0 \text{ or } k = 1 \\
1-i & \text{otherwise}.
\end{cases}$$
Notice that $[\phi]$ and $[\theta]$ are infinite and are disjoint. Now define $y \in 2^G$ by
$$y(g) = \begin{cases}
x(g) & \text{if } g \not\in \langle P \rangle \\
\theta(k) & \text{if } \exists s \neq 0 \ \exists k \in \ZZ \ g = a_s v^k \\
\phi(k) & \text{if } \exists k \in \ZZ \ g = v^k = a_0 v^k. \\
\end{cases}$$
Then $y =^* x$. Suppose $h \in G$ and $h \cdot y = y$. Then $h \in P \cup \{1_G\}$. Let $s$ and $k$ be such that $h^{-1} = a_s v^k$. Then the function
$$m \mapsto y(h^{-1} v^m) = y(a_s v^k v^m) = y(a_s v^{k+m})$$
lies in $[\phi] \cup [\theta]$. However,
$$\phi(m) = y(v^m) = (h \cdot y)(v^m) = y(h^{-1} v^m).$$
So by the pairwise disjointness of $[\phi]$ and $[\theta]$ we must have that $s = 0$ and $h^{-1} = a_0 v^k = v^k$. Since the restriction of $y$ to $\langle v \rangle$ is (essentially) $\phi$ and $v^k \cdot y = y$, we must have $k = 0$ and therefore $h = 1_G$. Thus $y =^* x$ and $y$ is aperiodic.

The three remaining cases deal with the scenario where $\langle P \rangle$ is not finitely generated. Let $H_0 \leq H_1 \leq \cdots$ be an increasing sequence of finitely generated subgroups of $\langle P \rangle$ with $\bigcup H_n = \langle P \rangle$. Now each $H_n$ has either $0$, $1$, or $2$ ends. By passing to a subsequence if necessary, we may suppose that the number of ends of $H_n$ is independent of $n \in \N$. We proceed by cases on the number of ends of the $H_n$'s.

\underline{Case 4:} Each $H_n$ has $0$ ends. Then each $H_n$ is finite and $\langle P \rangle$ is locally finite. Let $y \in 2^G$ be such that $y(g) = x(g)$ for nonidentity $g \in G$ and $y(1_G) = 1 - x(1_G)$. It will suffice to show that either $x$ or $y$ is aperiodic. Towards a contradiction, suppose $x$ and $y$ are both periodic. Let $a, b \in G$ be nonidentity group elements with $a \cdot x = x$ and $b \cdot y = y$. Then $a, b \in P$ and therefore $H = \langle a, b \rangle$ is finite. Let $\Gamma$ be the graph with vertex set $H$ and edge relation $\{(h, a^{\pm 1} h) \: h \in H\} \cup \{(h, b^{\pm 1} h) \: h \in H\}$. Since $x$ and $y$ agree on $G - \{1_G\}$ and $a \cdot x = x$ and $b \cdot y = y$, we have that $x$ must be constant on the connected components of $\Gamma - \{1_G\}$. We have $x(a) = x(1_G) \neq y(1_G) = y(b) = x(b)$. So $a$ and $b$ do not lie in the same connected component of $\Gamma - \{1_G\}$. Thus $\Gamma - \{1_G\}$ is not connected. Since $a$ and $b$ have finite order, $a$ is in the same connected component as $a^{-1}$ and similarly $b$ is in the same connected component as $b^{-1}$. Every connected component of $\Gamma - \{1_G\}$ must contain a point adjacent to $1_G$. Therefore there must be precisely two connected components of $\Gamma - \{1_G\}$, one containing $a$, call it $A$, and the other containing $b$, call it $B$. Multiplication on the right induces an automorphism of $\Gamma$, so for any $h \in H$ the connected components of $\Gamma - \{h\}$ are $A h$ and $B h$. Since $H$ is finite so are $A$ and $B$. Suppose that $|A| \leq |B|$ (the other case is identical). We have $1_G \in A a$ and therefore $1_G \not\in B a$. However $B a \cup \{a\}$ is connected, hence connected in $\Gamma - \{1_G\}$. Since $a \in B a \cup \{a\}$, we must have that $B a \cup \{a\} \subseteq A$. Therefore $|B| + 1 \leq |A|$, contradicting $|A| \leq |B| < \infty$. Thus either $x$ or $y$ is aperiodic. In any case, $x$ is almost equal to an aperiodic element of $2^G$.

\underline{Case 5:} Each $H_n$ has $1$ end. Since each $H_n$ has $1$ end, each $H_n$ must be infinite. By the Key Observation, we have that for every $n \in \N$ and every right coset $H_n a$ of $H_n$, $x$ is constant on a cofinite subset of $H_n a$. We claim that $x$ is constant on a cofinite subset of $\langle P \rangle$. Notice that $x$ is constant on all but finitely many of the right cosets of $H_0$ (although this constant value may change from coset to coset). So if $x$ is not constant on any cofinite subset of $\langle P \rangle$, then there must be right cosets $H_0 a, H_0 b \subseteq \langle P \rangle$ such that $x$ takes the value $0$ infinitely many times on $H_0 a$ and $x$ takes the value $1$ infinitely many times on $H_0 b$. Since $\bigcup H_m = \langle P \rangle$, there is $n \in \N$ with $H_0 a, H_0 b \subseteq H_n$. But then $x$ takes the value $0$ and the value $1$ infinitely many times on $H_n$, contradicting $x$ being constant on a cofinite subset of $H_n$. Thus $x$ must be constant on a cofinite subset of $\langle P \rangle$. Let $i \in \{0, 1\}$ be this constant value. Define $y \in 2^G$ by $y(1_G) = 1 - i$, $y(p) = i$ for $1_G \neq p \in \langle P \rangle$, and $y(g) = x(g)$ for $g \not\in \langle P \rangle$. Then $y =^* x$. If $h \in G$ and $h \cdot y = y$ then by definition we have $h \in P \cup \{1_G\}$. So $y(h^{-1}) = (h \cdot y)(1_G) = y(1_G)$ and therefore $h = 1_G$ due to how $y$ was defined on $\langle P \rangle$. So $y$ is aperiodic and is almost equal to $x$.

In Case 6 below, we consider the final possible scenario in which each $H_n$ has two ends. Before handling this case, we make the general claim that $H_0$ must be of finite index within each $H_n$. To see this, for each $n$ let $K_n$ be a normal finite index subgroup of $H_n$ isomorphic to $\ZZ$ (see Theorem \ref{BH THM}). Then $(H_0 K_n) / K_n$ is a subgroup of the finite group $H_n / K_n$, and thus $(H_0 K_n) / K_n$ is finite. By the isomorphism theorems of group theory, we have that $H_0 / (H_0 \cap K_n) \cong (H_0 K_n) / K_n$ is finite. So $H_0 \cap K_n$ has finite index in $H_0$. However, $H_0 \cap K_n$ is a subgroup of $K_n \cong \ZZ$ and therefore must have finite index within $K_n$. Thus $H_0 \cap K_n$ has finite index in $H_n$, so in particular $H_0$ has finite index in $H_n$. This proves the claim. We now continue to Case 6.

\underline{Case 6:} Each $H_n$ has $2$ ends. Since each $H_n$ has $2$ ends, each $H_n$ must be infinite. As in Case 3, by Theorem \ref{BH THM} and the claim above there is $v \in H_0$ of infinite order such that $\langle v \rangle$ has finite index in each $H_n$. For each $n \in \N$ there are infinite sets $L_n, R_n \subseteq H_n$ (in some sense corresponding to the two ends of $H_n$) with $L_n \cap R_n = \varnothing$, $L_n \cup R_n$ cofinite within $H_n$, and $x$ constant on each of $L_n$ and $R_n$. By swapping $L_n$ and $R_n$ if necessary, we can assume that for each $h \in H_n$ there is $m \in \N$ such that $h v^{-k} \in L_n$ and $h v^k \in R_n$ whenever $k \geq m$ (see Case 3). It follows from this last property that $L_n \cap L_{n+1}$ and $R_n \cap R_{n+1}$ are nonempty. Since $x$ is constant on each $L_n$ and $L_n \cap L_{n+1} \neq \varnothing$, it follows that $x$ is constant on $L_\infty = \bigcup L_n$ and similarly $x$ is constant on $R_\infty = \bigcup R_n$. We first show that $x$ is constant on $L_\infty \cup R_\infty$. For this it suffices to show that $x$ is constant on $L_n \cup R_n$ for some $n \in \N$. By the Key Observation we know that $x$ is constant on all but finitely many of the right cosets of $H_0$ in $G$. Since $H_0$ is finitely generated and $\langle P \rangle$ is not, we must have that $H_0$ is of infinite index in $\langle P \rangle$. So there is a right coset $H_0 a$ of $H_0$ in $\langle P \rangle$ on which $x$ is constant. Since $\langle P \rangle$ is the increasing union of the $H_n$'s, there is $n \in \N$ with $H_0 a \subseteq H_n$. Now $\langle v \rangle$ has finite index in $H_n$, so there is $k \in \N$ with $\langle v^k \rangle \lhd H_n$ (Lemma \ref{lem:Poincare}). So $H_0 a \langle v^k \rangle = H_0 \langle v^k \rangle a = H_0 a$ and therefore we have $L_n \cap H_0 a \neq \varnothing$ and $R_n \cap H_0 a \neq \varnothing$. Since $x$ is constant on $L_n$, $H_0 a$, and $R_n$, $x$ must be constant on $L_n \cup R_n$. So for each $n \in \N$ $x$ is constant on $L_n \cup R_n$, and $L_n \cup R_n$ is cofinite in $H_n$. Since every right coset of $H_n$ in $\langle P \rangle$ is contained in some $H_m$, it follows that $x$ is constant on a cofinite subset of each right coset of each $H_n$ in $\langle P \rangle$. An argument identical to that appearing in Case 5 now shows that $x$ is constant on a cofinite subset of $\langle P \rangle$. The construction in Case 5 then shows that there is an aperiodic $y \in 2^G$ with $y =^* x$.

To finish the proof, we show that (iv) implies (i). Suppose that $G$ contains a free nonabelian subgroup. Since every free nonabelian group contains a free nonabelian group of countably infinite rank, there is a free nonabelian subgroup $H \leq G$ of countably infinite rank. Let $a_0, a_1, a_2, \ldots$ be a free generating set for $H$. Let $C \subseteq G$ be a complete set of representatives of the right cosets of $H$ in $G$. Specifically, $G = H C$ and for $c \neq d \in C$ we have $H c \cap H d = \varnothing$. Let $F$ be the set of all functions $y: A \rightarrow 2$ where $A \subseteq G$ is finite. Then $F$ is countable, and we can find an injective function $\sigma: F \rightarrow \N$ such that for every $y \in F$ $\dom(y) \subseteq \langle a_0, a_1, \ldots, a_{\sigma(y)} \rangle C$. Now recursively define $x \in 2^G$ as follows. First set $x(c) = 0$ for all $c \in C$. Now fix $g \in G$, and let $h \in H$ be such that $g \in h C$. Let $i \in \N$ be such that the reduced word representation of $h$ begins on the left with $a_i^k$, where $k \neq 0$ and $|k|$ is maximized. If there is $y \in F$ with $a_i^{-k} g \in \dom(y)$ and $\sigma(y) = i-1$, then set $x(g) = y(a_i^{-k} g)$ (this is well defined since $\sigma$ is injective). Otherwise, set $x(g) = x(a_i^{-k} g)$ (this case is where the recursion occurs). Clearly $a_0 \cdot x = x$. So $x$ is periodic. Now suppose $z$ is almost equal to $x$. Then there is $y \in F$ such that $z(g) = y(g)$ for $g \in \dom(y)$ and $z(g) = x(g)$ for $g \not\in \dom(y)$. Let $i = \sigma(y)+1$. We claim that $a_i \cdot z = z$. It suffices to show that for every $k \in \ZZ$, every $c \in C$, and every $h \in H$ not beginning with $a_i$ or $a_i^{-1}$ we have $z(a_i^k h c) = z(h c)$. Fix such $k$, $c$, and $h$. Then $a_i^k h c \not\in \dom(y)$ since $\sigma(y) = i - 1$. If $hc \not\in \dom(y)$ then $z(hc) = x(hc) = x(a_i^k hc) = z(a_i^k hc)$ (the second equality follows from the definition of $x$, for which the second case applies). Similarly, if $hc \in \dom(y)$ then $z(hc) = y(hc) = x(a_i^k hc) = z(a_i^k hc)$ (by the definition of $x$ in the first case). Thus $a_i \cdot z = z$.
\end{proof}

A famous problem in group theory was the von Neumann Conjecture which states that a group is nonamenable if and only if it contains a free nonabelian subgroup. The von Neumann Conjecture was disproven by Olshanskii \cite{O} in 1981. However, the above theorem provides us with a dynamical characterization for which groups contain free nonabelian subgroups.

In practice, it may be useful to know for which groups the equivalent properties listed in the previous theorem fail. The following corollary addresses this issue.

\begin{cor} \label{COR ALMOST PERIOD}
Let $\mathcal{C}$ be the class of countable groups $G$ which have the property that for every periodic $x \in 2^G$ there is an aperiodic $y \in 2^G$ with $y =^* x$. Then $\mathcal{C}$ coincides with the class of countable groups which do not contain free nonabelian subgroups. Moreover, $\mathcal{C}$ contains all countable amenable groups and all countable torsion groups and is closed under the operations of: extensions; increasing unions; passing to subgroups; and taking quotients.
\end{cor}

By the operation of extension, we mean that if $G$ is a countable group, $K \lhd G$, and both $K, G / K \in \mathcal{C}$, then $G \in \mathcal{C}$.

\begin{proof}
The fact that $\mathcal{C}$ coincides with the class of countable groups which do not contain free nonabelian subgroups follows immediately from the previous theorem. Clearly torsion groups cannot contain free subgroups, and it is well known that amenable groups cannot contain any free nonabelian subgroups. Thus these groups are members of $\mathcal{C}$. It is simple to see that $\mathcal{C}$ is closed under the operations of increasing unions and passing to subgroups. It is also closed under taking quotients, since if a quotient $G / K$ of $G$ contains a free group, then one can pick any two elements $g_1 K, g_2 K \in G / K$ which generate a rank two free group and see that $g_1$ and $g_2$ generate a rank two free group in $G$. Finally, suppose that $G$ is a countable group and $K \lhd G$ satisfies $K, G / K \in \mathcal{C}$. Towards a contradiction, suppose that $G \not\in \mathcal{C}$. Then $G$ contains a free nonabelian subgroup $H$. Since $H \cap K$ is a normal subgroup of the free group $H$, it must be either trivial or free and nonabelian (if it is free and abelian then it cannot be normal). Since $K \in \mathcal{C}$, we must have that $H \cap K$ is trivial. It follows then that $H$ embeds into $G / K$, contradicting $G / K \in \mathcal{C}$. We conclude that $G \in \mathcal{C}$ and $\mathcal{C}$ is closed under extensions.
\end{proof}

Our argument within Case 4 of the previous theorem suggests defining and studying a new, but related, notion. For $x, y \in 2^G$ we write $y =^{**} x$ \index{$=^{**}$} if $x$ and $y$ agree everywhere on $G$ except at \emph{precisely} one point (so $x \neq^{**} x$). In the theorem below, we study how much of the previous theorem still holds when we work with $=^{**}$ instead of almost equality.

\begin{theorem} \label{STRONG ALMOST PERIOD}
Let $G$ be a countable group. The following are equivalent:
\begin{enumerate}
\item[\rm (i)]  if $x \in 2^G$ is periodic, then every $y =^{**} x$ is aperiodic;
\item[\rm (ii)] $G$ does not contain any subgroup which is a free product of nontrivial groups.
\end{enumerate}
\end{theorem}

Notice the slight difference between clause (ii) of this theorem and clause (iii) of Theorem \ref{ALMOST PERIOD}.

\begin{proof}
$\neg$(ii) $\Rightarrow$ $\neg$(i). By the previous theorem, we are done if $G$ contains a subgroup which is a free product of two nontrivial groups one of which has more than two elements. So all that remains is to handle the case where $\ZZ_2 * \ZZ_2 \leq G$. Let $a, b \in G$ be involutions (meaning $a^2 = b^2 = 1_G$) with $\langle a, b \rangle \cong \ZZ_2 * \ZZ_2$. Set $H = \langle a, b \rangle$. Notice that every nonidentity element of $H$ can be written uniquely in one of the following four forms:
$$ababab \cdots aba, \ \ bababa \cdots baba, \ \ ababab \cdots abab, \ \ bababa \cdots bab.$$
We say $h \in H$ ends with $a$ if $h$ can be written in one of the two forms on the left. Otherwise, we say $h$ ends with $b$. Define $x \in 2^G$ by
$$x(g) = \begin{cases}
0 & \text{if } g \not\in H \\
0 & \text{if } g \in H \text{ and } g \text{ ends with } a \\
1 & \text{otherwise}.
\end{cases}$$
So if $x(g) = 1$ then $g \in H$ and either $g = 1_G$ or else $g$ ends with $b$. It is easy to see that $b \cdot x = x$. Now let $y \in 2^G$ be equal to $x$ everywhere except have the opposite value at $1_G$. Then it is again easy to check that $a \cdot y = y$.

$\neg$(i) $\Rightarrow$ $\neg$(ii). Now suppose that $x, y \in 2^G$ are both periodic and $x =^{**} y$. We must show that $G$ contains a subgroup which is the free product of two nontrivial groups. By replacing $(x,y)$ with $(g \cdot x, g \cdot y)$ if necessary, we can assume that $1_G$ is the single point at which $x$ and $y$ disagree. Let $A, B \leq G$ be the stabilizer subgroups of $x$ and $y$, respectively. Set $H = \langle A \cup B \rangle$. We claim that $H \cong A * B$. Fix generating sets (possibly infinite) $S_A$ and $S_B$ for $A$ and $B$, respectively. We also require that $1_G \not\in S_A \cup S_B$, $S_A = S_A^{-1}$, and $S_B = S_B^{-1}$. Then $S_A \cup S_B$ generates $H$. Let $\Gamma$ be the left Cayley graph of $H$ associated to the generating set $S_A \cup S_B$. Since $A$ and $B$ are the stabilizers of $x$ and $y$ and $x$ and $y$ differ only at $1_G$, it follows that $x$ and $y$ are constant and agree on each of the connected components of $\Gamma - \{1_G\}$. Let $C_A$ be the union of those connected components of $\Gamma - \{1_G\}$ which contain an element of $A$. Similarly define $C_B$. If $a \in A - \{1_G\}$ and $b \in B - \{1_G\}$ then
$$x(a) = x(1_G) \neq y(1_G) = y(b) = x(b).$$
From the above inequality we have that $A \cap B = \{1_G\}$ and $S_A \cap S_B = \varnothing$. More importantly, since $x$ is constant on the connected components of $\Gamma - \{1_G\}$, it follows that $C_A \cap C_B = \varnothing$. In particular, $B \cap C_A = \varnothing$ and $A \cap C_B = \varnothing$. Also notice that every connected component of $\Gamma - \{1_G\}$ must contain a vertex adjacent to $1_G$ and therefore $C_A \cup C_B = H - \{1_G\}$.

Since $H = \langle A \cup B \rangle$, there is a surjective homomorphism $\phi: A * B \rightarrow H$ with $\phi(a) = a$ and $\phi(b) = b$ for all $a \in A$ and $b \in B$. We claim that $\phi$ is an isomorphism. Define the length of $k \in A * B$ to be $0$ if $k$ is the identity and to be
$$\min\{n \: k = c_1 c_2 \cdots c_n, \ \forall j \ c_j \in S_A \cup S_B\}$$
otherwise. Towards a contradiction, suppose the kernel of $\phi$ is nontrivial. Let $k$ be a nontrivial element of the kernel with minimum length. Clearly $k \not\in A \cup B$. Let
$$k = c_1 c_2 \cdots c_n$$
be a minimal length representation of $k$. Then the sequence
$$1_G, \ \phi(c_n), \ \phi(c_{n-1} c_n), \ \ldots, \ \phi(c_1 c_2 \cdots c_n)$$
is a closed path in $\Gamma$ beginning and ending at $1_G$. Therefore (by minimality of the length of the expression), the non-endpoint vertices traversed by this path lie in a single connected component of $\Gamma - \{1_G\}$. Therefore we must have either $c_1, c_n \in S_A$ or $c_1, c_n \in S_B$ (but not both since $S_A \cap S_B = \varnothing$). We will consider the case $c_1, c_n \in S_A$. The argument for the other case is identical. Let $m \leq n$ be minimal with $c_m, c_{m+1}, \ldots, c_n \in S_A$. Then $m > 1$ since $k \not\in A$. By conjugating $k$ by $c_m c_{m+1} \cdots c_n$ we get
$$k' = c_m c_{m+1} \cdots c_n c_1 c_2 \cdots c_{m-1}.$$
Since $k'$ is a conjugate of $k$ and the kernel is normal, $k'$ must also lie in the kernel. Also, since $k'$ is a conjugate of $k$ it must be nontrivial. As $k$ was chosen to have minimal length, it must be that the length of $k'$ is greater than or equal to the length of $k$. However, the above representation of $k'$ shows that the length of $k'$ is less than or equal to the length of $k$. So $k$ and $k'$ must have the same length, and the above representation of $k'$ must be of minimal length. Now $k'$ and its representation above have all of the same properties which we assumed of $k$ and its represenation. Therefore, arguing just as we did before, we must have that either $c_m, c_{m-1} \in S_A$ or $c_m, c_{m-1} \in S_B$. Since $c_m \in S_A$ and $S_A \cap S_B = \varnothing$, we must have $c_{m-1} \in S_A$. This contradicts the definition of $m$. We conclude that $\phi$ is injective and is thus an isomorphism. So $A * B \cong H \leq G$ as claimed.
\end{proof}

In the corollary below, $\Stab(x)$ denotes the stabilizer subgroup of $x \in 2^G$.

\begin{cor} \label{COR FREE PROD}
Let $G$ be a countable group and let $x, y \in 2^G$. If $x =^{**} y$ then $\langle \Stab(x) \cup \Stab(y) \rangle \cong \Stab(x) * \Stab(y)$.
\end{cor}

\begin{proof}
Set $A = \Stab(x)$ and $B = \Stab(y)$. The claim is trivial if either $A = \{1_G\}$ or $B = \{1_G\}$. So suppose $A$ and $B$ are nontrivial. Let $g \in G$ be the unique group element with $x(g) \neq y(g)$. Set $x' = g^{-1} \cdot x$ and $y' = g^{-1} \cdot y$. Then $x'(1_G) \neq y'(1_G)$, $g^{-1} A g = \Stab(x')$, and $g^{-1} B g = \Stab(y')$. Then
$$\langle A \cup B \rangle = g \langle g^{-1} A g \cup g^{-1} B g \rangle g^{-1} \cong \langle g^{-1} A g \cup g^{-1} B g \rangle \cong g^{-1} A g * g^{-1} B g \cong A * B.$$
where the second $\cong$ follows from the proof of the previous theorem.
\end{proof}

The following corollary summarizes the two previous theorems of this section.

\begin{cor}
Let $G$ be a countable group.
\begin{enumerate}
\item[\rm (i)] If $G$ does not contain any subgroup which is a free product of nontrivial groups, then for every periodic $x \in 2^G$ every $y =^{**} x$ is aperiodic.
\item[\rm (ii)] If $G$ contains $\ZZ_2 * \ZZ_2$ as a subgroup and if every subgroup of $G$ which is the free product of two nontrivial groups is isomorphic to $\ZZ_2 * \ZZ_2$, then for every periodic $x \in 2^G$ there is an aperiodic $y \in 2^G$ with $y =^* x$ but there are periodic $w, z \in 2^G$ with $w =^{**} z$.
\item[\rm (iii)] If $G$ contains a subgroup which is the free product of two nontrivial groups one of which has more than two elements, then there is a periodic $x \in 2^G$ such that every $y =^*x$ is also periodic.
\end{enumerate}
\end{cor}

\begin{cor}
Let $G$ be a countable group. For every periodic $x \in 2^G$ every $y =^{**} x$ is aperiodic if $G$ is:
\begin{enumerate}
\item[\rm (i)] finite;
\item[\rm (ii)] locally finite;
\item[\rm (iii)] a torsion group;
\item[\rm (iv)] amenable and does not contain $\ZZ_2 * \ZZ_2$.
\end{enumerate}
\end{cor}

\begin{proof}
By the previous theorem, it suffices to show that $G$ does not contain any subgroups which are free products of nontrivial groups. Free products of nontrivial groups are infinite, have finitely generated subgroups which are infinite, and have elements of infinite order. Therefore finite groups, locally finite groups, and torsion groups must have the stated property. As stated previously, it is well known that amenable groups do not contain free nonabelian subgroups. Therefore, by Theorem \ref{ALMOST PERIOD}, the only possible subgroup of an amenable group which is a free product is $\ZZ_2 * \ZZ_2$. Therefore any amenable groups not containing $\ZZ_2 * \ZZ_2$ must have the stated property.
\end{proof}

\begin{cor}
If $G$ is a finite group and $k > 1$ is an integer then $k^G$ contains at least $(k-1) k^{|G|-1}$ many aperiodic points.
\end{cor}

\begin{proof}
We first point out as we did at the beginning of this paper that all of our results for $2^G$ immediately generalize to $k^G$ for all $k > 1$ (only minor changes need to be made to all of the proofs in this paper). So if $G$ is finite, $k > 1$ is an integer, and $x, y \in k^G$ differ at precisely one coordinate, then $x$ and $y$ cannot both be periodic. Set $A = G - \{1_G\}$. Then $|k^A| = k^{|G|-1}$. Each element of $k^A$ can be extended in $k$ different ways to an element of $k^G$. If $x \in k^A$, then at most one extension of $x$ to $k^G$ is periodic. Therefore at least $k-1$ extensions of $x$ are aperiodic. Since the extensions of $x, y \in k^A$ are distinct if $x \neq y$, we have that $k^G$ contains at least $(k-1)k^{|G|-1}$ many aperiodic elements.
\end{proof}

Now that we better understand how periodic points behave under almost equality, we can finally prove that every near $2$-coloring is an almost $2$-coloring. In fact, we prove something much stronger.

\begin{theorem} \label{NEAR IS ALMOST}
Let $G$ be a countable group and let $x$ be a near $2$-coloring. Then either $x$ is a $2$-coloring or else every $y =^{**} x$ is a $2$-coloring.
\end{theorem}

\begin{proof}
Suppose $x$ is not a $2$-coloring. Fix $y =^{**} x$ and towards a contradiction suppose that $y$ is not a $2$-coloring. Since $x$ and $y$ are near $2$-colorings but not $2$-colorings, they must be periodic by clause (d) of Lemma \ref{lem:almostcoloringlemma}. Let $a, b \in G$ be nonidentity elements with $a \cdot x = x$ and $b \cdot y = y$. Set $H = \langle a, b \rangle$. Then $H \cong \langle a \rangle * \langle b \rangle$ by Corollary \ref{COR FREE PROD}. Set $h = ab \in H$ and notice that $h$ has infinite order. Let $p \in G$ be the unique element with $x(p) \neq y(p)$. Set
$$B = \{p, b^{-1} p\}.$$
Let $A, T \subseteq G$ be finite with the property that
$$\forall g \in G - A \ \exists t \in T \ x(ght) \neq x(gt).$$
Now pick $g \in \langle h \rangle$ with
$$g \not\in A \cup B T^{-1}.$$
Fix $t \in T$. Then $g t \not\in B$. Since $b \cdot y = y$ and $gt \neq p \neq b gt$, we have
$$x(b g t) = y(b g t) = y(g t) = x(g t).$$
Since $a \cdot x = x$ we have
$$x(h g t) = x(a b g t) = x(b g t) = x(g t).$$
However, $g \in \langle h \rangle$ so
$$x(g h t) = x(h g t) = x(g t).$$
Since $t \in T$ was arbitrary and $g \not\in A$ this contradicts $x$ being a near $2$-coloring.
\end{proof}

\begin{cor}
For every countable group $G$, every near $2$-coloring is an almost $2$-coloring.
\end{cor}

\begin{cor} \label{NEAR IS ALMOST SMRY}
For a countable group $G$ and $x \in 2^G$, the following are equivalent:
\begin{enumerate}
\item[\rm (i)] either $x$ is a $2$-coloring or else every $y \in 2^G$ differing from $x$ on precisely one coordinate is a $2$-coloring;
\item[\rm (ii)] there is a $2$-coloring $y \in 2^G$ which differs from $x$ on at most one coordinate;
\item[\rm (iii)] there is a $2$-coloring $y \in 2^G$ which differs from $x$ on finitely many coordinates;
\item[\rm (iv)] $x$ is an almost $2$-coloring;
\item[\rm (v)] $x$ is a near $2$-coloring;
\item[\rm (vi)] for every nonidentity $s \in G$ there are finite sets $A, T \subseteq G$ so that for all $g \in G - A$ there is $t \in T$ with $x(gt) \neq x(gst)$;
\item[\rm (vii)] every limit point of $[x]$ is aperiodic.
\end{enumerate}
\end{cor}

\begin{proof}
The implications (i) $\Rightarrow$ (ii) and (ii) $\Rightarrow$ (iii) are obvious. (iii) $\Rightarrow$ (iv) is by definition. (iv) $\Rightarrow$ (v) is clause (c) of Lemma \ref{lem:almostcoloringlemma}. (v), (vi), and (vii) are equivalent by Lemma \ref{NEAR EQUIV}. (v) $\Rightarrow$ (i) is Theorem \ref{NEAR IS ALMOST}.
\end{proof}

These results lead to an alternative proof of the density of $2$-colorings (the original proof was given in Theorem \ref{thm:density}).

\begin{cor}
For every countably infinite group $G$, the collection of $2$-colorings on $G$ is dense in $2^G$.
\end{cor}

\begin{proof}
By Theorem \ref{GEN COL}, there exists a $2$-coloring $x$ on $G$. Let $y \in 2^G$ and $\epsilon > 0$. Let $r \in \N$ be such that $2^{-r} < \epsilon$, and let $g_0, g_1, \ldots$ be the enumeration of $G$ used in defining the metric $d$ on $2^G$. Define $x'$ by
$$x'(g) = \begin{cases}
y(g) & \text{if } g = g_i \text{ for } 0 \leq i \leq r \\
x(g) & \text{otherwise}.
\end{cases}$$
Then $d(x', y) < 2^{-r} < \epsilon$. If $x'$ is a $2$-coloring, then we are done. If $x'$ is not a $2$-coloring, then the function $x'' \in 2^G$ defined by
$$x''(g) = \begin{cases}
x'(g) & \text{if } g \neq g_{r+1} \\
1 - x'(g) & \text{if } g = g_{r+1}.
\end{cases}$$
is a $2$-coloring by Theorem \ref{NEAR IS ALMOST} (since $x'$ is an almost $2$-coloring, in particular a near $2$-coloring). Also, $d(x'', y) < \epsilon$.
\end{proof}

Now we can further characterize extendability of partial functions to $2$-colorings.

\begin{cor} \label{COR COFINITE}
Let $y \in 2^{\subseteq G}$ be a partial function with cofinite domain. The following are equivalent:
\begin{enumerate}
\item[\rm (i)] there is a $2$-coloring $x \in 2^G$ extending $y$;
\item[\rm (ii)] for every $x_0, x_1 \in 2^G$ extending $y$ with $x_0 =^{**} x_1$, either $x_0$ or $x_1$ is a $2$-coloring;
\item[\rm (iii)] every element of $\overline{[y]} \cap 2^G$ is aperiodic.
\end{enumerate}
\end{cor}

\begin{proof}
(i) $\Rightarrow$ (ii). Let $x \in 2^G$ be a $2$-coloring extending $y$, and let $x_0, x_1 \in 2^G$ extend $y$ with $x_0 =^{**} x_1$. Since $y$ has cofinite domain, $x_0$ almost equals $x$ and hence is an almost $2$-coloring. If $x_0$ is not a $2$-coloring, then by the previous theorem $x_1$ is a $2$-coloring.

(ii) $\Rightarrow$ (iii). Let $x_0, x_1 \in 2^G$ extend $y$ with $x_0 =^{**} x_1$. Without loss of generality, we may suppose that $x_0$ is a $2$-coloring on $G$. Then clearly $\overline{[y]} \cap 2^G \subseteq \overline{[x_0]}$ since $x_0$ extends $y$. Since $x_0$ is a $2$-coloring, every element of $\overline{[x_0]}$ is aperiodic.

(iii) $\Rightarrow$ (i). If $x \in 2^G$ extends $y$, then every limit point of $[x]$ lies in $\overline{[y]} \cap 2^G$. Thus every limit point of $[x]$ is aperiodic, and so $x$ is a near $2$-coloring by Lemma \ref{NEAR EQUIV} (or by Corollary \ref{NEAR IS ALMOST SMRY} above). By the previous theorem, if $x$ is not a $2$-coloring then we can change the value of $x$ at a point $g \not\in \dom(y)$ to get a $2$-coloring $x'$ extending $y$.
\end{proof}

So far in this chapter we have studied, among other things, under what conditions on $A \subseteq G$ and $y: A \rightarrow 2$ we can extend $y$ to a $2$-coloring on $G$. Although we were unable to answer this question in fullest generality, we were able to answer it for certain subsets $A \subseteq G$ (specifically for sets $A$ which are either slender or cofinite). In addressing the general question of which partial functions can be extended to $2$-colorings, we make the following conjecture. 

\begin{conj}
Let $G$ be a countable group, let $A \subseteq G$, let $k > 1$ be an integer, and let $y: A \rightarrow k$. Then $y$ can be extended to a $k$-coloring if and only if $\overline{[y]} \cap k^G$ consists of aperiodic points.
\end{conj}

If $y$ can be extended to a $k$-coloring, then it is easy to see that $\overline{[y]} \cap k^G$ consists of aperiodic points. The difficult question to resolve is if this condition is sufficient. Clearly this conjecture implies Corollary \ref{COR COFINITE}. Also, if $A$ is slender and $y$ is as above, then $\overline{[y]} \cap k^G$ must be empty. Thus the implication (i) $\Rightarrow$ (iii) appearing in Theorem \ref{THM EXT SLENDER} also follows from the above conjecture. We would like to emphasize that in all of the results of this paper, the obvious necessary conditions have always been sufficent. This is the main reason that we formally make this conjecture.

\section{\label{sec:autoext}Automatic extendability}

We have shown in the last section that any partial function on a proper subgroup of a countably infinite group can be extended to a $2$-coloring of the full group. In this section we consider a curious question: when is it the case that any extension of
a $2$-coloring on a subgroup is automatically a $2$-coloring of the full group?

We know from early on that this happens to $\mathbb{Z}$ (see the discussion following Definition~\ref{def:block}). The goal of this section is to determine all countable groups with this automatic extendability property. It will turn out that $\mathbb{Z}$ is the only group with this property.

\begin{prop}
Let $G$ be a countably infinite group and let $H \leq G$ be nontrivial. The following are equivalent:
\begin{enumerate}
\item[\rm (i)] If $x \in 2^H$ is a 2-coloring on $H$ and $y \in 2^G$ extends $x$, then $y$ is a 2-coloring on $G$;
\item[\rm (ii)] $|G:H|$ is finite and for every nonidentity $g \in G$, $\langle g \rangle \cap H \neq \{1_G\}$.
\end{enumerate}
\end{prop}

\begin{proof}
(i) $\Rightarrow$ (ii). Towards a contradiction, suppose $|G:H|$ is infinite. Let $H, a_1 H, a_2 H, \ldots$ be an enumeration of the left cosets of $H$ in $G$, and let $x \in 2^H$ be a 2-coloring. Extend $x$ to $y \in 2^G$ by defining $y(g) = 0$ for all $g \in G - \dom(x)$. Then $\lim a_n^{-1} \cdot y = 0$, a contradiction. We conclude $|G:H|$ is finite.

Again, towards a contradiction suppose $a \in G - \{1_G\}$ satisfies $\langle a \rangle \cap H = \{1_G\}$. Let $x \in 2^H$ be a 2-coloring on $H$. Extend $x$ to $y \in 2^G$ by defining $y(a^n h) = x(h)$ for all $n \in \mathbb{Z}$ and $h \in H$, and set $y(g) = 0$ for all other $g \in G$. Then it easy to see that $y$ is well-defined and periodic: $a \cdot y = y$. This is a contradiction.

(ii) $\Rightarrow$ (i). Let $x \in 2^H$ be a 2-coloring on $H$, and let $y \in 2^G$ extend $x$. It is enough to show that every $w \in \overline{[y]}$ is not periodic. Fix $w \in \overline{[y]}$. Let $a_0 H = H, a_1 H, \ldots, a_n H$ be an enumeration of all left cosets of $H$ in $G$. Let $(g_m)_{m \in \N}$ be a sequence of elements of $G$ with $w = \lim g_m \cdot y$. By passing to a subsequence if necessary, we can assume that each $g_m$ lies in the same left coset of $H$, say $a_i H$. For each $m \in \N$ let $h_m \in H$ be such that $g_m = a_i h_m$. Then
$$w = \lim_{m \rightarrow \infty} g_m \cdot y = \lim_{m \rightarrow \infty} (a_i h_m) \cdot y$$
$$= \lim_{m \rightarrow \infty} a_i \cdot (h_m \cdot y) = a_i \cdot (\lim_{m \rightarrow \infty} h_m \cdot y) = a_i \cdot z$$
where $z = \lim h_m \cdot y$. We must show that $w = a_i \cdot z$ is not periodic, so it will suffice to show that $z$ is not periodic. Suppose $g \in G$ satisfies $g \cdot z = z$. Clearly $\lim h_m \cdot x \subseteq z$, so that $z \upharpoonright H$ is a 2-coloring on $H$. In particular, $z \upharpoonright H \in 2^H$ is not periodic. Thus $\langle g \rangle$ must intersect $H$ trivially, from which it follows that $g = 1_G$. We conclude $y$ is a 2-coloring on $G$.
\end{proof}

\begin{theorem}
Let $G$ be a countably infinite group. The following are equivalent:
\begin{enumerate}
\item[\rm (i)] If $H \leq G$ is any nontrivial subgroup, $x \in 2^H$ is any 2-coloring on $H$, and $y \in 2^G$ is any extension of $x$, then $y$ is a 2-coloring on $G$;
\item[\rm (ii)] $G = \mathbb{Z}$.
\end{enumerate}
\end{theorem}

\begin{proof}
(ii) $\Rightarrow$ (i) is clear from the previous proposition.

(i) $\Rightarrow$ (ii). By the previous proposition, every nontrivial subgroup of $G$ has finite index. The result now follows from the next proposition.
\end{proof}

Recall that $\Z(G)$ denotes the center of $G$.

\begin{prop}
If $G$ is an infinite group and every nontrivial subgroup of $G$ has finite index, then $G = \ZZ$.
\end{prop}

\begin{proof}
We first show that $G$ posseses the following properties:
\begin{enumerate}
\item [\rm (i)] $G$ is countable;
\item [\rm (ii)] every nonidentity element of $G$ has infinite order;
\item [\rm (iii)] for every $g \in G - \{1_G\}$, there is $k > 0$ with $\langle g^k \rangle \lhd G$;
\item [\rm (iv)] every normal cyclic subgroup of $G$ is contained in $\Z(G)$;
\item [\rm (v)] $\Z(G)$ is isomorphic to $\mathbb{Z}$;
\item [\rm (vi)] $G / \Z(G)$ contains no nontrivial abelian normal subgroups;
\item [\rm (vii)] if $G$ is solvable then $G = \mathbb{Z}$.
\end{enumerate}
After establishing these properties we will use them to complete the proof of the proposition. We now proceed to prove each of the clauses above.

(i). This is immediate from considering the cosets of $\langle g \rangle$ for any $g \in G - \{1_G\}$.

(ii). If $g \in G - \{1_G\}$, then $\langle g \rangle$ has finite index in the infinite group $G$. Hence $g$ has infinite order.

(iii). Let $g \in G - \{1_G\}$, and consider the action of $G$ on the left cosets of $\langle g \rangle$ by left multiplication. This action induces a homomorphism $G \rightarrow S_n$ with kernel $K$, where $n = [G : \langle g \rangle]$. So $K$ is normal and has finite index. As elements of $K$ fix the left coset $1_G \cdot \langle g \rangle$, we have $K \subseteq \langle g \rangle$. Thus, $K = \langle g^k \rangle$ for some $k > 0$.

(iv). Suppose $\{1_G\} \neq \langle g \rangle \lhd G$ and let $h \in G$. We will show $h g = g h$. Since $\langle g \rangle$ is normal, conjugation by $h$ induces an automorphism of $\langle g \rangle$. If $h g h^{-1} = g$, then there is nothing to show. Towards a contradiction, suppose $h g h^{-1} = g^{-1}$. Then $h \neq 1_G$, so $h$ has infinite order and $\langle g \rangle$ is normal of finite index, so there are nonzero $k, m \in \mathbb{Z}$ with $h^k = g^m$. We then have
$$h^{-k} g^m = 1_G = h 1_G h^{-1} = h h^{-k} g^m h^{-1} = h^{-k} h g^m h^{-1} = h^{-k} g^{-m}.$$
Thus, $g^m = g^{-m}$ with $m \neq 0$, contradicting (ii).

(v). By (ii), (iii), and (iv) we have $\Z(G) \neq \{1_G\}$. Let $H = \langle h \rangle$ be a maximal cyclic subgroup of $\Z(G)$. Towards a contradiction, suppose $H \neq \Z(G)$. Let $a \in \Z(G) - H$. Then $\langle a, h \rangle \leq \Z(G)$, so $\langle a, h \rangle$ is abelian. By (ii) $\langle a, h \rangle$ is isomorphic to a subgroup of $\mathbb{Z}^2$. However, $\mathbb{Z}^2$ has a nontrivial subgroup of infinite index, so $\langle a, h \rangle \not\cong \mathbb{Z}^2$. As $h \neq 1_G$, we must have $\langle a, h \rangle \cong \mathbb{Z}$. This contradicts the maximality of $H$.

(vi). Let $K \lhd G / \Z(G)$ be abelian, and let $\pi : G \rightarrow G / \Z(G)$ be the quotient map. Now $\Z(G) \leq \pi^{-1}(K) \lhd G$, and we want to show that $K$ is trivial. Hence, by (iv) it will suffice to show that $\pi^{-1}(K)$ is cyclic. Let $H = \langle h \rangle$ be a maximal cyclic subgroup of $\pi^{-1}(K)$ containing $\Z(G)$. Towards a contradiction, suppose $H \neq \pi^{-1}(K)$. Let $a \in \pi^{-1}(K) - H$, and let $k > 0$ be such that $\langle h^k \rangle = \Z(G)$. We have $\pi(a h a^{-1} h^{-1}) = 1_{G / \Z(G)}$ so for some $m \in \mathbb{Z}$
$$a h a^{-1} h^{-1} = h^{mk} \in \Z(G) \text{ and } a h a^{-1} = h^{mk +1}.$$
However, $h^k \in \Z(G)$, so
$$h^k = a h^k a^{-1} = (h^{mk+1})^k = h^{mk^2+k}.$$
As $h$ has infinite order, we must have $m = 0$. Thus $a h a^{-1} = h$ so $a$ and $h$ commute. Then $\mathbb{Z} \cong \langle a, h \rangle \leq \pi^{-1}(K)$, contradicting the maximality of $H$.

(vii). If $G$ is solvable then so is $G / \Z(G)$. If $G = \Z(G)$ then we are done by (v). Towards a contradiction, suppose $G \neq \Z(G)$. As $G / \Z(G)$ is solvable, its derived series terminates after finitely many steps. The last nontrivial group appearing in the derived series for $G / \Z(G)$ is normal and abelian, contradicting (vi).

Now we use the properties we have established of $G$ and complete the proof. By clause (v), $\Z(G)$ is nontrivial so $G / \Z(G)$ is a finite group. Let $\pi: G \rightarrow G / \Z(G)$ be the quotient map and let $P$ be any Sylow subgroup of $G / \Z(G)$. Every nontrivial subgroup of $N = \pi^{-1}(P)$ has finite index in $N$. However, $N$ is nilpotent, in particular solvable. Thus $N \cong \mathbb{Z}$ by clause (vii). In particular, $P = \pi(N)$ is cyclic. We now apply the following theorem of group theory (Theorem 10.1.10 of \cite{RobinsonBook}).

\begin{theorem}[H\"{o}lder, Burnside, Zassenhaus]
If $K$ is a finite group, then all of its Sylow subgroups are cyclic if and only if $K$ has a presentation
$$K = \langle a, b \: a^m = 1_K = b^n, \ b^{-1} a b = a^r \rangle$$
where $r^n \equiv 1 \mod m$, $m$ is odd, $0 \leq r < m$, and $m$ and $n(r-1)$ are coprime.
\end{theorem}

Let $a, b \in G / \Z(G)$ be as in the presentation above for $K = G / \Z(G)$. Then $\langle a \rangle$ is a normal abelian subgroup of $G / \Z(G)$. Hence by clause (vi) $\langle a \rangle$ is trivial and $a = 1_{G / \Z(G)}$. But then $\langle b \rangle = G / \Z(G)$ is a normal abelian (improper) subgroup. Again we conclude $b = 1_{G / \Z(G)}$. Thus $G / \Z(G)$ is trivial, so $G = \Z(G)$. Now $G \cong \mathbb{Z}$ by clause (v).
\end{proof}


\chapter{Further Questions}

Throughout the paper we have mentioned a number of interesting further questions that we do not have answers to. In this final chapter we collect them together and mention some general directions for further studies. The problems will be listed according to their nature, not in the order of the chapters covering them. To make this chapter a useful reference for the reader, we repeat some definitions, recall some proven facts, and include some remarks on the listed problems.

\section{Group structures}

Much of what we did in the paper was to provide a combinatorial, and in fact almost geometric, analysis of the structure of a general countable group. The results we obtained were sufficient for our purposes. But the analysis is to a large extent incomplete.

The strongest sense of geometric and combinatorial regularity of group structures is given by the notion of a ccc group. Recall that a countable group is ccc if it admits a coherent, cofinal, and centered sequence of tilings.  This notion is closely related to the study of monotileable amenable groups by Weiss. In particular, the following problem raised by Weiss was explicit in \cite{W}.

\begin{problem}[Weiss \cite{W}]
Is every countable group an MT group? That is, does every countable group admit a cofinal sequence of tilings?
\end{problem}

Weiss proved that the class of all MT groups is closed under group extensions and contains all residually finite groups and all solvable groups. In contrast we showed that the class of all ccc groups contains all residually finite groups, free products, nilpotent groups and polycyclic groups. But we only know that the class is closed under products and finite index group extensions.
We do not know of any countable group which fails to be ccc. This prompts the following questions.

\begin{problem}
Is every countable group ccc? In particular, is every solvable group ccc?
\end{problem}

Turning to a different subject, we also introduced the purely group theoretic notion of a flecc group and showed the curious property that a countably infinite group $G$ being flecc corresponds to the class of all $2$-colorings of $G$ forming a $\bf\Sigma^0_2$-complete subset of $2^G$. Recall that a countable group $G$ is flecc if there is a finite set $A\subseteq G-\{1_G\}$ such that for every non-identity $g\in G$ there is $n\in \mathbb{Z}$ and $h\in G$ with $hgh^{-1}\in A$.
 The following basic questions about flecc groups are open.

\begin{problem} Is a quotient of a flecc group flecc?
\end{problem}

\begin{problem} Is the product of two flecc groups flecc?
\end{problem}

Partial results are given in Section \ref{sec:flecc}. In particular, we showed that a normal subgroup of a flecc group is flecc. Also, the product of two flecc groups, where at least one of them is torsion, is flecc. We also completely characterized the abelian flecc groups. It turns out that this class coincides with the class of all abelian groups with the minimal condition.

\section{$2$-colorings}

In this paper we succeeded in constructing numerous group colorings with various properties. However, several construction-type problems remain open. In this section we summarize some typical ones.

We have shown that any countably infinite group $G$ admits perfectly many pairwise orthogonal $2$-colorings (this even occurs within any given open neighborhood of $2^G$). It is a trivial consequence that the same holds for $k$-colorings for any $k\geq 2$. A curious problem is to consider finite groups. Here the group $G$ and the parameter $k$ both matter in determining the maximal number of pairwise orthogonal $k$-colorings on $G$.

\begin{problem} For each finite group $G$ compute the maximal number of pairwise orthogonal $k$-colorings on $G$ for all $k\geq 2$.
\end{problem}

We obtained quite a few results involving almost equality $=^*$, especially in our proof that all near $2$-colorings are almost $2$-colorings. The fact that a $2$-coloring can be almost equal to a periodic element (which is equivalent to saying that there are groups without the ACP) is still a bit counter-intuitive and surprising. Even if we have given a complete and satisfactory characterization for all groups with the ACP, there are questions which appear to be only slightly more demanding than the ACP and which we do not know the answers to. For instance, the simplest group without the ACP is the meta-abelian group $\mathbb{Z}_2*\mathbb{Z}_2$ (which can also be expressed as a semidirect product of $\mathbb{Z}$ with $\mathbb{Z}_2$). One can directly construct a $2$-coloring $x$ on $\mathbb{Z}_2*\mathbb{Z}_2$ so that the result of turning $x(1_G)$ to $1-x(1_G)$ is a periodic element. However, we do not know if the $2$-coloring $x$ can be minimal. In general, we have not been able to construct any almost-periodic $2$-coloring that turns out to be minimal.

\begin{problem} Is there a minimal $2$-coloring on ${\mathbb Z}_2*{\mathbb Z}_2$ that is almost equal to a periodic element?
\end{problem}

In the last chapter we considered some extension problems about $2$-colorings. We showed that $\mathbb{Z}$ is the only countable group with the property that any extension of a $2$-coloring on a subgroup is a $2$-coloring on the whole group. We also completely characterized all subsets $A\subseteq G$ for which an arbitrary function $c:A\to 2$ can be extended to a $2$-coloring on $G$. In this direction the ultimate question seems to be: which partial functions on $G$ can be extended to $2$-colorings on $G$? Here we mention a necessary condition that might be sufficient.

Given a partial function $c:G\rightharpoonup 2$ define $\overline{[c]}$ as follows. Let $c^*:G\rightarrow 3$ be defined as $c^*(g)=c(g)$ if $g\in\dom(c)$ and $c^*(g)=2$ otherwise. Then let $\overline{[c]}=\overline{[c^*]}\cap 2^G$.

\begin{problem}
Given a partial function $c$ on $G$, are the following equivalent:
\begin{itemize}
\item[(i)] $c$ can be extended to $2$-colorings on $G$;
\item[(ii)] $\overline{[c]}\subseteq F(G)$?
\end{itemize}
\end{problem}

It is easy to see that (i)$\Rightarrow$(ii). So the real question is whether the converse holds.

We have seen that the set of $2$-colorings always has measure zero and is always meager. Thus, in some sense, the set of $2$-colorings is very small. However, on the other hand, we have shown that every non-empty open subset of $2^G$ contains continuum-many $2$-colorings with the closure of their orbits pairwise disjoint. This was even further strengthened in Section \ref{sec:charext}. So under certain viewpoints, the set of $2$-colorings is large. The following two questions address other notions of largeness.

\begin{problem}[Juan Souto]
For groups $G$ in which a notion of entropy exists, what is the largest possible entropy of a free subflow of $2^G$?
\end{problem}

\begin{problem}[Juan Souto]
For a given group $G$, what is the largest possible Hausdorff dimension of a free subflow of $2^G$?
\end{problem}

Finally, the fundamental method has been seen to be a tremendously useful tool for the constructive study of Bernoulli shifts. In Chapter \ref{CHAP STUDY} we developed specialized tools which work in conjunction with the fundamental method in order to produce minimal elements of $2^G$ and also pairs of points of $2^G$ whose orbit closures display some rigidity with respect to topological conjugacy. There are likely other general constructions which combine with the fundamental method in order to produce more specialized elements and subflows of $2^G$. The constructive methods in this paper would be of much more interest to the ergodic theory community if the following question were to have a positive answer.

\begin{problem}[Ralf Spatzier]
Can the fundamental method be improved in order to construct a variety subflows of $2^G$ which support ergodic probability measures?
\end{problem}

\section{Generalizations}

One of the most intriguing questions for us is: to what extent can the results of this paper be generalized? This takes many forms and can be probed in many directions. The most important direction, it seems to us, is to generalize results about Bernoulli flows to more general dynamical systems.

\begin{problem} Let $G$ be a countable group acting continuously on a Polish space $X$. Suppose there is at least one aperiodic element in $X$. Does there exist a hyper aperiodic element?
\end{problem}

We do not know the answer even when $X$ is assumed to be compact. In addition, for dynamical systems in which hyper aperiodic elements do exist, one can inquire about their density, orthogonality, etc.

Recall that $(2^\mathbb{N})^G$ is a universal Borel $G$-space. If $X$ is a compact, zero-dimensional Polish space on which $G$ acts continuously, then there is a continuous $G$-embedding (which is necessarily a homeomorphic embedding preserving $G$-actions) from $X$ into $(2^\mathbb{N})^G$. Therefore, studying hyper aperiodic elements in $(2^\mathbb{N})^G$ might be relevant to the above general problem, at least for the case when the phase space is compact and zero-dimensional.

We are fairly certain that our methods used in this paper can be used to answer many questions about the space $(2^\mathbb{N})^G$, although we have not worked out their details.

Also throughout the paper we have considered some variations of $2$-colorings whenever it is convenient. For instance, for the concept of two-sided $2$-colorings, it is trivial to note that for abelian groups they are identical to the concept of $2$-colorings. We also constructed examples of two-sided $2$-colorings for the free groups, and examples of $2$-colorings on free groups that are not two-sided. We have not attempted to systematically study this concept in conjunction with minimality, orthogonality, etc. Any such question is likely open.

Yet another direction of generalization is to consider the concept of $2$-colorings on semigroups (with the definition given by the combinatorial formulation of $2$-colorings). Other than the results we mentioned about $\mathbb{N}$ nothing is known about their general existence and properties.

Finally, the notion of hyper aperiodicity or $2$-coloring can each be generalized to the context of uncountable groups and their actions. We do not know of any interesting connection between different formulations and any significant consequence they might entail.

\section{Descriptive complexity}

The most significant application of the fundamental method in this paper is the determination of the descriptive complexity of the topological conjugacy relation for free Bernoulli subflows. Working in conjunction with the method Clemens invented in \cite{JC}, we showed that, as long as the group is not locally finite, the topological conjugacy relation for free subflows is always universal for all countable Borel equivalence relations (this result was also independently obtained by Clemens). However, understanding the complexity of this relation restricted on minimal free subflows has met significant challenges. The following very concrete problem is still open.

\begin{problem}
What is the complexity of the topological conjugacy relation for minimal free subflows of $2^{\mathbb{Z}}$?
\end{problem}

We also showed that, for locally finite groups $G$, the topological conjugacy relation for all subflows of $2^G$ is Borel bireducible with $E_0$. Moreover, the same is true when this relation is restricted on free subflows or minimal free subflows. For general groups $G$, we only know that the conjugacy relation for minimal free subflows is always at least as complex as $E_0$ in the Borel reducibility hierarchy. The general problem of determining their complexity is wide open.

\begin{problem}
For an arbitrary countable group $G$ that is not locally finite, what is the complexity of the topological conjugacy relation for minimal free subflows of $2^G$?
\end{problem}  

\backmatter
%

\bibliographystyle{amsalpha}


\printindex
\end{document}